\documentclass[10pt,reqno]{amsart}
\usepackage[T1]{fontenc}
\usepackage[margin=1in]{geometry}
\usepackage{bbm}
\usepackage{amsfonts}
\usepackage{latexsym, amssymb, amsmath, amscd, amsthm, amsxtra}
\usepackage{mathtools}
\usepackage{enumerate}
\usepackage[all]{xy}
\usepackage{mathrsfs}
\usepackage{fancyhdr}
\usepackage{listings}
\usepackage{hyperref}
\usepackage{soul}
\usepackage{xcolor}
\usepackage{dsfont}
\usepackage{stmaryrd}
\usepackage{bm}
\usepackage{comment}
\usepackage{diagbox}
\usepackage{verbatim}
\usepackage{footmisc}
\usepackage{aliascnt}

\allowdisplaybreaks[4]
\hypersetup{colorlinks=true}
\usepackage{tikz}
\usetikzlibrary{decorations.pathmorphing, decorations.markings, decorations.pathreplacing}
\usetikzlibrary{arrows.meta}

\usepackage[rightcaption]{sidecap}

\usepackage{graphicx, color}

\numberwithin{equation}{section}
\numberwithin{figure}{section}

\makeatletter
\newcommand{\rmnum}[1]{\uppercase{{\expandafter{\romannumeral #1}}}}
\makeatother

\theoremstyle{plain} \newtheorem{theorem}{Theorem}[section]
\newtheorem*{theorem*}{Theorem}
\newaliascnt{lemma}{theorem}
\newtheorem{lemma}[lemma]{Lemma}
\aliascntresetthe{lemma}
\newtheorem*{lemma*}{Lemma}
\newaliascnt{corollary}{theorem}
\newtheorem{corollary}[corollary]{Corollary}
\aliascntresetthe{corollary}
\newtheorem*{corollary*}{Corollary}
\newaliascnt{proposition}{theorem}

\aliascntresetthe{proposition}
\newtheorem*{proposition*}{Proposition}
\newaliascnt{definition}{theorem}
\newtheorem{definition}[definition]{Definition}
\aliascntresetthe{definition}
\newtheorem*{definition*}{Definition}
\newaliascnt{conjecture}{theorem}

\aliascntresetthe{conjecture}
\newtheorem*{conjecture*}{Conjecture}

\newtheorem{assumption}{Assumption}

\theoremstyle{definition}
\newaliascnt{example}{theorem}
\newtheorem{example}[example]{Example}
\aliascntresetthe{example}
\newtheorem*{example*}{Example}
\newaliascnt{remark}{theorem}
\newtheorem{remark}[remark]{Remark}
\aliascntresetthe{remark}
\newtheorem*{remark*}{Remark}
\newtheorem{claim}{Claim}

\usepackage{cleveref}
\crefname{theorem}{Theorem}{Theorems}
\Crefname{theorem}{Theorem}{Theorems}
\crefname{lemma}{Lemma}{Lemmas}
\Crefname{lemma}{Lemma}{Lemmas}
\crefname{corollary}{Corollary}{Corollaries}
\Crefname{corollary}{Corollary}{Corollaries}
\crefname{proposition}{Proposition}{Propositions}
\Crefname{proposition}{Proposition}{Propositions}
\crefname{definition}{Definition}{Definitions}
\Crefname{definition}{Definition}{Definitions}
\crefname{conjecture}{Conjecture}{Conjectures}
\Crefname{conjecture}{Conjecture}{Conjectures}
\crefname{assumption}{Assumption}{Assumptions}
\Crefname{assumption}{Assumption}{Assumptions}
\crefname{example}{Example}{Examples}
\Crefname{example}{Example}{Examples}
\crefname{remark}{Remark}{Remarks}
\Crefname{remark}{Remark}{Remarks}
\crefname{claim}{Claim}{Claims}
\Crefname{claim}{Claim}{Claims}

  \renewcommand{\r}{\mathrm}   \newcommand{\cal}{\mathcal} 
\newcommand{\scr}{\mathscr} 
 
 \newcommand{\ol}[1]{\overline{#1} \!\,} \newcommand{\wh}{\widehat}
\newcommand{\wt}{\widetilde}
\renewcommand{\txt}[1]{\text{\rm{#1}}}

\definecolor{darkred}{rgb}{0.9,0,0.3}
\definecolor{darkblue}{rgb}{0,0.3,0.9}

\newcommand{\nc}{\normalcolor}

\usepackage{ifthen}
\def\comment#1{\ifthenelse{\isodd{\value{page}}}{\marginpar{\raggedright\scriptsize{\textcolor{darkred}{#1}}}}{\marginpar{\raggedleft\scriptsize{\textcolor{darkred}{#1}}}}}

\newcommand{\bal}{{\boldsymbol{\alpha}}}
\renewcommand{\P}{\mathbb{P}}
\newcommand{\E}{\mathbb{E}}
\newcommand{\R}{\mathbb{R}}
\newcommand{\C}{\mathbb{C}}
\newcommand{\N}{\mathbb{N}}
\newcommand{\Z}{\mathbb{Z}}

\newcommand{\cG}{{\mathcal G}}
\newcommand{\cM}{{\mathcal M}}
\newcommand{\cW}{{\mathcal W}}

\newcommand{\cI}{\mathcal{I}}
\newcommand{\cC}{\mathcal{C}}
\newcommand{\cR}{{\mathcal R}}
\newcommand{\cO}{{\mathcal{O}}}
\newcommand{\cS}{{\mathcal{S}}}
\newcommand{\cB}{{\mathcal{B}}}
\newcommand{\cX}{{\mathcal{X}}}
\newcommand{\cD}{{\mathcal{D}}}
\newcommand{\cL}{{\mathcal{L}}}
\newcommand{\cK}{{\mathcal K}}
\newcommand{\cQ}{{\mathcal{Q}}}
\newcommand{\cJ}{{\mathcal{J}}}
\newcommand{\cZ}{\mathcal{Z}}
\newcommand{\cA}{\mathcal{A}}
\newcommand{\cU}{\mathcal{U}}
\newcommand{\cP}{\mathcal{P}}
\newcommand{\cE}{\mathcal{E}}

\newcommand{\fc}{\mathfrak c}

\newcommand{\fP}{\mathfrak P}

\newcommand{\fS}{\mathfrak S}

\newcommand{\sC}{\mathsf C}

\newcommand{\sE}{\mathsf E}
\newcommand{\bU}{\mathbf U}
\newcommand{\bTheta}{\boldsymbol{\Theta}}

\newcommand{\sx}{\mathsf x}

\newcommand{\e}{{\varepsilon}}
\newcommand{\fd}{\mathfrak{d}}
\newcommand{\rd}{\mathrm{d}}

\newcommand{\ii}{\mathrm{i}}
\newcommand{\dd}{\mathrm{d}}

 \renewcommand{\leq}{\le}
\renewcommand{\epsilon}{\varepsilon}

\newcommand{\contract}{\text{loop-contraction}}

\newcommand{\floor}[1] {\lfloor {#1} \rfloor}
\newcommand{\ceil}[1]  {\lceil  {#1} \rceil}

\newcommand{\qq}[1]{[\![{#1}]\!]}

\newcommand{\p}[1]{({#1})}
\newcommand{\pb}[1]{\bigl({#1}\bigr)}
\newcommand{\pB}[1]{\Bigl({#1}\Bigr)}
\newcommand{\pbb}[1]{\biggl({#1}\biggr)}

\newcommand{\pa}[1]{\left({#1}\right)}

\newcommand{\q}[1]{[{#1}]}
\newcommand{\qb}[1]{\bigl[{#1}\bigr]}
\newcommand{\qB}[1]{\Bigl[{#1}\Bigr]}
\newcommand{\qbb}[1]{\biggl[{#1}\biggr]}

\newcommand{\qa}[1]{\left[{#1}\right]}

\newcommand{\h}[1]{\{{#1}\}}

\newcommand{\ha}[1]{\left\{{#1}\right\}}

\newcommand{\abs}[1]{\lvert #1 \rvert}
\newcommand{\absb}[1]{\bigl\lvert #1 \bigr\rvert}

\newcommand{\absbb}[1]{\biggl\lvert #1 \biggr\rvert}

\newcommand{\absa}[1]{\left\lvert #1 \right\rvert}

\newcommand{\norm}[1]{\lVert #1 \rVert}
\newcommand{\normb}[1]{\bigl\lVert #1 \bigr\rVert}
\newcommand{\normB}[1]{\Bigl\lVert #1 \Bigr\rVert}

\newcommand{\norma}[1]{\left\lVert #1 \right\rVert}

\newcommand{\braket}[1]{\langle #1 \rangle}

\newcommand{\avg}[1]{\mathrm{Tr}\,( #1 )}
\newcommand{\avgb}[1]{\mathrm{Tr}\,\big( #1 \big)}

\newcommand{\avga}[1]{\mathrm{Tr}\left( #1 \right)}

\DeclareMathOperator{\re}{Re}
\DeclareMathOperator{\im}{Im}

\DeclareMathOperator{\sgn}{sgn}

\DeclareMathOperator{\spec}{spec}

\newcommand{\ba}{{\bf{a}}}
\newcommand{\bx}{{\bf{x}}}
\newcommand{\by}{{\bf{y}}}
\newcommand{\bu}{{\bf{u}}}
\newcommand{\bv}{{\bf{v}}}

\newcommand{\bp}{{p}}

\newcommand{\al}{\alpha}

\newcommand{\be}{\begin{equation}}
	\newcommand{\ee}{\end{equation}}

\newcommand{\oo}{\mathrm{o}}
\newcommand{\OO}{\mathrm{O}}

\newcommand{\HS}{\txt{HS}}

\newcommand{\bsigma}{\bm{\sigma}}

\newcommand{\ZL}{\mathbb{Z}_N}

\newcommand{\Gc}{\mathring{G}}

\newcommand{\NCT}{\txt{NCT}}

\newcommand{\bpsi}{\boldsymbol{\psi}}

\newcommand{\upcirc}[1]{\mathring{#1}}

\newcommand{\bbT}{\mathbb{T}}

\newcommand{\tf}{t_{\txt{f}}}

\newcommand{\opr}[1]{\mathrm{O}_\prec\left({#1}\right)}
\allowdisplaybreaks

\numberwithin{equation}{section}

\begin{document}

	\title{Localization Lengths of Power-Law Random Band Matrices}
	
	\author{Jiaqi Fan$^\star$}	
	\author{Fan Yang$^\dagger$} 
	\author{Jun Yin$^\S$}

	\thanks{\hspace{-9.3pt}$^\star$Qiuzhen College, Tsinghua University, Beijing, China, \href{mailto:fanjq24@mails.tsinghua.edu.cn}{fanjq24@mails.tsinghua.edu.cn}}
	\thanks{$^\dagger$Yau Mathematical Sciences Center, Tsinghua University, and Beijing Institute of Mathematical Sciences and Applications, Beijing, China,
		\href{mailto:fyangmath@mail.tsinghua.edu.cn}{fyangmath@mail.tsinghua.edu.cn}}
	\thanks{$^\S$Department of Mathematics, University of California, Los Angeles, Los Angeles, CA, USA,
		\href{mailto:jyin@math.ucla.edu}{jyin@math.ucla.edu}}
	
	\begin{abstract}
		We study large $N\times N$ power-law random band matrices $H=(H_{ij})$ with centered complex Gaussian entries, where the variances satisfy a power-law decay $\mathbb{E}|H_{ij}|^2\propto (|i-j|/W+1)^{-1-\alpha}$, for some exponent $\alpha>-1$ and bandwidth $W\gg 1$. We establish the following lower bounds, with high probability, on the localization length $\ell$ of bulk eigenvectors in the different regimes of $\alpha$: (1) $\ell=N$ if $-1<\alpha<0$; (2) $\ell \ge W^{C}$ for any large constant $C>0$ if $0 \le \alpha \le 1$; (3) $\ell \ge W^{\alpha/(\alpha-1)}$ if $1 < \al <2$; (4) $\ell \ge W^{2}$ if $ \alpha \ge 2$. These results verify the physical conjecture of \cite{Power-RBM} on the delocalized side. The main difficulty in the proof lies in handling the interplay between the non-mean-field nature of the model and the slow decay of the variance profile. To address this issue, a key technical ingredient is a new dynamical analysis of $T$-variables formed from pairs of resolvent entries of $H$. In contrast to the fundamental works on regular random band matrices with fast-decaying variances \cite{Band1D,erdos2025zigzagstrategyrandomband}, this approach does not rely on higher-order resolvent loops.
		
	\end{abstract}
	\maketitle
	
	\tableofcontents
	
	\section{Introduction}

	The study of the celebrated Anderson localization-delocalization phenomenon originates from the seminal work of Anderson in 1958 \cite{Anderson1958Absence}. To theoretically understand electron transport and metal-insulator transitions in disordered crystals, Anderson introduced the discrete random Schr\"odinger operator, now widely known as Anderson's tight-binding model. On the $d$-dimensional integer lattice $\mathbb{Z}^d$, the Hamiltonian of the Anderson model is typically defined as $H = -\Delta + \lambda V$, where $\Delta$ is the discrete Laplacian, $V = \{v_x\}_{x\in\mathbb{Z}^d}$ represents a disordered potential field consisting of i.i.d.~random variables, and $\lambda$ denotes the disorder strength. A central quantity of interest is the \emph{localization length}, denoted by $\ell$. Roughly speaking, the localization length characterizes the spatial decay scale of the eigenvectors: if a wavefunction $\psi$ exhibits exponential decay, for instance $\psi(x) \sim \exp(-|x|/\ell)$, it is said to be localized; conversely, if the wavefunction extends over the entire system (conceptually corresponding to $\ell = \infty$), it is considered extended. In the physics literature, a combination of theoretical approaches---most notably the one-parameter scaling theory \cite{PRL_Anderson}---has led to detailed predictions concerning the behavior of the localization length across different dimensions. In low dimensions ($d\in\{1, 2\}$), it is widely believed that arbitrarily weak disorder $\lambda > 0$ leads to localization of all eigenstates. Specifically, in the weak disorder regime $\lambda \to 0$, the localization length is predicted to scale as $\ell \propto \lambda^{-2}$ in $d=1$, and to diverge exponentially as $\log \ell \propto \lambda^{-2}$ in $d=2$. In higher dimensions $d \ge 3$, the system is predicted to undergo an Anderson transition at a critical disorder strength $\lambda_c$: states are localized at the spectral edges or under strong disorder $\lambda > \lambda_c$, and extended in the bulk of the spectrum under weak disorder $\lambda < \lambda_c$. For a more detailed overview and further references, we refer the reader to \cite{mott1961theory,ishii1973localization,Borland1963TheNature,Aizenman_book,Kirsch2007}.

	Beyond the Anderson model with diagonal random potentials, another important class of models for studying the Anderson localization-delocalization phenomenon is provided by the celebrated random band matrix (RBM) model \cite{scalingabndCGMLIF1990PRL,ConJ-Ref2,ScalingPropertyBandMatrixFYMA1991PRL}. RBMs offer a natural interpolation between mean-field models, such as Wigner matrices \cite{Wigner}, and random Schr\"odinger operators. Roughly speaking, the off-diagonal elements $H_{xy}$ of an RBM $H$ are independent centered random variables (up to the Hermitian symmetry), whose variances $S_{xy}:=\mathbb{E}|H_{xy}|^2$ depend on the distance $|x-y|$ and a parameter $W$, known as the bandwidth. Usually, it is assumed that the variance decays exponentially or has compact support when $|x-y| \gg W$. Based on the nonlinear $\sigma$-model \cite{ScalingPropertyBandMatrixFYMA1991PRL} and numerical simulations \cite{scalingabndCGMLIF1990PRL,ConJ-Ref2,PhysRevLett.66.986,MWilkinson_1991}, physicists predict that the localization length $\ell$ of an RBM on the discrete torus $\ZL^d:=\Z^d/(N\Z)^d$ depends strongly on the dimension $d$ and the bandwidth $W$: in $d=1$, the localization length is expected to scale as $\ell \propto W^2\wedge N$; in $d=2$, $\ell \propto  [\exp\pa{cW^2}] \wedge N$; in $d \ge 3$, a delocalization transition occurs when $W$ exceeds a critical constant $W_{c}$.

	Although short-range models such as the standard Anderson model and short-range RBMs---which we will also refer to as regular RBMs---have been extensively studied in both the physics and mathematics literature, a closer look at Anderson's seminal work \cite{Anderson1958Absence} reveals an interesting nuance. In his original paper, Anderson considered a slightly more general model in dimension $d=3$ in which the discrete Laplacian $-\Delta$ is replaced by a more general hopping operator $\Gamma$,\footnote{It is denoted by $V$ in \cite{Anderson1958Absence}. To avoid confusion with the tight-binding model, we use a different notation here.} whose matrix elements decay faster than $r^{-d}$. This framework includes the nearest-neighbor Laplacian $-\Delta$ as a special case, but also allows for long-range power-law interactions.
	Motivated by the study of such long-range effects and regular RBMs, Mirlin et al.~\cite{Power-RBM} later introduced the \emph{power-law random band matrix} (PRBM) model as a paradigm for investigating Anderson localization-delocalization transitions in systems with long-range hopping. More precisely, a PRBM $H=(H_{xy})$ is an $N\times N$ random Hermitian matrix indexed by the discrete torus $\ZL=\mathbb{Z}/(N\mathbb{Z})$, whose entries are independent centered Gaussian random variables (up to the Hermitian symmetry). Given a bandwidth parameter $W\gg1$, the variances $S_{xy}:=\mathbb{E}|H_{xy}|^2$ decay polynomially when $|x-y|_N\gg W$ (where $|\cdot|_N$ denotes the periodic distance on the discrete circle \smash{$\ZL$}), according to a fixed exponent $\alpha\in(-1,\infty)$:\footnote{Our exponent $\alpha$ corresponds to the exponent $2\alpha-1$ in \cite{Power-RBM}. We choose the decay exponent $\alpha+1$ so that it matches the power-law decay exponent of the $\alpha$-stable density.}
	\begin{equation}\label{eq:decay_exponent}
		S_{xy} \propto \pa{{|x-y|_N}/{W} + 1}^{-1-\alpha}. \end{equation}
	One can normalize the variances so that $\sum_y S_{xy}\equiv1$ for all $x$.

	By mapping this model to a nonlinear $\sigma$-model with nonlocal interaction, Mirlin et al.~\cite{Power-RBM} derived several physical predictions regarding the behavior of PRBMs in different regimes of the exponent $\alpha$:
	\begin{itemize}
		\item[(C1)] For $\alpha<0$, the model behaves like a mean-field Wigner ensemble.
		
		\item[(C2)] For $\alpha>2$, the model behaves similarly to a one-dimensional (1D) regular RBM.
		
		\item[(C3)] For $0<\alpha<2$, the model effectively behaves like a regular RBM in dimension $2/\alpha$. In particular, a phase transition is expected at $\alpha=1$, corresponding to the critical dimension $2$ for regular RBMs.
	\end{itemize}
	Based on these predictions, the characteristic spatial scale of wave-packet propagation for $H$ at time $t$ is conjectured to be \smash{$\ell(t)=Wt^{\frac{1}{\alpha\wedge 2}}$} for $\alpha>0$. In particular, the dynamics are expected to be diffusive for $\alpha\ge2$, superdiffusive for $\alpha\in(1,2)$, and for $\alpha\le1$ the spreading scale becomes \smash{$\ell(t)=W(t^{\frac{1}{2}})^{\frac{2}{\alpha}}$}, which coincides with the characteristic diffusive scale of a regular RBM in effective dimension $2/\alpha$. 
	Combining these heuristics with the Thouless scaling theory, one can predict the occurrence of localization or delocalization by comparing $\ell(N)$---the propagation scale at time $N$ (where $N$ corresponds to the reciprocal of the typical eigenvalue spacing in the spectral bulk)---with the system size $N$. This leads to the following conjecture for the localization length $\ell_{\alpha,W,N}$:
	\begin{equation}\label{definition_ell_P}
		\ell_{\al,W,N}\asymp \begin{cases}
			N,\ \ &\txt{if}\ \ \alpha\in (-1,1),\\
			W^{\alpha/(\alpha-1)}\wedge N,\ \ &\txt{if}\ \ \alpha\in (1,2),\\
			W^2\wedge N,\ \ &\txt{if}\ \ \alpha\in (2,\infty).
		\end{cases}
	\end{equation}
	Owing to its rich phase diagram and intriguing physical behavior, the PRBM model has attracted considerable attention in the physics community. Extensive studies have been devoted to the model itself \cite{Cuevas2001Anomalously,varga2000critical,evers2000fluctuations,mirlin2000multifractality,cuevas2001fluctuations,varga2002fluctuation,yevtushenko2003virial,PhysRevLett.79.1913,ndawana2003energy,PhysRevB.72.064108,PhysRevE.73.026213,paley2005statistical,mendez2006scattering,Mildenberger2006BoundaryMI,monthus2009statistical,mendez2010scattering,mendez2014generalized,bogomolny2018power,carrera2021multifractal,martinez2023scattering}, as well as its applications in quantum chaos \cite{ALTSHULER1997487,BORGONOVI1999317,PhysRevE.104.064202,hopjan2023scale,lydzba2024normal,santra2025complexitytransitionschaoticquantum}, and in various other contexts \cite{Pono1997Coherent,Kravtsov2000Energy,varga2008entanglement,monthus2009statistics,PhysRevB.98.134205,thomson2020dynamics,lydzba2020eigenstate,rao2022power,zhang2024magnetic,cohen2024complexity,buijsman2026power,zarateherrada2025dynamicaldetectionextendednonergodic,faez2011critical,martinez2023coherent,vallejo2024reducing,sen2025dissipativespectralformfactor}. 
    The PRBM model also has a hierarchical analogue known as the \emph{ultrametric ensemble}, and the connection between these two models has been discussed in the physics literature \cite{Fyodorov_2009,bogomolny2018power}. The localization-delocalization transition in the ultrametric ensemble has also been studied mathematically \cite{10.1214/18-EJP197,vonSoostenWarzel2019}. Furthermore, other long-range models, such as long-range Anderson models, have also attracted significant interest in recent years \cite{PhysRevB.57.10232,C.Yeung_1987,PhysRevLett.64.547,L.S.Levitov_1989,levitov1999critical,GARCIAGARCIA2004361,vega2019multifractality,CLG2023critical}. We refer the reader to \cite[Section III]{evers2008anderson} for a detailed review of PRBMs in the physics literature.

	\subsection{Overview of the main results}
	
	The main objective of this work is to provide a mathematically rigorous justification of the conjecture \eqref{definition_ell_P} from the delocalized side by establishing (almost) sharp lower bounds on the localization length $\ell_{\alpha,W,N}$ over the full non-mean-field regime $\alpha\in(-1,\infty)$. We consider a slightly more general setting than \eqref{eq:decay_exponent}, assuming only the following upper bound on the variances:
	\begin{equation*}
		0 \leq S_{xy} \lesssim \cS_{xy},\qquad \txt{with}\quad\cS_{xy} := \frac{1}{\cZ_{\alpha}} \bigg( \frac{|x-y|_N}{W} + 1 \bigg)^{-1-\alpha},
	\end{equation*}
	where $\mathcal{Z}_\alpha$ is the normalization constant of $\mathcal{S}$. In addition, we impose certain regularity conditions on the resolvent of the variance matrix $S=(S_{xy})$, as detailed in Assumptions \ref{assumption_input_bound} and \ref{assumption_input_bound<0} below.
	These assumptions can be verified by a delicate Fourier analysis of the resolvent of $S$ in the case where $S$ is translation invariant and induced by a profile function whose Fourier transform satisfies suitable regularity conditions; see Example \ref{example_profile_function}. In particular, this framework includes the 1D regular RBMs for any large exponent $\alpha>2$. For definiteness, we define the \emph{localization length} of the eigenvectors $\{\bpsi_k\}$ of $H$ by
	\be\label{eq:localization_length}
    \ell(\bpsi_k):= \inf\Big\{0\le \ell\le N: \max_{x_0\in \ZL} \sum_{|x-x_0|\le \ell}|\bpsi_k(x)|^2 \ge 1/2 \Big\}.
	\ee 
	We note that it is equally natural to replace the constant $1/2$ in \eqref{eq:localization_length} by any fixed value $a\in(0,1)$. For the results obtained in this work, all such choices lead to localization lengths of the same order.\footnote{However, we note that, unlike the Anderson model or short-range RBMs with exponentially decaying variance profiles, the eigenvectors of PRBMs are no longer exponentially localized. Instead, their entries are expected to exhibit power-law decay, as suggested by the physics literature.}

	Under the above assumptions and the definition \eqref{eq:localization_length}, and assuming $W\ge N^{\varepsilon}$ for an arbitrarily small constant $\varepsilon>0$,\footnote{This condition can be relaxed to $W\ge1$ in the regime $\alpha\in(-1,0)$.} we prove in Theorem \ref{theorem_delocalization} that \eqref{definition_ell_P} holds as a lower bound for the localization length of the bulk eigenvectors. The proof follows from an optimal local law for the Green's function $G(z)=(H-z)^{-1}$ for spectral parameters $z=E+\mathrm{i}\eta$ with $\eta>0$ down to the optimal scale $\eta\gg \ell_{\alpha,W,N}^{-1}$ (Theorem \ref{theorem_local_law_entrywise}). Moreover, we establish several characteristic signatures of the bulk eigenvectors and eigenvalues in the delocalized phase (i.e., when $\ell_{\alpha,W,N}=N$), including quantum unique ergodicity (QUE) of eigenvectors (Theorem \ref{theorem_que} and Corollary \ref{cor_que}) and the universality of local bulk spectral statistics (Theorem \ref{theorem_local_statistics}). These results are obtained as consequences of quantum diffusion estimates for the Green's function (Theorem \ref{theorem_local_law_loops}).

	We remark that our results can be extended to power-law random band matrices in higher dimensions $d\ge 2$. More precisely, for an exponent $\alpha\in(-d,\infty)$, consider a $d$-dimensional PRBM $H$ with variances $S_{xy} \propto \p{{|x-y|_N}/{W} + 1}^{-d-\alpha}$ for $x,y \in \Z^d/(N\Z)^d$. We conjecture that this model exhibits three distinct regimes depending on whether $\alpha<0$, $0<\alpha<d$, or $\alpha>d$. In particular, when $\alpha>d$, the model is expected to behave similarly to a $d$-dimensional regular RBM. For $\alpha\in(-d,d)$, the corresponding localization length $\ell_{\alpha,d,W,N}$ is conjectured to be of order $N$ for all $d\ge2$.
	We believe that our methods can be extended to this higher-dimensional setting to show that $\ell_{\alpha,d,W,N}\gtrsim N$ for all $d\ge2$ and $\alpha\in(-d,\infty)$ under the assumption $W\ge N^{\varepsilon}$. However, for simplicity of presentation, and noting that the phase diagram in higher dimensions is \emph{less rich} than that of the 1D PRBM, we do not pursue this direction in the present paper.

	We now briefly review some related work, beginning with results on regular RBMs. 
	There has been a series of breakthroughs regarding the delocalization of regular RBMs, and many of the techniques and ideas developed in these works provide a fundamental framework for the present paper. In one dimension, delocalization was established progressively under the assumption $W\ge N^{a}$ for various exponents $a>1/2$ in a sequence of works \cite{erdHos2011quantum,erdHos2011quantum1,delocal,BaoErd2015,HeMa2018,bourgade2017universality,bourgade2019random,bourgade2020random,Band1D_III,DY}, and was finally proved under the sharp condition $W\gg N^{1/2}$ in \cite{Band1D} using a dynamical approach known as the \emph{analysis of loop hierarchy}. With further technical refinements, this strategy was subsequently extended to 2D (two-dimensional) RBMs in \cite{Band2D}, establishing delocalization under the assumption $W\ge N^{\varepsilon}$ for an arbitrarily small constant $\varepsilon>0$, and later to the spectral edge in \cite{Bandedge}. Delocalization was also extended to 1D RBMs with general variance profiles and non-Gaussian entries in \cite{erdos2025zigzagstrategyrandomband} (see the next paragraph for more details).
	In higher dimensions, a diagrammatic expansion method was developed to establish delocalization for dimensions $d\ge7$ under the assumption $W\ge N^{\varepsilon}$ \cite{yang2021delocalization,YYYTexpansion2022CMP,xu2024bulk}. More recently, in \cite{dubova2025delocalizationnonmeanfieldrandommatrices}, the dynamical method and the diagrammatic expansion were combined to extend delocalization results for regular RBMs to all dimensions $d\ge3$ under the condition $W\ge N^{\varepsilon}$. On the localized side, localization has been established for 1D regular RBMs under the condition $W\ll N^{a}$ for various exponents $a<1/2$ \cite{Sch2009,PelSchShaSod,Cipolloni2024,Chen2022}, and under the sharp condition $W\ll \sqrt{N}$ in \cite{Localization1_2}. These results complete the phase diagram for the localization-delocalization transition in 1D RBMs.

	Compared with the extensive studies of regular RBMs in all dimensions in the mathematics community and the widespread interest in power-law RBMs in the physics literature, rigorous results on power-law RBMs remain relatively scarce in the mathematical literature. For PRBMs with relatively fast decay (i.e., for large decay exponent $\alpha$), Erd\H{o}s and Riabov \cite{erdos2025zigzagstrategyrandomband} recently established delocalization and QUE of bulk eigenvectors, as well as the universality of local bulk spectral statistics, under the sharp condition $W \gg N^{1/2}$.
	Strictly speaking, the results for PRBMs in \cite{erdos2025zigzagstrategyrandomband} were explicitly stated and proved for $\alpha\ge 7$. However, by carefully comparing the admissibility conditions in Definition 2.4 therein with our Assumption \ref{assumption_input_bound} below, we believe that their arguments can be extended directly to all $\alpha\ge3$. When $\alpha<3$, however, their assumption (2.25) inevitably breaks down, thereby invalidating some key steps of their proof.
	On the other hand, near the spectral edges, Liu and Zou \cite{liu2025edgestatisticsrandomband,liu2025edgeuniversalityinhomogeneousrandom} generalized the moment method pioneered by Sodin \cite{Sodin2010Spectral} to establish a phase transition in the extreme eigenvalue statistics for 1D PRBMs when $\alpha>1$, together with partial results in higher dimensions. Notably, none of these recent developments addresses the intermediate regime $\alpha\in(0,1)$, in which the PRBM exhibits behavior analogous to that of high-dimensional regular RBMs.

	\subsection{New ideas} \label{subsec:ideas}

	In the following discussion, it is helpful to keep the conjectures (C1)--(C3) in mind. Although these conjectures concerning the comparison with regular RBMs are heuristic in nature, they provide a useful criterion for assessing the relative difficulty of different regimes of $\alpha$. Moreover, they serve as a guiding principle for extracting insights from corresponding studies of regular RBMs in analogous dimensions.
	
	We begin by introducing some basic notation. Let $G(z):=(H-z)^{-1}$ denote the Green's function (or resolvent) of $H$, where $z\in\mathbb{C}$ with $\eta=\im z>0$. Our analysis in all regimes is based on the classical observation that the off-diagonal entries $|G_{xy}|^2$ are stochastically dominated by the so-called $T$-variables and $\mathcal{L}$-loops (see, e.g., Lemma \ref{2_loop_to_1_chain}), defined by
	\begin{equation}\label{eq:TcL}
		T_{xy}=\sum_a S_{xa}|G_{a y}|^2,\quad \cL_{xy}=\sum_{a,b} S_{xa}|G_{a b}|^2S_{b y}.
	\end{equation}
	These quantities are expected to exhibit a \emph{quantum diffusion} behavior, in the sense that they are well approximated by a classical diffusion profile given by the resolvent of the Markov chain on $\ZL$ with transition matrix $S$, namely
	\(\Theta_{t}:=S/(1-t S),\) 
	which we refer to as the \emph{$\Theta$-propagator}. One may think of $t\in[0,1)$ as a time parameter satisfying $1-t\asymp \eta$. More precisely, we adopt the dynamical framework introduced in \cite{DY,Band1D}, in which the random matrix evolves according to a matrix-valued Brownian motion:
	\[
	\dd (H_{t})_{xy} = \sqrt{S_{xy}} \, \dd (\mathbf{B}_{t})_{xy}, \quad \text{with}\quad H_0 = 0,
	\]
	where \( (\mathbf{B}_{t})_{xy} \) are independent standard complex Brownian motions for \( x, y \in  \ZL\), subject to the Hermitian symmetry constraint. Following \cite{10.1214/19-ECP278, Sooster2019, DY}, we consider the Green's function flow \(G_t := (H_t - z_t)^{-1},\) with a suitably chosen time-dependent spectral parameter $z_t$ (see \eqref{eq:zt}), whose dynamics are naturally renormalized at leading order.

	We first consider the almost mean-field regime $\alpha\in(-1,0)$. In this regime, iterating the $T$-variable representation yields that $|G_{xy}|^2$ can be bounded by $\sum_{\alpha} (S^k)_{xa} |G_{a y}|^2$ for any fixed $k\ge1$. Since there exists a fixed $k$ (depending on $\alpha$) such that $\max_{x,y} S^k_{xy} = O(N^{-1})$, we can bound $|G_{xy}|^2$ with $N^{-1}\sum_\al |G_{\al y}|^2=\im G_{yy}/N\eta$ by Ward's identity. Thus, the model exhibits strong mean-field behavior, and techniques developed for Wigner matrices apply in this regime.
	In the critical case $\alpha=0$, the model is conjectured to behave similarly to an infinite-dimensional regular RBM. In this case, we adapt ideas from high-dimensional RBMs \cite{YYYTexpansion2022CMP} to control high moments of the $T$-variables via a $T$-expansion. The remaining challenging regimes are:
	\begin{itemize}
		\item $\alpha\in[1,\infty)$, which is expected to exhibit behavior analogous to regular RBMs in dimensions $1\le d\le2$ (with $\alpha=1$ corresponding to the critical dimension $d=2$). Techniques developed for 1D and 2D RBMs \cite{Band1D,Band2D,erdos2025zigzagstrategyrandomband} are therefore relevant.

		\item $\alpha\in(0,1)$, which is expected to behave similarly to regular RBMs in dimensions $d>2$, suggesting the relevance of techniques for higher-dimensional RBMs \cite{dubova2025delocalizationnonmeanfieldrandommatrices,yang2021delocalization}.
		
	\end{itemize}
	We now discuss the main difficulties in these regimes and the new ideas developed to overcome them. A key reason why the above heuristics are effective is that the short-range decay and local regularity properties of the $\Theta$-propagator for power-law decaying $S$ resemble those in short-range models of the corresponding dimensions. However, a fundamental challenge in the power-law setting is that the long-range decay of the $\Theta$-propagator---being polynomial with a fixed exponent---is significantly slower than the exponential decay observed in regular RBMs. This is the primary source of the technical difficulties encountered in our analysis.

	\medskip
	\noindent
	{\bf The regime $\alpha \in [1, \infty)$}.
	In this regime, our model exhibits properties analogous to those of 1D and 2D regular RBMs. In particular, for $t\in[0,1)$, the upper bound on the $\Theta$-propagator takes a form closely resembling its counterparts in the 1D and 2D regular RBM settings:
	\begin{equation}
		\Theta_t(x,y)\lesssim \q{\ell_t(1-t)}^{-1}\q{|x-y|_N/\ell_t+1}^{-1-\alpha},\quad \text{with}\quad \ell_t:=\min\q{W(1-t)^{-\frac{1}{\alpha\wedge 2}},N}.\label{eq_slow_theta} 
	\end{equation}
	The key distinction is that, when $|x-y|_N$ exceeds the \emph{characteristic length} $\ell_t$, the exponential decay observed in regular RBMs is replaced here by a much slower polynomial decay.
	
	Let $G_t^+\equiv G_t$, $G_t^-\equiv G_t^*$, and for each $x\in \ZL$, let $S^{(x)}$ denote the diagonal matrix with entries \smash{$S^{(x)}_{ab} =S_{xa}\delta_{ab}$}. We define the resolvent loops $\cL^{(n)}$ and resolvent chains $\mathcal C^{(n)}$ of order (or length) $n\in \N$ as:
	\be\label{eq:nGloops}
	{\cal L}^{(n)}_{t,\bsigma,\bx}:= \mathrm{Tr} \prod_{i=1}^n \p{G_t^{\sigma_i} S^{(x_i)}},\quad {\cal C}^{(n)}_{t,\bsigma,\bx'} (x,y):=  \p{G^{\sigma_1}_tS^{\pa{x_1}}G^{\sigma_2}_t\cdots S^{\pa{x_{n-1}}}G^{\sigma_n}_t}_{xy}, 
	\ee
	where $\bsigma=(\sigma_1,\ldots,\sigma_n)\in\{+,-\}^n$, $\bx=(x_1,\ldots, x_n)\in \ZL^n,$ and $ \bx'=(x_1,\ldots, x_{n-1})\in \ZL^{n-1}.$ The analysis of the hierarchy of such resolvent loops is referred to as \emph{analysis of the loop hierarchy} in \cite{Band1D}, and constitutes a cornerstone of a series of recent works on regular RBMs, including \cite{Band1D,Band2D,erdos2025zigzagstrategyrandomband,Bandedge,fan2025blockreductionmethodrandom,dubova2025delocalizationnonmeanfieldrandommatrices}. In this framework, the evolution of loops of length $n$ depends on estimates for loops of lengths $1,2,\ldots,n+1$, and $2n+2$. In particular, the analysis of short loops requires control over loops of significantly larger lengths, reflecting a structure analogous to the classical BBGKY hierarchy in many-body dynamics. This difficulty is resolved in the foundational works \cite{Band1D, erdos2025zigzagstrategyrandomband} via related but distinct approaches:
	\begin{itemize}
		\item In \cite{Band1D}, each $(2n+2)$-loop arising in the evolution equation for $n$-loops is decomposed into two shorter loops of lengths $2\lceil(n+1)/2\rceil$ and $2\lfloor(n+1)/2\rfloor$, respectively. These, together with the other loop contributions from the equation, can be bounded directly, yielding improved bootstrap estimates for the $n$-loops (i.e., estimates sharper than the a priori bounds). Iterating this procedure finitely many times produces sufficiently strong bounds on the maximum norms of loops of arbitrary length.

		\item In \cite{erdos2025zigzagstrategyrandomband}, a stopping-time argument is developed for loops and chains of length up to a large fixed even integer \(K \in \mathbb{N}\), with precise decay estimates tracked. During the analysis, all loops and chains of length exceeding \(K\) are decomposed into objects of length at most \(K\), at the cost of controllable small factors $W^\e$.
	\end{itemize}
	
	However, both approaches break down in the present power-law setting, particularly in the regime $\al\in[1,3)$. In the argument of \cite{Band1D}, the iterative scheme relies crucially on sufficiently strong initial estimates. These a priori bounds are obtained via certain \emph{continuity estimates} formulated solely in the max-norm. Such control is adequate in short-range models due to the rapid exponential decay  $\exp(-c|x-y|_N/\ell_t)$ of the $\Theta$-propagator, which constitutes a key convenience in their analysis. In particular, this allows one to track loop bounds entirely through their max-norms, together with the observation that resolvent entries $(G_t)_{xy}$ become negligible whenever $|x-y|_N\ge W^\e\ell_t$, for arbitrarily small $\e>0$. Both properties are stable under the flow and persist throughout the evolution.
	In contrast, this simplification is no longer available in the power-law setting, due to the slow decay $\q{|x-y|_N/\ell_t+1}^{-1-\alpha}$ in \eqref{eq_slow_theta}. As a result, the scale $\ell_t$ no longer provides an effective cutoff for the decay of resolvent entries, and it becomes necessary to track their precise spatial decay. However, once the polynomial decay factors are divided, the a priori bounds derived from the continuity estimates deteriorate significantly, thereby invalidating the bootstrap scheme of \cite{Band1D}.

	The need to control the precise decay of resolvent loops has already been recognized in \cite{erdos2025zigzagstrategyrandomband}, where such estimates are successfully established. Nevertheless, their method encounters additional difficulties in the regime $\al\in[1,3)$. More specifically, the losses of factors $W^\e$, arising from cutting long loops and chains, necessitate sharper propagation estimates than those required in \cite{Band1D}. In \cite{erdos2025zigzagstrategyrandomband}, these refined estimates are achieved through the so-called \emph{two regularizations} (see Lemma 7.4 therein), which involve taking finite differences of two $\Theta$-propagators. By comparison, the \emph{sum-zero operator} introduced in \cite[Definition 5.12]{Band1D} acts on a single propagator and may thus be viewed as a \emph{one-regularization} in the terminology of \cite{erdos2025zigzagstrategyrandomband}. However, the proof of the two-regularization estimates relies on assumption (2.25) in \cite{erdos2025zigzagstrategyrandomband}, which in our setting is only valid for $\al\ge 3$, as mentioned above.

	In this work, our strategy for the regime $\alpha\in[1,3)$ is to combine the framework of \cite{Band1D} with an extensive use of the so-called \emph{loop-contraction inequalities}, introduced in \cite{dubova2025delocalizationnonmeanfieldrandommatrices} for the study of regular RBMs in dimensions $d\ge 3$.\footnote{We expect that the method of \cite{erdos2025zigzagstrategyrandomband} can also be extended to the regime $\alpha\in[1,3)$ with additional technical effort. Our choice to adopt the approach of \cite{Band1D,dubova2025delocalizationnonmeanfieldrandommatrices} is intended to demonstrate that the techniques developed therein---originally for regular RBMs with a specific block variance profile and without tracking the precise decay of higher-order resolvent loops (of length $n\ge 3$)---can be adapted to models with general variance profiles and refined to capture precise decay estimates for resolvent loops of length $\ge 3$. The ideas developed here may prove useful in future studies of non-mean-field random matrices.}
	Broadly speaking, loop-contraction inequalities form a class of techniques that exploit Ward's identities to control partial sums of long loops over subsets of indices in terms of shorter loops. Compared with the loop-cutting arguments in \cite{Band1D, erdos2025zigzagstrategyrandomband}, this approach has the key advantage of yielding an additional factor $\q{\ell_t(1-t)}^{-1/2}$, which arises from the gain of an extra off-diagonal resolvent entry in the contraction procedure. This gain is crucial for our analysis.     
	With this strategy in mind, adapting the framework of \cite{Band1D} to our setting reduces to two main tasks: (i) establishing sufficiently strong (albeit non-optimal) a priori bounds to initiate the bootstrap argument, and (ii) developing a sequence of bootstrap estimates that strictly improve upon these initial bounds. 
	To this end, we first prove a weak local law for resolvent loops of length 2, which yields sharp upper bounds for these loops and, consequently, sharp decay estimates for off-diagonal resolvent entries, namely $|(G_{t})_{xy}|\lesssim \sqrt{\Theta_t(x,y)}$ with high probability, where $\Theta_t(x,y)$ satisfies \eqref{eq_slow_theta}. The proof of this weak local law relies on a stopping-time argument applied simultaneously to resolvent loops of lengths 2,3,4, thereby involving loops of length up to 10. To truncate this hierarchy and close the argument, we reduce the length of long loops via loop-contraction inequalities. 

    Equipped with the sharp resolvent estimates, we derive precise decay bounds for longer loops in terms of the long-range decay factors $\q{|x-y|_N/\ell_t+1}^{-(1+\alpha)/2}$ appearing in $\sqrt{\Theta_t(x,y)}$. However, this procedure does not yield sufficiently many factors of \smash{$\q{\ell_t(1-t)}^{-1/2}$} to initiate the bootstrap procedure. As a consequence, the argument of \cite{Band1D} breaks down, in particular in the treatment of the $(2n+2)$-loops arising in the martingale estimates. At this stage, loop-contraction inequalities once again play a crucial role: they provide the missing factor \smash{$\q{\ell_t(1-t)}^{-1/2}$}, thereby compensating for the loss in the naive loop estimates and yielding bootstrap bounds strong enough for iteration. 
    Furthermore, unlike \cite{Band1D}, where the block variance profile of the model enables a self-consistent analysis of resolvent loops, the general variance profile considered in this work also requires control of the resolvent chains arising from the evolution equations for the loops. In \cite{erdos2025zigzagstrategyrandomband}, this difficulty is addressed by incorporating the dynamics of the resolvent chains, thereby making the coupled system of evolution equations for loops and chains self-consistent. Inspired by \cite[Lemma A.2]{Band1D}, we instead develop a static method that combines resolvent identities with standard large-deviation bounds to obtain sharp estimates for chains in terms of bounds on loops; see \Cref{lemma_Xi_reshape}. This approach allows us to decouple the loop dynamics from the chain dynamics and thereby restore self-consistency in the dynamical analysis of the loops. For further details, we refer to Sections \ref{sec_step_3_supercritical}--\ref{sec_step_5_supercritical}.

	Before concluding the discussion of the regime $\alpha\in[1,\infty)$, we highlight another technical challenge that arises in our analysis, namely, the derivation of sufficiently strong upper bounds for the so-called $\cK$-loops, which represent the deterministic limits of resolvent loops. Such bounds have previously been established in \cite{Band1D,Band2D,truong2025localizationlengthfinitevolumerandom,dubova2025delocalizationnonmeanfieldrandommatrices} via the static approach based on the \emph{tree-representation formula}, and in \cite{erdos2025zigzagstrategyrandomband,fan2025blockreductionmethodrandom} through two distinct dynamical approaches.
    Motivated by the decay properties of resolvent entries, one expects that sharp upper bounds on resolvent loops exhibit the so-called \emph{loop decay} (see \Cref{def_decay_factors}), in which each edge of the loop contributes a factor \smash{$\sqrt{\Theta_t(\cdot,\cdot)}$}. This type of decay also appears in the estimates obtained in \cite{erdos2025zigzagstrategyrandomband}. It is therefore natural to conjecture that this decay is likewise optimal for the $\cK$-loops. 
    However, from a technical perspective, this decay structure leads to difficulties in the regime $\al\in[1,3)$. The underlying issue is the absence of the two-regularization mechanism mentioned above, which plays a crucial role in both the tree-representation approach and the dynamical methods. More precisely, the long-range decay factor between adjacent indices, $(|\cdot|_N/\ell_t+1)^{-(1+\alpha)/2}$, is insufficient to close the required estimates.
	Guided by the tree-representation formula for $\cK$-loops, we instead establish an upper bound---conjecturally optimal---with a stronger decay structure, which we term the \emph{tree-shaped decay}. In this formulation, each pair of adjacent indices contributes a full decay factor $(|\cdot|_N/\ell_t+1)^{-1-\alpha}$ (see \Cref{def_tree_shaped_decay} and \Cref{lemma_n_cK_tree_bound}). Crucially, this tree-shaped decay is stable under two-regularization-type estimates, which enables us to establish it by adapting the dynamical method developed in \cite{fan2025blockreductionmethodrandom}.\footnote{There is no doubt that the same bound can also be obtained using the tree-representation method of \cite{Band1D,truong2025localizationlengthfinitevolumerandom} or the dynamical approach of \cite{erdos2025zigzagstrategyrandomband}.}

	\medskip
	\noindent
	{\bf The regime $\alpha \in (0, 1)$}. In this regime, the model exhibits features analogous to those of short-range random band matrices in dimensions $d\ge 3$. This analogy is partially reflected in the behavior of the $\Theta$-propagator, which satisfies the bound
	\begin{equation}\label{introduction_Theta_t_B_t}
		\Theta_t(x,y)\lesssim B_t(|x-y|):=\qa{W^{-1}(|x-y|_N/W+1)^{\alpha-1}+(N(1-t))^{-1}}\cdot (|x-y|_N/\ell_t+1)^{-2\alpha},
	\end{equation}
	where $\ell_t$ is defined in \eqref{eq_slow_theta}. The structure of the \emph{shape parameter} $B_t$ largely parallels that of short-range RBMs in dimension $d\ge 3$ (see \Cref{assumption_input_bound} and \cite[Lemma 2.19]{dubova2025delocalizationnonmeanfieldrandommatrices}), in the sense that
	\[(|x-y|_N/W+1)^{\alpha-1}=\qb{(|x-y|_N/W+1)^{-(d_{\txt{eff}}-2)}}^{1/d_{\txt{eff}}},\quad \text{with}\quad d_{\txt{eff}}:=2/\al.\]
	Here, $(|x-y|_N/W+1)^{-(d_{\txt{eff}}-2)}$ coincides with the characteristic \emph{short-range decay} of the $\Theta$-propagator for $d_{\txt{eff}}$-dimensional regular RBMs. This correspondence suggests that the regime $\al\in(0,1)$ presents greater analytical challenges than other regimes, and indicates that techniques developed for high-dimensional random band matrices in \cite{dubova2025delocalizationnonmeanfieldrandommatrices,yang2021delocalization,YYYTexpansion2022CMP} may be relevant.
	
	Motivated by \cite{dubova2025delocalizationnonmeanfieldrandommatrices}, a natural approach is to extend the loop hierarchy analysis developed there to the present setting. However, several essential obstacles arise. As in the regime $\al\in[1,\infty)$, the precise long-range spatial decay is not tracked in \cite{dubova2025delocalizationnonmeanfieldrandommatrices} for resolvent loops of order $n\ge 3$. In fact, an interesting observation in that work is that, even though the $\cK$-loops exhibit short-range decay in $|x-y|_N/W$, the corresponding loop bounds in the max-norm can be propagated along the flow without loss (after suitable one-regularization), and only the sharp decay of length-2 resolvent loops is required. This mechanism again relies on the effective cutoff of resolvent entries beyond the scale $\ell_t$.
	In the present setting, however, it is necessary to track the full long-range decay of resolvent loops of arbitrary length. A natural next step would therefore be to combine the ideas of \cite{dubova2025delocalizationnonmeanfieldrandommatrices} with those of \cite{erdos2025zigzagstrategyrandomband} in order to recover such decay in a manner stable under the flow. Unfortunately, these two approaches are not directly compatible. The main obstruction is the failure of the convolution inequality \begin{equation}
		B_t^{1/2}(0)\sum_{a\in\ZL}B_t^{1/2}(|x-a|)B_t^{1/2}(|a-y|)\lesssim \frac{1}{1-t}B_t^{1/2}(|x-y|),\label{eq:half-half}
	\end{equation}
	which corresponds to assumption (2.24) in \cite{erdos2025zigzagstrategyrandomband}. Heuristically, this inequality expresses that summing over the decay factors associated with two resolvent entries $(G_{t})_{xa}$ and $(G_{t})_{ay}$ reproduces the decay of the single entry $(G_{t})_{xy}$, in accordance with Ward's identity. 
    Its failure in the present setting, as well as in regular RBMs due to the short-range component \smash{$(|\cdot |_N/W+1)^{-(d_{\txt{eff}}-2)}$} of $B_t$ in dimensions $d_{\txt{eff}}>2$, invalidates the stopping-time argument of \cite{erdos2025zigzagstrategyrandomband}. 
    To overcome this difficulty, \cite{dubova2025delocalizationnonmeanfieldrandommatrices} develops a refined diagrammatic expansion technique combined with a \emph{nested structure} of graphs---an approach originating from \cite{yang2021delocalization,YYYTexpansion2022CMP}---to handle the resulting technical complications.

	In attempting to adapt the approaches of \cite{dubova2025delocalizationnonmeanfieldrandommatrices,erdos2025zigzagstrategyrandomband} to our setting with $\al\in[0,1)$, a substantially more severe and fundamental difficulty arises. Namely, the propagation of sharp loop decay factors---represented by loops whose edges each contribute a factor \smash{$B^{1/2}_t$}---for resolvent loops of length $n\ge 3$ produces a singular term that, in contrast to other regimes, cannot be removed by the sum-zero operator of \cite{Band1D} or by the one- or two-regularization procedures of \cite{erdos2025zigzagstrategyrandomband}.
	This obstruction originates from the extremely slow long-range decay of $B_t(|x-y|)$ in \eqref{introduction_Theta_t_B_t}. More precisely, for $|x-y|_N\ge \ell_t$, the decay behaves as  \smash{$|x-y|_N^{-(1+\al)}$}, so that the corresponding half-decay factor  \smash{$|x-y|_N^{-(1+\al)/2}$} is no longer summable in either variable when $\al\in[0,1)$. This lack of summability invalidates the regularization mechanisms that are effective in all other cases, including regular RBMs and power-law models with $\al\ge 1$. As a consequence, even at a formal level---disregarding error terms in the evolution equation for loops (see \eqref{evolution_n_cL_cK_integrated})---the sharp loop decay structure is not preserved under the evolution kernel (defined in \eqref{eq:bUst}), and fails to propagate to later times; see \Cref{rmk_failureofLoop} for a more formal discussion.

	Interestingly, although the loop decay structure breaks down, the stronger \emph{tree-shaped decay}---in which each edge of a tree graph carries a full $B_t$ factor---remains stable under the evolution. This decay structure, established for the $\cK$-loops,\footnote{In this work, we establish this bound only for $\al\ge 1$, although the argument extends to $\al\in(0,1)$ as well.} is compatible with the one-regularization procedure and allows for the removal of the singular terms. While such a strong decay cannot be expected for general resolvent loops of length $n\ge 3$, it coincides with the loop decay in the simplest case of 2-loops, where both notions reduce to the shape parameter $B_{t}(|x-y|)$.
	This observation motivates a shift in strategy. To retain the advantages of the dynamical framework, we restrict our analysis to \emph{resolvent loops and chains of length at most 2}, as defined in \eqref{eq:nGloops} with $n\le 2$. Following the RBM literature, we refer to chains of length 2 as \emph{$T$-variables}, noting that diagonal 2-chains correspond to the first definition in \eqref{eq:TcL}.
	We are thus led to a rather constrained setting in which only loops and chains of length at most 2 are available as analytical tools for power-law models with $\al\in(0,1)$. This poses a significant challenge, as longer loops play a crucial role in prior works \cite{Band1D,erdos2025zigzagstrategyrandomband,Band2D,dubova2025delocalizationnonmeanfieldrandommatrices} in truncating and controlling the loop hierarchy that naturally arises in the evolution of 2-loops. 
    To overcome this difficulty, we develop a new mechanism to truncate the hierarchy at order 2. Roughly speaking, we incorporate the dynamics of the $T$-variables into our analysis, where the control of off-diagonal $T$-variables partially replaces the role played by higher-order loops in previous approaches.\footnote{We remark that this new methodology---without relying on higher-order resolvent loops---can also be applied to regular RBMs in dimensions $d\ge 2$. It yields delocalization, QUE, and bulk universality results analogous to those in \cite{Band2D,dubova2025delocalizationnonmeanfieldrandommatrices}, but with a significantly simpler argument.  However, this simplification comes at the cost of weaker quantum diffusion estimates for 2-loops and the loss of higher-order loop control.}

	As in \cite{dubova2025delocalizationnonmeanfieldrandommatrices}, our analysis begins by addressing the failure of the convolution inequality \eqref{eq:half-half}. A key step in this direction is to obtain sufficiently strong bounds on the \emph{tadpole resolvent diagram}
	\begin{equation}
		f_{xy}(G):=\sum_{a,b\in\ZL}G_{xa}G_{ay}S_{ab}[G(z)-m(z)]_{bb},\label{tadpole_intro}
	\end{equation}
	where $m$ denotes the Stieltjes transform of the Wigner semicircle law and represents the deterministic limit of $G$. In \cite[Section 7]{dubova2025delocalizationnonmeanfieldrandommatrices}, such bounds are established via a diagrammatic expansion method, as discussed above. However, in the present power-law setting with slowly decaying variance profiles, the graphs generated by this expansion cannot, in general, be controlled in the desired manner (see \Cref{sec:graph_ideas} for details about this issue).
	At a heuristic level, this difficulty stems from the fact that, in regular RBMs, vertices connected by $S$-edges form local structures---referred to as \emph{molecules} in \cite{yang2021delocalization,dubova2025delocalizationnonmeanfieldrandommatrices}---with spatial extent of order at most $\OO(W)$. This allows for a clean separation between local (intra-molecular) and global (inter-molecular) structures. In contrast, in the power-law setting, the slow decay of the $S$-edges entangles local and global structures, preventing such a decomposition. Moreover, the expansion scheme in \cite{dubova2025delocalizationnonmeanfieldrandommatrices}, which employs all three resolvent expansion formulas from \cite[Lemmas 7.11--7.13]{dubova2025delocalizationnonmeanfieldrandommatrices}, proceeds in a largely unrestricted manner and does not explicitly track local structures. As a result, it may generate unfavorable configurations that cannot be bounded directly. 
	To overcome this issue, we develop a new expansion strategy that is both more economical and more structured. Specifically, we rely exclusively on the $GG$-expansion formula from \cite[Lemma 7.13]{dubova2025delocalizationnonmeanfieldrandommatrices} and perform significantly fewer expansion steps while preserving detailed control over the local structures within molecules. Within this framework, the graph associated with $|f_{xy}(G)|^{p}$, for fixed even $p\in 2\N$, can be expanded into a sum of $\OO(1)$ well-structured diagrams, which we term \emph{snake-like graphs} (see Definition \ref{def_graph_snakes}). These graphs can then be bounded effectively by applying the Cauchy-Schwarz inequality and Ward's identity in a carefully chosen order. We refer to \Cref{sec:graph_ideas} for a detailed exposition of this diagrammatic method and the underlying ideas.

	We now explain how the dynamics of 2-loops and $T$-variables are analyzed. In the study of 2-loops, the longer loops and chains arising in the evolution can be controlled via loop-contraction inequalities (see the proof of \Cref{lemma_martingale_term}), following the ideas in \cite{dubova2025delocalizationnonmeanfieldrandommatrices}, together with the bounds on the tadpole diagram discussed above.
	To close the argument, it remains to control the diagonal $T$-variables appearing in the evolution equations. However, due to the absence of an averaging effect induced by $S^{(x)}$ in \eqref{eq:nGloops}, the loop-contraction inequality applied to $T$-variables incurs a loss in the decay factor, which prevents a simultaneous treatment of 2-loops and $T$-variables (see Remark \ref{rmk:MG_for_T_fails}). To overcome this difficulty, we perform an additional graph expansion that rewrites $T$-variables in terms of 2-loops, up to some errors, thereby allowing the analysis of the dynamics to be closed at the level of 2-loops; see \Cref{lemma_double_difference_decomposition}. The resulting estimates for 2-loops, although non-optimal, are sufficient to establish sharp local laws for individual resolvent entries, which in turn imply eigenvector delocalization.
	
	The proofs of QUE and bulk universality, however, require significantly stronger bounds, in particular on the expectations of 2-loops. In previous works \cite{Band1D,Band2D,dubova2025delocalizationnonmeanfieldrandommatrices}, such estimates are obtained from sharp bounds on both 2-loops and 3-loops. In our setting, bounds on 3-loops can be partially recovered via diagrammatic expansion techniques. By contrast, while the input 2-loop estimates may not be sharp, they must still be sufficiently strong to support the analysis of their expectations. Unfortunately, the available 2-loop bounds fall short of this requirement. This necessitates the derivation of sharper 2-loop estimates, which is, however, a delicate task due to the absence of adequately controlled higher-order loops in their evolution equations.

	Observing that the QUE estimate follows once sufficiently strong bounds on 2-loops are obtained at spectral parameters $z_t$ approaching the optimal scale $\im z_t\gg N^{-1}$, we divide the dynamical analysis into two time regimes: the \emph{non-flat regime} $I_0:=[0,1-W/N]$ and the \emph{flat regime} $I_1:=[1-W/N,1)$. In the latter regime, the shape parameter $B_t$ defined in \eqref{introduction_Theta_t_B_t} is dominated by the flat contribution $(N(1-t))^{-1}$.
	It therefore suffices to improve the 2-loop estimates in the flat regime. This refinement is achieved in three steps. In the first step, we establish bounds on the off-diagonal $T$-variables via a dynamical analysis in the regime $t\in I_0$. This yields an additional small factor $[N(1-t)]^{-c}$, for some constant $c>0$, at the expense of a loss in spatial decay; see \Cref{lemma_weak_local_law_T_along_flow}. This step relies on the previously established weak a priori bounds for 2-loops, combined with the application of loop-contraction inequalities to the martingale terms. Although this procedure of handling the martingale terms incurs a loss of decay, such a loss is immaterial for the subsequent analysis in the flat regime $I_1$, where the shape parameter $B_t$ is already essentially flat. In contrast, the gain of the additional factor $[N(1-t)]^{-c}$ plays a crucial role in the next step.
	In the second step, we improve the 2-loop estimates in the flat regime via a stopping-time argument that simultaneously tracks 2-loops and $T$-variables. A key observation is that the martingale terms in the evolution equations for both quantities can be controlled directly through the off-diagonal $T$-variables. This mechanism bypasses the limitations of the loop-contraction approach and leads to a stronger bound on the 2-loops in the flat regime.
	In the final step, we combine the improved 2-loop estimates with the off-diagonal $T$-variable bounds obtained in the previous steps to derive sufficiently strong estimates for the expectations of 2-loops. These bounds imply the desired QUE results for bulk eigenvectors in the delocalized phase. The delocalization and QUE estimates can then be used to establish bulk universality of local eigenvalue statistics by following the strategy of \cite{xu2024bulk}.

	\nc

	\medskip
	\noindent
	{\bf Organization of the remaining text.} The remainder of the paper is organized as follows. In \Cref{sec_model_and_main_results}, we introduce the power-law random band matrix model studied in this work, formulate the regularity assumptions on the $\Theta$-propagators, and state our main results.
	The proofs of the main results in the regime $\alpha\in(0,1)$ are presented in Sections \ref{sec:pfal01}--\ref{sec_proof_of_que}. Specifically, in \Cref{sec:pfal01} we develop the underlying dynamical framework of the proof and reformulate the main results in a dynamical form, given by Lemmas \ref{lemma_weak_local_law_T_along_flow}--\ref{lemma_bootstrap}. Section \ref{sec_step_1_to_3} is devoted to the proof of \Cref{lemma_bootstrap}, based on a key high-moment estimate for the tadpole diagram \eqref{tadpole_intro}, stated in Lemma \ref{claim_graph_bound_light_weight_u}. The proof of this lemma is given in \Cref{sec_graphical_argument}, where we develop a new set of graphical tools. The remaining results---namely, Lemmas \ref{lemma_weak_local_law_T_along_flow}--\ref{theorem_que_along_flow}---are established in \Cref{sec_proof_of_que} via an analysis of the dynamics of the $T$-variables.
	The proofs of the main results in the remaining regimes $\alpha\in(-1,0)$, $\alpha=0$, and $\alpha\in[1,\infty)$ are given in Sections \ref{sec_proof_of_remaining_cases_2}, \ref{sec_proof_alpha_0}, and \ref{sec_proof_of_remaining_cases_1}, respectively. 

To streamline the presentation, we defer the proofs of several auxiliary results to the appendix. Some of these results are included primarily for illustrative purposes, others are relatively straightforward, and still others follow arguments similar to those in the existing literature. Moreover, they are all largely independent of the main proofs.

	\medskip
	{\noindent {\bf Notation.}}
	To facilitate the presentation, we introduce some necessary notation that will be used throughout this paper. We will use the set of natural numbers $\N=\{1,2,3,\ldots\}$ and the complex upper half-plane $\C_+:=\{z\in \C:\im z>0\}$.  
	In this paper, we are interested in the asymptotic regime with $N\to \infty$. When we refer to a constant, it will not depend on $N$ or $W$. Unless otherwise noted, we will use $C$, $D$ etc.~to denote large positive constants, whose values may change from line to line. Similarly, we will use $\e$, $\delta$, $\tau$, $c$, $\fc$, $\fd$ etc.~to denote small positive constants. 
	For any two (possibly complex) sequences $\xi_N$ and $\zeta_N$ depending on $N$, $\xi_N = \OO(\zeta_N)$, $\zeta_N=\Omega(\xi_N)$, or $\xi_N \lesssim \zeta_N$ means that $|\xi_N| \le C|\zeta_N|$ for some constant $C>0$, whereas $\xi_N=\oo(\zeta_N)$ or $|\xi_N|\ll |\zeta_N|$ means that $|\xi_N| /|\zeta_N| \to 0$ as $N\to \infty$. We say that $\xi_N \asymp \zeta_N$ if $\xi_N = \OO(\zeta_N)$ and $\zeta_N = \OO(\xi_N)$. For any $\al,\beta\in\R$, we denote $\llbracket \al, \beta\rrbracket: = [\al,\beta]\cap {\mathbb Z}$, $\qq{\al}:=\qq{1,\al}$, $\al\vee \beta:=\max\{\al, \beta\}$, and $\al\wedge \beta:=\min\{\al, \beta\}$. 
	Given a vector $\mathbf v$, $|\mathbf v|\equiv \|\mathbf v\|_2$ denotes the Euclidean norm and $\|\mathbf v\|_p$ denotes the $\ell^p$-norm. 
	Given a matrix $\cal A = (\cal A_{ij})$, $\|\cal A\|$, $\|\cal A\|_{p\to p}$, and $\|\cal A\|_{\infty}\equiv \|\cal A\|_{\max}:=\max_{i,j}|\cal A_{ij}|$ denote the operator (i.e., $\ell^2\to \ell^2$) norm,  $\ell^p\to \ell^p$ norm (where we allow $p=\infty$), and maximum (i.e., $\ell^\infty$) norm, respectively. We will use $\cal A_{ij}$ and $ \cal A(i,j)$ interchangeably in this paper. We will use $I_n$ to denote an $n\times n$ identity matrix.
	
	Given an event $\Xi$, let $\mathbf 1_\Xi$ or $\mathbf 1(\Xi)$ denote its indicator function. We will say an event $\Xi$ holds with high probability (w.h.p.) if for any constant $D>0$, $\mathbb P(\Xi)\ge 1- N^{-D}$ for large enough $N$. For clarity of presentation, we will use the following notion of stochastic domination introduced in \cite{Average_fluc}. Let \(\xi=\left(\xi^{(N)}(u):N\in\mathbb N, u\in U^{(N)}\right)\) and \(\zeta=\left(\zeta^{(N)}(u):N\in\mathbb N, u\in U^{(N)}\right)\)
	be two families of non-negative random variables, where $U^{(N)}$ is a possibly $N$-dependent parameter set. We say $\xi$ is stochastically dominated by $\zeta$, uniformly in $u$, if for any fixed (small) $\tau>0$ and (large) $D>0$, 
	\be\label{stoch_domination}\mathbb P\bigg(\bigcup_{u\in U^{(N)}}\left\{\xi^{(N)}(u)>N^\tau\zeta^{(N)}(u)\right\}\bigg)\le N^{-D}\ee
	for large enough $N\ge N_0(\tau, D)$, and we will use the notation $\xi\prec\zeta$. If for some complex family $\xi$ we have $|\xi|\prec\zeta$, then we will also write $\xi \prec \zeta$ or $\xi=\OO_\prec(\zeta)$. As a convention, for two \emph{deterministic} quantities $\xi$ and $\zeta$, we will write $\xi\prec\zeta$ if and only if $|\xi|\le N^\tau |\zeta|$ for any constant $\tau>0$ for large enough $N\ge N_0(\tau)$. 
	
	\subsection*{Acknowledgement}  
	Fan Yang is supported in part by the National Key R\&D Program of China (No.~2023YFA1010400) and NSFC (No.~12526201). We would like to thank Guangyi Zou for valuable discussions. AI tools were used only for polishing the language of this paper. All mathematical results, arguments, and computations were developed and written by the authors.

	\section{The model and main results}\label{sec_model_and_main_results}

	\subsection{The power-law random band matrix model}

Throughout the paper, we identify $\ZL=\Z/(N\Z)$ as a discrete circle. Accordingly, for $x,y\in \ZL$, we define the periodic distance by 
	\[ |x-y|\equiv |x-y|_{N}:=\min_{k\in\Z} | x- y + kN|. \]
We now state the precise assumptions of the power-law random band matrix model.  
	
	\begin{definition}[Power-law random band matrices]\label{considered_model}
		We consider a class of $N\times N$ complex Hermitian random matrices $H:=(h_{xy}:x,y\in \ZL)$ defined on the discrete circle $\ZL$, where the entries $h_{xy}$ are independent Gaussian random variables up to the Hermitian symmetry $h_{xy}=\overline h_{yx}$. More precisely, given an $N\times N$ symmetric doubly stochastic variance matrix $S=(S_{xy}:x,y\in \ZL)$, the diagonal entries of $H$ are real Gaussian random variables, while the off-diagonal entries are complex Gaussian random variables, distributed as follows: 
		\be\label{bandcw0}
		h_{xy}\sim \begin{cases}
			\mathcal{N}_{\R}(0, S_{xy}),\ & \text{if}\ x=y,\\ 
			\mathcal{N}_{\C}(0, S_{xy}),\ & \text{if}\ x\ne y.  
		\end{cases}
		\ee
		We assume that the entries of $S$ exhibit a power-law decay on a scale $1\le W\le N/2$ with a fixed exponent $\alpha\in\R$. Specifically, there exists a constant $C_{\txt{RBM}}\geq 1$ such that 
		\begin{equation}\label{alpha_decay}
			\begin{aligned}
				0\le S_{xy}\leq C_{\txt{RBM}}\cS_{xy},\quad \text{with}\quad \cS_{xy}:=\frac{1}{\cZ_{\alpha}}\pa{\frac{|x-y|}{W}+1}^{-1-\alpha},
			\end{aligned}
		\end{equation}
		where $\cZ_{\alpha}$ is the normalization constant
		\begin{equation*}
			\begin{aligned}            \cZ_{\alpha}\equiv\cZ_{\alpha,N,W}:=\sum_{x\in\ZL}\pa{\frac{|x|}{W}+1}^{-1-\alpha}.
			\end{aligned}
		\end{equation*}
		Elementary calculus classifies the order of $\cZ_\alpha$ as follows:
		\begin{equation}
			\cZ_\alpha \asymp 
			\begin{cases} 
				W(W/N)^\alpha, & \text{if } \alpha < 0, \\
				W \log(N/W + 1), & \text{if } \alpha = 0, \\
				W, & \text{if } \alpha > 0.
			\end{cases}
		\end{equation}
		We refer to $H$ defined above as a \emph{power-law random band matrix with bandwidth $W$ and exponent $\alpha$}.
	\end{definition}
	
	\begin{remark}
		We remark that the focus of this work is on the regime $\alpha\in(-1,\infty)$. Indeed, when $\alpha\leq -1$, the assumption \eqref{alpha_decay} in \Cref{considered_model} implies that $\max_{x,y\in \ZL}S_{xy}\lesssim N^{-1}$. 
		Random matrix models satisfying such mean-field-type conditions fall into the well-studied class of generalized Wigner matrices and have been analyzed in a series of works under quite general assumptions; see, for example, 		\cite{EKYYEJP18,QuadraticAOEK,ajanki2017universality,erdHos2024eigenstate}. In particular, for $\alpha\leq -1$, under mild regularity assumptions, the analogues of our main results follow as straightforward corollaries of the main results established in
		\cite{ajanki2017universality} and \cite{erdHos2024eigenstate}.
	\end{remark}
	
	We denote the eigenvalues of $H$ by $\lambda_1\leq \lambda_2\leq \cdots \leq \lambda_{N}$, and the corresponding normalized eigenvectors by $(\boldsymbol{\psi}_k)_{k=1}^{N}$. It is well known that, under the normalizing condition $\sum_y S_{xy}\equiv 1$, the empirical spectral measure \smash{$N^{-1}\sum_{k=1}^{N}\delta_{\lambda_k}$} converges almost surely to the famous Wigner semicircle law, whose density is \(\rho_{\txt{sc}}\pa{x}=\sqrt{(4-x^2)_+}/2\pi\) with spectral edges at $\pm 2$. 
	Define the Green's function (or resolvent) of $H$ by
	\begin{equation}\label{def_Green}
		\begin{aligned}
			G\p{z}:=(H-z)^{-1},\quad z\in\C_+.
		\end{aligned}
	\end{equation}
	As $N\to \infty$, the Green's function $G\p{z}$ converges entrywise to the scalar matrix $m(z)I_N$, where $m(z)$ is the Stieltjes transform of $\rho_{\txt{sc}}$, defined by
	\begin{equation}\label{def_m_sc}
		\begin{aligned}
			m\p{z}\equiv m_{\txt{sc}}\pa{z}:=\frac{-z+\sqrt{z^2-4}}{2}=\int_{\R}\frac{\rho_{\txt{sc}}\pa{x}}{x-z}\,\rd x.
		\end{aligned}
	\end{equation}
	One readily verifies that $m(z)$ satisfies the self-consistent equation 
	\begin{equation}\label{eq:self_cons_m}    
		m\p{z}=-\pa{m\p{z}+z}^{-1}.
	\end{equation} 
	In this work, we focus on the bulk regime of the spectrum. Specifically, we always consider spectral parameters $z$ lying in the domain
	\begin{equation}\label{def_spectral_domain}
		\begin{aligned}
			\mathbf{D}_{\kappa,\fc}:=\ha{z=E+\ii\eta\in\C_+: |E|\leq 2-\kappa, N^{\fc}\eta_*\leq \eta \leq \fc^{-1}}
		\end{aligned}
	\end{equation}
	for fixed small constants $\kappa,\fc>0$, where the scale $\eta_*$ is defined in \eqref{def_eta_star_and_W_c} below.

	In addition to \Cref{considered_model}, we impose further technical assumptions on the ``diffusion profile'' of the model, namely on the Green's function of the classical random walk on $\ZL$ with transition probabilities given by the variance matrix $S$. For notational consistency, we introduce a collection of shape parameters that provide upper bounds on these diffusion profiles in different regimes of $\alpha\in(-1,\infty)$.

	\begin{definition}[Shape parameters]\label{definition_shape_parameter}
		For any $\alpha\in (-1,\infty)$, spectral scale $\eta>0$, and length scale $r\in[0,\infty)$, we define the \emph{shape parameter} $B(\eta,r)$ by
		\begin{equation}\label{def_B_t}
			B(\eta,r):=\begin{dcases}
				W^{-1}(N/W)^{\alpha}(r/W+1)^{-1-\alpha}+(N\eta)^{-1}, \ & \text{if}\ \ \ \alpha\in(-1,0],\\
				[W^{-1}(r/W+1)^{\alpha-1}+(N\eta)^{-1}]\cdot (r/\ell(\eta)+1)^{-2\alpha},\ & \text{if}\ \ \  \alpha\in(0,1),\\
				[\eta\ell(\eta)]^{-1}\cdot (r/\ell(\eta)+1)^{-1-\alpha},\ & \text{if}\ \ \  \alpha\in[1,\infty).
			\end{dcases}
		\end{equation}
		Here, $\ell(\eta)$ denotes the characteristic length scale
		\begin{equation}\label{eq:elleta}
			\ell(\eta):=\begin{dcases}
				N,\ & \text{if}\ \ \ \alpha\in(-1,0],\\
				\min(W\eta^{-\frac{1}{\alpha\wedge 2}},N),\ & \text{if}\ \ \ \alpha\in(0,\infty).
			\end{dcases}
		\end{equation}
		In addition, we introduce the \emph{zero-mode-removed shape parameter} for $\al\ge 0$, which controls the diffusion profile after removing the zero Fourier mode:
		\begin{equation}\label{def_c_B_t}
			\begin{aligned}
				\upcirc{B}\p{\eta,r}:=\begin{dcases}
					W^{-1}\pa{r/W+1}^{\alpha-1},\ & \text{if}\ \ \ \alpha\in[0,1),\\
					N^{-1}(N/W)^{\alpha\wedge 2},\ & \text{if}\ \ \ \alpha\in[1,\infty).
				\end{dcases}
			\end{aligned}
		\end{equation}
		Finally, we define the \emph{difference parameter}, which controls the first-order difference of the diffusion profile for $\al>0$:
		\begin{equation}\label{def_R_t}
			\begin{aligned}
				R\p{\eta,r}:=\begin{dcases}
					W^{-1}(r/W+1)^{-1},\ & \text{if}\ \ \ \alpha\in (0,1),\\
					W^{\alpha-2}[\ell(\eta)]^{1-\alpha}(r/W+1)^{\alpha-2}(r/\ell(\eta)+1)^{1-\alpha},\ & \text{if}\ \ \ \alpha\in[1,2),\\
					\ell(\eta)^{-1}(r/\ell(\eta)+1)^{-1},\ & \text{if}\ \ \ \alpha\in[2,\infty).
				\end{dcases}
			\end{aligned}
		\end{equation}
	\end{definition}
	
	\begin{remark}
		We briefly discuss the physical meaning of the key parameter $\ell(\eta)$ in our model. Heuristically, the Green's function $G(z)$ with $\im z=\eta$ corresponds to the time-evolution operator $e^{\ii Ht}$ under the identification $t=\eta^{-1}$. Thus, for $\al>0$ and in the large-system limit $N\to \infty$, \smash{$\ell(\eta)=\ell(t^{-1})=Wt^{\frac{1}{\alpha\wedge 2}}$} describes the characteristic spatial scale of the wave-packet propagation at time $t$. This interpretation is consistent with the observations in \cite{Power-RBM}. When $\al\ge 2$, the dynamics are diffusive.
		For $1<\al<2$, since the exponent $\al^{-1}>1/2$, the wave packet spreads superdiffusively. 
		For $\al\le 1$, we have \smash{$\ell(t^{-1})=W(t^{\frac{1}{2}})^{\frac{2}{\alpha}}$}, which coincides with the characteristic spatial scale of a $(2/\alpha)$-dimensional short-range random band matrix. \end{remark}

	\begin{assumption}\label{assumption_input_bound}
		Consider a power-law random band matrix $H$ with variance profile $S$ from \Cref{considered_model}, with a fixed exponent $\alpha\ge 0$. Let $c_0\in(0,1)$ be a small constant. Suppose that the following estimates hold for some deterministic parameter $C_{\alpha}\equiv C_{\alpha}(N)>0$, uniformly in $t\in\qa{0,1-N^{-1}}$ and $\xi\in \C$ with $|\xi|=1$ and $|\xi-1|\ge c_0$. 
		For notational convenience, for $t\in[0,1)$, we introduce the abbreviations
		\begin{equation}\label{eq:abbrvB}
			B_t(r)\equiv B(1-t,r),\quad  \upcirc{B}_t\p{r}\equiv \upcirc{B}\p{1-t,r},\quad R_t\p{r}\equiv R\p{1-t,r},\quad \ell_t\equiv \ell (1-t).
		\end{equation}
		\begin{itemize}
			\item {\bf Upper bounds:} For any $x,y\in\ZL$, we have
			\begin{equation}\label{Theta_bound_opposite_charge}
				\begin{aligned}
					\absa{\pa{\frac{S}{1-t S}}_{xy}}\leq C_{\alpha}B_t\p{|x-y|},  
				\end{aligned}
			\end{equation}
			\begin{equation}\label{Theta_bound_equaled_charge}
\absa{\pa{\frac{S}{1-t\xi S}}_{xy}}\leq C_{\alpha}B_0\p{|x-y|}. \end{equation}

			\item {\bf First-order difference bound:} For $t\in\qa{0,1-N^{-1}}$ with $\ell_t < N$, we have for any $x,y,z\in\ZL$,
			\begin{equation}\label{Theta_difference_bound}
				\begin{aligned}
					\quad    \absa{\pa{\frac{S}{1-tS}}_{xy}-\pa{\frac{S}{1-tS}}_{xz}}\prec \braket{y-z}_W R_t(|x-y|\vee |x-z|) {B_t\p{|x-y|\wedge|x-z|}},
				\end{aligned}
			\end{equation}
			where we denote $\braket{y-z}_W:=|y-z|+W$ for any $y,z\in\ZL$.
			\item {\bf Zero-mode-removed upper bound:}  For $t\in\qa{0,1-N^{-1}}$ with $\ell_t=N$, we have for any $x,y\in\ZL$,
			\begin{equation}\label{zero_mode_removed_bound}
				\begin{aligned}
					\absa{\pa{\frac{S}{1-tS}}_{xy}-\frac{1}{N(1-t)}}\prec \upcirc{B_t}\p{|x-y|}.
				\end{aligned}
			\end{equation}
		\end{itemize}
        For $\al\in (0,1)$, we assume that \eqref{Theta_bound_opposite_charge}--\eqref{zero_mode_removed_bound} hold with $C_\alpha\lesssim 1$. For $\al\ge 1$, we assume that \eqref{Theta_bound_opposite_charge}--\eqref{zero_mode_removed_bound} hold with $C_\alpha\prec 1$. Finally, for $\al=0$, we assume that \eqref{Theta_bound_opposite_charge}, \eqref{Theta_bound_equaled_charge}, and \eqref{zero_mode_removed_bound} hold with $C_\alpha\prec 1$. The weaker assumption $C_\alpha\prec 1$ is imposed mainly for the two critical cases $\al\in\{0,1\}$, where additional logarithmic corrections may appear in the estimate \eqref{Theta_bound_opposite_charge}. Nevertheless, the argument for $\al=1$ extends to all $\al\ge 1$ under the assumption $C_\alpha\prec 1$. 
	\end{assumption}

	For $\alpha\ge 0$, we expect the bounds in \Cref{assumption_input_bound} to hold under assumption~\eqref{alpha_decay}, together with some mild additional regularity conditions on the variance profile. A rigorous verification of these bounds, however, would require a rather involved analysis of the resolvent-type operator $S/\p{1-tS}$.
	For instance, under the additional assumption that $S$ is translation invariant, one may exploit its Fourier series representation or invoke local CLTs for the associated random walk on $\ZL$ with transition matrix $S$.\footnote{Depending on the decay exponent, this walk exhibits Gaussian behavior for $\al\ge 2$ and Lévy $\al$-stable behavior for $0<\al<2$.} Developing such an analysis in full generality would substantially deviate from the main objective of the present paper, namely, the derivation of optimal lower bounds on localization lengths for PRBMs with arbitrary decay exponent $\al>0$. We therefore do not pursue this direction here. Instead, we provide a concrete example demonstrating that \Cref{assumption_input_bound} is indeed satisfied by a canonical class of power-law random band matrix models, namely those with a precise variance profile function. These models constitute the primary motivation for this work.

	\begin{example}\label{example_profile_function}
		Fix a parameter $\alpha>0$, and let $f_{\alpha}$ be the continuous density of a symmetric probability distribution with characteristic function $\phi_{\alpha}(t)$. Suppose that the following conditions hold.
		\begin{enumerate}
			
			\item {\bf Power-law decay:} For some constant $C_{\alpha}>1$,
			\begin{equation}\label{f_alpha_decay}
				\begin{aligned}
					0\leq f_{\alpha}(x)\leq C_{\alpha} (1+|x|)^{-1-\alpha}.
				\end{aligned}
			\end{equation}
			\item {\bf Small $t$ behavior:} For some constant $0<c_{\alpha}<1$,
			\begin{equation}\label{psi_spectral_gap}
				\begin{aligned}
					1-\phi_{\alpha}(t)\geq c_{\alpha}\pa{ |t|^{\alpha\wedge 2}\wedge 1}.
				\end{aligned}
			\end{equation}
			\item {\bf Regularity of the characteristic function:} The characteristic function $\phi_{\alpha}(t)$ is $(\ceil{\alpha}+2)$-times differentiable on $\R\setminus\ha{0}$ and $(\ceil{\alpha}-1)$-times differentiable at $t=0$. Moreover, given a constant $\delta>0$, for any sufficiently small \(\varepsilon>0\) there exists a constant $C\equiv C(\alpha,\varepsilon)>1$ such that
			\begin{equation}\label{phi_alpha_regularity}
				\begin{aligned}
					\absa{\phi_{\alpha}^{\pa{\ceil{\alpha}+k}}(t)} \leq C|t|^{\alpha-\ceil{\alpha}-k}(1+\mathbf{1}_{\alpha\geq 1}\cdot |t|^{-\varepsilon}),\quad \forall k\in\ha{0,1,2}\txt{ and } t\in(-\delta,\delta)\setminus\ha{0}.
				\end{aligned}
			\end{equation}
			The additional singular factor \(|t|^{-\varepsilon}\) is included to control possible logarithmic corrections when \(\alpha\) is an integer. For instance, logarithmic factors in \(|t|\) appear for the \(1\)-stable density \(f_1\).

			\item {\bf Decay of the characteristic function:} There exist constants $\varepsilon,\delta>0$ and $C\equiv C(\alpha,\varepsilon)>1$ such that 
			\begin{equation}\label{phi_alpha_decay}
				\begin{aligned}
					\max_{k\in\ha{0,1,\ldots,\ceil{\alpha}+2}}\absa{\phi_{\alpha}^{(k)}(t)}\leq C\pa{1+|t|}^{-2-\varepsilon},\quad \forall |t|\geq \delta.
				\end{aligned}
			\end{equation}
			Without loss of generality, we assume that the constant \(\delta\) here coincides with that in \eqref{phi_alpha_regularity}.
			
		\end{enumerate}
		This setting includes a broad class of power-law densities. In particular, one can readily verify that $\alpha$-stable distributions for $\alpha\in(0,2)$ and Student's $t$-distributions satisfy the above assumptions.     
		We now define the variance profile induced by the profile function $f_{\alpha}$ as
		\begin{equation}
			\begin{aligned}
				S_{xy}:=\frac{1}{Z_{\alpha}}\sum_{n\in\Z}f_{\alpha}\pa{\frac{x-y+nN}{W}},\quad \forall x,y\in\ZL,
			\end{aligned}
		\end{equation}
		where the normalization constant is given by
		\begin{equation}
			\begin{aligned}
				Z_\alpha:=\sum_{x\in\ZL}\sum_{n\in\Z}f_{\alpha}\pa{\frac{x+nN}{W}}=\sum_{x\in\Z}f_{\alpha}\pa{\frac{x}{W}}.
			\end{aligned}
		\end{equation}
		Assuming $W\geq N^{c}$ for some constant $c>0$, the properties (i)--(iv) imply that the variance profile $S$ satisfies the assumptions of \Cref{assumption_input_bound} for any $\al>0$. The detailed proof is postponed to \Cref{appendix_proof_of_example_profile_function}.
		Finally, we remark that in the critical regime $\alpha=0$, precise estimates for $S/\p{1-tS}$, in particular the correction factors, are model dependent and highly sensitive to the detailed structure of $f_0$. For this reason, we do not consider this case here.
		
	\end{example}

	For $\alpha<0$, the corresponding estimates \eqref{Theta_bound_opposite_charge}, \eqref{Theta_bound_equaled_charge}, and \eqref{zero_mode_removed_bound} remain valid under the following spectral gap condition; see \Cref{lemma_input_bound_subcritical}.
	\begin{assumption}\label{assumption_input_bound<0}
		Consider a power-law random band matrix $H$ with variance profile $S$ from \Cref{considered_model}, with a fixed exponent $\alpha \in (-1,0)$. Assume there exists a constant $c_\alpha>0$ such that $\spec(S)\subset [-1+c_\al,1-c_\al]\cup \{1\}$ and the Perron-Frobenius eigenvalue $1$ is simple. 
	\end{assumption}

    This spectral gap condition can be readily guaranteed by imposing suitable irreducibility assumptions. For example, it is sufficient to assume that there exists a fixed \(k\in\mathbb{N}\) such that
\(
\min_{x,y}(S^k)_{xy}\gtrsim N^{-1}.
\)
In particular, this condition is satisfied if, in addition, one assumes that \(S_{xy}\ge c\mathcal{S}_{xy}\) for some constant \(c>0\).

	\subsection{Main results}
	
	To state the main results of this work, we first introduce two key parameters, $\eta_*$ and $W_c$, depending on $\alpha$, defined by
	\begin{equation}\label{def_eta_star_and_W_c}
		\eta_* := 
		\begin{cases} 
			N^{-1},                                  &\text{for} \ \ -1< \alpha \le 1 \\
			W^{-\frac{\alpha}{\alpha-1}} + N^{-1},   &\text{for} \ \ 1 < \alpha < 2, \\
			W^{-2} + N^{-1}                         &\text{for} \ \ \alpha \ge 2
		\end{cases},
		\quad \text{and} \quad 
		W_c := 
		\begin{cases} 
			1,                 &\text{for} \ \ -1< \alpha \le 1 \\
			N^{1-1/\alpha},    &\text{for} \ \ 1 < \alpha < 2 \\
			N^{1/2},           &\text{for} \ \ \alpha \ge 2
		\end{cases}.
	\end{equation}
	Here, $\eta_*$ represents the minimal scale of $\im z$ at which a precise local law for $G(z)$ can be established. The quantity $\eta_*^{-1}$ also characterizes the localization length of bulk eigenvectors (recall \eqref{definition_ell_P}). The parameter $W_c$ denotes the critical bandwidth at which the localization length $\eta_*^{-1}$ becomes comparable to the system size $N$. Thus, as $W$ crosses $W_c$, the system is conjectured to undergo an Anderson localization-delocalization transition. Our first main result establishes an (almost) sharp upper bound on the $\ell^\infty$-norm of bulk eigenvectors for power-law random band matrices across all regimes of $\alpha$. In particular, it yields the correct lower bound on the localization length of bulk eigenvectors.
	
	\begin{theorem}[Localization length]\label{theorem_delocalization}
		Consider the power-law random band matrix $H$ introduced in \Cref{considered_model}. Assume that its variance profile matrix $S$ satisfies \Cref{assumption_input_bound} when $\alpha\ge 0$, and \Cref{assumption_input_bound<0} when $\alpha\in(-1,0)$.
		Let $\kappa,\varepsilon\in(0,1)$ be arbitrarily small constants. We further assume that the bandwidth satisfies $W\ge N^{\varepsilon}$ in the case $\alpha\ge 0$. Then, the following estimate holds for any constants $\tau,D>0$ when $N$ is sufficiently large (recall that $\{\lambda_k\}$ and $\{\boldsymbol{\psi}_k\}$ denote the eigenvalues and eigenvectors of $H$):
		\begin{equation}\label{eq:delocal}
			\begin{aligned}
				\P\pa{\max_{k:|\lambda_k|\leq 2-\kappa}\|\boldsymbol{\psi}_k\|_{\infty}^2\leq N^{\tau}\eta_*}\geq 1-N^{-D}.
			\end{aligned}
		\end{equation}
		As a consequence, recalling the definition of the localization length in \eqref{eq:localization_length}, we obtain that for any constants $\tau,D>0$ and sufficiently large $N$,
		\begin{align}
			\P\pa{\min_{k: |\lambda_k | \leq 2 - \kappa} \ell(\boldsymbol{\psi}_k) \ge N^{-\tau}\eta_*^{-1}}\ge 1-N^{-D}.\label{eq:length}  
		\end{align}
		In particular, if $W\geq \mathbf{1}_{\alpha\geq 0}\cdot N^{\varepsilon}  W_c +1$, then the bulk eigenvectors are completely delocalized with high probability, i.e., for any constants $\tau,D>0$ and sufficiently large $N$,
		\begin{equation}\label{eq:delocal2}
			\begin{aligned}
				\P\pa{\max_{k:|\lambda_k|\leq 2-\kappa}\|\boldsymbol{\psi}_k\|_{\infty}^2\leq N^{-1+\tau}}\geq 1-N^{-D}.
			\end{aligned}
		\end{equation}
	\end{theorem}
	
	\Cref{theorem_delocalization} is an immediate consequence of the following local law for the resolvent defined in \eqref{def_Green}.

	\begin{theorem}[Local semicircle law]\label{theorem_local_law_entrywise}
		Under the assumptions of \Cref{theorem_delocalization}, for any constants $\fc,\tau,D>0$, the following events hold with probability at least $1-N^{-D}$ for sufficiently large $N$:
		\begin{align}
			&\bigcap_{z\in\mathbf{D}_{\kappa,\fc}}\bigcap_{x,y\in\ZL}\ha{\absa{G_{xy}(z)-m(z)\delta_{xy}}^2\leq N^{\tau}B(\im z, |x-y|)},\label{eq_entrywise_local_law}\\
			&\bigcap_{z\in\mathbf{D}_{\kappa,\fc}}\bigcap_{x\in\ZL}\ha{\absbb{\sum_{y\in\ZL}S_{xy}G_{yy}(z)-m(z)}\leq N^{\tau}B(\im z,0)},\label{eq_average_local_law}
		\end{align}
		where the spectral domain $\mathbf{D}_{\kappa,\fc}$ is defined in \eqref{def_spectral_domain}, and $B(\im z,\cdot)$ is introduced in \Cref{definition_shape_parameter}.
	\end{theorem}

	\begin{proof}[\bf Proof of \Cref{theorem_delocalization}]
		The estimate \eqref{eq:delocal} follows directly from the entrywise local law \eqref{eq_entrywise_local_law} via the bound 
		\be\label{eq:ukx}|\boldsymbol{\psi}_k(x)|^2 \le \eta\im G_{xx}(\lambda_k + \ii \eta),\quad \forall \eta>0\, .\ee 
		Applying \eqref{eq_entrywise_local_law} to $G_{xx}(\lambda_k+\mathrm{i}\eta)$ with $\eta = N^{\tau}\eta_*$, we obtain $\im G_{xx}(\lambda_k+\mathrm{i}\eta)=\OO(1)$ with high probability. Combined with \eqref{eq:ukx}, this yields \eqref{eq:delocal}. The bounds \eqref{eq:length} and \eqref{eq:delocal2} follow immediately from \eqref{eq:delocal}.
	\end{proof}

	In the complete delocalization regime, namely when $W\geq \mathbf{1}_{\alpha\geq 0}\cdot N^{\varepsilon}  W_c +1$,
	one can further establish a \emph{quantum unique ergodicity} (QUE) estimate for bulk eigenvectors, as stated below. 
	
	\begin{theorem}[Quantum unique ergodicity]\label{theorem_que}
		In the setting of \Cref{theorem_delocalization}, assume further that $W\geq \mathbf{1}_{\alpha\geq 0}\cdot N^{\varepsilon}  W_c +1$.
		Given any $E\in\qa{-2+\kappa,2-\kappa}$, we define the interval
		\[\cI_{E}\equiv\cI_E(\varepsilon_0):=\ha{x\in\R:|x-E|\leq \eta(\varepsilon_0)},\qquad \forall \varepsilon_0\in(0,\mathfrak{c}_0),\]
		where the parameter $\eta(\varepsilon_0)$ and the constant $\mathfrak{c}_0$ are defined as
		\begin{equation*}
			\begin{aligned}
				\eta(\varepsilon_0):=\begin{cases}            
					N^{-1-\varepsilon_0}\cdot W^{(1+\alpha)/2}N^{-\alpha/2},\ &\txt{if}\ \  \alpha\in(-1,0),\\
					N^{-1-\varepsilon_0}\cdot W^{1/2},\ &\txt{if}\ \  \alpha\in[0,1)\\
					N^{-1/2-\varepsilon_0}(W/N)^{(\alpha\wedge 2)/2}\ &\txt{if}\ \  \alpha\in[1,\infty)
				\end{cases},
				\quad
				\fc_0:=\begin{cases}
					-\alpha/2,\ &\text{if} \ \ \alpha\in(-1,0)\\
					\varepsilon/2,\ &\text{if} \ \ \alpha\in[0,1)\\
					(\alpha\wedge 2)\cdot (\varepsilon/2),\ &\text{if} \ \ \alpha \in [1,\infty)
				\end{cases}.
			\end{aligned}
		\end{equation*}
		Note that $\fc_0$ is chosen such that $\eta(\e_0)\ge N^{-1+\fc_0-\varepsilon_0}$. Define $\fc_1>0$ as
		\begin{equation*}
			\begin{aligned}
				\fc_1:=\begin{cases}
					\pa{2\varepsilon_0}\wedge\qa{(\fc_0-\varepsilon_0)/4},\ &\text{if} \ \ \alpha\in(-1,1)\\
					(2\varepsilon_0)\wedge \qa{(\fc_0-\varepsilon_0)/2},\ &\text{if} \ \ \alpha \in [1,\infty)
				\end{cases}.
			\end{aligned}
		\end{equation*}
		Then, for any constant $c\in(0,\fc_1/2)$, the following estimate holds for any constant $\tau>0$:
		\begin{equation}\label{eq_que}
			\begin{aligned}
				\sup_{|E|\leq 2-\kappa}\max_{x\in\ZL}\P\pa{\max_{i,j:\lambda_i,\lambda_j\in\cI_E}\absa{\sum_{a\in\ZL}S_{xa} \ol{\bpsi}_{i}(a)\bpsi_j(a)-\delta_{ij}N^{-1}}\geq N^{-1-c}}\leq N^{-\fc_1+2c+\tau}.
			\end{aligned}
		\end{equation}
	\end{theorem}
	
	When $\alpha>0$, Theorem \ref{theorem_que} further implies a more conventional form of QUE for bulk eigenvectors that is commonly used in the literature. Roughly speaking, it asserts that, with probability $1-\oo(1)$, every bulk eigenvector of $H$ is asymptotically equidistributed in $\ell^2$-mass on all scales larger than $W$. In particular, in the complete delocalization regime, this implies that the localization length of every bulk eigenvector satisfies $\ell(\boldsymbol{\psi}_k)=\Omega(N)$ with probability $1-\oo(1)$. 
	\begin{corollary}\label{cor_que}
In the setting of \Cref{theorem_que}, suppose that a non-empty subset $I\equiv I_N\subseteq \ZL$ satisfies
		\begin{equation}\label{eq_que_I_assumption}
			\frac{1}{|I|}\sum_{a\in I}\sum_{b\in \ZL\setminus I}S_{ab}\leq N^{-\fc}
		\end{equation}
		for some constant $\fc>0$. Then there exists a constant $c>0$, depending on $\fc$ and $\fc_1$, such that
		\begin{equation}\label{eq_flat_que}
			\sup_{|E|\leq 2-\kappa}\P\pa{\max_{k:\lambda_k\in\cI_E}\absa{\frac{1}{|I|}\sum_{a\in I}\absa{\bpsi_k(a)}^2-N^{-1}}\geq N^{-1-c}}\leq N^{-c}.
		\end{equation}
	\end{corollary}

	Note that when $\alpha>0$, any interval $I=\qq{x,y}\subseteq \ZL$ with $|I|\geq W^{1+\delta}$ for some constant $\delta>0$ satisfies the assumption \eqref{eq_que_I_assumption} for a suitable constant $\fc>0$. In contrast, when $\alpha\leq 0$, the variance profile is sufficiently flat that the assumption \eqref{eq_que_I_assumption} cannot hold unless the subset $I$ satisfies $|\ZL\setminus I|\le N^{1-\delta}$ for some constant $\delta>0$. Establishing sharp results in the regime $\alpha\le 0$ would require proving a refined two-resolvent local law for quantities of the form $\avg{G(z_1)AG(z_2)A^*}$ with $z_1,z_2\in \{z,\overline z\}$ for a general deterministic diagonal matrix $A$. Although we believe that such an extension could be obtained by combining the methods of the present paper with those developed in \cite{CEDS_EJP,CES_Forum,erdHos2024eigenstate}, we do not pursue this direction here due to space limitations.

	As another important consequence of the QUE estimates in Theorem \ref{theorem_que}, combined with the Green's function comparison argument developed in \cite{xu2024bulk}, we obtain \emph{universality of the local bulk eigenvalue statistics} in the complete delocalization regime. Specifically, for power-law random band matrices, the local eigenvalue statistics near any fixed bulk energy $E$ coincide asymptotically with those of Wigner matrices.
	Let $p_H(\lambda_1,\ldots,\lambda_N)$ denote the joint symmetrized probability density of the (unordered) eigenvalues of $H$. For $1\le n\le N$, define the $n$-point correlation function by
	$$
	p_{H}^{(n)}\left(\lambda_1, \cdots, \lambda_n\right)
	:=\int_{\R^{N-n}} p_{H}\left(\lambda_1, \cdots, \lambda_{N}\right) \mathrm{d} \lambda_{n+1} \cdots \mathrm{d} \lambda_{N}.
	$$
	Let $p_{\mathrm{GUE}}^{(n)}$ denote the corresponding $n$-point correlation function for an $N\times N$ GUE matrix.

	\begin{theorem}[Bulk universality]\label{theorem_local_statistics}
		In the setting of Theorem \ref{theorem_que}, let $O \in C_c^{\infty}\left(\mathbb{R}^n\right)$ be a smooth, compactly supported test function. Then, for any $|E|\le 2-\kappa$ and fixed $n\in \N$, we have
		\begin{equation}\label{eq:universality}
			\lim_{N\to \infty}	\int_{\mathbb{R}^n} \mathrm{~d} \boldsymbol{\alpha}\; O(\boldsymbol{\alpha}) \left[p_{H}^{(n)}-p_{\rm{GUE}}^{(n)}\right]\left(E+\frac{\alpha_1}{N  }, \ldots, E+\frac{\alpha_n}{N}\right)   =0, 
		\end{equation} 
		where $\bal$ denotes $\boldsymbol{\alpha}=\left(\alpha_1, \ldots, \alpha_n\right)$. 
	\end{theorem}

	A key ingredient in the proof of the QUE estimate in \Cref{theorem_que} is the following \emph{quantum diffusion estimate}, stated in \Cref{theorem_local_law_loops}. It roughly asserts that, for $|x-y|\gg W$ and $\eta\ll 1$, the quantity $\E|G_{xy}|^2$ can be approximated by the Green's function of an associated classical random walk. To formulate this result, we introduce the so-called \emph{$\mathcal L$-loops} and their deterministic approximations. 
	For $\boldsymbol{\sigma}=(\sigma,\sigma')\in\{-,+\}^2$ and $x,y\in\mathbb Z_N$, define
	\begin{equation}\label{eq:def_Lloop}
		\cL_{xy}^{\bsigma}(z):=\sum_{a,b\in\ZL}S_{xa}S_{yb}G(z_{\sigma})_{ba}G(z_{\sigma'})_{ab},\quad \cK_{xy}^{\bsigma}(z):= \pa{\frac{m(z_{\sigma})m(z_{\sigma'})S^2}{1-m(z_{\sigma})m(z_{\sigma'})S}}_{xy},
	\end{equation}
	where we set $z_{+}:=z$ and $z_{-}=\ol{z}$.
	Note that the \emph{$\mathcal L$-loops} represent local averages of $G(z_{\sigma})_{ba}G(z_{\sigma'})_{ab}$, weighted by the entries of the variance matrix $S$. On the other hand, the $\mathcal K$-loops describe the effective “diffusion profile” of the random walk on $\ZL$ whose transition probabilities are given by $S$.
	In addition to the spectral domain defined in \eqref{def_spectral_domain}, we introduce the \emph{flat spectral regime}
	\begin{equation}\label{eq:Dflat}
		\begin{aligned}
			\mathbf{D}_{\kappa,\fc}^{\txt{flat}}:=\ha{z=E+\ii \eta\in \mathbf{D}_{\kappa,\fc}: |E|\leq 2-\kappa,\, \eta\leq \eta_{\txt{flat}}},
		\end{aligned}
	\end{equation}
	where we denote
	\begin{equation}\label{eq:etaflat}
		\eta_{\txt{flat}}:=\begin{cases}
			\pa{W/N}^{1+\al}, \ &\text{if}\ \alpha\in(-1,0)\\
			W/N, \ &\text{if}\ \alpha\in [0,1)\\
			\pa{W/N}^{\al\wedge 2}, \ &\text{if}\ \alpha\ge 1
		\end{cases}. \end{equation}
	For $z=E+\mathrm i\eta\in\mathbf D_{\kappa,\fc}^{\mathrm{flat}}$, the associated shape parameter is flat: since $1-|m(z)|^2\asymp \eta \leq \eta_{\txt{flat}}$, we obtain
	\begin{equation*}
		\begin{aligned}
			B_{|m(z)|^2}(r)\prec \p{N\eta}^{-1},\quad \forall 0\le r\le N.
		\end{aligned}
	\end{equation*}
	
	\begin{theorem}[Quantum diffusion]\label{theorem_local_law_loops}
		
		Under the assumptions of \Cref{theorem_delocalization}, for any fixed constants $\tau,D>0$, there exists $N_0(\tau,D)$ such that, for all $N\ge N_0$, the corresponding events occur with probability $\ge 1-N^{-D}$.
		\begin{itemize}
			\item[(i)] {\bf The case $\al\ge 1$.} We have
			\begin{align}
				&\bigcap_{z\in\mathbf{D}_{\kappa,\fc}}\bigcap_{x,y\in\ZL}\ha{\max_{\bsigma\in\ha{-,+}^2}\absa{(\cL-\cK)_{xy}^{\bsigma}(z)}\leq N^{\tau}B_{|m|^2}(0)B_{|m|^2}(|x-y|)},\label{eq_loop_local_law_1}
			\end{align}
			where we abbreviate $m\equiv m(z)$, and $B_{|m|^2}(\cdot )=B(1-|m|^2,\cdot)$ is defined in the sense of \eqref{eq:abbrvB}. 
			
			\item[(ii)]{ \bf The case $\alpha\in[0,1)$.} We have
			\begin{align}
				&\bigcap_{z\in\mathbf{D}_{\kappa,\fc}}\bigcap_{x,y\in\ZL}\ha{\max_{\bsigma\in\ha{-,+}^2}\absa{(\cL-\cK)_{xy}^{\bsigma}(z)}\leq N^{\tau}[B_{|m|^2}(0)]^{1/5}B_{|m|^2}(|x-y|)}.\label{eq_loop_local_law_2}
			\end{align}
			Moreover, in the flat regime, this estimate improves to
			\begin{align}    	&\bigcap_{z\in\mathbf{D}_{\kappa,\fc}^{\txt{flat}}}\bigcap_{x,y\in\ZL}\ha{\max_{\bsigma\in\ha{-,+}^2}\absa{(\cL-\cK)_{xy}^{\bsigma}(z)}\leq N^{\tau}\pa{W^{-6/5}+\pa{N\im z}^{-7/4}}}.\label{eq_loop_local_law_2_improved}
			\end{align}    
			\item[(iii)]{ \bf The case $\alpha\in(-1,0)$.} In the flat regime, we have\footnote{A more detailed inspection of the proof in the regime $\alpha\in(-1,0)$ shows that these bounds can be extended to the full domain $\mathbf D_{\kappa,\fc}$. However, outside the flat regime the resulting estimates take a more intricate form depending explicitly on $\alpha$, without yielding any further improvement for the QUE statement. For brevity, we therefore restrict to the current formulation.}     \begin{align}    	&\bigcap_{z\in\mathbf{D}_{\kappa,\fc}^{\txt{flat}}}\bigcap_{x,y\in\ZL}\ha{\max_{\bsigma\in\ha{-,+}^2}\absa{(\cL-\cK)_{xy}^{\bsigma}(z)}\leq N^{\tau}\pa{N\im z}^{-7/4}}.\label{eq_loop_local_law_3}    	    \end{align}

		\end{itemize}

		The preceding high-probability estimates admit stronger bounds at the level of expectations. More precisely, for any small constant $\tau>0$ and sufficiently large $N$, the following estimates hold. If $\alpha\ge 1$, then
		\begin{equation}\label{eq_expected_local_law_1}
			\begin{aligned}
				\max_{\bsigma\in\ha{-,+}^2}\absa{\E(\cL-\cK)_{xy}^{\bsigma}(z)}\leq N^{\tau}[B_{|m|^2}(0)]^{3/2}B_{|m|^2}(|x-y|)
			\end{aligned}
		\end{equation}
		for each $z\in\mathbf D_{\kappa,\fc}$ and $x,y\in \ZL$.
		If $\alpha\in(-1,1)$, then for each $z\in\mathbf D_{\kappa,\fc}^{\mathrm{flat}}$ and $x,y\in \ZL$,
        \begin{equation}\label{eq_expected_loop_local_law}     
			\max_{\bsigma\in\ha{-,+}^2}\absa{\E(\cL-\cK)_{xy}^{\bsigma}(z)}\leq N^{\tau}\times \begin{cases}W^{-7/4}(N/W)^{7\alpha/4} +\pa{N\im z}^{-9/4},\ & \text{if}\ \alpha\in(-1,0)\\
			    W^{-6/5}+\pa{N\im z}^{-9/4},\ & \text{if}\ \alpha\in [0,1)
			\end{cases}.
		\end{equation}        
\end{theorem}

	A byproduct of the proof of the above quantum diffusion estimates is a local law for the \emph{$T$-variables}, together with their deterministic approximations $\Theta$. For $\boldsymbol{\sigma}=(\sigma,\sigma')\in\{-,+\}^2$ and $x,y,y'\in\mathbb Z_N$, define
	\begin{equation}\label{eq:def_Tchain}
		T_{x,yy'}^{\bsigma}(z):=\sum_{a\in\ZL}S_{xa}G(z_{\sigma})_{ya}G(z_{\sigma'})_{ay'},\quad \Theta_{x,yy'}^{\bsigma}(z):= \delta_{yy'}\pa{\frac{m(z_{\sigma})m(z_{\sigma'})S}{1-m(z_{\sigma})m(z_{\sigma'})S}}_{xy}.
	\end{equation}
	When $y\neq y'$, we refer to $T_{x,yy'}^{\boldsymbol{\sigma}}$ as an \emph{off-diagonal} $T$-variable; otherwise, it is called a \emph{diagonal} $T$-variable. We also refer to $\Theta^{\boldsymbol{\sigma}}$ as the \emph{$\Theta$-propagators}.

	\begin{lemma}[$T$-lemma]\label{theorem_local_law_T}
		Under the assumptions of \Cref{theorem_delocalization}, the following statements hold. For any fixed constants $\tau,D>0$, there exists $N_0(\tau,D)$ such that, for all $N\ge N_0$, the corresponding events occur with probability at least $1-N^{-D}$.
		\begin{itemize}
			\item[(i)] {\bf The case $\al\ge 1$.} We have 
			\begin{align}        &\bigcap_{z\in\mathbf{D}_{\kappa,\fc}}\bigcap_{x,y,y'\in\ZL}\ha{\max_{\bsigma\in\ha{-,+}^2}\absa{(T-\Theta)_{x,yy'}^{\bsigma}(z)}\leq N^{\tau}\qa{B_{|m|^2}(0)}^{1/2}\qa{B_{|m|^2}(|x-y|)B_{|m|^2}(|x-y'|)}^{1/2}}.\label{eq_chain_local_law_1}
			\end{align}
			
			\item[(ii)] {\bf The case $\alpha\in[0,1)$.} We have
			\begin{align}    	&\bigcap_{z\in\mathbf{D}_{\kappa,\fc}}\bigcap_{x,y,y'\in\ZL}\ha{\max_{\bsigma\in\ha{-,+}^2}\absa{(T-\Theta)_{x,yy'}^{\bsigma}(z)}\leq N^{\tau}\qa{B_{|m|^2}(0)}^{7/10}\qa{B_{|m|^2}(|x-y|\wedge |x-y'|)}^{1/2}}.\label{eq_chain_local_law_2}
			\end{align}
			Moreover, in the flat regime, this estimate improves to
			\begin{align}    	&\bigcap_{z\in\mathbf{D}_{\kappa,\fc}^{\txt{flat}}}\bigcap_{x,y,y'\in\ZL}\ha{\max_{\bsigma\in\ha{-,+}^2}\absa{(T-\Theta)_{x,yy'}^{\bsigma}(z)}\leq N^{\tau}\pa{W^{-6/5}+\pa{N\im z}^{-3/2}}}.\label{eq_chain_local_law_2_improved}
			\end{align}
			
			\item[(iii)] {\bf The case $\alpha\in(-1,0)$.} In the flat regime, we have
			\begin{align}
				&\bigcap_{z\in\mathbf{D}_{\kappa,\fc}^{\txt{flat}}}\bigcap_{x,y,y'\in\ZL}\ha{\max_{\bsigma\in\ha{-,+}^2}\absa{(T-\Theta)_{x,yy'}^{\bsigma}(z)}\leq N^{\tau}\pa{N\im z}^{-3/2}}.\label{eq_chain_local_law_3}
			\end{align}
		\end{itemize}
	\end{lemma}
	
	In addition to being of independent interest, the above $T$-variable estimate constitutes a key and novel technical ingredient in the proofs of the local laws in \Cref{theorem_local_law_entrywise} and the quantum diffusion estimates in \Cref{theorem_local_law_loops}, as discussed in \Cref{subsec:ideas}. In particular, it partially replaces the role of the high-order $G$-loop analysis (see \Cref{def_n_loops_and_chains} for its definition) employed in earlier works \cite{Band1D,Band2D,erdos2025zigzagstrategyrandomband,truong2025localizationlengthfinitevolumerandom,dubova2025delocalizationnonmeanfieldrandommatrices,Bandedge}. 
	
	\subsection{Quantum unique ergodicity and bulk universality}\label{sec_proof_of_local_statistics}
	
	In this subsection, we prove the QUE estimates using the quantum diffusion estimates in Theorem \ref{theorem_local_law_loops}. These estimates further yield bulk universality of the local eigenvalue statistics by adopting the strategy of \cite{xu2024bulk}. Such arguments are standard and have been used for regular random band matrix models; see, for example, \cite{Band1D,Band2D,erdos2025zigzagstrategyrandomband,dubova2025delocalizationnonmeanfieldrandommatrices}. Therefore, we only outline the proof.
	
	First, given the QUE estimate \eqref{eq_que}, the proof of Theorem \ref{theorem_local_statistics} is almost identical to that in \cite[Section 2.3]{Band1D} (which assumes a block variance profile) and \cite[Section 15]{erdos2025zigzagstrategyrandomband} for 1D regular random band matrices, provided that we can establish the following lemma.
	
	\begin{lemma}
		Under the assumptions of \Cref{theorem_local_statistics}, let $\lambda\in[0,N^{-1+\tau}]$ for an arbitrarily small constant $\tau\in (0,1/100)$. Define the matrix $H^{\lambda}$ by
		\begin{equation*}
			\begin{aligned}
				H^{\lambda}:=\sqrt{1-\lambda}\, H+\sqrt{\lambda}\, H_{\txt{GUE}},
			\end{aligned}
		\end{equation*}
		where $H_{\txt{GUE}}$ is an $N\times N$ GUE matrix independent of $H$, whose entries have mean zero and variance $N^{-1}$. Then the main results, namely Theorems \ref{theorem_delocalization}--\ref{theorem_que} and \ref{theorem_local_law_loops}, continue to hold for $H^{\lambda}$.
	\end{lemma}
	\begin{proof} 
		Note that the matrix $H^{\lambda}$ satisfies essentially the same assumptions as $H$, except that its variance profile $S^\lambda$ is now given by \( S^\lambda_{xy} =\p{1-\lambda}S_{xy}+\lambda/N\), where $S_{xy}$ satisfies \eqref{alpha_decay}. Moreover, after a slight modification of the characteristic lengths and shape parameters, all the estimates required in \Cref{assumption_input_bound} continue to hold with $S$ replaced by $S^\lambda$. This follows from the identity
		\begin{equation*}
			\frac{1}{1-t\xi S^\lambda}=(1-P)\frac{1}{1-t(1-\lambda)\xi S}(1-P)+\frac{1}{1-t\xi}P,  
		\end{equation*}
		where $P:=\mathbf e\mathbf e^*$ denotes the projection onto the vector $\mathbf e=N^{-1/2}(1,\ldots,1)^*$. Using these modified estimates as input and repeating the arguments in the subsequent sections therefore yields Theorems \ref{theorem_delocalization}--\ref{theorem_que} and \ref{theorem_local_law_loops} for each fixed $\lambda\in[0,N^{-1+\tau}]$. Finally, a standard $N^{-C}$-net and perturbation argument extends these estimates uniformly to all $\lambda\in[0,N^{-1+\tau}]$. We omit the details.
	\end{proof}

	To prove \Cref{theorem_que}, we set $z=E+\ii \eta$ with $\eta= \eta(\varepsilon_0)$. By the spectral decomposition of $H$, we obtain
	\begin{align}        \E\sum_{i,j:\lambda_i,\lambda_j\in\cI_E}\absa{\bpsi_i^*\pb{S^{(x)}-N^{-1}}\bpsi_j}^2&\lesssim \eta^2\E\mathrm{Tr}\qb{\im G(z)\pb{S^{(x)}-N^{-1}}\im G(z)\pb{S^{(x)}-N^{-1}}}\nonumber\\
		&= \frac{\eta^2}{N^2}\sum_{y,y'\in\ZL}\E\mathrm{Tr}\qb{\im G(z)\pb{S^{(x)}-S^{(y)}}\im G(z)\pb{S^{(x)}-S^{(y')}}},\label{proof_que_spectral_decomposition_bound}
	\end{align}
	where we recall that $S^{(x)}$ denotes the diagonal matrix defined by 
	\begin{equation}\label{eq:def_Sx}
		S^{(x)}_{ab}:=\delta_{ab}S_{xa},\quad \forall a,b\in \ZL.
	\end{equation} 
	Since $\eta\ge N^{-1+\fc_0-\varepsilon_0}$, for any small constant $\tau>0$, applying the expected quantum diffusion estimates \eqref{eq_expected_local_law_1} and \eqref{eq_expected_loop_local_law}, we can bound the right-hand side (RHS) of \eqref{proof_que_spectral_decomposition_bound} by
	\begin{equation}\label{proof_que_spectral_decomposition_bound2}
		N^{\tau/2}\cdot\eta^2\cE_{\E}+\eta^2\max_{\bsigma \in\ha{-,+}^2}\max_{x,y,y'\in\ZL}\absa{\cK_{xy}^{\bsigma}(z)-\cK_{xy'}^{\bsigma}(z)}\leq N^{\tau}\eta^2\pa{\cE_{\E}+\cE_{\cK}}.
	\end{equation}
	Here, $\cE_{\E}$ and $\cE_{\cK}$ denote the errors arising from the expected quantum diffusion estimates and from the flatness of $\cK$, respectively:
	\begin{equation*}
		\begin{aligned}
			\cE_{\E}:=\begin{cases}
				W^{-7(1+\alpha)/4}N^{7\alpha/4}+(N\eta)^{-9/4},\ &\text{if}\ \alpha\in (-1,0)\\
				W^{-6/5}+(N\eta)^{-9/4},\ &\text{if}\ \alpha\in [0,1)\\
				(N\eta)^{-5/2},\ &\text{if}\ \alpha\in[1,\infty)
			\end{cases} ,
			\quad
			\cE_{\cK}:=\begin{cases}
				W^{-1-\alpha}N^{\alpha},\ &\text{if}\ \alpha\in (-1,0)\\
				W^{-1},\ &\text{if}\ \alpha\in [0,1)\\
				N^{-1}(N/W)^{\alpha\wedge 2},\ &\text{if}\ \alpha\in[1,\infty)
			\end{cases}.
		\end{aligned}
	\end{equation*}
	In deriving the error $\cE_{\cK}$, we used the assumptions \eqref{Theta_bound_equaled_charge} and \eqref{zero_mode_removed_bound} for $\al\ge 0$ together with the bound \eqref{zero_mode_removed_bound_subcritical} below for $\al\in(-1,0)$. From \eqref{proof_que_spectral_decomposition_bound} and \eqref{proof_que_spectral_decomposition_bound2}, it follows that for any small constant $\tau>0$, \begin{equation}\label{eq_expectation_que}
		\E\sum_{i,j:\lambda_i,\lambda_j\in\cI_E}\absa{\bpsi_i^*\pb{S^{(x)}-N^{-1}}\bpsi_j}^2\leq N^{-2-\fc_1+\tau}.
	\end{equation}
	The estimate \eqref{eq_que} then follows immediately from Markov's inequality.
	
	Finally, we prove Corollary \ref{cor_que}. For any $\bpsi_k$ with $\lambda_k\in\cI_{E}$ and $|E|\leq 2-\kappa$, we write
	\begin{align}
		&\frac{1}{|I|}\sum_{a\in I}\absa{\bpsi_k(a)}^2=\frac{1}{|I|}\sum_{a\in I}\sum_{b\in\ZL}S_{ab}\absa{\bpsi_k(a)}^2\nonumber\\
		=&\frac{1}{|I|}\sum_{b\in I}\sum_{a\in \ZL}S_{ab}\absa{\bpsi_k(a)}^2+\frac{1}{|I|}\sum_{a\in I}\sum_{b\in \ZL\setminus I}S_{ab}\absa{\bpsi_k(a)}^2-\frac{1}{|I|}\sum_{a\in \ZL\setminus I}\sum_{b\in I}S_{ab}\absa{\bpsi_k(a)}^2.\label{eq_decomposition_I_que}
	\end{align}
	On the event $\ha{\lambda_k\in\cI_{E}}$, we bound the last two terms on the right-hand side of \eqref{eq_decomposition_I_que} using the delocalization estimate \eqref{eq:delocal2} together with the assumption \eqref{eq_que_I_assumption}. This yields, for some constant $c>0$,
	\begin{equation}\label{eq:Ireduce}
		\frac{1}{|I|}\sum_{a\in I}\absa{\bpsi_k(a)}^2=\frac{1}{|I|}\sum_{b\in I}\sum_{a\in \ZL}S_{ab}\absa{\bpsi_k(a)}^2+\opr{N^{-1-c}}.
	\end{equation}
	From \eqref{eq_expectation_que}, we have
	\[ \E\max_{\lambda_k\in\cI_E}\frac{1}{|I|}\sum_{b\in I}\absa{\sum_{a\in \ZL}S_{ab}\absa{\bpsi_k(a)}^2-N^{-1}}^2\leq N^{-2-\fc_1+\tau}. \]
	Combining this bound with Markov's inequality, we obtain the following estimate for some constant $c>0$:
	\begin{equation*}
		\begin{aligned}
			\P\pa{\max_{\lambda_k\in \cI_{E}}\absa{\frac{1}{|I|}\sum_{b\in I}\sum_{a\in \ZL}S_{ab}\absa{\bpsi_k(a)}^2-N^{-1}}\geq N^{-1-c}}\leq N^{-c}.
		\end{aligned}
	\end{equation*}
	Together with \eqref{eq:Ireduce}, this completes the proof of \eqref{eq_flat_que}.

	\section{Proof of the main results for \texorpdfstring{$\al\in(0,1)$}{alpha(0,1)}}\label{sec:pfal01}
	
	In this section, we prove the main results, \Cref{theorem_local_law_entrywise,theorem_local_law_loops} and \Cref{theorem_local_law_T}, in the most delicate regime $\alpha\in(0,1)$. The argument relies on several estimates that will be established in subsequent sections.

	\subsection{Notation and tools}\label{sec:notations}

	In this subsection, we collect several tools and notational conventions that will be used throughout the proof. We begin by introducing the flow framework. \begin{definition}[Matrix Brownian motion]\label{def_H_t_flow}
		In the setting of \Cref{theorem_delocalization}, let $S$ denote the variance matrix of $H$. Let $\pa{\mathbf{B}_{xy}(t):x,y\in\ZL}$ be a family of independent complex Brownian motions subject to the Hermitian symmetry \smash{$\mathbf{B}_{xy}(t)=\overline{\mathbf{B}_{yx}(t)}$}. Then for each $t>0$, $t^{-1/2}\mathbf{B}_t$ is an $N\times N$ GUE matrix with zero mean and unit variance entries. We define the matrix Brownian motion $H_t$ by
		\begin{equation}\label{MBM}
			\begin{aligned}
				\rd (H_t)_{xy}=\sqrt{S_{xy}}\,\rd \mathbf{B}_{xy}\p{t}, \quad \forall x,y\in\ZL.
			\end{aligned}
		\end{equation}
		We denote the resolvent of $H_{t}$ by $G_t(z):=(H_t-z)^{-1}$ for $z\in \C_+$.
	\end{definition}

	As in \cite{10.1214/19-ECP278, Sooster2019, DY, Band1D}, it is convenient to study the resolvent of $H_t$ at a suitably chosen time-dependent spectral parameter $z_t$, whose dynamics are naturally normalized at leading order.
	
	\begin{definition}[Renormalized flow]\label{def_z_t_flow}
		In the setting of \Cref{theorem_delocalization}, recall the definition of $m$ from \eqref{def_m_sc}. Fix a flow parameter $\mathsf{E}\in  (-2,2)$. The \emph{renormalized spectral parameter flow} is defined by
		\begin{equation}\label{eq:zt}
			\begin{aligned}
				z_t\equiv z_t(\sE):=\mathsf{E}+(1-t)m\p{\mathsf{E}},\quad t\in\qa{0,1},
			\end{aligned}
		\end{equation}
		where $m(\mathsf{E})$ denotes the boundary value $m(\mathsf{E}+\mathrm{i}0_+)$ obtained by continuously extending $m$ from $\mathbb{C}_+$ to $\overline{\mathbb{C}}_+$. We write $z_t = E_t + \mathrm{i}\eta_t$, where
		\begin{equation}\label{re_im_z_t}
			\begin{aligned}
				E_t\equiv E_t\p{\mathsf{E}}:=\mathsf{E}+(1-t)\txt{Re}\,m\p{\mathsf{E}},\quad \eta_t\equiv \eta_t\p{\mathsf{E}}:=(1-t)\im m\p{\mathsf{E}}.
			\end{aligned}
		\end{equation} 
		Then, the corresponding \emph{renormalized resolvent flow} is defined by 
		\begin{equation}\label{self_Gt}
			\begin{aligned}
				G_t\equiv G_{t,\mathsf{E}}:=G_t(z_t(\sE))=\pa{H_t-z_t\p{\mathsf{E}}}^{-1}.
			\end{aligned}
		\end{equation}
	\end{definition}
	
	To illustrate a key advantage of the flow \eqref{eq:zt}, let $m_t(z)$ denote the unique solution to the equation 
	\begin{equation}\label{eq:mtz}
		m_t(z)=-(z+tm_t(z))^{-1},\quad  \text{with}\quad  \im m_t(z)>0\ \ \text{for}\ \ z\in \C_+.
	\end{equation} 
	It is known that $m_t(z) I_N$ describes the deterministic entrywise limit of the resolvent $G_t(z)$. Combining this equation with \eqref{eq:self_cons_m}, one verifies that
	\begin{equation}\label{eq:mE}
		m_t(z_t(\sE))\equiv m(\sE),\quad \forall t\in [0,1], 
	\end{equation} 
	that is, the deterministic limit of $G_{t,\mathsf{E}}$ remains invariant along the renormalized flow. 
	
	For a given target spectral parameter $z\in\mathbf{D}_{\kappa,\mathfrak c}$, our goal is to analyze the original resolvent $G(z)=(H-z)^{-1}$. This can be achieved via the renormalized flow by an appropriate choice of the flow parameter $\mathsf{E}$, as stated in the following lemma.
	
	\begin{lemma}[Lemma 2.8 of \cite{Band1D}]\label{zztE}
		Under the assumptions of \Cref{theorem_delocalization}, fix any $z\in \mathbf{D}_{\kappa,\fc}$. Define
		\be\label{eq:t0E0}\tf\equiv \tf(z)=|m(z)|^2=\frac{\im m(z)}{\im m(z)+ \im z},\quad \sE\equiv \sE(z)=-2\frac{\re m(z)}{|m(z)|}\, .\ee
		Then, for the power-law RBM model, we have 
		\begin{equation}\label{eq:zztE}
			\sqrt{\tf}m(\sE)=m(z), \quad     z_{\tf}(\sE) =\sqrt{\tf} z  , \quad G(z) \stackrel{d}{=} \sqrt{\tf} G_{\tf,\sE} ,
		\end{equation} 
		where ``$\stackrel{d}{=}$'' means equality in distribution.
	\end{lemma}

	In the subsequent analysis, we fix a target spectral parameter $z=E+\mathrm{i}\eta\in \mathbf{D}_{\kappa,\mathfrak c}$. Accordingly, we choose the time parameter $\tf$ and the flow parameter $\mathsf{E}$ as in \eqref{eq:t0E0}. From the second identity in \eqref{re_im_z_t}, we observe that $1-t\asymp \eta_t$ uniformly in $t\in [0,\tf]$. Hence, along the flow from $t=0$ to $\tf$, the imaginary part $\eta_t$ decreases from $\eta_0\asymp 1$ to $\eta_{\tf}\asymp 1-\tf \asymp \eta \ge N^{\fc}\eta_*$.  
	For brevity, unless the dependence on $\mathsf{E}$ needs to be emphasized, we suppress this parameter in the notation and write simply $z_t$, $E_t$, $\eta_t$, $m$, and, in particular, $G_t \equiv G_{t,\mathsf{E}}$. Our primary focus is on the evolution of $G_t$ and the associated $\mathcal L$-loops and $T$-variables (see \eqref{eq:def_Lloop} and \eqref{eq:def_Tchain}).

	\begin{definition}[$\cL$-loops and $T$-variables]\label{def_loops_and_T_variables}
		Within the flow framework of Definition \ref{def_z_t_flow}, for any $t\in[0,1)$, $x,y,y'\in\ZL$, and charge $\boldsymbol\sigma=(\sigma_1,\sigma_2)\in\{-,+\}^2$, we define
		\begin{equation}\label{eq:cL_and_T}            \cL_{t,xy}^{\bsigma}:=\avg{G_t\p{\sigma_1}S^{\pa{x}}G_t\p{\sigma_2}S^{\pa{y}}},\quad T_{t,x,yy'}^{\bsigma}:=\p{G_t\p{\sigma_1}S^{\pa{x}}G_t\p{\sigma_2}}_{yy'},
		\end{equation}
		where $G_t(+)\equiv G_t^+:=G_t$, $G_t(-)\equiv G_t^-:=G_t^*$, and $S^{(x)}$ is defined in \eqref{eq:def_Sx}.
		When no ambiguity arises, we use the abbreviations
		\begin{equation}\label{abbreviation_cL_T}
			\begin{aligned}
			\cL_{t,xy}\equiv \cL_{t,xy}^{(-,+)},\quad T_{t,xy}^{\bsigma}\equiv T_{t,x,yy}^{\bsigma},\quad T_{t,x,yy'}\equiv T_{t,x,yy'}^{(-,+)},\quad T_{t,xy}\equiv T_{t,xy}^{(-,+)}.
			\end{aligned}
		\end{equation}
		As deterministic approximations to the $\mathcal L$-loops and $T$-variables, we introduce the corresponding $\mathcal K$-loops and $\Theta$-propagators:
		\begin{equation}\label{eq:KTheta}
			\begin{aligned}
				\cK_{t,xy}^{\bsigma}:=\pa{\frac{m\p{\sigma_1}m\p{\sigma_2}S^2}{1-tm(\sigma_1)m(\sigma_2)S}}_{xy},\quad \Theta_{t,x,yy'}^{\bsigma}:=\delta_{yy'}\pa{\frac{m\p{\sigma_1}m\p{\sigma_2}S}{1-tm(\sigma_1)m(\sigma_2)S}}_{xy},
			\end{aligned}
		\end{equation}
		where $m(+)\equiv m^+:=m(\mathsf{E})$ and $m(-)\equiv m^-:=\overline m(\mathsf{E})$. Similar abbreviations for $\cK$ and $\Theta$ to those introduced in \eqref{abbreviation_cL_T} will be used when no confusion arises. The charge $\boldsymbol\sigma$ will be written either as a superscript or as an argument, depending on convenience.
	\end{definition}

	Let $\cal A$ be a Hermitian matrix, and define its resolvent by $R(z):=(\cal A-z)^{-1}$ for $z= E+ \ii \eta\in \C_+$. Using the algebraic identity $\im R:=(R-R^*)/(2\ii)= \eta RR^*= \eta R^*R$, we get the well-known Ward's identity:
	\be\label{eq_Ward0}
	\begin{split}
		\sum_x \overline {R_{xy}}  R_{xy'} = \frac{(\im R)_{yy'}}{\eta},\quad
		\sum_x R_{yx} \overline {R_{y'x}}  = \frac{(\im R)_{yy'}}{\eta} . \end{split}
	\ee
	Since $S$ is doubly stochastic, it follows from \eqref{eq_Ward0} that, for $\sigma_1\ne \sigma_2$, the following Ward's identities hold for the $\cL$-loops and the $T$-variables:
	\be\label{eq_Ward}
	\sum_y \cL^{\bsigma}_{t,xy} = \sum_y \cL^{\bsigma}_{t,yx} = \frac{\im \avg{G_t S^{(x)}} }{ \eta_t},\quad \sum_x T_{t,x,yy'}^{\bsigma}=\frac{(\im G_t)_{yy'}}{\eta_t},\quad \sum_y T_{t,xy}^{\bsigma}=\frac{\im \avg{G_t S^{(x)}} }{ \eta_t}.
	\ee
    We will use the identities \eqref{eq_Ward0} and \eqref{eq_Ward} tacitly throughout the sequel.

	Using \eqref{eq:zt}, \eqref{self_Gt}, \eqref{eq:mE}, and \eqref{eq:KTheta}, and applying Itô's formula, we obtain the evolution equations summarized in the following lemma. Since the derivation is standard, we omit the details.

	\begin{lemma}\label{lemma_evolution_equations}
		Under the notation of \Cref{def_loops_and_T_variables}, the following evolution equations hold for $\boldsymbol{\sigma}=(\sigma_1,\sigma_2)\in\{-,+\}^2$, $t\in[0,1)$, and $x,y,y'\in\mathbb Z_N$.
		\begin{itemize}
			\item For $\cL$-loops, we have
			\begin{equation}\label{cL_evolution}
				\begin{aligned}
					\rd \cL_{t,xy}^{\bsigma}=\sum_{a\in\ZL}T_{t,xa}^{\bsigma}\cL_{t,a y}^{\bsigma}\dd t+\cW_{t,xy}^{\cL,\bsigma}\, \rd t+\rd \cB_{t,xy}^{\cL,\bsigma},
				\end{aligned}
			\end{equation}
			where the light-weight\footnote{We refer to a factor of the form $G_{xx}^\sigma-m^\sigma$, for $x\in\ZL$ and $\sigma\in\ha{-,+}$, as a light-weight.} term is given by
			\begin{align}
				\cW_{t,xy}^{\cL,\bsigma}:=&~\sum_{a\in\ZL}\pa{G_t^{\sigma_2}-m^{\sigma_2}}_{aa}\avg{G_t^{\sigma_2}S^{\pa{a}}G_t^{\sigma_2}S^{\pa{y}}G_t^{\sigma_1}S^{\pa{x}}}\nonumber\\
				+&~\sum_{a\in\ZL}\pa{G_t^{\sigma_1}-m^{\sigma_1}}_{aa}\avg{G_t^{\sigma_1}S^{\pa{a}}G_t^{\sigma_1}S^{\pa{x}}G_t^{\sigma_2}S^{\pa{y}}},\label{def_light_weight_term}
			\end{align}
			and the martingale term is given by
			\begin{equation}\label{eq:BLt}
				\rd\cB_{t,xy}^{\cL,\bsigma}=\sum_{a,b\in\ZL}(\partial_{ab}\cL_{t,xy}^{\bsigma})\sqrt{S_{ab}}\,\rd \mathbf{B}_{ab}\p{t}.
			\end{equation}
			Here, for $(a,b)\in \ZL^2$, we write $\partial_{ab}\equiv \partial_{(H_t)_{ab}}$.
			
			\item For $T$-variables, we have
			\begin{equation}\label{T_evolution}
				\begin{aligned}
					\rd T_{t,x,yy'}^{\bsigma}=\sum_{a\in\ZL}T_{t,xa}^{\bsigma}T_{t,a, yy'}^{\bsigma}\,\rd t+\cW_{t,x,yy'}^{T,\bsigma}\, \rd t+\rd \cB_{t,x,yy'}^{T,\bsigma},
				\end{aligned}
			\end{equation}
			where the light-weight term is given by
			\begin{align}
				\cW_{t,x,yy'}^{T,\bsigma}:=&~\sum_{a\in\ZL}\pa{G_t^{\sigma_2}-m^{\sigma_2}}_{aa}\p{G_t^{\sigma_1}S^{\pa{x}}G_t^{\sigma_2}S^{\pa{a}}G_t^{\sigma_2}}_{yy'}\nonumber\\
				+&~\sum_{a\in\ZL}\pa{G_t^{\sigma_1}-m^{\sigma_1}}_{aa}\p{G_t^{\sigma_1}S^{\pa{a}}G_t^{\sigma_1}S^{\pa{x}}G_t^{\sigma_2}}_{yy'},\label{def_light_weight_term_T}
			\end{align}
			and the martingale term is given by
			\begin{equation}\label{eq:BTt}
				\rd\cB_{t,x,yy'}^{T,\bsigma}=\sum_{a,b\in\ZL}(\partial_{ab}T_{t,x,yy'}^{\bsigma})\sqrt{S_{ab}}\,\rd \mathbf{B}_{ab}\p{t}.
			\end{equation}
			\item For $\cK$-loops, we have
			\begin{equation}\label{cK_evolution}
				\partial_t \cK_{t,xy}^{\bsigma}=\sum_{a\in\ZL}\Theta_{t,xa}^{\bsigma}\cK_{t,a y}^{\bsigma}=\sum_{a\in\ZL}\cK_{t,xa}^{\bsigma}\Theta_{t,a y}^{\bsigma}.
			\end{equation}
			\item For $\Theta$-propagators, we have
			\begin{equation}\label{Theta_evolution}
				\partial_t \Theta_{t,xy}^{\bsigma}=\sum_{a\in\ZL}\Theta_{t,xa}^{\bsigma}\Theta_{t,a y}^{\bsigma}.
			\end{equation}
		\end{itemize}
	\end{lemma}
	
	Subtracting \eqref{cK_evolution} from \eqref{cL_evolution}, and using the identity
	\[ \sum_a (\cL-\cK)_{t,xa}^{\bsigma}\Theta_{t,a y}^{\bsigma} = \sum_{a,b} (T-\Theta)_{t,xb}^{\bsigma}S_{ba}\Theta_{t,a y}^{\bsigma}= \sum_{b} (T-\Theta)_{t,xb}^{\bsigma}\cK_{t,b y}^{\bsigma},\]
	we obtain the following evolution equation for the difference $(\mathcal L-\mathcal K)$:
	\begin{equation}\label{cL_cK_evolution}
		\rd (\cL-\cK)_{t,xy}^{\bsigma}=\sum_{a}\Theta_{t,xa}^{\bsigma}(\cL-\cK)_{t,a y}^{\bsigma}\,\rd t+\sum_a (\cL-\cK)_{t,xa}^{\bsigma}\Theta_{t,a y}^{\bsigma}\,\rd t+\cW_{t,xy}^{\cL,\bsigma}\, \rd t+\mathcal{E}_{t,xy}^{\cL,\bsigma}\, \rd t+\rd \cB_{t,xy}^{\cL,\bsigma},
	\end{equation}
	where the quadratic error term is given by
	\begin{equation}\label{eq:quadL}            \mathcal{E}_{t,xy}^{\cL,\bsigma}:=\sum_{a}(T-\Theta)_{t,xa}^{\bsigma}(\cL-\cK)_{t,a y}^{\bsigma}.
	\end{equation}
	Similarly, subtracting \eqref{Theta_evolution} from \eqref{T_evolution} yields 
	\begin{align}
		\rd (T-\Theta)_{t,x,yy'}^{\bsigma}=&~\sum_{a}\Theta_{t,xa}^{\bsigma}(T-\Theta)_{t,a, yy'}^{\bsigma}\,\rd t+\sum_a (T-\Theta)_{t,xa}^{\bsigma}\Theta_{t,a, yy'}^{\bsigma}\,\rd t \nonumber\\
		&~+\cW_{t,x,yy'}^{T,\bsigma}\, \rd t+\mathcal{E}_{t,x,yy'}^{T,\bsigma}\, \rd t+\rd \cB_{t,x,yy'}^{T,\bsigma},\label{T_Theta_evolution}
	\end{align}
	where the quadratic error term is given by
	\begin{equation}\label{eq:quadT}
		\mathcal{E}_{t,x,yy'}^{T,\bsigma}:=\sum_{a}(T-\Theta)_{t,xa}^{\bsigma}(T-\Theta)_{t,a, yy'}^{\bsigma}.
	\end{equation}

	We define the evolution kernel from time $s$ to time $t$, for any $0\leq s\leq t< 1$, by
	\begin{equation}\label{U_s_t_equal_1_Theta}
		\begin{aligned}
			\cU_{s,t}^{\bsigma}:=\frac{1-sm(\sigma_1)m(\sigma_2)S}{1-tm(\sigma_1)m(\sigma_2)S}=1+\pa{t-s}\Theta_t^{\bsigma}.
		\end{aligned}
	\end{equation}
	This kernel satisfies the matrix evolution equation
	\begin{equation*}
		\begin{aligned}
			\frac{\rd }{\rd t}\,\cU_{s,t}^{\bsigma}=\Theta_t^{\bsigma}\cdot\cU_{s,t}^{\bsigma}=\cU_{s,t}^{\bsigma}\cdot\Theta_t^{\bsigma},\quad \text{with}\quad \cU_{s,s}^{\bsigma}=I.
		\end{aligned}
	\end{equation*}
	Using this kernel, the SDEs \eqref{cL_cK_evolution} and \eqref{T_Theta_evolution} can be rewritten in their corresponding integral formulations.

	\begin{lemma}
		Suppose $s\in [0,1)$ and let $s\le \tau<1$ be a stopping time with respect to the matrix Brownian motion $\mathbf{B}(t)$. For any $\bsigma=(\sigma_1,\sigma_2)\in\ha{-,+}^2$ and $x,y\in\ZL$, equation \eqref{cL_cK_evolution} can be rewritten as
		\begin{equation}\label{cL_cK_evolution_step_2}
			\begin{aligned}
				(\cL-\cK)_{\tau,xy}^{\bsigma}=(\cL-\cK)_{s,xy}^{\bsigma}&+\int_{s}^{\tau}\qa{\Theta_t^{\bsigma}\circ (\cL-\cK)_t^{\bsigma}}_{xy}\,\rd t +\int_{s}^{\tau}\cW_{t,xy}^{\cL,\bsigma}\, \rd t+\int_{s}^{\tau}\mathcal{E}_{t,xy}^{\cL,\bsigma}\, \rd t+\int_{s}^{\tau}\rd \cB_{t,xy}^{\cL,\bsigma}.
			\end{aligned}
		\end{equation}
		Applying Duhamel's principle yields the integrated representation 
		\begin{align}
			(\cL-\cK)_{\tau,xy}^{\bsigma}=&~\qa{\cU_{s,\tau}^{\bsigma}\circ(\cL-\cK)_{s}^{\bsigma}}_{xy}+\int_{s}^{\tau}\qb{\cU_{t,\tau}^{\bsigma}\circ\cW_{t}^{\cL,\bsigma}}_{xy}\, \rd t \nonumber\\
			&~+\int_{s}^{\tau}\qb{\cU_{t,\tau}^{\bsigma}\circ\mathcal{E}_{t}^{\cL,\bsigma}}_{xy}\, \rd t+\int_{s}^{\tau}\qb{\cU_{t,\tau}^{\bsigma}\circ\rd \cB_t^{\cL,\bsigma}}_{xy}.\label{integrated_cL_cK_evolution}
		\end{align}
		Here, for any $0\leq s\leq t< 1$, we define the actions of $\Theta_t$ and $\cU_{s,t}$ on a matrix $X\in\C^{\ZL\times\ZL}$ as
		\begin{equation}\label{def_action_Theta_cU}
			\begin{aligned}
				\Theta_t^{\bsigma}\circ X:=\Theta_t^{\bsigma}\cdot X+X\cdot\Theta_t^{\bsigma},\quad \cU_{s,t}^{\bsigma}\circ X:=\cU_{s,t}^{\bsigma}\cdot X\cdot\cU_{s,t}^{\bsigma}.
			\end{aligned}
		\end{equation}
		Similarly, equation \eqref{T_Theta_evolution} implies that for diagonal $T$-variables,
		\begin{align}
			(T-\Theta)_{\tau,xy}^{\bsigma}=&~\qa{\cU_{s,\tau}^{\bsigma}\circ(T-\Theta)_{s}^{\bsigma}}_{xy}+\int_{s}^{\tau}\qb{\cU_{t,\tau}^{\bsigma}\circ\cW_{t}^{T,\bsigma}}_{xy}\, \rd t \nonumber\\
			&~+\int_{s}^{\tau}\qb{\cU_{t,\tau}^{\bsigma}\circ\mathcal{E}_{t}^{T,\bsigma}}_{xy}\, \rd t+\int_{s}^{\tau}\qb{\cU_{t,\tau}^{\bsigma}\circ\rd \cB_t^{T,\bsigma}}_{xy}.\label{integrated_T_Theta_evolution_d}
		\end{align}
		For off-diagonal $T$-variables with $y\ne y'$, we obtain (recall that $\Theta_{t,x,yy'}^{\bsigma}=0$ in this case)
		\begin{equation}\label{integrated_T_Theta_evolution_od}
			\begin{aligned}
				T_{\tau,x,yy'}^{\bsigma}=\qa{\cU_{s,\tau}^{\bsigma}\cdot T_{s,\cdot,yy'}^{\bsigma}}_{x}&+\int_{s}^{\tau}\qb{\cU_{t,\tau}^{\bsigma}\cdot\cW_{t,\cdot,yy'}^{T,\bsigma}}_{x}\, \rd t +\int_{s}^{\tau}\qb{\cU_{t,\tau}^{\bsigma}\cdot\mathcal{E}_{t,\cdot,yy'}^{T,\bsigma}}_{x}\, \rd t+\int_{s}^{\tau}\qb{\cU_{t,\tau}^{\bsigma}\cdot\rd \cB_{t,\cdot,yy'}^{T,\bsigma}}_{x}.
			\end{aligned}
		\end{equation}
	\end{lemma}

    \begin{remark}\label{remark_non_measurability}
        Rigorously speaking, the martingale terms in the integrated representations \eqref{integrated_cL_cK_evolution}, \eqref{integrated_T_Theta_evolution_d}, and \eqref{integrated_T_Theta_evolution_od} are ill-defined, since the integrand processes therein may not be progressively measurable, due to the dependence on the stopping time $\tau$. Therefore, taking \eqref{integrated_cL_cK_evolution} as an example, we should regard the martingale term therein as an abbreviation
        \begin{equation}
            \int_{s}^{\tau}\cU_{t,\tau}^{\bsigma}\circ\rd \cB_t^{\cL,\bsigma}:=\cU_{0,\tau}^{\bsigma}\circ\int_{s}^{\tau}(\cU_{0,t}^{\bsigma})^{-1}\circ\rd \cB_t^{\cL,\bsigma}.
        \end{equation}
        The integrated martingale terms in \eqref{integrated_T_Theta_evolution_d} and \eqref{integrated_T_Theta_evolution_od} should be understood in a similar sense.
    \end{remark}

	Applying the Burkholder-Davis-Gundy (BDG) inequality, we obtain the following lemma, which provides high-moment bounds for the martingale terms in \eqref{integrated_cL_cK_evolution}, \eqref{integrated_T_Theta_evolution_d}, and \eqref{integrated_T_Theta_evolution_od}.
	
	\begin{lemma}\label{lemma_martingale_to_cB}
		Suppose $s\in [0,1)$ and let $s\le \tau<1$ be a stopping time with respect to the matrix Brownian motion $\mathbf{B}(t)$. Let $p\in2\N$ and let $A_t,B_t\in\C^{\ZL\times\ZL}$ be deterministic matrices. Then, there exists a constant $C_p>0$, depending only on $p$, such that
		\begin{equation}\label{bound_martingale_to_cB}
			\begin{aligned}
				\E\absbb{\int_{s}^{\tau}\pa{A_t\cdot\rd \cB_{t}^{\star,\bsigma}\cdot B_t}_{xy}}^p\leq C_p \E\pa{\int_s^{\tau}\qa{  \pa{A_t\otimes \ol{A_t}\otimes B_t \otimes \ol{B_t} } \circ\pa{\cB\otimes\cB}_{t}^{\star,\bsigma}}_{xxyy}\,\rd t}^{p/2}
			\end{aligned}
		\end{equation}
		for all $\star\in\ha{\cL,T}$, $\bsigma=(\sigma_1,\sigma_2)\in\ha{-,+}^2$, and $x,y\in\ZL$.    Here, the 4-tensor $\pa{\cB\otimes\cB}_{t;x_1x_2y_1y_2}$ is defined as
		\begin{align}
			\pa{\cB\otimes\cB}_{t;x_1x_2y_1y_2}^{\cL,\bsigma}:=&\sum_{a}(G_t^{\sigma_1}S^{\pa{x_1}}G_t^{\sigma_2}S^{\pa{y_1}}G_t^{\sigma_1}S^{\pa{a}}G_t^{-\sigma_1}S^{\pa{y_2}}G_t^{-\sigma_2}S^{\pa{x_2}}G_t^{-\sigma_1})_{aa}\label{def_tensor_cB}\\
			+&\sum_a (G_t^{-\sigma_2}S^{\pa{x_2}}G_t^{-\sigma_1}S^{\pa{y_2}}G_t^{-\sigma_2}S^{\pa{a}}G_t^{\sigma_2}S^{\pa{y_1}}G_t^{\sigma_1}S^{\pa{x_1}}G_t^{\sigma_2})_{aa},\nonumber\\
			\pa{\cB\otimes\cB}_{t;x_1x_2y_1y_2}^{T,\bsigma}:=&\sum_{a}(G_t^{\sigma_1}S^{\pa{x_1}}G_t^{\sigma_2})_{ay_1}(G_t^{\sigma_1}S^{\pa{a}}G_t^{-\sigma_1})_{y_1y_2}\p{G_t^{-\sigma_2}S^{\pa{x_2}}G_t^{-\sigma_1}}_{y_2a}\label{def_tensor_cB2}\\
			+&\sum_a (G_t^{-\sigma_2}S^{\pa{x_2}}G_t^{-\sigma_1})_{ay_2}(G_t^{-\sigma_2}S^{\pa{a}}G_t^{\sigma_2})_{y_2y_1}\p{G_t^{\sigma_1}S^{\pa{x_1}}G_t^{\sigma_2}}_{y_1a}.\nonumber
		\end{align}
		The tensor product \(  A_t \otimes \ol{A_t} \otimes B_t \otimes \ol{B_t} \) acts on a 4-tensor \( \cX \) by
		\begin{equation}\label{eq:action4tensor}
		 \qa{ \pa{A_t\otimes \ol{A_t}\otimes B_t \otimes \ol{B_t} } \circ \cX}_{x_1x_2y_1y_2}=\sum_{a_1,a_2,b_1,b_2\in\ZL}A_{t,x_1a_1}\ol{A}_{t,x_2a_2}B_{t,b_1y_1}\ol{B}_{t,b_2y_2}\cdot\cX_{a_1a_2b_1b_2}.
		\end{equation}
		For the martingale terms arising from the evolution of the off-diagonal \(T\)-variables, we similarly have
		\begin{equation}\label{bound_martingale_to_cB_T_od}
			\begin{aligned}
				\E\absbb{\int_{s}^{\tau}\qa{A_t\cdot\rd \cB_{t,\cdot,yy'}^{T,\bsigma}}_{x}}^p\leq C_p\E\pa{\int_{s}^\tau\qa{  (A_t\otimes \ol{A_t} )\circ(\cB\otimes\cB)_{t,yy'}^{T,\bsigma}}_{xx}\,\rd t}^{p/2}
			\end{aligned}
		\end{equation}
		for any $\bsigma=(\sigma_1,\sigma_2)\in\ha{-,+}^2$ and $y, y'\in \ZL$, where 
		\begin{align}
			(\cB\otimes\cB)_{t;yy'}^{T,\bsigma}\pa{x_1,x_2}:=&\sum_{a}(G_t^{-\sigma_2}S^{\pa{x_2}}G_t^{-\sigma_1})_{y'a}(G_t^{\sigma_1}S^{\pa{a}}G_t^{-\sigma_1})_{yy}(G_t^{\sigma_1}S^{\pa{x_1}}G_t^{\sigma_2})_{ay'}\label{def_2tensor_cB2}\\
			+&\sum_{a}(G_t^{\sigma_1}S^{\pa{x_1}}G_t^{\sigma_2})_{ya}(G_t^{-\sigma_2}S^{\pa{a}}G_t^{\sigma_2})_{y'y'}(G_t^{-\sigma_2}S^{\pa{x_2}}G_t^{-\sigma_1})_{ay}.\nonumber
		\end{align}
		The action of $A_t\otimes \ol{A_t} $ on 2-tensors is defined analogously to \eqref{eq:action4tensor}, with $B_t=I$.
		We shall refer to the tensors $\cB\otimes\cB$ defined above as the \emph{quadratic variation tensors}.
	\end{lemma}

	\subsection{Proofs of \Cref{theorem_local_law_entrywise,theorem_local_law_loops} and \Cref{theorem_local_law_T} for $\al\in(0,1)$}\label{sec:strat_al01}
	
	In this subsection, we present the proofs of the main results in the regime $\al\in(0,1)$. The cornerstone of our approach is a dynamical analysis of the $\cL$-loops and the $T$-variables along the flows in Definitions \ref{def_H_t_flow} and \ref{def_z_t_flow}. More precisely, we establish \Cref{theorem_local_laws_along_flow}--\ref{theorem_que_along_flow}, from which \Cref{theorem_local_law_entrywise,theorem_local_law_loops} and \Cref{theorem_local_law_T} for $\al\in(0,1)$ follow as direct consequences.

	\begin{lemma}\label{theorem_local_laws_along_flow}
		Under the assumptions in \Cref{theorem_delocalization}, fix $\al\in(0,1)$, small constants $\kappa,\fc\in(0,1)$, and a spectral parameter $z=E+\ii \eta\in\mathbf{D}_{\kappa,\fc}$. Consider the flow framework in \Cref{def_H_t_flow,def_z_t_flow}, with the flow parameter chosen as in \eqref{eq:t0E0}. Then the following estimates hold uniformly for $t\in[0,\tf]$.
		\begin{itemize}
			\item {\bf Local law:} For any $x,y\in\ZL$, the following entrywise and averaged local laws hold:
			\begin{equation}\label{flow_local_law}
				\begin{aligned}
					\absa{(G_t-m)_{xy}}^2\prec B_t\p{|x-y|},\quad |\avg{(G_t-m)S^{\pa{x}}}|\prec B_t\p{0}.
				\end{aligned}
			\end{equation}
			\item {\bf $\cal L$-loop estimate:} For any $\bsigma\in\ha{-,+}^2$ and $x,y\in\ZL$, the following $\cal L$-loop estimate holds:
			\begin{equation}\label{flow_loop_local_law}
				\begin{aligned}
					(\cL_t-\cK_t)_{xy}^{\bsigma}\prec \qa{B_t\p{0}}^{1/5}B_t\p{|x-y|}.
				\end{aligned}
			\end{equation}
			Here, the exponent \( 1/5 \) may be replaced by any positive constant not exceeding \( 1/4 \).
			
		\end{itemize}
	\end{lemma}

	\begin{lemma}\label{lemma_weak_local_law_T_along_flow}
		In the setting of \Cref{theorem_local_laws_along_flow}, the following estimate for the $T$-variables holds uniformly for $t\in[0,1-W/N]$ (equivalently, $1-t\ge \eta_{\txt{flat}}$, see \eqref{eq:etaflat}). For any $\bsigma\in\h{-,+}^2$ and $x,y,y'\in\ZL$, we have
		\begin{equation}\label{weak_local_law_T_along_flow}
			\begin{aligned}
				\pa{T-\Theta}_{t,x,yy'}^{\bsigma}\prec \qa{B_t(0)}^{7/10}\qa{B_t(|x-y|\wedge |x-y'|)}^{1/2}.
			\end{aligned}
		\end{equation}
	\end{lemma}
	\begin{lemma}\label{lemma_flat_local_law_along_flow}
		In the setting of \Cref{theorem_local_laws_along_flow}, the following refined $\cL$-loop and $T$-variable estimates hold uniformly for $t\in[1-W/N,\tf]$:
		\begin{align}
			\max_{\bsigma\in\h{-,+}^2}\max_{x,y\in\ZL}\absa{(\cL-\cK)_{t,xy}^{\bsigma}}&\prec W^{-6/5}+\pa{{N(1-t)}}^{-7/4},\label{flat_local_law_along_flow}\\
			\max_{\bsigma\in\h{-,+}^2}\max_{x,y,y'\in\ZL}\absa{(T-\Theta)_{t,x,yy'}^{\bsigma}}&\prec W^{-6/5}+\pa{{N(1-t)}}^{-3/2}.\label{flat_local_law_along_flow2}
		\end{align}
	\end{lemma}
	
	\begin{lemma}\label{theorem_que_along_flow}
		In the setting of \Cref{theorem_local_laws_along_flow}, the following improved estimates for expected $\cL$-loops hold uniformly for $t\in[1-W/N,\tf]$:
		\begin{equation}\label{exp_flat_local_law}        \max_{\bsigma\in\h{-,+}^2}\max_{x,y\in\ZL}\absa{\E(\cL-\cK)_{t,xy}^{\bsigma}}\prec W^{-6/5}+\pa{N(1-t)}^{-9/4}.
		\end{equation}
	\end{lemma}

	\begin{proof}[\bf Proof of \Cref{theorem_local_law_entrywise,theorem_local_law_loops} and \Cref{theorem_local_law_T} for $\al\in(0,1)$]
		By \eqref{eq:zztE}, we have $G\pa{z}\stackrel{d}{=}\sqrt{\tf}G_{\tf,\sE}$. Moreover, in view of \eqref{eq:t0E0} and \eqref{re_im_z_t}, we know that $1-\tf\asymp \im z \asymp 1-|m(z)|^2$, which implies that 
		\begin{equation}
			B(\im z,r)\asymp B_{\tf}(r)\asymp B_{|m(z)|^2}(r),\quad \forall 0\le r\le N.\label{eq:Betar}  
		\end{equation}
		Then, for each fixed $z\in\mathbf{D}_{\kappa,\fc}$, \Cref{theorem_local_laws_along_flow} yields the estimates \eqref{eq_entrywise_local_law}, \eqref{eq_average_local_law}, and \eqref{eq_loop_local_law_2}, whereas \eqref{weak_local_law_T_along_flow} in \Cref{lemma_weak_local_law_T_along_flow} and \eqref{flat_local_law_along_flow2} in \Cref{lemma_flat_local_law_along_flow} together imply \eqref{eq_chain_local_law_2}. 
		For each fixed $z\in\mathbf{D}_{\kappa,\fc}^{\txt{flat}}$, \Cref{lemma_flat_local_law_along_flow} yields \eqref{eq_loop_local_law_2_improved} and \eqref{eq_chain_local_law_2_improved}, and \Cref{theorem_que_along_flow} gives \eqref{eq_expected_loop_local_law}. 
		Finally, a standard $N^{-C}$-net and perturbation argument extends these estimates uniformly to all $z$.
	\end{proof}
	
	Note that at $t=0$, we have $G_{0}^\sigma=m^\sigma I$ for any $\sigma\in\ha{-,+}$. By \Cref{def_loops_and_T_variables}, it follows immediately that for all $x,y,y'\in\ZL$ and $\boldsymbol\sigma=(\sigma_1,\sigma_2)\in\{-,+\}^2$,
	\begin{equation}\label{eq:initialL}
		{\cL}^{\bsigma}_{0, xy}= {\cK}^{\bsigma}_{0, xy},\quad T^{\bsigma}_{0, x,yy'}= \Theta^{\bsigma}_{0, x,yy'}.   
	\end{equation}
	Thus, \Cref{theorem_local_laws_along_flow} follows directly from the following bootstrap statement, which propagates the $\cL$-loop estimate progressively along the stochastic flow starting from $t=0$.
	
	\begin{lemma}\label{lemma_bootstrap}
		In the setting of \Cref{theorem_local_laws_along_flow}, suppose that the estimates \eqref{flow_local_law} and \eqref{flow_loop_local_law} hold at some fixed time $s\in\qa{0,\tf}$. Then, there exists a sufficiently small constant $c_\alpha\in(0,1)$, depending only on $\kappa,\ \fc$, and $\alpha$, such that for any $t\in\qa{s,\tf}$ satisfying
		\begin{equation}\label{bootstrap_assumption}
			\begin{aligned}
				\qa{B_t\p{0}}^{c_\alpha}\leq \frac{1-t}{1-s}\leq 1,
			\end{aligned}
		\end{equation}
		the estimates \eqref{flow_local_law} and \eqref{flow_loop_local_law} also hold at time $t$. \end{lemma}
	
	\begin{proof}[\bf Proof of \Cref{theorem_local_laws_along_flow}]
	For any $t\in[0,\tf]$, we iterate \Cref{lemma_bootstrap} from time $0$ up to time $t$ along an increasing sequence of times $\{t_k:0\le k \le n\}$ with $t_0=0$ and $t_n=t$, chosen such that $n=\OO(1)$ and $(1-t_k)/(1-t_{k-1})\ge [B_{t_k}(0)]^{c_\al}$. After $n$ induction steps, we obtain \Cref{theorem_local_laws_along_flow} at each fixed time $t\in[0,\tf]$. Finally, a standard $N^{-C}$-net argument combined with perturbation estimates extends these bounds uniformly to all $t\in[0,\tf]$.
	\end{proof}

	The proof of \Cref{lemma_bootstrap} is organized into three steps.

	\medskip
	\noindent
	{\bf Step 1} (Weak bounds for $\cL$-loops and weak local law). The $\cL$-loops satisfy the following a priori bound:
	\begin{equation}\label{step_1_priori_bound}
		\begin{aligned}
			\max_{\bsigma\in\h{-,+}^2}\max_{x,y\in\ZL}\absa{\cL_{u,xy}^{\bsigma}}\prec \frac{1-s}{1-u} \qa{B_s\p{0}}^{1/2},\quad \forall u\in\qa{s,t}.
		\end{aligned}
	\end{equation}
	In addition, the following weak entrywise local law holds:
	\begin{equation}\label{step_1_weak_local_law}
		\begin{aligned}
			\|G_u-m\|_{\max}\prec \qa{B_u\p{0}}^{1/5},\quad \forall u\in\qa{s,t}.
		\end{aligned}
	\end{equation}

	\noindent
	{\bf Step 2} (Sharp bounds for $\cal L$-loops and sharp local law). There exists a sufficiently large constant $C_\alpha>1$, depending only on $\kappa$, $\fc$, and $\alpha$, such that for any $x,y\in\ZL$,
	\begin{equation}\label{step_2_weak_estimate_loops}
		\begin{aligned}
			\max_{\bsigma\in\h{-,+}^2}\absa{(\cL_u^{\bsigma}-\cK_u^{\bsigma})_{xy}}\prec \pa{\frac{1-s}{1-u}}^{C_\alpha}\qa{B_u\p{0}}^{1/5}B_u\p{|x-y|},\quad \forall u\in[s,t].
		\end{aligned}
	\end{equation}
	As a consequence, if the constant $c_\alpha$ in \eqref{bootstrap_assumption} is chosen sufficiently small, we obtain the following sharp bound for the $\cL$-loops:
	\begin{equation}\label{step_2_sharp_bound_loops}
		\begin{aligned}        \max_{\bsigma\in\h{-,+}^2}\absa{\cL_{u,xy}^{\bsigma}}\prec B_u\p{|x-y|},\quad \forall u\in\qa{s,t}.
		\end{aligned}
	\end{equation}
	Moreover, for any $x,y\in\ZL$, the following sharp local laws hold:
	\begin{equation}\label{step_2_sharp_local_law}
		\begin{aligned}
			\absa{(G_u-m)_{xy}}^2\prec B_u\p{|x-y|},\quad \absa{\avg{(G_u-m)S^{\pa{x}}}}\prec B_u\p{0},\quad \forall u\in\qa{s,t}.
		\end{aligned}
	\end{equation}

	\noindent
	{\bf Step 3} ($\cL$-loop estimate). For any $x,y\in\ZL$, we establish the refined estimate
	\begin{equation}\label{step_3_estimate_loops}
		\begin{aligned}
			\max_{\bsigma\in\h{-,+}^2}\absa{(\cL_u^{\bsigma}-\cK_u^{\bsigma})_{xy}}\prec \qa{B_u\p{0}}^{1/5}B_u\p{|x-y|},\quad \forall u\in\qa{s,t}.
		\end{aligned}
	\end{equation}
	
	Finally, we note that the estimates \eqref{step_1_priori_bound}--\eqref{step_3_estimate_loops} can be upgraded to high-probability bounds that hold uniformly for all $u\in[s,t]$ once again via a standard $N^{-C}$-net argument. For brevity, we shall not repeatedly emphasize this uniformity in the arguments that follow.

	The proof of \Cref{lemma_bootstrap} constitutes the main focus of the next section. More precisely, fix a time $t$ satisfying \eqref{bootstrap_assumption}. The above three-step strategy will be implemented in \Cref{sec_step_1_to_3}\footnote{  Since we restrict to the analysis of $\cL$-loops of length 2, our three steps correspond to the first, second, and fifth steps in the previous works \cite{Band1D,Band2D,truong2025localizationlengthfinitevolumerandom,dubova2025delocalizationnonmeanfieldrandommatrices}.}, based on a sequence of estimates for the error terms appearing in the evolution equation \eqref{cL_cK_evolution}.
	The derivation of these estimates requires several new ideas. In particular, we must treat the general variance profile structure using only the $T$-variables and $\cL$-loops of length 2. This is in contrast to the previous works \cite{Band1D,Band2D,truong2025localizationlengthfinitevolumerandom,dubova2025delocalizationnonmeanfieldrandommatrices}, which assume a block variance structure, as well as to the approach of \cite{erdos2025zigzagstrategyrandomband}, where the analysis relies on both $\cL$-loops and $\cal C$-chains (see \Cref{def_n_loops_and_chains} below) of arbitrarily large lengths.
	These new arguments will be developed in detail in \Cref{sec_estimates_error_terms}.
	Moreover, as discussed in the introduction, the treatment of the light-weight term presents a substantially greater challenge than in the case of regular RBMs in dimensions $d\ge 3$ considered in \cite{dubova2025delocalizationnonmeanfieldrandommatrices}. 
    To address this difficulty, we develop new graphical techniques in \Cref{sec_graphical_argument}, extending the methods of \cite{dubova2025delocalizationnonmeanfieldrandommatrices}, which allow us to recover the optimal decay of the $\cal L$-loops in Steps 2 and 3. We refer to \Cref{sec:graph_ideas} for a more detailed discussion of the main difficulties and the key new ideas.

	Finally, once \Cref{lemma_bootstrap}, and hence \Cref{theorem_local_laws_along_flow}, have been established, we turn to the proofs of Lemmas \ref{lemma_weak_local_law_T_along_flow}--\ref{theorem_que_along_flow} in \Cref{sec_proof_of_que}. These results are derived from \eqref{flow_local_law} and \eqref{flow_loop_local_law}, the evolution equations \eqref{integrated_cL_cK_evolution}--\eqref{integrated_T_Theta_evolution_od}, as well as certain estimates established in the course of the three-step proof of \Cref{lemma_bootstrap}.
	
	\begin{remark}
		Note that the $\cL$-loop defined in \eqref{eq:cL_and_T} is a weighted average of the diagonal $T$-variables. Therefore, any bound on the $(T-\Theta)$ variables immediately yields the corresponding bound for the $(\cL-\cK)$ loops. Moreover, the evolution equation \eqref{T_Theta_evolution} for the $T$-variables has a structural advantage over \eqref{cL_cK_evolution} for the $\cL$-loops: the quadratic error term in \eqref{eq:quadT} is self-consistent, in the sense that it depends only on the $(T-\Theta)$ variables. In contrast, the quadratic error in \eqref{eq:quadL} depends on both $(T-\Theta)$ variables and $(\cL-\cK)$ loops, which complicates the analysis. 
		However, there is a serious obstacle in applying a direct dynamical analysis to \eqref{T_Theta_evolution}. Its martingale term cannot be controlled with sufficient precision, unlike the martingale term in \eqref{cL_cK_evolution}; see \Cref{rmk:MG_for_T_fails}. This reflects a well-known phenomenon in the study of RBMs: the additional averaging in the definition of the $\cL$-loops produces an extra fluctuation-averaging mechanism, leading to significantly improved estimates. For this reason, our proof is based on a dynamical analysis of the $\cL$-loops as in \Cref{theorem_local_laws_along_flow}. Additional work is required to handle the quadratic error term \eqref{eq:quadL} and to turn it into a self-consistent structure depending only on $(\cL-\cK)$ loops; see \Cref{lemma_double_difference_decomposition}.
	\end{remark}

	\section{Dynamical analysis of \texorpdfstring{$\cL$}{cL}-loops: Proof of \texorpdfstring{\Cref{lemma_bootstrap}}{Lemma 3.12}}\label{sec_step_1_to_3}
	
	In this section, we prove \Cref{lemma_bootstrap} by implementing the three-step strategy outlined above, namely, by successively establishing the estimates \eqref{step_1_priori_bound}--\eqref{step_3_estimate_loops}.

	\subsection{Step 1: Weak bounds for $\cL$-loops and weak local law}
	
	Step 1 in the proof of \Cref{lemma_bootstrap} follows closely the argument in \cite[Section 5.1]{Band1D} for 1D regular RBMs. Accordingly, we provide only a brief outline.

	\begin{lemma}\label{2_loop_to_1_chain}
		In the setting of \Cref{lemma_bootstrap}, fix any constant $c>0$ and, for any $u\in[s,t]$, define the event
		\begin{equation*}
			\begin{aligned}
				\Omega_{u,c}:=\ha{\|G_u-m\|_{\max}\leq W^{-c}}.
			\end{aligned}
		\end{equation*}
		Then, uniformly for $u\in[s,t]$, the off-diagonal entries of $G_u$ satisfy
		\begin{equation}\label{upper_bound_2_to_1}
			\mathbf{1}\p{\Omega_{u,c}}\cdot |\pa{G_u}_{xy}|^2\prec \mathbf{1}\p{\Omega_{u,c}}\cdot
			T^{(-,+)}_{u,xy}\wedge T^{(+,-)}_{u,yx} +\cS_{xy} \prec \cL_{u,xy}+\cS_{xy},\quad \forall  x\neq y \in\ZL, 
		\end{equation}
		where we use the abbreviations introduced in \eqref{abbreviation_cL_T}. For the diagonal entries, we have
		\begin{equation}\label{entrywise_bound_2_to_1}            \mathbf{1}\p{\Omega_{u,c}}\cdot\max_{x\in\ZL}\absa{(G_u)_{xx}-m}^2\prec \max_{x,y\in\ZL}\pa{\cL_{u,xy}+\cS_{xy}}.
		\end{equation}
		Furthermore, suppose that for some constant $c>0$ and a deterministic control parameter $\Psi_u$ satisfying $W^{-1/2}\le \Psi_u\le W^{-c}$, we have
		\begin{equation*}
			\begin{aligned}
				\|G_u-m\|_{\max}\prec W^{-c},\quad \text{and}\quad     \|\cL_{u}\|_{\max}\prec \Psi_{u}^2.
			\end{aligned}
		\end{equation*}
		Then the following averaged local law holds:
		\begin{equation}\label{average_bound_2_to_1}
			\max_{x\in\ZL}|\avg{(G_u-m)S^{\pa{x}}}|\prec \Psi_u^2 .
		\end{equation}
		The above estimates also hold in the regime $\alpha\in[1,\infty)$.
	\end{lemma}
	\begin{proof}
		This statement was established as Lemma 4.1 in \cite{Band1D} for the 1D regular RBM. The proof given there is model-independent and therefore applies equally to the present power-law RBM. We omit the details.
	\end{proof}

	\begin{lemma}\label{lemma_continuity_estimate}
		In the setting of \Cref{lemma_bootstrap}, assume that $\varepsilon\leq s\leq t\leq \tf$ for some fixed constant $\varepsilon>0$.\footnote{For $t\leq \varepsilon$, the $\max$-norm bound on $\cL_{u}^{\bsigma}$ is directly provided by the $\OO(1)$ operator norm bound of $G_u$.} Then, for any fixed large constant $C>0$, on the event $\Omega_u:=\ha{\|G_u\|_{\max}\leq C}$, the following continuity estimate holds uniformly for all $u\in[s,t]$:
		\begin{equation}     \label{step_1_priori_bound_pf}      \mathbf{1}\p{\Omega_u}\cdot\max_{\bsigma\in\h{-,+}^2} \|\cL_{u}^{\bsigma}\|_{\max}\prec \frac{1-s}{1-u} \qa{B_s\p{0}}^{1/2}.
		\end{equation}
	\end{lemma}
	\begin{proof}
		This statement is essentially contained as a special case of the proof of \cite[Lemma 5.1]{Band1D}. More precisely, by the local law \eqref{flow_local_law} and $\cL$-loop estimate \eqref{flow_loop_local_law} at time $s$, we have
		\[\norma{G_s-m}_{\max}\prec \qa{B_s\p{0}}^{1/2},\quad 
		\max_{\bsigma\in\{-,+\}^2}\max_{x,y\in\ZL}\absa{\cL_{s,xy}^{\bsigma}}\prec B_s\p{0}.\]
		Combining these bounds with the argument in \cite[Lemma 5.1]{Band1D} for $\cL$-loops of length 2 yields \eqref{step_1_priori_bound_pf}.\footnote{In fact, \cite[Lemma 5.1]{Band1D} assumes bounds for $\cL$-loops of arbitrary length and derives a weak estimate for any $n$-$\cal L$-loop with fixed $n$ (see \Cref{def_n_loops_and_chains}). Here, only the upper bounds of $2$-loops are available, which results in the exponent $1/2$ in the RHS of \eqref{step_1_priori_bound_pf}.} Alternatively, one may also follow the argument in the proof of equation (5.7) in \cite{yang2021delocalization}. We omit the details.
	\end{proof}
	
	Combining Lemmas \ref{2_loop_to_1_chain} and \ref{lemma_continuity_estimate}, Step 1 of the proof of \Cref{lemma_bootstrap}---that is, the derivation of \eqref{step_1_priori_bound} and \eqref{step_1_weak_local_law}---proceeds exactly as in Section 5.1 of \cite{Band1D}. We therefore omit further details.

	\subsection{Estimates for the terms in the evolution equation}

	Before turning to the proofs of Steps 2 and 3, we collect in this subsection several estimates for the terms appearing on the right-hand side of the evolution equation \eqref{cL_cK_evolution_step_2}, as well as some basic bounds on the shape parameters and the associated evolution kernels. We begin with a convolution inequality for the shape parameters.
	
	\begin{lemma}\label{lemma_B_t}
		In the setting of \Cref{assumption_input_bound}, fix $\alpha\in (0,1)$. Then, there exists a constant $C_\alpha>0$, depending only on $\alpha$, such that the following estimate holds for all $0\leq s\leq t\leq 1- N^{-1}$:
		\begin{equation}\label{convolution_2_B}
			\qa{B_s(|\cdot|)*B_t(|\cdot|)}\pa{x-y}\equiv \sum_{a\in\ZL}B_s\pa{|x-a|}B_t\pa{|a-y|}\le \frac{C_\alpha}{1-s}B_t\pa{|x-y|}.
		\end{equation}
	\end{lemma}
	
	Next, for the evolution kernel introduced in \eqref{U_s_t_equal_1_Theta}, we establish the following bound on its $\ell^\infty \to \ell^\infty$ operator norm with respect to the action defined in \eqref{def_action_Theta_cU}.

	\begin{lemma}\label{lemma_U_s_t}
		In the setting of \Cref{assumption_input_bound}, fix $\alpha\in (0,1)$. We regard the shape parameter $B_t$ as a matrix defined by $B_{t,xy}:=B_t(|x-y|)$. Then the following estimates hold for $\cU_{s,t}$, uniformly in $\bsigma=(\sigma_1,\sigma_2)\in\{-,+\}^2$ and $0\le s\le t\le 1-N^{-1}$.
		\begin{itemize}
			\item For any matrix $X\in\C^{\ZL\times \ZL}$, we have
			\begin{equation}\label{bound_cU_s_t_1}
				\begin{aligned}
					\norm{(\mathcal U_{s,t}^{\bsigma}\circ X)/B_{t}}_{\max}\prec \frac{1-s}{1-t} \|X/B_{s}\|_{\max},
				\end{aligned}
			\end{equation}
			where, with a slight abuse of notation, $A/B$ denotes the entrywise quotient of two matrices; for instance, $(X/B_s)_{xy}=X_{xy}/B_{s,xy}$.

			\item Suppose one of the following two conditions holds: (i) $\sigma_1=\sigma_2$; (ii) $\sigma_1\neq \sigma_2$ and $X\in\C^{\ZL\times\ZL}$ satisfies either the \emph{right sum-zero property}
			\begin{equation}\label{right_sum_zero}
				\sum_{y\in\ZL}X_{xy}=0,\quad \forall x\in\ZL,
			\end{equation}
			or the \emph{left sum-zero property}
			\begin{equation}\label{left_sum_zero}
				\sum_{x\in\ZL}X_{xy}=0,\quad \forall y\in\ZL.
			\end{equation}
			Then, we have the improved bound
			\begin{equation}\label{bound_cU_s_t_2}
				\begin{aligned}
					\norm{(\mathcal U_{s,t}^{\bsigma}\circ X)/B_{t}}_{\max}\prec \|X/B_{s}\|_{\max}.
				\end{aligned}
			\end{equation}
			\item Suppose that $\ell_s=\ell_t=N$ (i.e., $1-t\le 1-s\le (W/N)^\alpha$), and that one of the following conditions is satisfied: (i) $\sigma_1=\sigma_2$; (ii) $\sigma_1\neq \sigma_2$ and the vector $\bv\in\C^{\ZL}$ satisfies the sum-zero property $\sum_{a}\bv_a=0$. Then, we have the estimate         
            \begin{equation}\label{bound_cU_s_t_3}
				\begin{aligned}
					\norm{\mathcal U_{s,t}^{\bsigma}\cdot \bv}_{\max}\prec \|\bv\|_{\max}.
				\end{aligned}
			\end{equation}
		\end{itemize}
		
	\end{lemma}

	The proofs of \Cref{lemma_B_t,lemma_U_s_t} are postponed to \Cref{appendix_additional_proof_B_t_and_U_s_t,appendix_additional_proof_B_t_and_U_s_t2}, respectively.
    We now state the estimates for the terms appearing in the evolution equation \eqref{cL_cK_evolution_step_2}, which constitute the main inputs for Steps 2 and 3 in the proof of \Cref{lemma_bootstrap}. 
	To this end, for any length parameter $\ell\ge 0$ and $u\in[s,t]$, we introduce the truncated shape parameter
	\begin{equation*}
		B_u^{\pa{\ell}}\p{r}:=B_u\pa{r\wedge \ell}, \quad \forall r\geq 0.
	\end{equation*}
	In addition, we remove the spatial decay factor from the $(\cL-\cK)$-loops and define, for $u\in[s,t]$,
	\begin{equation}\label{eq:JUl}
		J_u^{\pa{\ell}}:=\max_{\bsigma\in\h{-,+}^2}\max_{x,y\in\ZL}\absa{(\cL-\cK)_{u,xy}^{\bsigma}}/B_u^{\pa{\ell}}\p{|x-y|}.
	\end{equation}
	With these definitions in place, we first estimate the second term on the RHS of \eqref{cL_cK_evolution_step_2} as follows. 
	
	\begin{lemma}\label{bound_leading_term}
		Under the assumptions of \Cref{lemma_bootstrap}, for any length parameter $\ell\geq 0$ and $u\in\qa{s,t}$, we have \begin{equation}\label{eq_bound_leading_term}          \max_{\bsigma\in\h{-,+}^2}\max_{x,y\in\ZL}\absa{\qa{\Theta_u^{\bsigma}\circ(\cL-\cK)_u^{\bsigma}}_{xy}}\Big/B_u^{\pa{\ell}}\p{|x-y|}\leq \frac{C}{1-u}\cdot J_u^{\pa{\ell}},
		\end{equation}
		where the constant $C>0$ depends only on $\kappa$, $\fc$, and $\alpha$.
	\end{lemma}
	\begin{proof}
		For any $\bsigma \in\h{-,+}^2$ and $x,y\in\ZL$, by the definition \eqref{def_action_Theta_cU} and the assumptions \eqref{Theta_bound_opposite_charge} and \eqref{Theta_bound_equaled_charge}, we have 
		\begin{equation*}
			\begin{aligned}
				\absa{\qa{\Theta_u^{\bsigma}\circ(\cL-\cK)_u^{\bsigma}}_{xy}}\leq C_{\alpha} J_u^{\pa{\ell}} \sum_{a\in\ZL} B_u\pa{|x-a|}B_u^{\pa{\ell}}\pa{|a-y|}+ C_{\alpha} J_u^{\pa{\ell}} \sum_{a\in\ZL} B_u^{\pa{\ell}}\pa{|x-a|}B_u\pa{|a-y|}.
			\end{aligned}
		\end{equation*}
		By \eqref{def_B_t}, we have 
		\begin{equation}
			\sum_{a\in\ZL}B_u\pa{|a|}\asymp (1-u)^{-1}.\label{eq:L1Bu}
		\end{equation} 
		Using \eqref{eq:L1Bu}, together with the convolution bound \eqref{convolution_2_B}, we bound the first term by
        \begin{align}
            \sum_{a\in\ZL} B_u\pa{|x-a|}B_u^{\pa{\ell}}\pa{|a-y|}=&~ \sum_{|a-y|\leq \ell} B_u\pa{|x-a|}B_u\pa{|a-y|}+  \sum_{|a-y|> \ell} B_u\pa{|x-a|}B_u\pa{\ell}\nonumber\\
			\lesssim&~\frac{1}{1-u} \qa{B_u\pa{|x-y|}+B_u\pa{\ell}}\asymp \frac{1}{1-u}  B_u^{\pa{\ell}}\pa{|x-y|}.\label{eq:convolution_BBell}
        \end{align}
The second term is estimated in the same way. Dividing by $B_u^{(\ell)}(|x-y|)$ yields \eqref{eq_bound_leading_term}.
	\end{proof}
	
	Next, we estimate the quadratic error term, i.e., the fourth term on the RHS of \eqref{cL_cK_evolution_step_2}. This is analogous to estimate (3.26) in \cite{dubova2025delocalizationnonmeanfieldrandommatrices}. However, in the absence of a block structure in the variance profile matrix $S$, additional ideas and more delicate analysis are required; these will be developed in the proof of \Cref{lemma_double_difference_decomposition}.
	
	\begin{definition}[Admissible control parameter]\label{assm_admissible}
		For any $u\in[s,t]$, let $\Psi_u:[0,\infty)\to(0,\infty)$ be a deterministic, non-increasing function. We say that $\{\Psi_u(\ell):\ell\ge 0\}_{u\in[s,t]}$ is a family of \emph{admissible control parameters} if the following conditions hold uniformly: 
		\begin{itemize}
			\item [(i)] ({\bf Polynomial decay}) For any constant $c>0$, there exists a constant $C\equiv C(c)>0$ such that
			\begin{equation}\label{polydecay_Psiu}
				\Psi_u(0)\le C\Psi_u(W),\quad \text{and}\quad \Psi_u(c\ell)\le C\Psi_u\pa{\ell}, \quad \forall \ell\ge 0.   
			\end{equation}
			
			\item[(ii)] ({\bf Lower bound}) $\Psi_u(0)$ satisfies $W^{-1/2}\le \Psi_u(0)\leq W^{-\e}$ for some constant $\e>0$. Moreover, we have the pointwise lower bound
			\begin{equation}\label{polydecay_Psiu2}
				\cS_{xy}\prec \Psi_u^2(|x-y|),\quad \forall x,y\in \ZL.
			\end{equation}
			
			\item[(iii)] ({\bf Convolution property}) For all $x,y\in \ZL$, the following estimate holds: 
			\begin{equation}\label{polydecay_Psiu3}
				\sum_a \cS_{xa}\Psi_u^2(|a-y|)\prec \Psi_u^2(|x-y|).
			\end{equation}
		\end{itemize}
		
	\end{definition}
	
	\begin{lemma}\label{bound_double_difference_term}
		Under the assumptions of \Cref{lemma_bootstrap}, let $\{\Psi_{u}\}_{u\in[s,t]}$ be a family of admissible control parameters. Suppose that for any $u\in[s,t]$, the following bounds hold  uniformly:
		\begin{equation}\label{double_difference_input}
			\max_{\bsigma\in\h{(-,+),(+,-)}}\cL_{u,xy}^{\bsigma}\prec \Psi_{u}^{2}(|x-y|),\quad B_{u}(|x-y|)\prec \Psi_{u}^{2}(|x-y|),\quad  \forall x,y\in\ZL.
		\end{equation}
		Then for any deterministic length parameter $\ell\geq 0$, there exists a constant $C>0$ depending only on $\kappa$, $\fc$, and $\alpha$ such that the quadratic error term satisfies, with high probability and uniformly in $u\in[s,t]$ and $x,y\in\ZL$,
		\begin{equation}\label{eq_bound_double_difference_term}
			\begin{aligned}
				\max_{\bsigma\in\h{-,+}^2} \absa{\mathcal{E}_{u,xy}^{\cL,\bsigma}}\leq \frac{C}{1-u}\qB{ J_u^{\pa{\ell}}+\mathbf 1_{\ell\ne 0}\pb{J_u^{\pa{\ell}}}^2}\cdot B_u^{\pa{\ell}}\pa{|x-y|}+\OO_\prec\pbb{\frac{\Psi_u^{1/2}(0)}{1-u}\cdot \Psi_u^2(|x-y|)}.
			\end{aligned}
		\end{equation}
	\end{lemma}
	\begin{proof}
		A main difficulty in the proof is that it is not straightforward to extract an effective bound on $T_u-\Theta_u$ from a bound on the $(\cL_u-\cK_u)$-loops. This contrasts with the case of block random band matrices studied in \cite{dubova2025delocalizationnonmeanfieldrandommatrices}, where $\mathcal{E}_{u}^{\cL,\bsigma}$ involves only products of two $(\cL_u-\cK_u)$-loops. To overcome this issue, we use the following decomposition lemma, which rewrites $\mathcal{E}_{u}^{\cL,\bsigma}$ as a product of two $(\cL_u-\cK_u)$-loops up to a controllable error term. Its proof is deferred to \Cref{sec_estimates_error_terms}.
		
		\begin{lemma}\label{lemma_double_difference_decomposition}
			For any $\bsigma=(\sigma_1,\sigma_2)\in\h{-,+}^2$ and $x,y\in\ZL$, we have
			\begin{equation}\label{eq_decomposition_double_difference}
				\begin{aligned}                \mathcal{E}_{u,xy}^{\cL,\bsigma}=um\p{\sigma_1}m\p{\sigma_2}\sum_{a\in\ZL}(\cL-\cK)_{u,xa}^{\bsigma}(\cL-\cK)_{u,ay}^{\bsigma}+\OO_\prec\pbb{\frac{\Psi_u^{1/2}(0)}{1-u}\cdot \Psi_u^2(|x-y|)}.
				\end{aligned}
			\end{equation}
		\end{lemma}
		
		It remains to control the first term on the RHS of \eqref{eq_decomposition_double_difference}. By the definition of $J_u^{(\ell)}$ in \eqref{eq:JUl}, we have
		\begin{align}
			\absbb{\sum_{a\in\ZL}(\cL-\cK)_{u,xa}^{\bsigma}(\cL-\cK)_{u,ay}^{\bsigma}} &\leq \sum_{|a-x| < \ell}|(\cL-\cK)_{u,xa}^{\bsigma}|\cdot|(\cL-\cK)_{u,ay}^{\bsigma}|+\sum_{|a-x|\ge \ell}|(\cL-\cK)_{u,xa}^{\bsigma}|\cdot|(\cL-\cK)_{u,ay}^{\bsigma}|\nonumber\\
			\leq \mathbf 1_{\ell\ne 0}\pb{J_u^{\pa{\ell}}}^2\sum_{a\in\ZL}&B_u\p{|x-a|} B_u^{(\ell)}\p{|a-y|} +J_u^{\pa{\ell}}\sum_{a\in\ZL}B_u\p{\ell}\pa{\absa{\cL_{u,ay}^{\bsigma}}+\absa{\cK_{u,ay}^{\bsigma}}}.\label{eq:quad0}
		\end{align}
		Using the convolution inequality \eqref{eq:convolution_BBell}, we bound the first term on the RHS as \begin{align} 
			\mathbf 1_{\ell\ne 0}\pb{J_u^{\pa{\ell}}}^2\sum_{a\in\ZL}B_u\p{|x-a|} B_u^{(\ell)}\p{|a-y|}  
\lesssim \mathbf 1_{\ell\ne 0}\frac{\pb{J_u^{\pa{\ell}}}^2}{1-u} B_u^{(\ell)}\p{|x-y|}.\label{eq:quad1}\end{align}
		For the second term in \eqref{eq:quad0}, we use Ward's identity (together with Cauchy-Schwarz when $\sigma_1=\sigma_2$) and the local law \eqref{step_1_weak_local_law} to conclude that with high probability, 
		\begin{align}           
			J_u^{\pa{\ell}} B_u\p{\ell}\sum_{a\in\ZL}\pa{\absa{\cL_{u,ay}^{\bsigma}}+\absa{\cK_{u,ay}^{\bsigma}}}&\leq J_u^{\pa{\ell}} B_u\p{\ell} \max_{\bsigma\in\h{(-,+),(+,-)}}\sum_{a\in\ZL}\pa{\cL_{u,ay}^{\bsigma}+\cK_{u,ay}^{\bsigma}} \nonumber\\
			& =J_u^{\pa{\ell}} B_u\p{\ell}\pa{\frac{\avga{(\im G_u)S^{(y)}}}{\eta_u}+\frac{1}{1-u}}\lesssim \frac{J_u^{\pa{\ell}}}{1-u} B_u\p{\ell}.\label{eq:quad2}
		\end{align}
		Combining \eqref{eq:quad1} and \eqref{eq:quad2} with \eqref{eq:quad0}, and inserting the resulting bound into \eqref{eq_decomposition_double_difference}, we obtain \eqref{eq_bound_double_difference_term}.    
	\end{proof}

	The martingale term, i.e., the last term on the RHS of \eqref{cL_cK_evolution_step_2}, can be controlled using \Cref{lemma_martingale_to_cB} together with an extension of the loop-contraction inequality in Lemma 3.16 of \cite{dubova2025delocalizationnonmeanfieldrandommatrices}. The detailed proof is deferred to \Cref{sec_estimates_error_terms2}.
	
	\begin{lemma}\label{lemma_martingale_term}
		In the setting of \Cref{bound_double_difference_term}, the quadratic variation tensor defined in \eqref{def_tensor_cB} satisfies, with high probability and uniformly in $u\in[s,t]$ and $x,y\in\ZL$, 
		\begin{equation}\label{martingale_term_bound}
			\begin{aligned}
				\max_{\bsigma\in\h{-,+}^2}  \absa{(\cB\otimes\cB)_{u;xxyy}^{\cL,\bsigma}}\prec \frac{\Psi_u(0)}{1-u}\cdot \Psi_u^4\pa{|x-y|}.
			\end{aligned}
		\end{equation}
	\end{lemma}
	
	Finally, we estimate the light-weight term, i.e., the third term on the RHS of \eqref{cL_cK_evolution_step_2}. 
	
	\begin{lemma}\label{lemma_light_weight_term}
		In the setting of \Cref{bound_double_difference_term}, the light-weight term satisfies, with high probability and uniformly in $u\in[s,t]$ and $x,y\in\ZL$, 
		\begin{equation}\label{light_weight_bound}           \max_{\bsigma\in\h{-,+}^2}\absa{\cW_{u,xy}^{\cL,\bsigma}}\prec \frac{\Psi_u(0)}{1-u} \cdot \Psi_u^2\pa{|x-y|}.
		\end{equation}
	\end{lemma}
	\begin{proof}
		The proof relies on the following lemma, whose proof is deferred to \Cref{sec_graphical_argument}, where we develop new graphical tools tailored to this purpose.
		
		\begin{lemma}\label{claim_graph_bound_light_weight_u}
			For any $\sigma\in\h{-,+}$ and $x,y\in\ZL$, we have
			\begin{equation}\label{graph_bound_light_weight_u}
				\sum_{a,b\in\ZL}(G_u^{\sigma})_{xa}(G_u^{\sigma})_{ay}S_{ab}(G_u^{\sigma}-m^{\sigma})_{bb}=\opr{\frac{\Psi_u(0)}{1-u}\cdot \Psi_u\pa{|x-y|}}.
			\end{equation}
		\end{lemma}

		We rewrite the first term in \eqref{def_light_weight_term} as
		\begin{equation}\label{eq:LW1}
			\sum_{a,b,c,d\in\ZL}(G_u^{\sigma_2}-m^{\sigma_2})_{aa}S_{ab}(G_u^{\sigma_2})_{db}(G_u^{\sigma_2})_{bc}\cdot S_{cy}(G_u^{\sigma_1})_{cd}S_{dx}.
		\end{equation}
		Applying \eqref{graph_bound_light_weight_u}, we obtain
		\begin{align*}
			(\ref{eq:LW1}) &\prec\frac{\Psi_u(0)}{1-u}\sum_{c,d\in\ZL}\cS_{xd}\cS_{yc}[\delta_{cd}+\Psi_u(|c-d|)]\Psi_u(|c-d|)\prec \frac{\Psi_u(0)}{1-u}\Psi^2_u(|x-y|),
		\end{align*}
		where, in the first step, we used the assumption \eqref{double_difference_input} and \Cref{2_loop_to_1_chain} to bound $(G_u^{\sigma_1})_{cd}$ by $\delta_{cd}+\Psi_u(|c-d|)$, and in the second step, we used the bound \eqref{polydecay_Psiu3} and \(\sum_{c,d}\cS_{xd}\delta_{cd}\Psi_u(|c-d|)\cS_{c y}\leq\sum_{c}\cS_{xc}\cS_{c y} \lesssim \cS_{xy}. \) The second term in \eqref{def_light_weight_term} can be treated analogously. This completes the proof of \Cref{lemma_light_weight_term}.
	\end{proof}

	\subsection{Step 2: Sharp bounds for \texorpdfstring{$\cal L$}{L}-loops and sharp local law}\label{sec:Gronwall}
	
	In this step, we prove the estimates \eqref{step_2_weak_estimate_loops} and \eqref{step_2_sharp_local_law} by analyzing the equation \eqref{cL_cK_evolution_step_2} together with Lemmas \ref{bound_leading_term}, \ref{bound_double_difference_term}, \ref{lemma_martingale_term}, and \ref{lemma_light_weight_term}. The sharp local law \eqref{step_2_sharp_local_law} follows immediately from \eqref{step_2_sharp_bound_loops} (which itself is a consequence of \eqref{step_2_weak_estimate_loops}) and \Cref{2_loop_to_1_chain}. Therefore, it suffices to prove \eqref{step_2_weak_estimate_loops}. We follow the strategy of Section 3.3 in \cite{dubova2025delocalizationnonmeanfieldrandommatrices}. We first establish the $\max$-estimate
	\begin{equation}\label{step_2_L_infty_bound}
		\begin{aligned}
			\max_{\bsigma\in\h{-,+}^2}\|\cL_u^{\bsigma}-\cK_u^{\bsigma}\|_{\max}\prec \pa{\frac{1-s}{1-u}}^{C}B_u^{6/5}(0),\quad \forall u\in\qa{s,t}.
		\end{aligned}
	\end{equation}
	It suffices to establish the following estimate for any fixed integer $k\geq 0$: there exist a sufficiently large constant $C\equiv C(k)>0$ and a sufficiently small constant $\fc_{\alpha}(k)>0$ such that, if the constant $c_{\alpha}$ in \eqref{bootstrap_assumption} satisfies $c_{\alpha}\le \fc_{\alpha}(k)$, then
	\begin{equation}\label{step_2_L_infty_bound_induction}
		\begin{aligned}
			\max_{\bsigma\in\h{-,+}^2}\|\cL_u^{\bsigma}-\cK_u^{\bsigma}\|_{\max}\prec \pa{\frac{1-s}{1-u}}^{C}B_u^{6/5}(0)+\pa{\frac{1-s}{1-u}}^{C}\qa{B_s\p{0}}^{(5/4)^k/2},\quad \forall u\in\qa{s,t}.
		\end{aligned}
	\end{equation}
	Given this bound, taking $k=4$ and noting that $(5/4)^4/2>6/5$, the second term on the RHS of \eqref{step_2_L_infty_bound_induction} is dominated by the first term, provided that the constant $c_{\alpha}$ in \eqref{bootstrap_assumption} is chosen sufficiently small. 
	The case $k=0$ of \eqref{step_2_L_infty_bound_induction}, with $C=1$, follows directly from the bound \eqref{step_1_priori_bound}. Next, assume that \eqref{step_2_L_infty_bound_induction} holds for some fixed integer $k\geq 0$ and a constant $C\equiv C(k)>0$. Choose $\fc_{\alpha}(k+1)$ sufficiently small so that, for any $0<c_{\alpha}\le\fc_{\alpha}(k+1)$, we have \smash{$[(1-s)/(1-u)]^CB_u^{1/5}(0)\leq 1$} for all $u\in\qa{s,t}$. Then we obtain 
	\begin{equation}\label{continuity_max_bound_loops}        \max_{\bsigma\in\h{-,+}^2}\|\cL_u^{\bsigma}\|_{\max}\prec\Psi_u^2(0), \quad \text{with}\quad  \Psi_u^2(\cdot)\equiv\Psi_u^2(0):=B_u(0)+\pa{\frac{1-s}{1-u}}^{C}\qa{B_s\p{0}}^{(5/4)^k/2}. 
	\end{equation}
	Applying the induction hypothesis in \Cref{lemma_bootstrap} to the first term on the RHS of \eqref{cL_cK_evolution_step_2}, \Cref{bound_leading_term} (with $\ell=0$) to the second term, \Cref{lemma_light_weight_term} to the third term, \Cref{bound_double_difference_term} (with $\ell=0$) to the fourth term, and \Cref{lemma_martingale_term} together with \eqref{bound_martingale_to_cB} to the martingale term, we obtain that for any $\bsigma\in\{-,+\}^2$, $x,y\in\ZL$, and $r\in[s,t]$, the following estimate holds with high probability:
	\begin{equation*}
		\begin{aligned}
			\absa{(\cL-\cK)_{r,xy}^{\bsigma}}\le &~        \int_{s}^{r}\frac{C_0}{1-u}J_u^{\pa{0}}B_u(0)\,\rd u +\opr{B_s^{1/5}(0)B_s(|x-y|)+\int_{s}^r\frac{\Psi_u^{5/2}(0)}{1-u}\,\rd u+\pa{\int_s^r\frac{\Psi_u^5(0)}{1-u}\,\rd u}^{1/2}}\\
			\leq&~B_r(0)\int_s^r\frac{C_0}{1-u}J_u^{\pa{0}}\,\rd u+\opr{B_r^{6/5}(0)+\pa{\frac{1-s}{1-r}}^{5C/4}\qa{B_s\p{0}}^{(5/4)^{k+1}/2}},
		\end{aligned}
	\end{equation*}
	for some constant $C_0>0$ depending only on $\kappa$, $\fc$, and $\alpha$. Here we used $\Psi_u(0)$ from \eqref{continuity_max_bound_loops} and the elementary estimate  $\int_s^{r}\p{1-u}^{-1}\rd u\lesssim \log N\prec 1.$ 
	Dividing both sides by $B_r(0)$ and taking the supremum over $\bsigma$ and $(x,y)$, we obtain
	\begin{equation}\label{Gronwall_inequality0}
		J_r^{\pa{0}}\leq \int_s^{r} \frac{C_0}{1-u}J_u^{\pa{0}}\dd u +\opr{B_r^{1/5}(0)+\pa{\frac{1-s}{1-r}}^{5C/4}\qa{B_s\p{0}}^{(5/4)^{k+1}/2-1}},\quad \forall r\in\qa{s,t}.
	\end{equation}
	Recall Gr{\"o}nwall's inequality: for any non-decreasing function $a$ and non-negative function $b$,
	\begin{equation}\label{Gronwall_inequality}
		\begin{aligned}
			f(r)\leq a(r)+\int_s^rb(u)f(u)\,\rd u,\quad \forall r\in\qa{s,t}\implies f(r)\leq a(r)\exp\pa{\int_s^r b(u)\,\rd u} .
		\end{aligned}
	\end{equation}
	Applying \eqref{Gronwall_inequality} to \eqref{Gronwall_inequality0}, we obtain
	\begin{equation*}
		\begin{aligned}
			J_r^{\pa{0}}\prec \pa{\frac{1-s}{1-r}}^{C_0}\qa{B_r^{1/5}(0)+\pa{\frac{1-s}{1-r}}^{5C/4}\qa{B_s\p{0}}^{(5/4)^{k+1}/2-1}},\quad \forall r\in\qa{s,t},
		\end{aligned}
	\end{equation*}
	which, together with \eqref{eq:JUl}, establishes the estimate \eqref{step_2_L_infty_bound_induction} for a sufficiently large constant $C(k+1)>0$.
	
	  Using \eqref{step_2_L_infty_bound} and choosing $c_\alpha$ in \eqref{bootstrap_assumption} sufficiently small, we obtain
	\begin{equation*}        \max_{\bsigma\in\h{-,+}^2}\absa{\cL_{u,xy}^{\bsigma}}\prec B_u^{\pa{0}}\pa{|x-y|},\quad \forall u\in\qa{s,t},
	\end{equation*}
	where we also used that $\max_{\bsigma\in\h{-,+}^2}|\cK_{u,xy}^{\bsigma}|\prec B_u(|x-y|)$, which follows from $\cK_u=S\Theta_u$ together with the bounds \eqref{Theta_bound_opposite_charge} and \eqref{Theta_bound_equaled_charge}.
	We now perform an induction on the scale of length parameters. Suppose that, for a collection of length parameters $\{K_u:u\in[s,t]\}$ depending continuously on $u\in\qa{s,t}$, we have
	\begin{equation}\label{step_2_cL_bound_input}
		\max_{\bsigma\in\h{(-,+),(+,-)}}\cL_{u,xy}^{\bsigma}\prec B_u^{\pa{K_u}}\pa{|x-y|}=:\Psi_u^2\pa{|x-y|},\quad \forall u\in\qa{s,t},
	\end{equation}
	and that $B_u^{\pa{K_u}}\pa{|x-y|}$ is \emph{non-decreasing in $u$} for every $x,y\in\ZL$. Under \eqref{step_2_cL_bound_input}, we have $\Psi_u(0)=B_u^{1/2}\pa{0}$. Define the stopping time
	\begin{equation}\label{eq:stoppingtime_step2}
		\tau=t\wedge T,\quad \txt{with} \quad T:=\inf\ha{u\geq s: J_u^{\pa{K_u}}\geq B_u^{1/6}(0)},
	\end{equation}
   where we adopt the convention $\inf\emptyset=\infty$.
	By the induction hypothesis in \Cref{lemma_bootstrap} at time $s$, we know that \smash{$J_s^{\pa{K_s}}\prec B_s^{1/5}(0)$}, hence $\tau>s$ with high probability (unless $t=s$). Combining \eqref{step_2_cL_bound_input} with Lemmas
	\ref{bound_leading_term},
	\ref{bound_double_difference_term},
	\ref{lemma_martingale_term} (together with \eqref{bound_martingale_to_cB}),
	and \ref{lemma_light_weight_term}, we obtain from the evolution equation
	\eqref{cL_cK_evolution_step_2} that for any
	$\bsigma\in\{-,+\}^2$, $x,y\in\ZL$, and $r\in[s,\tau]$, the following bound holds with high probability for a constant $C>0$: \begin{equation*}
		\begin{aligned}
			\frac{\absa{(\cL-\cK)_{r,xy}^{\bsigma}}}{B_r^{(K_r)}\pa{|x-y|}}&\leq{ \int_s^r \frac{C/2}{1-u} \pa{J_u^{\pa{K_u}}+\pb{J_u^{\pa{K_u}}}^2}\rd u+\opr{B_r^{1/5}(0)+\int_{s}^r\frac{B_r^{1/4}(0)}{1-u}\,\rd u+\pa{\int_s^r\frac{B_r^{1/2}(0)}{1-u}\,\rd u}^{1/2}}}\\
			&\leq \int_s^r \frac{C}{1-u} J_u^{\pa{K_u}}\,\rd u+\opr{ B_r^{1/5}\pa{0}},
		\end{aligned}
	\end{equation*}
	where the second step uses $J_u^{\pa{K_u}}\le B_u^{1/6}(0)$ on $[s,\tau]$ by the definition \eqref{eq:stoppingtime_step2}. 
	Applying Gr{\"o}nwall's inequality \eqref{Gronwall_inequality} yields for any fixed constant $\e>0$, with high probability,
	\begin{equation}
		\begin{aligned}
			J_r^{\pa{K_r}}\leq \pa{\frac{1-s}{1-r}}^{C}\cdot W^{\varepsilon} B_r^{1/5}\pa{0},\quad \forall r\in\qa{s,\tau}.
		\end{aligned}
	\end{equation}
	A standard continuity argument then implies that $\tau=t$ with high probability, provided $\e$ is chosen sufficiently small. Consequently, for any $\bsigma\in\{-,+\}^2$ and $x,y\in\ZL$,
	\begin{equation}\label{step_2_iterate_result}
		\begin{aligned}
			(\cL-\cK)_{u,xy}^{\bsigma}\prec \pa{\frac{1-s}{1-u}}^{C}B_u^{1/5}(0)\cdot B_u^{\pa{K_u}}\pa{|x-y|},\quad \forall u\in\qa{s,t}.
		\end{aligned}
	\end{equation}
	To complete the induction, define a new collection of length parameters $\{K_u':u\in[s,t]\}$ as the unique positive solution of
	\begin{equation*}
		\begin{aligned}
			B_u(K_u')=B_u^{1/6}(0)B_u(K_u).
		\end{aligned}
	\end{equation*}
	With this choice, and taking $c_\alpha$ in \eqref{bootstrap_assumption} sufficiently small, \eqref{step_2_iterate_result} implies
	\begin{align}        \absa{\cL_{u,xy}^{\bsigma}}&\leq\absa{\cK_{u,xy}^{\bsigma}}+\absa{(\cL-\cK)_{u,xy}^{\bsigma}}\prec B_u\pa{|x-y|}+\pa{\frac{1-s}{1-u}}^{C} B_u^{1/5}(0)\cdot\qa{B_u\pa{|x-y|}+B_u(K_u)}\nonumber\\
		&\lesssim B_u\pa{|x-y|}+B_u^{1/6}(0)B_u(K_u)\asymp B_u^{\p{K_u'}}\pa{|x-y|},\quad \forall u\in\qa{s,t}.\label{step_2_cL_bound_output}
	\end{align}
	Fix any large constant $D>0$. Iterating the improvement from \eqref{step_2_cL_bound_input} to \eqref{step_2_cL_bound_output} $\OO(1)$ times, we eventually obtain \eqref{step_2_cL_bound_input} for a collection $\{K_u:u\in[s,t]\}$ such that $B_u(K_u)\le W^{-D}$. For $D$ sufficiently large, this yields
	\begin{equation*}
		\begin{aligned}        \max_{\bsigma\in\h{-,+}^2}\absa{\cL_{u,xy}^{\bsigma}}\prec B_u^{\pa{K_u}}\pa{|x-y|}\lesssim B_u(|x-y|)+W^{-D}\lesssim B_u(|x-y|).     \end{aligned}
	\end{equation*}
	Plugging this back into \eqref{step_2_iterate_result} completes the proof of \eqref{step_2_weak_estimate_loops}.

	\subsection{Step 3: $\cL$-loop estimate}
	
	In this step, we have already established the sharp $\cL$-loop bound \eqref{step_2_sharp_bound_loops}. Hence, we may choose the control parameter $\Psi_u^2(|x-y|)=B_u\pa{|x-y|}$ for $x,y\in \ZL$ and $u\in\qa{s,t}$. 
	To prove \eqref{step_3_estimate_loops}, we use the integrated evolution equation \eqref{integrated_cL_cK_evolution} with $\tau=r$ for $r\in[s,t]$. We bound the last three terms on the RHS of \eqref{integrated_cL_cK_evolution} by inserting the sharp bound \eqref{step_2_sharp_bound_loops} into Lemmas \ref{bound_double_difference_term}, \ref{lemma_martingale_term}, and \ref{lemma_light_weight_term}, and then applying the evolution kernel bound \eqref{bound_cU_s_t_1}. 
	Specifically, for the light-weight term, by \Cref{lemma_light_weight_term} and \eqref{bound_cU_s_t_1}, the following bound holds for any $x,y\in\ZL$ and $s\le u\le r\le t$,
	\begin{equation}        \qa{\cU_{u,r}^{\bsigma}\circ\cW_{u}^{\cL,\bsigma}}_{xy}\prec \frac{B_r^{1/2}(0)}{1-r}\cdot B_r(|x-y|).\label{eq:LW_step3}
	\end{equation}
	Integrating over $u$ gives
	\begin{equation}\label{step_3_light_weight}
		\begin{aligned}
			\int_{s}^r\qa{\cU_{u,r}^{\bsigma}\circ\cW_{u}^{\cL,\bsigma}}_{xy}\,\rd u\prec \frac{1-s}{1-r}\cdot B_r^{1/2}(0)B_r(|x-y|)\leq B_r^{1/5}(0)B_r\pa{|x-y|},
		\end{aligned}
	\end{equation}
	provided $c_\alpha$ in \eqref{bootstrap_assumption} is chosen sufficiently small. For the quadratic error term, applying \eqref{eq_decomposition_double_difference} with $\Psi_u^2(|x-y|)=B_u(|x-y|)$ and using the weak estimate \eqref{step_2_weak_estimate_loops}, we obtain
	\begin{equation*}
		\begin{aligned}
			\mathcal{E}_{u,xy}^{\cL,\bsigma}&\prec B_u^{2/5}(0)\pa{\frac{1-s}{1-u}}^{2C_\alpha}\sum_{a}B_u\pa{|x-a|}B_u\pa{|a-y|}+\frac{B_u^{1/4}(0)}{1-u}\cdot B_u\pa{|x-y|}\\
			&\lesssim\frac{B_u^{1/4}(0)}{1-u}\cdot B_u\pa{|x-y|},
		\end{aligned}
	\end{equation*}
	where we used the convolution inequality \eqref{convolution_2_B} and again chose $c_\alpha$ sufficiently small. Applying \eqref{bound_cU_s_t_1} and integrating in $u$, we get
	\begin{equation}\label{step_3_double_difference}
		\begin{aligned}
			\int_s^r \qa{\cU_{u,r}^{\bsigma}\circ \mathcal{E}_u^{\cL,\bsigma}}_{xy}\,\rd u\prec B_r^{1/5}(0)B_r(|x-y|).
		\end{aligned}
	\end{equation}
	To control the martingale term, we apply \Cref{lemma_martingale_to_cB}. 
	By Cauchy-Schwarz (CS), we have
	\begin{equation*}
		\begin{aligned}
			\absa{(\cB\otimes\cB)_{u;x_1x_2y_1y_2}^{\cL,\bsigma}}\lesssim \qa{(\cB\otimes\cB)_{u;x_1x_1y_1y_1}^{\cL,\bsigma}}^{1/2}\qa{(\cB\otimes\cB)_{u;x_2x_2y_2y_2}^{\cL,\bsigma}}^{1/2} .
		\end{aligned}
	\end{equation*}
	Let $X_{u,\cB}^{\cL,\bsigma}(a,b):=\qb{(\cB\otimes\cB)_{u;aabb}^{\cL,\bsigma}}^{1/2}$. Since $\absa{\cU_{u,r}^{\bsigma}(a,b)}\leq  \cU_{u,r}^{(-,+)}(a,b)$, we obtain from
	\eqref{bound_cU_s_t_1} and \eqref{martingale_term_bound} that 
	\begin{align}        \ha{\qa{\cU_{u,r}^{\bsigma}\otimes\pa{\cU_{u,r}^{\bsigma}}^*\otimes\cU_{u,r}^{\bsigma}\otimes\pa{\cU_{u,r}^{\bsigma}}^*}\circ(\cB\otimes\cB)_{u}^{\cL,\bsigma}}_{xxyy}&\lesssim \qa{\cU^{(-,+)}_{u,r}\circ X_{u,\cB}^{\cL,\bsigma}}^2_{xy} \nonumber\\
		&\prec \pa{\frac{1-s}{1-r}}^2\frac{B_r^{1/2}(0)}{1-u}B_r^2(|x-y|).\label{eq:UUUUBB}
	\end{align}
	Plugging this into \eqref{bound_martingale_to_cB} (with $A_u=B_u=\cU_{u,r}^{\bsigma}$) and applying Markov's inequality yields
	\begin{equation}\label{step_3_martingale}
		\int_{s}^{r}\qa{\cU_{u,r}^{\bsigma}\circ\rd \cB_{u}^{\cL,\bsigma}}_{xy}\prec \frac{1-s}{1-r}B_r^{1/4}(0)B_r(|x-y|)\leq B_r^{1/5}(0)B_r(|x-y|),
	\end{equation}
	again for sufficiently small $c_\alpha$.
	
	It remains to bound the first term in \eqref{integrated_cL_cK_evolution} with the contribution from the initial condition. If $\sigma_1=\sigma_2$, then by the induction hypothesis \eqref{flow_loop_local_law} and \eqref{bound_cU_s_t_2},
	\begin{equation*}
		\begin{aligned}
			\qa{\cU_{s,r}^{\bsigma}\circ(\cL-\cK)_{s}^{\bsigma}}_{xy}\prec B_r^{1/5}(0)B_r(|x-y|).
		\end{aligned}
	\end{equation*}
	If $\sigma_1\neq\sigma_2$, we decompose
	\begin{align}
		(\cL-\cK)_{s,xy}&=\cA_1(x,y)+\cA_2(x,y),\label{step_3_decomposition_leading}
	\end{align}
	where, for brevity, we suppress the superscript $\bsigma$. The term $\cal A_2$ is defined by
	\begin{equation}
		\cA_2(x,y):=\sum_w(\cL-\cK)_{s,xw}\cdot (1-s)\Theta_{s,xy}=\frac{\im \avg{(G_s-m)S^{\pa{x}}}}{\eta_s}(1-s)\Theta_{s,xy},\label{eq:A2_Ward}
	\end{equation}
	and $\cA_1$ is the remainder. In the second equality above, we used Ward's identity and the identity $\im m/\eta_s=(1-s)^{-1}$ by \eqref{re_im_z_t}.
	Using the assumption \eqref{Theta_bound_opposite_charge}, along with \eqref{flow_local_law} and \eqref{flow_loop_local_law} at time $s$, we obtain 
	\begin{align}
		&|\cA_1\pa{x,y}|\prec B_s^{1/5}(0)B_s(|x-y|)+B_s(0)B_s(|x-y|)\lesssim B_s^{1/5}(0)B_s\pa{|x-y|},\label{eq:A1xy}\\
		&|\cA_2\pa{x,y}|\prec B_s(0)B_s(|x-y|).\label{eq:A2xy}
	\end{align}
	Moreover, the definition \eqref{eq:A2_Ward} and $(1-s)\sum_y\Theta_{s,xy}=1$ imply $\sum_{y}\cA_1\pa{x,y}=0$ for all $x\in\ZL$. Applying \eqref{bound_cU_s_t_1} to $\cU_{s,r}^{\bsigma}\circ\cA_2$
	and \eqref{bound_cU_s_t_2} to $\cU_{s,r}^{\bsigma}\circ\cA_1$, we obtain
	\begin{equation}\label{step_3_leading}
		\begin{aligned}
			\qa{\cU_{s,r}^{\bsigma}\circ(\cL-\cK)_s^{\bsigma}}_{xy}\prec B_r^{1/5}(0)B_r(|x-y|)+{\frac{1-s}{1-r}}B_r(0)B_r(|x-y|)\lesssim B_r^{1/5}(0)B_r(|x-y|),
		\end{aligned}
	\end{equation}
	again choosing $c_\alpha$ sufficiently small.
	
	Finally, combining \eqref{step_3_light_weight}, \eqref{step_3_double_difference}, \eqref{step_3_martingale}, and \eqref{step_3_leading}, the integrated evolution equation \eqref{integrated_cL_cK_evolution} yields 
	\begin{equation*}
		\begin{aligned}
			(\cL-\cK)_{r,xy}^{\bsigma}\prec B_r^{1/5}(0)B_r(|x-y|),\quad \forall r\in\qa{s,t},
		\end{aligned}
	\end{equation*}
	which completes Step 3 in the proof of \Cref{lemma_bootstrap}.

	\subsection{Proof of \Cref{lemma_double_difference_decomposition}}\label{sec_estimates_error_terms}
	For notational simplicity, we present the proof only for the most delicate case $\bsigma = (-,+)$; the argument for general charges $\bsigma \in \{-,+\}^2$ is identical. Throughout the proof, we suppress the charge parameter $\bsigma$ in the superscript and, whenever no confusion arises, also omit the time parameter $u$ from the notation. For instance, we write $G \equiv G_u$, $\cL \equiv \cL_u^{\bsigma}$, $\cK \equiv \cK_u^{\bsigma}$, $T \equiv T_u^{\bsigma}$, and $\Theta \equiv \Theta_u^{\bsigma}$.

	To establish the estimate \eqref{eq_decomposition_double_difference}, we decompose
	\begin{equation*}
		\begin{aligned}
			\mathcal{E}_{xy}^{\cL,\bsigma}=\sum_{a\in\ZL}(T-\Theta)_{xa}(\cL-\cK)_{ay}=\sum_{a\in\ZL}(T-\Theta)_{xa}\cL_{ay}-\sum_{a\in\ZL}(T-\Theta)_{xa}\cK_{ay}.
		\end{aligned}
	\end{equation*}
	It therefore suffices to prove that    
	\begin{align}
		\sum_{a\in\ZL}(T-\Theta)_{xa}\cL_{ay}=u\sum_{a\in\ZL}(\cL-\cK)_{xa}\cL_{ay}+\opr{\frac{\Psi_u^{1/2}(0)}{1-u}\cdot\Psi_u^2\pa{|x-y|}},\label{eq_first_decomposition_double_diffenrence}\\
		\sum_{a\in\ZL}(T-\Theta)_{xa}\cK_{ay}=u\sum_{a\in\ZL}(\cL-\cK)_{xa}\cK_{ay}+\opr{\frac{\Psi_u^{1/2}(0)}{1-u}\cdot\Psi_u^2\pa{|x-y|}}.\label{eq_2nd_decomposition_double_diffenrence}
	\end{align}
	We present only the proof of \eqref{eq_first_decomposition_double_diffenrence}, since \eqref{eq_2nd_decomposition_double_diffenrence} follows from a similar, and much simpler, argument. To this end, define
	\begin{equation}
		\begin{aligned}
			\Xi_{xy}\equiv \Xi_{u,xy}:=\sum_{a\in\ZL}(T-\Theta)_{xa}\cL_{ay}-u\sum_{a\in\ZL}(\cL-\cK)_{xa}\cL_{ay}.
		\end{aligned}
	\end{equation}
	We prove \eqref{eq_first_decomposition_double_diffenrence} by establishing the following high-moment bound: for any fixed $p\in\N$,
	\begin{equation}\label{high_moment_bound_Xi_double_difference}
		\begin{aligned}
			\E|\Xi_{xy}|^{2p}\prec \qbb{\frac{\Psi_u^{1/2}(0)}{1-u}\cdot\Psi_u^2\pa{|x-y|}}^{2p}.
		\end{aligned}
	\end{equation}
	
	Using the definition of $G_u$ and the self-consistent equation \( m=m_u(z_u)=-(z_u+um_u(z_u))^{-1}\), we verify the identity
	\begin{equation}
		G_u=m-G_u(H_u+um)m=m-m(H_u+um)G_u.\label{eq:Gu-m}
	\end{equation} 
	Invoking this relation, we rewrite
	\begin{align}
		\E|\Xi_{xy}|^{2p}&=\E\sum_{i,j\in\ZL}S_{xi}G_{ij}G^*_{ji}\cL_{jy}\Xi_{xy}^{p-1}\ol{\Xi}_{xy}^{p} -\E\sum_{a\in\ZL}\Theta_{xa}\cL_{ay}\Xi_{xy}^{p-1}\ol{\Xi}_{xy}^{p}-u\E\sum_{a\in\ZL}(\cL-\cK)_{xa}\cL_{ay}\Xi_{xy}^{p-1}\ol{\Xi}_{xy}^{p}\nonumber\\
		&=-m\E\sum_{i,j,\beta}S_{xi}G_{i\beta}(H_u)_{\beta j}G^*_{ji}\cL_{jy}\Xi_{xy}^{p-1}\ol{\Xi}_{xy}^{p} -um^2\E\sum_{i,j}S_{xi}G_{ij}G^*_{ji}\cL_{jy}\Xi_{xy}^{p-1}\ol{\Xi}_{xy}^{p}\nonumber\\
		&\quad+m\E\sum_{i}S_{xi}G^*_{ii}\cL_{iy}\Xi_{xy}^{p-1}\ol{\Xi}_{xy}^{p}  -\E\sum_{a}\Theta_{xa}\cL_{ay}\Xi_{xy}^{p-1}\ol{\Xi}_{xy}^{p}-u\E\sum_{a}(\cL-\cK)_{xa}\cL_{ay}\Xi_{xy}^{p-1}\ol{\Xi}_{xy}^{p}.\label{eq_Xi_H}
	\end{align}
	Applying Gaussian integration by parts with respect to the entries of $H_u$, we expand the first term as  
	\begin{equation*}
		\begin{aligned}
			&\quad -m\E\sum_{i,j,\beta\in\ZL}S_{xi}G_{i\beta}(H_u)_{\beta j}G^*_{ji}\cL_{jy}\Xi_{xy}^{p-1}\ol{\Xi}_{xy}^{p}\\
			&= um\E\sum_{i,j,\beta\in\ZL}S_{xi}S_{\beta j}\pa{G_{ij}G_{\beta\beta}G^*_{ji}+G_{i\beta}G^*_{jj}G^*_{\beta i}}\cL_{jy}\Xi_{xy}^{p-1}\ol{\Xi}_{xy}^{p}\\
			& -um\E\sum_{i,j,\beta\in\ZL}S_{xi}S_{\beta j}G_{i\beta}G^*_{ji}\partial_{j\beta}\pa{\cL_{jy}\Xi_{xy}^{p-1}\ol{\Xi}_{xy}^{p}}.
		\end{aligned}
	\end{equation*}
	Substituting this expression into \eqref{eq_Xi_H}, and using $\Theta-u\cK=S$ together with $|m|=1$, we obtain   \begin{equation}\label{high_moment_decomposition_double_difference}
		\begin{aligned}
			&\E|\Xi_{xy}|^{2p}=um\E\sum_{j,\beta\in\ZL}(G-m)_{\beta\beta}S_{\beta j}T_{xj}\cL_{jy}\Xi_{xy}^{p-1}\ol{\Xi}_{xy}^{p}\\
			&+um\E\sum_{j\in\ZL}\pa{G^*-\ol{m}}_{jj}\cL_{xj}\cL_{jy}\Xi_{xy}^{p-1}\ol{\Xi}_{xy}^{p}+m\E\sum_{i\in\ZL}S_{xi}(G^*-\ol{m})_{ii}\cL_{iy}\Xi_{xy}^{p-1}\ol{\Xi}_{xy}^{p}\\
			&-um\E\sum_{i,j,\beta\in\ZL}S_{xi}S_{\beta j}G_{i\beta}G^*_{ji}\partial_{j\beta}\pa{\cL_{jy}\Xi_{xy}^{p-1}\ol{\Xi}_{xy}^{p}}.
		\end{aligned}
	\end{equation}
By \Cref{2_loop_to_1_chain}, together with the assumptions \eqref{double_difference_input} and \eqref{polydecay_Psiu}--\eqref{polydecay_Psiu3}, we obtain the following estimates:
	\begin{equation}
		\|G-m\|_{\max}\prec \Psi_u(0),\quad \mathbf 1_{x\ne y}|G_{xy}|\prec \Psi_{u}(|x-y|),\quad T_{xy}\prec \Psi_{u}^2(|x-y|).\label{eq:localPsiu}
	\end{equation}
	In particular, \eqref{eq:localPsiu} implies that
	\begin{equation*}
		\begin{aligned}
			\sum_{j,\beta\in\ZL}(G-m)_{\beta\beta}S_{\beta j}T_{xj}\cL_{jy}&\prec \Psi_u(0)\cdot \sum_{|j-x|\leq |x-y|/2}T_{xj}\cL_{jy}+\Psi_u(0)\cdot\sum_{|j-x|>|x-y|/2}T_{xj}\cL_{jy} \\
			&\prec \Psi_u(0)\cdot\Psi_u^2\pa{|x-y|}\sum_{j\in\ZL}\pa{T_{xj}+\cL_{jy}}\prec \frac{\Psi_u(0)}{1-u}\cdot \Psi_u^2\pa{|x-y|},
		\end{aligned}
	\end{equation*}
	where in the last step we used Ward's identity. Consequently, the first term on the RHS of \eqref{high_moment_decomposition_double_difference} can be bounded as 
	\begin{equation}\label{eq:boundfirstterm}
		um\E\sum_{j,\beta\in\ZL}(G-m)_{\beta\beta}S_{\beta j}T_{xj}\cL_{jy}\Xi_{xy}^{p-1}\ol{\Xi}_{xy}^{p}\prec \frac{\Psi_u(0)}{1-u}\Psi_u^2\pa{|x-y|}\cdot \pa{\E|\Xi_{xy}|^{2p}}^{(2p-1)/2p}.
	\end{equation}
	The second and third terms on the RHS of \eqref{high_moment_decomposition_double_difference} admit the same bound by an analogous argument.
	For the last term in
	\eqref{high_moment_decomposition_double_difference},
	we establish the bounds
	\begin{align}
		&\sum_{i,j,\beta\in\ZL}S_{xi}S_{\beta j}\absa{G_{i\beta}G^*_{ji}\partial_{j\beta}\pa{\cL_{jy}}}\prec \frac{\Psi_u(0)}{1-u}\cdot \Psi_u^2\pa{|x-y|},\label{bound_pull_in_edge_term_1}\\
		&\sum_{i,j,\beta\in\ZL}S_{xi}S_{\beta j}\qa{\absa{G_{i\beta}G^*_{ji}\cL_{jy}\partial_{j\beta}\pa{\Xi_{xy}}}+\absa{G_{i\beta}G^*_{ji}\cL_{jy}\partial_{j\beta}\pa{\ol{\Xi}_{xy}}}}\prec \frac{\Psi_u(0)}{(1-u)^2} \cdot \Psi_u^4\pa{|x-y|}. \label{bound_pull_in_edge_term_2}
	\end{align}
	Although these estimates can be verified directly by expanding the derivatives $\partial_{j\beta}(\cL_{jy})$, $\partial_{j\beta}(\Xi_{xy})$, and $\partial_{j\beta}(\ol{\Xi}_{xy})$, we provide a more concise argument in \Cref{sec:defer_proof} using the graphical tools developed in \Cref{sec_graphical_argument}. 
	Combining \eqref{bound_pull_in_edge_term_1} and \eqref{bound_pull_in_edge_term_2} with Hölder's inequality yields
	\begin{equation}\label{eq:boundlastterm}
		um\E\sum_{i,j,\beta\in\ZL}S_{xi}S_{\beta j}G_{i\beta}G^*_{ji}\partial_{j\beta}\pa{\cL_{jy}\Xi_{xy}^{p-1}\ol{\Xi}_{xy}^{p}}\prec \sum_{k=1,2}\qa{\frac{\Psi_u^{1/2}(0)}{1-u}\Psi_u^2\pa{|x-y|}}^k\cdot \pa{\E|\Xi|^{2p}}^{(2p-k)/2p}.
	\end{equation}
	Substituting \eqref{eq:boundfirstterm} and \eqref{eq:boundlastterm} into
	\eqref{high_moment_decomposition_double_difference}, and applying Young's inequality, we obtain \eqref{high_moment_bound_Xi_double_difference}. By Markov's inequality, this further implies
	\begin{equation*}
		\begin{aligned}
			\Xi_{xy}\prec{\frac{1}{1-u}\Psi_u^{1/2}(0)\Psi_u^2\pa{|x-y|}}.
		\end{aligned}
	\end{equation*}
	This completes the proof of \eqref{eq_decomposition_double_difference}.

	\subsection{Proof of \Cref{lemma_martingale_term}}\label{sec_estimates_error_terms2}
	
	It suffices to bound the first term in \eqref{def_tensor_cB}, since the second term can be treated analogously by symmetry. As before, we present the proof only for the most delicate case $\bsigma=(+,-)$; the argument for general charges $\bsigma\in\{-,+\}^2$ is identical. In this case, we decompose the first term in \eqref{def_tensor_cB} (with $t$ replaced by $u$) as
	\begin{equation*}
		\begin{aligned}
			\sum_{a\in\ZL}(G_uS^{\pa{x}}G_u^*S^{\pa{y}}G_uS^{\pa{a}}G_u^*S^{\pa{y}}G_uS^{\pa{x}}G_u^*)_{aa}=I_{\leq}+I_{>},
		\end{aligned}
	\end{equation*}
	where $I_\leq$ is defined as
	\begin{equation*}
		\begin{aligned}
			I_\leq:=\sum_{a:|x-a|\leq |x-y|/2}(G_uS^{\pa{x}}G_u^*S^{\pa{y}}G_uS^{\pa{a}}G_u^*S^{\pa{y}}G_uS^{\pa{x}}G_u^*)_{aa},
		\end{aligned}
	\end{equation*}
	and $I_{>}$ denotes the complementary contribution over the region \(\{a:|x-a|>|x-y|/2\}\). We bound $I_\leq$ as
	\begin{equation}\label{bound_I_leq}
		I_{\leq}=\sum_{a:|x-a|\leq |x-y|/2}\psi_{a,\leq}A_{a,\leq}\psi_{a,\leq}^*\leq \sum_{a:|x-a|\leq |x-y|/2}\|A_{a,\leq}\|\cdot\|\psi_{a,\leq}\|_2^2,
	\end{equation}
	where the (row) vector $\psi_{a,\leq}\in \C^N$ and the matrix $A_{a,\leq}\in \C^{N\times N}$ are defined as 
	\begin{equation*}
		\begin{aligned}
			\psi_{a,\leq}\pa{i}:=(G_uS^{\pa{x}}G_u^*)_{ai}\sqrt{S_{iy}}, \quad A_{a,\leq}\pa{i,j}:=\sqrt{S_{iy}}(G_uS^{\pa{a}}G_u^*)_{ij}\sqrt{S_{jy}},\quad  \forall i,j\in\ZL.
		\end{aligned}
	\end{equation*}
	Bounding the operator norm of $A_{a,\leq}$ by the Hilbert-Schmidt (HS) norm yields
	\begin{align}
		\|A_{a,\leq}\|^2&\leq \|A_{a,\leq}\|_{\HS}^2=\sum_{i,j}S_{iy}(G_uS^{\pa{a}}G_u^*)_{ij}S_{jy}(G_uS^{\pa{a}}G_u^*)_{ji} =\avga{G_uS^{\pa{a}}G_u^*S^{(y)}G_uS^{\pa{a}}G_u^*S^{(y)}} .\label{4_loop_term_Ward_inequality}
	\end{align}
	Applying Cauchy-Schwarz together with \eqref{double_difference_input}, we obtain
	\begin{equation}\label{4_loop_term_Ward_inequalityCS}
		\|A_{a,\leq}\|^2\leq \pbb{\sum_{i,k}S_{yi}|(G_u)_{ik}|^2S_{k a}}^2=\pa{\cL_{u,ya}^{(-,+)}}^2\prec \Psi_u^4\pa{|a-y|}\lesssim \Psi_u^4\pa{|x-y|},
	\end{equation}
	where in the last step, we used that $|a-y|\geq |x-y|-|x-a|\geq |x-y|/2$ whenever $|x-a|\leq |x-y|/2$, together with the polynomial decay property \eqref{polydecay_Psiu}. Substituting \eqref{4_loop_term_Ward_inequalityCS} into \eqref{bound_I_leq} and invoking Ward's identity \eqref{eq_Ward0}, we deduce
	\begin{equation}\label{3_loop_term_Ward_inequality}
		\begin{aligned}
			I_{\leq}&\prec \Psi_u^2\pa{|x-y|} \sum_{a}\|\psi_{a,\leq}\|^2= \Psi_u^2\pa{|x-y|} \avga{G_uS^{\pa{x}}G_u^*S^{(y)}G_uS^{\pa{x}}G_u^*}\\
			&    =\frac{\Psi_u^2\pa{|x-y|}}{2\ii\eta_u}\pa{\avga{G_uS^{\pa{x}}G_u^*S^{\pa{y}}G_uS^{\pa{x}}}-\avga{G_u^*S^{\pa{x}}G_u^*S^{\pa{y}}G_uS^{\pa{x}}}}.
		\end{aligned}
	\end{equation}
	Applying Cauchy-Schwarz once more to each term, and using the local law \eqref{eq:localPsiu} together with the argument leading to \eqref{4_loop_term_Ward_inequalityCS}, we conclude
	\begin{equation}\label{eq:Ileq}
		I_\leq\prec \frac{\Psi_u^2\pa{|x-y|}}{\eta_u}\cdot \Psi_u(0)\cL_{u,xy}^{(+,-)} \prec \frac{ \Psi_u(0)}{1-u}\cdot\Psi_u^4\pa{|x-y|}.
	\end{equation}
	
	For $I_{>}$, we employ a slightly different quadratic form:
	\begin{equation}\label{bound_I_g}            I_{>}=\sum_{a'\in\ZL}\psi_{a',>}A_{a',>}\psi_{a',>}^*\leq \sum_{a'\in\ZL} \|A_{a',>}\|\cdot\|\psi_{a',>}\|^2,
	\end{equation}
	where the vector \smash{$\psi_{a',>}\in \C^N$ and the matrix $A_{a',>}\in \C^{N\times N}$} are defined as 
	\begin{equation*}
		\begin{aligned}
			\psi_{a',>}\pa{i}:=(G_u^*S^{\pa{y}}G_u)_{a' i}\sqrt{S_{i x}},\quad A_{a',>}\pa{i,j}:=\sum_{a:|x-a|>|x-y|/2}\sqrt{S_{xi}}(G_u^*)_{ia}S_{aa'}(G_u)_{aj}\sqrt{S_{jx}}.
		\end{aligned}
	\end{equation*}
	Bounding again by the Hilbert-Schmidt norm gives
	\begin{equation}\label{eq:A>HS}
		\|A_{a',>}\|^2\leq \|A_{a',>}\|_{\HS}^2=\sum_{a_1:|x-a_1|>|x-y|/2}\sum_{a_2:|x-a_2|>|x-y|/2}(G_u S^{(x)}G_u^*)_{a_2a_1}S_{a_1a'}(G_uS^{(x)}G_u^*)_{a_1a_2}S_{a_2a'}.
	\end{equation}    
	Using the local law \eqref{eq:localPsiu} together with
	\eqref{polydecay_Psiu}--\eqref{polydecay_Psiu3}, we estimate
	\begin{align*}
		\absa{(G_uS^{(x)}G_u^*)_{a_1a_2}} &\le \sum_i S_{xi}\absa{(G_u)_{a_1 i}(G_u^*)_{i a_2}} \le \pB{\sum_i S_{xi}|(G_u)_{a_1 i}|^2}^{1/2} \pB{\sum_i S_{xi}|(G_u)_{a_2i }|^2}^{1/2} \\
		&\prec  \sqrt{[S_{xa_1}+\Psi^2_u(|x-a_1|)][S_{xa_2}+\Psi_u^2(|x-a_2|)]} \prec  \Psi_u^2(|x-y|)
	\end{align*}
	uniformly over the summation region. Inserting this bound into
	\eqref{eq:A>HS} yields \smash{$\|A_{a',>}\|^2\prec \Psi_u^4\pa{|x-y|}$}. Applying this estimate in \eqref{bound_I_g} and arguing as in \eqref{3_loop_term_Ward_inequality}, we conclude
	\[I_>\prec \frac{\Psi_u(0)}{1-u}\cdot\Psi_u^4\pa{|x-y|} .\]
	Together with \eqref{eq:Ileq}, this completes the proof of \Cref{lemma_martingale_term}.

	\begin{remark}\label{rmk:MG_for_T_fails}
		It is tempting to extend the idea of the loop-contraction inequality to the martingale term in the evolution equation \eqref{T_evolution} for the diagonal $T$-variables. If this were possible, one could restrict the dynamical analysis to the $T$-variables and thereby avoid some of the technical difficulties arising from the analysis of loops, such as those encountered in \Cref{lemma_double_difference_decomposition}. However, this extension fails due to the lack of averaging in the $T$-variables compared with the $\cL$-loops.
		More precisely, when bounding $I_>$ in the proof above, we estimate \smash{$\sum_{a'\in\ZL}\|\psi_{a',>}\|^2$} by certain three-resolvent loops (analogous to \eqref{3_loop_term_Ward_inequality}), where the averaging is introduced through $S^{(x)}$, $S^{(y)}$, and $S^{(x)}$. These terms can then be further bounded by $\Psi_u(0)\Psi^2_u(|x-y|)$ via a CS argument.
		In contrast, for the quadratic variation tensor of the martingale terms in the evolution equation of $T_{t,xy}$, the analogous step leads to a bound of the form $(G_uS^{(x)}G_u^*)_{yy}\im (G_u)_{yy}$ when estimating the counterpart of \smash{$\sum_{a'\in\ZL}\|\psi_{a',>}\|^2$}. Since the diagonal resolvent entry is of constant order and does not yield any small factor, this expression can only be further bounded by $\Psi^2_u(|x-y|)$. The resulting bound therefore lacks the additional $\Psi_u(0)$ factor, which invalidates the arguments in Steps 2 and 3 of the proof of Lemma \ref{lemma_bootstrap}.
	\end{remark}

	\section{Graphical tools: Proof of \texorpdfstring{\Cref{claim_graph_bound_light_weight_u}}{Claim 2}}\label{sec_graphical_argument}

	This section is devoted to the proof of \Cref{claim_graph_bound_light_weight_u} by developing several new graphical tools. Recall that $H$ is a power-law random band matrix defined in Definition \ref{considered_model}, with $\alpha\in (0,1)$, and satisfying Assumption \ref{assumption_input_bound}. For notational convenience, we work under assumptions slightly more general than those in Lemma~\ref{claim_graph_bound_light_weight_u}.
 More precisely, fix $t\in[0,1-N^{-1}]$ and $z=E+\ii \eta$, and denote $G\equiv G_t(z)=(\sqrt{t}H-z)^{-1}$ and $m\equiv m_t(z)$ (recall \eqref{eq:mtz}). 
	Let $\Psi:[0,\infty)\to (0,\infty)$ be an admissible control parameter in the sense of \Cref{assm_admissible}. By monotonicity of $\Psi$, we have $\Psi(0)=\max_{x,y}\Psi\p{|x-y|}$. Suppose that $\|G_t(z)-m_t(z)\|_{\max}\prec \Psi(0) \le W^{-\e}$, $|m_t(z)|=\OO(1)$, and  \begin{equation}\label{double_difference_input_2}
		\cL_{xy}\equiv\cL_{xy}\pa{z}:=\sum_{a,b\in\ZL}S_{xa}|G_{ab}\p{z}|^2S_{by}\prec \Psi^2\p{|x-y|}.
	\end{equation}
	In addition, assume that the following deterministic bounds hold:
	\begin{equation}\label{graph:Spm}
		S_{xy}^+:=\pa{\frac{tS}{1-tm^2S}}_{xy}\prec \cS_{xy},\quad S_{xy}^{-}:=\pa{\frac{tS}{1-t\ol{m}^2S}}_{xy}\prec \cS_{xy}.
	\end{equation}	
	Note that Lemma~\ref{claim_graph_bound_light_weight_u} is recovered by taking $t=u$, $z=z_u$, and $\cL=\cL_u$. 
	We define
	\begin{equation}\label{fxyG}
		f_{xy}\pa{G}:=\sum_{a,b}S_{ab}(G-m)_{bb}G_{xa}G_{a y}.
	\end{equation}
	By Markov's inequality, in order to prove \eqref{graph_bound_light_weight_u}, it suffices to show that for any fixed $p\in 2\N$, 
	\begin{equation}\label{light_weight_term_moment_bound}
		\E\absa{f_{xy}\p{G}}^{p}\prec \eta^{-p}\qa{\Psi\p{0}}^p \qa{\Psi\p{|x-y|}}^p.
	\end{equation}

	We first derive a simple bound on $f_{xy}(G)$ and explain the main difficulty in improving it to \eqref{light_weight_term_moment_bound}.
	Under the assumption \eqref{double_difference_input_2}, the estimates \eqref{upper_bound_2_to_1} and \eqref{entrywise_bound_2_to_1} give
	\begin{equation} \label{eq:Gxixy}       \absa{G_{xy}}^2\prec\delta_{xy}+\xi_{xy}^2,\quad \text{with}\quad \xi_{xy}^2:=\cS_{xy}+\cL_{xy}(z)\prec \Psi^2(|x-y|).
	\end{equation}
	Together with the averaged local law \eqref{average_bound_2_to_1}, this implies
	\begin{equation*}
		\absa{f_{xy}(G)}\prec\Psi^2(0)\pB{\delta_{xy}+\Psi(|x-y|)+\sum_{a\in\ZL}\Psi(|x-a|)\Psi(|a-y|)}.
	\end{equation*}
	However, when the control parameter $\Psi(\cdot)$ does not decay sufficiently fast, the last term on the RHS can be highly singular. Indeed, even if we assume the sharp control parameter $\Psi^2(|x-y|)=B(\eta,|x-y|)$ defined in \eqref{def_B_t}, this term yields only a factor of $\eta^{-1}$ and therefore provides no improvement over the following estimate based on Ward's identity:
	\begin{equation*}
		\absa{f_{xy}(G)}\prec \Psi^2(0)\sum_{a\in\ZL}|G_{xa}G_{ay}| \le \eta^{-1}\Psi^2(0)\sqrt{\im G_{xx}\cdot \im G_{yy}} \prec \eta^{-1}\Psi^2(0).
	\end{equation*}
	Compared with \eqref{light_weight_term_moment_bound}, this estimate lacks the decay factor $\Psi(|x-y|)$ and is therefore too weak when $|x-y|\gg W$. 
	The key issue is that there are not enough resolvent entries within $f_{xy}(G)$. More precisely, to obtain both the decay factor $\Psi(|x-y|)$ and the small factor $\Psi(0)$, one needs two $G$-entries (including light-weights), at least one of which must have length $\gtrsim |x-y|$\footnote{By the length of a $G_{xy}$ entry, we mean $|x-y|$, which can be interpreted as the length of the solid edge representing $G_{xy}$; see \Cref{def_graph1}.} to produce the required spatial decay factor $\Psi(|x-y|)$. In addition, two more $G$-entries are needed for an application of Ward's identity, which yields the factor $\eta^{-1}$. By contrast, within $f_{xy}(G)$, there are only three $G$-entries, including the light-weight.

	To overcome this difficulty, a key observation in \cite{dubova2025delocalizationnonmeanfieldrandommatrices} is that, by representing resolvent expressions as graphs and performing a sequence of systematic diagrammatic expansions of the graph $|f_{xy}(G)|^p$, one can generate $p$ additional $G$-entries of length $\gtrsim|x-y|$---namely, one extra $G$-entry for each factor $f_{xy}(G)$ or $\overline f_{xy}(G)$. We adopt this idea by introducing appropriate graphical notation and bounding $\E\absa{f_{xy}\p{G}}^{p}$ via a carefully designed diagrammatic expansion strategy. However, in contrast to \cite{dubova2025delocalizationnonmeanfieldrandommatrices}, the slow power-law decay of the variance profile prevents the strategy therein from yielding the bound \eqref{light_weight_term_moment_bound} in our setting, and new ideas and graphical tools are therefore required.
	In \Cref{sec:graph_ideas}, we provide a heuristic discussion of the obstacles to adopting the strategy of \cite{dubova2025delocalizationnonmeanfieldrandommatrices}, explain the key ideas used to overcome these obstacles, and outline our new expansion strategy. We first introduce the graphical notation used in the proof.

	\subsection{Graphical notation}
	
	In this subsection, we introduce the basic graphical notation used throughout our arguments. We begin with the definition of $G$-graphs.
	
	\begin{definition}[$G$-graphs] \label{def_graph1} 
		Given a graph with vertices and edges, we assign the following structures and refer to the resulting object as a $G$-graph.  
		\begin{itemize}
			\item {\bf Vertices:} Vertices represent matrix indices appearing in our expressions. Each vertex is designated as either an external or an internal vertex: external vertices correspond to indices whose values are fixed, whereas internal vertices correspond to summation indices that will be summed over.
			
			\item{\bf  Solid edges:} We use $(x,y)$ to denote a solid edge from vertex $x$ to vertex $y$. Each solid edge represents a resolvent entry. More precisely, a blue (resp.~red) oriented solid edge from $x$ to $y$ represents a factor $G_{xy}$ (resp.~$\overline G_{xy}$).

			\item {\bf Light-weights:} A factor $(G-m)_{xx}$ (resp.~$\overline {(G-m)}_{xx}$) is represented by a blue (resp.~red) self-loop at $x$. Following the convention in \cite{yang2021delocalization}, we refer to these factors as blue and red light-weights, respectively.
            Moreover, if the only edges incident to a vertex $x$ are a light-weight and an $S$-waved edge (defined below), then we call this light-weight a \emph{leaf light-weight}. 
            To distinguish light-weights from the self-loops representing factors $G_{xx}$ or $\overline {G}_{xx}$, which are referred to as regular weights in \cite{dubova2025delocalizationnonmeanfieldrandommatrices}, we represent the latter by self-loops labeled with ``R''.

			\item {\bf Waved edges:}
			\begin{itemize}
				\item A black waved edge (also referred to as $S$-waved edge) between $x$ and $y$ represents a factor $S_{xy}$. 
				
				\item A blue (resp.~red) waved edge between $x$ and $y$ represents a factor $S^+_{xy}$ (resp.~$S^-_{xy}$).
			\end{itemize}
			
			\item {\bf Dotted edges:} A dotted edge between $x$ and $y$ represents the factor $\mathbf 1_{x=y}$, while a $\times$-dotted edge represents the factor $\mathbf 1_{x\ne y}$.

			\item{\bf Coefficient:} Each graph is associated with a coefficient of order $\OO(1)$, which is a polynomial in $m$, $\overline m$, and $t$.
		\end{itemize}

        For convenience, we collectively refer to solid edges and light-weights as {\bf $G$-edges}. To each $G$-graph $\mathcal G$, we assign a \emph{value} as follows. First, we take the product of all edge factors together with the coefficient associated with the graph. We then sum over all internal indices corresponding to internal vertices, while keeping the external indices fixed at their prescribed values.
		For a linear combination of graphs $\sum_i c_i \cal G_i$, we define its value in the natural way as the linear combination of the values of the individual graphs $\cal G_i$. For simplicity, throughout the following proof, we will abuse notation by identifying a graph (a geometric object) with its value (the corresponding analytic expression).
	\end{definition}

Following \cite{yang2021delocalization}, given a $G$-graph, we define the corresponding \emph{molecular graph} by taking the quotient graph with respect to the equivalence classes induced by connectivity through waved edges.
	
	\begin{definition}[Molecules and molecular graphs]\label{mole_graph}
		For any $G$-graph $\cG$, we say that two vertices belong to the same molecule if they are connected by a path consisting entirely of waved edges. A molecule is called an \emph{external molecule} if it contains at least one external vertex; otherwise, it is called an \emph{internal molecule}. The \emph{molecular graph} of $\cG$, denoted by $\cG_{\cM}$, is obtained from $\cG$ through the following two operations:
		\begin{itemize}
			\item merging all vertices within the same molecule into a single vertex representing that molecule;
			
			\item discarding all other components except the solid edges and dotted edges between different molecules (in particular, edges within a single molecule and all weights are discarded).
			
		\end{itemize}
	\end{definition}
	
	We remark that our definition of molecules and molecular graphs differs slightly from those in \cite{yang2021delocalization,dubova2025delocalizationnonmeanfieldrandommatrices}, where two vertices belong to the same molecule whenever they are connected by a path each of whose edges is \emph{dotted} or waved. To facilitate the estimation of $G$-graphs, we introduce their scaling orders as follows. The scaling order will be used to control the $\max$-norm bound of the corresponding $G$-graphs.
	
	\begin{definition}[Scaling order]\label{def_scaling_orer}
		For a $G$-graph $\cG$, we define its \emph{scaling order}, denoted by $\mathsf{ord}\pa{\cG}$, as
		\begin{equation*}
			\begin{aligned}
				\mathsf{ord}\pa{\cG}:=\#\ha{\txt{solid edges in }\cG}+\#\ha{\txt{light-weights in }\cG}-2\#\ha{\txt{internal molecules of }\cG}.
			\end{aligned}
		\end{equation*}
	\end{definition}

Heuristically, the scaling order predicts a max-norm bound for the graph through the following two steps: (1) summing over each internal molecule by Ward's identity consumes two solid edges; (2) the estimates \eqref{upper_bound_2_to_1} and \eqref{entrywise_bound_2_to_1} bound each of the remaining $\mathsf{ord}\pa{\cG}$ $G$-edges by $\Psi(0)$, resulting in the factor $\Psi(0)^{\mathsf{ord}\pa{\cG}}$.   
    In our proof, we will repeatedly apply the following $GG$-expansion formulas to $G$-graphs. These formulas can be derived using Gaussian integration by parts; see \cite[Lemma 3.14]{yang2021delocalization}.
	
	\begin{lemma}[$GG$-expansions]\label{lemma:GGexpansion}
		Let $f(G)$ be a differentiable function of $G$. Consider the graph $\cG=G_{xy}G_{y'x}f\pa{G}$, which contains the \emph{$GG$-pair $G_{xy}G_{y'x}$ centered at $x$}. Then
		\begin{align}
			\cG& =_{\E} m\delta_{xy}G_{y' x}f(G)+  m^3  S^{+}_{xy} G_{y' y} f(G)  + t m \sum_a  S_{xa} \Gc_{a a} \cal G + t m^3 \sum_{a,b}  S^{+}_{xa}  S_{ab} \Gc_{b b} G_{a y}   G_{y'a} f(G) \nonumber\\
			&  + t m\Gc_{xx }    \sum_a S_{xa}G_{a y}   G_{y'a} f(G) + t m^3 \sum_{a,b}  S^+_{xa} S_{ab}  \Gc_{aa }  G_{b y}   G_{y'b} f(G)\label{out_GG_expansion}\\
			& - t m \sum_{ a}  S_{xa}G_{a y} G_{y' x} \partial_{ax}f(G) - t m^3 \sum_{a,b} S^{+}_{xa} S_{ab}G_{b y} G_{y' a} \partial_{ba}f(G),\nonumber
		\end{align}
		where ``$=_{\E}$'' means ``equal in expectation'', and we write $\Gc:=(G-m)$ and $\partial_{ax}\equiv \partial_{(H_t)_{ax}}$ with $H_t\equiv \sqrt{t}H$.  
		We refer to \eqref{out_GG_expansion} as the \emph{out-edge} $GG$-expansion of the $GG$-pair $G_{xy}G_{y'x}$ with respect to the edge $G_{xy}$. 
		Similarly, we have the following \emph{in-edge} $GG$-expansion of the same $GG$-pair with respect to the edge $G_{y'x}$: 
		\begin{equation}\label{in_GG_expansion}
			\begin{aligned}
				\cG& =_{\E} mG_{xy}\delta_{y'x}f(G)+  m^3  S^{+}_{xy'} G_{y' y} f(G)  + t m \sum_a  S_{xa} \Gc_{aa} \cal G + t m^3 \sum_{a,b}  S^{+}_{xa}  S_{ab} \Gc_{b b} G_{a y}   G_{y'a} f(G) \\
				&  + t m\Gc_{xx }    \sum_a S_{xa}G_{a y}   G_{y'a} f(G) + t m^3 \sum_{a,b}  S^+_{xa} S_{ab}  \Gc_{aa }  G_{b y}   G_{y'b} f(G)\\
				& - t m \sum_{ a}  S_{xa}G_{x y} G_{y' a} \partial_{ x a}f(G) - t m^3 \sum_{a,b} S^{+}_{xa} S_{ab}G_{a y} G_{y' b} \partial_{ a b}f(G).
			\end{aligned}
		\end{equation}
		Analogous $\ol{GG}$-expansions hold for $\ol{GG}$-pairs $\overline G_{xy}\overline G_{y'x}$. 
	\end{lemma}

	When no confusion arises, we will refer to both $GG$-pairs and $\ol{GG}$-pairs simply as $GG$-pairs. Furthermore, for clarity of presentation, we classify the terms on the RHS of \eqref{out_GG_expansion} and \eqref{in_GG_expansion} as follows. The first two terms will be referred to as the \emph{reshaped terms}, since the $GG$-pair is reshaped into a structure containing either a dotted or a waved edge, which we refer to as a dotted structure or a waved structure, respectively.
	The third through sixth terms will be called the \emph{growing terms}, as a light-weight is grown in the graph. The final two terms will be referred to as the \emph{$\partial$-terms}. 
	In addition, the subgraphs $S_{xa}G_{ay}G_{y' x}$, with inner vertices $x$ and $a$, and $S_{xa}^+S_{ab}G_{by}G_{y'a}$, with inner vertices $a$ and $b$, arising from the $\partial$-terms in \eqref{out_GG_expansion}, will be referred to as \emph{extended $GG$-pairs} with the corresponding inner vertices (which are precisely the vertices denoted by $\star_1$ and $\star_2$ in the corresponding derivative $\partial_{\star_1\star_2}$). The extended $GG$-pairs arising from the $\partial$-terms in \eqref{in_GG_expansion} are defined analogously.
    For reference, the complete out-edge $GG$-expansion is represented schematically in \Cref{tikz_out_edge_expansion}.

    \begin{figure}[htbp]
  \centering
  \begin{tikzpicture}[
    x=0.70cm,
    y=0.70cm,
    fermion/.style={
      thick,
      postaction={decorate},
      decoration={markings,
        mark=at position 0.56 with {\arrow[scale=1.1]{Latex}}}
    },
    photon/.style={
      thick,
      decorate,
      decoration={snake,segment length=4.5pt,amplitude=1.1pt}
    },
    dot/.style={circle,fill=black,inner sep=0.9pt},
    fbox/.style={
      draw=gray!65,
      densely dashed,
      rounded corners=1.5pt,
      fill=gray!6,
      inner xsep=4pt,
      inner ysep=2.5pt,
      font=\scriptsize
    },
    potential/.style={
      gray!70,
      densely dashed,
      line width=0.55pt,
      shorten >=1.0pt
    },
    index/.style={font=\scriptsize,inner sep=1pt},
    panel/.style={font=\scriptsize},
    ident/.style={thick,densely dotted}
  ]

\newcommand{\panelbox}{\path[use as bounding box] (-3.15,-1.75) rectangle (3.15,2.05);}

\begin{scope}[shift={(0,0)}]
    \panelbox
    \node[dot,label={[index]above left:$x$}]  (x)  at (-1.15,0.35) {};
    \node[dot,label={[index]above right:$y$}] (y)  at (0.85,0.35) {};
    \node[dot,label={[index]left:$y'$}]        (yp) at (-0.60,-0.85) {};
    \node[fbox] (F) at (2.05,1.20) {$f(G)$};
    \draw[fermion,blue] (x) -- (y);
    \draw[fermion,blue] (yp) -- (x);
    \node[panel] at (0,-1.48) {(a) $\mathcal G$};
  \end{scope}

\begin{scope}[shift={(7.15,0)}]
    \panelbox
    \node[dot,label={[index]above left:$x$}]  (x)  at (-1.15,0.35) {};
    \node[dot,label={[index]above right:$y$}] (y)  at (0.85,0.35) {};
    \node[dot,label={[index]left:$y'$}]        (yp) at (-0.60,-0.85) {};
    \node[fbox] (F) at (2.05,1.20) {$f(G)$};
    \draw[ident] (x) -- (y);
    \draw[fermion,blue] (yp) -- (x);
    \node[panel] at (0,-1.48) {(b) term 1 $(m)$};
  \end{scope}

\begin{scope}[shift={(14.30,0)}]
    \panelbox
    \node[dot,label={[index]above left:$x$}]  (x)  at (-1.15,0.35) {};
    \node[dot,label={[index]above right:$y$}] (y)  at (0.85,0.35) {};
    \node[dot,label={[index]left:$y'$}]        (yp) at (-0.60,-0.85) {};
    \node[fbox] (F) at (2.05,1.20) {$f(G)$};
    \draw[photon,blue] (x) -- (y);
    \draw[fermion,blue] (yp) -- (y);
    \node[panel] at (0,-1.48) {(c) term 2 $(m^3)$};
  \end{scope}

\begin{scope}[shift={(0,-4.55)}]
    \panelbox
    \node[dot,label={[index]left:$x$}]        (x)  at (-1.10,0.05) {};
    \node[dot,label={[index]right:$y$}]       (y)  at (0.85,0.05) {};
    \node[dot,label={[index]left:$y'$}]       (yp) at (-1.60,-0.80) {};
    \node[dot,label={[index]below right:$a$}] (a)  at (-0.65,1.15) {};
    \node[fbox] (F) at (2.05,1.20) {$f(G)$};
    \draw[fermion,blue] (x) -- (y);
    \draw[fermion,blue] (yp) -- (x);
    \draw[photon] (x) -- (a);
    \draw[thick,blue] (a) arc (-90:270:0.40);
    \node[panel] at (0,-1.48) {(d) term 3 $(tm)$};
  \end{scope}

\begin{scope}[shift={(7.15,-4.55)}]
    \panelbox
    \node[dot,label={[index]above left:$x$}]  (x)  at (-1.75,0.75) {};
    \node[dot,label={[index]below left:$a$}]  (a)  at (-0.55,0.25) {};
    \node[dot,label={[index]right:$b$}]       (b)  at (0.25,1.05) {};
    \node[dot,label={[index]right:$y$}]       (y)  at (1.20,0.10) {};
    \node[dot,label={[index]left:$y'$}]       (yp) at (-0.35,-0.90) {};
    \node[fbox] (F) at (2.10,1.25) {$f(G)$};
    \draw[photon,blue] (x) -- (a);
    \draw[photon] (a) -- (b);
    \draw[fermion,blue] (a) -- (y);
    \draw[fermion,blue] (yp) -- (a);
    \draw[thick,blue] (b) arc (-90:270:0.40);
    \node[panel] at (0,-1.48) {(e) term 4 $(tm^3)$};
  \end{scope}

\begin{scope}[shift={(14.30,-4.55)}]
    \panelbox
    \node[dot,label={[index]left:$x$}]        (x)  at (-1.45,0.75) {};
    \node[dot,label={[index]above:$a$}]       (a)  at (-0.20,0.00) {};
    \node[dot,label={[index]right:$y$}]       (y)  at (1.20,0.00) {};
    \node[dot,label={[index]left:$y'$}]       (yp) at (-0.45,-0.90) {};
    \node[fbox] (F) at (2.10,1.25) {$f(G)$};
    \draw[photon] (x) -- (a);
    \draw[fermion,blue] (a) -- (y);
    \draw[fermion,blue] (yp) -- (a);
    \draw[thick,blue] (x) arc (-90:270:0.40);
    \node[panel] at (0,-1.48) {(f) term 5 $(tm)$};
  \end{scope}

\begin{scope}[shift={(0,-9.10)}]
    \panelbox
    \node[dot,label={[index]above left:$x$}]  (x)  at (-1.75,0.70) {};
    \node[dot,label={[index]below left:$a$}]  (a)  at (-0.55,0.45) {};
    \node[dot,label={[index]above:$b$}]       (b)  at (0.25,-0.05) {};
    \node[dot,label={[index]right:$y$}]       (y)  at (1.30,-0.05) {};
    \node[dot,label={[index]left:$y'$}]       (yp) at (-0.20,-0.90) {};
    \node[fbox] (F) at (2.10,1.25) {$f(G)$};
    \draw[photon,blue] (x) -- (a);
    \draw[photon] (a) -- (b);
    \draw[fermion,blue] (b) -- (y);
    \draw[fermion,blue] (yp) -- (b);
    \draw[thick,blue] (a) arc (-90:270:0.40);
    \node[panel] at (0,-1.48) {(g) term 6 $(tm^3)$};
  \end{scope}

\begin{scope}[shift={(7.15,-9.10)}]
    \panelbox
    \node[dot,label={[index]above left:$x$}]  (x)  at (-1.35,0.40) {};
    \node[dot,label={[index]above:$a$}]       (a)  at (-0.10,0.40) {};
    \node[dot,label={[index]right:$y$}]       (y)  at (1.20,0.40) {};
    \node[dot,label={[index]left:$y'$}]       (yp) at (-1.40,-0.80) {};
    \node[fbox] (F) at (2.10,1.25) {$f(G)$};
    \node[index,text=gray!75] at (2.10,1.82) {$\partial_{ax}$};
    \draw[photon] (x) -- (a);
    \draw[fermion,blue] (a) -- (y);
    \draw[fermion,blue] (yp) -- (x);
    \draw[potential] (F.west) .. controls (1.15,1.25) and (0.35,0.85) .. (a);
    \draw[potential] (F.south west) .. controls (0.90,0.75) and (-0.45,1.05) .. (x);
    \node[panel] at (0,-1.48) {(h) term 7 $(-{t}m)$};
  \end{scope}

\begin{scope}[shift={(14.30,-9.10)}]
    \panelbox
    \node[dot,label={[index]above left:$x$}]  (x)  at (-1.75,0.70) {};
    \node[dot,label={[index]above:$a$}]       (a)  at (-0.55,0.35) {};
    \node[dot,label={[index]above:$b$}]       (b)  at (0.25,-0.10) {};
    \node[dot,label={[index]right:$y$}]       (y)  at (1.35,-0.10) {};
    \node[dot,label={[index]left:$y'$}]       (yp) at (-0.45,-0.85) {};
    \node[fbox] (F) at (2.10,1.25) {$f(G)$};
    \node[index,text=gray!75] at (2.10,1.82) {$\partial_{ba}$};
    \draw[photon,blue] (x) -- (a);
    \draw[photon] (a) -- (b);
    \draw[fermion,blue] (b) -- (y);
    \draw[fermion,blue] (yp) -- (a);
    \draw[potential] (F.west) .. controls (1.15,1.25) and (0.25,0.90) .. (a);
    \draw[potential] (F.south west) .. controls (1.30,0.65) and (0.75,0.25) .. (b);
    \node[panel] at (0,-1.48) {(i) term 8 $(-{t}m^3)$};
  \end{scope}
  \end{tikzpicture}
  \caption{Schematic representation of the out-edge $GG$-expansion in
  \eqref{out_GG_expansion}. Panel~(a) represents the input graph
  $\mathcal G=G_{xy}G_{y'x}f(G)$, while panels~(b)--(i) represent, in
  order, the eight terms on the right-hand side. The scalar coefficient
  of each graph is displayed in parentheses, and the vertices labeled
  $a$ and $b$ are internal vertices and are therefore summed over. The
  gray dashed box is a schematic placeholder for the unspecified factor
  $f(G)$; it is neither a vertex nor an edge and may share indices with
  the displayed subgraph. In panels~(h) and~(i), the two gray dashed
  connectors join this box to the two indices of the displayed derivative:
  $a,x$ in panel~(h) and $b,a$ in panel~(i). They schematically represent
  the two potential $G$-edges pulled out of $f(G)$ when the derivative acts
  on a $G$-edge inside $f(G)$. Since the differentiated edge is unspecified,
  these connectors carry neither color nor orientation and are not graph
  edges before the derivative is expanded.}\label{tikz_out_edge_expansion}
\end{figure}
	
\begin{figure}[htbp]
	\centering
	\vspace{-0.2cm}
    \hspace{-0.5cm}
	\begin{minipage}[c]{0.4\textwidth}
		\caption{Graphs induced by the seventh term on the RHS of \eqref{out_GG_expansion}:}
		\label{tikz_glued_two_tadpoles}
	\end{minipage}
    \hfill    
	\begin{minipage}[c]{0.6\textwidth}
		\centering
		\begin{tikzpicture}[
			scale=0.4,
			fermion/.style={postaction={decorate}, decoration={markings,
				mark=at position 0.55 with {\arrow[scale=1.2]{Latex}}}},
			photon/.style={decorate, decoration={snake, segment length=5pt,
				amplitude=1.2pt}},
			dot/.style={circle, fill=black, inner sep=1.0pt}
			]

\node[dot, label=left:{$x_1$}] (wx1) at (-3.0, 0.5) {};
				\node[dot, label=left:{$x_2$}] (wx2) at (-2.5, -2.0) {};
				\node[dot, label=right:{$y_1$}] (wy1) at (3.5, 0.5) {};
				\node[dot, label=right:{$y_2$}] (wy2) at (3.5, -2.0) {};
				\node[dot, label=left:{$a_1$}] (wL) at (-1.0, 2.5) {};
				\node[dot] (wT) at (-1.0, 4.5) {};
				\node[dot, label=right:{$a_1'$}] (wR1) at (1.5, 2.5) {};
				\node[dot, label=right:{$b_2$}] (wB) at (0.0, 0.0) {};
				\node[dot, label=below:{$a_2$}] (wR2) at (0.0, -1.5) {};

				\draw[thick, fermion, blue] (wx1) -- (wL);
				\draw[thick, photon] (wL) -- (wT);
				\draw[thick, blue] (wT) arc (-90:270:0.6);
				\draw[thick, photon] (wL) -- (wR1);
				\draw[thick, fermion, blue] (wR1) -- (wy1);
				\draw[thick, fermion, red] (wB) -- (wL);
				\draw[thick, fermion, red] (wR1) -- (wB);
				\draw[thick, fermion, red] (wx2) -- (wR2);
				\draw[thick, fermion, red] (wR2) -- (wy2);
				\draw[thick, photon] (wR2) -- (wB);
				\node at (0.25,-3.35) {\scriptsize (a) $\partial$ on the light-weight};

\begin{scope}[shift={(10.5,0)}]
			\node[dot, label=left:{$x_1$}] (x1) at (-3.0, 0.5) {};
			\node[dot, label=left:{$x_2$}] (x2) at (-2.5, -2.0) {};
			\node[dot, label=right:{$y_1$}] (y1) at (3.5, 0.5) {};
			\node[dot, label=right:{$y_2$}] (y2) at (3.5, -2.0) {};
			\node[dot, label=left:{$a_1$}] (vL) at (-1.0, 2.5) {};
			\node[dot] (vT) at (-1.0, 4.5) {};
			\node[dot, label=right:{$a_1'$}] (vR1) at (1.5, 2.5) {};
			\node[dot, label=right:{$a_2$}] (vR2) at (0.0, 0.0) {};
			\node[dot, label=right:{$b_2$}] (vB) at (0.0, -1.8) {};

			\draw[thick, fermion, blue] (x1) -- (vL);
			\draw[thick, fermion, red] (x2) -- (vL);
			\draw[thick, photon] (vL) -- (vT);
			\draw[thick, blue] (vT) arc (-90:270:0.6);
			\draw[thick, photon] (vL) -- (vR1);
			\draw[thick, fermion, red] (vR1) -- (vR2);
			\draw[thick, fermion, blue] (vR1) -- (y1);
			\draw[thick, fermion, red] (vR2) -- (y2);
			\draw[thick, photon] (vR2) -- (vB);
			\draw[thick, red] (vB) arc (90:450:0.6);
			\node at (0.25,-3.35) {\scriptsize (b) $\partial$ on the in-edge};
			\end{scope}
		\end{tikzpicture}
	\end{minipage}
	\vspace{-0.3cm}
\end{figure}
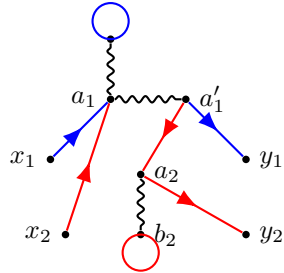

	For any fixed $p\in 2\N$, we estimate the graph $\big|f_{xy}(G)\big|^{p} = [f_{xy}(G)]^{p/2}[\overline{f_{xy}(G)}]^{p/2}$ by expanding it using the $GG$-expansions introduced above. First, $\absa{f_{xy}(G)}^p$ is represented by the graph in \eqref{graph_tadpole_g}:	
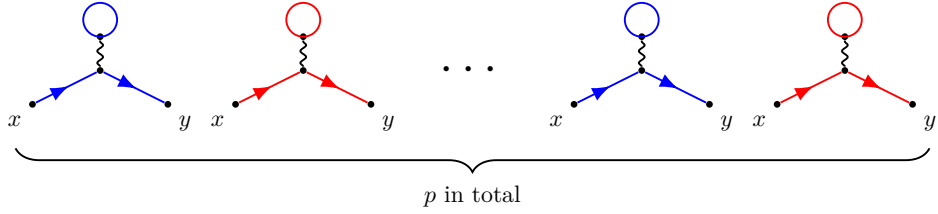
    \begin{equation}  \label{graph_tadpole_g} 
\parbox[c]{0.8\linewidth}{\begin{center}\scalebox{0.7}{
			\begin{tikzpicture}[
				scale=0.5,         fermion/.style={postaction={decorate}, decoration={markings, mark=at position 0.55 with {\arrow[scale=1.2]{Latex}}}},
				photon/.style={decorate, decoration={snake, segment length=5pt, amplitude=1.2pt}},
				dot/.style={circle, fill=black, inner sep=1.0pt}
				]
				
				\def\drawDiagram#1#2{
					\begin{scope}[shift={(#1, 0)}]
						\node[dot, label=below left:{$x$}] (x) at (-2.0, 0) {};
						\node[dot, label=below right:{$y$}] (y) at (2.0, 0) {};
						\node[dot] (v1) at (0, 1) {};
						\node[dot] (v2) at (0, 2) {};
						
						\draw[thick, fermion, #2] (x) -- (v1);
						\draw[thick, fermion, #2] (v1) -- (y);
						\draw[thick, photon] (v1) -- (v2);
						\draw[thick, #2] (v2) arc (-90:270:0.5);
					\end{scope}
				}

				\drawDiagram{0}{blue}
				
				\drawDiagram{6}{red}
				
				\node at (11, 1) {\Huge $\cdots$};
				
				\drawDiagram{16}{blue}
				
				\drawDiagram{22}{red}

				\draw[thick, decorate, decoration={brace, amplitude=10pt, mirror}] 
				(-2.5, -1.3) -- (24.5, -1.3) 
				node[midway, below=12pt] {$p \text{ in total}$};
				
		\end{tikzpicture}}
	\end{center}}
    \end{equation}
	Here, a blue tadpole diagram represents a factor $f_{xy}(G)$, whereas a red tadpole diagram represents a factor \smash{$\ol{f_{xy}(G)}$}. In our expansion strategy, we first select a blue tadpole diagram $f_{x_1y_1}(G)$\footnote{Here and in the following proof, we use $(x_i,y_i)$ to distinguish the endpoints of the $p$ tadpole diagrams, although we keep in mind that they all represent the same pair $(x,y)$.} from the graph \smash{$\absa{f_{xy}(G)}^p$}, and apply the $GG$-expansion in either \eqref{out_GG_expansion} or \eqref{in_GG_expansion} with respect to the $GG$-pair $ G_{x_1a_1}G_{a_1y_1}$. 
	For instance, suppose that the out-edge $GG$-expansion  \eqref{out_GG_expansion} is applied. This expansion generates $\OO(1)$ terms. The terms arising from the first six terms on the RHS of \eqref{out_GG_expansion} can be viewed as \emph{local} deformations of the chosen tadpole diagram; in particular, they do not ``pull in'' edges from the remaining $(p-1)$ tadpole diagrams.
	Among the terms generated by the last two terms on the RHS of \eqref{out_GG_expansion}, some still correspond to local deformations if the $\partial$-derivative acts on the selected tadpole diagram $f_{x_1y_1}(G)$. For these graphs, we say that the number of \emph{effective sections} remains one.\footnote{Note that the tadpole diagrams may transform during the expansion procedure. To track these changes, we introduce the notion of \emph{sections}, defined as the subgraphs induced by the corresponding molecules.} Two examples are illustrated in Figure~\ref{tikz_deformed_tadpoles}.
	On the other hand, if a new graph is generated from the last two terms on the RHS of \eqref{out_GG_expansion} and the $\partial$-operator acts on another tadpole diagram, say $f_{x_2y_2}(G)$ or \smash{$\ol{f_{x_2y_2}(G)}$}, then this diagram is ``pulled in'' and ``glued'' to the selected diagram. In this case, the number of effective sections increases to two. Two examples of this situation are shown in Figure~\ref{tikz_glued_two_tadpoles}.
    
	\begin{figure}[htbp]
		\vspace{-0cm}
		\centering
		\hspace{-2.0cm}
		\begin{minipage}[c]{0.4\textwidth}
			\vspace{-1.0cm}
			\caption{Graphs induced by the second and the seventh terms on the RHS of (\eqref{out_GG_expansion}): }
			\label{tikz_deformed_tadpoles}
		\end{minipage}
		\hspace{-1.0cm}
		\begin{minipage}[c]{0.3\textwidth}
			\centering
			\begin{tikzpicture}[scale=0.6, 
				fermion/.style={
					postaction={decorate},
					decoration={markings, mark=at position 0.55 with {\arrow[scale=1.2]{Latex}}} 
				},
				photon/.style={
					decorate,
					decoration={snake, segment length=5pt, amplitude=1.2pt}
				},
				dot/.style={
					circle, fill=black, inner sep=1.0pt
				}
				]
				
				\node[dot, label=below left:{$x$}] (x) at (-2.0, 0) {};

				\node[dot, label=below right:{$y$}] (y) at (2.0, 0) {};

				\node[dot] (v1) at (0, 1) {};

				\node[dot] (v2) at (0, 2) {};

				\draw[thick, photon, blue] (v1) -- (y);

				\draw[thick, photon] (v1) -- (v2);
				
				\draw[thick, fermion, blue] (x) -- (y);

				\draw[thick,blue] (v2) arc (-90:270:0.5);
			\end{tikzpicture}
		\end{minipage}\hspace{-0.5cm}
		\begin{minipage}[c]{0.3\textwidth}
			\centering
			\vspace{-0.0cm}
			\begin{tikzpicture}[scale=0.6, 
				fermion/.style={
					postaction={decorate},
					decoration={markings, mark=at position 0.55 with {\arrow[scale=1.2]{Latex}}} 
				},
				photon/.style={
					decorate,
					decoration={snake, segment length=5pt, amplitude=1.2pt}
				},
				dot/.style={
					circle, fill=black, inner sep=1.0pt
				}
				]
				
				\node[dot, label=below left:{$x$}] (x) at (-2.0, 0) {};

				\node[dot, label=below right:{$y$}] (y) at (2.0, 0) {};

				\node[dot] (v1) at (0, 1.25) {};

				\node[dot] (v2) at (0, 3.0) {};
				
				\node[dot] (v3) at (1,1) {};

				\draw[thick, fermion, blue] (x) -- (v1);

				\draw[thick, photon] (v1) -- (v3);
				
				\draw[thick, fermion, blue] (v2) -- (v3);

				\draw[thick, photon] (v1) -- (v2);
				
				\draw[thick, fermion, blue] (v1) to[out=130, in=-120] (v2);
				
				\draw[thick, fermion, blue] (v3) -- (y);

			\end{tikzpicture}
		\end{minipage}    \vspace{-0.5cm}
	\end{figure}
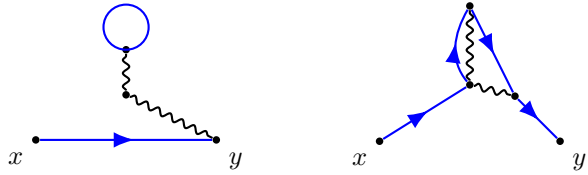

	For the graphs in Figures~\ref{tikz_deformed_tadpoles} and \ref{tikz_glued_two_tadpoles}, one can readily verify, using \eqref{eq:Gxixy} together with the averaged local law \eqref{average_bound_2_to_1}, that for $x\neq y$ the left graph in \Cref{tikz_deformed_tadpoles} is bounded by $\opr{\Psi^2(0)\Psi(|x-y|)}$, and the left graph in \Cref{tikz_glued_two_tadpoles} is bounded by $\opr{\eta^{-2}\Psi^2(0)\Psi^2(|x-y|)}$. However, the remaining two graphs do not admit strong enough estimates and therefore require further expansion. These expansions must be carried out in a carefully designed manner so that the resulting graphs can be bounded as required in \eqref{light_weight_term_moment_bound}. This procedure involves several nontrivial operations and new graphical concepts. Therefore, before presenting the general strategy, we provide a heuristic discussion of the key ideas and the expansion procedure.

	\subsection{New ideas and snake-like graphs}\label{sec:graph_ideas}
	
	The strategy in \cite[Section 7]{dubova2025delocalizationnonmeanfieldrandommatrices} roughly consists of three steps:
\begin{itemize}
	\item Expand the graph $|f_{xy}(G)|^p$ into a sum of $\OO(1)$ 
	\emph{locally standard graphs} by repeatedly applying suitable light-weight and edge expansions, together with the $GG$-expansions from \Cref{lemma:GGexpansion}.
	Here, locally standard graphs roughly refer to graphs without weights or $GG$-pairs,
	i.e., every $G$-edge is paired with a unique $\overline G$-edge.
    \item Bound the resulting locally standard graphs in terms of their
	corresponding molecular graphs.

	\item Control the molecular graphs by applying \eqref{eq:Gxixy} together with Ward's identities.
\end{itemize}
The first and the third steps (provided that the second step can be carried out) are largely model-independent and therefore apply equally well in our setting. The main purpose of the first step is to generate additional $G$-edges, which play two crucial roles: 
\begin{enumerate}
	\item {\bf Sufficient $\Psi$-factors:}
	They provide the necessary factors of both $\Psi(0)$ and $\Psi(|x-y|)$.

	\item {\bf Nested property:}
	They complete the graphical structure in such a way that, even after removing the edges used to extract the $\Psi$-factors, the remaining graph retains sufficient
	connectivity to allow repeated applications of Ward's identity when summing
	over each vertex.
\end{enumerate} 
The nested property then enables the implementation of the third step through
\eqref{eq:Gxixy} and Ward's identities.

On the other hand, the second step is designed to simplify the complicated structure of the $G$-graphs, so that the key nested property becomes apparent in the corresponding molecular graphs after ignoring the detailed structure within each molecule. For regular RBMs, this step is relatively routine, dating back at least to \cite{Band1D_III} and subsequently adopted in \cite{yang2021delocalization,YYYTexpansion2022CMP,dubova2025delocalizationnonmeanfieldrandommatrices}. However, this step breaks down completely in our setting due to the slow power-law decay of the waved edges, characterized by $B_0(|x-y|)\asymp\cal S_{xy}$.
Indeed, for regular RBMs, the reason why the internal structure of each molecule can be neglected is that the waved edges decay exponentially beyond the scale $W$. As a consequence, the diameter of any molecule is at most $\OO_\prec\p{W}$, up to an error of order $\OO\p{W^{-D}}$ for any large constant $D>0$. On the other hand, the control parameter $\Psi(\cdot)$ is essentially flat on the scale $W$, so that each molecule can effectively be treated as a single point when tracking the decay of the solid edges. This mechanism breaks down in the present power-law setting. Here, the decay of the edges $S_{xy}$ and $S^{\pm}_{xy}$ remains effective on the scale of the entire system, so that molecules are no longer ``local structures'' on the scale $W$, in contrast to the situation for regular RBMs.

One might attempt to bypass the second step and work directly with the $G$-graphs without introducing molecular graphs. However, in such an approach, the molecular graph structure becomes intertwined with the local structure inside molecules, which destroys the nested property and therefore invalidates the third step. Consequently, we still need a mechanism to shrink the internal structure of each molecule by summing over the vertices within the molecule in a suitable way. A key tool for this purpose is the following inequality based on the CS inequality:
\begin{equation}\label{eq:reduction_xi}
	\sum_{b\in\ZL}\xi_{xb}\xi_{yb} \cS_{ba}\leq \xi'_{xa}\xi'_{ya},\quad \text{with}\quad (\xi'_{xa})^2:=\sum_{b\in\ZL}\xi^2_{xb}\cS_{ba},
\end{equation}
where we recall from \eqref{eq:Gxixy} that $\xi_{xy}$ controls the solid edges between different molecules. This inequality effectively reduces the number of vertices within a molecule. Moreover, the new edge weight $\xi'$ inherits two key properties from $\xi$: the bound $ \sum_{a}(\xi'_{ax})^2\prec \eta^{-1}$, which follows from Ward's identity, and the bound  $\absa{\xi'_{xy}}\prec \Psi(|x-y|)$, which follows from \eqref{double_difference_input_2} and \eqref{polydecay_Psiu3}. 
However, this approach does not always work. Intuitively, each waved edge can ``carry'' at most two solid edges, which either provide decay factors or are used in an application of Ward's identity. During the reduction of graphs using \eqref{eq:reduction_xi}, the method fails to yield satisfactory bounds when an endpoint of a waved edge is connected to more than two solid edges. To illustrate this issue, consider the following graphs in Figure \ref{tikz_bad_graphs}: 
	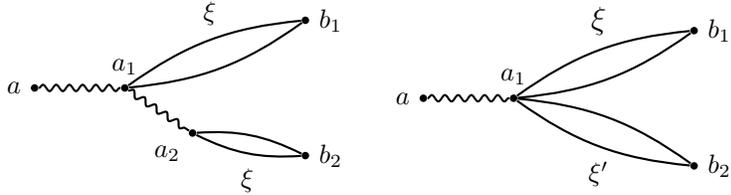
\begin{figure}[htbp]
		\vspace{-0cm}
		\centering
		\hspace{-1.0cm}
		\begin{minipage}[c]{0.4\textwidth}
			\vspace{-1.0cm}
			\caption{Examples of bad graphs: }
			\label{tikz_bad_graphs}
		\end{minipage}
		\hspace{-2.0cm}
		\begin{minipage}[c]{0.44\textwidth}
			\centering
			\begin{tikzpicture}[
				scale=0.6,         fermion/.style={postaction={decorate}, decoration={markings, mark=at position 0.55 with {\arrow[scale=1.2]{Latex}}}},
				photon/.style={decorate, decoration={snake, segment length=5pt, amplitude=1.2pt}},
				dot/.style={circle, fill=black, inner sep=1.0pt}
				]
				
				\node[dot, label=left:{ $a$}] (a) at (0, 0) {};
				
				\node[dot, label=above:{ $a_1$}] (a1) at (2.0, 0) {};
				\node[dot, label=right:{ $b_1$}] (b1) at (6.0, 1.5) {};
				
				\node[dot, label=below left:{ $a_2$}] (a2) at (3.5, -1.0) {};
				\node[dot, label=right:{ $b_2$}] (b2) at (6.0, -1.5) {};

				\draw[thick, photon] (a) -- (a1);               \draw[thick, photon] (a1) -- (a2);              
				\draw[thick] (a1) to[bend left=15] node[midway, above=2pt] { $\xi$} (b1);
				\draw[thick] (a1) to[bend right=15] (b1);
				
				\draw[thick] (a2) to[bend left=15] (b2);
				\draw[thick] (a2) to[bend right=15] node[midway, below=2pt] { $\xi$} (b2);
				
			\end{tikzpicture}
		\end{minipage}
		\hspace{-2.0cm}
		\begin{minipage}[c]{0.4\textwidth}
			\centering
			\begin{tikzpicture}[
				scale=0.6,         fermion/.style={postaction={decorate}, decoration={markings, mark=at position 0.55 with {\arrow[scale=1.2]{Latex}}}},
				photon/.style={decorate, decoration={snake, segment length=5pt, amplitude=1.2pt}},
				dot/.style={circle, fill=black, inner sep=1.0pt}
				]
				
				\node[dot, label=left:{ $a$}] (a) at (0, 0) {};
				
				\node[dot, label=above:{ $a_1$}] (a1) at (2.0, 0) {};
				\node[dot, label=right:{ $b_1$}] (b1) at (6.0, 1.5) {};
				
				\node[dot, label=right:{ $b_2$}] (b2) at (6.0, -1.5) {};

				\draw[thick, photon] (a) -- (a1);               
				\draw[thick] (a1) to[bend left=15] node[midway, above=2pt] {\Large $\xi$} (b1);
				\draw[thick] (a1) to[bend right=15] (b1);
				
				\draw[thick] (a1) to[bend left=15] (b2);
				\draw[thick] (a1) to[bend right=15] node[midway, below=2pt] { $\xi'$} (b2);
				
			\end{tikzpicture}
		\end{minipage}
		\hspace{-1.0cm}
		\vspace{-0.5cm}
	\end{figure}
    
	If we sum over the $\xi^2_{a_2 b_2}$ edges attached to $a_2$, we obtain two new solid edges \smash{$(\xi'_{a_1 b_2})^2$}. The vertex $a_1$ is then incident to four solid edges. In some situations, the distance between $a$ and $b_1$ may be of order $\gtrsim |x-y|$, so that the waved edge $(a,a_1)$ together with the solid edges $\xi^2_{a_1 b_1}$ must supply the decay factor $\Psi^2(|x-y|)$. Hence, we need to extract such a factor using \eqref{eq:reduction_xi}.
	However, since the waved edge $(a,a_1)$ can be long, the distance $|a_1-b_2|$ may be as small as $\OO(W)$. In this case, the factor \smash{$(\xi'_{a_1 b_2})^2$} can only be bounded by the non-decaying quantity $\Psi^2(0)$. On the other hand, it may also happen that the distance between $a$ and $b_2$ is of order $\gtrsim |x-y|$, in which case one should extract another factor $\Psi^2(|x-y|)$ from the two solid edges incident to $b_2$. This contribution is missed in the above argument. 
	In other situations, we may need to obtain an $\eta^{-1}$ factor by applying Ward's identity to the summation over the two solid edges $(\xi'_{a_1 b_2})^2$. However, the above argument yields only the bound $N\Psi^2(0)$, which is too weak for our purposes.
	
	Unfortunately, configurations such as those in \Cref{tikz_bad_graphs}, or even more severe ones in which a waved edge carries more than four solid edges, cannot be avoided within the expansion strategy of \cite{dubova2025delocalizationnonmeanfieldrandommatrices}. As a result, the locally standard graphs produced by that strategy cannot, in general, be bounded by the RHS of \eqref{light_weight_term_moment_bound}, as required. 
	Moreover, that approach does not track the detailed structure within each molecule, making it impossible to detect when such unfavorable configurations arise. Even if we stop the expansion earlier to avoid them, the resulting graphs may still fail to satisfy sufficiently strong bounds. 
	These limitations necessitate the development of a new expansion strategy. In our approach, the expansions are carried out more carefully, so as to avoid generating such unfavorable graphs; more efficiently, by eliminating unnecessary operations present in the previous strategy; and more transparently, by explicitly tracking the internal structure of each molecule throughout the expansion.
	
\medskip

	To explain the basic idea of our strategy, recall from the discussion after \eqref{eq:Gxixy} that our goal is roughly to generate $p$ additional $G$-edges, one for each tadpole diagram, in order to ensure both the \emph{sufficient $\Psi$-factors} condition and the \emph{nested property} described above. In particular, these additional edges should be ``\emph{long}'', in the molecular sense that they connect different molecules. By contrast, $G$-edges generated within a molecule (see, e.g., the second graph in Figure \ref{tikz_deformed_tadpoles}) yield only additional $\Psi(0)$-factors. Such long edges arise when one tadpole diagram pulls in an edge from another tadpole diagram, as illustrated by the graphs in Figure \ref{tikz_glued_two_tadpoles}.    
    Starting from a single tadpole diagram, which serves as the ``tail'', we repeatedly pull in long edges from other tadpoles, which collectively form the ``body''. 
    Figuratively, this procedure connects a chain of tadpoles to form a ``\emph{headless snake}''. The procedure stops once the headless snake either acquires a ``head'' and becomes a ``snake'', or bites one of its own body sections (or its tail) and becomes an ``\emph{ouroboros}''; see Definitions \ref{def_graph_types} and \ref{def_graph_snakes} for the formal definitions.\footnote{In mythology, the Ouroboros bites only its own tail; here we also allow it to bite its own body.} From a technical perspective, both termination mechanisms establish the nested property. Moreover, the edge-pulling operations automatically generate sufficiently many new edges to ensure the sufficient \(\Psi\)-factors condition.

Explicitly, in the first step of the expansion strategy, we select a tadpole diagram $f_{x_1y_1}(G)$ from the graph $|f_{xy}(G)|^p$ as the ``quasi-tail'' section.\footnote{Here, the word ``tail'' means ``the first section''. The prefix ``quasi'' indicates that this first section has not yet been connected to another section and hence is not yet an actual tail; see \Cref{def_graph_types}.} 
We then apply the $GG$-expansion \eqref{out_GG_expansion} or
\eqref{in_GG_expansion}, which produces eight new graphs corresponding to the
eight terms on the RHS of these expansions.
For the graphs produced by the first two terms (the reshaped terms), the expansion procedure terminates immediately, since these graphs already satisfy the nested property and the sufficient $\Psi$-factors condition. They can therefore be bounded directly, as illustrated by the first graph in \Cref{tikz_deformed_tadpoles}. The remaining $(p-1)$ unaffected tadpole diagrams can be handled by a simple application of Young's inequality.\footnote{That is, if we can show \smash{\(\E\absa{f_{xy}\p{G}}^{p}\prec \qa{\eta^{-1}\Psi\p{0}\Psi\p{|x-y|}}^k \E \absa{f_{xy}\p{G}}^{p-k}\)}, where $(p-k)$ denotes the number of unaffected sections, then applying Hölder's and Young's inequalities yields \smash{\(\E\absa{f_{xy}\p{G}}^{p}\prec N^\e \qa{\eta^{-1}\Psi\p{0}\Psi\p{|x-y|}}^p +N^{-\e}\E \absa{f_{xy}\p{G}}^{p}\)}. Solving for $\E\absa{f_{xy}\p{G}}^{p}$ and using the arbitrariness of $\e$, we obtain \eqref{light_weight_term_moment_bound}.} These two graphs provide the basic examples of a ``quasi-tail'' section undergoing a
``$\mathsf{reshape}$'' operation and thereby becoming a
``\emph{tail-head}'', which turns the graph into a ``snake''. See
Definitions \ref{def_graph_types}--\ref{def_graph_operations} for the formal
definitions of the terminology used here. 
On the other hand, the remaining six graphs require further expansion and fall into the following two classes.
    \begin{itemize}
        \item[(A)] The new graphs corresponding to the third through sixth terms (the growing terms) and the last two terms (the $\partial$-terms) on the RHS of \eqref{out_GG_expansion} or \eqref{in_GG_expansion}, where the $\partial$-derivative acts on a $G$-edge within the current section, still contain only the initial quasi-tail section. They are not connected to any other section and therefore require further expansion, since they do not satisfy the nested property; see the second graph in \Cref{tikz_deformed_tadpoles} for an example. The operations transforming the graph $f_{x_1y_1}(G)$ into these new graphs belong to the class of ``$\mathsf{grow}$'' operations. As illustrated by graphs (d)--(i) in \Cref{tikz_out_edge_expansion}, these new graphs satisfy the following four crucial properties:
        \begin{itemize}
            \item The number of internal molecules remains one.
            \item At least one of the external vertices $x_1,y_1$ remains incident to a $GG$-pair. Therefore, inductively, as long as a quasi-tail section is acted on by a $\mathsf{grow}$ operation, the resulting graphs can be expanded again via a $GG$-expansion.\footnote{At each step, the expansion must be performed with respect to an edge incident to an external vertex. This ensures the proper molecular-graph structure and preserves the nested property when the expansion induces a $\mathsf{reshape}$ operation.} \item  The number of $G$-edges, and hence the scaling order of the graph, increases by at least one. Therefore, heuristically, a $\mathsf{grow}$ operation reduces the size of the graph by at least a factor $\Psi(0)\leq W^{-c}$ for some constant $c>0$.
            \item The $\mathsf{grow}$ operation preserves the simple internal structure of the molecule, even under repeated applications (see the formal definition of quasi-tail sections in \Cref{def_graph_types}).  This structural simplicity allows us to extract the required $\Psi$-factors from the generated $G$-edges throughout the expansion.
        \end{itemize}

        \item[(B)] For the graphs generated from the last two terms (the $\partial$-terms), the $\partial$-derivative may act on a $G$-edge belonging to a new tadpole diagram. In this case, the new tadpole is pulled in as a new effective section, while the original quasi-tail is ``completed'' and becomes a genuine ``tail'' section. See the blue part of the graphs in \Cref{tikz_glued_two_tadpoles} for an example of such a tail section, and \Cref{def_graph_types} $\txt{(iii)}$ for the formal definition of the completion of a quasi-tail.
    \end{itemize}

     Based on the above discussion of the first step of the expansion strategy, the possible outcomes can be roughly divided into the following three cases:
    \begin{itemize}
        \item The tadpole diagram is successively acted on by a sequence of $\mathsf{grow}$ operations, followed by a $\mathsf{reshape}$ operation. In this case, the expansion procedure terminates.
        
        \item The tadpole diagram is acted on by sufficiently many $\mathsf{grow}$ operations so that its contribution becomes $\opr{W^{-D}}$ for a large constant $D\geq 0$. In this case, we also stop expanding the graph.
        
        \item The tadpole diagram is successively acted on by a sequence of $\mathsf{grow}$ operations, after which a new effective section is pulled in. The number of effective sections then increases to two.
    \end{itemize}
    In the last case, the two graphs in \Cref{tikz_glued_two_tadpoles} represent, roughly speaking, the two possible types of outcomes. Each graph contains two effective sections, and hence, heuristically, requires two additional long $G$-edges for control. However, for the first graph, a direct computation shows that the required nested property has already been established. Therefore, only one additional $G$-edge is needed to provide sufficient $\Psi$-factors. The missing factor is supplied by the \emph{leaf light-weight} in the blue section through the averaged local law \eqref{average_bound_2_to_1}, which yields one additional $\Psi(0)$-factor compared with the entrywise bound \eqref{entrywise_bound_2_to_1}. In this case, the new section is called a ``weight-head'' section, and the operation that inserts this new section into the preceding one is called $\mathsf{insert}$. This operation turns the graph into a genuine snake. 
    
    In contrast, the second graph in \Cref{tikz_glued_two_tadpoles} requires further expansion. The operation that extends the preceding graph to the current one is called $\mathsf{extend}$, and the new section (the red part) is called a ``quasi-body'' section. In general, a quasi-body section has almost identical structure to the one in this example: it consists of a \emph{``free'' solid edge} and a disjoint component, called the \emph{quasi-tail component} (in the second graph of \Cref{tikz_glued_two_tadpoles}, this is the tadpole diagram component containing $a_2$ and $b_2$). Therefore, we can apply the same expansion strategy as the one used in the first step for $f_{x_1y_1}(G)$. Ignoring the free solid edge, the resulting graphs have exactly the same local structure as those obtained above, as long as only the growing terms are generated or the $\partial$-derivative acts only on $G$-edges within the \emph{quasi-tail component}.\footnote{We again emphasize that the $GG$-expansion should always be applied with respect to an edge incident to an external vertex. Here, in the second graph of \Cref{tikz_glued_two_tadpoles}, this edge is $\overline{G}_{a_2y_2}$.}
    Otherwise, if the graph is generated by the $\partial$-terms in \eqref{out_GG_expansion} or \eqref{in_GG_expansion} with the $\partial$-derivative acting on an edge outside the quasi-tail component, the resulting graphs can be roughly divided into the following three cases, which are again illustrated using the second graph in \Cref{tikz_glued_two_tadpoles}:    
\begin{itemize}
        \item[(1)] If the $\partial$-derivative acts precisely on the free solid edge $\ol{G}_{x_2a_1}$ of the quasi-body section, the resulting graph has a similar graphical structure to the first graph in \Cref{tikz_glued_two_tadpoles}. Hence, it satisfies both the nested property and the sufficient $\Psi$-factors condition, and becomes a genuine snake. In this case, the operation and the resulting section are referred to as a $\mathsf{reshape}$ of the second type and a ``body-head II'' section, respectively.\footnote{A ``Body-head I'' section is obtained by applying a $\mathsf{reshape}$ operation to the quasi-tail component in the quasi-body.}
        \item[(2)] If the $\partial$-derivative acts on a $G$-edge in the tail section, the headless snake pulls two edges from the tail section into the quasi-tail component of the quasi-body section. The resulting graph then immediately acquires both the nested property and the sufficient $\Psi$-factors condition, and becomes an ouroboros. In this case, the operation is referred to as a $\mathsf{bite}$ operation.
        \item[(3)] If the $\partial$-derivative acts on a new tadpole diagram, this tadpole is pulled in and becomes a new quasi-body section. We then move to this new section and repeat the above expansion strategy for the quasi-body discussed above.
    \end{itemize}

Following the above expansion strategy, we generate a sum of well-structured \emph{snake-like graphs}, which, as the name suggests, consist of a ``tail'', several intermediate ``body'' sections, and a ``head'' section, unless they are headless snakes or ouroboroi. A key feature of snake-like graphs is that their molecular graphs belong to a constrained family with simple (linear) structures. This structure allows us to extract the required $\Psi$-factors and apply Ward's identities successively according to an explicit order; see \Cref{fig_snake_mole} and \Cref{lemma_section_structure_aux}.
Moreover, for any fixed constant $D>0$, the expansion terminates after $\OO(1)$ steps, up to an error of order $\opr{W^{-D}}$, and produces a sum of snake-like graphs that are not headless snakes (see \Cref{lemma_terminate}). Indeed, every $\mathsf{grow}$ operation reduces the size of the graph by at least a factor $W^{-c}$ for some constant $c>0$, whereas every other operation either terminates the expansion immediately or pulls in a new tadpole diagram. Since there are at most $p$ tadpole diagrams, the total
number of such operations is bounded.

	We conclude this subsection by comparing our strategy with that of the previous work \cite{dubova2025delocalizationnonmeanfieldrandommatrices,yang2021delocalization}. Compared with \cite{dubova2025delocalizationnonmeanfieldrandommatrices}, our strategy is more streamlined in the following sense. The approach in \cite{dubova2025delocalizationnonmeanfieldrandommatrices}, building on \cite{yang2021delocalization}, applies both weight expansions and $GG$-pair expansions in order to remove \emph{all} weights and $GG$-pairs. The objective there is to retain only $G\overline G$-pairs in the resulting graphs, which then allows the application of the so-called $T$-expansions of $G\overline G$-pairs.
	In contrast, our goal is more modest: we only need to generate one additional long $G$-edge for each tadpole diagram, and we do not require the graph to be locally standard. Consequently, our expansion procedure is more efficient. In particular, we no longer perform weight expansions, nor do we attempt to eliminate all $GG$-pairs in the headless snakes. For instance, in the second graph of Figure \ref{tikz_glued_two_tadpoles}, if the section corresponding to $x_1$ and $y_1$ pulls in a blue solid edge and forms two $GG$-pairs, we do not expand these pairs further; instead, we move to the new section and continue the $GG$-expansion there.
	As a result, the structures of the resulting graphs become more transparent. Moreover, to avoid the unfavorable configurations described in \Cref{tikz_bad_graphs}, we always expand a $GG$-pair in the newest section of the headless snake that is incident to at least one external solid edge (i.e., a solid edge connected to an external vertex $x$ or $y$). This ensures that each molecule is incident to at most five long $G$-edges, and that each vertex carries at most two $G$-edges throughout the graph reduction procedure based on \eqref{eq:reduction_xi}.

	\subsection{Expansion strategy}
	
	We first define the constituent components of the snake-like graphs.
	
	\begin{definition}[Sections]\label{def_graph_types}
		A snake-like graph is composed of single-colored sections of the following types. 
		\begin{itemize}
			
			\item[(i)] A \emph{quasi-tail} section with external vertices $x,y$ is a graph satisfying the following properties. The subgraph formed by all internal vertices together with the waved edges is a tree. Each internal vertex either (1) is incident to exactly two solid edges pointing in opposite directions (i.e., one is an in-edge and the other an out-edge at that vertex), (2) is attached to exactly one light-weight, or (3) is incident to no $G$-edge.
			Each of the vertices $x$ and $y$ is connected to exactly one internal vertex by a solid edge, referred to as the \emph{$x$-solid edge} and \emph{$y$-solid edge}, respectively. Moreover, these two external solid edges point in opposite directions: for example, if the $x$-solid edge is an in-edge at $x$, then the $y$-solid edge is an out-edge at $y$, and vice versa. In Figure~\ref{tikz_quasitail}, we illustrate several examples of quasi-tail sections with external vertices $x,y$.

			\begin{figure}[htbp]\label{tikz_example_quasi_tails}
				\vspace{-0.4cm}
				\centering
				\begin{minipage}{0.4\textwidth}
					\centering
					\begin{tikzpicture}[scale=0.4, 
						fermion/.style={
							postaction={decorate},
							decoration={markings, mark=at position 0.55 with {\arrow[scale=1.2]{Latex}}} 
						},
						photon/.style={
							decorate,
							decoration={snake, segment length=5pt, amplitude=1.2pt}
						},
						dot/.style={
							circle, fill=black, inner sep=1.0pt
						}
						]
						
						\node [dot, label=left:{\Large $x$}] (x)  at (-3.0, -0.5) {};
						\node [dot] (v1) at (-1.5, -2.5) {};
						\node [dot] (v2) at (-0.2,  2.0) {};
						\node [dot] (v3) at (-0.2, -0.2) {};
						\node [dot] (v4) at ( 2.0, -2.5) {};
						\node [dot] (v5) at ( 2.2, -0.5) {};
						\node [dot] (v6) at ( 3.2,  1.5) {};
						\node [dot, label=right:{\Large $y$}] (y)  at ( 4.5, -1.0) {};
						
						\draw [fermion] (x) to[bend right=10] (v1);
						\draw [fermion] (v1) to[bend left=10] (v2);
						\draw [fermion] (v2) to[bend left=5] (v5);
						\draw [fermion] (v5) to[bend left=15] (y);
						
						\draw [fermion] (v3) to[bend left=15] (v4);
						\draw [fermion] (v4) to[bend left=15] (v3);
						
						\draw [photon] (v1) -- (v3);
						\draw [photon] (v1) to[bend right=10] (v4);
						\draw [photon] (v2) to[bend right=10] (v3);
						\draw [photon] (v3) to[bend left=10] (v5);
						\draw [photon] (v5) -- (v6);
						
						\draw (v6) to[out=45, in=135, loop, min distance=1.5cm] (v6);
						
					\end{tikzpicture}
				\end{minipage}    \hspace{-2.0cm}
				\begin{minipage}{0.4\textwidth}
					\centering
					\begin{tikzpicture}[scale=0.4, 
						fermion/.style={
							postaction={decorate},
							decoration={markings, mark=at position 0.55 with {\arrow[scale=1.2]{Latex}}}
						},
						photon/.style={
							decorate,
							decoration={snake, segment length=5pt, amplitude=1.2pt}
						},
						dot/.style={
							circle, fill=black, inner sep=1.0pt     }
						]
						
						\node [dot, label=below left:{\Large $x$}] (x)  at (-2.5, -3.0) {};
						\node [dot, label=below right:{\Large $y$}] (y)  at ( 2.0, -3.0) {};
						
						\node [dot] (d) at ( 0.0, -0.5) {};   \node [dot] (a) at (-2.5,  1.5) {};   \node [dot] (b) at ( 0.5,  3.0) {};   \node [dot] (c) at ( 3.0,  0.5) {}; 
						
						\draw [fermion] (x) -- (d);
						\draw [fermion] (d) -- (y);
						
						\draw [fermion] (a) -- (c);
						
						\draw [fermion] (c) to[bend right=30] (b);
						\draw [fermion] (b) to[bend right=20] (a);
						
						\draw [photon] (a) to[bend right=30] (d);
						\draw [photon] (d) to[bend right=15] (c);
						
						\draw [photon] (b) to[bend right=5]  (c);
						
					\end{tikzpicture}
				\end{minipage}
				\hspace{-2.0cm}
				\begin{minipage}{0.4\textwidth}
					\centering
					\begin{tikzpicture}[scale=0.4, 
						fermion/.style={
							postaction={decorate},
							decoration={markings, mark=at position 0.55 with {\arrow[scale=1.2]{Latex}}}
						},
						photon/.style={
							decorate,
							decoration={snake, segment length=5pt, amplitude=1.2pt}
						},
						dot/.style={
							circle, fill=black, inner sep=1.0pt     }
						]
						
						\node [dot, label=below left:{\Large $x$}] (x)  at (-2.5, -3.0) {};
						\node [dot, label=below right:{\Large $y$}] (y)  at ( 2.0, -3.0) {};
						
						\node [dot] (d) at ( 0.0, -0.5) {};   \node [dot] (a) at (-2.5,  1.5) {};   \node [dot] (b) at ( 0.5,  3.0) {};   \node [dot] (c) at ( 3.0,  0.5) {}; 
						
						\draw [fermion] (x) -- (d);
						\draw [fermion] (d) -- (y);

						\draw [fermion] (b) to[bend right=20] (a);
						\draw [fermion] (a) to[bend right=20] (b);
						
						\draw [photon] (a) to[bend right=30] (d);
						\draw [photon] (d) to[bend right=15] (c);
						
						\draw [photon] (b) to[bend right=5]  (c);
						
					\end{tikzpicture}
				\end{minipage}
				\caption{Examples of quasi-tail sections. For simplicity, the colors are not shown.} \label{tikz_quasitail}
			\end{figure}

			Roughly speaking, quasi-tail sections are introduced to represent the possible subgraphs generated from a tadpole diagram by iteratively applying $GG$-expansions and selecting either a growing term or a $\partial$-term in which the $\partial$-derivative acts on one of the $G$-edges within the section. In other words, this corresponds to case (A) in the first step of expansion discussed in the previous section. For instance, the second graph in Figure~\ref{tikz_deformed_tadpoles} is a blue quasi-tail section obtained from such an operation applied to a tadpole diagram.
            
			\item[(ii)] 
			A \emph{quasi-body} section with external vertices $x,y$ and outer vertices $a,b$ is a graph consisting of a free solid edge between $x$ (resp.~$y$) and $a$, referred to as the $x$-solid edge (resp.~$y$-solid edge), together with a quasi-tail with external vertices $b$ and $y$ (resp.~$b$ and $x$). Moreover, we assume that the $x$- and $y$-solid edges point in opposite directions. An example of a quasi-body section is illustrated in Figure~\ref{tikz_quasibody}, where, for simplicity, the colors of the solid edges are not shown.
			
			\begin{figure}[htbp]\label{tikz_example_quasi_body}
				\centering
				\vspace{-0.4cm}
				\begin{minipage}{0.45\textwidth}
					\centering
					\begin{tikzpicture}[
						scale=0.4, 
						fermion/.style={
							postaction={decorate},
							decoration={markings, mark=at position 0.55 with {\arrow[scale=1.2]{Latex}}}
						},
						photon/.style={
							decorate,
							decoration={snake, segment length=5pt, amplitude=1.2pt}
						},
						dot/.style={
							circle, fill=black, inner sep=1.0pt
						}
						]
						
						\node [dot, label=below:{\Large $x$}] (x)  at (-7.0, 0.0) {};
						\node [dot, label=above:{\Large $a$}] (a)  at (-1.0,  3.0) {};
						
						\node [dot, label=above:{\Large $b$}] (b)  at ( 1.0,  3.0) {};
						\node [dot, label=above:{\large $v_3$}] (v1) at ( 3.5,  2.0) {};         \node [dot, label=right:{\large $v_2$}] (v2) at ( 6.0,  2.0) {};         \node [dot, label=below:{\large $v_1$}] (v3) at ( 6.0, 0.0) {};         \node [dot, label=right:{\Large $y$}] (y)  at ( 10.0, 0.0) {};

						\draw [fermion] (x) -- (a);
						
						\draw [fermion] (b) -- (v1);
						\draw [fermion] (v1) -- (v2);
						
						\draw [fermion] (v2) to[bend left=45] (v3);
						
						\draw [fermion] (v3) -- (y);
						
						\draw [photon] (v1) to[bend right=10] (v3);
						
						\draw [photon] (v2) to[bend right=15] (v3);
						
					\end{tikzpicture}
				\end{minipage}
				\caption{An example of a quasi-body section.}\label{tikz_quasibody}
			\end{figure}
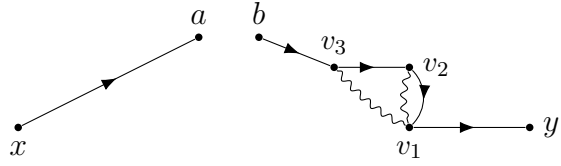
			
			In our expansion scheme, the first quasi-body section arises from a tadpole diagram, denoted by $\cal T_k$, as follows. We apply a $GG$-expansion to a $GG$-pair in another section and select the $\partial$-term in which the $\partial$-derivative acts on a solid edge $e$ inside $\cal T_k$. The edge $e$ is then transformed into two edges $e_1$ and $e_2$, which, together with the remaining edges in $\cal T_k$, form a quasi-body section.
			This type of quasi-body section will be referred to as a \emph{tadpole-quasi-body}, which consists of a free solid edge and a tadpole diagram. For instance, the red solid edges together with the waved edge $(a_2,b_2)$ in the second graph of Figure~\ref{tikz_glued_two_tadpoles} form a tadpole-quasi-body section with outer vertices $a_1$ and $a'_1$.       
			More generally, similarly to the quasi-tail sections, the notion of a quasi-body section is introduced to represent the possible subgraphs generated from a tadpole-quasi-body section by iteratively applying $GG$-expansions with respect to a $GG$-pair containing an external solid edge and selecting either a growing term or a $\partial$-term in which the $\partial$-derivative acts on one of the $G$-edges within the quasi-tail component.  \end{itemize}   
		Both quasi-tail and quasi-body sections will collectively be referred to as \emph{quasi-type sections}. A crucial feature of quasi-type sections is that at least one external edge belongs to a $GG$-pair to which a $GG$-expansion can be applied. In this expansion, if a $\partial$-term is selected, the $GG$-pair is replaced by an extended $GG$-pair. In this case, we say that the quasi-type section is completed and becomes a regular \emph{tail} or \emph{body} section.
		\begin{itemize}
			\item[(iii)] The set of (regular) \emph{tail} sections with external vertices $x,y$ and inner vertices $a,b$ is constructed as follows. First, take any quasi-tail section $\fS$ with external vertices $x,y$, and choose any $\star\in\ha{x,y}$. Next, identify the unique $GG$-pair containing the external edge incident to $\star$. Finally, replace this $GG$-pair by any of the corresponding extended $GG$-pairs with inner vertices $a,b$. We refer to this procedure as the \emph{completion} of the quasi-tail $\fS$ with respect to $\star$.

			\item[(iv)] The set of (regular) \emph{body} sections with external vertices $x,y$, outer vertices $a,b$, and inner vertices $a',b'$ is constructed analogously by completing quasi-body sections, via a completion of the quasi-tail component as described above. The only difference is that the vertex $\star$ chosen for the completion must be an external vertex of the quasi-body section, rather than an outer vertex.
			For example, Figure~\ref{tikz_body_section} shows a completion of the quasi-body section in Figure~\ref{tikz_quasibody} with respect to the vertex $y$, rather than the vertex $b$, which is also an external vertex of the quasi-tail component but serves as an outer vertex of the quasi-body section.         \begin{figure}[htbp]\label{tikz_example_body}
				\centering
				\vspace{-0.4cm}
				\begin{minipage}{0.45\textwidth}
					\centering
					\begin{tikzpicture}[
						scale=0.4, 
						fermion/.style={
							postaction={decorate},
							decoration={markings, mark=at position 0.55 with {\arrow[scale=1.2]{Latex}}}
						},
						photon/.style={
							decorate,
							decoration={snake, segment length=5pt, amplitude=1.2pt}
						},
						dot/.style={
							circle, fill=black, inner sep=1.0pt
						}
						]
						
						\node [dot, label=below:{$x$}] (x)  at (-7.0, 0.0) {};
						\node [dot, label=above:{$a$}] (a)  at (-1.0,  3.0) {};
						
						\node [dot, label=above:{$b$}] (b)  at ( 1.0,  3.0) {};
						\node [dot] (v1) at ( 3.5,  2.0) {};         \node [dot] (v2) at ( 6.0,  2.0) {};         \node [dot, label=below:{$a'$}] (v3) at ( 6.0, 0.0) {};         \node [dot, label=below:{$b'$}] (v4) at ( 8.0, 0.0) {};
						\node [dot, label=right:{$y$}] (y)  at ( 10.0, 0.0) {};

						\draw [fermion, blue] (x) -- (a);
						
						\draw [fermion, blue] (b) -- (v1);
						\draw [fermion, blue] (v1) -- (v2);
						
						\draw [fermion, blue] (v2) to[bend left=45] (v3);
						
						\draw [fermion, blue] (v4) -- (y);
						
						\draw [photon] (v1) to[bend right=10] (v3);
						
						\draw [photon] (v2) to[bend right=15] (v3);
						\draw [photon] (v3) to (v4);
						
					\end{tikzpicture}
				\end{minipage}
				\caption{A body section obtained by the completion of Figure \ref{tikz_quasibody}.}\label{tikz_body_section}
				\vspace{-0.4cm}
			\end{figure}
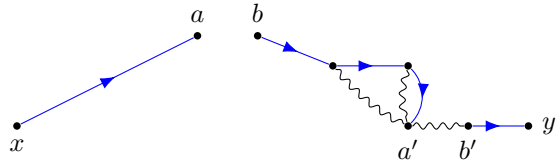
		\end{itemize}
		Finally, we introduce four types of \emph{head} sections with external vertices $x,y$.         \begin{itemize}
			\item[(v)] {\bf Tail-head}: Take any quasi-tail with external vertices $x,y$ and choose one of its external vertices $\star$. We then replace the $GG$-pair containing the $\star$-solid edge with the corresponding dotted structure or waved structure generated from the $GG$-expansion with respect to the solid edge incident to $\star$ (recall the notation introduced below Lemma \ref{lemma:GGexpansion}), identifying the two vertices connected by a dotted edge when present. In other words, a tail-head section is generated from the reshaped terms in the expansion of a $GG$-pair inside a quasi-tail.
			
			\item[(vi)] {\bf Body-head I} with outer vertices $a,b$: Take any quasi-body section with external vertices $x,y$ and outer vertices $a,b$, and choose an external vertex $\star\in\{x,y\}$ belonging to the quasi-tail component. We then replace the $GG$-pair containing the $\star$-solid edge with the corresponding dotted structure or waved structure generated from the $GG$-expansion with respect to the solid edge incident to $\star$, again identifying the two vertices connected by a dotted edge if such an edge appears. In other words, a body-head of the first type is generated from the reshaped terms in the expansion of a $GG$-pair inside a quasi-body.
			
			\item[(vii)] {\bf Body-head II} with outer vertices $a,b$: Take any body section with external vertices $x,y$, outer vertices $a,b$, and inner vertices $c,d$, and pull the free solid edge in the section to the inner vertices in the following sense. In a blue body section, given an edge $G_{uv}$ and inner vertices $c,d$, where $c$ is attached to a unique $G$-edge $G_{cc'}$ and $d$ is attached to a unique $G$-edge $G_{d'd}$, pulling $G_{uv}$ to vertices $c,d$ means replacing the factor $G_{uv}G_{cc'}G_{d'd}$ with $G_{uc}G_{dv}G_{cc'}G_{d'd}$. 
			Similar operations can be defined for red body sections. In other words, a body-head of the second type is generated from the $\partial$-terms in the expansion of a $GG$-pair inside the section, where the $\partial$-derivative acts on its own free solid edge.
			
			\item[(viii)] {\bf Weight-head} with outer vertices $a,b$: Take a tadpole diagram \smash{$f_{xy}\pa{G}$} or \smash{$\ol{f_{xy}\pa{G}}$}, and replace the light-weight at vertex $w$ in this diagram by a pair of solid edges $(a,w)$ and $(w,b)$.
			In other words, a weight-head arises from the $\partial$-terms in the expansion of a $GG$-pair inside a quasi-tail or quasi-body section, where the $\partial$-derivative acts on the light-weight in another tadpole diagram.
		\end{itemize}
		
		We will refer to a quasi-tail, tail, or tail-head section as a \emph{tail-type section}, and to a quasi-body or body section as a \emph{body-type section}. A head section that is \emph{not} a tail-head will be called a \emph{regular head}. For any section $\fS$ with inner or outer vertices, we denote the corresponding vertex sets by $V_{\txt{in}}\pa{\fS}$ and $V_{\txt{out}}\pa{\fS}$, respectively. These vertices serve as ``connectors'' between different sections, along which the sections are joined together.
	\end{definition}
	
	We now define three types of snake-like graphs constructed from the sections introduced above.
	
	\begin{definition}[Snake-like graphs]\label{def_graph_snakes}
		Consider a family of graphs $\cG$ formed by $q$ sections, i.e., $\cG=\fS_1\fS_2\cdots\fS_q$, for some $q\in\qq{p}$. When $q=1$, we let $\fS_1$ be either a quasi-tail or a tail-head section. For $q\ge 2$, suppose that the following properties hold:
		\begin{itemize}
			\item $\fS_1$ is a tail section.
			
			\item $\fS_i$ is a body section for $1<i<q$.
			
			\item $\fS_q$ is either a quasi-body section or a regular head.
			
			\item For any $i\in\qq{q-1}$, we have $V_{\txt{in}}\pa{\fS_i}=V_{\txt{out}}\pa{\fS_{i+1}}$.
			
			\item For any $i\in\qq{q}$, all $G$-edges in $\fS_i$ have the same color.
			
		\end{itemize}
		We call $\cG$ a \emph{\bf headless snake} if $\fS_q$ is a quasi-type section, and a \emph{\bf snake} if $\fS_q$ is a head section. 
		
		An \emph{\bf ouroboros} is a graph constructed from a headless snake $\cG=\fS_1\fS_2\cdots\fS_q$ with $q\geq 2$ as follows. First, we complete the quasi-body section $\fS_q$ with respect to one of its external vertices $x$ or $y$, and denote the resulting section by $\overline\fS_q$. This completion produces two inner vertices $a$ and $b$. Next, we pull a $G$-edge (possibly a light-weight) from $\fS_1\cdots\fS_{q-1}$ to $a,b$ in the following sense. 
		If $\overline\fS_q$ is a blue body section where $a$ is attached to a unique $G$-edge $G_{aa'}$ and $b$ is attached to a unique $G$-edge $G_{b'b}$, then pulling an edge $G_{uv}$ to the vertices $a,b$ means replacing the factor $G_{uv}G_{aa'}G_{b'b}$ with $G_{ua}G_{bv}G_{aa'}G_{b'b}$. Similarly, pulling an edge $\overline G_{uv}$ to $a,b$ means replacing the factor $\overline G_{uv}G_{aa'}G_{b'b}$  with $\overline G_{ub}\overline G_{av}G_{aa'}G_{b'b}$. Analogous operations can be defined when $\overline\fS_q$ is a red body section. In other words, an ouroboros is generated from the $\partial$-terms in the expansion of a $GG$-pair within the section $\fS_q$ of a headless snake, where the $\partial$-derivative acts on the preceding sections $\fS_1\cdots\fS_{q-1}$. If the pulled $G$-edge in the above completion procedure is a solid edge connecting two different molecules, we call the resulting graph an ouroboros of type I. Otherwise, if the pulled edge is a $G$-edge within a single molecule, we call the graph an ouroboros of type II. Moreover, the two edges produced by this pulling procedure are referred to as the distinguished edges.
		
	\end{definition}

	\begin{figure}[htbp]
		\centering
		\vspace{-0.15cm}
		\begin{tikzpicture}[
			x=0.42cm,
			y=0.47cm,
			mol/.style={circle, fill=black, inner sep=0.9pt},
			ext/.style={circle, fill=black, inner sep=0.75pt},
			edge/.style={line width=0.55pt},
			omit/.style={line width=0.45pt, gray, dash pattern=on 1.2pt off 1.2pt},
			dist/.style={line width=1.1pt},
			panel/.style={font=\scriptsize, align=center, text width=2.5cm},
			every label/.style={font=\tiny, inner sep=0.6pt}
			]

\begin{scope}[xshift=0cm]
				\node[mol] (h-a1) at (0,5.0) {};
				\node[mol] (h-a2) at (0,3.75) {};
				\node[mol] (h-a3) at (0,2.50) {};
				\node[mol] (h-aqm) at (0,0.25) {};
				\node[mol] (h-aq) at (0,-1.0) {};
				\node[ext,label=left:$x_1$] (h-x1) at (-1.65,4.8) {};
				\node[ext,label=right:$y_1$] (h-y1) at (1.65,4.8) {};
				\node[ext,label=left:$x_2$] (h-x2) at (-1.65,3.55) {};
				\node[ext,label=right:$y_2$] (h-y2) at (1.65,3.55) {};
				\node[ext,label=left:$x_3$] (h-x3) at (-1.65,2.30) {};
				\node[ext,label=right:$y_3$] (h-y3) at (1.65,2.30) {};
				\node[ext,label=left:$x_q$] (h-xq) at (-1.65,-1.2) {};
				\node[ext,label=right:$y_q$] (h-yq) at (1.65,-1.2) {};
				\draw[edge] (h-x1)--(h-a1)--(h-y1);
				\draw[edge] (h-x2)--(h-a1)--(h-a2)--(h-y2);
				\draw[edge] (h-x3)--(h-a2)--(h-a3)--(h-y3);
				\draw[omit] (h-a3)--(h-aqm);
				\draw[edge] (h-xq)--(h-aqm)--(h-aq)--(h-yq);
				\node[panel] at (0,-2.35) {(a) headless};
			\end{scope}

\begin{scope}[xshift=2.70cm]
				\node[ext,label=left:$x_1$] (t-x1) at (-1.35,1.80) {};
				\node[ext,label=right:$y_1$] (t-y1) at (1.35,1.80) {};
				\draw[edge] (t-x1)--(t-y1);
				\node[panel] at (0,-2.35) {(b) tail-head\\($q=1$)};
			\end{scope}

\begin{scope}[xshift=5.40cm]
				\node[mol] (bi-a1) at (0,5.0) {};
				\node[mol] (bi-a2) at (0,3.75) {};
				\node[mol] (bi-a3) at (0,2.50) {};
                \node[mol] (bi-bq-1) at (0,1.00) {};
				\node[mol] (bi-bq) at (0,0.00) {};
				\node[ext,label=left:$x_1$] (bi-x1) at (-1.65,4.8) {};
				\node[ext,label=right:$y_1$] (bi-y1) at (1.65,4.8) {};
				\node[ext,label=left:$x_2$] (bi-x2) at (-1.65,3.55) {};
				\node[ext,label=right:$y_2$] (bi-y2) at (1.65,3.55) {};
				\node[ext,label=left:$x_3$] (bi-x3) at (-1.65,2.30) {};
				\node[ext,label=right:$y_3$] (bi-y3) at (1.65,2.30) {};
				\node[ext,label=left:$x_{q-1}$] (bi-xq-1) at (-1.65,-0.20) {};
				\node[ext,label=right:$y_{q-1}$] (bi-yq-1) at (1.65,-0.20) {};
                \node[ext,label=left:$x_q$] (bi-xq) at (-1.65,-1.2) {};
				\node[ext,label=right:$y_q$] (bi-yq) at (1.65,-1.2) {};
				\draw[edge] (bi-x1)--(bi-a1)--(bi-y1);
				\draw[edge] (bi-x2)--(bi-a1)--(bi-a2)--(bi-y2);
				\draw[edge] (bi-x3)--(bi-a2)--(bi-a3)--(bi-y3);
				\draw[omit] (bi-a3)--(bi-bq);
                \draw[edge] (bi-xq-1)--(bi-bq-1)--(bi-bq)--(bi-yq-1);
				\draw[edge] (bi-xq)--(bi-bq)--(bi-yq);
				\node[panel] at (0,-2.35) {(c) body-head I};
			\end{scope}

\begin{scope}[xshift=8.10cm]
				\node[mol] (bii-a1) at (0,5.0) {};
				\node[mol] (bii-a2) at (0,3.75) {};
				\node[mol] (bii-a3) at (0,2.50) {};
				\node[mol] (bii-bq) at (0,0.25) {};
				\node[mol] (bii-aq) at (0,-1.0) {};
				\node[ext,label=left:$x_1$] (bii-x1) at (-1.65,4.8) {};
				\node[ext,label=right:$y_1$] (bii-y1) at (1.65,4.8) {};
				\node[ext,label=left:$x_2$] (bii-x2) at (-1.65,3.55) {};
				\node[ext,label=right:$y_2$] (bii-y2) at (1.65,3.55) {};
				\node[ext,label=left:$x_3$] (bii-x3) at (-1.65,2.30) {};
				\node[ext,label=right:$y_3$] (bii-y3) at (1.65,2.30) {};
				\node[ext,label=left:$x_q$] (bii-xq) at (-1.65,-1.2) {};
				\node[ext,label=right:$y_q$] (bii-yq) at (1.65,-1.2) {};
				\draw[edge] (bii-x1)--(bii-a1)--(bii-y1);
				\draw[edge] (bii-x2)--(bii-a1)--(bii-a2)--(bii-y2);
				\draw[edge] (bii-x3)--(bii-a2)--(bii-a3)--(bii-y3);
				\draw[omit] (bii-a3)--(bii-bq);
				\draw[edge] (bii-xq)--(bii-aq)--(bii-yq);
				\draw[edge] (bii-bq) to[bend left=18] (bii-aq);
				\draw[edge] (bii-bq) to[bend right=18] (bii-aq);
				\node[panel] at (0,-2.35) {(d) body-head II\\or weight-head};
			\end{scope}

\begin{scope}[xshift=10.80cm]
				\node[mol] (oi-a1) at (0,5.0) {};
				\node[mol] (oi-a2) at (0,3.75) {};
				\node[mol] (oi-a3) at (0,2.50) {};
				\node[mol] (oi-aqm) at (0,0.25) {};
				\node[mol] (oi-aq) at (0,-1.0) {};
				\node[ext,label=left:$x_1$] (oi-x1) at (-1.65,4.8) {};
				\node[ext,label=right:$y_1$] (oi-y1) at (1.65,4.8) {};
				\node[ext,label=left:$x_2$] (oi-x2) at (-1.65,3.55) {};
				\node[ext,label=right:$y_2$] (oi-y2) at (1.65,3.55) {};
				\node[ext,label=left:$x_3$] (oi-x3) at (-1.65,2.30) {};
				\node[ext,label=right:$y_3$] (oi-y3) at (1.65,2.30) {};
				\node[ext,label=left:$x_q$] (oi-xq) at (-1.65,-1.2) {};
				\node[ext,label=right:$y_q$] (oi-yq) at (1.65,-1.2) {};
				\draw[edge] (oi-x1)--(oi-a1)--(oi-y1);
				\draw[edge] (oi-x2)--(oi-a1)--(oi-a2)--(oi-y2);
				\draw[edge] (oi-x3)--(oi-a2);
				\draw[edge] (oi-a3)--(oi-y3);
				\draw[omit] (oi-a3)--(oi-aqm);
				\draw[edge] (oi-xq)--(oi-aqm)--(oi-aq)--(oi-yq);
				\draw[dist] (oi-a2) to[bend right=18] (oi-aq);
				\draw[dist] (oi-a3) to[bend left=18] (oi-aq);
				\node[panel] at (0,-2.35) {(e) ouroboros I};
			\end{scope}

\begin{scope}[xshift=13.50cm]
				\node[mol,label=above:$a_1$] (oii-a1) at (0,5.0) {};
				\node[mol] (oii-a2) at (0,3.75) {};
				\node[mol] (oii-a3) at (0,2.50) {};
				\node[mol] (oii-aqm) at (0,0.25) {};
				\node[mol] (oii-aq) at (0,-1.0) {};
				\node[ext,label=left:$x_1$] (oii-x1) at (-1.65,4.8) {};
				\node[ext,label=right:$y_1$] (oii-y1) at (1.65,4.8) {};
				\node[ext,label=left:$x_2$] (oii-x2) at (-1.65,3.55) {};
				\node[ext,label=right:$y_2$] (oii-y2) at (1.65,3.55) {};
				\node[ext,label=left:$x_3$] (oii-x3) at (-1.65,2.30) {};
				\node[ext,label=right:$y_3$] (oii-y3) at (1.65,2.30) {};
				\node[ext,label=left:$x_q$] (oii-xq) at (-1.65,-1.2) {};
				\node[ext,label=right:$y_q$] (oii-yq) at (1.65,-1.2) {};
				\draw[edge] (oii-x1)--(oii-a1)--(oii-y1);
				\draw[edge] (oii-x2)--(oii-a1)--(oii-a2)--(oii-y2);
				\draw[edge] (oii-x3)--(oii-a2)--(oii-a3)--(oii-y3);
				\draw[omit] (oii-a3)--(oii-aqm);
				\draw[edge] (oii-xq)--(oii-aqm)--(oii-aq)--(oii-yq);
				\draw[dist] (oii-a2) to[bend left=18] (oii-aq);
				\draw[dist] (oii-a2) to[bend right=18] (oii-aq);
				\node[panel] at (0,-2.35) {(f) ouroboros II};
			\end{scope}
		\end{tikzpicture}
		\vspace{-0.15cm}
		\caption{Representative molecular graphs of all types of snake-like graphs. Panels (a)--(f) show, respectively, a headless snake, a tail-head snake, a body-head-I snake, a body-head-II or weight-head snake, and ouroboroi of types I and II. The body-head-II and weight-head cases have identical molecular topology. In panels (e) and (f), the thicker edges denote the distinguished edges. The dashed segments indicate omitted intermediate body sections.}\label{fig_snake_mole}
	\end{figure}
	
As can be seen from Figure \ref{fig_snake_mole}, the molecular graphs of all types of snake-like graphs exhibit particularly simple and explicit structures (see also \Cref{lemma_section_structure_aux} below). In particular, the number of internal molecules is $q$ for panels (a), (d), (e), and (f), and $q-1$ for panels (b) and (c). We note that snake and ouroboros graphs can be bounded directly, whereas headless snakes require further expansion using the $GG$-expansion. We therefore introduce the corresponding \emph{graph operations} in the following definition.

	\begin{definition}[Graph operations]\label{def_graph_operations}
		Consider a graph $\cG=\mathscr{S}_q\cdot F_{q+1}\cdots F_{p}$, where each $F_i\in\h{f_{x_iy_i}\p{G},\ol{f_{x_iy_i}\p{G}}}$, $i\in\qq{q+1,p}$, is called a \emph{free tadpole diagram}. The graph $\mathscr{S}_q=\fS_1\cdots\fS_q$ is a headless snake with $q$ \emph{effective sections}, where the external vertices of $\fS_i$, $i\in\qq{q}$, are also denoted by $x_i,y_i$. (Recall that all pairs $(x_i,y_i)$ represent the same pair $(x,y)$.) 
		We first choose an external vertex $\star$ in the quasi-tail component of $\fS_q$ (where $\fS_q$ is a quasi-tail when $q=1$ and a quasi-body when $q\ge 2$). If two such external vertices exist, we select the one whose $\star$-solid edge points toward $\star$, for definiteness (recall that the two external solid edges are assumed to point in opposite directions). Next, we identify the unique $GG$-pair containing $\star$. For definiteness, we assume that this $GG$-pair is $G_{a\star}G_{ba}$; the case $\overline G_{a\star}\overline G_{ba}$ can be treated analogously.
		We then apply the $GG$-expansion with respect to the out-edge $G_{a\star}$ of the pair $G_{a\star}G_{ba}$ in the graph $\cal G$, yielding a sum of graphs corresponding to the eight terms in \eqref{out_GG_expansion}. The resulting operations are classified as follows: \begin{itemize}
			\item For graphs generated by the reshaped terms in \eqref{out_GG_expansion}, we identify the two vertices connected by a dotted edge in $\fS_q$, whenever such an edge is present.\footnote{This identification ensures that the scaling order increases by one.} The resulting operation is called a reshape operation of the first type.

            \item If $\fS_q$ is a quasi-body section, then for graphs induced by the $\partial$-terms in \eqref{out_GG_expansion} where the $\partial$ operator acts on the free solid edge, the operation is called $\mathsf{reshape}$ of the second type.         
			\item For graphs induced by the $\partial$-terms in \eqref{out_GG_expansion} where the $\partial$-operator acts on a light-weight factor in some $F_i$, the operation is called $\mathsf{insert}$.
			\item For graphs induced by the $\partial$-terms in \eqref{out_GG_expansion} where the $\partial$ operator acts on a $G$-edge in $\fS_1\cdots\fS_{q-1}$, the operation is called $\mathsf{bite}$. 
			\item For graphs induced by the $\partial$-terms in \eqref{out_GG_expansion} where the $\partial$-operator acts on an external $G$-edge in some $F_i$, the operation is called $\mathsf{extend}$.
			\item For the remaining graphs---namely, those induced by the growing terms in \eqref{out_GG_expansion} and those induced by the $\partial$-terms where the $\partial$-operator acts on a $G$-edge in the quasi-tail component of $\fS_q$---the operation is called $\mathsf{grow}$.
		\end{itemize}
	\end{definition}
	
	It is straightforward to verify that each of these operations increases the scaling order of the graph by at least $1$. Moreover, the $\mathsf{reshape}$ operations transform $\fS_q$ into a tail-head or body-head section, so that $\mathscr{S}_q$ becomes a snake. The $\mathsf{insert}$ operation converts a free tadpole diagram into a weight-head and inserts it into $\scr{S}_q$, producing a snake with $(q+1)$ effective sections. The $\mathsf{bite}$ operation transforms $\mathscr{S}_q$ into an ouroboros. The $\mathsf{extend}$ operation converts a free tadpole diagram into a quasi-body section and attaches it to $\scr{S}_q$, yielding a headless snake with $(q+1)$ effective sections. Finally, the $\mathsf{grow}$ operation preserves the headless-snake structure of $\mathscr{S}_q$ and does not modify the free tadpole diagrams.
	In summary, the operations $\mathsf{extend}$ and $\mathsf{grow}$ will also be referred to as \emph{intermediate operations}, as they produce headless snakes, whereas the remaining operations will be called \emph{terminal operations}, since they produce snakes or ouroboros graphs.

	We can now describe our expansion strategy. Starting from the initial graph $|f_{xy}\p{G}|^p=f_{x_1y_1}(G)\cdot f_{xy}(G)^{p/2-1}\ol{f_{xy}(G)}^{p/2}$ with the simplest headless snake $f_{x_1y_1}(G)$, we iteratively apply the $GG$-expansion to each headless snake arising in the expansion, as specified in \Cref{def_graph_operations}. The expansion is terminated if either (i) the resulting graph $\cG$ is a snake or an ouroboros, or (ii) the scaling order (recall Definition \ref{def_scaling_orer}) of $\cG$ is sufficiently large so that $\Psi^{\mathsf{ord}(\cG)}(0)\leq N^{-D}$. As discussed in \Cref{sec:graph_ideas}, this procedure terminates after $\OO(1)$ steps. This yields the following lemma.

	\begin{lemma}\label{lemma_terminate}
		Given any large constant $M>0$, we can expand $\big|f_{xy}(G)\big|^{p}$ into a linear combination of $\OO(1)$ graphs $\cal G_{\mu}$ in the sense of equal expectation:
		\be\label{eq:local_Gs}
		\big|f_{xy}(G)\big|^{p} =_{\E} \sum_{\mu}  \cal G_{\mu} + \mathcal Err.
		\ee
		Here, $\mathcal{E}rr$ denotes a sum of $\OO(1)$ graphs. Each $\cG_{\mu}$ and each graph in $\mathcal Err$ is of the form $\mathscr{S}_q\cdot F_{q+1}\cdots F_{p}$ for some $q\in\qq{p}$, where $F_i$, $i\in\qq{q+1,p}$, are free tadpole diagrams. For any graph in $\mathcal{E}rr$, $\mathscr{S}_q$ is a snake-like graph of scaling order $\ge M$, while, for each $\cG_\mu$, $\mathscr{S}_q=\fS_1\cdots\fS_q$ is either a snake or an ouroboros with $q$ \emph{effective sections}. 
	\end{lemma}

	\subsection{Upper bound of graphs}\label{sec_upper_bound_of_graphs}
	
	With Lemma \ref{lemma_terminate}, to establish \eqref{light_weight_term_moment_bound} it suffices to control the snake and ouroboros graphs $\mathscr{S}_q$. To this end, we have the following lemma.
	
	\begin{lemma}\label{snake_bound}
		Given a snake-like graph $\mathscr{S}_q=\fS_1\cdots\fS_q$ with $q$ effective sections, where the external vertices of $\fS_i$, $i\in\qq{q}$, are denoted by $x_i,y_i$ with $x_i\ne y_i$, we have the rough bound
		\begin{equation}\label{snake_bound_rough}
			\begin{aligned}
				\absa{\mathscr{S}_q}\prec N^{q}\qa{\Psi\p{0}}^{\mathsf{ord}\pa{\mathscr{S}_q}-q-1}.
			\end{aligned}
		\end{equation}
		Moreover, if $\mathscr{S}_q$ is a snake or an ouroboros, then the following improved bound holds:
		\begin{equation}\label{snake_bound_improved}
			\begin{aligned}
				\absa{\mathscr{S}_q}\prec \frac{1}{\eta^q}\qa{\Psi\p{0}}^{\mathsf{ord}\pa{\mathscr{S}_q}-\mathsf{ord}\pa{\txt{Aux}\pa{\mathscr{S}_q}}+n_{\txt{llw}}\pa{\mathscr{S}_q}}\prod_{i=1}^{q}\Psi\pa{|x_i-y_i|},
			\end{aligned}
		\end{equation}
		where the auxiliary graph $\txt{Aux}\pa{\mathscr{S}_q}$ and its scaling order $\mathsf{ord}\pa{\txt{Aux}\pa{\mathscr{S}_q}}$ are defined in \Cref{lemma_val_to_aux} below, and $n_{\txt{llw}}\pa{\mathscr{S}_q}$ denotes the number of leaf light-weights in the graph. In addition, if $\mathscr{S}_q$ is a snake or an ouroboros generated from the expansion in \Cref{lemma_terminate}, we have 
		\begin{equation}\label{eq:ordSq}
			\mathsf{ord}\pa{\mathscr{S}_q}-\mathsf{ord}\pa{\txt{Aux}\pa{\mathscr{S}_q}}+n_{\txt{llw}}\pa{\mathscr{S}_q}\ge q.
		\end{equation}
	\end{lemma}
	
	\begin{proof}[\bf Proof of \eqref{light_weight_term_moment_bound}]
		When $x=y$, the bound \eqref{light_weight_term_moment_bound} follows directly from the averaged local law \eqref{average_bound_2_to_1} and Ward's identity. Hence, we assume $x\neq y$ in the following. 
		Using Lemma \ref{lemma_terminate}, we expand $\E|f_{xy}(G)|^{p}$ as in \eqref{eq:local_Gs}. We choose $M$ sufficiently large such that $N^p\q{\Psi(0)}^{M-p-1}\le \eta^{-p}\q{\Psi(0)}^p\q{\Psi(N)}^p$. Then, by \eqref{snake_bound_rough}, the contribution $\E\mathcal Err$ can be bounded by the RHS of \eqref{light_weight_term_moment_bound}. 
		Next, for each $\cal G_\mu$ of the form $\mathscr{S}_q\cdot F_{q+1}\cdots F_{p}$, applying \eqref{snake_bound_improved} and \eqref{eq:ordSq} yields
		\begin{equation}\label{S_q_q_power}
			\E\absa{\cal G_\mu}\prec \eta^{-q}\qa{\Psi(0)}^{q}\qa{\Psi\pa{|x-y|}}^q\cdot \qa{\E|f_{xy}\pa{G}|^{p}}^{(p-q)/p}.
		\end{equation}
		Combining this bound with Hölder's inequality, we obtain
		\begin{equation*}
			\begin{aligned}
				\E|f_{xy}\pa{G}|^p\prec \sum_{q=1}^{p}\frac{1}{\eta^{q}}\qa{\Psi(0)}^{q}\qa{\Psi\pa{|x-y|}}^q\cdot\pa{\E|f_{xy}\pa{G}|^{p}}^{(p-q)/p}+ \eta^{-p}\q{\Psi(0)}^p\q{\Psi(|x-y|)}^p.
			\end{aligned}
		\end{equation*}
		Applying Young's inequality then yields \eqref{light_weight_term_moment_bound}. 
	\end{proof}

	For the proof of \Cref{snake_bound}, we first provide a heuristic explanation under the simplifying assumption that all solid edges in our graphs correspond to off-diagonal entries of $G$. Roughly speaking, the $\mathsf{ord}\pa{\mathscr{S}_q}+2n_{M}=n_S\pa{\mathscr{S}_q}+n_{\txt{lw}}\pa{\mathscr{S}_q}$ edges in $\mathscr{S}_q$ contribute the $\Psi(0)$ factors in \eqref{snake_bound_rough}, where $n_{M},n_S,n_{\txt{lw}}$ denote the number of internal molecules, solid edges, and light-weights, respectively. Moreover, recall from Figure \ref{fig_snake_mole} that $n_{M}\in\ha{q-1,q}$. The factor $N^q$ arises from the summations over at most $q$ internal molecules.\footnote{The loss of $\Psi(0)$-factors is due to the disappearance of solid edges resulting from possible identifications of different vertices connected by dotted edges.}
	For the improved bound \eqref{snake_bound_improved}, the quantity $\mathsf{ord}\pa{\txt{Aux}\pa{\mathscr{S}_q}}$ roughly captures the molecular graph structure (recall Definition \ref{mole_graph} and Figure \ref{fig_snake_mole}). In this auxiliary graph, the $q$ long solid edges lying on the $q$ paths connecting $x_i$ and $y_i$ yield the $q$ decaying factors $\Psi(|x_i-y_i|)$, while the application of Ward's identities to the summations over at most $q$ internal molecules produces the factor $\eta^{-q}$. The remaining $G$-edges within molecules contribute $\mathsf{ord}\pa{\mathscr{S}_q}-\mathsf{ord}\pa{\txt{Aux}\pa{\mathscr{S}_q}}$ factors of $\Psi(0)$. Finally, each leaf light-weight in $\mathscr{S}_q$ improves the bound by an additional factor $\Psi(0)$. This gain arises because the bound encoded by the scaling order uses the entrywise local law estimate $\Psi(0)$, whereas the actual estimate available from the averaged local law \eqref{average_bound_2_to_1} yields the stronger factor $\Psi^2(0)$.

	However, the actual proof is not so direct, since each solid edge may correspond either to a diagonal entry, in which case it is controlled by a dotted edge, or to an off-diagonal entry. To account for this, we assign to each pair of vertices $a,b$ connected by a solid edge the decomposition $1=\mathbf{1}_{a=b}+\mathbf{1}_{a\neq b}$. This induces the following decomposition of $\mathscr{S}_q$:
	\begin{equation}\label{decomposition_S_q}
		\mathscr{S}_q=\sum_{A\subseteq \mathsf{SP}\pa{\mathscr{S}_q}}\chi_{A}\cdot \mathscr{S}_q.
	\end{equation}
	Here, for any graph $\cG$, we denote by $\mathsf{SP}\p{\cG}:=\ha{\{a,b\}\subseteq V\pa{\cG}: a,b \txt{ are connected by a solid edge}}$ the set of unordered pairs of vertices connected by a solid edge. Moreover, we define
	\begin{equation*}
		\begin{aligned}
			\chi_A:=\prod_{\ha{a,b}\in A}\mathbf{1}_{a=b}\cdot \prod_{\ha{a,b}\in\mathsf{SP}\pa{\mathscr{S}_q}\setminus A}\mathbf{1}_{a\neq b}.
		\end{aligned}
	\end{equation*}
	To simplify the graph reduction procedure based on \eqref{eq:reduction_xi}, we introduce the notion of \emph{value graphs}, which provide upper bounds for the corresponding $G$-graphs.

	\begin{definition}[Value graphs]\label{def_value_graphs}
		Given an undirected graph consisting of external and internal vertices together with the edges between them, we assign the following structures to it and refer to the resulting object as a \emph{value graph}.
		\begin{itemize}
			
			\item{\bf  Solid edges:} We use $(x,y)$ to denote a solid edge between vertices $x$ and $y$. Every solid edge represents an entry of a $\Psi$-matrix. Here, a $\Psi$-matrix is defined as an $N\times N$ (random)  symmetric matrix $\xi$ with nonnegative entries satisfying, for all $x,y\in \ZL$,
			\begin{equation}\label{xi_matrix_upper_bound}
				{\xi_{xy}}\prec \Psi(|x-y|),\quad \sum_{a\in\ZL}\xi_{xa}^2= \sum_{a\in\ZL}\xi_{ax}^2 \prec \eta^{-1}. \end{equation}
			Different solid edges may correspond to $\xi$-factors arising from different $\Psi$-matrices.

            \item {\bf Dotted edges:} A dotted edge between $x$ and $y$ represents the factor $\mathbf{1}_{x=y}$, while a $\times$-dotted edge represents the factor $\mathbf{1}_{x\neq y}$.

			\item {\bf Waved edges:} A waved edge between $x$ and $y$ represents a factor $\cS_{xy}$.
			
			\item {\bf Weight factors:} Each graph is associated a factor that is a power of $\Psi(0)$.
			
			\item {\bf Dashed edges:} For convenience, we occasionally introduce dashed edges in the graph, each representing the trivial factor $1$.
		\end{itemize}
		The relevant notions---such as scaling orders and molecular value graphs---are defined analogously to those for $G$-graphs, where, besides solid and dotted edges, dashed edges are also retained in the molecular graphs, and the scaling orders are defined as $\#\ha{\txt{solid edges}}-2\#\ha{\txt{internal molecules}}$ since value graphs have no light-weights.
	\end{definition}
	
The primary example of a $\Psi$-matrix is  
	\begin{equation}\label{def_xi_edges}    \xi_{xy}^2:=\cS_{xy}+\cL_{xy}+\cL_{yx}=\cS_{xy}+\sum_{a,b\in\ZL}S_{xa}(|G_{ab}|^2+|G_{ba}|^2)S_{b y}, \quad \forall x,y\in\ZL.
	\end{equation}
	Moreover, given any $\Psi$-matrix $\xi$, the admissibility of the $\Psi$-parameter (see Definition \ref{assm_admissible}) implies that the following matrix is again a $\Psi$-matrix (recall \eqref{eq:reduction_xi}):
	\begin{equation}\label{xi_matrix_convolution}
		\pb{\xi'_{xy}}^2:= \sum_{a\in\ZL}\pa{\xi_{xa}^2\cS_{ay}+\xi_{ya}^2\cS_{ax}},\qquad  \forall x,y\in\ZL.
	\end{equation}
	For any graph $\chi_A\cdot \mathscr{S}_q$ appearing in the decomposition \eqref{decomposition_S_q}, we may apply the local laws \eqref{entrywise_bound_2_to_1} and \eqref{average_bound_2_to_1} to bound it by its associated value graph, as stated in the following lemma.

	\begin{lemma}\label{lemma_graphs_to_value_graphs}
	Given a term $\chi_A\cdot \mathscr{S}_q$ from \eqref{decomposition_S_q}, we define its associated value graph $\txt{Val}_A\pa{\mathscr{S}_q}$ through the following procedure.
		\begin{enumerate}
			\item[(1)] Starting from $\chi_{A}\cdot \mathscr{S}_q$, remove all light-weights and all waved edges attached to leaf light-weights (together with the vertices incident to leaf light-weights), and add a factor \smash{$\q{\Psi(0)}^{n_{\txt{lw}}\pa{\mathscr{S}_q}+n_{\txt{llw}}\pa{\mathscr{S}_q}}$} to the graph, where $n_{\txt{lw}}$ and $n_{\txt{llw}}$ denote the number of light-weights and leaf light-weights, respectively.

			\item[(2)] In the resulting graph from Step (1), consider any pair of vertices $a$ and $b$ connected by a dotted edge within the same molecule. Merge these two vertices and remove all $G$-solid edges between them. The subgraph formed by the waved edges and vertices within this molecule then contains a loop (which may be a self-loop at a vertex). Eliminate one waved edge within this molecule to restore the tree structure of the waved edges, and add a weight factor $\Psi^2(0)\ge W^{-1}$ to the graph. Repeat this step until no such pair of vertices remains.
			
			\item[(3)] In the resulting graph from the previous steps, replace each $G$-solid edge whose endpoints are connected by $\times$-dotted edges with the corresponding $\xi$-solid edge defined in \eqref{def_xi_edges}. Moreover, replace each $G$-solid edge whose endpoints are connected by a dotted edge with a dashed edge. (Note that, after Step (2), such $G$-solid edges must connect vertices belonging to different molecules.) Keep all dotted and $\times$-dotted edges. Replace each waved edge by an $\cS$-waved edge. 
		\end{enumerate}
		The value graph $\txt{Val}_A\pa{\mathscr{S}_q}$ obtained through this procedure provides an upper bound for the graph $\chi_A\cdot \mathscr{S}_q$:
		\begin{equation}\label{graphs_to_value_graphs}
			\begin{aligned}
				\absa{\chi_A\cdot \mathscr{S}_q}\prec \txt{Val}_A\pa{\mathscr{S}_q}.
			\end{aligned}
		\end{equation}
		Moreover, $\txt{Val}_A\pa{\mathscr{S}_q}$ inherits the following structural properties from $\mathscr{S}_q$.
		\begin{enumerate}
			\item The molecular graph of $\txt{Val}_A\pa{\mathscr{S}_q}$ is graphically isomorphic to the molecular graph of $\mathscr{S}_q$. \item The waved edges within each molecule form a tree.
			\item Each vertex is incident to at most two solid or dashed edges.
			
			\item The following inequality holds: 
			\begin{equation}
				n_{\Psi}\pa{\txt{Val}_A\pa{\mathscr{S}_q}}+n_{\txt{dash}}\pa{\txt{Val}_A\pa{\mathscr{S}_q}}+\mathsf{ord}\pa{\txt{Val}_A\pa{\mathscr{S}_q}} \geq \mathsf{ord}\pa{\mathscr{S}_q}+n_{\txt{llw}}\pa{\mathscr{S}_q} ,\label{eq:relation_ord}    
			\end{equation} 
			where $n_{\Psi}$ and $n_{\txt{dash}}$ denote the number of $\Psi(0)$-factors and dashed edges, respectively, in $\txt{Val}_A\pa{\mathscr{S}_q}$.
		\end{enumerate}
	\end{lemma}
	
	\begin{proof}
		The upper bound \eqref{graphs_to_value_graphs} follows immediately from the local laws \eqref{upper_bound_2_to_1} and \eqref{average_bound_2_to_1}, together with the fact that merging two vertices results in a loss of at most two $G$-solid edges, since every vertex in $\mathscr{S}_q$ has degree at most 2. This loss is compensated by the gain of a weight factor resulting from the removal of a waved edge, as ensured by the bounds \eqref{alpha_decay}, \eqref{graph:Spm}, and the assumption $\Psi^2(0)\geq W^{-1}$ from Definition \ref{assm_admissible}.
        Properties (i) and (ii) follow directly from the definition. For property (iii), note that the original graph $\mathscr{S}_q$ satisfies this property, and it is straightforward to verify that it is preserved throughout each step in the construction of $\txt{Val}_A\pa{\mathscr{S}_q}$.
		It remains to prove property (iv) by tracking the three steps in the construction of $\txt{Val}_A\pa{\mathscr{S}_q}$. Starting from the graph $\chi_A\cdot \mathscr{S}_q$, which has the same scaling order as $\mathscr{S}_q$ by definition, we observe that Step (1) preserves the quantity $n_{\Psi}+\mathsf{ord}+n_{\txt{llw}}$, which is exactly the RHS of \eqref{eq:relation_ord} for $\chi_{A}\cdot\mathscr{S}_q$. 
		Next, suppose that Step (2) is iterated $k$ times. In this case, at most $2k$ solid edges are removed, while a factor $\Psi^{2k}(0)$ is added to the graph. Hence, the relation \eqref{eq:relation_ord} remains valid. Finally, Step (3) in the construction of $\txt{Val}_A\pa{\mathscr{S}_q}$ preserves the quantity $n_{\txt{dash}}+\mathsf{ord}$. Combining the above observations, we obtain property (iv) for $\txt{Val}_A\pa{\mathscr{S}_q}$.
	\end{proof}

	With \eqref{graphs_to_value_graphs} in hand, it remains to bound the associated value graph $\txt{Val}_A\pa{\mathscr{S}_q}$. Following the strategy in \cite{dubova2025delocalizationnonmeanfieldrandommatrices}, our aim is to bound it in terms of its molecular graph. To this end, we employ CS inequalities such as \eqref{eq:reduction_xi} to reduce the number of vertices within each molecule to one.
	This approach works well in the case $A=\emptyset$, where one can exploit the fact that only a limited number of solid edges are incident to each molecule in the molecular graph, due to the snake-like structure of $\mathscr{S}_q$. However, when $A\neq \emptyset$, the situation becomes more delicate. If one were to merge all vertices connected by dotted edges, as in \cite{dubova2025delocalizationnonmeanfieldrandommatrices}, a single vertex might end up incident to more than two edges, which prevents the direct application of \eqref{eq:reduction_xi}. In such cases, we need to exploit the fact that the number of internal vertices is reduced, thereby producing additional waved edges whose decay can be utilized. For this purpose, we adopt a slightly less natural but technically convenient device: rather than merging vertices connected by dotted edges, we apply the following inequalities to convert a dotted edge between molecules into a $\xi$-type edge: 
	\begin{equation}\label{eq:convert-delta}
		\sum_{a\in\ZL}\delta_{xa}\cS_{ay}\lesssim (\cS_{xy})^{1/2},\qquad \sum_{a\in\ZL}\delta_{xa}\xi_{x'a}\cS_{ay}\leq \pB{\sum_{a\in\ZL}\delta_{xa}\cS_{ay}}^{1/2}\pB{\sum_{a\in\ZL}\xi_{x'a}^2\cS_{ay}}^{1/2}.
	\end{equation}
	Although these inequalities are not sharp, they have the advantage that they preserve the graphical structure, unlike the vertex-merging operation, which simplifies the subsequent analysis. With this preparation, we are ready to prove the following reduction lemma, which bounds value graphs in terms of their \emph{auxiliary graphs}.
	
	\begin{lemma}\label{lemma_val_to_aux}
		In the setting of \Cref{snake_bound}, for any $A\subseteq \mathsf{SP}\pa{\mathscr{S}_q}$, the associated value graph $\txt{Val}_A\pa{\mathscr{S}_q}$ can be bounded by its auxiliary graph, denoted by $\txt{Aux}\pa{\mathscr{S}_q}$, as follows:
		\begin{equation}
			\begin{aligned}
				\txt{Val}_A\pa{\mathscr{S}_q}\prec \qa{\Psi(0)}^{n_{\Psi}\pa{\txt{Val}_A\pa{\mathscr{S}_q}}+n_{\txt{dash}}\pa{\txt{Val}_A\pa{\mathscr{S}_q}}+\mathsf{ord}\pa{\txt{Val}_A\pa{\mathscr{S}_q}}-\mathsf{ord}\pa{\txt{Aux}\pa{\mathscr{S}_q}}}\cdot\txt{Aux}\pa{\mathscr{S}_q}.
			\end{aligned}
		\end{equation}
		Here, the auxiliary graph of $\mathscr{S}_q$ and its scaling order are defined as follows.
		\begin{enumerate}
			\item             
            If \(\mathscr{S}_q\) is a headless snake, a snake, or an ouroboros of type I (see Definition \ref{def_graph_snakes}), then the auxiliary graph \(\txt{Aux}(\mathscr{S}_q)\) is defined as the molecular graph associated with \(\mathscr{S}_q\), with each solid or dashed edge between molecular vertices replaced by a \smash{\(\wt\xi\)}-edge. Here, a \smash{\(\wt\xi\)}-edge between \(a\) and \(b\) represents the factor \smash{\(\wt\xi_{ab}=\xi_{ab}+\delta_{ab}\)} for some \(\Psi\)-matrix \(\xi\). 
            The dotted and $\times$-dotted edges are discarded. External vertices are fixed as the corresponding $x_i$ or $y_i$, while internal vertices are summed over $\ZL$. In this case, we define the scaling order of $\txt{Aux}\pa{\mathscr{S}_q}$ as 
			\[\mathsf{ord}\pa{\txt{Aux}\pa{\mathscr{S}_q}}:=\#\h{\wt\xi\txt{-edges in }\txt{Aux}\pa{\mathscr{S}_q}}-2\#\h{\txt{internal vertices in }\txt{Aux}\pa{\mathscr{S}_q}}.\]
			
			\item If $\mathscr{S}_q$ is an ouroboros of type II, the auxiliary graph $\txt{Aux}\pa{\mathscr{S}_q}$ is defined as in case (i), except that one of the two distinguished edges is removed from the graph.\footnote{We remove one of the edges to ensure that every vertex in the auxiliary graph has degree at most five. For example, in panel (f) of \Cref{fig_snake_mole}, if both distinguished edges are connected to \(a_1\), then \(a_1\) has degree six.} In this case, we define
			\[\mathsf{ord}\pa{\txt{Aux}\pa{\mathscr{S}_q}}:=\#\h{\wt\xi\txt{-edges in }\txt{Aux}\pa{\mathscr{S}_q}}-2\#\h{\txt{internal vertices in }\txt{Aux}\pa{\mathscr{S}_q}}+1,\]
			where the $+1$ accounts for the loss of a distinguished edge in the construction.
		\end{enumerate}
	\end{lemma}
	\begin{proof}
		The proof relies on the following claim.
		
		\begin{claim}\label{shrink_to_star}
			Consider a graph $\cG$ with external vertices $x_1,x_2,\ldots,x_r$ and internal vertices $y_1,y_2,\ldots,y_k$, in which the waved edges form a tree on the vertices $y_1,y_2,\ldots,y_k$. Let the remaining edges of $\cG$ be denoted by $\ha{e_{i}}_{i=1}^r$. Assume that each $e_{i}$ is either a $\xi$-solid edge or a dotted edge connecting $x_i$ to some internal vertex. If $r\leq 5$, then there exist some $\Psi$-matrices $\h{\xi^{(i)}}_{i=1}^r$ such that \begin{equation}
				\cG\prec \sum_{y\in\ZL}\prod_{i=1}^{r}\wt\xi_{x_i y}^{(i)},\quad \text{with}\quad \wt \xi^{(i)}_{ab}:=\delta_{ab}+\xi^{(i)}_{ab}. \end{equation}
			Graphically, this inequality shows that $\cG$ can be bounded by a star graph obtained by collapsing all internal vertices of $\cG$ into a single internal vertex $y$, and replacing all external edges with \smash{$\wt\xi$}-edges. 
		\end{claim}

		With \Cref{shrink_to_star}, \Cref{lemma_val_to_aux} follows readily by bounding the solid edges within each molecule by $\Psi(0)$-factors. This procedure generates		$\mathsf{ord}\pa{\txt{Val}_A\pa{\mathscr{S}_q}}+n_{\txt{dash}}\pa{\txt{Val}_A\pa{\mathscr{S}_q}}-\mathsf{ord}\pa{\txt{Aux}\pa{\mathscr{S}_q}}$ additional $\Psi(0)$-factors. Next, we apply \Cref{shrink_to_star} molecule by molecule to collapse each internal molecule into a single vertex. Note that the degree of each vertex in $\txt{Aux}\pa{\mathscr{S}_q}$ is at most five, and this shrinking procedure does not affect the scaling order of the graph.

		It remains to prove \Cref{shrink_to_star}. We argue by induction on $k\ge 1$. The case $k=1$ is trivial. Assume that $k\ge 2$ and that \Cref{shrink_to_star} holds for $1,\ldots,k-1$. Using the assumption $r\le 5$ and the graphical fact that every tree with at least two vertices has at least two leaves, we deduce that there exists a leaf internal vertex $y_i$ that is incident to exactly one waved edge and at most two $e_i$-edges. Without loss of generality, we may assume that this vertex is $y_1$. Then $y_1$ is connected either to (i) no external vertex, (ii) one external vertex $x_1$, or (iii) two external vertices $x_1$ and $x_2$, while the other endpoint of the waved edge incident to $y_1$ is $y_2$. 
		In case (i), summing over $y_1$ directly reduces the number of internal vertices by one. In case (ii), applying the CS inequality yields
		\begin{align}\label{eq_e_to_e_example}
			\sum_{y_1}e_{x_1y_1}\cS_{y_1y_2}\le \pbb{\sum_{y_1}e_{x_1y_1}^2\cS_{y_1y_2}}^{1/2}=:e'_{x_1y_2}.
		\end{align}
		By \eqref{xi_matrix_convolution} and \eqref{eq:convert-delta}, the symmetrized $e'_{x_1y_2}+e'_{y_2x_1}$ represents a $\xi$-solid edge between $x_1$ and $y_2$, and the internal vertex $y_1$ is removed from the resulting graph. Case (iii) can be handled similarly. Applying the CS inequality and symmetrization again, we obtain two new $\xi$-solid edges while reducing the number of internal vertices by one:
		\begin{align*}
			\sum_{y_1}e_{x_1y_1}e_{x_2y_1}\cS_{y_1y_2}\leq \pbb{ \sum_{y_1}e_{x_1y_1}^2\cS_{y_1y_2}}^{1/2}\pbb{\sum_{y_1}e_{x_2y_1}^2\cS_{y_1y_2}}^{1/2}.
		\end{align*}
		This completes the induction and hence the proof of \Cref{shrink_to_star}.
	\end{proof}

	Now, to prove \Cref{snake_bound}, it remains to bound the auxiliary graph. To this end, we first state the following lemma, which describes the simple sectional structure of the auxiliary graph.
	
	\begin{lemma}[Sectional structure]\label{lemma_section_structure_aux}
    In the setting of \Cref{snake_bound}, the auxiliary graphs consist of the following three types of sections, which are induced from the original sections defined in \Cref{def_graph_types}.

		\begin{enumerate}
			\item A \emph{tail} section (corresponding to the quasi-tail and tail sections defined in Definition \ref{def_graph_types}) with external vertices $x,y$ and inner vertex $a$ represents the graph \smash{$\wt\xi_{xa}\wt\xi_{a y}$}.
			
			\item A \emph{body} section (corresponding to the quasi-body and body sections defined in Definition \ref{def_graph_types}) with external vertices $x,y$, inner vertex $a$, and outer vertex $b$ represents the graph \smash{$\wt\xi_{xa}\wt\xi_{ab}\wt\xi_{b y}$} or  \smash{$\wt\xi_{ya}\wt\xi_{ab}\wt\xi_{b x}$}.
			
			\item The three types of \emph{head} sections take the following forms.
			\begin{itemize}
				\item Tail-head with external vertices $x,y$: \smash{$\wt\xi_{xy}$}.
				\item Body-head I with external vertices $x,y$ and outer vertex $b$: \smash{$\wt\xi_{xb}\wt\xi_{b y}$}.
				\item Body-head II (corresponding to the body-head II and weight-head sections defined in Definition \ref{def_graph_types}) with external vertices $x,y$, free vertex $a$, and outer vertex $b$: \smash{$\wt\xi_{xa}\wt\xi_{ab}^2\wt\xi_{a y}$}.
			\end{itemize}
		\end{enumerate}
		Here, the notation $\wt\xi$ does not refer to any specific $\wt\xi$-matrix (defined as $I_{N}+\xi$ for some $\Psi$-matrix $\xi$), but merely indicates that the corresponding edge is a $\wt\xi$-edge. Using these notations, we can classify the structure of $\txt{Aux}\pa{\mathscr{S}_q}$ as follows; see Figure \ref{fig_snake_mole} for an illustration.
		\begin{itemize}
			\item[(a)] If $\mathscr{S}_q$ is a headless snake, then we have the section decomposition $\txt{Aux}\pa{\mathscr{S}_q}=\wt\fS_1\wt\fS_2\cdots\wt\fS_q$, where each \smash{$\wt\fS_i$} is a section with external vertices $x_i,y_i$. The section \smash{$\wt\fS_1$} is a tail with inner vertex $a_1$, while $\wt\fS_i$ is a body with inner vertex $a_i$ and outer vertex $b_i$ for $i\in\qq{2,q}$. Moreover, the consistency condition $a_i=b_{i+1}$ holds for $i\in\qq{q-1}$.
			
			\item[(b)] If $\mathscr{S}_q$ is a snake, then $\txt{Aux}\pa{\mathscr{S}_q}=\wt\fS_1\wt\fS_2\cdots\wt\fS_q$, where each $\wt\fS_i$ is a section with external vertices $x_i,y_i$. If $q=1$, then $\fS_1$ is a tail-head. If $q\ge2$, then \smash{$\wt\fS_1$} is a tail with inner vertex $a_1$, \smash{$\wt\fS_q$} is either a body-head I or a body-head II with outer vertex $b_q$, and each intermediate section $\fS_i$, $i\in\qq{2,q-1}$, is a body with inner vertex $a_i$ and outer vertex $b_i$. Moreover, the consistency condition $a_i=b_{i+1}$ holds for $i\in\qq{q-1}$.

			\item[(c)] If $\mathscr{S}_q$ is an ouroboros of type I (recall Definition \ref{def_graph_snakes}), then $\txt{Aux}(\mathscr{S}_q)$ is obtained from a headless snake \smash{$\wt\fS_1\wt\fS_2\cdots\wt\fS_q$} of the form described in (a) by pulling a \smash{$\smash{\wt\xi}$}-edge from \smash{$\wt\fS_1\cdots\wt\fS_{q-1}$} to $a_q$. If $\mathscr{S}_q$ is an ouroboros of type II, then $\txt{Aux}(\mathscr{S}_q)$ is obtained from a headless snake by connecting $a_q$ to some $a_i$, $i\in\qq{q-1}$, by a \smash{$\smash{\wt\xi}$}-edge. In both cases, the consistency condition $a_i=b_{i+1}$ holds for $i\in\qq{q-1}$.    
			
		\end{itemize}
		From the above construction, we can directly verify the following scaling orders:
		\begin{itemize}
			\item $\mathsf{ord}\pa{\txt{Aux}\pa{\mathscr{S}_q}}=q-1$ if $\mathscr{S}_q$ is a headless snake;
			\item $\mathsf{ord}\pa{\txt{Aux}\pa{\mathscr{S}_q}}=q+1$ if $\mathscr{S}_q$ is an ouroboros of type II;
            \item $\mathsf{ord}\pa{\txt{Aux}\pa{\mathscr{S}_q}}=q$ in all other cases.
		\end{itemize}
	\end{lemma}
	\begin{proof}
		The lemma follows directly from \Cref{def_graph_snakes} together with the definitions of the value and auxiliary graphs in Definition \ref{def_value_graphs} and Lemma \ref{lemma_val_to_aux}. We omit the details.
	\end{proof}

	With the sectional structure described in Lemma \ref{lemma_section_structure_aux}, we are ready to provide the last piece of the proof for Lemma \ref{snake_bound}.
	
	\begin{lemma}\label{aux_bound}
		In the setting of \Cref{snake_bound}, the auxiliary graph constructed in Lemma \ref{lemma_val_to_aux} satisfies the bound
		\begin{equation}\label{aux_bound_rough}
			\begin{aligned}
				\txt{Aux}\pa{\mathscr{S}_q}\prec N^q.
			\end{aligned}
		\end{equation}
		Moreover, if $\mathscr{S}_q$ is a snake or an ouroboros, then the following improved bound holds:
		\begin{equation}\label{aux_bound_improved}
			\txt{Aux}\pa{\mathscr{S}_q}\prec\frac{1}{\eta^q}\prod_{i=1}^q\Psi\pa{|x_i-y_i|}.
		\end{equation}
	\end{lemma}
	\begin{proof}  
		The bound \eqref{aux_bound_rough} is immediate: we bound each solid edge by $1$ and sum over the at most $q$ internal vertices. For the bound \eqref{aux_bound_improved}, the spatial decay factors therein arise from the $q$ paths $\ha{\mathfrak{P}_i}_{i=1}^{q}$ consisting of \smash{$\wt \xi$}-edges within the graph $\txt{Aux}(\mathscr{S}_q)$, where the path $\mathfrak{P}_i$ connects $x_i$ and $y_i$. To see the existence of these paths, let \smash{$\wt \fS_1\wt \fS_2\cdots \wt \fS_q$} be the sectional representation given in \Cref{lemma_section_structure_aux}. If \(\mathscr{S}_q\) is a snake or an ouroboros of type II, then each section \smash{\(\wt \fS_i\)} yields the required path \(\mathfrak{P}_i\) for \(i\in\qq{q}\) (if a section contains more than one such path, we choose the shortest one); see \Cref{fig_snake_mole}. Otherwise, if \(\mathscr{S}_q\) is an ouroboros of type I, then an edge in \smash{$\wt \fS_1\wt \fS_2\cdots \wt \fS_{q-1}$} is pulled, but the path is clearly retained.

        Next, we explain the approach to extract the spatial decay factor. The case $q=1$ is trivial since $x_1\ne y_1$. For $q\geq 2$, the cases where $\mathscr{S}_q$ is an ouroboros of type I or a snake with $\fS_q$ being a body-head I were already treated in Lemma 7.23 of \cite{dubova2025delocalizationnonmeanfieldrandommatrices}. 
        More precisely, we expand the auxiliary graph $\txt{Aux}(\mathscr{S}_q)$ into a sum of $\OO(1)$ terms by writing each \smash{$\wt\xi$}-edge as \smash{$\wt\xi_{ab}=\xi_{ab}+\delta_{ab}$}. Each graph generated by this expansion is obtained by replacing every \smash{$\wt\xi$}-edge in $\txt{Aux}(\mathscr{S}_q)$ with either a $\xi$-edge or a dotted edge. Let $B$ denote the set of edges replaced by dotted edges, and denote the resulting graph by $\txt{Aux}_B(\mathscr{S}_q)$. Since $x_i\neq y_i$ for all $i\in\qq{q}$, the contribution $\txt{Aux}_B(\mathscr{S}_q)$ vanishes whenever there exists a path of dotted edges connecting $x_i$ to $y_i$ for some \(i\). In particular, each section $\fS_i$ must contain at least one $\xi$-edge. We then merge each pair of vertices connected by a dotted edge in \(\txt{Aux}_B(\mathscr{S}_q)\). In each resulting graph, there exist \(q\) edge-disjoint \(\xi\)-paths connecting \(x_i\) and \(y_i\), one for each \(i\in\qq{q}\).
Moreover, the resulting graph inherits the property that each internal vertex is traversed by at least two distinct paths inherited from the original graph; see \Cref{fig_snake_mole}. This verifies the key path properties required by \cite[Lemma 7.23]{dubova2025delocalizationnonmeanfieldrandommatrices}. The other assumptions therein can be verified similarly. The bound \eqref{aux_bound_improved} then follows from that result. It suffices to prove \eqref{aux_bound_improved} when $\mathscr{S}_q$ is an ouroboros of type II or a snake with $\fS_q$ being a body-head II or weight-head.

        To this end, we prove the required bound for a slightly more general class of graphs. More precisely, we consider the graph $\txt{Aux}^{g}(\mathscr{S}_q)$ obtained from $\txt{Aux}(\mathscr{S}_q)$ by replacing an arbitrary \smash{$\wt \xi$}-edge in each path $\fP_i$ with a \emph{ghost edge},\footnote{We remark that when $\fS_q$ is a body-head II or weight-head, the two edges connecting $a_{q-1}$ and $a_q$ (where \(a_q\) denotes the central vertex of the body-head II or weight-head section) are not included in any path.} representing the trivial factor $1$, for every $i\in\qq{q}$. We refer to $\txt{Aux}^g(\mathscr{S}_q)$ as a \emph{ghostification} of $\txt{Aux}(\mathscr{S}_q)$. 
		We claim that all such graphs satisfy  
		\begin{equation}\label{remaining_graph_bound}
			\txt{Aux}^g(\mathscr{S}_q)\prec \eta^{-\#\{\txt{internal molecules of }\mathscr{S}_q\}}.
		\end{equation} 
		Observe that this estimate implies \eqref{aux_bound_improved}. Indeed, in each path $\fP_i$ we bound the longest $\wt\xi$-edge along a path from $x_i$ to $y_i$, say the edge $(a_i,b_i)$, by the corresponding decay factor $\Psi(|a_i-b_i|)+\delta_{a_ib_i}$. By the triangle inequality, we have $|a_i-b_i|\gtrsim |x_i-y_i|>0$. Hence, $\delta_{a_ib_i}=0$ and any graph $\txt{Aux}(\mathscr{S}_q)$ can be bounded by $\prod_{i=1}^q\Psi\pa{|x_i-y_i|}$ times a ghostification graph $\txt{Aux}^g(\mathscr{S}_q)$, which yields the remaining factor $\OO_{\prec}(\eta^{-q})$ by \eqref{remaining_graph_bound}, since the number of internal molecules is at most $q$. 

        It remains to prove \eqref{remaining_graph_bound}. We first consider the case where $\mathscr{S}_q$ is a snake with body-head II or weight-head. By the consistency condition $a_i=b_{i+1}$ for $\qq{q-1}$, we can represent the graph $\txt{Aux}^g(\mathscr{S}_q)$ as
        \begin{equation*}
            \txt{Aux}^g(\mathscr{S}_q)=\bu^*K_2\cdots K_{q-1}D_q\bv.
        \end{equation*}
        Here, we denote the edges of $\txt{Aux}^g(\mathscr{S}_q)$ by $e$, and define the vectors $\bu,\bv\in\C^{\ZL}$, the matrices $K_i\in\C^{\ZL\times\ZL}$, $i\in\qq{2,q-1}$, and the matrix $D_q$ by
        \begin{equation}\label{def_bu_bv_K_i}
        \begin{aligned}            
                \bu(a_1)=e_{x_1a_1}e_{a_1y_1},\quad & \bv(a_{q})=e_{x_qa_{q}}e_{a_q y_q},\\ 
                K_{i}(a_{i-1},a_i)=e_{x_ia_{i-1}}e_{a_{i-1}a_{i}}e_{a_i y_i},\quad & D_{q}(a_{q-1},a_{q})=e_{a_{q-1}a_{q}}e_{a_{q}a_{q-1}},
        \end{aligned}
        \end{equation}
        where, with a slight abuse of notation, the symbol \(e\) denotes both an edge and the value associated with that edge. 
        By the definition of ghostification and the bound \smash{$\sum_{a}(\wt \xi_{xa}^2+\wt \xi_{ax}^2)\prec \eta^{-1}$}, we obtain
        \begin{equation}\label{eq:KHS}
                \|\bu\|_{2}\prec \eta^{-1/2},\quad \|\bv\|_{2}\prec \eta^{-1/2},\quad \|K_i\|_{2\to 2}\leq \|K_i\|_{\HS}\prec \eta^{-1}.
        \end{equation}
        Moreover, by Schur's test,
        $$\|D_q\|_{2\to 2}\leq \sqrt{\|D_q\|_{\infty\to \infty}\|D_q\|_{1\to 1}}=\pa{\max_{a}\sum_{b}|D_q(a,b)|}^{1/2}\pa{\max_{b}\sum_{a}|D_q(a,b)|}^{1/2}\prec \eta^{-1}.$$
        Therefore, we obtain
        \begin{equation*}
            \begin{aligned}
                \txt{Aux}^g(\mathscr{S}_q) \leq \|\bu\|_2\|\bv\|_2\|D_q\|_{2\to 2}\prod_{i=2}^{q-1}\|K_i\|_{2\to 2}\prec \eta^{-q}.
            \end{aligned}
        \end{equation*}

        For the case where $\mathscr{S}_q$ is an ouroboros of type II, suppose that the distinguished edge is $e_{a_ia_q}$ for some $i\in\qq{q-1}$. We define $\bu,K_2,\ldots,K_{q}$ as in \eqref{def_bu_bv_K_i} and set
        \begin{equation*}
            \bv(a_i):=(\bu^*K_2\cdots K_{i})(a_i),\quad \mathfrak{E}_{a_qa_i}:=e_{a_ia_q}.
        \end{equation*}
        With these notations, we can write \(           \txt{Aux}^g(\mathscr{S}_q)=\avga{\mathfrak{D}_{\bv}K_{i+1}\cdots K_{q}\mathfrak{E}},\)
where $\mathfrak{D}_{\bv}$ is the diagonal matrix defined by $\mathfrak{D}_{\bv}(a,b):=\delta_{ab}\bv(a)$. Then, by the CS inequality, we obtain
        \begin{equation*}
            \begin{aligned}
                \txt{Aux}^g(\mathscr{S}_q)\leq \|K_{i+1}\cdots K_{q}\|_{\HS}\|\mathfrak{E}\mathfrak{D}_{\bv}\|_{\HS} \leq \prod_{i+1\leq j\leq q}\|K_j\|_{\HS}\cdot \|\mathfrak{E}\mathfrak{D}_{\bv}\|_{\HS}\prec \eta^{-q},
            \end{aligned}
        \end{equation*}
        where we used \eqref{eq:KHS} and the bound
        $$ \|\mathfrak{E}\mathfrak{D}_{\bv}\|_{\HS}^2=\sum_{a,b}|\mathfrak{E}_{ab}|^2|\bv(b)|^2\prec \eta^{-1}\|\bv\|_2^2\leq \eta^{-1}\|\bu\|_2^2\prod_{j=2}^i\|K_j\|_{2\to 2}^2\prec \eta^{-2i}.
        $$
        This proves \eqref{remaining_graph_bound} and completes the proof of
\eqref{aux_bound_improved}.
  	\end{proof}

	We now collect the technical ingredients developed above to complete the proof of \Cref{snake_bound}.

	\begin{proof}[\bf Proof of \Cref{snake_bound}]
	Combining the decomposition \eqref{decomposition_S_q}, the bound \eqref{graphs_to_value_graphs}, and Lemmas \ref{lemma_val_to_aux} and \ref{aux_bound}, we immediately obtain the bounds \eqref{snake_bound_rough} and \eqref{snake_bound_improved}. It remains to prove \eqref{eq:ordSq}. From the expansion strategy, we know that $\mathscr{S}_q$ is generated by at least $(q-1)$ operations if $\mathscr{S}_q$ is a snake with a weight-head, and by at least $q$ operations otherwise. We consider the following cases.
		
		First, consider a snake with a weight-head. Suppose that $\mathscr{S}_q$ is generated by at least $q$ operations from Definition \ref{def_graph_operations}. Note that the original graph $\prod_{i=1}^{q}f_{x_iy_i}(G)$ has scaling order $q$, while each operation increases the scaling order by at least $1$. Hence $\mathsf{ord}\pa{\mathscr{S}_q}\ge 2q$. Together with the fact that $\mathsf{ord}\pa{\txt{Aux}\pa{\mathscr{S}_q}}=q$ by \Cref{lemma_section_structure_aux}, this implies \eqref{eq:ordSq}.
		Otherwise, if $\mathscr{S}_q$ is generated by exactly $(q-1)$ operations, then $\mathsf{ord}\pa{\mathscr{S}_q}-\mathsf{ord}\pa{\txt{Aux}\pa{\mathscr{S}_q}}\ge q-1$. In this case, in order to obtain a snake with $q$ sections, the final operation must be an $\mathsf{insert}$, while all preceding operations must be $\mathsf{extend}$. This implies that the leaf light-weight in $\fS_1$ is never affected during the expansion procedure, which yields $n_{\txt{llw}}\pa{\mathscr{S}_q}\ge 1$ and hence establishes \eqref{eq:ordSq}.
		
		In all other cases, the $q$ operations increase the scaling order by $q$, which implies $\mathsf{ord}\pa{\mathscr{S}_q}\ge 2q$. If $\mathscr{S}_q$ is not an ouroboros of type II, then $\mathsf{ord}\pa{\txt{Aux}\pa{\mathscr{S}_q}}=q$ by Lemma \ref{lemma_section_structure_aux}, which again yields \eqref{eq:ordSq}. 
		Finally, it remains to consider the case where $\mathscr{S}_q$ is an ouroboros of type II. In this case, Lemma \ref{lemma_section_structure_aux} implies that $\mathsf{ord}\pa{\txt{Aux}\pa{\mathscr{S}_q}}=q+1$. This already establishes \eqref{eq:ordSq} if $\mathsf{ord}\pa{\mathscr{S}_q}\ge 2q+1$. In the case $\mathsf{ord}\pa{\mathscr{S}_q}=2q$, a similar argument to above shows that at least one of $\fS_1$ and $\fS_q$ contains a leaf light-weight. Consequently $n_{\txt{llw}}\pa{\mathscr{S}_q}\ge 1$, which again implies \eqref{eq:ordSq}.
	\end{proof}

	\subsection{Proofs of \eqref{bound_pull_in_edge_term_1} and \eqref{bound_pull_in_edge_term_2}}\label{sec:defer_proof}
	
	We now briefly outline the proofs of the bounds \eqref{bound_pull_in_edge_term_1} and \eqref{bound_pull_in_edge_term_2} using the graphical tools developed in this section, omitting the full technical details. We take \eqref{bound_pull_in_edge_term_2} as an example; the proof of \eqref{bound_pull_in_edge_term_1} is similar but slightly simpler.
	
	For any term generated from the derivatives of $\Xi_{xy}$ or $\overline{\Xi}_{xy}$ on the left-hand side (LHS) of \eqref{bound_pull_in_edge_term_2}, we bound it as follows. For each $|G_{ab}|$-factor in the term, we bound it by a \smash{$\wt \xi^B$}-factor defined by
	\[\wt \xi^B_{ab}:=\delta_{ab}+\xi_{ab}+\qa{B_u(|a-b|)}^{1/2},\]
	where $\xi_{ij}$ is defined in \eqref{def_xi_edges}. For each $S_{ab}$-factor, we bound it by $\cS_{ab}$. For each $\Theta_{ab}$-factor, we bound it by two \smash{$\xi^B_{ab}$} factors, where \smash{$\xi^{B}_{ij}:=\xi_{ij}+[B_u(|i-j|)]^{1/2}$}. For each $\cK_{ab}=(S\Theta)_{ab}$ factor, we bound it by a waved edge together with two \smash{$\xi^B$}-edges. We interpret every such term as a value graph in which \smash{$\wt \xi^B$}-solid edges represent the \smash{$\wt \xi^B$}-factors, \smash{$\xi^B$}-solid edges represent the $\xi^B$-factors, and waved edges represent the $\cS$-factors.
	Note that each such graph contains two internal molecules and two external molecules containing $x$ and $y$, respectively. Moreover, there are four edge-disjoint paths of (\smash{$\wt \xi^B$} or $\xi^B$) solid edges on the molecular graph connecting the two external molecules. In addition, the partial derivative pulls a solid edge between molecules into the molecule containing $\beta$ and $j$, so there are at least 9 solid edges in the molecular graph.

	We next perform the molecular graph reduction as in Lemma \ref{lemma_val_to_aux}. More precisely, we have
	\begin{equation}\label{eq:newxi2}
      \sum_{a}\qa{\cS_{ia}(\wt \xi_{aj}^{B})^2+(\wt \xi_{ia}^{B})^2\cS_{aj}} \prec \ha{\xi_{ij}'+\qa{B_{u}(|i-j|)}^{1/2}}^2, 
\end{equation}
    and an analogous bound holds with the \(\wt\xi^B\)-factors on the LHS replaced by \(\xi^B\)-factors. Here, we again interpret the RHS as two \(\xi^B\)-edges, although the corresponding \(\xi\)-factors originate from different \(\Psi\)-matrices.
    We then apply the CS inequality, together with estimates such as \eqref{eq:newxi2}, to all relevant internal vertices---following an appropriate order of summation---to eliminate all waved edges and the $\delta$-terms from the \smash{$\wt \xi^{B}$}-edges. In the resulting graph, each of the four paths contributes a factor $\Psi_u(|x-y|)$ by selecting the longest $\xi^B$-edge along that path. The summation over the two internal vertices consumes four solid edges and produces a factor $\eta_u^{-2}\asymp (1-u)^{-2}$ by Ward's identity. The remaining extra solid edge contributes a factor $\Psi_u(0)$. Altogether, this yields the bound \eqref{bound_pull_in_edge_term_2}.

	\section{Dynamical analysis of \texorpdfstring{$T$}{T}-variables: Proofs of Lemmas \texorpdfstring{\ref{lemma_weak_local_law_T_along_flow}-\ref{theorem_que_along_flow}}{Lemmas}}\label{sec_proof_of_que}

	In this section, we establish the local laws for the $\cL$-loops and $T$-variables stated in Lemmas \ref{lemma_weak_local_law_T_along_flow}-\ref{theorem_que_along_flow}. The proof relies on the local laws and $\cL$-loop estimates obtained in \Cref{theorem_local_laws_along_flow}, together with a dynamical analysis of the evolution equations \eqref{T_Theta_evolution}, \eqref{integrated_T_Theta_evolution_d}, and \eqref{integrated_T_Theta_evolution_od} for the $T$-variables.

	\subsection{Proof of \Cref{lemma_weak_local_law_T_along_flow}}\label{sec_proof_of_weak_local_law_T_along_flow}
	
	From equation \eqref{T_Theta_evolution}, we deduce
	\begin{equation*}
		\begin{aligned}
			\rd (T-\Theta)_{t,x,yy'}^{\bsigma}=\sum_{a\in\ZL}\Theta_{t,xa}^{\bsigma}(T-\Theta)_{t,a, yy'}^{\bsigma}\,\rd t+\mathscr{E}_{t,x,yy'}^{T,\bsigma}\,\rd t+\cW_{t,x,yy'}^{T,\bsigma}\, \rd t+\rd \cB_{t,x,yy'}^{T,\bsigma},
		\end{aligned}
	\end{equation*}
	where the new term $\mathscr{E}_{t,x,yy'}^{T,\bsigma}$ combines  $\sum_a (T-\Theta)_{t,xa}^{\bsigma}\Theta_{t,a, yy'}^{\bsigma}$ with the quadratic error term $\mathcal{E}_{t,x,yy'}^{T,\bsigma}$: \begin{equation*}
		\begin{aligned}
			\mathscr{E}_{t,x,yy'}^{T,\bsigma}:=\sum_{a\in\ZL}(T-\Theta)_{t,xa}^{\bsigma}T_{t,a,yy'}^{\bsigma}=\sum_{a\in\ZL}(\cL-\cK)_{t,xa}^{\bsigma}\p{G_{t}^{\sigma_1}}_{ya}\p{G_{t}^{\sigma_2}}_{ay'}.
		\end{aligned}
	\end{equation*}
	Applying Duhamel's principle, we obtain that for any $\bsigma\in\h{-,+}^2$ and $x,y,y'\in\ZL$,
	\begin{equation}\label{integrated_T_Theta_evolution_scrE}
		\begin{aligned}
			(T-\Theta)_{t,x,yy'}^{\bsigma}=\qa{\cU_{s,t}^{\bsigma}\cdot(T-\Theta)_{s,\cdot,yy'}^{\bsigma}}_{x}&+\int_{s}^{t}\qb{\cU_{u,t}^{\bsigma}\cdot\cW_{u,\cdot,yy'}^{T,\bsigma}}_{x}\, \rd u\\
			&+\int_{s}^{t}\qb{\cU_{u,t}^{\bsigma}\cdot\mathscr{E}_{u,\cdot,yy'}^{T,\bsigma}}_{x}\, \rd u+\int_{s}^{t}\qb{\cU_{u,t}^{\bsigma}\cdot\rd \cB_{u,\cdot,yy'}^{T,\bsigma}}_{x}.
		\end{aligned}
	\end{equation}
	When $y\ne y'$, this equation reduces to \eqref{integrated_T_Theta_evolution_od}. For the diagonal case, however, we work with \eqref{integrated_T_Theta_evolution_scrE} instead of \eqref{integrated_T_Theta_evolution_d}, since the map $X\mapsto \cU_{s,t}^{\bsigma}\circ X$ does not propagate the weak decay $\smash{B_t^{1/2}}$ appropriately. More precisely, the decay of \smash{$ (\cU_{s,t}^{\bsigma}\cdot B_s^{1/2}\cdot \cU_{s,t}^{\bsigma})_{xy}$} may be substantially slower than that of \smash{$B_t^{1/2}$}. In contrast, the matrix product $X\mapsto \cU_{s,t}^{\bsigma}\cdot X$ preserves the desired $\smash{B_t^{1/2}}$ decay. (We emphasize that $\smash{B_t^{1/2}}$ denotes the entrywise square root, rather than the matrix square root.) This follows from
	\eqref{U_s_t_equal_1_Theta}, \eqref{Theta_bound_opposite_charge}, \eqref{Theta_bound_equaled_charge}, and the next lemma. The proof of \Cref{convolution_B_half_B} is deferred to \Cref{sec_proof_of_convolution_B_half_B}.
	
	\begin{lemma}\label{convolution_B_half_B}
		For any $0\leq s\leq t<1$ with $\ell_t<N$, we have
		\begin{equation}\label{convolution_half_B1}
			\sum_{a\in\ZL}B_t(|x-a|)B_s^{1/2}(|a-y|)\prec \frac{1}{1-s} B_t^{1/2}(|x-y|),
		\end{equation}
		while, for any $0\leq s\leq t\leq 1-W/N$ with $\ell_s=\ell_t=N$, we have
		\begin{equation}\label{convolution_half_B2}
			\sum_{a\in\ZL}\upcirc{B_t}(|x-a|)B_s^{1/2}(|a-y|)\prec \frac{1}{1-s}B_t^{1/2}(|x-y|).
		\end{equation}
	\end{lemma}

	As in the proof of \Cref{lemma_light_weight_term}, using the entrywise local law in \eqref{flow_local_law} together with \Cref{claim_graph_bound_light_weight_u}, we can bound the light-weight term as follows: for any $\bsigma\in\h{-,+}^2$ and $x,y,y'\in\ZL$,
	\begin{equation}\label{eq:LW_T1/2}
		\cW_{u,x,yy'}^{T,\bsigma}\prec \frac{B_u^{1/2}(0)}{1-u}B_u^{1/2}(|x-y|)B_u^{1/2}(|x-y'|).
	\end{equation}
	For the $\mathscr{E}$-term, applying \eqref{flow_local_law}, \eqref{flow_loop_local_law}, \eqref{convolution_2_B}, and Cauchy-Schwarz, we obtain 
	\begin{align}
		\mathscr{E}_{u,x,yy'}^{T,\bsigma}&\prec B_u^{1/5}(0) \sum_{a\in\ZL}B_u(|x-a|)\qb{\delta_{ay}+B_u^{1/2}(|a-y|)}\qb{\delta_{ay'}+B_u^{1/2}(|a-y'|)} \nonumber\\
		&\prec\frac{B_u^{1/5}(0)}{1-u}B_u^{1/2}(|x-y|)B_u^{1/2}(|x-y'|).\label{eq:scrE_T1/2}
	\end{align}
	For the martingale term, we bound the tensor in \eqref{def_2tensor_cB2} (with $x_1=x_2=x$) as
	\begin{align}
		(\cB\otimes\cB)_{u;yy'}^{T,\bsigma}(x,x)&=\sum_{a\in\ZL}\pa{\absb{(G_u^{\sigma_1}S^{\pa{x}}G_u^{\sigma_2})_{ay'}}^2(G_u^{\sigma_1}S^{\pa{a}}G_u^{-\sigma_1})_{yy} +\absb{(G_u^{\sigma_1}S^{\pa{x}}G_u^{\sigma_2})_{ya}}^2(G_u^{-\sigma_2}S^{\pa{a}}G_u^{\sigma_2})_{y'y'}} \nonumber\\
		&\prec B_u(0)\qa{(G_u^{-\sigma_2}S^{\pa{x}}G_u^{-\sigma_1}G_u^{\sigma_1}S^{\pa{x}}G_u^{\sigma_2})_{y'y'}+(G_u^{\sigma_1}S^{\pa{x}}G_u^{\sigma_2}G_u^{-\sigma_2}S^{\pa{x}}G_u^{-\sigma_1})_{yy}}\nonumber\\
		&=\frac{B_u(0)}{\eta_u}\qa{(G_u^{-\sigma_2}S^{\pa{x}}(\im G_u)S^{\pa{x}}G_u^{\sigma_2})_{y'y'}+(G_u^{\sigma_1}S^{\pa{x}}(\im G_u)S^{\pa{x}}G_u^{-\sigma_1})_{yy}}\nonumber\\
		& \prec \frac{B_u^{3/2}(0)}{1-u}\qa{B_u(|x-y|)+B_u(|x-y'|)}.\label{martingale_bound_T_od}
	\end{align}
	Here, in the second step we used the local law \eqref{flow_local_law} to bound $(G_u^{\sigma_1}S^{\pa{a}}G_u^{-\sigma_1})_{yy}$ and $(G_u^{-\sigma_2}S^{\pa{a}}G_u^{\sigma_2})_{y'y'}$ by $\OO_\prec(B_u(0))$, in the third step we applied Ward's identity, and in the last step we used the local law \eqref{flow_local_law} again together with the convolution condition \eqref{polydecay_Psiu3} for the control parameter $B_u$, which follows directly from its definition.
	
	\begin{remark}\label{remark_weak_decay_T_local_law}
		Note that the estimate \eqref{martingale_bound_T_od} loses one factor of $B_u$ decay---namely, ${B_u(|x-y|)+B_u(|x-y'|)}$ instead of ${B_u(|x-y|)\cdot B_u(|x-y'|)}$. One could instead use the argument from the proof of \Cref{lemma_martingale_term}, as in \Cref{sec_estimates_error_terms2}. Although that approach yields the optimal decay $(1-u)^{-1}\cdot B_u(|x-y|)B_u(|x-y'|)$, it misses a factor of $B_u^{1/2}(0)$ (recall \Cref{rmk:MG_for_T_fails}) and therefore only gives
		\[\pa{T-\Theta}_{t,x,yy'}^{\bsigma}\prec B_t^{1/2}(|x-y|)B_t^{1/2}(|x-y'|),\] 
		which is in fact a trivial consequence of the local law \eqref{flow_local_law}. In particular, this leads to the weaker bound $(N\eta)^{-1}$ instead of $(N\eta)^{-6/5}$ in \eqref{weak_local_law_T_along_flow} at $1-t=\eta_{\txt{flat}}$. The fact that the exponent $6/5$ is strictly larger than $1$ is crucial for proving \Cref{lemma_flat_local_law_along_flow,theorem_que_along_flow}.  
	\end{remark}

	We now set $s=0$ and take $t\in[0,1-(W/N)^\alpha]$ in \eqref{integrated_T_Theta_evolution_scrE}. In this case, the first term vanishes. Combining \eqref{eq:LW_T1/2} and \eqref{eq:scrE_T1/2} with \eqref{U_s_t_equal_1_Theta}, \eqref{Theta_bound_opposite_charge}, and \eqref{convolution_2_B}, and performing the integration in $u$, we obtain the bound
	\begin{equation}\label{eq:LWQuad_Teq}       \int_{0}^{t}\qb{\cU_{u,t}^{\bsigma}\cdot\cW_{u,\cdot,yy'}^{T,\bsigma}}_{x}\, \rd u+\int_{0}^{t}\qb{\cU_{u,t}^{\bsigma}\cdot\mathscr{E}_{u,\cdot,yy'}^{T,\bsigma}}_{x}\, \rd u \prec B_t^{1/5}(0)B_t^{1/2}(|x-y|)B_t^{1/2}(|x-y'|).
	\end{equation}
	For the martingale term, applying \eqref{bound_martingale_to_cB_T_od} and arguing as in \eqref{eq:UUUUBB} via the CS inequality, we deduce that for any fixed $p\in 2\N$,
	\begin{equation*}
		\begin{aligned}        \E\absa{\int_{0}^{t}\qb{\cU_{u,t}^{\bsigma}\cdot\rd \cB_{u,\cdot,yy'}^{T,\bsigma}}_{x}}^{p}\lesssim \E\pa{\int_0^{t}\qa{\cU_{u,t}^{(-,+)}\cdot v_{u,\cB}^{T,\bsigma}}^2_{x}\, \rd u}^{p/2},
		\end{aligned}
	\end{equation*}
	where the vector $v_{u,\cB}^{T,\bsigma}$ is defined by $v_{u,\cB}^{T,\bsigma}(x):=\qb{(\cB\otimes\cB)_{u;yy'}^{T,\bsigma}(x,x)}^{1/2}$ for $x\in\ZL$. 
	Invoking \Cref{convolution_B_half_B} and \eqref{martingale_bound_T_od}, and integrating in $u$, we obtain
	\begin{equation*}
			\begin{aligned}        \E\absa{\int_{0}^{t}\qb{\cU_{u,t}^{\bsigma}\cdot\rd \cB_{u,\cdot,yy'}^{T,\bsigma}}_{x}}^{p} \prec \qa{B_t^{3/2}(0)B_t(|x-y|\wedge |x-y'|)}^{p/2}.
		\end{aligned}
	\end{equation*}
	By Markov's inequality, this implies
	\begin{equation}\label{eq:MG_Teq}
		\int_{0}^{t}\qb{\cU_{u,t}^{\bsigma}\cdot\rd \cB_{u,\cdot,yy'}^{T,\bsigma}}_{x}\prec B_t^{3/4}(0)B_t^{1/2}(|x-y|\wedge |x-y'|)\le B_t^{7/10}(0)B_t^{1/2}(|x-y|\wedge |x-y'|).
	\end{equation}
	Substituting \eqref{eq:LWQuad_Teq} and \eqref{eq:MG_Teq} into \eqref{integrated_T_Theta_evolution_scrE}, we conclude \eqref{weak_local_law_T_along_flow} for all $t\in[0,1-(W/N)^\alpha]$.
	
	Next, taking $s=1-(W/N)^{\alpha}$ and $t\in[1-(W/N)^{\alpha},1-W/N]$ in \eqref{integrated_T_Theta_evolution_scrE}, we may bound the second and third terms exactly as in \eqref{eq:LWQuad_Teq}. It therefore remains to estimate the first and fourth terms in \eqref{integrated_T_Theta_evolution_scrE}. 
	For the first term, using \eqref{U_s_t_equal_1_Theta} together with the initial estimate \eqref{weak_local_law_T_along_flow} established at time $s$, we obtain
	\begin{align}
		\qa{\cU_{s,t}^{\bsigma}\cdot(T-\Theta)_{s,\cdot,yy'}^{\bsigma}}_{x}&=(T-\Theta)_{s,x,yy'}^{\bsigma}+(t-s)\qa{\Theta_t^{\bsigma}\cdot(T-\Theta)_{s,\cdot,yy'}^{\bsigma}}_x \nonumber\\
		&=(t-s)\qa{\Theta_t^{\bsigma}\cdot(T-\Theta)_{s,\cdot,yy'}^{\bsigma}}_x+\opr{B_s^{7/10}(0)B_s^{1/2}(|x-y|\wedge |x-y'|)}. \label{eq:UST_initial}
	\end{align}
	It remains to control the first term on the RHS. If $\bsigma=(\sigma_1,\sigma_2)$ satisfies $\sigma_1=\sigma_2$, then by \eqref{Theta_bound_equaled_charge}, we have
	\begin{equation}\label{eq:BBcirc}
		\absa{\Theta_{t,ab}^{\bsigma}}\prec B_0(|a-b|)\lesssim \upcirc{B_t}(|a-b|).
	\end{equation}
	Combining this with \eqref{convolution_half_B2}, we deduce
	\begin{align}
		(t-s)&\qa{\Theta_t^{\bsigma}\cdot(T-\Theta)_{s,\cdot,yy'}^{\bsigma}}_x\prec (t-s) B_s^{7/10}(0) \sum_{a\in\ZL}\mathring B_t(|x-a|)B_s^{1/2}(|a-y|\wedge|a-y'|)\nonumber\\
		&\prec \frac{t-s}{1-s}B_t^{7/10}(0)\qa{B_t^{1/2}(|x-y|)+B_t^{1/2}(|x-y'|)}\lesssim B_t^{7/10}(0)B_t^{1/2}(|x-y|\wedge|x-y'|).\label{eq:UST_initial_equal}
	\end{align}
	It remains to treat the case $\sigma_1\neq\sigma_2$. To this end, we introduce the following \emph{zero-mode-removing operator}.
	
	\begin{definition}[Zero-mode-removing operator]\label{def_sum_zero_operators}
		For $\bv\in\C^{\ZL}$, we define the zero-mode-removing operator by
		\begin{equation*}
			\begin{aligned}
				(\cQ\bv)_{x}:=\bv(x)-\frac{1}{N}\sum_{y\in\ZL}\bv(y),\quad \forall  x\in\ZL .
			\end{aligned}
		\end{equation*}
		Equivalently, $\cQ$ can be represented by the matrix $\cQ = 1 - \cP$, where $\cP$ is the rank-one projection given by $\cP_{xy} := N^{-1}$ for all $x,y\in\ZL$. Similarly, for a matrix $X\in\C^{\ZL\times\ZL}$, we define
		\begin{equation*}
			\begin{aligned}
				(\cO\circ X)_{xy}:=(\cQ X\cQ)_{xy}=X_{xy}-\frac{1}{N}\sum_{a\in\ZL}X_{ay}-\frac{1}{N}\sum_{b\in\ZL}X_{xb}+\frac{1}{N^2}\sum_{a,b\in\ZL}X_{ab}.
			\end{aligned}
		\end{equation*}
		By construction, for any $\bv\in\C^{\ZL}$, $X\in\C^{\ZL\times\ZL}$, and $x,y\in\ZL$, we have
		\begin{equation*}
			\begin{aligned}
				\sum_{a\in\ZL}(\cQ\bv)_a=0,\quad \sum_{a\in\ZL}(\cO\circ X)_{ay}=\sum_{b\in\ZL}(\cO\circ X)_{xb}=0.
			\end{aligned}
		\end{equation*}
		Note that the matrices $\cQ$ and $\cP$ are orthogonal projections. In particular, they commute with $S$, and hence also with $\cU_{u,r}^{\bsigma}$ and $\Theta_u^{\bsigma}$ for all $\bsigma\in\h{-,+}^2$ and $0\le u\le r\le \tf$.
	\end{definition}
	
	We now consider the case $\sigma_1\neq\sigma_2$. With the above notation, we decompose
	\begin{align}
		(t-s)&\qa{\Theta_t^{\bsigma}\cdot(T-\Theta)_{s,\cdot,yy'}^{\bsigma}}_x=(t-s)\qa{\cP\cdot\Theta_t^{\bsigma}\cdot(T-\Theta)_{s,\cdot,yy'}^{\bsigma}}_x+(t-s)\qa{{\cQ\cdot\Theta_t^{\bsigma}}\cdot(T-\Theta)_{s,\cdot,yy'}^{\bsigma}}_x.\nonumber
	\end{align}
	Using $(\cP\cdot\Theta_t^{\bsigma})_{xy} \equiv [N(1-t)]^{-1}$, the identity \smash{$\sum_x \Theta_{s,xy}^{\bsigma}\equiv(1-s)^{-1}=\im m/\eta_s$} from \eqref{re_im_z_t}, and Ward's identity \eqref{eq_Ward}, we obtain
	\begin{align*}
		(t-s)\qa{\cP\cdot\Theta_t^{\bsigma}\cdot(T-\Theta)_{s,\cdot,yy'}^{\bsigma}}_x=\frac{t-s}{1-t}\frac{\qa{\im(G_s-m)}_{yy'}}{N\eta_s}\prec \frac{B_s^{1/2}(0)}{N(1-t)}\lesssim B_t^{7/10}(0)B_t^{1/2}(|x-y|\wedge|x-y'|),
	\end{align*}
	where we used \eqref{flow_local_law} in the second step. For the second term, using the bound $\absa{\pa{\cQ\cdot\Theta_t^{\bsigma}}_{ab}}\prec \upcirc{B_t}(|a-b|)$ from \eqref{zero_mode_removed_bound}, together with the same argument leading to \eqref{eq:UST_initial_equal}, we obtain
	\begin{align*}
		(t-s)\qa{{\cQ\cdot\Theta_t^{\bsigma}}\cdot(T-\Theta)_{s,\cdot,yy'}^{\bsigma}}_x\prec B_t^{7/10}(0)B_t^{1/2}(|x-y|\wedge|x-y'|).
	\end{align*}
	Combining the above two estimates yields
	\begin{align}
		(t-s)\qa{\Theta_t^{\bsigma}\cdot(T-\Theta)_{s,\cdot,yy'}^{\bsigma}}_x
		\prec B_t^{7/10}(0)B_t^{1/2}(|x-y|\wedge|x-y'|).\label{eq:UST_initial_opposite}
	\end{align}
	Substituting \eqref{eq:UST_initial_equal} and \eqref{eq:UST_initial_opposite} into \eqref{eq:UST_initial}, we see that the first term in \eqref{integrated_T_Theta_evolution_scrE} is bounded by the RHS of \eqref{weak_local_law_T_along_flow}.

	It remains to control the fourth term in \eqref{integrated_T_Theta_evolution_scrE}. First consider the case $\bsigma=(\sigma_1,\sigma_2)$ with $\sigma_1=\sigma_2$. By \eqref{bound_martingale_to_cB_T_od}, for any fixed $p\in2\N$, we have
	\begin{equation}\label{eq:martigale_Tp}
		\E\absbb{\int_{s}^{t}\qa{\cU_{u,t}^{\bsigma}\cdot\rd \cB_{u,\cdot,yy'}^{T,\bsigma}}_{x}}^p\lesssim\E\pa{\int_{s}^t\qa{(\cU_{u,t}^{\bsigma}\otimes \cU_{u,t}^{\ol{\bsigma}})\circ(\cB\otimes\cB)_{u,yy'}^{T,\bsigma}}_{xx}\,\rd u}^{p/2},
	\end{equation}
	where $\ol{\bsigma}:=(-\sigma_2,-\sigma_1)$. Arguing as in \eqref{martingale_bound_T_od}, we obtain
	\begin{align}
		&\qa{(\cU_{u,t}^{\bsigma}\otimes \cU_{u,t}^{\ol{\bsigma}})\circ(\cB\otimes\cB)_{u,yy'}^{T,\bsigma}}_{xx}\nonumber\\
		&=\sum_{a,x_1,x_2\in\ZL}\cU_{u,t}^{\bsigma}(x,x_1)\cU_{u,t}^{\ol{\bsigma}}(x,x_2)(G_u^{-\sigma_2}S^{\pa{x_2}}G_u^{-\sigma_1})_{y'a}(G_u^{\sigma_1}S^{\pa{a}}G_u^{-\sigma_1})_{yy}(G_u^{\sigma_1}S^{\pa{x_1}}G_u^{\sigma_2})_{ay'}\nonumber\\
		&+\sum_{a,x_1,x_2\in\ZL}\cU_{u,t}^{\bsigma}(x,x_1)\cU_{u,t}^{\ol{\bsigma}}(x,x_2)(G_u^{\sigma_1}S^{\pa{x_1}}G_u^{\sigma_2})_{ya}(G_u^{-\sigma_2}S^{\pa{a}}G_u^{\sigma_2})_{y'y'}(G_u^{-\sigma_2}S^{\pa{x_2}}G_u^{-\sigma_1})_{ay}\nonumber\\
		&\prec \frac{B_u(0)}{1-u}\sum_{x_1,x_2\in\ZL}\absa{\cU_{u,t}^{\bsigma}(x,x_1)\cU_{u,t}^{\ol{\bsigma}}(x,x_2)}(G_u^{-\sigma_2}S^{\pa{x_2}}(\im G_u) S^{\pa{x_1}}G_u^{\sigma_2})_{y'y'}\nonumber\\
		&+\frac{B_u(0)}{1-u}\sum_{x_1,x_2\in\ZL}\absa{\cU_{u,t}^{\bsigma}(x,x_1)\cU_{u,t}^{\ol{\bsigma}}(x,x_2)}(G_u^{\sigma_1}S^{\pa{x_1}}(\im G_u)S^{\pa{x_2}}G_u^{-\sigma_1})_{yy}.\label{martingale_bound_T_od_equal_charge}
	\end{align}
	Using \eqref{U_s_t_equal_1_Theta}, \eqref{Theta_bound_equaled_charge}, the local law \eqref{flow_local_law}, and the bound \eqref{eq:BBcirc}, we further estimate the RHS by 
	\begin{equation*}
		\begin{aligned}
			(\ref{martingale_bound_T_od_equal_charge})\prec &~\frac{B_u^{3/2}(0)}{1-u}\sum_{x_1,x_2\in\ZL}[\delta_{xx_1}+(t-u)\upcirc{B_u}(|x-x_1|)][\delta_{xx_2}+(t-u)\upcirc{B_u}(|x-x_2|)]\\
			&~\times \qa{B_u^{1/2}(|x_1-y'|) B_u^{1/2}(|x_2-y'|)+ B_u^{1/2}(|x_1-y|) B_u^{1/2}(|x_2-y|)}\\
			\prec&~        \frac{B_u^{3/2}(0)}{1-u}\qa{B_u(|x-y'|)+B_u(|x-y|)},
		\end{aligned}
	\end{equation*}
	where in the second step, we applied the same argument as in \eqref{eq:UST_initial_equal}, based on \eqref{convolution_half_B2}. Substituting this bound into \eqref{eq:martigale_Tp}, performing the integration in $u$, and applying Markov's inequality, we obtain
	\begin{equation}\label{eq:Tmartingale_equal}
		\int_{s}^{t}\qa{\cU_{u,t}^{\bsigma}\cdot\rd \cB_{u,\cdot,yy'}^{T,\bsigma}}_{x}\prec B_t^{3/4}(0)B_t^{1/2}(|x-y|\wedge |x-y'|)\lesssim B_t^{7/10}(0)B_t^{1/2}(|x-y|\wedge |x-y'|).
	\end{equation}
	We next consider the case $\sigma_1\neq\sigma_2$. Using $(\cP\cdot\Theta_t^{\bsigma})_{xy}\equiv [N(1-t)]^{-1}$, we decompose
	\begin{equation*}
		\begin{aligned}
			\int_{s}^{t}\qa{\cU_{u,t}^{\bsigma}\cdot\rd \cB_{u,\cdot,yy'}^{T,\bsigma}}_{x}=\int_{s}^{t}\rd \cB_{u,x,yy'}^{T,\bsigma}&+\int_{s}^{t}(t-u)\qa{(\cQ\cdot\Theta_t^{\bsigma})\cdot\rd \cB_{u,\cdot,yy'}^{T,\bsigma}}_{x} +\int_{s}^{t}\frac{t-u}{1-t}\qa{\cP\cdot \rd \cB_{u,\cdot,yy'}^{T,\bsigma}}_{x}.
		\end{aligned}
	\end{equation*}
	The first term is bounded directly by \eqref{bound_martingale_to_cB_T_od} and \eqref{martingale_bound_T_od}. 
	The second term is handled exactly as in the case $\sigma_1=\sigma_2$, except that we now use the estimate \smash{$\pa{\cQ\cdot\Theta_{t}^{\bsigma}}_{ab}\prec \upcirc{B_t}(|a-b|)$} from \eqref{zero_mode_removed_bound} in place of \eqref{eq:BBcirc}. For the third term, we again apply \eqref{bound_martingale_to_cB_T_od} together with
	\begin{equation*}
		\begin{aligned}
			\qa{(\cP\otimes \cP)\circ(\cB\otimes\cB)_{u,yy'}^{T,\bsigma}}_{xx}\prec\frac{B_u(0)}{1-u}\pa{\frac{1}{N(1-u)}}^{2},
		\end{aligned}
	\end{equation*}
	which follows readily from Ward's identity and the local law \eqref{flow_local_law}. 
	Combining these bounds yields \eqref{eq:Tmartingale_equal} also in the case $\sigma_1\neq\sigma_2$, thereby completing the proof of \Cref{lemma_weak_local_law_T_along_flow}.

	\subsection{Proof of \Cref{lemma_flat_local_law_along_flow}}\label{sec_proof_of_flat_local_law_along_flow}
	Given any small constant $c>0$, we introduce the stopping times
	\begin{align}
\tau^{\cL}\equiv\tau^{\cL}(c):=\tf\wedge\inf\ha{t\geq  1-\frac{W}{N}:\max_{\bsigma\in\h{-,+}^2}\norma{(\cL-\cK)_{t}^{\bsigma}}_{\max}\geq \frac{W^{c}}{W^{6/5}}+\frac{W^{c}}{\q{N(1-t)}^{7/4}} },\label{tauLt}\\
 \tau^{T}\equiv\tau^{T}(c):=\tf\wedge \inf\ha{t\geq 1-\frac{W}{N}:\max_{\bsigma\in\h{-,+}^2}\max_{x,y,y'\in\ZL}\absa{(T-\Theta)_{t,x,yy'}^{\bsigma}}\geq \frac{W^{c}}{W^{6/5}}+ \frac{W^{c}}{\q{N(1-t)}^{3/2}}},\label{tauTt}
	\end{align}
 where we adopt the convention $\inf\emptyset=\infty$.
	By \eqref{flow_loop_local_law} and \eqref{weak_local_law_T_along_flow} at $t=1-W/N$, for any small constant $c>0$, the following holds with high probability
	\begin{equation}\label{eq:initialW/N}
		\max_{\bsigma}\norma{(\cL-\cK)_{t}^{\bsigma}}_{\max}\leq W^{-6/5+c},\qquad \max_{\bsigma}\max_{x,y,y'\in\ZL}\absa{(T-\Theta)_{t,x,yy'}^{\bsigma}}\leq W^{-6/5+c}.
	\end{equation}
	By continuity in $t$, it follows that $\tau^{\cL}(c)>1-W/N$ and $\tau^{T}(c)>1-W/N$ with high probability. 
	In the sequel, we fix two sufficiently small constants $c_{\cL},c_T\in(0,1)$ such that $c_{\cL}\ge 100 c_T$ and $W^{100c_{\cL}}B_{\tf}(0)\le 1$. For $\star\in\{\cL,T\}$, we write $\tau^\star=\tau^\star(c_\star)$ and set $\tau:=\tau^{\cL}\wedge \tau^{T}$. We prove \eqref{flat_local_law_along_flow} and \eqref{flat_local_law_along_flow2} by establishing a self-improving estimate for the $\cL$-loops and $T$-variables, based on the evolution equations \eqref{integrated_cL_cK_evolution}, \eqref{integrated_T_Theta_evolution_d}, and \eqref{integrated_T_Theta_evolution_od} with initial time $s=1-W/N$. We begin by estimating the error terms appearing in these equations.

	We first control the light-weight terms in \eqref{cL_cK_evolution} and \eqref{T_Theta_evolution}. Uniformly for $u\in[s,\tau]$, we claim
	\begin{equation}\label{flat_light_weight_bound}
		\begin{aligned}        \max_{\bsigma\in\h{-,+}^2}\max_{x,y\in\ZL}\absa{\cW_{u,xy}^{\cL,\bsigma}}+\max_{\bsigma\in\h{-,+}^2}\max_{x,y,y'\in\ZL}\absa{\cW_{u,x,yy'}^{T,\bsigma}}\prec \frac{W^{c_T}}{1-u}\qa{\frac{W^{-6/5}}{\sqrt{N(1-u)}}+ \frac{1}{\q{N(1-u)}^{2}}} .
		\end{aligned}
	\end{equation}
	We prove this bound for $\cW_{u,xy}^{\cL,\bsigma}$; the term $\cW_{u,x,yy'}^{T,\bsigma}$ is handled analogously. By the definition of $\tau$ and the local laws in  \eqref{flow_local_law}, the first term in \eqref{def_light_weight_term} satisfies
	\begin{equation*}
		\begin{aligned}
			&\sum_{a}\pa{G_u^{\sigma_2}-m^{\sigma_2}}_{aa}\avg{G_u^{\sigma_2}S^{\pa{a}}G_u^{\sigma_2}S^{\pa{y}}G_u^{\sigma_1}S^{\pa{x}}}=\sum_{b,c}\pb{G_u^{\sigma_1}S^{\pa{x}}G_u^{\sigma_2}}_{cb}\avgb{(G_u^{\sigma_2}-m^{\sigma_2})S^{\pa{b}}}\pb{G_u^{\sigma_2}}_{bc}S_{cy}\\
			&\prec\sum_{b\neq c}\qa{\frac{W^{c_T}}{W^{6/5}}+\frac{W^{c_T}}{\q{N(1-u)}^{3/2}}}\frac{1}{\q{N(1-u)}^{3/2}}S_{cy}+\sum_{c}\frac{1}{\q{N(1-u)}^2}S_{cy} \\
			&\lesssim \frac{W^{c_T}}{1-u}\qa{\frac{W^{-6/5}}{\sqrt{N(1-u)}}+ \frac{1}{\q{N(1-u)}^{2}}}.     \end{aligned}
	\end{equation*}
	By symmetry, the second term in \eqref{def_light_weight_term} obeys the same bound, which establishes \eqref{flat_light_weight_bound} for $\cW_{u,xy}^{\cL,\bsigma}$.
	Next, using the definitions \eqref{tauLt} and \eqref{tauTt}, the quadratic error term in \eqref{cL_cK_evolution} satisfies, uniformly for $u\in[s,\tau]$,
	\begin{equation}\label{flat_double_difference_bound_cL}
		\begin{aligned}
			\max_{\bsigma\in\h{-,+}^2}\max_{x,y\in\ZL}\absa{\cE_{u,xy}^{\cL,\bsigma}}&\le   \sum_{a}\qa{\frac{W^{c_T}}{W^{6/5}}+ \frac{W^{c_T}}{\q{N(1-u)}^{3/2}}}\qa{\frac{W^{c_{\cL}}}{W^{6/5}}+\frac{W^{c_{\cL}}}{\q{N(1-u)}^{7/4}}}\\
			&\lesssim \frac{ W^{c_T+c_{\cL}}}{1-u}\frac{1}{\q{N(1-u)}^{1/5}}\qa{\frac{1}{W^{6/5}}+\frac{1}{\q{N(1-u)}^{7/4}}},
		\end{aligned}
	\end{equation}
	where we used that $W^{-1}\le [N(1-u)]^{-1}$ for $u\ge 1-W/N$ in the second step. Similarly, the quadratic error term in \eqref{T_Theta_evolution} satisfies, uniformly for $u\in[s,\tau]$,
	\begin{equation}\label{flat_double_difference_bound_T}
		\begin{aligned}        \max_{\bsigma\in\h{-,+}^2}\max_{x,y,y'\in\ZL}\absa{\mathcal{E}_{u,x,yy'}^{T,\bsigma}}\prec \frac{W^{2c_T}}{1-u} \frac{1}{\q{N(1-u)}^{1/5}}\qa{\frac{1}{W^{6/5}}+ \frac{1}{\q{N(1-u)}^{3/2}}}.
		\end{aligned}
	\end{equation}

	Finally, to bound the martingale terms in \eqref{integrated_cL_cK_evolution}, \eqref{integrated_T_Theta_evolution_d}, and \eqref{integrated_T_Theta_evolution_od}, we estimate the tensors in \eqref{def_tensor_cB}, \eqref{def_tensor_cB2}, and \eqref{def_2tensor_cB2} as follows uniformly in $u\in\q{s,\tau}$:
	\begin{align}
		\max_{x_1,x_2,y_1,y_2,\bsigma} \pa{   \absa{(\cB\otimes \cB)_{u,y_1y_2}^{T,\bsigma}(x_1,x_2)}+\absa{(\cB\otimes \cB)_{u,x_1x_2y_1y_2}^{T,\bsigma}}}\prec \frac{W^{c_{T}}}{1-u} \qa{\frac{1}{W^{6/5}}+\frac{1}{\q{N(1-u)}^{3/2}}}^2,\label{martingale_offT}\\
		\max_{x_1,x_2,y_1,y_2,\bsigma}\absa{(\cB\otimes \cB)_{u,x_1x_2y_1y_2}^{\cL,\bsigma}}\prec  \frac{W^{2c_{T}}}{1-u} \frac{1}{\sqrt{N(1-u)}} \qa{\frac{1}{W^{6/5}}+ \frac{1}{\q{N(1-u)}^{3/2}}}^2.\label{martingale_cL}
	\end{align}
	These bounds are proved by arguments analogous to those used in \Cref{lemma_martingale_term} and in the derivation of \eqref{martingale_bound_T_od}; we therefore only sketch the main steps. First, the estimate for $(\cB\otimes \cB)_{u,y_1y_2}^{T,\bsigma}(x_1,x_2)$ follows the same strategy as in \eqref{martingale_bound_T_od}. The only modification is that, in the final step, the resulting three-resolvent chain contains a $T$-variable; this term is controlled by combining the local law \eqref{flow_local_law} with the a priori bound on $T$ from \eqref{tauTt}.
	Next, the bound for $(\cB\otimes \cB)_{u,x_1x_2y_1y_2}^{T,\bsigma}$ follows immediately from the CS inequality,
	\begin{equation*}
		\begin{aligned}        \absa{(\cB\otimes\cB)_{u;x_1x_2y_1y_2}^{T,\bsigma}}^2\leq (\cB\otimes\cB)_{u;x_1x_1y_1y_1}^{T,\bsigma}(\cB\otimes\cB)_{u;x_2x_2y_2y_2}^{T,\bsigma}=(\cB\otimes\cB)_{u;y_1y_1}^{T,\bsigma}(x_1,x_1)(\cB\otimes\cB)_{u;y_2y_2}^{T,\bsigma}(x_2,x_2),
		\end{aligned}
	\end{equation*}
	and the previously established bound for the diagonal terms. Finally, \eqref{martingale_cL} is obtained by an argument parallel to the proof of \Cref{lemma_martingale_term} in \Cref{sec_estimates_error_terms2}. More precisely, we bound the relevant expression by \eqref{bound_I_leq}, now summing over the full region $a\in\ZL$. In the subsequent analysis, the four-resolvent loop in \eqref{4_loop_term_Ward_inequality}, which contains two $T$-variables, is controlled using the stopping-time bound from \eqref{tauTt}. Similarly, the three-resolvent loop in \eqref{3_loop_term_Ward_inequality}, containing one $T$-variable, is estimated by combining \eqref{tauTt} with the local law \eqref{flow_local_law} for the remaining resolvent entry.

	We now combine the bounds \eqref{flat_light_weight_bound}--\eqref{martingale_cL} to complete the proof of \eqref{flat_local_law_along_flow} and \eqref{flat_local_law_along_flow2}.
	First, if $\bsigma=(\sigma_1,\sigma_2)$ satisfies $\sigma_1=\sigma_2$, substituting the estimates \eqref{flat_light_weight_bound}--\eqref{martingale_cL}, together with the initial condition \eqref{eq:initialW/N}, into \eqref{integrated_cL_cK_evolution}, \eqref{integrated_T_Theta_evolution_d}, and \eqref{integrated_T_Theta_evolution_od}, and applying the evolution kernel bound \eqref{bound_cU_s_t_3} for case (i), we integrate over $u$ (using \Cref{lemma_martingale_to_cB} for the martingale term\footnote{{In fact, \Cref{lemma_martingale_to_cB} does not directly apply here, since the coefficient matrices of the martingale term depends on a random stopping time and the BDG inequality therefore fails. However, the desired estimates can be recovered by first applying \Cref{lemma_martingale_to_cB} with $\tau$ replaced by any fixed time, and then employing a standard $N^{-C}$-net argument to extend the estimates to the random time $\tau$. Moreover, this procedure introduces an additional term of order $\OO\pa{N^{-D}}$ for any large constant $D\geq 1$, which could be absorbed by other much larger terms in the following steps. Therefore, here and below, we omit this subtlety whenever such application of BDG inequality is required.}}) 
Exploiting the assumptions $W^{100c_{\cL}}B_{\tf}(0)\le 1$ and $c_{\cL}\ge 100c_T$, we obtain that for all $t\in[s,\tau]$ and $x,y,y'\in\ZL$,
	\begin{equation}\label{self_improved_flat_bound}
		\begin{aligned}
			(\cL-\cK)_{t,xy}^{\bsigma}\prec \frac{W^{c_T}}{W^{6/5}}+\frac{W^{c_T}}{\q{N(1-t)}^{7/4}},\quad (T-\Theta)_{t,x,yy'}^{\bsigma}\prec \frac{W^{c_T/2}}{W^{6/5}}+ \frac{W^{c_T/2}}{\q{N(1-t)}^{3/2}}.
		\end{aligned}
	\end{equation}
	In the case $\sigma_1\neq \sigma_2$, we employ the zero-mode-removing operators from \Cref{def_sum_zero_operators}. Applying $\cQ$ to both sides of \eqref{integrated_T_Theta_evolution_od}, and $\cO$ to both sides of \eqref{integrated_cL_cK_evolution} and \eqref{integrated_T_Theta_evolution_d}, we proceed as in the equal-sign case, now using the evolution kernel bound \eqref{bound_cU_s_t_3} corresponding to case (ii). This yields, for all $t\in[s,\tau]$ and $x,y,y'\in\ZL$ with $y\neq y'$,
	\begin{equation}\label{self_improved_flat_bound2}
		\begin{aligned}
			&\qa{\cO\circ(\cL-\cK)_{t}^{\bsigma}}_{xy}\prec \frac{W^{c_T}}{W^{6/5}}+\frac{W^{c_T}}{\q{N(1-t)}^{7/4}},\\ 
			&\absa{\qa{\cO\circ(T-\Theta)_{t}^{\bsigma}}_{xy}}+\absb{\q{\cQ\cdot(T-\Theta)_{t,\cdot,yy'}^{\bsigma}}_{x}}\prec \frac{W^{c_T/2}}{W^{6/5}}+ \frac{W^{c_T/2}}{\q{N(1-t)}^{3/2}}.
		\end{aligned}
	\end{equation}
	On the other hand, using Ward's identity together with the averaged local law in \eqref{flow_local_law}, we obtain 
	\begin{align}
		\qa{(\cL-\cK)_{t}^{\bsigma}-\cO\circ(\cL-\cK)_{t}^{\bsigma}}_{xy}&=\frac{\im \q{\avg{G_tS^{(x)}}-m}}{N\eta_{t}}+\frac{\im \q{\avg{G_tS^{(y)}}-m}}{N\eta_{t}}-\frac{\im\q{N^{-1}\avg{G_t}-m}}{N\eta_{t}}\nonumber\\
		&\prec {\q{N(1-t)}^{-2}}.\label{eq:cL-cK_flat}
	\end{align}
	Similarly, using Ward's identity and the entrywise local law in \eqref{flow_local_law}, we obtain  
	\begin{align*}
		&\absa{\qa{(T-\Theta)_{t}^{\bsigma}-\cO\circ(T-\Theta)_{t}^{\bsigma}}_{xy}}+\absb{\q{(T-\Theta)_{t,\cdot,yy'}^{\bsigma}-\cQ\cdot(T-\Theta)_{t,\cdot,yy'}^{\bsigma}}_{x}}\prec  {\q{N(1-t)}^{-3/2}}.
	\end{align*}
	Combining these bounds with \eqref{self_improved_flat_bound2}, we conclude that \eqref{self_improved_flat_bound} remains valid when $\sigma_1\neq\sigma_2$.
	
	We observe that the estimates in \eqref{self_improved_flat_bound} improve upon the bounds assumed in the definitions of $\tau^{\cL}(c_{\cL})$ and $\tau^{T}(c_T)$ by a factor of $W^{-c_T/2}$, since $c_{\cL}\ge 100c_T$. Together with the fact that $\tau=\tau^{\cL}(c_{\cL})\wedge\tau^{T}(c_T)>s$ with high probability and a standard continuity argument, this implies that $\tau\ge \tf$ with high probability.
	This completes the proof of \Cref{lemma_flat_local_law_along_flow}, as $c_{\cL}$ and $c_T$ may be chosen arbitrarily small.

	\subsection{Proof of \texorpdfstring{\Cref{theorem_que_along_flow}}{Theorem 7.1}}\label{sec_proof_of_theorem_que_along_flow}

	Taking expectations in \eqref{integrated_cL_cK_evolution}, we obtain that for $1-W/N\le s\le t\le \tf$,
	\begin{equation}\label{E_cL_cK_evolution_integrated}
		\begin{aligned}
			\E(\cL-\cK)_{t,xy}^{\bsigma}=\qa{\cU_{s,t}^{\bsigma}\circ\E(\cL-\cK)_{s}^{\bsigma}}_{xy}+\int_{s}^{t}\qb{\cU_{u,t}^{\bsigma}\circ\E\cW_{u}^{\cL,\bsigma}}_{xy}\, \rd u+\int_{s}^{t}\qb{\cU_{u,t}^{\bsigma}\circ\E\mathcal{E}_{u}^{\cL,\bsigma}}_{xy}\, \rd u.
		\end{aligned}
	\end{equation}
	To exploit the refined evolution kernel bound \eqref{bound_cU_s_t_2} in the case $\sigma_1\neq\sigma_2$, we apply $\cO$ to both sides of \eqref{E_cL_cK_evolution_integrated}. Since $\cQ$ commutes with $\cU_{u,t}^{\bsigma}$, we obtain
	\begin{equation}\label{cO_E_cL_cK_evolution_integrated}
		\begin{aligned}
			\E\qa{\cO\circ(\cL-\cK)_{t}^{\bsigma}}_{xy}&=\qa{\cU_{s,t}^{\bsigma}\circ\E\qa{\cO\circ(\cL-\cK)_{s}^{\bsigma}}}_{xy}\\
			&+\int_{s}^{t}\qb{\cU_{u,t}^{\bsigma}\circ\E\p{\cO\circ\cW_u^{\cL,\bsigma}}}_{xy}\, \rd u+\int_{s}^{t}\qb{\cU_{u,t}^{\bsigma}\circ\E\p{\cO\circ\mathcal{E}_u^{\cL,\bsigma}}}_{xy}\, \rd u.
		\end{aligned}
	\end{equation}
	We analyze \eqref{E_cL_cK_evolution_integrated} and \eqref{cO_E_cL_cK_evolution_integrated} using \Cref{lemma_weak_local_law_T_along_flow,lemma_flat_local_law_along_flow}, along with the following local law in expectation.
	
	\begin{lemma}\label{lemma_expected_single_resolvent_average_local_law}
		In the setting of \Cref{theorem_que_along_flow}, for any $t\in[1-W/N,\tf]$, we have 
		\begin{equation*}
			\begin{aligned}            \max_{x\in\ZL}\absa{\E\avg{(G_t -m)S^{\pa{x}}}} \prec \qa{N(1-t)}^{-2}.
			\end{aligned}
		\end{equation*}
	\end{lemma}
	\begin{proof}
		Applying \eqref{eq:Gu-m} to $\avg{(G_t-m)S^{\pa{x}}}$ and using Gaussian integration by parts, we obtain
		\begin{equation*}
			\begin{aligned}
				\sum_{y}(1-tm^2S)_{xy}\E\avg{(G_t-m)S^{\pa{y}}}=tm\sum_{a,b}\E S_{ab}S_{xa}(G_{t}-m)_{bb}(G_t-m)_{aa},\quad \forall x\in\ZL.
			\end{aligned}
		\end{equation*}
		Applying \eqref{eq:Gu-m} once more to $(G_t-m)_{aa}$ and performing Gaussian integration by parts again yields
		\begin{equation*}
			\begin{aligned}
				\sum_{y}(1-tm^2S)_{xy}\E\avg{(G_t-m)S^{\pa{y}}}&=t^2m^2\sum_{a,b,c}\E S_{ab}S_{xa}S_{ca}(G_{t})_{ba}(G_t)_{ac}(G_t)_{cb}\\
				&+t^2m^2\sum_{a}\E S_{xa} (G_t)_{aa}\qb{\avg{(G_t-m)S^{\pa{a}}}}^2.
			\end{aligned}
		\end{equation*}
		The second term on the RHS is bounded by $\opr{\q{N(1-t)}^{-2}}$ using the local laws in \eqref{flow_local_law}. For the first term, we apply the $GG$-expansion from \eqref{out_GG_expansion} to $(G_t)_{ba}(G_t)_{ac}$, which produces a finite sum of graphs. Each term can be estimated individually by $\opr{\q{N(1-t)}^{-2}}$ using again the local laws in \eqref{flow_local_law}. We omit the routine details.
		Finally, invoking the bound $\|(1-tm^2S)^{-1}\|_{\max\to\max}\prec 1$, which follows from \eqref{Theta_bound_equaled_charge}, we conclude the proof.
	\end{proof}

	The light-weight term in \eqref{E_cL_cK_evolution_integrated} can be estimated using the same argument as in \cite[Lemma 6.2]{dubova2025delocalizationnonmeanfieldrandommatrices}.
	
	\begin{lemma}\label{lemma_expected_light_weight_term_bound}
		In the setting of \Cref{theorem_que_along_flow}, for any $t\in[1-W/N,\tf]$, we have 
		\begin{equation}\label{eq:Exp_LW}            \max_{\bsigma\in\h{-,+}^2}\max_{x,y\in\ZL}\absa{\E\cW_{t,xy}^{\cL,\bsigma}} \prec \frac{1}{1-t}\frac{1}{N(1-t)}\qa{\frac{1}{W^{6/5}}+ \frac{1}{\q{N(1-t)}^{3/2}}}.
		\end{equation}
	\end{lemma}
	\begin{proof}
		The proof follows the same argument as that of Lemma 6.2 in \cite{dubova2025delocalizationnonmeanfieldrandommatrices}. The only differences are that Lemma 6.1 therein is replaced by our \Cref{lemma_expected_single_resolvent_average_local_law}, and the estimate for $(\cL-\cK)$ used there is replaced by our bounds \eqref{flat_local_law_along_flow} and \eqref{flat_local_law_along_flow2}. Therefore, we omit the details.
	\end{proof}

	Finally, using \eqref{flat_local_law_along_flow}, we bound the quadratic error term as follows.
	\begin{lemma}\label{lemma_expected_double_difference_bound}
		In the setting of \Cref{theorem_que_along_flow}, for any $t\in[1-W/N,\tf]$, we have
		\begin{equation}\label{eq:Exp_quadraticerror}            \max_{\bsigma\in\h{-,+}^2}\max_{x,y\in\ZL}\absa{\E\cE_{t,xy}^{\cL,\bsigma}} \prec \frac{1}{1-t}\qa{\frac{1}{W^{6/5}}+\frac{1}{\q{N(1-t)}^{9/4}}}.
		\end{equation}
	\end{lemma}
	\begin{proof}
		For notational simplicity, we assume $\bsigma=(-,+)$; the general case is identical up to minor notational changes. Proceeding as in the proof of \Cref{lemma_double_difference_decomposition} in \Cref{sec_estimates_error_terms}, and using \eqref{eq:Gu-m} together with a Gaussian integration by parts argument analogous (and slightly simpler) to that below \eqref{eq_Xi_H}, we obtain
		\begin{align}
			&\E\cE_{t,xy}^{\cL,\bsigma}=t\E\sum_{j}(\cL-\cK)_{t,xj}(\cL-\cK)_{t,jy}+tm\E\sum_{j,\beta}(G-m)_{\beta\beta}S_{\beta j}T_{t,xj}(\cL-\cK)_{t,jy}\nonumber\\
			&+tm\E\sum_{j}(G_t^*-\ol{m})_{jj}\cL_{t,xj}(\cL-\cK)_{t,jy}+m\E\sum_{i}S_{xi}(G_t^*-\ol{m})_{ii}(\cL-\cK)_{t,iy}\nonumber\\
			&-tm\E\sum_{i,j,\beta}S_{xi}S_{\beta j}(G_t)_{i\beta}(G_t^*)_{ji}\partial_{j\beta}\p{\cL_{t,jy}}.\label{eq:quaderror_E}
		\end{align}
		Using \eqref{flat_local_law_along_flow} and the fact that $W^{-1}\le [N(1-t)]^{-1}$, the first term on the RHS is bounded by
		\begin{equation*}
			\begin{aligned}
				\OO_\prec \pa{ \frac{1}{1-t}\qa{\frac{1}{W^{7/5}}+\frac{1}{\q{N(1-t)}^{5/2}}}}.
			\end{aligned}
		\end{equation*}
		Applying \eqref{flat_local_law_along_flow} together with the local law \eqref{flow_local_law}, the second, third, and fourth terms are bounded by
		\begin{equation*}
			\begin{aligned}
				\OO_\prec \pa{ \frac{1}{1-t} \frac{1}{\sqrt{N(1-t)}} \qa{\frac{1}{W^{6/5}}+\frac{1}{\q{N(1-t)}^{7/4}}}}.
			\end{aligned}
		\end{equation*}
		For the last term in \eqref{eq:quaderror_E}, we expand $\partial_{j\beta}(\cL_{t,jy})$ and apply the local law \eqref{flow_local_law} together with the bound on off-diagonal $T$-variables from \eqref{flat_local_law_along_flow2}. This yields
		\begin{align*}
			\E\sum_{i,j,\beta}S_{xi}S_{\beta j}(G_t)_{i\beta}(G_t^*)_{ji}\partial_{j\beta}\p{\cL_{t,jy}}&\prec  \frac{N}{\sqrt{N(1-t)}}\qa{\frac{1}{W^{6/5}}+\frac{1}{\q{N(1-t)}^{3/2}} }^2 \\
            &\lesssim \frac{\p{1-t}^{-1}}{\sqrt{N(1-t)}}\qa{\frac{1}{W^{7/5}}+\frac{1}{\q{N(1-t)}^{2}}}.
		\end{align*}
		Combining the above bounds yields \eqref{eq:Exp_quadraticerror}.
	\end{proof}

	With the above inputs, we are ready to prove the estimate \eqref{exp_flat_local_law}.
	When $\sigma_1=\sigma_2$, we insert the initial condition \eqref{eq:initialW/N} at time $s=1-W/N$, together with the bounds \eqref{eq:Exp_LW} and \eqref{eq:Exp_quadraticerror}, into \eqref{E_cL_cK_evolution_integrated}. Using the evolution kernel estimate \eqref{bound_cU_s_t_3} for case (i) and performing the integral over $u$, we obtain \eqref{exp_flat_local_law}.
	Next, consider the case $\sigma_1\ne\sigma_2$. In this case we work with \eqref{cO_E_cL_cK_evolution_integrated}. By \Cref{lemma_expected_single_resolvent_average_local_law}, after taking expectation the estimate \eqref{eq:cL-cK_flat} improves to
	\begin{align}
		\E \qa{(\cL-\cK)_{t}^{\bsigma}-\cO\circ(\cL-\cK)_{t}^{\bsigma}}_{xy}     &\prec {\q{N(1-t)}^{-3}} .\label{eq:cL-cK_flat_exp}
	\end{align}
	Moreover, by \Cref{def_sum_zero_operators} we have $\|\cO\circ X\|_{\max}\lesssim \|X\|_{\max}$. Hence, the quantities $\cO\circ\E(\cL-\cK)_{s}^{\bsigma}$ with $s=1-W/N$, $\cO\circ\E \cW_u^{\cL,\bsigma} $, and $\cO\circ\E \mathcal{E}_u^{\cL,\bsigma}$ satisfy the same bounds as in \eqref{eq:initialW/N}, \eqref{eq:Exp_LW}, and \eqref{eq:Exp_quadraticerror}. Plugging these bounds into \eqref{cO_E_cL_cK_evolution_integrated}, applying the evolution kernel estimate \eqref{bound_cU_s_t_3} for case (ii), and performing the integral over $u$, we obtain
	\begin{equation}\label{exp_flat_local_law_cO}        \max_{\bsigma\in\h{(-,+),(+,-)}}\max_{x,y\in\ZL}\absa{\E \qa{\cO\circ (\cL-\cK)_{t}^{\bsigma}}_{xy}}\prec W^{-6/5}+\pa{N(1-t)}^{-9/4}.
	\end{equation}
	Combining \eqref{exp_flat_local_law_cO} with \eqref{eq:cL-cK_flat_exp} yields \eqref{exp_flat_local_law} for the case $\sigma_1\ne\sigma_2$.

	\section{Proof of the main results for \texorpdfstring{$\alpha\in (-1,0)$}{alpha in (-1,0)}}\label{sec_proof_of_remaining_cases_2}

	In this section, we present the proof of the main results for $\alpha\in(-1,0)$. For consistency of notation, we again adopt the flow framework and notation introduced in \Cref{sec:notations}. In fact, the dynamical analysis can be replaced by a standard bootstrap argument along a sequence of spectral parameters whose imaginary parts decrease multiplicatively. Such an approach also extends to general entry distributions under certain moment assumptions; see, e.g., \cite{QuadraticAOEK,ajanki2017universality,EKYYEJP18}. However, for simplicity of presentation, we still decrease the imaginary part of the spectral parameter along the flow  \eqref{eq:zt}, and assume that the entries of $H$ are Gaussian as in \eqref{bandcw0}. The main results then follow directly from the following estimates established along the flow.

	\begin{theorem}\label{local_laws_along_flow_subcritical}
		Under the assumptions of \Cref{theorem_delocalization}, fix $\alpha\in (-1,0)$, small constants $\kappa,\fc\in(0,1)$, and a spectral parameter $z=E+\ii \eta\in\mathbf{D}_{\kappa,\fc}$. Consider the flow framework in \Cref{def_H_t_flow,def_z_t_flow}, with the flow parameter chosen as in \eqref{eq:t0E0}. Then the following estimates hold.
		\begin{itemize}
			\item {\bf Local law:} For any $x,y\in\ZL$, the following entrywise and averaged local laws hold uniformly for $t\in[0,\tf]$:
			\begin{equation}\label{flow_local_law_subcritical}
				\begin{aligned}
					\absa{(G_t-m)_{xy}}^2\prec B_t(|x-y|),\quad \big|\avg{(G_t-m)S^{\pa{x}}}\big|\prec B_t\pa{0},
				\end{aligned}
			\end{equation}
			where we recall the shape parameter $B_t(\cdot )=B(1-t,\cdot)$ defined in \eqref{def_B_t}.
			
			\item {\bf $\cal L$-loop and $T$-variable estimates:} For any $\bsigma\in\ha{-,+}^2$ and $x,y,y'\in\ZL$, we have uniformly for $t\in[1-(W/N)^{1+\alpha},\tf]$ that
			\begin{equation}\label{flow_loop_local_law_subcritical}
				\begin{aligned}
					(T_{t}-\Theta_{t})_{x,yy'}^{\bsigma}\prec \qa{N(1-t)}^{-3/2},\quad(\cL_t-\cK_t)_{xy}^{\bsigma}\prec \qa{N(1-t)}^{-7/4}.
				\end{aligned}
			\end{equation}
			
			\item {\bf Expected $\cal L$-loop estimate:} For any $\bsigma\in\ha{-,+}^2$ and $x,y\in\ZL$, we have uniformly for $t\in[1-(W/N)^{1+\alpha},\tf]$ that       
			\begin{equation}\label{flow_expected_local_law_subcritical}
				\max_{\bsigma\in\h{-,+}^2}\max_{x,y\in\ZL}\absa{\E(\cL-\cK)_{t,xy}^{\bsigma}}\prec W^{-7/4}\pa{{N}/{W}}^{7\alpha/4}+\qa{{N(1-t)}}^{-9/4}.
			\end{equation}
		\end{itemize}
		
	\end{theorem}

	\begin{proof}[\bf Proof of \Cref{theorem_local_law_entrywise,theorem_local_law_loops} and \Cref{theorem_local_law_T} for $\al\in (-1,0)$]
		By \eqref{eq:zztE}, we have \smash{$G(z)\stackrel{d}{=}\sqrt{\tf}\,G_{\tf,\sE}$}. Hence, by \eqref{eq:Betar}, for each fixed $z\in\mathbf{D}_{\kappa,\fc}$ (or $z\in\mathbf{D}_{\kappa,\fc}^{\txt{flat}}$) the estimates obtained along the flow directly imply the corresponding bounds for the original resolvent. More precisely, \eqref{flow_local_law_subcritical} yields the entrywise and averaged local laws \eqref{eq_entrywise_local_law} and \eqref{eq_average_local_law}, \eqref{flow_loop_local_law_subcritical} implies \eqref{eq_loop_local_law_3} and \eqref{eq_chain_local_law_3}, and \eqref{flow_expected_local_law_subcritical} gives the expected bound \eqref{eq_expected_loop_local_law}. Finally, a standard $N^{-C}$-net argument extends these estimates uniformly to all $z\in\mathbf{D}_{\kappa,\fc}$ (or $z\in\mathbf{D}_{\kappa,\fc}^{\txt{flat}}$).
	\end{proof}
	
	\subsection{Proof of \Cref{local_laws_along_flow_subcritical}}

	We first establish the analogues of \eqref{Theta_bound_opposite_charge}, \eqref{Theta_bound_equaled_charge}, and \eqref{zero_mode_removed_bound} for the regime $\alpha\in(-1,0)$. These estimates are collected in the following lemma, whose proof is deferred to \Cref{sec_proof_of_lemma_input_bound_subcritical}.
	
	\begin{lemma}\label{lemma_input_bound_subcritical}
		Consider a power-law random band matrix $H$ from \Cref{considered_model} with variance profile matrix $S$ satisfying \eqref{alpha_decay} for a fixed parameter $\alpha\in(-1,0)$. Let $\sE$ be a flow parameter satisfying $|\sE|\le 2-\e$ for some constant $\e>0$. Then, under \Cref{assumption_input_bound<0}, the $\Theta$-propagators defined in \Cref{def_loops_and_T_variables} satisfy the following uniform bounds for $0\le t\le 1-N^{-1}$, with a constant $C_\alpha$ depending only on $\alpha$, $c_\al$, and $\e$:
		\begin{itemize}
			\item {\bf Upper bounds:} For $\bsigma=(\sigma_1,\sigma_2)$ with $\sigma_1\neq \sigma_2$, we have        \begin{equation}\label{Theta_upper_bound_opposite_charge_subcritical}
				\begin{aligned}
					\Theta_{t,xy}^{\bsigma}\prec B_t(|x-y|),\quad  \forall x,y\in\ZL,
				\end{aligned}
			\end{equation}
			while for $\bsigma=(\sigma_1,\sigma_2)$ with $\sigma_1=\sigma_2$, we have       \begin{equation}\label{Theta_upper_bound_equal_charge_subcritical}
				\begin{aligned}
					\Theta_{t,xy}^{\bsigma}\prec \cS_{xy}\asymp B_0(|x-y|),\quad  \forall x,y\in\ZL.
				\end{aligned}
			\end{equation}
			\item {\bf Zero-mode-removed upper bound:} For $\bsigma=(\sigma_1,\sigma_2)$ with $\sigma_1\neq \sigma_2$, we have         \begin{equation}\label{zero_mode_removed_bound_subcritical}
				\begin{aligned}
					\absa{\Theta_{t,xy}^{\bsigma}-\frac{1}{N(1-t)}}\prec  \cS_{xy}\asymp B_0(|x-y|),\quad  \forall x,y\in\ZL.
				\end{aligned}
			\end{equation}
		\end{itemize}
	\end{lemma}

    The key ingredient in the proof of this lemma is the following mixing property of the Markov chain with transition matrix \(S\), whose proof is direct and is therefore also deferred to \Cref{sec_proof_of_lemma_input_bound_subcritical}. 
	
	\begin{lemma}\label{lemma_subcritical_S_mix}
		Define $k_{\txt{mix}}:=\inf\ha{k\in \N: 1+k\alpha\leq 0}$. Then there exists a constant $C_\al>0$, depending only on $\al$, such that the $k$-step probabilities of the Markov chain with transition matrix $S$ satisfy
		\begin{equation}\label{bound_S_mix}
			S^{k}_{xy}\le \begin{cases}
				C_\al W^{-1}(N/W)^{k\alpha}(|x-y|/W+1)^{-1-k\alpha},\ & \text{if}\ k<k_{\txt{mix}},\\
				C_\al N^{-1}\log(N/W+1),\ & \text{if}\ k\geq k_{\txt{mix}}.
			\end{cases}
		\end{equation}
		Moreover, there exists a constant $c_\al>0$, depending only on $\al$, such that
		\begin{equation}\label{S_mix_exponential}
			\begin{aligned}
				\absa{S_{xy}^k-{N}^{-1}}\le c_\al^{-1} \exp(-c_\al k),\quad \forall k\in \N.
			\end{aligned}
		\end{equation}
	\end{lemma}
	
	Due to the slow decay of $S$ in the regime $\alpha\in(-1,0)$, \Cref{2_loop_to_1_chain} takes the following modified form.
	
	\begin{lemma}\label{2_loop_to_1_chain_subcritical}
		In the setting of \Cref{local_laws_along_flow_subcritical}, the bounds \eqref{upper_bound_2_to_1} and \eqref{entrywise_bound_2_to_1} remain valid with $W^{-c}$ in the definition of $\Omega_{u,c}$ there replaced by $N^{-c}$. Furthermore, suppose that $\|G_u-m\|_{\max}\prec N^{-c}$ for some constant $c>0$, and \smash{$\|\cS\|_{\max}+\|\cL_{u}\|_{\max}\prec \Psi_{u}^2$}
		for a deterministic control parameter $\Psi_u$ satisfying $0< \Psi_u\leq N^{-\varepsilon}$, where $\varepsilon>0$ is a constant. Then we have the averaged local law
		\begin{equation}\label{average_bound_2_to_1_subcritical}
			\begin{aligned}
				\big|\avgb{(G_u-m)\cal D}\big|\prec \Psi_u^2
			\end{aligned}        
		\end{equation}
		for any deterministic diagonal matrix $\cal D$ satisfying $\|\cal D\|_{\max}\prec \Psi_u^2$ and $\sum_i |\cal D_{ii}|\prec 1$.     
	\end{lemma}
	
	\begin{proof}
		This statement was established as Lemma 4.1 in \cite{Band1D} for 1D RBMs, and we omit the proof here.
	\end{proof}

	The continuity estimate in \Cref{lemma_continuity_estimate} remains valid in the regime $\alpha\in(-1,0)$. With these preparations, \Cref{local_laws_along_flow_subcritical} follows directly from the following lemma.

	\begin{lemma}\label{main_lemma_subcritical}
		In the setting of \Cref{local_laws_along_flow_subcritical}, fix a time $t\in[0,\tf]$ and assume that $\|G_t-m\|_{\max}\prec N^{-c}$ for some constant $c>0$.
		Then, at time $t$, the upper bound
		\begin{equation}\label{upper_bound_T_subcritial}
			\begin{aligned}
				\max_{\bsigma\in\ha{-,+}^2} \absa{T_{t,xy}^{\bsigma}}\prec B_t(|x-y|)
			\end{aligned}
		\end{equation}
		and the local laws in \eqref{flow_local_law_subcritical} hold. Moreover, if $t\in[1-(W/N)^{1+\alpha},\tf]$, the estimate \eqref{flow_loop_local_law_subcritical} also holds at time $t$.

	\end{lemma}
	\begin{proof}[\bf Proof of \Cref{local_laws_along_flow_subcritical}]
		We proceed by induction in time $t\in[0,\tf]$. At time $t=0$, the estimate \eqref{flow_local_law_subcritical} holds trivially since $G_0=m$. Suppose that \eqref{flow_local_law_subcritical} holds at some time $s\in[0,\tf]$, and let $t\in[s,\tf]$ satisfy
		\begin{equation*}
			\begin{aligned}
				\frac{1-s}{1-t}B_t^{1/10}(0)\leq 1.
			\end{aligned}
		\end{equation*}
		By \Cref{main_lemma_subcritical}, we obtain the upper bound \eqref{upper_bound_T_subcritial} at time $s$. Then, using Lemmas \ref{2_loop_to_1_chain_subcritical} and \ref{lemma_continuity_estimate}, we proceed as in Step~1 of the proof of \Cref{lemma_bootstrap} (see the arguments in Section~5.1 of \cite{Band1D}) to obtain   
		\begin{equation*}
			\begin{aligned}
				\|G_t-m\|_{\max}\prec B_t^{1/5}(0).
			\end{aligned}
		\end{equation*}
		This estimate implies, by \Cref{main_lemma_subcritical}, that \eqref{flow_local_law_subcritical} and \eqref{flow_loop_local_law_subcritical} also hold at time $t$. Here we also use that
		\begin{equation*}
			\begin{aligned}
				\absa{\cL_{s,xy}^{\bsigma}}\leq \sum_{a\in\ZL}\absa{T_{s,xa}^{\bsigma}}S_{ay}\prec B_s(|x-y|).
			\end{aligned}
		\end{equation*}
		Since $B_u(0)\le B_{\tf}(0)\le N^{-c}$ for some constant $c>0$, the above argument can be iterated $\OO(1)$ times, yielding \eqref{flow_local_law_subcritical} and \eqref{flow_loop_local_law_subcritical} in \Cref{local_laws_along_flow_subcritical}.
		
		Finally, using \eqref{flow_local_law_subcritical} and \eqref{flow_loop_local_law_subcritical}, and applying the methods used in the proofs of Lemmas \ref{lemma_expected_single_resolvent_average_local_law}--\ref{lemma_expected_double_difference_bound}, we obtain the following bounds for $t\in[1-(W/N)^{1+\alpha},\tf]$:
		\begin{equation*}
			\max_{x\in\ZL}\absa{\E\avg{(G_t-m)S^{\pa{x}}}}\prec \qa{N(1-t)}^{-2},
		\end{equation*}
		\begin{equation*}        \max_{\bsigma\in\h{-,+}^2}\max_{x,y\in\ZL}\absa{\E\cW_{t,xy}^{\cL,\bsigma}}\prec \frac{1}{1-t}\frac{1}{\q{N(1-t)}^{5/2}},\quad\max_{\bsigma\in\h{-,+}^2}\max_{x,y\in\ZL}\absa{\E\cE_{t,xy}^{\cL,\bsigma}}\prec \frac{1}{1-t}\frac{1}{\q{N(1-t)}^{9/4}}.
		\end{equation*}
		Together with \eqref{flow_local_law_subcritical} and \eqref{flow_loop_local_law_subcritical}, these bounds yield \eqref{flow_expected_local_law_subcritical} by the same argument as in the proof of \Cref{theorem_que_along_flow} in \Cref{sec_proof_of_theorem_que_along_flow}. We omit the details.
	\end{proof}

	\subsection{Proof of \Cref{main_lemma_subcritical}}
	
	The key to the proof of \Cref{main_lemma_subcritical} is that sharp upper bounds on the $\cL$-loops and $T$-variables can be derived from the weak input $\|G_t-m\|_{\max}\prec N^{-c}$ as follows. 
	More precisely, for any $\bsigma\in\h{-,+}^2$, $x,y\in\ZL$, and $t\in[0,\tf]$, we can bound the $T$-variables and $\cL$-loops by
	\begin{align}
		T_{t,xy}^{\bsigma}&\prec \frac{1}{W}\pa{\frac{N}{W}}^{\alpha}\pa{\frac{|x-y|}{W}+1}^{-1-\alpha}+\frac{1}{N(1-t)}= B_t(|x-y|),\label{subcritical_bound_T_and_loops1}\\
		\cL_{t,xy}^{\bsigma}&\prec \frac{1}{W}\pa{\frac{N}{W}}^{2\alpha}\pa{\frac{|x-y|}{W}+1}^{-1-2\alpha}\cdot\mathbf{1}_{1+2\alpha>0}+\frac{1}{N(1-t)}=:B_t^{\cL}(|x-y|).\label{subcritical_bound_T_and_loops2}
	\end{align}
	To derive the bound \eqref{subcritical_bound_T_and_loops1}, we apply \eqref{upper_bound_2_to_1} iteratively. For $\bsigma=(-,+)$, we obtain
	\begin{align*}
		T_{t,xy}^{\bsigma}&=\sum_{a\in\ZL}S_{xa}|(G_t)_{ay}|^2\prec \cS_{xy}+\sum_{a\in\ZL}S^2_{xa}|(G_t)_{ay}|^2\prec \cS_{xy}+ \cS^2_{xy}+\sum_{a\in\ZL}S^3_{xa}|(G_t)_{ay}|^2 \\
		&\prec\cdots\prec \sum_{1\leq k<k_{\txt{mix}}}\cS^{k}_{xy}+\sum_{a\in\ZL}S^{k_{\txt{mix}}}_{xa}|(G_t)_{ay}|^2.
	\end{align*}
Using the bound \eqref{bound_S_mix} and Ward's identity, we obtain
	\begin{equation*}
		\begin{aligned}
			T_{t,xy}^{\bsigma}&\prec \sum_{1\leq k<k_{\txt{mix}}}\frac{1}{W}\pa{\frac{N}{W}}^{k\alpha}\pa{\frac{|x-y|}{W}+1}^{-1-k\alpha}+\frac{1}{N}\sum_{a\in\ZL}|(G_t)_{ay}|^2\\
			&\prec \frac{1}{W}\pa{\frac{N}{W}}^{\alpha}\pa{\frac{|x-y|}{W}+1}^{-1-\alpha}+\frac{1}{N(1-t)}.
		\end{aligned}
	\end{equation*}
	The proof of \eqref{subcritical_bound_T_and_loops1} for the case $\bsigma=(+,-)$ is analogous, and the general case $\bsigma\in\ha{-,+}^2$ follows immediately from the cases $\bsigma\in\{(-,+),(+,-)\}$ by the CS inequality. The bound \eqref{subcritical_bound_T_and_loops2} can be established by a similar argument. Finally, combining \eqref{subcritical_bound_T_and_loops1} with \Cref{2_loop_to_1_chain_subcritical} yields the local laws in \eqref{flow_local_law_subcritical}.

	It remains to prove the estimates in \eqref{flow_loop_local_law_subcritical} for $t\in[1-(W/N)^{1+\alpha},\tf]$. In this regime, the bounds \eqref{subcritical_bound_T_and_loops1} and \eqref{subcritical_bound_T_and_loops2} simplify to
	\begin{equation}\label{subcritical_bound_T_and_loops}
		\max_{\bsigma\in\ha{-,+}^2}\|T_{t}^{\bsigma}\|_{\max}+\max_{\bsigma\in\ha{-,+}^2}\|\cL_{t}^{\bsigma}\|_{\max}\prec \qa{N(1-t)}^{-1}.
	\end{equation}
	We first prove the $T$-variable estimate by establishing the following self-improving bound. Suppose that for some deterministic control parameter $\cJ$ satisfying $ N^{-3/2}(1-t)^{-3/2}\leq \cJ\leq N^{-1/2}(1-t)^{-1/2}$, we have 
	\begin{equation}\label{eq:assm_boundT}
		\max_{\bsigma\in\ha{-,+}^2}\max_{x,y,y'\in\ZL} \absa{(T-\Theta)_{t,x,yy'}^{\bsigma}}\prec \cJ . \end{equation}
	Then the following improved bound holds:
	\begin{equation}\label{self_improved_T_local_law_subcritical}
		\begin{aligned}
			\max_{\bsigma\in\ha{-,+}^2}\max_{x,y,y'\in\ZL} \absa{(T-\Theta)_{t,x,yy'}^{\bsigma}}\prec  \frac{1}{\q{N(1-t)}^{3/2}}+\frac{\cJ^{1/2}}{\q{N(1-t)}^{3/4}}+ \frac{\cJ}{\q{N(1-t)}^{1/2}}.
		\end{aligned}
	\end{equation}
	Iterating this estimate $\OO(1)$ times yields the first bound in \eqref{flow_loop_local_law_subcritical}. To prove \eqref{self_improved_T_local_law_subcritical}, for any $\bsigma=(\sigma_1,\sigma_2)\in\h{-,+}^2$, we introduce
	\begin{equation*}
		\begin{aligned}
			\Xi_{t,x,yy'}^{T,\bsigma}:=\qb{(1-tm(\sigma_1)m(\sigma_2)S)\cdot (T-\Theta)_{t,\cdot,yy'}^{\bsigma}}_{x}.
		\end{aligned}
	\end{equation*}
	Using a Gaussian integration-by-parts argument similar to that used for \eqref{high_moment_decomposition_double_difference}, together with the assumption \eqref{eq:assm_boundT} and the local laws in \eqref{flow_local_law_subcritical}, we obtain that for any $\bsigma\in\h{-,+}^2$ and $x,y,y'\in\ZL$,
	\begin{equation}\label{bound_Xi_subcritial_T}
		\begin{aligned}
			\Xi_{t,x,yy'}^{T,\bsigma}\prec \frac{1}{\q{N(1-t)}^{3/2}}+\frac{\cJ^{1/2}}{\q{N(1-t)}^{3/4}}+ \frac{\cJ}{\q{N(1-t)}^{1/2}}.
		\end{aligned}
	\end{equation}
	Since the derivation is analogous, and in fact simpler than that of \eqref{high_moment_decomposition_double_difference}, we omit the details.
	If $\sigma_1=\sigma_2$, applying $(1-tm(\sigma_1)m(\sigma_2)S)^{-1}$ to both sides of \eqref{bound_Xi_subcritial_T} and using \eqref{Theta_upper_bound_equal_charge_subcritical} yields the desired bound \eqref{self_improved_T_local_law_subcritical}.
	If $\sigma_1\ne\sigma_2$, using the operators $\cP$ and $\cQ$ introduced in \Cref{def_sum_zero_operators}, we write
	\begin{align}
		\qb{\cQ\cdot(T-\Theta)_{t,\cdot,yy'}^{\bsigma}}_{x}&=\qb{\cQ\cdot\pa{1-tS}^{-1}\cdot \Xi_{t,\cdot,yy'}^{T,\bsigma}}_{x}=\qb{\cQ\cdot \Xi_{t,\cdot,yy'}^{T,\bsigma}}_{x}+t\qb{\cQ\cdot\Theta_{t}^{\bsigma}\cdot \Xi_{t,\cdot,yy'}^{T,\bsigma}}_{x}\nonumber\\
		&\prec \frac{1}{\q{N(1-t)}^{3/2}}+\frac{\cJ^{1/2}}{\q{N(1-t)}^{3/4}}+ \frac{\cJ}{\q{N(1-t)}^{1/2}},\label{eq:QT-Theta}     \end{align}
	where in the second step we used \eqref{bound_Xi_subcritial_T} together with \eqref{zero_mode_removed_bound_subcritical} for $\cQ\cdot\Theta_{t}^{\bsigma}$. On the other hand, using Ward's identity \eqref{eq_Ward} and the local law \eqref{flow_local_law_subcritical}, we obtain 
	\begin{equation}\label{eq:PT-Theta}
		\qb{\cP\cdot(T-\Theta)_{t,\cdot,yy'}^{\bsigma}}_{x}\prec \qa{{N(1-t)}}^{-3/2}.
	\end{equation}
	Combining \eqref{eq:QT-Theta} and \eqref{eq:PT-Theta}, and using the identity $\cP+\cQ=1$, yields \eqref{self_improved_T_local_law_subcritical}.

	The second estimate in \eqref{flow_loop_local_law_subcritical} for the $\cL$-loops follows from a similar argument. More precisely, arguing analogously to the above derivation for the $T$-variables, it suffices to prove that
	\begin{equation}\label{bound_Xi_subcritial_cL}
		\begin{aligned}
			\Xi_{t,xy}^{\cL,\bsigma}:=\qa{(1-tm\p{\sigma_1}m\p{\sigma_2}S)\cdot (\cL-\cK)_{t}^{\bsigma}}_{xy}\prec \qa{{N(1-t)}}^{-7/4}.
		\end{aligned}
	\end{equation}
	Without loss of generality, we assume $\bsigma=(-,+)$. For notational simplicity, we suppress the subscript $t$ as well as the superscripts $\cL$ and $\bsigma$. Employing a Gaussian integration-by-parts argument similar to that used for \eqref{high_moment_decomposition_double_difference}, we obtain for any fixed $p\in\N$ that 
	\begin{align}
		&\E|\Xi_{xy}|^{2p}=tm\E\sum_{i,a\in\ZL}S_{xi}S_{ia}(G-m)_{aa}T_{yi}^{(+,-)}\cdot \Xi_{xy}^{p-1}\ol{\Xi}_{xy}^{p}\nonumber\\
		&+tm\E\sum_{i\in\ZL}S_{xi}(G^*-\ol{m})_{ii}(\cL-\cK)_{iy} \cdot \Xi_{xy}^{p-1}\ol{\Xi}_{xy}^{p}+m\E\sum_{i\in\ZL}S_{xi}(G^*-\ol{m})_{ii}\Theta_{t,iy}\cdot \Xi_{xy}^{p-1}\ol{\Xi}_{xy}^{p}\nonumber\\
		&-tm\E\sum_{i,j,a\in\ZL}S_{xi}S_{ia}S_{jy}G_{a j}G_{ji}^*\cdot \partial_{a i}\pb{\Xi_{xy}^{p-1}\ol{\Xi}_{xy}^{p}}.\label{eq:boundL-K2p}
	\end{align}
	Using the local laws in \eqref{flow_local_law_subcritical}, the bound \eqref{subcritical_bound_T_and_loops}, and the first estimate in \eqref{flow_loop_local_law_subcritical} for the $T$-variables (whose diagonal version also provides an a priori bound on $\cL$-loops), we can bound the first, second, and fourth terms on the RHS by
	\begin{equation}
		\qa{{N(1-t)}}^{-2}\cdot\pa{\E|\Xi_{xy}|^{2p}}^{(2p-1)/2p}+\qa{{N(1-t)}}^{-7/2}\cdot\pa{\E|\Xi_{xy}|^{2p}}^{(2p-2)/2p}.\label{eq;Xip1}
	\end{equation}
	For the third term on the RHS of \eqref{eq:boundL-K2p}, we apply the averaged local law \eqref{average_bound_2_to_1_subcritical} with the choice $\cal D_{ii}=S_{xi} \Theta_{t,iy}$. Using also the bound $\|\Theta_t\|_{\max}\prec (N(1-t))^{-1}$ from \eqref{Theta_upper_bound_opposite_charge_subcritical}, we obtain 
	\begin{equation*}
		\sum_{i\in\ZL}S_{xi}(G^*-\ol{m})_{ii}\Theta_{t,iy}\prec \qa{{N(1-t)}}^{-2}.
	\end{equation*}
	Applying Hölder's inequality, the third term in \eqref{eq:boundL-K2p} can therefore be bounded by
	\begin{equation}\label{eq;Xip2}
		\qa{{N(1-t)}}^{-2}\cdot\pa{\E|\Xi_{xy}|^{2p}}^{(2p-1)/2p}.
	\end{equation}
	Substituting \eqref{eq;Xip1} and \eqref{eq;Xip2} into \eqref{eq:boundL-K2p}, and applying Young's and Markov's inequalities, yields \eqref{bound_Xi_subcritial_cL}.

	\section{Proof of the main results for \texorpdfstring{$\al=0$}{alpha0}}\label{sec_proof_alpha_0}

	In this critical case, the methods developed for the regimes $\alpha\in(0,1)$ and $\alpha\in(-1,0)$ both fail. More precisely, the approach for $\alpha\in(-1,0)$ breaks down because the bound \eqref{bound_S_mix} no longer holds. The method for $\alpha\in(0,1)$ fails at the level of Lemma \ref{bound_leading_term}, where the bound \smash{$\OO((1-u)^{-1}J_u^{(\ell)})$} must be replaced by \smash{$\OO_{\prec}((1-u)^{-1}J_u^{(\ell)})$}. This modification arises because in the case $\alpha=0$ we assume only $C_{\alpha}\prec 1$ in Assumption \ref{assumption_input_bound}. As a consequence, the Grönwall argument in \Cref{sec:Gronwall} is no longer valid.
	Nevertheless, the method for $\alpha\in(0,1)$ almost resolves the problem. More precisely, we adopt the framework of \Cref{sec:pfal01} and aim to establish Lemmas \ref{theorem_local_laws_along_flow}--\ref{theorem_que_along_flow}. For \Cref{theorem_local_laws_along_flow}, the problem is then reduced to proving \Cref{lemma_bootstrap} for $\alpha=0$. Combining this with the arguments in \Cref{sec_proof_of_que} yields Lemmas \ref{lemma_weak_local_law_T_along_flow}--\ref{theorem_que_along_flow}, which further implies all main results as in the argument below \Cref{theorem_que_along_flow}. To extend the proof of Lemma \ref{lemma_bootstrap} in \Cref{sec_step_1_to_3} to the case $\alpha=0$, we develop an alternative approach to Steps 2 and 3. Our method is inspired by \cite[Section 8]{YYYTexpansion2022CMP} and is based on the so-called $T$-expansion, corresponding to the expansion \eqref{eq:boundL-K2p} above.

	\medskip
	\noindent
	{\bf Step $\mathbf{2^\prime}$} (Sharp bounds for $\cL$-loops and sharp local law). For any $x,y\in\ZL$, we have the following optimal upper bound for $\cL$-loops:
	\begin{equation}\label{step_2_cL_upper_bound_critical}
		\max_{\bsigma\in\ha{-,+}^2}|\cL_{u,xy}^{\bsigma}|\prec B_u(|x-y|),\quad \forall u\in\qa{s,t}.
	\end{equation}
	As a consequence of \Cref{2_loop_to_1_chain}, the following sharp local laws hold for any $x,y\in \ZL$:
	\begin{equation}\label{step_2_sharp_local_law_critical}
		\begin{aligned}
			\absa{(G_u-m)_{xy}}^2\prec B_u\p{|x-y|},\quad \absa{\avg{(G_u-m)S^{\pa{x}}}}\prec B_u\p{0},\quad \forall u\in\qa{s,t}.
		\end{aligned}
	\end{equation}
	
	\noindent
	{\bf Step $\mathbf{3^\prime}$} ($T$-variables and $\cL$-loops estimates). For all $\bsigma\in\ha{-,+}^2$, $x,y,y'\in\ZL$, and $u\in\qa{s,t}$, we have that
	\begin{equation}\label{step_3_estimate_loops_critical}
		(T_u-\Theta_u)_{x,yy'}^{\bsigma}\prec \qa{B_u\p{0}}^{1/2}\q{B_u\p{|x-y|}B_u\p{|x-y'|}}^{1/2},\quad (\cL_u-\cK_u)_{xy}^{\bsigma}\prec \qa{B_u\p{0}}^{3/4}B_u\p{|x-y|} .
	\end{equation}
	
	\medskip

	First, since the arguments for \eqref{step_1_priori_bound} and \eqref{step_1_weak_local_law} remain valid for $\alpha=0$, these two estimates continue to hold. Moreover, assuming the upper bound \eqref{step_2_cL_upper_bound_critical}, we immediately obtain the local laws in \eqref{step_2_sharp_local_law_critical} from \Cref{2_loop_to_1_chain}. The estimates in \eqref{step_3_estimate_loops_critical} then follow by an argument similar to that used to prove \eqref{flow_loop_local_law_subcritical}. We therefore omit the details.
	It remains to establish the upper bound \eqref{step_2_cL_upper_bound_critical}. By \eqref{step_1_priori_bound}, there exists a deterministic control parameter $0<\Psi_u<W^{-c}$, for some sufficiently small constant $c>0$, such that
	\begin{equation}\label{step_2_input_upper_bound_critical}        \max_{\bsigma\in\ha{-,+}^2}\absa{\cL_{u,xy}^{\bsigma}}\prec B_u(|x-y|)+\Psi_u^2,\qquad \forall x,y\in \ZL,\ u\in\qa{s,t}.
	\end{equation}
	We next apply the following lemma, which provides a self-improving bound for the $\cL$-loops.

	\begin{lemma}\label{lemma_step_2_self_improved_upper_bound_critical}

		In the setting of Lemma \ref{lemma_bootstrap} with $\alpha=0$, suppose that the upper bound \eqref{step_2_input_upper_bound_critical} holds for a deterministic control parameter $\Psi_u$ satisfying $0<\Psi_u<W^{-c}$ for some constant $c>0$. Then the following improved upper bound holds uniformly for $u\in[s,t]$:
		\begin{equation}\label{eq:self_improve_al0}
			\max_{\bsigma\in\ha{-,+}^2}\absa{\cL_{u,xy}^{\bsigma}}\prec B_u(|x-y|)+B_u^{1/10}(0)\cdot \Psi_u^2,\quad \forall x,y\in\ZL.
		\end{equation}
	\end{lemma}
	
	Using the a priori bound \eqref{step_2_input_upper_bound_critical} as input and iterating \Cref{lemma_step_2_self_improved_upper_bound_critical} $\OO(1)$ times yields 
	\begin{equation*}        \max_{\bsigma\in\ha{-,+}^2}\absa{\cL_{u,xy}^{\bsigma}}\prec B_u(|x-y|)+W^{-D}\lesssim B_u(|x-y|),\qquad \forall x,y\in \ZL,\ u\in\qa{s,t},
	\end{equation*}
	for a large enough constant $D>0$, which concludes the upper bound \eqref{step_2_cL_upper_bound_critical}. 
	It remains to prove \Cref{lemma_step_2_self_improved_upper_bound_critical}. 
	Using the notation of \eqref{bound_Xi_subcritial_cL}, for any $\bsigma\in\ha{-,+}^2$ and $x,y\in\ZL$, we claim that
	\begin{equation}\label{step_2_Xi_bound_critical}
		\Xi_{u,xy}^{\cL,\bsigma}:=\qa{(1-um\p{\sigma_1}m\p{\sigma_2}S)\cdot (\cL-\cK)_{u}^{\bsigma}}_{xy}\prec  B_u^{1/10}(0)\cdot B_u(|x-y|)+B_u^{1/10}(0)\cdot \Psi_u^2.
	\end{equation}
	If $\sigma_1=\sigma_2$, we apply $(1-um(\sigma_1)m(\sigma_2)S)^{-1}$ to both sides of \eqref{step_2_Xi_bound_critical}. Using \eqref{Theta_bound_equaled_charge}, we obtain
	\begin{equation}\label{eq:self_improve_al0_Xi}
		\max_{\bsigma=(\sigma_1,\sigma_2): \sigma_1=\sigma_2}\absa{(\cL_u-\cK_u)_{xy}^{\bsigma}}\prec B_u^{1/10}(0)\cdot B_u(|x-y|)+B_u^{1/10}(0)\cdot \Psi_u^2.
	\end{equation}
	If $\sigma_1\ne\sigma_2$, we use the operators $\cP$ and $\cQ$ introduced in \Cref{def_sum_zero_operators} and write
	\begin{align}
		\qb{\cQ\cdot(\cL-\cK)_u^{\bsigma}}_{xy}&=\pb{\cQ\cdot \Xi_{u}^{\cL,\bsigma}}_{xy}+u\pb{\cQ\cdot\Theta_{u}^{\bsigma}\cdot \Xi_{u}^{\cL,\bsigma}}_{xy}\prec B_u^{1/10}(0)\cdot B_u(|x-y|)+B_u^{1/10}(0)\cdot \Psi_u^2,\label{eq:QT-Theta0}     \end{align}
	where in the second step we used \eqref{step_2_Xi_bound_critical} together with the bound \eqref{zero_mode_removed_bound} for $\cQ\cdot\Theta_{u}^{\bsigma}$.
	On the other hand, using Ward's identity and the local law $\avg{(G_u-m)S^{\pa{x}}}\prec B_u\p{0}+\Psi_u^2$ by \eqref{step_2_input_upper_bound_critical} and Lemma \ref{2_loop_to_1_chain}, we obtain 
	\begin{equation}\label{eq:PT-Theta0}
		\qb{\cP\cdot(\cL-\cK)_{u}^{\bsigma}}_{xy}\prec \qa{{N(1-u)}}^{-1}\cdot \pa{B_{u}(0)+\Psi_u^2}.
	\end{equation}
	Combining \eqref{eq:QT-Theta0} and \eqref{eq:PT-Theta0}, and using the identity $\cP+\cQ=1$, yields \eqref{eq:self_improve_al0} for each fixed $u\in[s,t]$. A standard $N^{-C}$-net argument then extends the bound uniformly to all $u\in\qa{s,t}$.

	Finally, we prove \eqref{step_2_Xi_bound_critical} by estimating the moments of $\Xi_{u,xy}^{\cL,\bsigma}$ using the expansion formula \eqref{eq:boundL-K2p} (with $t$ replaced by $u$). Using the bound \eqref{step_1_weak_local_law} together with Lemma \ref{2_loop_to_1_chain} under the assumption \eqref{step_2_input_upper_bound_critical}, and arguing as in the proof of \eqref{eq_decomposition_double_difference} in \Cref{sec_estimates_error_terms}, but more simply, we bound the right-hand side of \eqref{eq:boundL-K2p} by
	\begin{equation*}
		B_u^{1/5}(0)\qa{B_u(|x-y|)+\Psi_u^2}\cdot(\E|\Xi_{u,xy}^{\cL,\bsigma}|^{2p})^{(2p-1)/2p}+B_u^{1/5}(0)\qa{B_u(|x-y|)+\Psi_u^2}^2\cdot(\E|\Xi_{u,xy}^{\cL,\bsigma}|^{2p})^{(2p-2)/2p}.
	\end{equation*}
	Since the derivation is straightforward, we omit the details. Applying Young's inequality then yields
	\begin{equation*}
		\E|\Xi_{u,xy}^{\cL,\bsigma}|^{2p}\prec B_u^{p/5}(0)\qa{B_u(|x-y|)+\Psi_u^2}^{2p}.
	\end{equation*}
	Finally, the bound \eqref{step_2_Xi_bound_critical} follows from Markov's inequality.

	\section{Proof of the main results for \texorpdfstring{$\alpha\in [1,\infty)$}{alpha in [1,infty)}}\label{sec_proof_of_remaining_cases_1}

	In this section, we present the proofs of the main results, \Cref{theorem_local_law_entrywise,theorem_local_law_loops} and \Cref{theorem_local_law_T}, in the regime $\alpha\in[1,\infty)$. For consistency of notation, we adopt the flow framework introduced in \Cref{sec:notations}. Within this framework, the main results follow directly from the following estimates on the $\cL$-loops and $T$-variables established along the flow.
	\begin{theorem}\label{local_laws_along_flow_supercritical}
		Under the assumptions of \Cref{theorem_delocalization}, fix $\alpha\ge 1$, small constants $\kappa,\fc\in(0,1)$, and a spectral parameter $z=E+\ii \eta\in\mathbf{D}_{\kappa,\fc}$. Consider the flow framework in \Cref{def_H_t_flow,def_z_t_flow}, with the flow parameter chosen as in \eqref{eq:t0E0}. Then the following estimates hold uniformly for $t\in[0,\tf]$.
		\begin{itemize}
			\item {\bf Local law:} For any $x,y\in\ZL$, the following entrywise and averaged local laws hold:        \begin{equation}\label{flow_local_law_supercritical}
				\begin{aligned}
					\absa{(G_t-m)_{xy}}^2\prec B_t\p{|x-y|},\quad \big|\avg{(G_t-m)S^{\pa{x}}}\big|\prec B_t\p{0}.
				\end{aligned}
			\end{equation}
			\item {\bf $\cal L$-loop and $T$-variable estimates:} For any $\bsigma\in\ha{-,+}^2$ and $x,y,y'\in\ZL$, we have        \begin{equation}\label{flow_loop_local_law_supercritical}
				\begin{aligned}
					(T_t-\Theta_t)_{x,yy'}^{\bsigma}\prec B_t^{1/2}(0)\qa{B_t(|x-y|)B_t(|x-y'|)}^{1/2},\quad (\cL_t-\cK_t)_{xy}^{\bsigma}\prec B_t\p{0}B_t\p{|x-y|}.
				\end{aligned}
			\end{equation}
			\item {\bf Expected $\cal L$-loop estimate:} For any $\bsigma\in\ha{-,+}^2$ and $x,y\in\ZL$, we have        \begin{equation}\label{flow_expected_local_law_supercritical}
				\begin{aligned}
					\E(\cL_t-\cK_t)_{xy}^{\bsigma}\prec B_t^{3/2}\p{0}B_t\p{|x-y|}.
				\end{aligned}
			\end{equation}
		\end{itemize}
	\end{theorem}
	
	\begin{proof}[\bf Proof of \Cref{theorem_local_law_entrywise,theorem_local_law_loops} and \Cref{theorem_local_law_T} for $\al\in [1,\infty)$]
		By \eqref{eq:zztE}, we have $G(z)\stackrel{d}{=}\sqrt{\tf}G_{\tf,\sE}$. Hence, by \eqref{eq:Betar}, for each fixed $z\in\mathbf{D}_{\kappa,\fc}$ the estimates obtained along the flow directly yield the desired estimates for the original resolvent. 
		More precisely, \eqref{flow_local_law_supercritical} implies the entrywise and averaged local laws \eqref{eq_entrywise_local_law} and \eqref{eq_average_local_law}, \eqref{flow_loop_local_law_supercritical} yields the estimates \eqref{eq_loop_local_law_1} and \eqref{eq_chain_local_law_1}, and \eqref{flow_expected_local_law_supercritical} gives the expected bound \eqref{eq_expected_local_law_1}. 
		Finally, a standard $N^{-C}$-net argument extends these bounds uniformly to all $z\in\mathbf{D}_{\kappa,\fc}$.
	\end{proof}

	\subsection{Loop hierarchy}\label{sec_loop_hierarchy_supercritical}

	The remainder of this section is devoted to the proof of \Cref{local_laws_along_flow_supercritical}. As in the previous regimes, our approach is based on a dynamical analysis of the $\cL$-loops and $T$-variables. However, in contrast to the case $\alpha\in(0,1)$, and following the ideas developed for regular RBMs \cite{Band1D,Band2D,erdos2025zigzagstrategyrandomband}, our analysis in the regime $\alpha\in[1,\infty)$ relies on the dynamics of resolvent loops and chains of arbitrary length, defined below.\footnote{\label{footnoteloop} One could alternatively attempt to extend the proof for the regime $\alpha\in(0,1)$ to $\alpha\in[1,\infty)$. While most of the arguments carry over with minor modifications, a technical obstacle arises in Step 3 of the proof of \Cref{lemma_bootstrap}, namely in establishing \eqref{step_3_estimate_loops}. In the regime $\alpha\in[1,\infty)$, the improved evolution kernel bound \eqref{bound_cU_s_t_2} no longer holds under the sum-zero condition when $\sigma_1\ne\sigma_2$, and an additional power of $(1-s)/(1-t)$ appears; see \Cref{rmk_failureofLoop}. Cancelling this factor would require exploiting the CLT-type cancellation mechanism developed in \cite[Section 5]{dubova2025delocalizationnonmeanfieldrandommatrices}. 
    However, in the absence of the exponential decay of resolvent entries available in \cite{dubova2025delocalizationnonmeanfieldrandommatrices}, implementing this mechanism in our model becomes substantially more involved, if possible at all. 
}

	\begin{definition}[Loops and chains]\label{def_n_loops_and_chains}
		In the setting of \Cref{local_laws_along_flow_supercritical}, let $n\in\N$. For any charge $\bsigma=\pa{\sigma_1,\ldots,\sigma_n}\in\ha{-,+}^n$ and $\bx=(x_1,\ldots,x_n)\in\ZL^n$, we define the $\cL$-loop of length $n$ (called an $n$-$\cL$-loop) by
		\begin{equation*}            \cL_{t,\bsigma,\bx}^{\pa{n}}:=\avgb{G_t(\sigma_1)S^{\pa{x_1}}\cdots G_t(\sigma_n)S^{\pa{x_n}}}, \quad \forall t\in\qa{0,\tf}.
		\end{equation*}
		Moreover, for any charge $\bsigma\in\ha{-,+}^n$ and $(\bx,x,y)=(x_1,\ldots,x_{n-1},x,y)\in\ZL^{n+1}$ with $n\in\N$, we define the $\cC$-chain of length $n$ (called an $n$-$\cC$-chain) by
		\begin{equation*}            \cC_{t,\bsigma,\bx}^{\pa{n}}\pa{x,y}:=\pb{G_t(\sigma_1)S^{\pa{x_1}}G_t(\sigma_2)\cdots S^{(x_{n-1})} G_t(\sigma_n)}_{xy},\quad \forall t\in\qa{0,\tf}.
		\end{equation*}
		Note that the $1$-$\cC$-chain \smash{$\cC_{t,(\sigma),\bx}^{\pa{1}}\pa{x,y}$} coincides with the resolvent entry $G_{t,xy}(\sigma)$. For notational convenience, we abbreviate the diagonal $\cC$-chain $\cC_{t,\bsigma,\bx}^{(n)}(x,x)$ as
		\begin{equation}\label{abbreviation_diagonal_cC_chain}
			\begin{aligned}
				\cC_{t,\bsigma,(\bx,x)}^{\pa{n}}:=\cC_{t,\bsigma,\bx}^{\pa{n}}(x)\equiv\cC_{t,\bsigma,\bx}^{\pa{n}}(x,x).
			\end{aligned}
		\end{equation}
		Clearly, the loops and chains defined above generalize the $\cL$-loops and $T$-variables introduced in \Cref{def_loops_and_T_variables}. More precisely, for any charge $\bsigma\in\ha{-,+}^2$ and $x,y,y'\in\ZL$,
		\begin{equation*}            \cL_{t,\bsigma,\bx}^{\pa{2}}=\cL_{t,\bx}^{\bsigma},\quad \cC_{t,\bsigma,(x)}^{\pa{2}}(y,y')=T_{t,x,yy'}^{\bsigma},\quad \forall t\in\qa{0,\tf}.
		\end{equation*}
	\end{definition}
	
	To describe the dynamics of $\cL$-loops and $\cC$-chains, we introduce the following ``cut-and-glue'' operations acting on the indices.
	
	\begin{definition}\label{def_loop_operations}
		Fix $n\in\N$. We define the following operators.\footnote{Our definitions of \smash{$\cG_{k}^{(x)}$ and $(\cG_L)_{k,l}^{(x)}$} differ slightly from those in \cite{Band1D} in order to simplify the presentation of the chain factors in the loop hierarchy in \Cref{lemma_n_loop_evolution}.}
		
		\medskip
		\noindent
		\txt{(1)}
		For $k\in\qq{n}$ and $x\in\ZL$, we define the ``cut-and-glue'' operator $\cG_{k}^{(x)}$ by
		\begin{equation}
			\begin{aligned}
				\cG_{k}^{\pa{x}}\circ\pa{\bsigma,\bx}:= \qa{\pa{\sigma_k,\ldots,\sigma_n,\sigma_1,\ldots,\sigma_k},\pa{x_k,\ldots,x_n,x_1,\ldots,x_{k-1},x}}
			\end{aligned}
		\end{equation}
		for any $\bsigma=\pa{\sigma_1,\ldots,\sigma_n}\in\ha{-,+}^n$ and $\bx=(x_1,\ldots,x_n)\in\ZL^n$. Graphically, if $(\bsigma,\bx)$ is viewed as a discrete circle $(x_1,\ldots,x_n)$ in which $\sigma_i$ labels the $i$-th edge $(x_{i-1},x_i)$ for $i\in\qq{n}$, then the operator \smash{$\cG_k^{(x)}$} cuts the $k$-th edge and glues the two resulting ends to a new vertex $x$.

		\medskip
		\noindent
		\txt{(2)} For $k<l\in\qq{n}$ and $x\in\ZL$, we define two additional ``cut-and-glue'' operators: \smash{$(\cG_L)_{k,l}^{(x)}$}, which cuts from the left (“L”) of $k$, and \smash{$(\cG_R)_{k,l}^{(x)}$}, which cuts from the right of $k$, by
		\begin{equation}
			\begin{aligned}
				&\pa{\cG_L}_{k,l}^{\pa{x}}\circ\pa{\bsigma,\bx}:= \qa{\pa{\sigma_l,\ldots,\sigma_n,\sigma_1,\ldots,\sigma_k},\pa{x_l,\ldots,x_n,x_1,\ldots,x_{k-1},x}},\\
				&\pa{\cG_R}_{k,l}^{\pa{x}}\circ\pa{\bsigma,\bx}:= \qa{\pa{\sigma_k,\ldots,\sigma_l},\pa{x_k,\ldots,x_{l-1},x}}.
			\end{aligned}
		\end{equation}
		Graphically, starting from the circle $(\bsigma,\bx)$ described above, the configurations \smash{$(\cG_L)_{k,l}^{(x)}\circ(\bsigma,\bx)$} and \smash{$(\cG_R)_{k,l}^{(x)}\circ(\bsigma,\bx)$} correspond to the two circles obtained by cutting the $k$-th and $l$-th edges and gluing the four resulting ends to the vertex $x$.
	\end{definition}

	With the above notation, given an object $\cA_{\bsigma,\bx}^{(n)}$ with $\cA\in\ha{\cL,\cC}$ (and similarly for other loop or chain quantities introduced later, such as the $\cK$-loops in \Cref{def_n_cK_loop} below), the cut-and-glue operators act on \smash{$\cA_{\bsigma,\bx}^{(n)}$} through the indices as
	\begin{equation}\label{eq;action_cut_glue}        \cG_{k}^{\pa{x}}\circ\cA_{\bsigma,\bx}^{\pa{n}}:=\cA_{\cG_{k}^{\pa{x}}\circ\pa{\bsigma,\bx}}^{\pa{n+1}},\quad \pa{\cG_L}_{k,l}^{\pa{x}}\circ\cA_{\bsigma,\bx}^{\pa{n}}=\cA_{\pa{\cG_L}_{k,l}^{\pa{x}}\circ\pa{\bsigma,\bx}}^{\pa{n-l+k+1}},\quad \pa{\cG_R}_{k,l}^{\pa{x}}\circ\cA_{\bsigma,\bx}^{\pa{n}}=\cA_{\pa{\cG_R}_{k,l}^{\pa{x}}\circ\pa{\bsigma,\bx}}^{\pa{l-k+1}}.
	\end{equation}
	As a generalization of \Cref{lemma_evolution_equations}, we derive the following evolution equations for $\cL$-loops using Itô's calculus. These equations are referred to as the \emph{loop hierarchy} in \cite{Band1D,Band2D,dubova2025delocalizationnonmeanfieldrandommatrices}.
	
	\begin{lemma}[Loop hierarchy]\label{lemma_n_loop_evolution}
		Adopting the notation of \Cref{def_loop_operations}, let $\bsigma\in\ha{-,+}^n$ and $\bx=(x_1,\ldots,x_n)\in\ZL^n$ with $n\in\N$. The $\cL$-loops defined in \Cref{def_n_loops_and_chains} satisfy the following SDE:
		\begin{equation}\label{eq_evolution_n_loops}
			\begin{aligned}
				\rd \cL_{t,\bsigma,\bx}^{\pa{n}}&=\sum_{1\leq k<l\leq n}\sum_{x\in\ZL}\pa{\pa{\cG_L}_{k,l}^{\pa{x}}\circ\cL_{t,\bsigma,\bx}^{\pa{n}}}\cdot \pa{ \pa{\cG_R}_{k,l}^{\pa{x}}\circ\cC_{t,\bsigma,\bx}^{\pa{n}}} \rd t+\cW_{t,\bsigma,\bx}^{\pa{n}}\, \rd t+ \rd \cB_{t,\bsigma,\bx}^{\pa{n}}\\
				&=\sum_{1\leq k<l\leq n}\sum_{x\in\ZL}\pa{\pa{\cG_L}_{k,l}^{\pa{x}}\circ\cC_{t,\bsigma,\bx}^{\pa{n}}}\cdot \pa{\pa{\cG_R}_{k,l}^{\pa{x}}\circ\cL_{t,\bsigma,\bx}^{\pa{n}}} \rd t+\cW_{t,\bsigma,\bx}^{\pa{n}}\, \rd t+ \rd \cB_{t,\bsigma,\bx}^{\pa{n}},
			\end{aligned}
		\end{equation}
		where $\cC_{t,\bsigma,\bx}^{(n)}$ is defined in the sense of \eqref{abbreviation_diagonal_cC_chain} by \smash{$\cC_{t,\bsigma,\bx}^{\pa{n}}=\cC_{t,\bsigma,\bx'}^{\pa{n}}(x_n)$} with $\bx'=(x_1,\ldots, x_{n-1})$, and the cut-and-glue operators act according to \eqref{eq;action_cut_glue}. The light-weight and martingale terms are given by
		\begin{align}            \cW_{t,\bsigma,\bx}^{\pa{n}}&:=\sum_{k=1}^n\sum_{x\in\ZL}\avgb{\qa{G_t(\sigma_k)-m\pa{\sigma_k}}S^{\pa{x}}}\cdot \pa{\cG_{k}^{\pa{x}}\circ\cC_{t,\bsigma,\bx}^{\pa{n}}},\label{eq:LW_long}\\
			\rd \cB_{t,\bsigma,\bx}^{\pa{n}}&:=\sum_{x,y\in\ZL}\pa{\partial_{xy}\cL_{t,\bsigma,\bx}^{\pa{n}}}\cdot \sqrt{S_{xy}}\,\rd \mathbf{B}_{xy}(t).\label{eq:MG_long}
		\end{align}
	\end{lemma}
	
	The $\cL$-loops and $\cC$-chains can be well approximated by deterministic objects called $\cK$-loops and $\cK$-chains, defined below. Following \cite{dubova2025delocalizationnonmeanfieldrandommatrices,Band1D}, we refer to these as \emph{primitive loops}, which admit explicit tree representations as shown in \cite{Band1D,truong2025localizationlengthfinitevolumerandom}.   

	\begin{definition}[Primitive loops]\label{def_n_cK_loop}
	We first introduce the $\wh\cK$-loops of length $1$ by
		\begin{equation*}
			\begin{aligned}
				\wh \cK_{t,\sigma,x}^{\pa{1}}:=m\p{\sigma},\quad \forall t\in\qa{0,\tf},\  \sigma\in\ha{-,+},\  x\in\ZL.
			\end{aligned}
		\end{equation*}
		For $n\geq 2$, we define the $\wh\cK$-loops inductively as the unique solution to the system of equations
		\begin{equation}\label{evolution_wh_cK_loops}
			\rd \wh \cK_{t,\bsigma,\bx}^{\pa{n}}=\sum_{x,y\in\ZL}\sum_{1\leq k<l\leq n}\pa{\pa{\cG_L}_{k,l}^{\pa{x}}\circ \wh \cK_{t,\bsigma,\bx}^{\pa{n}}}  S_{xy} \pa{ \pa{\cG_R}_{k,l}^{\pa{y}}\circ \wh \cK_{t,\bsigma,\bx}^{\pa{n}}} \rd t,
		\end{equation}
		with initial condition 
		\smash{\(\wh \cK_{0,\bsigma,\bx}^{\pa{n}}=m\p{\sigma_1}\cdots m\p{\sigma_n}\delta_{x_1\cdots x_n}.\)} 
		To approximate $\cL$-loops and $\cC$-chains, we define $\cK^{\cL}$-loops and $\cK^{\cC}$-chains as follows. For $\bsigma=(\sigma_1,\ldots,\sigma_n)\in\ha{-,+}^n$ and $\bx=(x_1,\ldots,x_n)\in\ZL^n$, set
		\begin{equation*}
			\cK_{t,\bsigma,\bx}^{\cL,\pa{n}}:=\sum_{\ba=(a_1,\ldots,a_n)\in\ZL^n}S_{x_1a_1}\cdots S_{x_n a_n} \cdot \wh \cK_{t,\bsigma,\ba}^{\pa{n}}.
		\end{equation*}
		For $\bsigma=(\sigma_1,\ldots,\sigma_n)\in\ha{-,+}^n$ and $(\bx,x,y)=(x_1,\ldots,x_{n-1},x,y)\in\ZL^{n+1}$, define
		\begin{equation}\label{whK_sigmax_cC}
			\cK_{t,\bsigma,\bx}^{\cC,\pa{n}}(x,y):=\delta_{xy}\sum_{\ba=(a_1,\ldots,a_{n-1})\in\ZL^{n-1}}S_{x_1a_1}\cdots S_{x_{n-1}a_{n-1}} \cdot \wh \cK_{t,\bsigma,(\ba,x)}^{\pa{n}}.
		\end{equation}
		For simplicity of notation, we write  $\cK_{t,\bsigma,\bx}^{\pa{n}}\equiv\cK_{t,\bsigma,\bx}^{\cL,\pa{n}}$. We also introduce the following abbreviation for diagonal chains, analogous to \eqref{abbreviation_diagonal_cC_chain}:
		\begin{equation}\label{abbreviation_diagonal_cK_chain}
			\begin{aligned}
				\cK_{t,\bsigma,(\bx,x)}^{\cC,\pa
					n}:=\cK_{t,\bsigma,\bx}^{\cC,\pa
					n}(x)\equiv\cK_{t,\bsigma,\bx}^{\cC,\pa{n}}(x,x).
			\end{aligned}
		\end{equation}
		Furthermore, from \eqref{evolution_wh_cK_loops} it follows readily that for $n\ge2$,
		\begin{equation}\label{eq_evolution_n_cK_loops}
			\begin{aligned}
				\rd \cK_{t,\bsigma,\bx}^{\pa{n}}=&\sum_{1\leq k<l\leq n}\sum_{x\in\ZL}\pa{\pa{\cG_L}_{k,l}^{\pa{x}}\circ \cK_{t,\bsigma,\bx}^{\cL,\pa{n}}} \cdot \pa{\pa{\cG_R}_{k,l}^{\pa{x}}\circ \cK_{t,\bsigma,\bx}^{\cC,\pa{n}}} \rd t\\
				=&\sum_{1\leq k<l\leq n}\sum_{x\in\ZL}\pa{\pa{\cG_L}_{k,l}^{\pa{x}}\circ \cK_{t,\bsigma,\bx}^{\cC,\pa{n}}} \cdot \pa{\pa{\cG_R}_{k,l}^{\pa{x}}\circ \cK_{t,\bsigma,\bx}^{\cL,\pa{n}}} \rd t,
			\end{aligned}
		\end{equation}
		which coincides exactly with the leading term in \eqref{eq_evolution_n_loops}.
	\end{definition}
	
	For $n=2$, it is easy to see that $ \wh \cK_{t,\bsigma,\bx}^{\pa{2}}=m(\sigma_1)m(\sigma_2)\pa{1-tm(\sigma_1)m(\sigma_2)S}^{-1}$. Consequently,  
	\begin{equation}\label{whK_sigmax}
		\cK_{t,\bsigma,(x,y)}^{\cL,\pa{2}}=\cK_{t,xy}^{\bsigma},\quad  \cK_{t,\bsigma,\bx}^{\cC,\pa{2}}=\Theta_{t,x,yy}^{\bsigma},
	\end{equation}
	where $\cK_{t}^{\bsigma}$ and $\Theta_{t}^{\bsigma}$ are defined in \eqref{eq:KTheta}. 
	Consider the loop $\cL_{t,\bsigma,\bx}^{\pa{n}}$ with $\bsigma=(\sigma_1,\ldots,\sigma_n)\in\ha{-,+}^n$ and $\bx=(x_1,\ldots,x_n)\in\ZL^n$. Suppose that $\sigma_1\neq \sigma_{n}$. A direct application of Ward's identity \eqref{eq_Ward0} yields \begin{equation}\label{WI_calL}
		\sum_{x_n\in\ZL}\cL_{t,\bsigma,\bx}^{\pa{n}}=\frac{1}{2\ii \eta_t}\left( {\cL}^{(n-1)}_{t,  \wh\bsigma^{(+,n)}, \wh\bx^{(n)}}- {\cL}^{(n-1)}_{t,  \wh\bsigma^{(-,n)}, \wh\bx^{(n)}}\right),
	\end{equation}
	where $\wh\bsigma^{(\pm,n)}$ is defined by removing $\sigma_n$ from $\boldsymbol{\sigma}$ and replacing $\sigma_1$ with $\pm$, i.e., $\wh\bsigma^{(\pm,n)}=(\pm, \sigma_2, \cdots, \sigma_{n-1})$, and  $\wh\bx^{(n)}$ is obtained by removing $x_n$ from $\bx$, i.e.,
	$\wh\bx^{(n)}=(x_1, x_2,\cdots, x_{n-1})$. For $\cK$-loops, Corollary 4.21 in \cite{truong2025localizationlengthfinitevolumerandom} shows that an analogous identity holds, leading to the following result.

	\begin{lemma}[Ward's identity for $\cK$-loops]\label{lemma_loops_Wards_identities}
		Let $n\ge 2$ and let $\bsigma=(\sigma_1,\ldots,\sigma_n)$ satisfy $\sigma_1\neq \sigma_n$. Then
		\begin{equation}\label{WI_calK}
			\sum_{x_n\in\ZL}\cK_{t,\bsigma,\bx}^{\pa{n}}=\frac{1}{2\ii \eta_t}\left( {\cK}^{(n-1)}_{t,  \wh\bsigma^{(+,n)}, \wh\bx^{(n)}}- {\cK}^{(n-1)}_{t,  \wh\bsigma^{(-,n)}, \wh\bx^{(n)}}\right).
		\end{equation}
	\end{lemma}
	
	Subtracting the evolution equation \eqref{eq_evolution_n_cK_loops} from \eqref{eq_evolution_n_loops}, we obtain the following evolution equation for $(\cL-\cK)$-loops, which generalizes the length-two equation \eqref{cL_cK_evolution}.
	
	\begin{lemma}
		For any $\bsigma\in\ha{-,+}^n$ and $\bx\in\ZL^n$, the $(\cL-\cK)$-loops satisfy     \begin{equation}\label{eq_evolution_n_cL_cK_loops}
			\rd(\mathcal{L} - \mathcal{K})^{(n)}_{t, \boldsymbol{\sigma}, \bx}
			=\sum_{l_{\cK}=2}^{n}\qa{\mathfrak{D}_{l_{\cK}}(\cL-\cK)}_{t,\bsigma,\bx}^{\pa{n}}\, \rd t+ \mathcal{E}^{(n)}_{t, \boldsymbol{\sigma}, \bx}\dd t 
			+\mathcal{W}^{(n)}_{t, \boldsymbol{\sigma}, \bx}
			\dd t +
			\dd\mathcal{B}^{(n)}_{t, \boldsymbol{\sigma}, \bx}.
		\end{equation}
		The linear term is defined by
		\begin{align}
			\qa{\mathfrak{D}_{l_{\cK}}(\cL-\cK)}_{t,\bsigma,\bx}^{\pa{n}}:=&\sum_{1 \leq k < l \leq n:l-k=l_{\cK}-1} \sum_{x\in\ZL}
			(\mathcal{L} - \mathcal{K}^{\cL})^{(n+k-l+1)}_{t, \pa{\cG_L}^{(x)}_{k, l}\left(\boldsymbol{\sigma},\, \bx\right)}
			\mathcal{K}^{\cC,(l-k+1)}_{t,\pa{\cG_R}^{(x)}_{k, l}\left(\boldsymbol{\sigma} ,\bx\right)} \nonumber\\
			&+\sum_{1 \leq k < l \leq n: l-k=n-l_{\cK}+1} \sum_{x\in\ZL} \mathcal{K}^{\cC,(n+k-l+1)}_{t, \pa{\cG_L}^{(x)}_{k, l}\left(\boldsymbol{\sigma},\bx\right)}
			\left(\cL-\cK^{\cL}\right)^{(l-k+1)}_{t,\pa{\cG_R}^{(x)}_{k, l}\left(\boldsymbol{\sigma},\bx\right)} .\label{eq:linear_error_n}
		\end{align}
		The quadratic error term is given by
		\begin{align}
			\mathcal{E}^{(n)}_{t, \boldsymbol{\sigma}, \bx}:=&\sum_{1 \leq k < l \leq n} \sum_{x\in\ZL}
			(\mathcal{L} - \mathcal{K}^{\cL})^{(n+k-l+1)}_{t, \pa{\cG_L}^{(x)}_{k, l}\left(\boldsymbol{\sigma},\, \bx\right)}
			(\cC-\mathcal{K}^{\cC})^{(l-k+1)}_{t,\pa{\cG_R}^{(x)}_{k, l}\left(\boldsymbol{\sigma} ,\bx\right)}\nonumber\\
			=&\sum_{1 \leq k < l \leq n} \sum_{x\in\ZL}
			(\mathcal{C} - \mathcal{K}^{\cC})^{(n+k-l+1)}_{t, \pa{\cG_L}^{(x)}_{k, l}\left(\boldsymbol{\sigma},\, \bx\right)}
			(\cL-\mathcal{K}^{\cL})^{(l-k+1)}_{t,\pa{\cG_R}^{(x)}_{k, l}\left(\boldsymbol{\sigma} ,\bx\right)}.\label{eq:quad_error_n}
		\end{align}
	\end{lemma}
	
	The leading term on the RHS of \eqref{eq_evolution_n_cL_cK_loops} is $\qa{\mathfrak{D}_{2}(\cL-\cK)}_{t,\bsigma,\bx}^{(n)}$, which induces the following \emph{evolution kernel}.
	
	\begin{definition}[Evolution kernel]
		For any charge $\bsigma\in\ha{-,+}^n$ and $n$-tensor $\cal A^{(n)}$, the term $\qa{\mathfrak{D}_{2}(\cL-\cK)}_{t,\bsigma,\bx}^{\pa{n}}$ induces the following action of the $\Theta$-propagator on $n$-tensors:
		\begin{equation}\label{eq:bTheta_n}
			\left[\bTheta_{t,\boldsymbol{\sigma}}^{(n)}\circ \mathcal{A}^{(n)}\right]_{\mathbf{x}}=\sum_{i=1}^{n}\sum_{a_i\in\ZL}\Theta_{t,x_ia_i}^{(\sigma_i,\sigma_{i+1})}\cdot\mathcal{A}_{\mathbf{x}^{(i)}(a_i)}^{(n)},
		\end{equation}
		where we adopt the cyclic convention $\sigma_{n+1}=\sigma_1$, and $\bx^{(i)}$ is defined as \(\bx^{(i)}(a_i):= (x_1, \ldots, x_{i-1}, a_i, x_{i+1}, \ldots, x_n).\) 
		In particular, with \eqref{whK_sigmax}, we notice that     \[\qa{\mathfrak{D}_{2}(\cL-\cK)}_{t,\bsigma,\bx}^{\pa{n}}=\qa{\bTheta^{(n)}_{t,\bsigma}\circ (\cL-\cK)_{t,\bsigma}^{\pa{n}}}_{\bx}.\]
		For $0\leq s< t\leq \tf$, the evolution kernel from $s$ to $t$ associated with the operator $\bTheta_{t,\bsigma}^{(n)}$ is defined by
		\begin{equation}\label{eq:bUst}
			\qa{\bU_{s,t,\bsigma}^{\pa{n}}\circ \cA^{\pa{n}}}_{\bx}:=\sum_{\ba=(a_1,\ldots,a_n)\in\ZL^n}\prod_{i=1}^{n}\cU_{s,t}^{(\sigma_i,\sigma_{i+1})}\pa{x_i,a_i}\cdot \cA^{\pa{n}}_{\ba},
		\end{equation}
		where $\cU_{s,t}^{\bsigma}$ is defined in \eqref{U_s_t_equal_1_Theta}. One readily checks that
		\begin{equation*}
			\begin{aligned}
				\frac{\rd }{\rd t}\,\bU_{s,t,\bsigma}^{\pa{n}}\circ \cA^{\pa{n}}=\bTheta_{t,\bsigma}^{\pa{n}}\circ\bU_{s,t,\bsigma}^{\pa{n}}\circ \cA^{\pa{n}}
			\end{aligned}.
		\end{equation*}
	\end{definition}

	With the evolution kernel defined in \eqref{eq:bUst}, applying Duhamel's principle, we can rewrite the evolution equation \eqref{eq_evolution_n_cL_cK_loops} as follows. For any $s\in\qa{0,\tf}$ and stopping time $s\leq \tau\leq t$, we have  
	\begin{equation}\label{evolution_n_cL_cK_integrated}
		\begin{aligned}
			(\mathcal{L} - \mathcal{K})^{(n)}_{\tau, \boldsymbol{\sigma}, \bx}=&\pa{\bU_{s,\tau,\bsigma}^{\pa{n}}\circ \pa{\cL-\cK}_{s,\bsigma}^{\pa{n}}}_{\bx}+\sum_{l_{\cK}=3}^{n}\int_{s}^{\tau}\pa{\bU_{u,\tau,\bsigma}^{\pa{n}}\circ\qa{\mathfrak{D}_{l_{\cK}}(\cL-\cK)}_{u,\bsigma}^{\pa{n}}}_{\bx}\, \rd u\\
			&+ \int_{s}^\tau\pa{\bU_{u,\tau,\bsigma}^{\pa{n}}\circ\mathcal{E}^{(n)}_{u, \boldsymbol{\sigma}}}_{\bx}\dd u+\int_{s}^{\tau}\pa{\bU_{u,\tau,\bsigma}^{\pa{n}}\circ\mathcal{W}^{(n)}_{u, \boldsymbol{\sigma}}}_{\bx}
			\dd u + \int_{s}^{\tau}\pa{\bU_{u,\tau,\bsigma}^{\pa{n}}\circ
				\dd\mathcal{B}^{(n)}_{u, \boldsymbol{\sigma}}}_{\bx}    .
		\end{aligned}
	\end{equation}
    Here, the illegal stochastic integral in the martingale term is understood in a similar sense to that in \Cref{remark_non_measurability}. 
	To capture the spatial decay of the $\cL$-loops and $\cC$-chains, we introduce the decay factors used in the upper bounds of loops and chains, together with the corresponding decay-removed norm.
	
	\begin{definition}[Decay factors and decay-removed norm]\label{def_decay_factors}
		Given $\bx=(x_1,\ldots,x_n)\in\ZL^n$, we define the decay factor for loops by
		\begin{equation*}
			\begin{aligned}
				\mathscr{D}_{t,\bx}^{\pa{n}}:=\prod_{i=1}^{n}\mathscr{D}_{t}(|x_{i}-x_{i+1}|),\quad \text{with}\quad \mathscr{D}_{t}(r):=\pa{r/\ell_t+1}^{-(1+\alpha)/2} \ \ \text{for}\ \ r\ge 0,
			\end{aligned}
		\end{equation*}
		where we adopt the cyclic convention $x_{n+1}=x_1$, and the decay factor $\mathscr{D}_{t}(r)$ originates from \eqref{def_B_t} for $\alpha\ge1$. For chains, given $(\bx,x,y)=(x_1,\ldots,x_{n-1},x,y)\in\ZL^{n+1}$, we define
		\begin{equation*}
			\begin{aligned}
				\mathscr{D}_{t,\bx}^{\pa{n}}(x,y):=\mathscr{D}_t(|x-x_1|)\mathscr{D}_t(|x_1-x_2|)\cdots \mathscr{D}_t(|x_{n-2}-x_{n-1}|)\mathscr{D}_t(|x_{n-1}-y|).
			\end{aligned}
		\end{equation*}
	 In particular, for $n=1$, the decay factor $\mathscr{D}_{t,\bx}^{\pa{n}}(x,y)$ reduces to $\mathscr{D}_{t}(|x-y|)$. We again adopt the abbreviation $\mathscr{D}_{t,\bx}^{\pa{n}}(x)\equiv\mathscr{D}_{t,\bx}^{\pa{n}}(x,x)$, under which
		\begin{equation*}
			\begin{aligned}
				\mathscr{D}_{t,(\bx,x)}^{\pa{n}}= \mathscr{D}_{t,\bx}^{\pa{n}}(x)\equiv\mathscr{D}_{t,\bx}^{\pa{n}}(x,x),
			\end{aligned}
		\end{equation*}
		which is consistent with the abbreviations \eqref{abbreviation_diagonal_cC_chain} and \eqref{abbreviation_diagonal_cK_chain}. 
		For any $n$-tensor $\cA_{\bx}^{\pa{n}}$ on $\ZL^n$ and $t\in\qa{0,\tf}$, we define the \emph{decay-removed norm} of $\cA^{\pa{n}}$ at time $t$ by     \begin{equation*}            \normb{\cA^{\pa{n}}}_{\mathscr{D}_{t}^{-1}}:=\max_{\bx\in\ZL^n}\absb{\cA_{\bx}^{\pa{n}}}/\mathscr{D}_{t,\bx}^{\pa{n}}.
		\end{equation*}
	\end{definition}

    In the analysis of the loop hierarchy, estimates for the primitive loops (i.e., the \(\cK\)-loops) serve as a key technical input.
	
	\begin{lemma}\label{lemma_n_cK_bound}
	For the $\cK$-loops defined in \Cref{def_n_cK_loop}, the following bound holds uniformly for any $t\in\qa{0,\tf}$, $\bsigma\in\ha{-,+}^n$, and $\bx\in\ZL^n$ with $n\in\N$:
		\begin{equation}   \label{eq:KL_KC}       \absa{\cK_{t,\bsigma,\bx}^{\cL,\pa{n}}}+\absa{\cK_{t,\bsigma,\bx}^{\cC,(n)}}\prec B_t^{n-1}(0)\cdot\mathscr{D}_{t,\bx}^{\pa{n}}.
		\end{equation}
	\end{lemma}

For $\alpha\geq 3$, this estimate can be established using either the tree representation for \(\cK\)-loops, originally introduced in \cite{Band1D} and further developed in \cite{truong2025localizationlengthfinitevolumerandom}, or a further development of the dynamical analysis introduced in \cite{erdos2025zigzagstrategyrandomband,fan2025blockreductionmethodrandom}. However, for $\alpha\in[1,3)$, the decay factor within the loop decay turns out to be too weak for the strategies therein. Therefore, we introduce and establish a sharper bound than \eqref{eq:KL_KC}, where the loop-decay factor \smash{$\mathscr{D}_{t,\bx}^{(n)}$} is replaced by the stronger \emph{tree-shaped decay factor} \smash{$\mathscr{T}_{t,\bx}^{(n)}$}, defined as follows.
	
	\begin{definition}[Tree-shaped decay factors]\label{def_tree_shaped_decay}
		Given ordered vertices $\bx=(x_1,\ldots,x_n)\in\ZL^n$, let $\txt{NCT}(\bx)$ denote the set of all \emph{non-crossing trees} on $\bx$, i.e., trees $T=(V,E)$ with $V=(x_1,\ldots,x_n)$ such that no two edges $(x_i,x_j),(x_k,x_l)\in E$ coexist whenever $1\leq i<k<j<l\leq n$.\footnote{Geometrically, if the vertices $(x_1, \ldots, x_n)$ are arranged cyclically along a circle, a tree is non-crossing if and only if no two of its edges cross inside the circle.} The \emph{tree-shaped decay factor} \smash{$\mathscr{T}_{t,\bx}^{(n)}$} is then defined for any $t\in\qa{0,\tf}$ by
		\begin{equation}\label{eq:mathcalDT}
			\mathscr{T}_{t,\bx}^{(n)}:=\sum_{T\in\txt{NCT}(\bx)}\mathcal{D}_t(T),\quad \text{with}\quad \mathcal{D}_t(T):=\prod_{(x_i,x_j)\in E(T)}\mathscr{D}_t^2(|x_i-x_j|).
		\end{equation}
	\end{definition}

	\begin{lemma}\label{lemma_n_cK_tree_bound}
		For the $\cK$-loops defined in Definition \ref{def_n_cK_loop}, the following bound holds uniformly for any $t\in\qa{0,\tf}$, $\bsigma\in\ha{-,+}^n$, and $\bx\in\ZL^n$ with $n\in\N$:
		\begin{equation}\label{eq_n_cK_tree_bound}
			\begin{aligned}            \absa{\cK_{t,\bsigma,\bx}^{\cL,\pa{n}}}+\absa{\cK_{t,\bsigma,\bx}^{\cC,\pa{n}}}\prec B_t^{n-1}(0)\cdot\mathscr{T}_{t,\bx}^{\pa{n}}.
			\end{aligned}
		\end{equation}
	\end{lemma}

	We believe that the bound \eqref{eq_n_cK_tree_bound} is sharp. Its proof is deferred to \Cref{sec_proof_of_lemma_n_cK_bound}. 
    It is easy to see that Lemma \ref{lemma_n_cK_bound} follows as an immediate corollary.
	\begin{proof}[\bf Proof of Lemma \ref{lemma_n_cK_bound}]
		It suffices to prove that for any $\bx=(x_1,\ldots,x_n)\in\ZL^n$,
		\begin{equation}\label{eq_tree_bound_by_loop}
			\mathscr{T}_{t,\bx}^{(n)}\lesssim \mathscr{D}_{t,\bx}^{(n)}.
		\end{equation}
		We prove this bound by induction on $n$. For $n=1$, we have $\mathscr{T}_{t,\bx}^{(1)}=\mathscr{D}_{t,\bx}^{(1)}=1$ under the cyclic convention $x_1=x_2$. For $n=2$, it is immediate that \smash{$\mathscr{T}_{t,\bx}^{(2)}=\mathscr{D}_{t,\bx}^{(2)}$}.
		Now suppose $n\geq 3$ and that \eqref{eq_tree_bound_by_loop} holds for all $1,2,\ldots,n-1$. Let $T\in\txt{NCT}(\bx)$ with $\bx=(x_1,\ldots,x_n)\in\ZL^n$. Relabeling the vertices if necessary, we may assume that $x_1$ is a leaf of $T$ and that $(x_1,x_a)\in E(T)$ for some $a\leq \floor{(n+1)/2}\leq n-1$. Let $T_1$ and $T_2$ denote the restrictions of $T$ to the vertex sets $(x_1,\ldots,x_a)$ and $(x_{a},\ldots,x_n)$ respectively. Then, $T_1$ and $T_2$ are both non-crossing trees with at most $(n-1)$ vertices. Applying the induction hypothesis, we obtain 
		\begin{equation*}
			\begin{aligned}
				\prod_{(x_i,x_j)\in E(T)}\mathscr{D}_t^2(|x_i-x_j|)&=\prod_{(x_i,x_j)\in E(T_1)}\mathscr{D}_t^2(|x_i-x_j|)\cdot\prod_{(x_i,x_j)\in E(T_2)}\mathscr{D}_t^2(|x_i-x_j|)\\
				&\lesssim \mathscr{D}_t(|x_1-x_2|)\cdots \mathscr{D}_t(|x_a-x_1|)\cdot \mathscr{D}_t(|x_a-x_{a+1}|)\cdots \mathscr{D}_t(|x_n-x_{a}|)\lesssim \mathscr{D}_{t,\bx}^{(n)},
			\end{aligned}
		\end{equation*}
		where in the last step we used $\mathscr{D}_t(|x_a-x_1|)\mathscr{D}_t(|x_n-x_{a}|)\lesssim \mathscr{D}_t(|x_n-x_1|)$. This completes the induction and proves \eqref{eq_tree_bound_by_loop}. 
	\end{proof}

	\begin{remark}
		The tree-shaped bound \eqref{eq_n_cK_tree_bound} for $\cK$-loops is in fact rather strong, in the sense that one cannot expect an analogous bound to hold for the $\cL$-loops. Indeed, according to the sharp local law \eqref{flow_local_law_supercritical}, the loop \smash{$\cL_{t,\bsigma,\bx}^{\pa{n}}$} should exhibit the decay \smash{$ \mathscr{D}_{t,\bx}^{(n)}$}, up to additional factors of $B_t(0)$, which is strictly weaker than the tree-shaped decay.
		We conjecture that \smash{$ \mathscr{T}_{t,\bx}^{(n)}$} captures the correct spatial decay of the expectation of \smash{$\cL_{t,\bsigma,\bx}^{\pa{n}}$}, whereas the loop decay \smash{$ \mathscr{D}_{t,\bx}^{(n)}$} primarily arises from fluctuations. Consequently, the essential issue in establishing a tree-shaped decay for the $\cL$-loops lies in controlling the martingale term in the evolution equation \eqref{eq_evolution_n_loops}, which is the principal source of the fluctuations of the $\cL$-loops.

		In a broader sense, the tree-shaped decay in \eqref{eq_n_cK_tree_bound} is also necessary for understanding the behavior of the $\cK$-loops when extending the analysis to the regime $\alpha\in(0,1)$. For example, when $1-t\leq (W/N)^{\alpha}$,
		\begin{equation}\label{eq:Dt_al01}
			\mathscr{D}_t(r):=(r/W+1)^{\frac12(\alpha-1)}(r/\ell_t+1)^{-\alpha},
		\end{equation}
		as suggested by the shape parameter in \eqref{def_B_t} with $\al\in(0,1)$.\footnote{We prove Lemma \ref{lemma_n_cK_tree_bound} only for $\alpha\in[1,\infty)$, but the argument can be adapted to treat $\alpha\in(0,1)$. We do not pursue this extension here, since our analysis does not rely on higher-order $\cL$-loops and hence does not require the corresponding higher-order $\cK$-loops as their deterministic limits in any sense.}
		In the regime $\al\ge 1$, the decay factor $\mathscr{D}_t$ is stable under convolution in the following sense:  
		\begin{equation}\label{example_triangular_2_edge_to_1}
			\sum_{x\in\ZL}\mathscr{D}_{t}(|a-x|)\mathscr{D}_{t}(|b-x|)\lesssim \mathscr{D}_{t}(|a-b|)\sum_{x\in\ZL}\mathscr{D}_{t}(|a-x|\wedge |b-x|)\prec \ell_t\mathscr{D}_{t}(|a-b|),\quad \forall \al \ge 1.
		\end{equation}
		{Here, we employ the $\prec$-notation to absorb the potential $\log N$-factors in the case $\alpha=1$.} This implies that the loop decay factor \smash{$ \mathscr{D}_{t,\bx}^{(n)}$} is stable under partial summation over one vertex:
		\begin{equation}\label{example_triangular_2_edge_to_1_loop}
			\sum_{x_n\in\ZL}B_{t}^{n-1}(0)\mathscr{D}_{t,(x_1,\ldots,x_{n})}^{(n)}\prec \ell_t\cdot B_{t}^{n-1}(0) \mathscr{D}_{t,(x_1,\ldots,x_{n-1})}^{(n-1)} = (1-t)^{-1}\cdot B_{t}^{n-2}(0) \mathscr{D}_{t,(x_1,\ldots,x_{n-1})}^{(n-1)} .
		\end{equation}
		This property is consistent with Ward's identity for the $\cK$-loops in \eqref{WI_calK}. However, the convolution bound \eqref{example_triangular_2_edge_to_1} fails for the decay factor \eqref{eq:Dt_al01} in the regime $\al\in(0,1)$. Consequently, the relation \eqref{example_triangular_2_edge_to_1_loop} no longer holds, which is inconsistent with Ward's identity for the $\cK$-loops. In contrast, the tree-shaped decay avoids this difficulty and yields a similar vertex-reduction relation as in \eqref{example_triangular_2_edge_to_1_loop}, using the bound
		\[\sum_{x\in \ZL}B_t(x) \prec (1-t)^{-1},\quad \forall \al\in(0,1).\]
		This is consistent with Ward's identity for the $\cK$-loops in the regime $\al\in(0,1)$. It suggests that, in order to control long $\cL$-loops for $\al\in(0,1)$, one would need to exploit the sharp $\cK$-loop decay established in Lemma \ref{lemma_n_cK_tree_bound}. However, this approach encounters other substantial difficulties, as discussed earlier in the introduction, and therefore we do not pursue it here.

		Nevertheless, this idea may prove useful in the study of regular RBMs  on the torus $(\Z/L\Z)^d$ in dimensions $d\ge 3$, which can be viewed heuristically as corresponding to power-law RBMs with $\al=2/d$. In this setting, the propagator exhibits the decay 
		\begin{equation}\label{eq:Dt_dge3}
			B_t(r):=\qa{W^{-d}(r/W+1)^{-(d-2)}+[N(1-t)]^{-1}}\exp(-cr/\ell_t),\quad \text{with}\quad \ell_t=W(1-t)^{-1/2},
		\end{equation}
		for some constant $c>0$, where $N=L^d$ is the volume of the torus. The tree-shaped bound \eqref{eq_n_cK_tree_bound} extends readily to this setting and therefore strengthens Lemma 2.16 in \cite{dubova2025delocalizationnonmeanfieldrandommatrices}, where only a $\max$-norm bound for the $\cK$-loops was obtained.  
		With such an improvement, combining the ideas developed in \cite{dubova2025delocalizationnonmeanfieldrandommatrices} with those of the present work would allow several refinements of the results in that paper, including establishing sharp decay estimates for $\cL$-loops of arbitrarily large order, treating more general variance profiles, and removing several technical assumptions.
	\end{remark}

	Another key ingredient is the following evolution kernel estimate, whose proof is deferred to \Cref{sec_proof_of_n_evolution_kernel_bound}.

	\begin{lemma}\label{lemma_n_evolution_kernel_bound}
		For any $n$-tensor $\cA^{\pa{n}}:\ZL^n\to \C$ with $n\geq 2$, $\bsigma\in\ha{-,+}^n$, and $0\leq s<t\leq \tf$, we have the evolution kernel bound
		\begin{equation}\label{evolution_kernel_bound_1}
			\begin{aligned}
				\norma{\bU_{s,t,\bsigma}^{\pa{n}}\circ \cA^{\pa{n}}}_{\mathscr{D}_t^{-1}}\leq \norma{\cU_{s,t}^{\otimes n}\circ |\cA|^{\pa{n}}}_{\mathscr{D}_t^{-1}} \prec \frac{\ell_t}{\ell_s}\qa{\frac{\ell_s(1-s)}{\ell_t(1-t)}}^{n} \normb{ \cA^{\pa{n}}}_{\mathscr{D}_s^{-1}},
			\end{aligned}
		\end{equation}
		where the absolute tensor $|\cA|$ is defined by  $|\cA|_{\bx}^{\pa{n}}:=\abs{\cA_{\bx}^{\pa{n}}}$, and $\cU_{s,t}^{\otimes n}$ denotes the $n$-fold tensor product of the matrix \smash{$\cU_{s,t}=\cU_{s,t}^{(-,+)}=(1-sS)/(1-tS)$} defined in \eqref{U_s_t_equal_1_Theta}. 
		
		Furthermore, the evolution kernel estimate can be improved in the following two cases.
		\begin{itemize}
			\item[(I)] If $\bsigma=(\sigma_1,\ldots,\sigma_n)$ is non-alternating, i.e., $\sigma_i=\sigma_{i+1}$ for some $i\in\qq{n}$, then  
			\begin{equation}\label{evolution_kernel_bound_2}
				\begin{aligned}
					\norma{\bU_{s,t,\bsigma}^{\pa{n}}\circ \cA^{\pa{n}}}_{\mathscr{D}_t^{-1}}\prec \qa{\frac{\ell_s(1-s)}{\ell_t(1-t)}}^{n-1} \normb{ \cA^{\pa{n}}}_{\mathscr{D}_s^{-1}}.
				\end{aligned}
			\end{equation}
			\item[(II)] If $\bsigma=(\sigma_1,\ldots,\sigma_n)$ is alternating, i.e., $\sigma_i\neq\sigma_{i+1}$ for all $i\in\qq{n}$, and the tensor $\cA^{(n)}$ satisfies the sum-zero property
			\begin{equation}\label{sum_zero_property_supercritical}
				\begin{aligned}
					\sum_{x_2,\ldots,x_n\in\ZL}\cA_{x_1x_2\cdots x_n}^{\pa{n}}=0,\quad \forall x_1\in\ZL,
				\end{aligned}
			\end{equation}
			then we have 
			\begin{equation}\label{evolution_kernel_bound_3}
				\begin{aligned}
					\norma{\bU_{s,t,\bsigma}^{\pa{n}}\circ \cA^{\pa{n}}}_{\mathscr{D}_t^{-1}}\prec \qa{\frac{\ell_s(1-s)}{\ell_t(1-t)}}^{n} \normb{ \cA^{\pa{n}}}_{\mathscr{D}_s^{-1}}.
				\end{aligned}
			\end{equation}
		\end{itemize}
	\end{lemma}

	\begin{remark}\label{rmk_failureofLoop}
		In the regime $\alpha\in(0,1)$, the failure of the evolution kernel estimate \eqref{evolution_kernel_bound_3} is the main reason for the breakdown of the current loop hierarchy method. Consequently, as discussed above, the proof in the regime $\alpha\in(0,1)$ can only rely on the analysis of $\cL$-loops of length $2$. On the other hand, for the method based solely on $2$-loops developed in Sections \ref{sec:pfal01}--\ref{sec_proof_of_que}, the evolution kernel bound for $2$-tensors must be at least as strong as
		\begin{equation}\label{eq:UstA}
			\norma{\bU_{s,t,\bsigma}^{\pa{2}}\circ \cA^{\pa{2}}}_{\mathscr{D}_t^{-1}}\prec \frac{B_t(0)}{B_s(0)} \normb{ \cA^{\pa{2}}}_{\mathscr{D}_s^{-1}}.
		\end{equation}
		In the regime $\alpha\in(0,1)$, such an estimate is available when $\cA^{(2)}$ satisfies the sum-zero property \eqref{sum_zero_property_supercritical}, as established in \Cref{lemma_U_s_t}. (The main reason is that $B_s(0)=B_t(0)$ for $s<t$ with $\ell_t<N$.) However, for $\alpha\in[1,\infty)$, even when the sum-zero property holds for $n=2$, we only have the weaker bound
		\begin{equation}
			\begin{aligned}
				\norma{\bU_{s,t,\bsigma}^{\pa{2}}\circ \cA^{\pa{2}}}_{\mathscr{D}_t^{-1}}\prec \pa{\frac{B_t(0)}{B_s(0)}}^2 \norma{ \cA^{\pa{2}}}_{\mathscr{D}_s^{-1}},
			\end{aligned}
		\end{equation}
		which is insufficient compared with \eqref{eq:UstA}. This prevents us from directly extending the argument for $\alpha\in(0,1)$ to the regime $\alpha\ge 1$, as already mentioned in \cref{footnoteloop}.
	\end{remark}

	\subsection{Proof of \Cref{local_laws_along_flow_supercritical}}
	
	\Cref{local_laws_along_flow_supercritical} follows as an immediate consequence of the following estimates on $\cL$-loops and $\cC$-chains of arbitrary length.
	
	\begin{lemma}\label{local_law_with_n_loops_along_flow_supercritical}
		In the setting of \Cref{local_laws_along_flow_supercritical}, the estimates \eqref{flow_local_law_supercritical}--\eqref{flow_expected_local_law_supercritical} hold. Moreover, for any fixed $n\geq 2$, the following estimates hold uniformly in $\bsigma\in\ha{-,+}^n$, $\bx\in\ZL^n$, and $t\in\qa{0,\tf}$:
		\begin{align}
			\pa{\cL-\cK}_{t,\bsigma,\bx}^{\pa{n}}&\prec B_t^{n}(0)\cdot \mathscr{D}_{t,\bx}^{\pa{n}},\label{flow_n_loops_local_law_supercritical}\\
			\pa{\cC-\cK^{\cC}}_{t,\bsigma,\bx}^{\pa{n}}&\prec B_t^{n-1/2}(0)\cdot \mathscr{D}_{t,\bx}^{\pa{n}},\label{flow_n_chains_local_law_supercritical}
		\end{align}
		where $\cK\equiv\cK^{\cL}$, and $\cC_{t,\bsigma,\bx}^{(n)}$ and $\cK^{\cC,\pa{n}}_{t,\bsigma,\bx}$ are defined in the sense of \eqref{abbreviation_diagonal_cC_chain} and \eqref{abbreviation_diagonal_cK_chain}, respectively.
	\end{lemma}
	As in the proof of \Cref{theorem_local_laws_along_flow}, we establish \Cref{local_law_with_n_loops_along_flow_supercritical} using a bootstrap argument along the flow. Note that at time $t=0$, the left-hand sides of \eqref{flow_local_law_supercritical}--\eqref{flow_expected_local_law_supercritical}, \eqref{flow_n_loops_local_law_supercritical}, and \eqref{flow_n_chains_local_law_supercritical} all vanish, since $G_0(\sigma)=m(\sigma)I$. Therefore, \Cref{local_law_with_n_loops_along_flow_supercritical} follows from the following lemma combined with a standard induction argument.

	\begin{lemma}\label{lemma_bootstrap_supercritical}
		In the setting of \Cref{local_law_with_n_loops_along_flow_supercritical}, suppose that the estimates \eqref{flow_local_law_supercritical}, \eqref{flow_loop_local_law_supercritical}, and \eqref{flow_n_loops_local_law_supercritical} hold at some time $s\in\qa{0,\tf}$. Then, for any $t\in\qa{s,\tf}$ satisfying
		\begin{equation}\label{bootstrap_assumption_supercritical}
			\begin{aligned}
				B_t^{1/100}(0)\leq \frac{1-t}{1-s}\leq 1,
			\end{aligned}
		\end{equation}
		the estimates \eqref{flow_local_law_supercritical}, \eqref{flow_loop_local_law_supercritical}, and \eqref{flow_n_loops_local_law_supercritical} also hold at time $t$.  Furthermore, if the chain estimate \eqref{flow_n_chains_local_law_supercritical} holds at time $s\in\qa{0,\tf}$, then it also holds at any time $t\in\qa{s,\tf}$ satisfying \eqref{bootstrap_assumption_supercritical}. Finally, if both \eqref{flow_n_chains_local_law_supercritical} and the expected estimate \eqref{flow_expected_local_law_supercritical} hold at time $s$, then they also hold at any time $t\in\qa{s,\tf}$ satisfying \eqref{bootstrap_assumption_supercritical}.
	\end{lemma}
	
	\begin{proof}[\bf Proof of \Cref{local_law_with_n_loops_along_flow_supercritical}]
		We iterate \Cref{lemma_bootstrap_supercritical} from $t=0$ up to $\tf$ along an increasing sequence of times $\{t_k:0\le k \le n\}$ with $t_0=0$ and $t_n=\tf$, chosen such that $n=\OO(1)$ and $(1-t_k)/(1-t_{k-1})\ge [B_{t_k}(0)]^{1/100}$. After $n$ induction steps, we obtain \Cref{local_law_with_n_loops_along_flow_supercritical} for any time $t\in[0,\tf]$. Finally, a standard $N^{-C}$-net argument combined with perturbation estimates extends these bounds uniformly to all $t\in[0,\tf]$.
	\end{proof}

	Fix a time $t$ satisfying \eqref{bootstrap_assumption_supercritical}. We prove \Cref{lemma_bootstrap_supercritical} through the following six steps, which originate from \cite{Band1D} and have been used extensively in the study of regular random band matrices and block Anderson models  \cite{Band2D,truong2025localizationlengthfinitevolumerandom,fan2025blockreductionmethodrandom,dubova2025delocalizationnonmeanfieldrandommatrices}. 
	
	\medskip
	\noindent
	{\bf Step 1} (Weak bounds for $\cL$-loops and weak local law): The $\cL$-loops satisfy the following a priori bound:
	\begin{equation}\label{step_1_supercritical_continuity_max_bound}
		\begin{aligned}
			\max_{\bsigma\in\ha{-,+}^n}\max_{\bx\in\ZL^n}\absa{\cL_{u,\bsigma,\bx}^{\pa{n}}}\prec \pa{\frac{1-s}{1-u}}^{n-1}B_s^{n-1}(0),\quad \forall u\in\qa{s,t}.
		\end{aligned}
	\end{equation}
	In addition, the following weak local laws hold: 
	\begin{equation}\label{step_1_supercritical_weak_local_law}
		\norm{G_u-m}_{\max}\prec B_u^{1/4}(0),\quad \max_{x\in\ZL}\,\absb{\avgb{(G_u-m)S^{(x)}}}\prec B_u^{3/4}(0),\quad \forall u\in\qa{s,t}.
	\end{equation}

	\noindent
	{\bf Step 2} (Sharp local law and weak estimates for short $\cL$-loops): For any $x,y\in\ZL$, we have
	\begin{equation}\label{step_2_sharp_local_law_supercritical}
		\begin{aligned}
			\absa{(G_u-m)_{xy}}^2\prec B_u\pa{|x-y|},\quad |\avgb{(G_u-m)S^{\pa{x}}}|\prec B_u\pa{0},\quad \forall u\in\qa{s,t}.
		\end{aligned}
	\end{equation}
	Moreover, we have the following weak estimates for $\bsigma\in\ha{-,+}^n$ and $\bx\in\ZL^n$ with $n\in\ha{2,3,4}$:
	\begin{equation}\label{step_2_loops_local_law_supercritical}
		\begin{aligned}
			\absa{(\cL-\cK)_{u,\bsigma,\bx}^{\pa{n}}}\prec \qa{B_u(0)}^{n-(n+2)/(n+3)}\cdot \mathscr{D}_{u,\bx}^{\pa{n}}, \qquad \forall u\in\qa{s,t}.
		\end{aligned}
	\end{equation}
	In particular, together with \Cref{lemma_n_cK_bound}, this implies for $n=2,3,4$ that
	\begin{equation}\label{step_2_loops_bound_supercritical}
		\begin{aligned}
			\max_{\bsigma\in\ha{-,+}^n}\absa{\cL_{u,\bsigma,\bx}^{\pa{n}}}\prec B_u^{n-1}(0)\cdot \mathscr{D}_{u,\bx}^{\pa{n}},\qquad \forall u\in\qa{s,t}.
		\end{aligned}
	\end{equation}

	\noindent
	{\bf Step 3} (Sharp upper bound for $\cL$-loops): For any fixed $n\in\N$, we have the bound 
	\begin{equation}\label{step_3_n_bound_cL_supercritical}
		\begin{aligned}
			\max_{\bsigma\in\ha{-,+}^n}\absa{\cL_{u,\bsigma,\bx}^{\pa{n}}}\prec B_u^{n-1}(0)\cdot \mathscr{D}_{u,\bx}^{\pa{n}},\qquad \forall u\in\qa{s,t},\ \bx\in\ZL^n.
		\end{aligned}
	\end{equation}

	\noindent
	{\bf Step 4} (Sharp estimates for $\cL$-loops and $T$-variables): For any fixed $n\in\N$ and all $\bx\in\ZL^n$, we have 
	\begin{equation}\label{step_4_supercritical_sharp_local_law_loops}
		\max_{\bsigma\in\ha{-,+}^n}\absa{(\cL-\cK)_{u,\bsigma,\bx}^{\pa{n}}}\prec B_u^{n}(0)\cdot \mathscr{D}_{u,\bx}^{\pa{n}},\qquad \forall u\in\qa{s,t}.
	\end{equation}
	Moreover, for any $x,y,y'\in \ZL$, we have
	\begin{equation}\label{step_4_local_law_T_supercritical}
		\max_{\bsigma\in\ha{-,+}^2}\absa{(T-\Theta)_{u,x,yy'}^{\bsigma}}\prec B_u^{3/2}(0)\cdot \mathscr{D}_{u,x}^{\pa{2}}(y,y'),\qquad \forall u\in\qa{s,t}.
	\end{equation}
	
	\noindent
	{\bf Step 5} (Sharp chain estimates): For any fixed $n\in\N$, the following estimate holds for all $(\bx,x,y)\in\ZL^{n+1}$:
	\begin{equation}\label{step_5_supercritical_sharp_local_law_chains}
		\begin{aligned}
			\max_{\bsigma\in\ha{-,+}^n}\absa{\pa{\cC-\cK^{\cC}}_{u,\bsigma,\bx}^{\pa{n}}(x,y)}\prec B_u^{n-1/2}(0)\cdot \mathscr{D}_{u,\bx}^{\pa{n}}(x,y),\qquad \forall u\in\qa{s,t}.
		\end{aligned}
	\end{equation}

	\noindent
	{\bf Step 6} (Expected estimate for $2$-$\cL$-loops): For all $\bx\in\ZL^2$, we have
	\begin{equation}\label{step_6_expectedloop}
		\max_{\bsigma\in\ha{-,+}^2}\absa{\E\pa{\cL-\cK}_{u,\bsigma,\bx}^{\pa{2}}}\prec B_u^{5/2}(0)\cdot \mathscr{D}_{u,\bx}^{\pa{2}},\qquad \forall u\in\qa{s,t}.
	\end{equation}

	We remark that in the previous works 	\cite{Band1D,Band2D,truong2025localizationlengthfinitevolumerandom,fan2025blockreductionmethodrandom,dubova2025delocalizationnonmeanfieldrandommatrices},
	the six-step scheme is slightly different from ours. In those works, Step 5 establishes sharp estimates for $2$-loops with optimal decay, whereas sharp estimates for chains, as in our Step 5, are not derived. 
	In our approach, the original Step 5 can be skipped because, unlike in the above references, our techniques already allow us to track sharp decay for $\cL$-loops of arbitrary length in Step 4. This idea bears some similarity to the approach in \cite{erdos2025zigzagstrategyrandomband}, which applies in the regime $\alpha\ge3$. However, our argument differs in several key aspects, as discussed in the introduction, in order to handle the regime $1\le\alpha<3$. In particular, we decouple the dynamics of the $\cL$-loops from that of the $\cC$-chains using \Cref{lemma_Xi_reshape}, which provides almost sharp upper bounds on $\cC$-chains in terms of $\cL$-loops. 
    The sharp estimates for the chains are then obtained from the sharp \(\cL\)-loop estimates through the reshaping formula \eqref{reshape_chains_into_loops}, which is derived using resolvent identities and standard large deviation bounds. 
    Finally, our analysis relies on the a priori bound \eqref{step_1_supercritical_continuity_max_bound}, together with an extended version of the \(\contract\) inequality established in \Cref{lemma_n_loop_contraction}.

	The proof of Step 1 is essentially identical to that in \cite[Section 5.1]{Band1D}, up to minor notational changes. We therefore omit it here and focus on Steps 2--6 in the remainder of this section. Note that by \eqref{bootstrap_assumption_supercritical},
	\begin{align*}
		B_t(0)\leq \pa{\frac{1-t}{1-s}}^{100}\leq \pa{\frac{1-t}{1-s}}^{100}\pa{\frac{\ell_t}{\ell_s}}^{100}= \qa{\frac{B_s(0)}{B_t(0)}}^{100},
	\end{align*}
	which implies
	\begin{equation}\label{bound_B_t_to_B_s}
		\begin{aligned}
			B_t(0)\le B_s^{1-1/101}(0).
		\end{aligned}
	\end{equation}
	This bound will be used tacitly in the analysis below.

	\subsection{$\Xi$-parameters}\label{sec_Xi_parameter} 
	Before proceeding to the proofs of Steps 2--6, for clarity of presentation we first introduce five types of $\Xi$-parameters.
	\begin{itemize}
		\item The \emph{upper bound control parameters} are defined for any $n\in\N$ and $u\in\qa{s,t}$ as
		\begin{equation}\label{eq:Xi_Cd_Cod}
			\begin{aligned}
				\Xi_{u,n}^{\cC_{\rd}}:=&1+\sup_{\bsigma\in\ha{-,+}^n}\sup_{\pa{\bx,x}\in\ZL^n}\absb{\cC_{u,\bsigma,\bx}^{\pa{n}}\pa{x}}\big/(B_u^{n-1}(0)\mathscr{D}_{u,(\bx,x)}^{\pa{n}}),\\
				\Xi_{u,n}^{\cC_{\txt{od}}}:=&1+\sup_{\bsigma\in\ha{-,+}^n}\sup_{\pa{\bx,x,y}\in\ZL^{n+1}: x\neq y}\absb{\cC_{u,\bsigma,\bx}^{\pa{n}}\pa{x,y}}\big/(B_u^{n-1/2}(0)\mathscr{D}_{u,\bx}^{\pa{n}}\pa{x,y}),\\
				\Xi_{u,n}^{\cL}:=&1+\sup_{\bsigma\in\ha{-,+}^n}\sup_{\bx\in\ZL^n}\absb{\cL_{u,\bsigma,\bx}^{\pa{n}}}\big/(B_u^{n-1}(0)\mathscr{D}_{u,\bx}^{\pa{n}}).
			\end{aligned}
		\end{equation}
		\item The \emph{stability control parameters} are defined for any $n\in\N$ and $u\in\qa{s,t}$ as
		\begin{align}
			\Xi_{u,n}^{\cL-\cK}&:=1+\sup_{\bsigma\in\ha{-,+}^n}\sup_{\bx\in\ZL^n}\absb{(\cL-\cK)_{u,\bsigma,\bx}^{\pa{n}}}\big /(B_u^{n}(0)\mathscr{D}_{u,\bx}^{\pa{n}}),\label{eq;XiL-Kn}\\
            \Xi_{u,n}^{\cC-\cK^{\cC}}&:=1+\sup_{\bsigma\in\ha{-,+}^n}\sup_{\bx\in\ZL^n}\absb{(\cC-\cK^{\cC})_{u,\bsigma,\bx}^{\pa{n}}}\big /(B_u^{n-1/2}(0)\mathscr{D}_{u,\bx}^{\pa{n}}).\label{eq;XiC-KCn}
		\end{align}
	\end{itemize}
	
	We next record several static relations among these $\Xi$-parameters at a fixed time $u$, which will play a crucial role in the subsequent proofs. First, by definition and the bound on $\cK$ from \Cref{lemma_n_cK_bound}, we have
	\begin{equation}\label{eq;Xi_relate1}
		\Xi_{u,n}^{\cL}\prec 1+ B_u(0)\Xi_{u,n}^{\cL-\cK}\prec \Xi_{u,n}^{\cL},\qquad \Xi_{u,n}^{\cC_{\rd}}\prec 1+ B_u^{1/2}(0)\Xi_{u,n}^{\cC-\cK^{\cC}}\prec \Xi_{u,n}^{\cC_{\rd}}.
	\end{equation}
	Second, consider any $n$-loop $\cL_{u,\bsigma,\bx}^{\pa{n}}$ with $\bx=(x_1,\ldots,x_n)$. Using the basic bound
	\begin{equation*}
		\begin{aligned}
			\sum_{x\in\ZL}S_{a x}\mathscr{D}_{u}(|b-x|)\mathscr{D}_{u}(|x-c|)\prec \mathscr{D}_{u}(|b-a|)\mathscr{D}_u(|a-c|),
		\end{aligned}
	\end{equation*}
	we obtain, for $\bx'=(x_1,\ldots,x_{n-1})$, that
	\begin{equation}\label{eq;Xi_relate2}
		\cL_{u,\bsigma,\bx}^{\pa{n}}=\sum_{x\in\ZL}S_{x_n x}\cC_{u,\bsigma,(\bx',x)}^{\pa{n}}\prec B_u^{n-1}(0)\mathscr{D}_{u,\bx}^{\pa{n}}\Xi_{u,n}^{\cC_{\rd}}\ \ \implies \ \ \Xi_{u,n}^{\cL}\prec \Xi_{u,n}^{\cC_{\rd}}.
	\end{equation}
    Moreover, the CS inequality yields the following basic relations between the upper-bound control parameters. 
    
	\begin{lemma}\label{lemma_loop_bound_cut_long_to_short}
		For any fixed integer $k\ge1$ and $\star\in\ha{\cC_{\rd},\cL}$, we have
		\begin{equation}\label{loop_bound_cut_long_to_short}
			\begin{aligned}
				&B_u(0)\Xi_{u,4k}^{\star}\lesssim \Xi_{u,2k}^{\star}\Xi_{u,2k}^{\cL}, \quad &B_u(0)\Xi_{u,4k+2}^{\star}\lesssim \min\pa{ \Xi_{u,2k}^{\star}\Xi_{u,2k+2}^{\cL},\Xi_{u,2k}^{\cL}\Xi_{u,2k+2}^{\star}},\\
				&\pa{\Xi_{u,2k+1}^{\star}}^2\lesssim \Xi_{u,2k}^{\star}\Xi_{u,2k+2}^{\star},\quad &B_u(0)\pa{\Xi_{u,2k+1}^{\star}}^2\lesssim \Xi_{u,2}^{\cL}\pa{\Xi_{u,2k}^{\star}}^2.
			\end{aligned}
		\end{equation}
		Since these bounds hold uniformly for $u\in\qa{s,t}$, they remain valid when $u$ is replaced by any stopping time $\tau$ satisfying $s\le\tau\le t$.
	\end{lemma}

\begin{proof}
		We present the proof only for the case $\star=\cC_{\rd}$; the case $\star=\cL$ can be handled in the same way. To prove the first bound in \eqref{loop_bound_cut_long_to_short}, consider any $4k$-chain $\smash{\cC_{u,\bsigma,\bx}^{\pa{4k}}(x)}$ with $\bsigma=(\sigma_1,\ldots,\sigma_{4k})\in\{-,+\}^{4k}$ and \smash{$\bx=(x_1,\ldots,x_{4k-1})\in\ZL^{4k-1}$}. Without loss of generality, we may assume that the chain is symmetric, i.e., $x_{p}=x_{q}$ whenever $p+q=4k$, and $\sigma_{p}=-\sigma_q$ whenever $p+q=4k+1$, in which case the chain is non-negative. Otherwise, we apply the CS inequality to bound it by two symmetric chains:
		\begin{equation}\label{eq:CS-chain}
			\absa{\cC_{u,\bsigma,\bx}^{\pa{4k}}(x)}^2\leq \cC_{u,\bsigma',\bx'}^{\pa{4k}}(x)\cC_{u,\bsigma'',\bx''}^{\pa{4k}}(x),
		\end{equation}
		where we denote
		\begin{equation*}
			\begin{aligned}
				\bsigma':=(\sigma_1,\ldots,\sigma_{2k},-\sigma_{2k},\ldots,-\sigma_1),\quad &\bx':=(x_1,\ldots,x_{2k-1},x_{2k},x_{2k-1},\ldots,x_1),\\
				\bsigma'':=(-\sigma_{4k},\ldots,-\sigma_{2k+1},\sigma_{2k+1},\ldots,\sigma_{4k}),\quad &\bx'':=(x_{4k-1},\ldots,x_{2k+1},x_{2k},x_{2k+1},\ldots,x_{4k-1}).
			\end{aligned}
		\end{equation*}
		Under the symmetry assumption, we write     
        \begin{equation}\label{eq:CS-to-bound-chain}
			\begin{aligned}
				\cC_{u,\bsigma,\bx}^{\pa{4k}}(x)&=\pa{\sC S^{\pa{x_k}}\mathsf{D}S^{\pa{x_k}}\sC^*}_{xx}\leq (\sC S^{\pa{x_k}}\sC^*)_{xx} \norma{\pb{S^{\pa{x_k}}}^{1/2}\mathsf{D}\pb{S^{\pa{x_k}}}^{1/2}}\\
				&\leq (\sC S^{\pa{x_k}}\sC^*)_{xx}\txt{Tr}\,\pb{\mathsf{D}S^{\pa{x_k}}} \lesssim B_u^{4k-2}(0)\Xi_{u,2k}^{\cC_{\rd}}\Xi_{u,2k}^{\cL}\cdot \mathscr{D}_{u,\bx}^{\pa{4k}}(x),
			\end{aligned}
		\end{equation}
		where we denote
		\begin{equation*}
			\begin{aligned}
				\sC:=G_u(\sigma_1)S^{\pa{x_1}}\cdots G_u(\sigma_{k-1})S^{(x_{k-1})} G_u(\sigma_k),\quad \mathsf{D}:=G_u(\sigma_{k+1})S^{\pa{x_{k+1}}}\cdots G_u(\sigma_{3k-1})S^{(x_{3k-1})}G_u(\sigma_{3k}).
			\end{aligned}
		\end{equation*}
		This proves the first bound in \eqref{loop_bound_cut_long_to_short}. The second bound can be established in the same manner, while the third bound follows immediately from a Cauchy-Schwarz argument analogous to \eqref{eq:CS-chain}.

		Finally, we prove the last bound in \eqref{loop_bound_cut_long_to_short}. Consider an arbitrary chain $\smash{\cC_{u,\bsigma,\bx}^{\pa{2k+1}}(x)}$ with $\bsigma=(\sigma_1,\ldots,\sigma_{2k+1})\in\{-,+\}^{2k+1}$ and $\bx=(x_1,\ldots,x_{2k})\in\ZL^{2k}$. We write it in the form $\p{\sC_1S^{(x_k)}G_u(\sigma_{k+1})S^{(x_{k+1})}\sC_2}_{xx}$, where $\sC_1$ and $\sC_2$ are two $k$-chains. Then, we have 
		\begin{equation*}
			\begin{aligned}
				\absa{\cC_{u,\bsigma,\bx}^{\pa{2k+1}}(x)}^2             &\leq \pb{\sC_1S^{(x_k)}\sC_1^*}_{xx}\pb{\sC_2^*S^{(x_{k+1})}\sC_2}_{xx}\norma{\pb{S^{(x_k)}}^{1/2}G_u(\sigma_{k+1})\pb{S^{(x_{k+1})}}^{1/2}}^2\\
				&\leq \pb{\sC_1S^{(x_k)}\sC_1^*}_{xx}\pb{\sC_2^*S^{(x_{k+1})}\sC_2}_{xx}\avg{G_u(\sigma_{k+1})S^{(x_{k+1})}G_u(-\sigma_{k+1})S^{(x_k)}}\\
				&\lesssim B_u^{4k-1}(0)\Xi_{u,2}^{\cL}\pb{\Xi_{u,2k}^{\cC_{\rd}}}^2\cdot \qa{\mathscr{D}_{u,\bx}^{(2k+1)}(x)}^2.
			\end{aligned}
		\end{equation*}
		This proves the final bound in \eqref{loop_bound_cut_long_to_short}.
	\end{proof}

	Adopting the idea of \cite[Lemma A.2]{Band1D}, we establish the following lemma, which allows us to bound diagonal and off-diagonal $\cC$-chains inductively in terms of $\cL$-loops. This result will play a central role in controlling the chain factors arising in the loop hierarchy.\footnote{It is worth noting that \cite[Lemma A.2]{Band1D} is not essentially used in the proof there, since a block variance profile was assumed in \cite{Band1D}. In contrast, our arguments rely crucially on such a result in \Cref{lemma_Xi_reshape}, due to the general variance profile assumption \eqref{alpha_decay}.} Similar to \cite[Lemma A.2]{Band1D}, the proof of Lemma \ref{lemma_Xi_reshape} relies on repeated applications of resolvent identities to decouple the randomness, together with standard large deviation bounds. We therefore defer the proof to \Cref{sec_proof_of_Xi_bounding_relations}.

	\begin{lemma}\label{lemma_Xi_reshape}
		Under the assumptions of \Cref{lemma_bootstrap_supercritical}, suppose that the weak local law \eqref{step_1_supercritical_weak_local_law} holds uniformly for $u\in\qa{s,t}$. Fix a constant $c\in(0,1/100)$ and define the stopping time
		\begin{equation}\label{def_tau_cL}
			\tau_{\cL}\equiv \tau_{\cL}(c):=t\wedge\inf\ha{u\geq s: B_t^{c}(0) \Xi_{u,2}^{\cL}> 1},
		\end{equation}
        where we adopt the convention $\inf\emptyset=\infty$. 
		Then the following bounds hold.
		\begin{itemize}
			
			\item ({\bf Diagonal chains}) For any even $n\in 2\N$, we have
			\begin{align}
				\qquad \Xi_{u\wedge \tau_{\cL},n}^{\cC_{\txt{d}}}\prec &~ \Xi_{u\wedge \tau_{\cL},n}^{\cL}+B_{u\wedge \tau_{\cL}}(0)\sum_{\ell\ge 2}\sum_{2\leq n_1,\ldots,n_\ell\leq n-2: n_1+\cdots +n_\ell=n+\ell}\Xi_{u\wedge \tau_{\cL},n_1}^{\cC_{\rd}}\cdots \Xi_{u\wedge \tau_{\cL},n_\ell}^{\cC_{\rd}}\label{eq_cC_d_to_loops}\\
				&~+\sum_{\ell\ge 2} \sum_{2\leq n_1,\ldots,n_\ell\leq n-2: n_1+\cdots +n_\ell=n+\ell-1}\Xi_{u\wedge \tau_{\cL},n_1}^{\cC_{\rd}}\cdots \Xi_{u\wedge \tau_{\cL},n_\ell}^{\cC_{\rd}}.\nonumber
			\end{align}

			\item ({\bf Off-diagonal chains}) For off-diagonal $n$-chains with $n\ge 2$, we have
			\begin{align}
				\qquad &\Xi_{u\wedge \tau_{\cL},n}^{\cC_{\txt{od}}}\prec \pa{\Xi_{u\wedge\tau_{\cL},2n}^{\cC_{\txt{d}}} }^{1/2} +\sum_{\ell\ge 2}\sum_{n_1\geq 1,n_2,\ldots,n_\ell\geq 2:n_1+\cdots+n_\ell=n+\ell-1}\Xi_{u\wedge\tau_{\cL},n_1}^{\cC_{\txt{od}}}\Xi_{u\wedge\tau_{\cL},n_2}^{\cC_{\rd}}\cdots \Xi_{u\wedge\tau_{\cL},n_\ell}^{\cC_{\rd}}\label{eq_cC_od_to_cC_d}\\
				&+B_{u\wedge\tau_{\cL}}^{1/2}(0)\sum_{\ell\ge 2}\sum_{1\le n_1,n_2\leq n-1,n_3,\ldots,n_\ell\geq 2: n_1+\cdots+n_\ell=n+\ell-1}\Xi_{u\wedge\tau_{\cL},n_1}^{\cC_{\txt{od}}}\Xi_{u\wedge\tau_{\cL},n_2}^{\cC_{\txt{od}}}\Xi_{u\wedge\tau_{\cL},n_3}^{\cC_{\rd}}\cdots \Xi_{u\wedge\tau_{\cL},n_\ell}^{\cC_{\rd}}.\nonumber
			\end{align}

			\item ({\bf Reshaping $T$-variables into loops}) For any $u\in\qa{s,t}$ and $\bsigma=(\sigma_1,\sigma_2)\in\ha{-,+}^2$, the variable $T_{u,xy}^{\bsigma}$ can be replaced by an $\cL$-loop up to an error term as follows:
			\begin{align}
				&T_{u,xy}^{\bsigma}-m\p{\sigma_1}m\p{\sigma_2}\pa{u\cL_{u,xy}^{\bsigma}+S_{xy}}\prec B_u^{2}(0)\mathscr{D}_{u}^2(|x-y|)\qa{\Xi_{u,3}^{\cC_{\rd}}+ \pb{\Xi_{u,2}^{\cC_{\rd}}}^2}\label{reshape_2_chains_to_loops}\\
				&+B_u^{3/2}(0)\mathscr{D}_{u}^2(|x-y|)\pbb{\pa{\Xi_{u,2}^{\cL}}^{3}+\Xi_{u,4}^{\cL}+B_u(0)\sum_{\ell\ge 1, n_1,\ldots,n_\ell\geq 2:n_1+\cdots+n_\ell=4+\ell}\Xi_{u,n_1}^{\cC_{\rd}}\cdots \Xi_{u,n_\ell}^{\cC_{\rd}}}^{1/2}.\nonumber
			\end{align}

            \item ({\bf Reshaping chains into loops})
Suppose that $\Xi_{u,k}^{\star}\prec 1$ for any fixed $k\geq 1$, $u\in[s,t]$, and $\star\in\ha{\cC_{\rd},\cL}$. Let $n\geq 3$, and assume that, for all $u\in[s,t]$, 
            \[\Xi_{u,k}^{\cL-\cK}\prec 1, \quad \forall 1\leq k\leq n,\quad \text{and}\quad \Xi_{u,k}^{\cC-\cK^{\cC}}\prec 1, \quad \forall 1\leq k\leq n-1.\]
Then we have
            \begin{equation}\label{reshape_chains_into_loops}
                \Xi_{u,n}^{\cC-\cK^{\cC}}\prec 1,\qquad \forall u\in[s,t].
            \end{equation}
		\end{itemize}
		Finally, for any fixed $n\in\N$, a standard $N^{-C}$-net argument shows that these estimates hold uniformly for all $u\in\qa{s,t}$. In particular, they remain valid if $u$ is replaced by any stopping time $\tau$ satisfying $s\le \tau\le t$.
	\end{lemma}

	\subsection{Step 2: Sharp local law and weak estimates for short $\cL$-loops}\label{sec_step_2_supercritical}
	
	The sharp local law \eqref{step_2_sharp_local_law_supercritical} follows immediately from the sharp bound \eqref{step_2_loops_bound_supercritical} (which itself follows from \eqref{step_2_loops_local_law_supercritical}) together with \Cref{2_loop_to_1_chain}. It therefore remains to establish \eqref{step_2_loops_local_law_supercritical} for $n\in\{2,3,4\}$. To this end, we work with $\cL$-loops of length at most $4$ and define, for $n=2,3,4$, the stopping times
	\begin{equation}\label{eq:tau234}
	\tau_{n}:=t\wedge \inf\ha{u\geq s: \q{B_u(0)}^{(n+2)/(n+3)}\Xi_{u,n}^{\cL-\cK}> 1},
	\end{equation}
 with the convention $\inf\emptyset=\infty$. 
	We then set $\tau:=\tau_2\wedge \tau_3\wedge \tau_4$. By \Cref{lemma_n_cK_bound}, for any $n\in\ha{2,3,4}$, we have 
	\begin{equation}\label{step_2_cL_bound_supercritical}
		\begin{aligned}
		 \Xi_{u\wedge\tau,n}^{\cL}\lesssim 1,
		\end{aligned}
	\end{equation}
	which implies that $\tau\leq \tau_{\cL}$, where $\tau_{\cL}$ is defined in \eqref{def_tau_cL}. Using \eqref{step_2_cL_bound_supercritical}, we apply \eqref{eq_cC_d_to_loops} inductively to bound $\Xi_{u\wedge\tau,2}^{\cC_{\rd}}$ and $\Xi_{u\wedge\tau,4}^{\cC_{\rd}}$, and use the third inequality in \eqref{loop_bound_cut_long_to_short} to control $\Xi_{u\wedge\tau,3}^{\cC_{\rd}}$. This yields
	\begin{equation}\label{step_2_cC_rd_bound_supercritical}
		\begin{aligned}
			\Xi_{u\wedge\tau,n}^{\cC_{\rd}}\prec 1,\qquad \forall n\in\{2,3,4 \}.
		\end{aligned}
	\end{equation}
	Applying \Cref{2_loop_to_1_chain} and \eqref{eq_cC_od_to_cC_d} inductively then gives
	\begin{equation}\label{step_2_cC_od_bound_supercritical}
		\begin{aligned}
			\Xi_{u\wedge \tau, n}^{\cC_{\txt{od}}}\prec 1,\qquad \forall n\in\{1,2\}.
		\end{aligned}
	\end{equation}
	In the following proof, we will also need a bound on $\Xi_{u\wedge\tau,5}^{\cC_{\rd}}$, which follows from the last inequality in \eqref{loop_bound_cut_long_to_short}:
	\begin{equation}\label{step_2_cC_rd_bound_5_supercritical}
		\begin{aligned}
			\Xi_{u\wedge\tau,5}^{\cC_{\rd}}\prec B_{u\wedge \tau}^{-1/2}(0).
		\end{aligned}
	\end{equation}
	With these bounds, the estimate \eqref{reshape_2_chains_to_loops} at time $u\wedge \tau$ reduces to
	\begin{equation*}
		\begin{aligned}
			T_{u\wedge \tau,xy}^{\bsigma}=& m\p{\sigma_1}m\p{\sigma_2}\pa{u\wedge \tau\cdot\cL_{u\wedge \tau,xy}^{\bsigma}+S_{xy}}+\opr{B_{u\wedge \tau}^{3/2}(0)\mathscr{D}_{u\wedge \tau}^2(|x-y|)}.
		\end{aligned}
	\end{equation*}
	Subtracting $\Theta_{u\wedge \tau,xy}^{\bsigma}$ from both sides yields
	\begin{equation}\label{step_2_T_Theta_supercritical}
		\begin{aligned}
			(T-\Theta)_{u\wedge \tau,xy}^{\bsigma}&= u\wedge \tau\cdot m\p{\sigma_1}m\p{\sigma_2}(\cL-\cK)_{u\wedge \tau,xy}^{\bsigma}+\opr{B_{u\wedge \tau}^{3/2}(0)\mathscr{D}_{u\wedge \tau}^2(|x-y|)}\\
			&\prec B_{u\wedge \tau}^{6/5}(0)\mathscr{D}_{u\wedge \tau}^2(|x-y|),
		\end{aligned}
	\end{equation}
	where the second step follows from the definitions \eqref{eq:tau234} (with $n=2$) and \eqref{eq;XiL-Kn}.
	
	Next, we bound $\Xi_{r\wedge \tau, n}^{\cL-\cK}$ for $n\in\{2,3,4\}$ and $r\in\qa{s,t}$, using the estimates \eqref{step_2_cL_bound_supercritical}--\eqref{step_2_T_Theta_supercritical} as inputs. For any time $r\in\qa{s,t}$, it follows from \eqref{evolution_n_cL_cK_integrated} that
	\begin{align}
		(\mathcal{L} - \mathcal{K})^{(n)}_{r\wedge \tau, \boldsymbol{\sigma}, \bx}=&\pa{\bU_{s,r\wedge \tau,\bsigma}^{\pa{n}}\circ \pa{\cL-\cK}_{s,\bsigma}^{\pa{n}}}_{\bx}+\sum_{l_{\cK}=3}^{n}\int_{s}^{r\wedge \tau}\pa{\bU_{u,r\wedge \tau,\bsigma}^{\pa{n}}\circ\qa{\mathfrak{D}_{l_{\cK}}(\cL-\cK)}_{u,\bsigma}^{\pa{n}}}_{\bx}\, \rd u \label{eq:evolution_n_rtau}\\
		&+ \int_{s}^{r\wedge \tau}\pa{\bU_{u,r\wedge \tau,\bsigma}^{\pa{n}}\circ\mathcal{E}^{(n)}_{u, \boldsymbol{\sigma}}}_{\bx}\dd u  +\int_{s}^{r\wedge \tau}\pa{\bU_{u,r\wedge \tau,\bsigma}^{\pa{n}}\circ\mathcal{W}^{(n)}_{u, \boldsymbol{\sigma}}}_{\bx}
		\dd u + \int_{s}^{r\wedge \tau}\pa{\bU_{u,r\wedge \tau,\bsigma}^{\pa{n}}\circ
			\dd\mathcal{B}^{(n)}_{u, \boldsymbol{\sigma}}}_{\bx}
		.\nonumber
	\end{align}
	By the induction hypothesis on $(\cL-\cK)_{s,\bsigma}^{\pa{n}}$, and using the evolution kernel bound \eqref{evolution_kernel_bound_1} together with \eqref{bootstrap_assumption_supercritical}, the first term can be estimated as
	\begin{equation}\label{step_2_error_leading_term_n_2}
		\begin{aligned}
			\pa{\bU_{s,r\wedge \tau,\bsigma}^{\pa{n}}\circ \pa{\cL-\cK}_{s,\bsigma}^{\pa{n}}}_{\bx}\prec         \p{{\ell_{r\wedge \tau}}/{\ell_s}}\cdot B_{r\wedge \tau}^{n}(0)\mathscr{D}_{r\wedge \tau,\bx}^{(n)}\le B_{r\wedge \tau}^{n-1/5}(0)\mathscr{D}_{r\wedge \tau,\bx}^{(n)}.
		\end{aligned}
	\end{equation}
	For the quadratic error term with $n=2$ and $\bx=(x_1,x_2)$, using \eqref{step_2_T_Theta_supercritical} and the definition of $\tau_2$, we obtain 
	\[\mathcal{E}^{(2)}_{u, \boldsymbol{\sigma},\bx} \prec B_u^{12/5}(0)\sum_a \mathscr{D}_u^2(|x_1-a|)\mathscr{D}_u^2(|a-x_2|) \prec \ell_u{B_u^{12/5}(0)}\mathscr{D}_u^2(|x_1-x_2|)=\frac{B_u^{7/5}(0)}{1-u}\cdot \mathscr{D}_u^2(|x_1-x_2|).\]
	Together with \eqref{evolution_kernel_bound_1} and \eqref{bootstrap_assumption_supercritical}, this implies
	\begin{equation}\label{step_2_error_double_difference_n_2}
		\begin{aligned}
			\int_{s}^{r\wedge \tau}\pa{\bU_{u,r\wedge \tau,\bsigma}^{\pa{2}}\circ\mathcal{E}^{(2)}_{u, \boldsymbol{\sigma}}}_{\bx}\dd u\prec {\frac{1-s}{1-r\wedge \tau}} {B_{r\wedge \tau}^{7/5}(0)\mathscr{D}_{r\wedge \tau,\bx}^{(2)}}\le B_{r\wedge \tau}^{6/5+1/10}(0)\mathscr{D}_{r\wedge \tau,\bx}^{(2)}.
		\end{aligned}
	\end{equation}
	When $n\in\ha{3,4}$, we write each term in the quadratic error term as a product of a $(\cL-\cK)$-loop and a $(\cC-\cK^{\cC})$-chain, arranged so that the length of the loop does not exceed that of the chain. We then apply the definitions of $\tau_2$ and $\tau_3$ to bound the loop, and use the bounds \eqref{step_2_cC_rd_bound_supercritical} and \eqref{eq:KL_KC} to control the chain. Combining these estimates with the convolution inequality
	\begin{equation}\label{example_triangular_loop_recovery}
		\begin{aligned}
			&~\sum_{x\in\ZL}\mathscr{D}_{u}(|a-x|)\mathscr{D}_{u}(|b-x|)\mathscr{D}_{u}(|c-x|)\mathscr{D}_{u}(|d-x|)\\
			\lesssim &~\mathscr{D}_{u}(|a-b|)\mathscr{D}_{u}(|c-d|)\sum_{x\in\ZL}\mathscr{D}_{u}(|a-x|\wedge |b-x|)\mathscr{D}_{u}(|c-x|\wedge |d-x|)\\
			\prec &~\ell_u\cdot\mathscr{D}_{u}(|a-b|)\mathscr{D}_{u}(|c-d|),
		\end{aligned}
	\end{equation}
	which recovers the loop structure of $\mathscr{D}_{u,\bx}^{(n)}$, we obtain
	\[\mathcal{E}^{(n)}_{u, \boldsymbol{\sigma},\bx} \prec \ell_u B_u^{n+1-\frac{n+1}{n+2}}(0)\mathscr{D}_{u,\bx}^{(n)}=(1-u)^{-1}\cdot B_u^{n-\frac{n+1}{n+2}}(0)\mathscr{D}_{u,\bx}^{(n)}.\]
	Together with \eqref{evolution_kernel_bound_1} and \eqref{bootstrap_assumption_supercritical}, this implies
	\begin{equation}\label{step_2_error_double_difference_n_34}
		\begin{aligned}
			\int_{s}^{r\wedge \tau}\pa{\bU_{u,r\wedge \tau,\bsigma}^{\pa{n}}\circ\mathcal{E}^{(n)}_{u, \boldsymbol{\sigma}}}_{\bx}\dd u\prec \frac{1-s}{1-r\wedge \tau} {B_{r\wedge \tau}^{n-\frac{n+1}{n+2}}(0)\mathscr{D}_{r\wedge \tau,\bx}^{(n)}}\le B_{r\wedge \tau}^{n-\frac{n+2}{n+3}+\frac1{1000}}(0)\mathscr{D}_{r\wedge \tau,\bx}^{(n)}.
		\end{aligned}
	\end{equation}
	A nearly identical argument, using \eqref{eq:KL_KC}, shows that for $n\in\{3,4\}$,
	\begin{equation}\label{step_2_error_single_difference_n_34}
		\sum_{l_{\cK}=3}^{n}\int_{s}^{r\wedge \tau}\pa{\bU_{u,r\wedge \tau,\bsigma}^{\pa{n}}\circ\qa{\mathfrak{D}_{l_{\cK}}(\cL-\cK)}_{u,\bsigma}^{\pa{n}}}_{\bx}\, \rd u\prec B_{r\wedge \tau}^{n-\frac{n+2}{n+3}+\frac1{1000}}(0)\mathscr{D}_{r\wedge \tau,\bx}^{(n)}.
	\end{equation}
	For the light-weight term, we apply the bound \eqref{step_2_cC_rd_bound_supercritical} when $n\in\{2,3\}$ and the bound \eqref{step_2_cC_rd_bound_5_supercritical} when $n=4$ to control the $(n+1)$-chain in \eqref{eq:LW_long}. We then apply the averaged local law in \eqref{step_1_supercritical_weak_local_law} to bound the light-weight factor. Combining these bounds with the convolution inequality \eqref{example_triangular_2_edge_to_1}
	to recover the factor \smash{$\mathscr{D}_{u,\bx}^{(n)}$}, we obtain 
	\[\mathcal{W}^{(n)}_{u, \boldsymbol{\sigma},\bx} \prec \ell_u B_u^{n+1/4}(0)\mathscr{D}_{u,\bx}^{(n)}=(1-u)^{-1}\cdot B_u^{n-3/4}(0)\mathscr{D}_{u,\bx}^{(n)}.\]
	Together with \eqref{evolution_kernel_bound_1} and \eqref{bootstrap_assumption_supercritical}, this implies
	\begin{equation}\label{step_2_error_LW_n_234}
		\int_{s}^{r\wedge \tau}\pa{\bU_{u,r\wedge \tau,\bsigma}^{\pa{n}}\circ\mathcal{W}^{(n)}_{u, \boldsymbol{\sigma}}}_{\bx}\dd u\prec \frac{1-s}{1-r\wedge \tau} B_{r\wedge \tau}^{n-3/4}(0)\mathscr{D}_{r\wedge \tau,\bx}^{\pa{n}}\le  B_{r\wedge \tau}^{n-\frac{n+2}{n+3}+\frac1{1000}}(0)\mathscr{D}_{r\wedge \tau,\bx}^{\pa{n}}.
	\end{equation}
	
	It remains to bound the martingale term in \eqref{eq:evolution_n_rtau}. To this end, we employ the following analogue of \Cref{lemma_martingale_to_cB}, which is again a consequence of the Burkholder-Davis-Gundy inequality. To state it, we define the tensor product of a sequence of matrices $A_1,\ldots,A_n \in \C^{\ZL\times\ZL}$ acting on an $n$-tensor $\cX^{\pa{n}}$ by
	\begin{equation}
		\begin{aligned}
			\pa{\bigotimes_{i=1}^{n}A_{i}\circ \cX^{\pa{n}}}_{\bx}=\sum_{\ba=(a_1,\ldots,a_n)\in\ZL^n}\prod_{i=1}^{n}(A_{i})_{x_ia_i}\cdot\cX_{\ba}^{\pa{n}},\qquad \forall \bx=(x_1,\ldots,x_n)\in\ZL^n.
		\end{aligned}
	\end{equation}
	\begin{lemma}\label{lemma_bound_martingale_to_cB_supercritical}
	In the setting of \Cref{local_law_with_n_loops_along_flow_supercritical}, fix any constant $p\in2\N$. Given deterministic coefficient matrices \smash{$A_u^{\pa{1}},\ldots,A_u^{\pa{n}}\in\C^{\ZL\times\ZL}$}, $\bx=(x_1,\ldots,x_n)\in\ZL^n$, $s\in[0,1)$, and a stopping time $\tau\in [s,\tf]$, we have
		\begin{equation*}        \begin{aligned}
				\E\absa{\int_{s}^{\tau}\pa{\bigotimes_{i=1}^{n}A_u^{\pa{i}}\circ\rd \cB_{u,\bsigma}^{\pa{n}}}_{\bx}}^p\leq C_{p,n} \E\pa{\int_s^{\tau}\ha{\qa{\bigotimes_{i=1}^{n}A_u^{\pa{i}}}\otimes\qa{\bigotimes_{i=1}^{n} \ol{A}_u^{\pa{i}}}\circ\pa{\cB\otimes\cB}_{u,\bsigma}^{(2n)}}_{(\bx,\bx)}\,\rd u}^{p/2}.
			\end{aligned}
		\end{equation*} 
		Here, $C_{p,n}>0$ is a constant depending only on $p,n$, and $\pa{\cB\otimes\cB}_{u,\bsigma,\ba}^{(2n)}$, with $\ba=(a_1,\ldots,a_n,a_1',\ldots,a_n')$ and $\bsigma=(\sigma_1,\ldots,\sigma_n)$, is defined by
		\begin{equation}\label{def_tensor_cB_supecritical}
			\begin{aligned}
				\pa{\cB\otimes\cB}_{u,\bsigma,\ba}^{(2n)}:=\sum_{k=1}^{n}\sum_{x}\cC_{u,(\bsigma\otimes\ol{\bsigma})^{(k)},\ba^{(k)}(x)}^{\pa{2n+2}},
			\end{aligned}
		\end{equation}
		where we adopt the abbreviation \eqref{abbreviation_diagonal_cC_chain} and denote 
		\begin{equation*}
			\begin{aligned}
				(\bsigma\otimes\ol{\bsigma})^{(k)}:=&(\sigma_k,\ldots,\sigma_n,\sigma_1,\ldots,\sigma_k,-\sigma_k,\ldots,-\sigma_1,-\sigma_n,\ldots,-\sigma_k),\\
				\ba^{(k)}(x):=&(a_k,\ldots,a_n,a_1,\ldots,a_{k-1},x,a_{k-1}',\ldots,a_1',a_n',\ldots,a_k',x).
			\end{aligned}
		\end{equation*}    
	\end{lemma}
	
	With this lemma, applying an argument similar to that used for \eqref{eq:UUUUBB}, based on the CS inequality, we obtain that for any fixed $p\in2\N$,
	\begin{equation}\label{step_2_cU_martingale_to_cU_X}
		\begin{aligned}
            \E\absa{\int_{s}^{r\wedge \tau}\pa{\bU_{u,r,\bsigma}^{\pa{n}}\circ
					\dd\mathcal{B}^{(n)}_{u, \boldsymbol{\sigma}}}_{\bx}}^p\lesssim \E\pa{\int_{s}^{r\wedge \tau} \pa{\cU_{u,r}^{\otimes n}\circ X_{u,\bsigma}^{\cB,\pa{n}}}_{\bx}^2\,\rd u }^{p/2},
		\end{aligned}
	\end{equation}
	where $\cU_{u,r}\equiv \cU^{(-,+)}_{u,r}$ and we define the tensor $$X_{u,\bsigma,\bx}^{\cB,\pa{n}}:=\qa{(\cB\otimes\cB)_{u,\bsigma,(x_1,\ldots,x_n,x_1,\ldots,x_n)}^{(2n)}}^{1/2},\quad \forall\bx=(x_1,\ldots,x_n)\in\ZL^n.$$ 
	To bound the expression on the RHS of \eqref{step_2_cU_martingale_to_cU_X} for $n\in\ha{2,3}$, we use the estimate 
	\begin{align}
		\pa{X_{u,\bsigma,\bx}^{\cB,\pa{n}}}^2&=(\cB\otimes \cB)_{u,\bsigma,(x_1,\ldots,x_n,x_1,\ldots, x_n)}^{\pa{2n}}\prec B_u^{2n+1}(0) \Xi_{u,2n+2}^{\cC_{\rd}} \sum_{k=1}^n\sum_{x}\qa{\mathscr{D}_{u,\bx_k}^{\pa{n+1}}(x)}^2\nonumber\\
		& \prec\pa{1-u}^{-1}\cdot B_u^{2n}(0)\Xi_{u,2n+2}^{\cC_{\rd}}\pb{\mathscr{D}_{u,\bx}^{\pa{n}}}^2 \prec\pa{1-u}^{-1}\cdot B_u^{2n-1}(0)\pb{\mathscr{D}_{u,\bx}^{\pa{n}}}^2,\label{eq:martingale_2n+2}
	\end{align}
	where $\bx_k:=(x_k,\ldots,x_n, x_1,\ldots,x_{k-1})$ and we adopt the convention $x_0:=x_n$. In the second step, we bound the summation over $x$ using
	\begin{equation}\label{eq:convolution-D^2D^2}
		\sum_{x}\mathscr{D}_{u}^2(|x_{k-1}-x|)\mathscr{D}_{u}^2(|x-x_k|)\lesssim \ell_u\mathscr{D}_{u}^2(|x_{k-1}-x_k|),
	\end{equation}
	and in the last step, we used $B_u(0)\Xi_{u,6}^{\cC_{\rd}}\prec \Xi_{u,2}^{\cC_{\rd}}\Xi_{u,4}^{\cL}\prec 1$ and $B_u(0)\Xi_{u,8}^{\cC_{\rd}}\prec \Xi_{u,4}^{\cC_{\rd}}\Xi_{u,4}^{\cL}\prec 1$, which follow from \eqref{loop_bound_cut_long_to_short}, \eqref{step_2_cL_bound_supercritical}, and \eqref{step_2_cC_rd_bound_supercritical}. Substituting this bound into \eqref{step_2_cU_martingale_to_cU_X}, using the evolution kernel estimate \eqref{evolution_kernel_bound_1}, integrating over $u$, and applying Markov's inequality, we obtain for $n\in\{2,3\}$,
	\begin{equation}\label{step_2_error_martingale_n_23}
		\begin{aligned} 
            \int_{s}^{r\wedge \tau}\pa{\bU_{u,r\wedge \tau,\bsigma}^{\pa{n}}\circ
		    \dd\mathcal{B}^{(n)}_{u, \boldsymbol{\sigma}}}_{\bx}\prec \frac{1-s}{1-r\wedge \tau} B_{r\wedge \tau}^{n-1/2}(0)\mathscr{D}_{r\wedge \tau,\bx}^{\pa{n}}\le B_{r\wedge \tau}^{n-3/4}(0)\mathscr{D}_{r\wedge \tau,\bx}^{\pa{n}}.
		\end{aligned}
	\end{equation}
However, the above argument does not extend to the case $n=4$, since the $10$-chains appearing in \smash{$(\cB\otimes\cB)_{u,\bsigma,\ba}^{(8)}$} are too long to be controlled directly using \eqref{loop_bound_cut_long_to_short}. To overcome this difficulty, we extend the $\contract$ inequality developed in \cite{dubova2025delocalizationnonmeanfieldrandommatrices} (see also \Cref{sec_estimates_error_terms2}) to bound such long chains, as stated in the following lemma.
	
	\begin{lemma}\label{lemma_n_loop_contraction}
		Fix any integer $n\ge 2$.  For any $\bsigma\in\ha{-,+}^n$ and $\bx=(x_1,\ldots,x_n)\in\ZL^n$, we have 
		\begin{equation}\label{eq_n_loop_contraction}
			\begin{aligned}
				(\cB\otimes\cB)_{u,\bsigma,(x_1,\ldots,x_n,x_1,\ldots, x_n)}^{\pa{2n}}\prec \frac{1}{1-u}\qa{\pa{\Xi_{u,4}^{\cL}}^{1/2}+\Xi_{u,2}^{\cC_{\txt{od}}}+\Xi_{u,2}^{\cC_{\rd}}}\Xi_{u,n-1}^{\cC_{\rd}}\Xi_{u,n}^{\cC_{\rd}}\cdot B_u^{2n-3/2}(0)\pb{\mathscr{D}_{u,\bx}^{\pa{n}}}^2.
			\end{aligned}
		\end{equation}
	\end{lemma}
Now, taking $n=4$ in \eqref{eq_n_loop_contraction}, we obtain
	\begin{equation*}
		\begin{aligned} 
        \pa{X_{u,\bsigma,\bx}^{\cB,\pa{n}}}^2=(\cB\otimes \cB)_{u,\bsigma,(x_1,\ldots,x_n,x_1,\ldots,x_n)}^{\pa{2n}}\prec \frac{1}{1-u}\cdot B_u^{2n-3/2}(0)\pb{\mathscr{D}_{u,\bx}^{\pa{n}}}^2.
		\end{aligned}
	\end{equation*}
	Substituting this estimate into \eqref{step_2_cU_martingale_to_cU_X}, using the evolution kernel bound \eqref{evolution_kernel_bound_1}, integrating over $u$, and applying Markov's inequality, we obtain
	\begin{equation}\label{step_2_error_martingale_n_4}
		\begin{aligned}
            \int_{s}^{r\wedge \tau}\pa{\bU_{u,r\wedge \tau,\bsigma}^{\pa{n}}\circ
				\dd\mathcal{B}^{(n)}_{u, \boldsymbol{\sigma}}}_{\bx}\prec \frac{1-s}{1-r\wedge \tau} B_{r\wedge \tau}^{n-3/4}(0)\mathscr{D}_{r\wedge \tau,\bx}^{\pa{n}}\le B_{r\wedge \tau}^{n-\frac{n+2}{n+3}+\frac1{1000}}(0) \mathscr{D}_{r\wedge \tau,\bx}^{\pa{n}}.
		\end{aligned}
	\end{equation}
Combining the bounds \eqref{step_2_error_leading_term_n_2}, \eqref{step_2_error_double_difference_n_2}, \eqref{step_2_error_double_difference_n_34}, \eqref{step_2_error_single_difference_n_34}, \eqref{step_2_error_LW_n_234}, \eqref{step_2_error_martingale_n_23}, and \eqref{step_2_error_martingale_n_4}, we obtain
	\begin{equation*}
		\begin{aligned}
            (\mathcal{L} - \mathcal{K})^{(n)}_{r\wedge \tau, \boldsymbol{\sigma}, \bx}\prec B_{r\wedge \tau}^{-\frac{n+2}{n+3}+\frac1{1000}}(0)\cdot B_{r\wedge \tau}^{n}(0)\mathscr{D}_{r\wedge \tau,\bx}^{\pa{n}},\qquad \forall n\in\ha{2,3,4}.
		\end{aligned}
	\end{equation*}
	By the definition \eqref{eq;XiL-Kn}, this implies
	\begin{equation*}
		\begin{aligned}
            B_{r\wedge \tau}^{\frac{n+2}{n+3}}(0)\Xi_{r\wedge \tau, n}^{\cL-\cK}\prec B_{r\wedge\tau}^{\frac1{1000}}(0),\qquad \forall n\in\ha{2,3,4}, \ \ r\in\qa{s,t}.
		\end{aligned}
	\end{equation*}
	Using this bound together with the definitions \eqref{eq:tau234} and \eqref{def_tau_cL}, and the relation \eqref{eq;Xi_relate1}, a standard continuity argument yields that $\tau=t$ with high probability. In particular, for any $\bsigma\in\{-,+\}^n$ and $\bx\in\ZL^n$ we obtain
	\begin{equation}
		\begin{aligned}
			(\cL-\cK)_{u,\bsigma,\bx}^{\pa{n}}\prec B_u^{n-(n+2)/(n+3)}(0)\cdot \mathscr{D}_{u,\bx}^{\pa{n}},\qquad \forall n\in\ha{2,3,4}, \ \ u\in\qa{s,t}.
		\end{aligned}
	\end{equation}
	This proves \eqref{step_2_loops_local_law_supercritical} and completes Step 2 in the proof of Lemma \ref{lemma_bootstrap_supercritical}. It remains to prove Lemma \ref{lemma_n_loop_contraction}.

	\begin{proof}[\bf Proof of Lemma \ref{lemma_n_loop_contraction}]
By the definition \eqref{def_tensor_cB_supecritical}, it suffices to prove that for any $k\in\qq{n}$, $\bsigma=(\sigma_1,\ldots,\sigma_n)\in\ha{-,+}^n$, and $\mathsf{x}=(x_1,\ldots,x_n,x_1,\ldots,x_n)\in\ZL^{2n}$,
		\begin{equation}\label{eq_n_loop_contraction_k1}
			\sum_{x}\cC_{u,(\bsigma\otimes \ol{\bsigma})^{(k)},\sx^{(k)}(x)}^{(2n+2)} \prec \frac{1}{1-u}\qa{\pa{\Xi_{u,4}^{\cL}}^{1/2}+\Xi_{u,2}^{\cC_{\txt{od}}}+\Xi_{u,2}^{\cC_{\rd}}}\Xi_{u,n-1}^{\cC_{\rd}}\Xi_{u,n}^{\cC_{\rd}}\cdot B_u^{2n-3/2}(0)\pb{\mathscr{D}_{u,\bx}^{\pa{n}}}^2.
		\end{equation}
		The proof is similar to that of Lemma \ref{lemma_martingale_term} in \Cref{sec_estimates_error_terms2}. For notational convenience, we assume $k=1$; the other cases follow by symmetry. We split the LHS into two parts:
		\begin{equation*}
			\begin{aligned}
				I_{\leq}+I_{>}:=\sum_{|x-x_n|\leq |x_1-x_n|/2}\cC_{u,(\bsigma\otimes \ol{\bsigma})^{(1)},\sx^{(1)}(x)}^{(2n+2)}+\sum_{|x-x_n|>|x_1-x_n|/2}\cC_{u,(\bsigma\otimes \ol{\bsigma})^{(1)},\sx^{(1)}(x)}^{(2n+2)}.
			\end{aligned}
		\end{equation*}
		
		We first consider $I_{>}$:
		\begin{equation}\label{bound_I_g_supercritical}
			I_{>}=\sum_{|x-x_n|>|x_1-x_n|/2}\psi_{x,>}A_{x,>}\psi_{x,>}^*\leq \sum_{|x-x_n|>|x_1-x_n|/2}\|A_{x,>}\|\cdot\|\psi_{x,>}\|^2,
		\end{equation}
		where the (row) vector $\psi_{x,>}\in \C^N$ and the matrix $A_{x,>}\in \C^{N\times N}$ are defined by 
		\begin{equation*}
			\begin{aligned}
				\psi_{x,> }(i):=\cC_{u,\bsigma,(x_1,\ldots,x_{n-1})}^{\pa{n}}(x,i)\sqrt{S_{i x_n}},\quad A_{x,>}(i,j):=\sqrt{S_{ix_n}}\cC_{u,(\sigma_1,-\sigma_1),x}^{(2)}(i,j)\sqrt{S_{j x_n}}.
			\end{aligned}
		\end{equation*}
		Bounding the operator norm by the HS norm and using the constraint $|x-x_n|>|x_1-x_n|/2$, we obtain 
		\begin{equation*}
			\begin{aligned}
				\|A_{x,>}\|^2\leq \|A_{x,>}\|_{\HS}^2=\cL_{u,(\sigma_1,-\sigma_1,\sigma_1,-\sigma_1),(x,x_n,x,x_n)}^{(4)}\prec \Xi_{u,4}^{\cL} \cdot B_u^{3}(0)\mathscr{D}_{u}^4(|x_1-x_n|).
			\end{aligned}
		\end{equation*}
		Substituting this bound into \eqref{bound_I_g_supercritical} yields 
		\begin{align}
			I_{>}&\prec \pa{\Xi_{u,4}^{\cL}}^{1/2}\cdot B_u^{3/2}(0)\mathscr{D}_{u}^2(|x_1-x_n|) \sum_{x}\|\psi_{x,>}\|^2 \nonumber\\
			&= \pa{\Xi_{u,4}^{\cL}}^{1/2}\cdot  B_u^{3/2}(0)\mathscr{D}_{u}^2(|x_1-x_n|)\sum_{x}\cC_{u,(\sigma_1,\ldots,\sigma_n,-\sigma_n,\ldots,-\sigma_1),(x_1,\ldots,x_{n-1},x_n,x_{n-1},\ldots,x_1)}^{(2n)}(x).\label{eq:I>>}
		\end{align}
		We now apply a Cauchy-Schwarz argument similar to \eqref{eq:CS-to-bound-chain} to the $(2n)$-chain in \eqref{eq:I>>}. Depending on whether $n$ is even or odd, this chain can be bounded either by the product of an $n$-loop and an $n$-chain of the form $\smash{\cC_{u,\bsigma',\bx'}^{(n)}(x)}$, or by the product of an $(n-1)$-loop and an $(n+1)$-chain of the form $\smash{\cC_{u,\bsigma'',\bx''}^{(n+1)}(x)}$. In both cases, the $\cL$-loop and the $\cC$-chain are non-negative. We then estimate the sum over $x$ of the $l$-$\cC$-chain (with $l=n$ or $l=n+1$) using Ward's identity \eqref{eq_Ward0} together with the parameter \smash{$\Xi_{u,l-1}^{\cC_{\rd}}$}, while the $\cL$-loop is bounded by \smash{$\Xi_{u,2n-l}^{\cL}$} and the relation \eqref{eq;Xi_relate2}. This yields
		\begin{equation}\label{eq:Psi>l2}
			\sum_{x\in\ZL}\|\psi_{x,>}\|^2\prec \frac{1}{1-u}\Xi_{u,l-1}^{\cC_{\rd}}\Xi_{u,2n-l}^{\cC_{\rd}}\cdot B_u^{2n-3}(0)\qa{\mathscr{D}_{u,\bx}^{\pa{n}}/\mathscr{D}_{u}(|x_1-x_n|)}^2,
		\end{equation}
		where the factor $(1-u)^{-1}$ arises from Ward's identity, $l=n$ if $n$ is even, and $l=n+1$ if $n$ is odd. Substituting \eqref{eq:Psi>l2} into \eqref{eq:I>>}, we obtain
		\begin{equation}\label{eq:I>>2}
			I_{>}\prec \frac{1}{1-u}\pa{\Xi_{u,4}^{\cL}}^{1/2}\Xi_{u,n-1}^{\cC_{\rd}}\Xi_{u,n}^{\cC_{\rd}}\cdot B_u^{2n-3/2}(0)\pb{\mathscr{D}_{u,\bx}^{\pa{n}}}^2.
		\end{equation}
		
		The estimate for $I_{\le}$ is analogous. We write
		\begin{equation*}
			\begin{aligned}        I_{\leq}=\sum_{x'\in\ZL}\psi_{x',\le}^* A_{x',\le}\psi_{x',\le},
			\end{aligned}
		\end{equation*}
		where the (column) vector $\psi_{x',\le}\in \C^N$ and the matrix $A_{x',\le}\in \C^{N\times N}$ are defined by 
		\begin{equation*}
			\begin{aligned}
				&\psi_{x',\le}(i):=\sqrt{S_{x_1i}}\cC_{u,(\sigma_2,\ldots,\sigma_n,\sigma_1),(x_2,\ldots, x_n)}^{(n)}(i,x'),\\
				&A_{x',\le}(i,j):=\sum_{|x-x_n|\leq |x_1-x_n|/2}\sqrt{S_{x_1i}}G_u(-\sigma_1)_{ix}S_{x'x}G_u(\sigma_1)_{xj}\sqrt{S_{jx_1}}.
			\end{aligned}
		\end{equation*}
		The remainder of the argument proceeds as for $I_{>}$. The only difference is that instead of the bound $ B_u^{3/2}(0)\p{\Xi_{u,4}^{\cL}}^{1/2}\mathscr{D}_{u}^2(|x_1-x_n|) $, we estimate the HS norm of $A_{x',\le}$ by 
		\[\qa{\pb{B_u^{3/2}(0)\Xi_{u,2}^{\cC_{\txt{od}}}}^2+W^{-1}\pb{B_u(0)\Xi_{u,2}^{\cC_{\rd}}}^2 }^{1/2}\mathscr{D}_{u}^2(|x_1-x_n|) \lesssim B_u^{3/2}(0)\pb{\Xi_{u,2}^{\cC_{\txt{od}}}+\Xi_{u,2}^{\cC_{\rd}}}\mathscr{D}_{u}^2(|x_1-x_n|). \] 
		We omit the remaining details for brevity. Together with \eqref{eq:I>>2}, this proves \eqref{eq_n_loop_contraction}. \end{proof}

	\subsection{Step 3: Sharp upper bound for $\cL$-loops}\label{sec_step_3_supercritical}

    For the proof of Step 3, we state several \emph{dynamical relations} between the \(\Xi\)-parameters, obtained from the analysis of the evolution equation \eqref{evolution_n_cL_cK_integrated}. These relations yield self-improving estimates for the \(\Xi\)-parameters and provide the main technical input for the bootstrap arguments in Steps 3 and 4 of the proof of Lemma \ref{lemma_bootstrap_supercritical}. Dynamical inequalities of the type stated in the following lemma have been established for regular RBMs in \cite{Band1D,Band2D,truong2025localizationlengthfinitevolumerandom,dubova2025delocalizationnonmeanfieldrandommatrices}, and our proof follows essentially the same strategy. Therefore, we defer the proof of Lemma \ref{lemma_dynamical_inequality} to \Cref{subsec:pf_lemma_dynamical_inequality}.

	\begin{lemma}\label{lemma_dynamical_inequality}
		In the setting of Lemma \ref{lemma_bootstrap_supercritical}, assume that the bounds \eqref{step_1_supercritical_weak_local_law}--\eqref{step_2_loops_bound_supercritical} hold. Suppose that \smash{$\h{\wh \Xi_{u,n}^{\cL-\cK}}_{n}$} and \smash{$\h{\wh \Xi_{u,n}^{\cC_{\rd}}}_{n}$} are two deterministic families of parameters satisfying
		\begin{equation}\label{eq:Xi_by_whXi}
			\sup_{u\in\qa{s,t}}\Xi_{u,n}^{\cL-\cK}/\wh \Xi_{u,n}^{\cL-\cK}+\sup_{u\in\qa{s,t}}\Xi_{u,n}^{\cC_{\rd}}/\wh \Xi_{u,n}^{\cC_{\rd}}\prec 1
		\end{equation}
		for any fixed $n\in \N$. Then, for any fixed integer $n\geq 3$ and all $r\in\qa{s,t}$, the dynamical inequality
		\begin{equation}\label{eq_dynamical_inequality}
			\Xi_{r,n}^{\cL-\cK}\prec \sup_{u\in\qa{s,r}}\ha{\max_{k=2}^{n-1}\wh \Xi_{u,k}^{\cL-\cK}+\max_{k=2}^{n-1}\wh \Xi_{u,k}^{\cL-\cK}\wh \Xi_{u,n+2-k}^{\cC_{\rd}}+\wh \Xi_{u,n+1}^{\cC_{\rd}}}+\Lambda_{r,n}^{1/2}
		\end{equation}
		holds, where $\Lambda_{r,n}$ denotes one of the following two families of deterministic control parameters:
		\begin{itemize}
			\item ({\bf Bounding martingale terms with long chains})
			\begin{equation}
				\Lambda_{r,n}=\sup_{u\in\qa{s,r}}\wh \Xi_{u,2n+2}^{\cC_{\rd}};\label{eq:choose_Lambda1}
			\end{equation}
			\item ({\bf Bounding martingale terms with $\contract$ inequalities})
			\begin{equation}\Lambda_{r,n}=\frac{(1-s)^2}{(1-r)^2} B_r^{-3/2}(0)\sup_{u\in\qa{s,r}}\pa{\wh \Xi_{u,n-1}^{\cC_{\rd}}\wh \Xi_{u,n}^{\cC_{\rd}}}.\label{eq:choose_Lambda2}
			\end{equation}
		\end{itemize}
       Moreover, for $n=2$, we have uniformly in $r\in[s,t]$,
\begin{equation}\label{eq_dynamical_inequality_n_2}
            \Xi_{r,2}^{\cL-\cK}\prec 1+ B_r^{1/2}(0) \sup_{u\in[s,r]}\Xi_{u,2}^{\cL-\cK}\Xi_{u,2}^{\cC-\cK^{\cC}}+\sup_{u\in\qa{s,r}}\pa{\Xi_{u,2n+2}^{\cC_{\rd}}}^{1/2}.
        \end{equation}
        
\end{lemma}

    We now apply the dynamical inequalities obtained in Lemma \ref{lemma_dynamical_inequality} to establish the sharp upper bound \eqref{step_3_n_bound_cL_supercritical} for $\cL$-loops. Specifically, we prove by induction on $n\geq 1$ that
    \begin{equation}\label{step_2_induction_target_supercritical}
            \sup_{u\in[s,t]}\Xi_{u,l}^{\cL-\cK}\prec B_s^{-4/5}(0),\qquad \forall 1\leq l\leq 2n.
    \end{equation}
    For $n=1$, this follows directly from \eqref{step_2_sharp_local_law_supercritical} and \eqref{step_2_loops_local_law_supercritical}. For $n=2$, it remains to bound the $\Xi$-parameters $\Xi_{u,3}^{\cL-\cK}$ and $\Xi_{u,4}^{\cL-\cK}$. By \eqref{step_2_loops_local_law_supercritical}, \eqref{eq_cC_d_to_loops}, and the third estimate in \eqref{loop_bound_cut_long_to_short} with $k=1$, we have
    \begin{equation}\label{initial_pair_input_bound}
            \max_{1\leq l\leq 4}\max_{\star\in\ha{\cL,\cC_{\rd}}}\sup_{u\in[s,t]}\Xi_{u,l}^{\star}\prec 1,\qquad \sup_{u\in[s,t]}\Xi_{u,2}^{\cL-\cK}\prec B_s^{-4/5}(0).
    \end{equation}
    Together with the fourth inequality in \eqref{loop_bound_cut_long_to_short}, this implies
    \begin{equation}\label{initial_pair_input_bound_Xi5Cd}
        \sup_{\star\in\ha{\cL,\cC_{\rd}}}\sup_{u\in[s,t]}\Xi_{u,5}^{\star}\prec B_s^{-1/2}(0).
    \end{equation}
To bound $\Xi_{u,3}^{\cL-\cK}$, applying the dynamical inequality \eqref{eq_dynamical_inequality} with martingale term bounded as in \eqref{eq:choose_Lambda2}, we have uniformly in $r\in[s,t]$
    \begin{equation}\label{step_3_Xi_3_L_K_weak_estimate}
\Xi_{r,3}^{\cL-\cK}\prec  \sup_{u\in\qa{s,r}}\wh \Xi_{u,2}^{\cL-\cK}\wh \Xi_{u,3}^{\cC_{\rd}}+\sup_{u\in\qa{s,r}}\wh \Xi_{u,4}^{\cC_{\rd}}+\frac{1-s}{1-r}B_r^{-3/4}(0)\sup_{u\in[s,r]}\pa{\wh\Xi_{u,2}^{\cC_{\rd}}\wh\Xi_{u,3}^{\cC_{\rd}}}^{1/2}\prec B_s^{-4/5}(0).
\end{equation}
    Here, we also used the bounds in \eqref{initial_pair_input_bound}, the bootstrap assumption \eqref{bootstrap_assumption_supercritical}, and \eqref{bound_B_t_to_B_s}. 
    By a completely analogous argument, we obtain the following bound for
\(\Xi_{r,4}^{\cL-\cK}\), uniformly in \(r\in[s,t]\):
    \begin{equation*}
        \begin{aligned}
            \Xi_{r,4}^{\cL-\cK}\prec & \sup_{u\in\qa{s,r}}\wh \Xi_{u,2}^{\cL-\cK}\wh \Xi_{u,4}^{\cC_{\rd}}+\sup_{u\in\qa{s,r}}\wh \Xi_{u,3}^{\cL-\cK}\wh \Xi_{u,3}^{\cC_{\rd}}+\sup_{u\in\qa{s,r}}\wh \Xi_{u,5}^{\cC_{\rd}}+\frac{1-s}{1-r}B_r^{-3/4}(0)\sup_{u\in[s,r]}\pa{\wh\Xi_{u,3}^{\cC_{\rd}}\wh\Xi_{u,4}^{\cC_{\rd}}}^{1/2}\prec B_s^{-4/5}(0),
        \end{aligned}
    \end{equation*}
    where, in the second step, we used \eqref{initial_pair_input_bound}--\eqref{step_3_Xi_3_L_K_weak_estimate}.

    Next, suppose that $n\geq 3$ is an integer and that
    \begin{equation}\label{step_2_induction_assumption_0_supercritical}
            \sup_{u\in[s,t]}\Xi_{u,l}^{\cL-\cK}\prec B_s^{-4/5}(0),\qquad \forall 1\leq l\leq 2n-2.
    \end{equation}
    Together with \eqref{eq;Xi_relate1}, \eqref{eq_cC_d_to_loops}, and the third bound in \eqref{loop_bound_cut_long_to_short}, this, by an induction argument, gives
    \begin{equation}\label{step_2_induction_assumption_supercritical}
\sup_{\star\in\ha{\cL,\cC_{\rd}}}\sup_{u\in[s,t]}\Xi_{u,l}^{\star}\prec 1,\qquad \forall 1\leq l\leq 2n-2.
    \end{equation}
    Our goal is to prove that the bound in \eqref{step_2_induction_assumption_0_supercritical} holds for \(l=2n-1,2n\). To this end, we use a self-improving iteration. Let \(\Phi\equiv \Phi_N\geq 1\) be a deterministic control parameter, and assume that
    \begin{equation}\label{step_2_iteration_assumption_supercritical}
        \sup_{u\in[s,t]}\Xi_{u,2n}^{\cL-\cK}\prec \Phi.
    \end{equation}
    By \eqref{eq;Xi_relate1}, we have uniformly in $u\in[s,t]$
    \begin{equation}\label{step_2_iteration_input_supercritical}
        \Xi_{u,2n}^{\cL}\prec 1+B_u(0)\Phi,
    \end{equation}
    which, together with \eqref{eq_cC_d_to_loops} and \eqref{step_2_induction_assumption_supercritical}, implies
    \begin{equation}\label{step_2_iteration_cC_d_input_supercritical}
        \Xi_{u,2n}^{\cC_{\rd}}\prec 1+B_u(0)\Phi.
    \end{equation}
    Moreover, applying the third and fourth inequalities in \eqref{loop_bound_cut_long_to_short}, together with the bounds \eqref{step_2_induction_assumption_supercritical}, \eqref{step_2_iteration_input_supercritical}, and \eqref{step_2_iteration_cC_d_input_supercritical}, we obtain uniformly in $u\in[s,t]$
    \begin{equation}\label{step_2_iteration_input_odd_supercritical}
        \Xi_{u,2n-1}^{\star}\prec (1+B_u(0)\Phi)^{1/2},\qquad \Xi_{u,2n+1}^{\star}\prec B_u^{-1/2}(0)(1+B_u(0)\Phi),\qquad \forall \star\in\ha{\cL,\cC_{\rd}}.
    \end{equation}
    With these preparations, we first bound the parameter $\Xi_{u,2n-1}^{\cL-\cK}$. Applying the dynamical inequality \eqref{eq_dynamical_inequality}, we have uniformly in $r\in[s,t]$
    \begin{equation}\label{step_2_2n_1_bound_supercritical}
        \begin{aligned}
            &\Xi_{r,2n-1}^{\cL-\cK}\prec \sup_{u\in\qa{s,r}}\max_{k=2}^{2n-2}\wh \Xi_{u,k}^{\cL-\cK}\wh \Xi_{u,2n+1-k}^{\cC_{\rd}}+\sup_{u\in\qa{s,r}}\wh \Xi_{u,2n}^{\cC_{\rd}}+\frac{1-s}{1-r}B_r^{-3/4}(0)\sup_{u\in[s,r]}\pa{\wh\Xi_{u,2n-2}^{\cC_{\rd}}\wh\Xi_{u,2n-1}^{\cC_{\rd}}}^{1/2}\\
            & \quad\lesssim B_s^{-4/5}(0)(1+B_r(0)\Phi)^{1/2}+1+B_r(0)\Phi+B_s^{-4/5}(0)(1+B_r(0)\Phi)^{1/4}\lesssim B_s^{-4/5}(0)(1+B_r(0)\Phi),
        \end{aligned}
    \end{equation}
    where, in the second step, we applied the estimates \eqref{step_2_induction_assumption_0_supercritical}--\eqref{step_2_iteration_input_odd_supercritical} to bound the $\wh\Xi$-parameters. The parameter $\Xi_{u,2n}^{\cL-\cK}$ can be estimated analogously. The dynamical inequality \eqref{eq_dynamical_inequality} gives
    \begin{equation}
        \Xi_{r,2n}^{\cL-\cK}\prec \sup_{u\in\qa{s,r}}\max_{k=2}^{2n-1}\wh \Xi_{u,k}^{\cL-\cK}\wh \Xi_{u,2n+2-k}^{\cC_{\rd}}+\sup_{u\in\qa{s,r}}\wh \Xi_{u,2n+1}^{\cC_{\rd}}+\frac{1-s}{1-r}B_r^{-3/4}(0)\sup_{u\in[s,r]}\pa{\wh\Xi_{u,2n-1}^{\cC_{\rd}}\wh\Xi_{u,2n}^{\cC_{\rd}}}^{1/2}.
    \end{equation}
    Applying the estimates \eqref{step_2_induction_assumption_0_supercritical}--\eqref{step_2_iteration_input_odd_supercritical}, we can bound the second and third terms, as well as the first term with \(2\leq k\leq 2n-2\), by \smash{\(B_s^{-4/5}(0)(1+B_t(0)\Phi)\)}. For the first term with \(k=2n-1\), using \eqref{step_2_2n_1_bound_supercritical} together with \smash{\(\Xi_{u,3}^{\cC_{\rd}}\prec 1\)}, we obtain the same bound. 

    In sum, we have shown the following self-improving estimate: for any deterministic control parameter $\Phi\equiv \Phi_N\geq 1$, 
    \begin{equation}\label{step_2_iteration_implication_supercritical}
            \sup_{u\in[s,t]}\Xi_{u,2n}^{\cL-\cK}\prec \Phi\qquad\implies \qquad \sup_{u\in[s,t]}\Xi_{u,2n}^{\cL-\cK}\prec B_s^{-4/5}(0)+B_s^{1/6}(0)\Phi,
    \end{equation}
    where we used \eqref{bound_B_t_to_B_s} to estimate the coefficient \(B_s^{-4/5}(0)B_t(0)\) by \(B_s^{1/6}(0)\).
Starting from the trivial bound \smash{$\sup_{u\in[s,t]}\Xi_{u,2n}^{\cL-\cK}\leq W^C$} for some sufficiently large constant \(C>1\), and iterating the implication \eqref{step_2_iteration_implication_supercritical} $\OO(1)$ times, we conclude that
    \begin{equation*}
        \begin{aligned}
            \sup_{u\in[s,t]}\Xi_{u,2n}^{\cL-\cK}\prec B_s^{-4/5}(0).
        \end{aligned}
    \end{equation*}
    This completes the induction and the proof of \eqref{step_2_induction_target_supercritical} for all \(n\geq 1\).
Substituting the obtained bound into \eqref{eq;Xi_relate1} yields \eqref{step_3_n_bound_cL_supercritical} and completes Step 3 in the proof of \Cref{lemma_bootstrap_supercritical}. Finally, combining
\eqref{step_3_n_bound_cL_supercritical} with \Cref{lemma_Xi_reshape}, we obtain the following estimate, which will be used in subsequent arguments:
\begin{equation}\label{step_4_Xi_star_bound}
            \sup_{u\in[s,t]}\Xi_{u,n}^{\star}\prec 1,\qquad \forall\star\in\ha{\cL,\cC_{\rd},\cC_{\txt{od}}}.
    \end{equation}

	\subsection{Step 4: Sharp estimates for $\cL$-loops and $T$-variables}

    In addition to \eqref{step_3_n_bound_cL_supercritical}, we have already established the bound \eqref{step_4_Xi_star_bound}. We now prove the loop estimate \eqref{step_4_supercritical_sharp_local_law_loops} by induction on $n$ using the dynamical inequality \eqref{eq_dynamical_inequality}.
The case $n=1$ follows from the averaged local law in \eqref{step_2_sharp_local_law_supercritical}, together with the observation that \smash{$\mathscr D^{(1)}_{t,\bx}\equiv 1$} (since $x_2= x_1$ under the cyclic convention). 
For $n=2$, take the stopping time
\begin{equation*}
    \tau\equiv\tau(c):=t\wedge\inf\ha{u\geq s:\Xi_{u,2}^{\cL-\cK}\geq N^{c}},
\end{equation*}
where $c>0$ is a sufficiently small constant, and we adopt the convention $\inf\emptyset=\infty$. 
A standard $N^{-C}$-net argument allows us to take $r=\tau$ in \eqref{eq_dynamical_inequality_n_2}, which, together with the reshaping estimate \eqref{reshape_2_chains_to_loops} and the sharp upper bound \eqref{step_4_Xi_star_bound}, immediately yields a self-improving estimate for \smash{$\Xi_{u,2}^{\cL-\cK}$}, and further shows $\tau\geq t$ with high probability. This provides the sharp 2-loop estimate \smash{$\sup_{u\in[s,t]}\Xi_{u,2}^{\cL-\cK}\prec 1$}. Next, using \eqref{reshape_2_chains_to_loops}, we obtain, for any \(\bsigma=(\sigma_1,\sigma_2)\in\ha{-,+}^2\), \(x,y\in\ZL\), and \(u\in\qa{s,t}\),
	\begin{equation}\label{step_4_T_Theta_bound}
		\begin{aligned}
			(T-\Theta)_{u, xy}^{\bsigma}&= u \cdot m\p{\sigma_1}m\p{\sigma_2}(\cL-\cK)_{u, xy}^{\bsigma}+\opr{B_{u}^{3/2}(0)\mathscr{D}_{u}^2(|x-y|)} \prec B_{u}^{3/2}(0)\mathscr{D}_{u}^2(|x-y|).
		\end{aligned}
	\end{equation}
From \eqref{step_4_T_Theta_bound}, we deduce the estimate \eqref{step_4_local_law_T_supercritical} for the diagonal $T$-variables. The corresponding estimate for the off-diagonal $T$-variables follows directly from \eqref{step_4_Xi_star_bound} applied to $\Xi_{u,2}^{\cC_{\txt{od}}}.$
 
	Next, let $n\ge 3$ and suppose that \smash{$\sup_{u\in\qa{s,t}}\Xi_{u, k}^{\cL-\cK}\prec 1$} for all $k=1,2,\ldots, n-1$.  
	Applying \eqref{step_4_Xi_star_bound} to \eqref{eq_dynamical_inequality}, using the bound $\Lambda_{r,n}\prec 1$ (under the choice \eqref{eq:choose_Lambda1}), and applying an $N^{-C}$-net argument, we obtain 
	\begin{equation*}
		\sup_{u\in[s,t]}\Xi_{u,n}^{\cL-\cK}\prec \sup_{u\in\qa{s,t}} \max_{k=2}^{n-1}\wh \Xi_{u,k}^{\cL-\cK}  + 1 \prec 1.
	\end{equation*}
	This completes the induction and establishes \eqref{step_4_supercritical_sharp_local_law_loops} for all fixed $n$, thereby concluding Step 4.

 \subsection{Step 5: Sharp chain estimates}\label{sec_step_5_supercritical}
 The estimate \eqref{step_5_supercritical_sharp_local_law_chains} for off-diagonal chains is an immediate consequence of \eqref{step_4_Xi_star_bound} and the definition \eqref{eq:Xi_Cd_Cod}, while the corresponding estimate for diagonal  $1$- and $2$-chains (i.e., the resolvent entries and $T$-variables) has already been established in \eqref{step_2_sharp_local_law_supercritical} and  \eqref{step_4_local_law_T_supercritical}. A direct induction argument based on the reshaping bound \eqref{reshape_chains_into_loops} then yields
    \begin{equation}
        \Xi_{u,n}^{\cC-\cK^{\cC}}\prec 1,\qquad \forall u\in[s,t]\quad \txt{and}\quad n\geq 1.
    \end{equation}
    This completes Step 5 in the proof of \Cref{lemma_bootstrap_supercritical}.

	\subsection{Step 6: Expected estimate for 2-$\cL$-loops}

	The proof of this step combines ideas from the case $\alpha\in(0,1)$ (see \Cref{sec_proof_of_theorem_que_along_flow}) with the arguments used to establish expected 2-loop estimates in previous works \cite{Band1D,Band2D,truong2025localizationlengthfinitevolumerandom,fan2025blockreductionmethodrandom}. Accordingly, we present only a brief argument.
	We begin with the following expected single-resolvent local law, which can be proved analogously to Lemma \ref{lemma_expected_single_resolvent_average_local_law}; we therefore omit the details.
	
	\begin{lemma}
		In the setting of Lemma \ref{lemma_bootstrap_supercritical}, for any $u\in\qa{s,t}$, we have     \begin{equation}\label{expected_single_resolvent_local_law_supercritical}
			\max_{x\in\ZL}\absa{\E\avgb{\p{G_u-m}S^{(x)}}}\prec B_u^2(0).
		\end{equation}
	\end{lemma}
	
	Next, taking expectations on both sides of \eqref{evolution_n_cL_cK_integrated} with $n=2$, we obtain for any $\bsigma=(\sigma_1,\sigma_2)\in\ha{-,+}^2$, $\bx=(x,y)\in\ZL^2$, and $r\in\qa{s,t}$ that
	\begin{equation}
		\E(\mathcal{L} - \mathcal{K})^{(2)}_{r, \boldsymbol{\sigma}, \bx}=\pa{\bU_{s,r,\bsigma}^{\pa{2}}\circ \E\pa{\cL-\cK}_{s,\bsigma}^{\pa{2}}}_{\bx}+ \int_{s}^r\pa{\bU_{u,r,\bsigma}^{\pa{2}}\circ\E\mathcal{E}^{(2)}_{u, \boldsymbol{\sigma}}}_{\bx}\dd u +\int_{s}^{r}\pa{\bU_{u,r,\bsigma}^{\pa{2}}\circ\E \mathcal{W}^{(2)}_{u, \boldsymbol{\sigma}}}_{\bx}
		\dd u.\label{eq:EL-K2}
	\end{equation}
	We first consider the case $\sigma_1=\sigma_2$. By the induction hypothesis \eqref{step_6_expectedloop} at time $s$ and the evolution kernel bound \eqref{evolution_kernel_bound_2}, the first term is bounded by
	\begin{equation*}
		\begin{aligned}
			\pa{\bU_{s,r,\bsigma}^{\pa{2}}\circ \E\pa{\cL-\cK}_{s,\bsigma}^{\pa{2}}}_{\bx}\prec \frac{\ell_s(1-s)}{\ell_r(1-r)}B_s^{5/2}(0)\mathscr{D}_{r,\bx}^{(2)}\lesssim B_r^{5/2}(0)\mathscr{D}_{r,\bx}^{(2)}.
		\end{aligned}
	\end{equation*}
	For the second term, using the $\cL$-loop estimate \eqref{step_4_supercritical_sharp_local_law_loops} with $n=2$, the $T$-variable bound \eqref{step_4_local_law_T_supercritical}, the convolution inequality \eqref{eq:convolution-D^2D^2}, and the evolution kernel bound \eqref{evolution_kernel_bound_2}, we obtain
	\begin{equation*}
		\begin{aligned}
			\int_{s}^r\pa{\bU_{u,r,\bsigma}^{\pa{2}}\circ\E\mathcal{E}^{(2)}_{u, \boldsymbol{\sigma}}}_{\bx}\dd u\prec \int_{s}^r\frac{1}{1-u}\frac{\ell_u(1-u)}{\ell_r(1-r)} B_u^{5/2}(0)\mathscr{D}_{r,\bx}^{(2)}\dd u\prec B_r^{5/2}(0)\mathscr{D}_{r,\bx}^{(2)}.
		\end{aligned}
	\end{equation*}
	For the third term, we write 
	\begin{align}
		\E \mathcal{W}^{(2)}_{u, (\sigma_1,\sigma_2),(x,y)}&=\sum_{a\in\ZL}\E\avgb{\qb{G_u(\sigma_1)-m\pa{\sigma_1}}S^{(a)}}\pa{\cC-\cK^{\cC}}_{u,(\sigma_1,\sigma_2,\sigma_1),(x,y,a)}^{(3)}\label{eq;E-LW}\\
		&+\sum_{a\in\ZL}\E\avgb{\qb{G_u(\sigma_1)-m\pa{\sigma_1}}S^{(a)}}\cK_{u,(\sigma_1,\sigma_2,\sigma_1),(x,y,a)}^{\cC,(3)}+\qa{(\sigma_1,x)\leftrightarrow (\sigma_2,y)},\nonumber
	\end{align}
	where $\qa{(\sigma_1,x)\leftrightarrow (\sigma_2,y)}$ denotes the expression obtained by interchanging $(\sigma_1,x)$ and $(\sigma_2,y)$ in the preceding terms.
	Applying \eqref{step_5_supercritical_sharp_local_law_chains} together with the averaged local law \eqref{step_2_sharp_local_law_supercritical} to the first term, and \eqref{expected_single_resolvent_local_law_supercritical} together with \eqref{eq:KL_KC} to the second, we obtain
	\begin{equation*}
		\begin{aligned}
			\E \mathcal{W}^{(2)}_{u, (\sigma_1,\sigma_2),(x,y)}&\prec \sum_{a\in\ZL} B_u(0)\cdot B_u^{5/2}(0)\mathscr{D}^{(3)}_{u,(x,y,a)}+\sum_{a\in\ZL} B_u^2(0)\cdot B_u^{2}(0)\mathscr{D}^{(3)}_{u,(x,y,a)}+\qa{x\leftrightarrow y}\\
			&\prec (1-u)^{-1} B_u^{5/2}(0)\mathscr{D}_{u,\bx}^{(2)},
		\end{aligned}
	\end{equation*}
	where we also applied the convolution inequality \eqref{example_triangular_2_edge_to_1} in the second step. Using the evolution kernel bound \eqref{evolution_kernel_bound_2} and integrating over $u$, we obtain
	\begin{equation*}
		\begin{aligned}
			\int_{s}^{r}\pa{\bU_{u,r,\bsigma}^{\pa{2}}\circ\E \mathcal{W}^{(2)}_{u, \boldsymbol{\sigma}}}_{\bx}\dd u\prec B_r^{5/2}(0)\mathscr{D}_{r,\bx}^{(2)}.
		\end{aligned}
	\end{equation*}
	Substituting the above estimates into \eqref{eq:EL-K2}, we conclude that
	\begin{equation}\label{eq:EL-K_samesigma}
		\E(\mathcal{L} - \mathcal{K})^{(2)}_{r, \boldsymbol{\sigma}, \bx} \prec B_r^{5/2}(0)\mathscr{D}_{r,\bx}^{(2)}
	\end{equation}
	for all $r\in[s,t]$ in the case $\sigma_1=\sigma_2$.

For the case \(\sigma_1\neq \sigma_2\), the proof is similar, except that we use the evolution kernel bound \eqref{evolution_kernel_bound_3} instead. To this end, we decompose the relevant terms into a sum-zero component using a \emph{sum-zero operator} and a remaining component that can be estimated by Ward's identity. Similar arguments are also used in the derivation of the dynamical inequalities in Lemma \ref{lemma_dynamical_inequality}.
Since the argument is almost the same as that in \cite[Section 5.8]{Band1D}, we postpone the details of the proof of \eqref{eq:EL-K_samesigma} for \(\sigma_1\neq \sigma_2\) to the end of \Cref{subsec:pf_lemma_dynamical_inequality}, after formally introducing the sum-zero operator in Definition \ref{def_sum_zero_operator_supercritical}. This completes Step 6 and hence the proof of Lemma \ref{lemma_bootstrap_supercritical}.

\newpage
	
	\appendix

	\section{Proofs of some deterministic estimates}\label{appendix_additional_proof}

	\subsection{Proof of the bounds in \Cref{assumption_input_bound} for Example \ref{example_profile_function}}\label{appendix_proof_of_example_profile_function}

In this proof, to distinguish the graph distance from the usual absolute value, we use $|\cdot|_N$ to denote the periodic distance on $\ZL$. Since $f_\alpha$ is an even function, it follows immediately that $S_{xy}$ depends only on $|x-y|_N$. In particular, the matrix $S$ is translation invariant on $\ZL$. Therefore, for any $|\xi|<1$, we may write the Fourier expansion
	\begin{equation}\label{fourier_expansion_Theta}
			\pa{\frac{S}{1-\xi S}}_{xy}=\frac{1}{N}\sum_{\bp\in\bbT_N}\frac{\psi(\bp)}{1-\xi \psi(\bp)}\exp\pa{\ii \bp (x-y)}=:\frac{1}{N}\sum_{\bp\in\bbT_N}\Lambda_{\xi}(\bp)\exp\pa{\ii \bp (x-y)},
	\end{equation}
	where $\bbT_N:=(2\pi/N)\ZL$ and $\ha{\psi(p)}_{p\in\bbT_N}$ are the eigenvalues of $S$, given by
	\begin{equation*}
		\begin{aligned}
			\psi(\bp):=\sum_{x\in\ZL}S_{0x}\exp\pa{\ii \bp x}=\frac{1}{Z_{\alpha}}\sum_{x\in\ZL}\sum_{n\in\Z}f_{\alpha}\pa{\frac{x+nN}{W}}\exp\pa{\ii \bp x}=\frac{1}{Z_\alpha}\sum_{x\in\Z}f_{\alpha}\pa{\frac{x}{W}}\exp\pa{\ii \bp x}.
		\end{aligned}
	\end{equation*}
	Using the Poisson summation formula, we obtain for $p\in\bbT_N$ that \begin{equation}\label{eq:defpsip}
			\psi(\bp)=\frac{1}{Z_\alpha}\sum_{x\in\Z}f_{\alpha}\pa{\frac{x}{W}}\exp\pa{\ii \bp x}=\frac{W}{Z_\alpha}\sum_{n\in\Z}\phi_{\alpha}\pa{W\pa{2\pi n-p}},
	\end{equation}
	where the absolute convergence of the series follows from the decay assumptions \eqref{f_alpha_decay} and \eqref{phi_alpha_decay}. In particular, this representation extends $\psi$ and $\Lambda_{\xi}$ to $2\pi$-periodic functions on $\R$.

	Let $s(x)=|x|\wedge 1$. We summarize several basic properties of $\psi$, which follow from \eqref{f_alpha_decay}--\eqref{phi_alpha_decay} by straightforward estimates.
	\begin{itemize}
		\item The following lower bound holds uniformly for $p\in\qa{-\pi,\pi}$,
		\begin{equation}\label{spectral_gap_psi_p}
			\begin{aligned}
				1-\psi(p)\gtrsim [s(Wp)]^{\alpha\wedge 2}.
			\end{aligned}
		\end{equation}
		\item There exists a small constant $\e>0$ such that uniformly for $p\in\qa{-\pi,\pi}\setminus(-W^{-1},W^{-1})$,
		\begin{equation}\label{decay_psi_p}
			\begin{aligned}
				\absa{\psi^{(k)}(p)}\lesssim W^{k} (|Wp|+1)^{-2-\e},\quad \forall k\in\{0,1,\ldots,\ceil{\alpha}+2\}.
			\end{aligned}
		\end{equation}
		\item For any small constant $\e>0$, the derivatives of $\psi$ are uniformly bounded in $p\in\qa{-\pi,\pi}$ as 
		\begin{equation}\label{singularity_psi_p}
			\begin{aligned}
				& \absa{\psi'(p)}\lesssim \mathbf{1}_{\alpha>2}\cdot Ws(Wp)+\mathbf{1}_{1\leq \alpha\leq 2}\cdot W[s(Wp)]^{\alpha-1-\e}+\mathbf{1}_{0<\alpha<1}\cdot W [s(Wp)]^{\alpha-1},\\
				& \absa{\psi^{\pa{k}}(p)}\lesssim W^{k},\quad \forall 2\leq k\leq \ceil{\alpha}-1,\\
				& \absa{\psi^{\pa{k}}\p{p}}\lesssim  W^{k}[s(Wp)]^{\alpha-k}\pa{1+\mathbf{1}_{\alpha\geq 1}\cdot [s(Wp)]^{-\varepsilon}},\quad \forall  \ceil{\alpha}\le k \le  \ceil{\alpha}+2 .
			\end{aligned}
		\end{equation}
	\end{itemize}
	We briefly outline the proof of these bounds. First, taking $p=0$ in \eqref{eq:defpsip}, we obtain 
	\begin{equation*}
		Z_{\alpha}=W\sum_{n\in\Z}\phi_{\alpha}(2\pi nW)= W+W\cdot\OO\pa{\sum_{n=1}^{\infty}(nW)^{-2}}\asymp W.
	\end{equation*}
	Combining this estimate with the assumptions \eqref{psi_spectral_gap} and \eqref{phi_alpha_decay}, and using that $f_{\alpha}$ is even, we find
	\begin{equation*}
		\begin{aligned}
			&1-\psi(p)=\frac{W}{Z_\alpha}\sum_{n\in\Z}\ha{\phi_{\alpha}\pa{2\pi nW}-\frac{1}{2}\qa{\phi_{\alpha}\pa{W\pa{2\pi n-p}}+\phi_{\alpha}\pa{W\pa{2\pi n+p}}}}\\
			&\qquad =\frac{W}{Z_\alpha}\qa{1-\phi_{\alpha}(Wp)}+(|Wp|^{2}\wedge 1)\cdot\OO\pa{\sum_{n=1}^{\infty}(nW)^{-2}}\gtrsim [s(Wp)]^{\alpha\wedge 2}+\OO(W^{-2}[s(Wp)]^{2}) \gtrsim [s(Wp)]^{\alpha\wedge 2}.
		\end{aligned}
	\end{equation*}
	This concludes the proof of \eqref{spectral_gap_psi_p}. 
	For the proof of \eqref{decay_psi_p} and the last estimate in \eqref{singularity_psi_p}, we differentiate the Poisson representation of $\psi(p)$ in terms of $\phi_\alpha$ and apply the assumptions \eqref{phi_alpha_regularity} and \eqref{phi_alpha_decay}. 
	For the first bound in \eqref{singularity_psi_p} when $0<\al< 2$, the same argument applies, relying on the precise behavior of \smash{$\phi_\al'(t)$} near $t=0$ given by \eqref{phi_alpha_regularity}. 
	When $\al\ge 2$, it is more convenient to use the representation
	\[\psi^{(k)}(\bp)=\frac{1}{Z_\alpha}\sum_{x\in\Z}f_{\alpha}\pa{\frac{x}{W}}(\ii x)^{k}[\exp\pa{\ii \bp x}-\mathbf 1_{k=1}], \quad \forall 1\le k \le \ceil{\al}-1,\]
	where the subtraction of $\mathbf 1_{k=1}$ exploits the symmetry of $f_\alpha$. The first two bounds in \eqref{singularity_psi_p} then follow directly from the decay property \eqref{f_alpha_decay} together with elementary estimates on the oscillatory factor.

	We now rewrite \eqref{fourier_expansion_Theta} as 
	\begin{align}
		\pa{\frac{S}{1-\xi S}}_{xy}&=\frac{1}{N}\sum_{\bp\in\bbT_N}\Lambda_{\xi}(\bp)\exp\pa{\ii \bp (x-y)}=\frac{1}{N}\sum_{\bp\in\bbT_N}\sum_{k\in\Z}\wh\Lambda_{\xi}(k)\exp\pa{\ii k \bp}\exp\pa{\ii \bp (x-y)}\nonumber\\
		&=\sum_{n\in\Z}\wh \Lambda_{\xi}(y-x+nN),\label{Poisson_expansion_Theta}
	\end{align}
	where the Fourier coefficients $\wh\Lambda_\xi(k)$ are defined by
	\begin{equation}\nonumber
		\begin{aligned}
			\wh \Lambda_{\xi}(k):=\frac{1}{2\pi}\int_{-\pi}^{\pi}\Lambda_{\xi}(p)\exp\pa{-\ii k p}\, \rd p=\frac{1}{2\pi}\int_{-\pi}^{\pi}\frac{\psi(p)}{1-\xi\psi(p)}\exp\pa{-\ii k p}\, \rd p.
		\end{aligned}
	\end{equation}
	In the third step of \eqref{Poisson_expansion_Theta}, we used the commutativity of the summations, which is justified by the absolute summability of the series \smash{$\sum_{n\in\Z}\wh\Lambda_\xi(n)$} established below. We now prove the bounds \eqref{Theta_bound_opposite_charge}--\eqref{zero_mode_removed_bound} in \Cref{assumption_input_bound} for the following three regimes: (i) $\alpha\in[2,\infty)$, (ii) $\alpha\in[1,2)$, and (iii) $\alpha\in(0,1)$. The proof relies on a standard integration-by-parts argument, which has been widely used in the study of diffusion profiles for random band matrices; see, for example, \cite{yang2025delocalizationgeneralclassrandom, truong2025localizationlengthfinitevolumerandom, erdos2025zigzagstrategyrandomband}.
	
	\medskip
	\noindent
	{\bf Case (i): $\alpha\in[2,\infty)$.} 
The differentiability of $\psi$ at $0$, given by \eqref{singularity_psi_p}, together with periodicity, allows us to apply integration by parts $\ceil{\alpha}$ times. For $n\neq 0$, this yields
	\begin{equation}\label{eq:whLambdaxi}
		\wh \Lambda_{\xi}(n)=(\ii n)^{-\ceil{\alpha}}\frac{1}{2\pi}\int_{-\pi}^{\pi}\pa{\frac{\psi(p)}{1-\xi\psi(p)}}^{(\ceil{\alpha})}e^{-\ii n p}\, \rd p.
	\end{equation}
	Using \eqref{Poisson_expansion_Theta} and \eqref{eq:whLambdaxi}, we first prove the upper bounds \eqref{Theta_bound_opposite_charge} and \eqref{Theta_bound_equaled_charge}.
	
	For \eqref{Theta_bound_opposite_charge}, using \eqref{spectral_gap_psi_p} and \eqref{decay_psi_p}, we estimate
	\begin{equation}\label{wh_Lambda_upper_bound}
		\begin{aligned}
			\wh \Lambda_{t}(n)&=\frac{1}{2\pi}\int_{-\pi}^{\pi}\frac{\psi(p)}{1-t\psi(p)}e^{-\ii n p}\, \rd p\lesssim \int_{-\pi}^{\pi}\frac{|\psi(p)|}{1-t+[s(Wp)]^2}\, \rd p\\
			&\lesssim\int_{0}^{W^{-1}}\frac{1}{1-t+|Wp|^2}\, \rd p+\int_{ W^{-1}}^{\pi}|Wp|^{-2-\varepsilon}\, \rd p\lesssim \frac{1}{W\sqrt{1-t}}\le  \frac{1}{\ell_t(1-t)},\quad \forall n\in \Z.
		\end{aligned}
	\end{equation}
	To derive the decay of $\wh \Lambda_{t}$ for $|n|\geq W(1-t)^{-1/2}$, we rewrite \eqref{eq:whLambdaxi} (with $\xi=t$) as
	\begin{align}
		\wh \Lambda_{t}(n) &=(\ii n)^{-\ceil{\alpha}}\frac{1}{2\pi}\int_{|p|\leq |n|^{-1}}\pa{\frac{\psi(p)}{1-t\psi(p)}}^{\pa{\ceil{\alpha}}}e^{-\ii n p}\, \rd p \nonumber\\
		&+(\ii n)^{-\ceil{\alpha}}\frac{1}{2\pi}\int_{|n|^{-1}<|p|\leq \pi}\pa{\frac{\psi(p)}{1-t\psi(p)}}^{\pa{\ceil{\alpha}}}e^{-\ii n p}\, \rd p=:I_1+I_2.\label{eq:Lambda=I1I2}
	\end{align}
	To estimate $I_1$, using \eqref{spectral_gap_psi_p} and \eqref{singularity_psi_p}, we obtain that for any small constant $0<\varepsilon<1+\al-\ceil{\al}$,
	\begin{align}
		\absa{I_1}&\lesssim  |n|^{-\ceil{\alpha}}\int_{0}^{|n|^{-1}} \frac{W^{\ceil{\alpha}}|Wp|^{\alpha-\ceil{\alpha}-\varepsilon}}{(1-t+|Wp|^2)^2}\,\rd p+|n|^{-\ceil{\alpha}}\int_{0}^{|n|^{-1}}\frac{W^{\ceil{\alpha}}}{(1-t+|Wp|^2)^{\ceil{\alpha}/2+1}}\,\rd p \nonumber\\
		&\lesssim W^{\alpha-\varepsilon}|n|^{-\alpha-1+\varepsilon}(1-t)^{-2}+W^{\ceil{\alpha}}|n|^{-\ceil{\alpha}-1}(1-t)^{-\ceil{\alpha}/2-1}\nonumber \\
		&\lesssim \frac{(|n|/W)^\e}{W\sqrt{1-t}}\pa{\frac{|n|}{W(1-t)^{-1/2}}}^{-\alpha-1}.\label{eq:I1}
	\end{align}
	In the first step above, we used the following bound for $|p|\le |n|^{-1} \le W^{-1}$, derived from \eqref{singularity_psi_p}: 
	\begin{equation}\label{derivative_bound_0}
		\begin{aligned}
			\pa{\frac{\psi(p)}{1-t\psi(p)}}^{\pa{\ceil{\alpha}}}&\lesssim  \sum_{k=1}^{\ceil{\alpha}}\sum_{a_1,\ldots ,a_k\geq 1,\,a_1+\cdots+a_k=\ceil{\alpha}}\frac{W^{\ceil{\alpha}}|Wp|^{\#\ha{i:a_i=1}}\p{1+\mathbf{1}_{k=1}\cdot\absa{Wp}^{\alpha-\ceil{\alpha}-\varepsilon}}}{(1-t+|Wp|^2)^{k+1}}\\
			&\lesssim\frac{W^{\ceil{\alpha}}|Wp|^{\alpha-\ceil{\alpha}-\varepsilon}}{(1-t+|Wp|^2)^2}+\sum_{k=2}^{\ceil{\alpha}}\sum_{l=0}^{k} \frac{W^{\ceil{\alpha}}|Wp|^{l}}{(1-t+|Wp|^2)^{(\ceil{\alpha}+l)/2+1}}\\
			&\lesssim \frac{W^{\ceil{\alpha}}|Wp|^{\alpha-\ceil{\alpha}-\varepsilon}}{(1-t+|Wp|^2)^2}+ \frac{W^{\ceil{\alpha}}}{(1-t+|Wp|^2)^{\ceil{\alpha}/2+1}},
		\end{aligned}
	\end{equation}
	where $l:=\#\ha{i:a_i=1}$, and we used the inequality $\ceil{\al}\ge l+2(k-l)$ in the second step. 
    
	To estimate $I_2$, we apply integration by parts once more and obtain
	\begin{align}
		I_2&=(\ii n)^{-\ceil{\alpha}-1}\frac{1}{2\pi}\int_{|n|^{-1}<|p|\leq \pi}\pa{\frac{\psi(p)}{1-t\psi(p)}}^{\pa{\ceil{\alpha}+1}}e^{-\ii n p}\, \rd p \nonumber\\
		&+(\ii n)^{-\ceil{\alpha}-1}\frac{1}{2\pi}\qa{\left.\frac{\rd^{\ceil{\alpha}} }{\rd p^{\ceil{\alpha}}}\frac{\psi(p)}{1-t\psi(p)}\right|_{p=|n|^{-1}}e^{-\ii \sgn(n)}-\left.\frac{\rd^{\ceil{\alpha}} }{\rd p^{\ceil{\alpha}}}\frac{\psi(p)}{1-t\psi(p)}\right|_{p=-|n|^{-1}}e^{\ii \sgn(n)}}.\label{eq:I2}
	\end{align}
	Using \eqref{derivative_bound_0}, we bound the boundary terms in \eqref{eq:I2} by
	\begin{align*}
		\OO\pa{\frac{1}{|n|^{\ceil{\alpha}+1}}\pa{ \frac{W^{\ceil{\alpha}}|Wn^{-1}|^{\alpha-\ceil{\alpha}-\varepsilon}}{(1-t)^2}+ \frac{W^{\ceil{\alpha}}}{(1-t)^{\ceil{\alpha}/2+1}}}}=\OO\bigg({\frac{(|n|/W)^\e}{W\sqrt{1-t}}\pa{\frac{|n|}{W(1-t)^{-1/2}}}^{-\alpha-1}}\bigg).
	\end{align*}
	For the first term in \eqref{eq:I2}, using \eqref{spectral_gap_psi_p}--\eqref{singularity_psi_p} and an argument analogous to \eqref{derivative_bound_0}, we obtain
	\begin{equation*}
		\begin{aligned}
			&~ (\ii n)^{-\ceil{\alpha}-1}\frac{1}{2\pi}\int_{|n|^{-1}<|p|\leq \pi}\pa{\frac{\psi(p)}{1-t\psi(p)}}^{\pa{\ceil{\alpha}+1}}e^{-\ii n p}\, \rd p\\
			\lesssim & ~\frac{1}{|n|^{\ceil{\alpha}+1}}\int_{|n|^{-1}}^{W^{-1}}\qa{\frac{W^{\ceil{\alpha}+1}|Wp|^{\alpha-\ceil{\alpha}-1-\varepsilon}}{(1-t+|Wp|^2)^2} +\frac{W^{\ceil{\alpha}+1}}{(1-t+|Wp|^2)^{(\ceil{\alpha}+1)/2+1}}}\rd p + \frac{1}{|n|^{\ceil{\alpha}+1}}\int_{W^{-1}}^{\pi}\frac{W^{\ceil{\al}+1}}{(1+|Wp|)^{2+\e}} \,\rd p\\
			\lesssim & ~|n|^{-\ceil{\alpha}-1}W^{\ceil{\alpha}}\pa{ |Wn^{-1}|^{\alpha-\ceil{\alpha}-\varepsilon}(1-t)^{-2}+ (1-t)^{-\ceil{\alpha}/2-1}+ 1}\\        \lesssim  &~ {\frac{(|n|/W)^\e}{W\sqrt{1-t}}\pa{\frac{|n|}{W(1-t)^{-1/2}}}^{-\alpha-1}}.  
		\end{aligned}
	\end{equation*}
	Combining the estimates for the boundary terms and the integral term, we conclude 
	\begin{equation}\label{eq:I20}
		\begin{aligned}
		 \absa{I_2}\lesssim \frac{(|n|/W)^\e}{W\sqrt{1-t}}\pa{\frac{|n|}{W\p{1-t}^{-1/2}}}^{-\alpha-1}. \nc
		\end{aligned}
	\end{equation}
    
Substituting \eqref{eq:I1} and \eqref{eq:I20} into \eqref{eq:Lambda=I1I2} and using \eqref{wh_Lambda_upper_bound}, we obtain
	\begin{equation}\label{eq:Lambdatn}
		\begin{aligned}
			\abs{\wh\Lambda_t(n)}\lesssim \frac{(|n|/W+1)^\e}{W\sqrt{1-t}}\pa{\frac{|n|}{W(1-t)^{-1/2}}+1}^{-\alpha-1},\quad \forall n\in \Z.
		\end{aligned}
	\end{equation}
	Then, substituting this estimate into \eqref{Poisson_expansion_Theta} yields
	\begin{equation*}
		\begin{aligned}
			\pa{\frac{S}{1-tS}}_{xy}&\lesssim \frac{1}{W\sqrt{1-t}}\sum_{n\in \Z} \pa{\frac{|y-x+nN|}{W}+1}^\e\pa{\frac{|y-x+nN|}{W(1-t)^{-1/2}}+1}^{-\alpha-1}\\
			&\lesssim \frac{(N/W)^\e}{\ell_t(1-t)}\pa{\frac{|y-x|_N}{\ell_t}+1}^{-\alpha-1},
		\end{aligned}
	\end{equation*}
	which proves \eqref{Theta_bound_opposite_charge}, since $\e>0$ is arbitrary.
	The proof of \eqref{Theta_bound_equaled_charge} proceeds along the same lines, but is in fact considerably simpler. Indeed, one analyzes \smash{$\wh \Lambda_{t\xi}$} using the same integration-by-parts argument as above, while observing that $\min_p |1-t\xi \psi(p)|\gtrsim 1$. This uniform lower bound eliminates the singular behavior near $p=0$, so no delicate splitting of the integration region is required. We therefore omit the details.

	It remains to prove the estimates \eqref{Theta_difference_bound} and \eqref{zero_mode_removed_bound}. We claim that for $t$ satisfying $\ell_t<N$,
	\begin{equation}\label{eq:diff_1step}
		\absa{\pa{\frac{S}{1-tS}}_{0x}-\pa{\frac{S}{1-tS}}_{0,x+1}}\prec \frac{1}{W^{2}}\pa{\frac{|x|_N}{W(1-t)^{-1/2}}+1}^{-\alpha-2}.
	\end{equation}
	Summing this bound over $x$ and invoking translation invariance immediately yields \eqref{Theta_difference_bound}.
	To prove \eqref{eq:diff_1step}, we again use the representation \eqref{Poisson_expansion_Theta} and compute
	\begin{equation}\label{eq:Lambdan+1-n}
		\wh\Lambda_t(n)-\wh\Lambda_t(n+1)=\frac{1}{2\pi}\int_{-\pi}^{\pi}\frac{\psi(p)(1-e^{-\ii p})}{1-t\psi(p)}e^{-\ii np}\,\rd p.
	\end{equation}
	Applying integration by parts to \eqref{eq:Lambdan+1-n} $(\ceil{\al}+1)$ times, we obtain for any $n\ne 0$, 
	\begin{align}
		\wh\Lambda_t(n)-\wh\Lambda_t(n+1) &=(\ii n)^{-\ceil{\alpha}-1}\frac{1}{2\pi}\int_{|p|\leq |n|^{-1}}\pa{\frac{\psi(p)(1-e^{-\ii p})}{1-t\psi(p)}}^{\pa{\ceil{\alpha}+1}}e^{-\ii n p}\, \rd p \nonumber\\
		&+(\ii n)^{-\ceil{\alpha}-1}\frac{1}{2\pi}\int_{|n|^{-1}<|p|\leq \pi}\pa{\frac{\psi(p)(1-e^{-\ii p})}{1-t\psi(p)}}^{\pa{\ceil{\alpha}+1}}e^{-\ii n p}\, \rd p=:I_1+I_2.\label{eq:diffLambda=I1I2}
	\end{align}
	Using this representation and applying arguments analogous to those above to estimate $I_1$ and $I_2$, we deduce 
	\begin{equation}\label{eq:Lambda_diff0}
		\begin{aligned}
			\absa{\wh\Lambda_t(n)-\wh\Lambda_t(n+1)}\lesssim\frac{(|n|/W+1)^\e}{W^2}\pa{\frac{|n|}{W(1-t)^{-1/2}}+1}^{-\alpha-2},\quad \forall n\in\Z,
		\end{aligned}
	\end{equation}
	for any sufficiently small constant $\varepsilon>0$. 
	Substituting this bound into \eqref{Poisson_expansion_Theta} yields \eqref{eq:diff_1step}. 
	To prove \eqref{zero_mode_removed_bound}, we use the representation
	\begin{equation*}
		\begin{aligned}
			\pa{\frac{S}{1-tS}}_{xy}=\frac1N\sum_{\bp\in\bbT_N}\Lambda_{t}(\bp)\exp(\ii\bp(x-y)).
		\end{aligned}
	\end{equation*}
	From this we deduce
	\begin{equation*}
		\begin{aligned}
			\absa{\pa{\frac{S}{1-tS}}_{xy}-\frac{1}{N(1-t)}}\leq \frac{1}{N}\sum_{\bp\in\bbT_N:0<|\bp|\leq W^{-1}}\frac{1}{|W\bp|^2}+ \frac{1}{N}\sum_{\bp\in\bbT_N:|\bp|> W^{-1}}|\psi(\bp)|\lesssim \frac{N}{W^2},
		\end{aligned}
	\end{equation*}
	where in the second term we used \eqref{decay_psi_p} to control the contribution from $|p|>W^{-1}$.

	\medskip
	\noindent
	{\bf Case (ii): $\alpha\in [1,2)$.} The proofs of \eqref{Theta_bound_opposite_charge}, \eqref{Theta_bound_equaled_charge}, and \eqref{zero_mode_removed_bound} for $\alpha\in [1,2)$ proceed along the same lines as in the case $\alpha\ge 2$, and we therefore omit the details. The essential difference arises in the proof of \eqref{Theta_difference_bound}. We claim that, for any constant $\varepsilon>0$,
	\begin{equation}\label{eq:Lambda_diff}
		\absa{\wh\Lambda_t(n)-\wh \Lambda_t(n+1)}\lesssim \frac{1}{W^2}\pa{\frac{|n|}{W}+1}^{\alpha-2+\e}\pa{\frac{|n|}{W(1-t)^{-1/\alpha}}+1}^{-2\alpha}.
	\end{equation}
	Substituting this estimate into \eqref{Poisson_expansion_Theta} yields
	\begin{equation}\nonumber
		\absa{\pa{\frac{S}{1-tS}}_{0x}-\pa{\frac{S}{1-tS}}_{0,x+1}}\prec \frac{1}{W^2}\pa{\frac{|x|_N}{W}+1}^{\alpha-2}\pa{\frac{|x|_N}{W(1-t)^{-1/\alpha}}+1}^{-2\alpha},
	\end{equation}
	since $\e>0$ is arbitrary.
	Summing this bound over $x$ and using translation invariance yields \eqref{Theta_difference_bound}.
	
	The proof of \eqref{eq:Lambda_diff} follows the same strategy as that of \eqref{eq:Lambda_diff0} for $\alpha\ge 2$, with the additional complication that one must establish decay on two distinct scales.
	We first derive a uniform bound without spatial decay. Using \eqref{spectral_gap_psi_p} and \eqref{decay_psi_p} in \eqref{eq:Lambdan+1-n}, we obtain
	\begin{equation}\label{eq:diffLambda_nodecay}
		\begin{aligned}
			\absa{\wh\Lambda_t(n)-\wh\Lambda_t(n+1)}\lesssim \int_{0}^{W^{-1}}\frac{|p|}{1-t+|Wp|^\alpha}\,\rd p+\int_{W^{-1}}^{\pi}|p||\psi(p)|\rd p\lesssim \frac{1}{W^2}.
		\end{aligned}
	\end{equation}
	We next establish decay on the $W$-scale. For $|n|\ge W$, we decompose
	\begin{equation*}
		\begin{aligned}
			\wh\Lambda_t(n)-\wh\Lambda_t(n+1)&=\frac{1}{2\pi}\int_{|p|\leq |n|^{-1}}\frac{\psi(p)(1-e^{-\ii p})}{1-t\psi(p)}e^{-\ii np}\,\rd p +\frac{1}{2\pi}\int_{|p|> |n|^{-1}}\frac{\psi(p)(1-e^{-\ii p})}{1-t\psi(p)}e^{-\ii np}\,\rd p=:I_1+I_2.
		\end{aligned}
	\end{equation*}
	The first term is estimated directly:
	\begin{equation*}
		\begin{aligned}
			\absa{I_1}\lesssim \int_{0}^{|n|^{-1}}\frac{|p|}{1-t+|Wp|^{\alpha}}\,\rd p\lesssim W^{-\alpha}|n|^{\alpha-2}.
		\end{aligned}
	\end{equation*}
	For $I_2$, we perform one integration by parts and obtain 
	\begin{equation*}
		\begin{aligned}
			I_2&=  \frac{(\ii n)^{-1}}{2\pi}\int_{|p|> |n|^{-1}}\pa{\frac{\psi(p)(1-e^{-\ii p})}{1-t\psi(p)}}' e^{-\ii np}\rd p  +\frac{(\ii n)^{-1}}{2\pi}{\left.\frac{\psi(p)(1-e^{-\ii p})e^{-\ii np}}{1-t\psi(p)}\right|_{-|n|^{-1}}^{|n|^{-1}}}.
		\end{aligned}
	\end{equation*}
	Using \eqref{spectral_gap_psi_p}--\eqref{singularity_psi_p}, we bound this term by
	\begin{equation*}
		\begin{aligned}
			| I_2 | \lesssim (|n|/W)^\e\cdot W^{-\alpha}|n|^{\alpha-2}
		\end{aligned}
	\end{equation*}
	for any constant $\e>0$. Combining the above estimates with \eqref{eq:diffLambda_nodecay}, we conclude that
	\begin{equation}\label{eq:diffLambda_Wdecay}
		\begin{aligned}
			\absa{\wh\Lambda_t(n)-\wh\Lambda_t(n+1)}\lesssim \frac{1}{W^2}\pa{\frac{|n|}{W}+1}^{\alpha-2+\e},\quad \forall n\in \Z.
		\end{aligned}
	\end{equation}
	
	It remains to establish the decay for $|n|>W(1-t)^{-1/\alpha}$. In this regime, we employ the representation \eqref{eq:diffLambda=I1I2}. The term $I_1$ is estimated in the same manner as in \eqref{eq:I1}, while the second term can be treated analogously to \eqref{eq:I2}-\eqref{eq:I20}, where we replace $\lceil\alpha\rceil$ with $\lceil\alpha\rceil+1$. This yields
	\begin{equation}\label{eq:diffLambda_ltdecay}
		\begin{aligned}
			\absa{\wh\Lambda_t(n)-\wh\Lambda_t(n+1)}\lesssim (|n|/W)^\e \cdot W^{\alpha}|n|^{-2-\alpha}(1-t)^{-2}.
		\end{aligned}
	\end{equation}
	Combining this estimate with \eqref{eq:diffLambda_Wdecay} completes the proof of \eqref{eq:Lambda_diff}.

	\medskip
	\noindent
	{\bf Case (iii): $\alpha\in(0,1)$.} 
	We first prove \eqref{zero_mode_removed_bound} using the Fourier series
	\begin{equation}\label{eq:nonzeromode}
		\begin{aligned}
			\pa{\frac{S}{1-tS}}_{xy}-\frac{1}{N(1-t)}=\frac{1}{N}\sum_{\bp\in\bbT_N\setminus\ha{0}}\Lambda_t(\bp)\exp\pa{\ii \bp (x-y)}.
		\end{aligned}
	\end{equation}
	A direct estimate based on \eqref{spectral_gap_psi_p} and \eqref{decay_psi_p} yields
	\begin{equation}\label{eq:nonzeromode1}
		\begin{aligned}
			\absa{\pa{\frac{S}{1-tS}}_{xy}-\frac{1}{N(1-t)}}\lesssim \frac{1}{W}.
		\end{aligned}
	\end{equation}
	To obtain the spatial decay, assume that $|x-y|_N\ge W$ and decompose
	\begin{equation*}
		\begin{aligned}
			\pa{\frac{S}{1-tS}}_{xy}-\frac{1}{N(1-t)}&=\frac{1}{N}\sum_{\bp\in\bbT_N\setminus\ha{0}:|\bp|\leq 100|x-y|_N^{-1}}\Lambda_t(\bp)\exp\pa{\ii \bp (x-y)}\\
			&+\frac{1}{N}\sum_{\bp\in\bbT_N\setminus\ha{0}:|\bp|> 100|x-y|_N^{-1}}\Lambda_t(\bp)\exp\pa{\ii \bp (x-y)}=:I_1+I_2.
		\end{aligned}
	\end{equation*}
	The first term is estimated directly using \eqref{spectral_gap_psi_p}:
	\begin{equation}\label{eq:I1al01}
		|I_1|\lesssim W^{-\alpha}|x-y|_N^{\alpha-1}.    
	\end{equation}
	For the second term, we apply summation by parts and obtain
	\begin{align}
		I_2=&~\frac{1}{N}\frac{\exp\pa{\ii \delta_N (x-y)}}{1-\exp\pa{\ii\delta_{N}(x-y)}}\sum_{\bp\in\bbT_N:|\bp|> 100|x-y|_N^{-1}}\qa{\Lambda_t(\bp+\delta_N)-\Lambda_t(\bp)}\exp\pa{\ii \bp (x-y)}\nonumber\\
		&~ +\frac{1}{N}\frac{1}{1-\exp\pa{\ii\delta_{N}(x-y)}}\Lambda_t(-\bp_{x-y})\exp\pa{-\ii\bp_{x-y} (x-y)}  \nonumber\\
		&~  -\frac{1}{N}\frac{1}{1-\exp\pa{\ii\delta_{N}(x-y)}}\Lambda_t(\bp_{x-y}+\delta_N)\exp\pa{\ii (\bp_{x-y}+\delta_N) (x-y)},  \label{eq:I2alpha}
	\end{align}
	where $\delta_N=2\pi/N$ denotes the mesh size of $\mathbb T_N$, and $\bp_{x-y}$ is the rightmost point in $(-\infty,0]\cap\mathbb T_N$ satisfying $|\bp|>100|x-y|_N^{-1}$.
	Using \eqref{spectral_gap_psi_p}, the boundary terms are bounded by
	\begin{equation*}
		\begin{aligned}
			\OO\pa{|x-y|_N^{-1}\pa{\frac{|x-y|_N}{W}}^{\alpha}}=\OO\pa{ W^{-\alpha}|x-y|_N^{\alpha-1}}.
		\end{aligned}
	\end{equation*}
	For the first term on the right-hand side of \eqref{eq:I2alpha}, we use \eqref{spectral_gap_psi_p}--\eqref{singularity_psi_p} to obtain the bound
	\begin{equation*}
		\begin{aligned}
			&\frac{1}{|x-y|_N}\int_{|x-y|_N^{-1}}^{W^{-1}}\frac{W|Wp|^{\alpha-1}}{(1-t+|Wp|^{\alpha})^2}\,\rd p+\frac{1}{|x-y|_N}\int_{W^{-1}}^{\pi}\frac{W}{(|Wp|+1)^{2+\varepsilon}}\,\rd p\lesssim W^{-\alpha}|x-y|_N^{\alpha-1}. 
		\end{aligned}
	\end{equation*}
	Combining the above two estimates yields  
	\begin{equation}\label{eq:I2al01}
		|I_2|\lesssim W^{-\alpha}|x-y|_N^{\alpha-1}. 
	\end{equation}
	Together with \eqref{eq:I1al01} and \eqref{eq:nonzeromode1}, this proves \eqref{zero_mode_removed_bound} for $\alpha\in(0,1)$.

	Note that when $\ell_t = N$, that is, when $1-t \le (W/N)^{\alpha}$, the bound \eqref{Theta_bound_opposite_charge} follows directly from \eqref{zero_mode_removed_bound}.
	We therefore consider the case $\ell_t < N$. In this regime, we prove \eqref{Theta_bound_opposite_charge} using \eqref{eq:Lambda=I1I2} with $\lceil\alpha\rceil=1$. Repeating the argument between \eqref{wh_Lambda_upper_bound} and \eqref{eq:Lambdatn}, and applying the estimates \eqref{spectral_gap_psi_p}--\eqref{singularity_psi_p}, we obtain
	\begin{equation*}
		\begin{aligned}
			\abs{\wh\Lambda_t(n)}\lesssim W^{-1} (|n|/W+1)^{\al-1}, \quad \forall n\in \Z,
		\end{aligned}
	\end{equation*}
	and, for $|n|\ge W(1-t)^{-1/\alpha}$,
	\begin{equation*}
		\begin{aligned}
			\abs{\wh\Lambda_t(n)}\lesssim W^{\al} |n|^{-1-\al} (1-t)^{-2} .
		\end{aligned}
	\end{equation*}
	Combining these two bounds yields
	\begin{equation}\label{eq:Lambdatn01_2}
		\begin{aligned}
			\abs{\wh\Lambda_t(n)}\lesssim \frac{1}{W}\pa{\frac{|n|}{W}+1}^{\al-1}\pa{\frac{|n|}{W(1-t)^{-1/\alpha}}+1}^{-2\alpha},\quad \forall n\in \Z.
		\end{aligned}
	\end{equation}
	Substituting this estimate into \eqref{Poisson_expansion_Theta} yields \eqref{Theta_bound_opposite_charge}.
	The bounds \eqref{Theta_bound_equaled_charge} and \eqref{Theta_difference_bound} can be established by analogous integration-by-parts arguments as in the previous two cases; we omit the details.

	\subsection{Proof of Lemma \ref{lemma_B_t}}\label{appendix_additional_proof_B_t_and_U_s_t}
	For $\alpha\in(0,1)$, the convolution inequality \eqref{convolution_2_B} follows immediately from the following normalized tail domination estimate:
	\begin{equation}\label{bound_n_t_d}
		\begin{aligned}
			{B_s\p{r}}/{\|B_s\p{|\cdot|}\|_1}\lesssim {B_t\p{r}}/{V_t\p{r}},\quad \forall 0\le r \le N,
		\end{aligned}
	\end{equation}
	where the local cumulative function $V_t\p{r}$ is defined by
	\begin{equation*}
		V_t\p{r}:=\sum_{ \ell\in \Z:\, 0\leq  \ell \leq r}B_t\p{\ell}.
	\end{equation*}
	Assuming \eqref{bound_n_t_d}, together with \eqref{eq:L1Bu} and the elementary observation that for any fixed constant $C>0$ one has $B_u(r)\asymp B_u(Cr)$ uniformly for $u\in[s,t]$ and $r\ge0$, we obtain
	\begin{equation}\label{bound_half_half}
		\begin{aligned}
			&~\sum_{a\in\ZL}B_s\p{|x-a|}B_t\p{|a-y|} \\
			\leq&~ B_t(|x-y|/2)\sum_{|a-y|\geq |x-y|/2}B_s\pa{|x-a|}+B_s\pa{|x-y|/2}\sum_{|a-y|<|x-y|/2}B_t\pa{|a-y|}\\
			\lesssim&~\frac{1}{1-s}B_t\pa{|x-y|}+B_s\pa{|x-y|}V_t\pa{|x-y|}\lesssim \frac{1}{1-s}B_t\pa{|x-y|}.
		\end{aligned}
	\end{equation}
	
	It therefore remains to prove \eqref{bound_n_t_d}. For $\alpha\in(0,1)$, the local cumulative function can be estimated as
	\begin{align}
		V_t\pa{r}\lesssim&\sum_{0\leq  \ell\leq r} \frac{1}{W}\pa{\frac{\ell}{W}+1}^{\alpha-1}\pa{\frac{\ell}{\ell_t}+1}^{-2\alpha}+\frac{r+1}{N(1-t)}\cdot \mathbf{1}_{\ell_t=N}\nonumber\\
		\lesssim& \begin{cases} 
			(r+1)/W+(r+1)/[N(1-t)]\cdot \mathbf{1}_{\ell_t=N},\ &\txt{for }0\leq r\leq W\\
			(r/W)^\alpha+r/[N(1-t)]\cdot \mathbf{1}_{\ell_t=N},\ &\txt{for }W<r\leq \ell_t\\
			(\ell_t/W)^{\alpha}+r/[N(1-t)]\cdot \mathbf{1}_{\ell_t=N},\ &\txt{for }r>\ell_t
		\end{cases},\label{bound_V_t_r_0_1}\\
		\lesssim&(r+1)\qa{\frac{1}{W}\pa{\frac{r}{W}+1}^{\alpha-1}+\frac{1}{N(1-t)}}\pa{\frac{r}{\ell_t}+1}^{-\alpha}=:\mathcal{V}_t\pa{r}.\nonumber
	\end{align}	
	Combining \eqref{bound_V_t_r_0_1} with the definition \eqref{def_B_t}, we obtain
	\begin{equation}
		\begin{aligned}
			{B_s\p{r}}/{\mathcal{V}_s\p{r}}\le {B_t\p{r}}/{\mathcal{V}_t\p{r}}.
		\end{aligned}
	\end{equation}
	Since $\mathcal V_s(r)\lesssim \|B_s(|\cdot|)\|_1\asymp (1-s)^{-1}$, it follows that
	\begin{equation}
		\begin{aligned}
			\frac{B_s\p{r}}{\|B_s\p{|\cdot|}\|_1}\lesssim \frac{B_s\p{r}}{\mathcal{V}_s\p{r}}\le \frac{B_t\p{r}}{\mathcal{V}_t\p{r}}\lesssim \frac{B_t\p{r}}{V_t\p{r}}.
		\end{aligned}
	\end{equation}
	This proves \eqref{bound_n_t_d}, and hence completes the proof of the convolution inequality \eqref{convolution_2_B}.

	\subsection{Proof of Lemma \ref{lemma_U_s_t}}\label{appendix_additional_proof_B_t_and_U_s_t2}
	The first estimate \eqref{bound_cU_s_t_1} follows directly from the decomposition \eqref{U_s_t_equal_1_Theta}, together with the assumptions \eqref{Theta_bound_opposite_charge}, \eqref{Theta_bound_equaled_charge}, and the convolution bound \eqref{convolution_2_B}.
	For the improved bound \eqref{bound_cU_s_t_2}, consider first the case $\sigma_1=\sigma_2$. In this situation, the desired bound follows immediately from the decomposition \eqref{U_s_t_equal_1_Theta}, the estimate \eqref{Theta_bound_equaled_charge}, and the elementary inequality
	\begin{equation}
		\sum_{a\in\ZL}B_0\pa{|x-a|}B_s\pa{|a-y|}\lesssim B_s\pa{|x-y|}. \label{eq:SBS}
	\end{equation} 
	
	Next, assume $\sigma_1\neq\sigma_2$. Without loss of generality, suppose that $X$ satisfies the right sum-zero property. For notational simplicity, write $\cU_{s,t}\equiv \cU_{s,t}^{\bsigma}$ and $\Theta_t\equiv \Theta_t^{\bsigma}$. For any $x,y\in\ZL$, using the decomposition \eqref{U_s_t_equal_1_Theta}, we expand $(\cU_{s,t}\circ X)_{xy}$ as
	\begin{equation}\label{decomposition_U_s_t_X}
		(\cU_{s,t}\circ X)_{xy}=X_{xy}+(t-s)\sum_{a\in\ZL}\Theta_{t,xa}X_{a y}+(t-s)\sum_{b\in\ZL}X_{xb}\Theta_{t,b y}+(t-s)^2\sum_{a,b\in\ZL}\Theta_{t,xa}X_{ab}\Theta_{t,b y}.
	\end{equation}
	By the monotonicity of $B_u(|x-y|)$ in $u$, together with the bound \eqref{Theta_bound_opposite_charge} and the convolution estimate \eqref{convolution_2_B}, the first three terms on the RHS satisfy
	\begin{equation}\label{eq:simple_conv}
		X_{xy}+(t-s)\sum_{a\in\ZL}\Theta_{t,xa}X_{a y}+(t-s)\sum_{b\in\ZL}X_{xb}\Theta_{t,b y} \prec      \|X/B_{s}\|_{\max}\cdot B_t\pa{|x-y|}.
	\end{equation}
	It remains to estimate the fourth term in \eqref{decomposition_U_s_t_X}. Using the right sum-zero property \eqref{right_sum_zero}, the bound \eqref{Theta_bound_opposite_charge}, and the difference estimate \eqref{Theta_difference_bound}, we obtain for $t$ with $\ell_t<N$ that
	\begin{align}
		&~\sum_{a,b\in\ZL}\Theta_{t,xa}X_{ab}\Theta_{t,b y}=\sum_{a,b\in\ZL}\Theta_{t,xa}X_{ab}\pa{\Theta_{t,b y}-\Theta_{t,a y}}\nonumber\\
		\prec&~ \frac{\|X/B_{s}\|_{\max}}{W}\sum_{a,b\in\ZL}\pa{\frac{|a-y|}{W}+1}^{-1}B_{t}\pa{|x-a|}\pa{\frac{|a-b|}{W}+1}^{\alpha}\pa{\frac{|a-b|}{\ell_s}+1}^{-2\alpha}B_t\pa{|b-y|}\nonumber\\
		+&~ \frac{\|X/B_{s}\|_{\max}}{W}\sum_{a,b\in\ZL}\pa{\frac{|b-y|}{W}+1}^{-1}B_{t}\pa{|x-a|}\pa{\frac{|a-b|}{W}+1}^{\alpha}\pa{\frac{|a-b|}{\ell_s}+1}^{-2\alpha}B_t\pa{|a-y|}\nonumber\\
		\lesssim&~ \frac{\|X/B_{s}\|_{\max}}{W}\pa{\frac{\ell_s}{W}}^{2\alpha}\sum_{a,b\in\ZL}\pa{\frac{|a-y|}{W}+1}^{-1}B_{t}\pa{|x-a|}\pa{\frac{|a-b|}{W}+1}^{-\alpha}B_t\pa{|b-y|}\nonumber\\
		+&~ \frac{\|X/B_{s}\|_{\max}}{W}\pa{\frac{\ell_s}{W}}^{2\alpha}\sum_{a,b\in\ZL}\pa{\frac{|b-y|}{W}+1}^{-1}B_{t}\pa{|x-a|}\pa{\frac{|a-b|}{W}+1}^{-\alpha}B_t\pa{|a-y|}=:I_1+I_2.\label{eq:ThetaXTheta}
	\end{align}
For the term $I_1$, we claim the following two convolution bounds:
	\begin{align}\label{triangular_convolution_example_1}
		&\sum_{b\in\ZL}\pa{\frac{|a-b|}{W}+1}^{-\alpha}\cdot W^{-1}\pa{\frac{|b-y|}{W}+1}^{\alpha-1}\pa{\frac{|b-y|}{\ell_t}+1}^{-2\alpha}\prec \pa{\frac{|a-y|}{\ell_t}+1}^{-\alpha},\\
		&\sum_{a\in\ZL}\pa{\frac{|a-y|}{W}+1}^{-1}\pa{\frac{|a-y|}{\ell_t}+1}^{-\alpha}\cdot W^{-1}\pa{\frac{|x-a|}{W}+1}^{\alpha-1}\pa{\frac{|x-a|}{\ell_t}+1}^{-2\alpha} \nonumber\\
		&\prec \pa{\frac{|x-y|}{W}+1}^{\alpha-1}\pa{\frac{|x-y|}{\ell_t}+1}^{-2\alpha}\lesssim WB_t\pa{|x-y|}.\label{triangular_convolution_example_2}
	\end{align}
	Combining \eqref{triangular_convolution_example_1} and \eqref{triangular_convolution_example_2}, we conclude that
	\begin{equation}\label{eq:I1_TXT}
		I_1 \prec  {\|X/B_s\|_{\max}}\pa{\frac{\ell_s}{W}}^{2\alpha} B_t\pa{|x-y|} \le \frac{\|X/B_s\|_{\max}}{(1-s)^2} B_t\pa{|x-y|}.
	\end{equation}
For the proof of \eqref{triangular_convolution_example_1}, if $|a-y|\le \ell_t$, we simply use the bound $(|b-y|/\ell_t+1)^{-2\alpha}\le 1$. The desired estimate then follows immediately from Young's inequality applied to the convolution in $b$.
If $|a-y|>\ell_t$, we argue as in \eqref{bound_half_half} by decomposing the summation region according to whether $|a-b|\ge |a-y|/2$ or $|a-b|< |a-y|/2$. This yields
	\begin{equation*}
		\begin{aligned}
			&\frac1W\sum_{b\in\ZL}\pa{\frac{|a-b|}{W}+1}^{-\alpha}\pa{\frac{|b-y|}{W}+1}^{\alpha-1}\pa{\frac{|b-y|}{\ell_t}+1}^{-2\alpha}\\
			\lesssim& \pa{\frac{\ell_t}{W}}^\alpha\pa{\frac{|a-y|}{W}+1}^{-\alpha}+\pa{\frac{|a-y|}{W}}^{1-\alpha}\pa{\frac{|a-y|}{W}+1}^{\alpha-1}\pa{\frac{|a-y|}{\ell_t}+1}^{-2\alpha}\lesssim \pa{\frac{|a-y|}{\ell_t}+1}^{-\alpha}.
		\end{aligned}
	\end{equation*}
This establishes \eqref{triangular_convolution_example_1}. 
We next prove \eqref{triangular_convolution_example_2}. If $|x-y|>\ell_t$, a similar decomposition argument based on whether $|a-y|\ge |x-y|/2$ or $|a-y|< |x-y|/2$ gives the desired bound. If $W<|x-y|\le \ell_t$, the claim follows directly from the estimate 
	\begin{equation*}
		\begin{aligned}
			\frac1W \sum_{a\in\ZL}\pa{\frac{|a-y|}{W}+1}^{-1}\pa{\frac{|x-a|}{W}+1}^{\alpha-1}\prec \pa{\frac{|x-y|}{W}+1}^{\alpha-1}.
		\end{aligned}
	\end{equation*}
Finally, the remaining case $|x-y|\le W$ is immediate.
	
	Next, we estimate the term $I_2$ in \eqref{eq:ThetaXTheta}. By the triangle inequality, 
	\begin{equation*}
		\begin{aligned}
			B_t\pa{|x-a|}B_t\pa{|a-y|}&=B_t\pa{|x-a|\vee|a-y|}B_t\pa{|x-a|\wedge|a-y|} \lesssim B_t\pa{|x-y|}\pa{B_t\pa{|x-a|}+B_t\pa{|a-y|}}.
		\end{aligned}
	\end{equation*}
	Moreover, for any $w,w'\in\ZL$, we have
	\begin{equation*}
		\sum_{a,b\in\ZL}B_t\p{|w-a|}\pa{{|a-b|}/{W}+1}^{-\alpha}\pa{{|b-w'|}/{W}+1}^{-1}\prec W,
	\end{equation*}
	which follows immediately from the elementary bound $B_t\pa{|w-a|}\lesssim W^{-1}(|w-a|/W+1)^{\alpha-1}$ for $t$ with $\ell_t<N$. Combining these two estimates, we deduce
	\begin{align*}
		I_2&\lesssim \frac{\|X/B_{s}\|_{\max}}{(1-s)^2} B_{t}(|x-y|) \cdot \frac1W\sum_{a,b\in\ZL}\pa{\frac{|b-y|}{W}+1}^{-1}\pa{\frac{|a-b|}{W}+1}^{-\alpha}\pa{B_t\pa{|x-a|}+B_t\pa{|a-y|}} \\
		&\prec\frac{\|X/B_{s}\|_{\max}}{(1-s)^2} B_{t}(|x-y|). 
	\end{align*}
	Substituting this bound together with \eqref{eq:I1_TXT} into \eqref{eq:ThetaXTheta}, we obtain
	\[(t-s)^2\sum_{a,b\in\ZL}\Theta_{t,xa}X_{ab}\Theta_{t,b y} \prec \|X/B_{s}\|_{\max}\cdot B_{t}(|x-y|). \]
	Combined with \eqref{eq:simple_conv}, this proves \eqref{bound_cU_s_t_2} in the case $\ell_t<N$.

	We next consider the regime $\ell_s=\ell_t=N$. Using \eqref{Theta_bound_opposite_charge} and \eqref{zero_mode_removed_bound}, we write
	\begin{equation*}
		\begin{aligned}
			&~\sum_{a,b\in\ZL}\Theta_{t,xa}X_{ab}\Theta_{t,by}=\sum_{a,b\in\ZL}\Theta_{t,xa}X_{ab}\pa{\Theta_{t,by}-\frac{1}{N(1-t)}}\\
			\prec&~\|X/B_s\|_{\max}\sum_{a,b\in\ZL}B_t\pa{|x-a|}B_s\pa{|a-b|}\upcirc{B_t}\pa{|b-y|}\lesssim\frac{\|X/B_s\|_{\max}}{1-s}\sum_{b\in\ZL}B_t\pa{|x-b|}\upcirc{B_t}\pa{|b-y|}\\
			\asymp&~\frac{\|X/B_s\|_{\max}}{1-s}\sum_{b\in\ZL}\bigg[\frac{1}{W}\pa{\frac{|x-b|}{W}+1}^{\alpha-1}+\frac{1}{N(1-t)}\bigg]\cdot \frac1W\pa{\frac{|b-y|}{W}+1}^{\alpha-1}\\
			\lesssim&~\frac{\|X/B_s\|_{\max}}{1-s}\pa{\frac{N}{W}}^{\alpha}\bigg[\frac{1}{W}\pa{{\frac{|x-y|}{W}+1}}^{\alpha-1}+\frac{1}{N(1-t)}\bigg]\lesssim \frac{\|X/B_s\|_{\max}}{(1-s)^2}B_t\pa{|x-y|},
		\end{aligned}
	\end{equation*}
	where the third step follows from the convolution inequality \eqref{convolution_2_B}, and the final step uses the relation $1-s\le (W/N)^{\alpha}$. Together with \eqref{eq:simple_conv}, this establishes \eqref{bound_cU_s_t_2} in the case $\ell_s=N$.     
	Finally, the intermediate regime $\ell_s<N=\ell_t$ follows readily from the two cases above by introducing an intermediate time $t_0=1-(W/N)^{\alpha}$ between $s$ and $t$. This completes the proof of the improved estimate \eqref{bound_cU_s_t_2}.

	The bound \eqref{bound_cU_s_t_3} can be derived by a similar argument. When $\sigma_1=\sigma_2$, the claim again follows from \eqref{convolution_2_B}, \eqref{Theta_bound_equaled_charge}, and \eqref{eq:SBS}. If $\sigma_1\neq\sigma_2$ and $\bv$ satisfies the sum-zero property, we use the decomposition \( \cU_{s,t}\bv = \bv + (t-s)\Theta_t\bv,\) so it suffices to bound $\Theta_t\bv$. Using \eqref{Theta_bound_opposite_charge} and \eqref{zero_mode_removed_bound}, we obtain
	\begin{equation*}        \sum_{a\in\ZL}\Theta_{t,xa}\bv_a=\sum_{a\in\ZL}\pa{\Theta_{t,xa}-\frac{1}{N(1-t)}}\bv_{a}\prec \|\bv\|_{\max} \sum_{a\in\ZL}\upcirc{B_t}\pa{|x-a|}\lesssim \frac{\|\bv\|_{\max}}{1-s}.
	\end{equation*}
	This yields the desired bound \eqref{bound_cU_s_t_3}.

	\subsection{Proof of \Cref{convolution_B_half_B}}\label{sec_proof_of_convolution_B_half_B}
	
	For \eqref{convolution_half_B1}, if $|x-y|\geq \ell_t$, we apply the triangle inequality to obtain
	\begin{equation*}
		\begin{aligned}
			B_t^{1/2}(|x-a|)B_s^{1/2}(|a-y|)\lesssim &~B_t^{1/2}(|x-y|)B_s^{1/2}(|a-y|)+B_t^{1/2}(|x-a|)B_s^{1/2}(|x-y|)\\
			\asymp&~ B_t^{1/2}(|x-y|)B_s^{1/2}(|a-y|)+\pa{\frac{\ell_s}{\ell_t}}^{\alpha}B_t^{1/2}(|x-a|)B_t^{1/2}(|x-y|)\\
			\asymp&~ B_t^{1/2}(|x-y|)B_s^{1/2}(|a-y|)+\frac{1-t}{1-s}B_t^{1/2}(|x-a|)B_t^{1/2}(|x-y|).
		\end{aligned}
	\end{equation*}
	Combining this with the bound, valid for any $u\in\{s,t\}$ and $w\in\{x,y\}$,
	\begin{equation*}
		\begin{aligned}
			\sum_{a\in\ZL}B_t^{1/2}(|x-a|)B_u^{1/2}(|a-w|)
			\lesssim&~  \frac{1}{W}\sum_{a\in\ZL}\pa{\frac{|x-a|}{W}+1}^{(\alpha-1)/2}\pa{\frac{|a-w|}{W}+1}^{(\alpha-1)/2}\pa{\frac{|a-w|}{\ell_u}+1}^{-\alpha}\\
			\prec&~\pa{\frac{\ell_u}{W}}^{\alpha}\asymp \frac{1}{1-u}, 
		\end{aligned}
	\end{equation*}
	we obtain the desired estimate \eqref{convolution_half_B1} for $|x-y|\ge \ell_t$. We now consider the case $|x-y|<\ell_t$. In this regime,
	\begin{equation*}
		\begin{aligned}
			\sum_{a\in\ZL}B_t(|x-a|)B_s^{1/2}(|a-y|)\prec&~ \frac{1}{W^{3/2}} \sum_{a\in\ZL}\pa{\frac{|x-a|}{W}+1}^{\alpha-1}\pa{\frac{|a-y|}{W}+1}^{(\alpha-1)/2}\pa{\frac{|a-y|}{\ell_s}+1}^{-\alpha}\\
			\lesssim&~ \frac{1}{W^{3/2}}\pa{\frac{\ell_s}{W}}^{\alpha} \sum_{a\in\ZL}\pa{\frac{|x-a|}{W}+1}^{\alpha-1}\pa{\frac{|a-y|}{W}+1}^{-(\alpha+1)/2}.
		\end{aligned}
	\end{equation*}
	Using $(\ell_s/W)^{\alpha}=(1-s)^{-1}$ together with the triangle inequality 
	\begin{equation*}
		\begin{aligned}
			&\pa{\frac{|x-a|}{W}+1}^{(\alpha-1)/2}\pa{\frac{|a-y|}{W}+1}^{(\alpha-1)/2}\\
			\lesssim& \pa{\frac{|x-y|}{W}+1}^{(\alpha-1)/2}\qa{\pa{\frac{|x-a|}{W}+1}^{(\alpha-1)/2}+\pa{\frac{|a-y|}{W}+1}^{(\alpha-1)/2}},
		\end{aligned}
	\end{equation*}
	we conclude that
	\begin{equation*}
		\begin{aligned}
			&\sum_{a\in\ZL}B_t(|x-a|)B_s^{1/2}(|a-y|)\prec \frac{1}{1-s}  B_t^{1/2}(|x-y|) \cdot \frac{1}{W}\sum_{a\in\ZL}\pa{\frac{|x-a|}{W}+1}^{(\alpha-1)/2}\\
			&\qquad \times\qa{\pa{\frac{|x-a|}{W}+1}^{(\alpha-1)/2}+\pa{\frac{|a-y|}{W}+1}^{(\alpha-1)/2}}\pa{\frac{|a-y|}{W}+1}^{-\alpha}\prec  \frac{1}{1-s}B_t^{1/2}(|x-y|).
		\end{aligned}
	\end{equation*}
	This establishes \eqref{convolution_half_B1}. We now turn to \eqref{convolution_half_B2}. Let $s_0:=1-(W/N)^{\alpha}$. Applying \eqref{convolution_half_B1}, we obtain
	\begin{equation*}
		\begin{aligned}
			\sum_{a\in\ZL}\upcirc{B_t}(|x-a|)B_s^{1/2}(|a-y|)\prec&~ \sum_{a\in\ZL} B_{s_0}(|x-a|)B_{s_0}^{1/2}(|a-y|)+ \sum_{a\in\ZL}B_{s_0}(|x-a|) \pa{\frac{1}{N(1-s)}}^{1/2}\\
			\prec&~ \frac{1}{1-s_0}\qa{B_{s_0}(|x-y|)+\frac{1}{N(1-s)}}^{1/2}\lesssim \frac{1}{1-s}B_s^{1/2}(|x-y|).
		\end{aligned}
	\end{equation*}
	This proves \eqref{convolution_half_B2} and completes the proof of \Cref{convolution_B_half_B}.

	\subsection{Proof of \Cref{lemma_input_bound_subcritical}}\label{sec_proof_of_lemma_input_bound_subcritical}

	We first establish \eqref{zero_mode_removed_bound_subcritical} assuming \Cref{lemma_subcritical_S_mix}. For $\bsigma=(\sigma_1,\sigma_2)$ with $\sigma_1\neq\sigma_2$, we write
	\begin{align}
		&\Theta_{t,xy}^{\bsigma}-\frac{1}{N(1-t)}=\sum_{k=1}^\infty t^{k-1} \pa{S^{k}_{xy}-\frac{1}{N}}  \nonumber\\
		&\prec \sum_{k=1}^{k_{\txt{mix}}-1}\frac{1}{W}\pa{\frac{N}{W}}^{k\alpha}\pa{\frac{|x-y|}{W}+1}^{-1-k\alpha}+\frac{1}{N}\sum_{k=k_{\txt{mix}}}^{(\log N)^2}t^{k-1}+\sum_{k\ge (\log N)^2}e^{-c_\al k}. \label{eq:derv_Thetazero}
	\end{align}
	The first term is controlled by the case $k<k_{\mathrm{mix}}$ in \eqref{bound_S_mix}, while the second and third terms follow from the bounds for $k\ge k_{\mathrm{mix}}$ and from the exponential mixing estimate \eqref{S_mix_exponential}. Altogether, this yields
	\begin{align}
		&\Theta_{t,xy}^{\bsigma}-\frac{1}{N(1-t)} \prec \frac{1}{W}\pa{\frac{N}{W}}^{\alpha}\pa{\frac{|x-y|}{W}+1}^{-1-\alpha}+\frac1N \lesssim \cS_{xy} \asymp B_0(|x-y|),\nonumber
	\end{align}
	which proves \eqref{zero_mode_removed_bound_subcritical}. The bound \eqref{Theta_upper_bound_opposite_charge_subcritical} then follows immediately.
	We next consider the case $\sigma_1=\sigma_2$. Without loss of generality, we assume $\sigma_1=\sigma_2=+$. Writing $\xi:=m(\sE)^2$, we decompose
	\begin{equation*}
		\begin{aligned}
			\Theta_{t,xy}^{\bsigma}=\sum_{k=1}^{k_{\txt{mix}}-1}t^{k-1}(\xi S)^{k}_{xy}+\sum_{k=k_{\txt{mix}}}^{(\log N)^2}t^{k-1}(\xi S)^{k}_{xy}+\sum_{k\geq (\log N)^2}t^{k-1}(\xi S)^{k}_{xy} =:I_1+I_2+I_3.
		\end{aligned}
	\end{equation*}
	The terms $I_1$ and $I_2$ are bounded by $\cal S_{xy}$ using \eqref{bound_S_mix}, in the same manner as in \eqref{eq:derv_Thetazero}.
	For the third term, we use the exponential mixing estimate \eqref{S_mix_exponential} to obtain
	\begin{equation*}
		\begin{aligned}
			I_3=&\sum_{k\geq (\log N)^2}t^{k-1}(\xi S)^{k}_{xy}=\frac{1}{N}\sum_{k\geq (\log N)^2}t^{k-1}\xi^k +\OO\pa{\r e^{-c_{\alpha}(\log N)^2}\frac{1}{1-t}}        \lesssim \frac{1}{N} \lesssim \cal S_{xy}.
		\end{aligned}
	\end{equation*}
	Combining the bounds on $I_1$, $I_2$, and $I_3$ yields \eqref{Theta_upper_bound_equal_charge_subcritical}. It remains to prove \Cref{lemma_subcritical_S_mix}. 
	
	\begin{proof}[\bf Proof of \Cref{lemma_subcritical_S_mix}]
For $k=1$, the bound \eqref{bound_S_mix} follows directly from \eqref{alpha_decay}. Assume that \eqref{bound_S_mix} holds for all $1\le k\le k_0-1$ with $2\le k_0 < k_{\txt{mix}}$. Then, using the induction hypothesis and that $1+k_0\alpha>0$, which implies that
		\begin{equation*}
			\begin{aligned}
				S_{xy}^{k_0}&=\sum_{a\in\ZL}S^{k_0-1}_{xa}S_{ay}\lesssim \frac{1}{W^2}\pa{\frac{N}{W}}^{k_0\alpha}\sum_{a\in\ZL}\pa{\frac{|x-a|}{W}+1}^{-1-(k_0-1)\alpha}\pa{\frac{|a-y|}{W}+1}^{-1-\alpha}\\
				&\lesssim \frac{1}{W^2}\pa{\frac{N}{W}}^{k_0\alpha}\pa{\frac{|x-y|}{W}+1}^{-1-(k_0-1)\alpha}\sum_{|a-y|\leq |x-y|/2}\pa{\frac{|a-y|}{W}+1}^{-1-\alpha}\\
				&\quad  + \frac{1}{W^2}\pa{\frac{N}{W}}^{k_0\alpha}\pa{\frac{|x-y|}{W}+1}^{-1-\alpha}\sum_{|x-a|\le |x-y|/2}\pa{\frac{|x-a|}{W}+1}^{-1-(k_0-1)\alpha}\\
				&\quad + \frac{1}{W^2}\pa{\frac{N}{W}}^{k_0\alpha}\sum_{|x-a|\wedge |a-y|\ge |x-y|/2}\pa{\frac{|x-a|}{W}+1}^{-1-(k_0-1)\alpha}\pa{\frac{|a-y|}{W}+1}^{-1-\alpha}\\         &\lesssim \frac{1}{W}\pa{\frac{N}{W}}^{k_0\alpha} \pa{\frac{|x-y|}{W}+1}^{-1-k_0\alpha}.
			\end{aligned}
		\end{equation*}
		This establishes \eqref{bound_S_mix} for $k_0 < k_{\txt{mix}}$. For $k_0=k_{\txt{mix}}$ with $1+k_0\alpha\le 0$, we have 
		\begin{align}
			S_{xy}^{k_0}&=\sum_{a\in\ZL}S^{k_0-1}_{xa}S_{ay}\lesssim \frac{1}{W^2}\pa{\frac{N}{W}}^{k_0\alpha}\sum_{a\in\ZL}\pa{\frac{|x-a|}{W}+1}^{-1-(k_0-1)\alpha}\pa{\frac{|a-y|}{W}+1}^{-1-\alpha}\nonumber\\
			&\le \frac{C}{W}\pa{\frac{N}{W}}^{k_0\alpha}\cdot \pa{\frac{N}{W}}^{-1-k_0\alpha}\log \pa{\frac N W+1}=\frac{C}{N}\log \pa{\frac N W+1},\label{eq:Skmix}
		\end{align}
		which proves \eqref{bound_S_mix} at the mixing scale. The logarithmic factor appears in the borderline case $1+k_0\al=0$. For $k_0> k_{\txt{mix}}$, the desired bound follows immediately from \eqref{eq:Skmix}, along with the fact that $\|S\|_{\infty\to\infty}\le 1$. 
		To prove \eqref{S_mix_exponential}, we decompose $ S=\cQ S \cQ+\cal P$, where the projection operators $\cal P$ and $\cQ$ are defined in \Cref{def_sum_zero_operators}. Thus, for any $k\in \N$, we have
		\begin{equation*}
			\begin{aligned}
				S^{k}-\cal P&= (\cQ S \cQ)^k.
			\end{aligned}
		\end{equation*}
		The bound \eqref{S_mix_exponential} then follows directly from the spectral gap assumption in \Cref{assumption_input_bound<0}, which implies that the spectral radius of $\cQ S \cQ$ is strictly smaller than 1.
	\end{proof}

	\subsection{Proof of Lemma \ref{lemma_n_evolution_kernel_bound}}\label{sec_proof_of_n_evolution_kernel_bound}
	
	Let $\bsigma=(\sigma_1,\ldots,\sigma_n)\in\ha{-,+}^n$ with $n\geq 2$. We begin with the bound \eqref{evolution_kernel_bound_1}. By the definition \eqref{U_s_t_equal_1_Theta}, we have
	\begin{equation*}
		\begin{aligned}
			\absa{\cU_{s,t}^{(\sigma,\sigma')}(a,b)}=\absa{\delta_{ab}+(t-s)\Theta_{t,ab}^{(\sigma,\sigma')}}\leq \delta_{ab}+(t-s)\Theta_{t,ab}=\cU_{s,t}(a,b),
		\end{aligned}
	\end{equation*}
	Hence the first inequality in \eqref{evolution_kernel_bound_1} follows immediately:
	\begin{equation*}
		\begin{aligned}        \absa{\pa{\bU_{s,t,\bsigma}^{(n)}\circ\cA^{(n)}}_{\bx}}         &\leq \sum_{\ba=(a_1,\ldots,a_n)\in\ZL^n}\prod_{i=1}^{n}\cU_{s,t}(x_i,a_i)\cdot\absb{\cA_{\ba}^{(n)}}=\pa{\cU_{s,t}^{\otimes n}\circ |\cA|^{(n)}}_{\bx}.
		\end{aligned}
	\end{equation*}
	To prove the second inequality in \eqref{evolution_kernel_bound_1}, we expand
	\begin{equation}\label{Theta_I_decomposition_no_charge}
		\pa{\cU_{s,t}^{\otimes n}\circ |\cA|^{(n)}}_{\bx}=\sum_{I\subseteq \qq{n}}(t-s)^{|I|}\pa{\boldsymbol{\Theta}_{t}^{I,(n)}\circ\absa{\cA}^{(n)}}_{\bx},
	\end{equation}
	where for any $I\subseteq \qq{n}$ we define the operator $\boldsymbol{\Theta}_{t}^{I,(n)}$ by
	\begin{equation*}
		\begin{aligned}
			\pa{\boldsymbol{\Theta}_{t}^{I,(n)}\circ \cX^{(n)}}_{\bx}:=\sum_{\ba=(a_1,\ldots,a_n)\in\ZL^n}\prod_{i\in I}\Theta_{t,x_ia_i}\prod_{i\notin I}\delta_{x_ia_i}\cdot \cX^{(n)}_{\ba},\quad \forall \bx=(x_1,\ldots,x_n)\in\ZL^n.
		\end{aligned}
	\end{equation*}
	{We now bound the operator $\boldsymbol{\Theta}_{t}^{I,(n)}$ as}
	\begin{equation}\label{Theta_I_decomposition_bound}        \pa{\boldsymbol{\Theta}_{t}^{I,(n)}\circ\absa{\cA}^{(n)}}_{\bx}\prec \normb{\cA^{(n)}}_{\mathscr{D}_s^{-1}} \sum_{\ba=(a_1,\ldots,a_n)\in\ZL^n}\prod_{i\in I}B_t(|x_i-a_i|)\prod_{i\notin I}\delta_{x_ia_i}\cdot \mathscr{D}_{s,\ba}^{(n)}.
	\end{equation}
	For any $r_1,r_2\in\qa{s,t}$ with $r_1\leq r_2$, we have that
	\begin{equation}\label{triangular_bound_B_t_scrD_r_1_r_2}
		\begin{aligned}
			&~\sum_{a_i\in\ZL}B_t(|x_i-a_i|)\mathscr{D}_{r_1}(|a_i-a|)\mathscr{D}_{r_2}(|a_i-b|)\\
			=&~ \frac{1}{\ell_t(1-t)}\sum_{a_i\in\ZL}\mathscr{D}_t(|x_i-a_i|)\mathscr{D}_{r_1}(|a_i-a|)\cdot \mathscr{D}_t(|x_i-a_i|)\mathscr{D}_{r_2}(|a_i-b|)\\
			\lesssim&~\frac{\mathscr{D}_{t}(|x_i-b|)}{\ell_t(1-t)}\sum_{a_i\in\ZL}\mathscr{D}_t(|x_i-a_i|)\mathscr{D}_{r_1}(|a_i-a|)\prec \frac{\ell_{r_1}}{\ell_t(1-t)}\mathscr{D}_{t}(|x_i-a|)\mathscr{D}_{t}(|x_i-b|).
		\end{aligned}
	\end{equation}
	Applying this bound successively to the RHS of \eqref{Theta_I_decomposition_bound} for $a_1,\ldots,a_n$ (skipping indices $i\notin I$) yields 
	\begin{equation}\label{Theta_I_bound_no_charge}
		\begin{aligned}
			\pa{\boldsymbol{\Theta}_{t}^{I,(n)}\circ\absa{\cA}^{(n)}}_{\bx}\prec \pa{\frac{\ell_t}{\ell_s}}^{\mathbf{1}\pa{|I|=n}} \qa{\frac{\ell_s}{\ell_t(1-t)}}^{|I|}\normb{\cA^{(n)}}_{\mathscr{D}_s^{-1}}\cdot \mathscr{D}_{t,\bx}^{(n)}.
		\end{aligned}
	\end{equation}
	Substituting this estimate into \eqref{Theta_I_decomposition_no_charge} and summing over all $I\subseteq \qq{n}$ yields the bound \eqref{evolution_kernel_bound_1}.

	To improve the evolution kernel bound, we adopt a decomposition similar to that in \eqref{Theta_I_decomposition_no_charge}, i.e., we write 
	\begin{equation}\label{first_decomposition_cU_to_Theta_I}
		\bU_{s,t,\bsigma}^{(n)}=\sum_{I\subseteq \qq{n}}(t-s)^{|I|}\boldsymbol{\Theta}_{t,\bsigma}^{I,(n)},
	\end{equation}
	where the operator \smash{$\boldsymbol{\Theta}_{t,\bsigma}^{I,(n)}$} is defined for any $I\subseteq \qq{n}$ by
	\begin{equation*}
		\begin{aligned}
			\pa{\boldsymbol{\Theta}_{t,\bsigma}^{I,(n)}\circ \cX^{(n)}}_{\bx}=\sum_{\ba=(a_1,\ldots,a_n)\in\ZL^n}\prod_{i\in I}\Theta_{t,x_ia_i}^{(\sigma_i,\sigma_{i+1})}\prod_{i\notin I}\delta_{x_ia_i}\cdot \cX_{\ba}^{(n)},\quad \forall \bx=(x_1,\ldots,x_n)\in\ZL^n.
		\end{aligned}
	\end{equation*}
	Using the bound $\absb{\pb{\boldsymbol{\Theta}_{t,\bsigma}^{I,(n)}\circ \cX^{(n)}}_{\bx}}\le \pb{\boldsymbol{\Theta}_{t}^{I,(n)}\circ |\cX|^{(n)}}_{\bx}$ together with \eqref{Theta_I_bound_no_charge}, we obtain
	\begin{equation}\label{cU_cA_Theta_n_error}
		\begin{aligned}
			\pa{\bU_{s,t,\bsigma}^{(n)}\circ \cA^{(n)}}_{\bx}=(t-s)^n\pa{\boldsymbol{\Theta}_{t,\bsigma}^{\qq{n},(n)}\circ \cA^{(n)}}_{\bx}+\opr{\normb{\cA^{(n)}}_{\mathscr{D}_s^{-1}}\cdot\qa{\frac{\ell_s(1-s)}{\ell_t(1-t)}}^{n-1}\mathscr{D}_{t,\bx}^{(n)}}.
		\end{aligned}
	\end{equation}
	We now bound the first term in the case (I), where we assume that $\sigma_{i_0}=\sigma_{i_0+1}$ for some $i_0\in\qq{n}$. In this case, we apply the bound \eqref{Theta_bound_equaled_charge} to \smash{$\Theta_{t}^{(\sigma_{i_0},\sigma_{i_0+1})}(x_{i_0},a_{i_0})$}, and the bound \eqref{Theta_bound_opposite_charge} or \eqref{Theta_bound_equaled_charge} to \smash{$\Theta_{t}^{(\sigma_i,\sigma_{i+1})}(x_i,a_i)$} for $i\neq i_0$. 
	For the latter terms we use the bound \eqref{triangular_bound_B_t_scrD_r_1_r_2}, while for the former we apply
	\begin{equation}\label{triangular_bound_B_t_scrD_r_1_r_2_t0}
		\begin{aligned}
			&~\sum_{a_{i_0}\in\ZL}B_0(|x_{i_0}-a_{i_0}|)\mathscr{D}_{r_1}(|a_{i_0}-a|)\mathscr{D}_{r_2}(|a_{i_0}-b|)\\
			=&~ \frac{1}{W}\sum_{a_{i_0}\in\ZL}\mathscr{D}_0(|x_{i_0}-a_{i_0}|)\mathscr{D}_{r_1}(|a_{i_0}-a|)\cdot \mathscr{D}_0(|x_{i_0}-a_{i_0}|)\mathscr{D}_{r_2}(|a_{i_0}-b|)\\
			\lesssim&~\frac{\mathscr{D}_{r_2}(|x_{i_0}-b|)}{W}\sum_{a_{i_0}\in\ZL}\mathscr{D}_0(|x_{i_0}-a_{i_0}|)\mathscr{D}_{r_1}(|a_{i_0}-a|)\prec \mathscr{D}_{r_1}(|x_{i_0}-a|)\mathscr{D}_{r_2}(|x_{i_0}-b|).
		\end{aligned}
	\end{equation}
	Summing successively over $a_1,\ldots,a_n$ and applying \eqref{triangular_bound_B_t_scrD_r_1_r_2} and \eqref{triangular_bound_B_t_scrD_r_1_r_2_t0} accordingly at each step, we obtain
	\begin{equation*}
		\begin{aligned}
			\pa{\boldsymbol{\Theta}_{t,\bsigma}^{\qq{n},(n)}\circ \cA^{(n)}}_{\bx}\prec \norma{\cA}_{\mathscr{D}_s^{-1}}\cdot \qa{\frac{\ell_s}{\ell_t(1-t)}}^{n-1}\mathscr{D}_{t,\bx}^{(n)}.
		\end{aligned}
	\end{equation*}
	Combining this estimate with \eqref{cU_cA_Theta_n_error} yields the improved bound \eqref{evolution_kernel_bound_2}.
	
	It remains to establish the bound \eqref{evolution_kernel_bound_3} in case (II). In this case, the operator $\boldsymbol{\Theta}_{t,\bsigma}^{\qq{n},(n)}$ reduces to $\Theta_t^{\otimes n}$. Without loss of generality, by inserting an intermediate time $t_0=1-(W/N)^{\alpha\wedge 2}$ if necessary, it suffices to prove \eqref{evolution_kernel_bound_3} in the following two situations: (i) $s\leq t$ with $\ell_t<N$, or (ii) $\ell_s=\ell_t=N$. The second case is already implied by \eqref{evolution_kernel_bound_1}. We therefore focus on case (i). For any $I\subseteq \qq{n}$, define the $\vartheta$-operator by
	\begin{equation}\label{eq:vartheta_tI}
		\pa{\vartheta_{t,I}^{(n)}\circ \cX^{(n)}}_{\bx}:=\sum_{\ba=(a_1,\ldots,a_n)\in\ZL^n}\prod_{i\in I}\Theta_{t,x_i a_1}\cdot \prod_{i\notin I}\Theta_{t,x_i a_i}\cdot \cX_{\ba}^{(n)},\quad \forall \bx=(x_1,\ldots,x_n)\in\ZL^n.
	\end{equation}
	Note that  \smash{$\vartheta_{t,\qq{1}}^{(n)}=\boldsymbol{\Theta}_{t,\bsigma}^{\qq{n},(n)}$}. Using these identities and the sum-zero property \eqref{sum_zero_property_supercritical} of $\cA^{(n)}$, we obtain
	\begin{equation}\label{decomposition_bTheta_to_vartheta}
		\begin{aligned}
			\boldsymbol{\Theta}_{t,\bsigma}^{\qq{n},(n)}\circ \cA^{(n)}=\sum_{i=1}^{n-1}\pa{\vartheta_{t,\qq{i}}^{(n)}-\vartheta_{t,\qq{i+1}}^{(n)}}\circ \cA^{(n)}.
		\end{aligned}
	\end{equation}
	Therefore, in order to prove \eqref{evolution_kernel_bound_3}, it suffices to show that for each $i\in\qq{n-1}$,
	\begin{equation}\label{eq:reduceto_Theta_diff}
		\begin{aligned}
			\qa{\pa{\vartheta_{t,\qq{i}}^{(n)}-\vartheta_{t,\qq{i+1}}^{(n)}}\circ \cA^{(n)}}_{\bx}\prec \normb{\cA^{(n)}}_{\mathscr{D}_s^{-1}}\cdot \qa{\frac{\ell_s}{\ell_t(1-t)}}^n\mathscr{D}_{t,\bx}^{(n)},\qquad \forall \bx\in\ZL^n.
		\end{aligned}
	\end{equation}
	Without loss of generality, assume that $\norm{\cA^{(n)}}_{\mathscr{D}_s^{-1}}\le 1$. Using the definition of the $\vartheta$-operators together with the assumptions \eqref{Theta_bound_opposite_charge} and \eqref{Theta_difference_bound}, we obtain for $\bx=(x_1,\ldots,x_n)\in\ZL^n$,
	\begin{align}
		&\qa{\pa{\vartheta_{t,\qq{i}}^{(n)}-\vartheta_{t,\qq{i+1}}^{(n)}}\circ \cA^{(n)}}_{\bx}=\sum_{\ba=(a_1,\ldots,a_n)\in\ZL^n}\prod_{1\leq j\leq i}\Theta_{t,x_j a_1}\cdot\pa{\Theta_{t,x_{i+1}a_{i+1}}-\Theta_{t,x_{i+1}a_1}}\prod_{i+1< j\leq n}\Theta_{t,x_j a_j}\cdot \cA_{\ba}^{(n)}\nonumber\\
		&\prec \sum_{\ba}\braket{a_{i+1}-a_1}_WR_t(|x_{i+1}-a_{i+1}|)\prod_{1\leq j\leq i+1}B_{t}(|x_j-a_1|)\prod_{i+1< j\leq n}B_{t}(|x_j-a_j|)\cdot \mathscr{D}_{s,\ba}^{(n)}\nonumber\\
		&\quad +\sum_{\ba}\braket{a_{i+1}-a_1}_WR_t(|x_{i+1}-a_{1}|)\prod_{1\leq j\leq i}B_{t}(|x_j-a_1|)\prod_{i+1\leq j\leq n}B_{t}(|x_j-a_j|)\cdot \mathscr{D}_{s,\ba}^{(n)}\nonumber\\
		&\prec \frac{\ell_s^{i-1}}{\q{\ell_t(1-t)}^{n}}\sum_{a_1,a_{i+1},\ldots,a_n}\braket{a_{i+1}-a_1}_WR_t(|x_{i+1}-a_{i+1}|)\mathscr{D}_{t,(x_1,\ldots,x_{i+1})}^{(i+2)}(a_1)\prod_{i+1< j\leq n}\mathscr{D}_{t}^2(|x_j-a_j|)\cdot \mathscr{D}_{s,(a_1,a_{i+1},\ldots,a_n)}^{(n-i+1)}\nonumber\\
		&\quad +\frac{\ell_s^{i-1}}{\q{\ell_t(1-t)}^{n}} \sum_{a_1,a_{i+1},\ldots,a_n}\braket{a_{i+1}-a_1}_WR_t(|x_{i+1}-a_{1}|)\mathscr{D}_{t,(x_1,\ldots,x_i)}^{(i+1)}(a_1)\prod_{i+1\leq j\leq n}\mathscr{D}_{t}^2(|x_j-a_j|)\cdot \mathscr{D}_{s,(a_1,a_{i+1},\ldots,a_n)}^{(n-i+1)}\nonumber\\
		&=:\frac{\ell_s^{i-1}}{\q{\ell_t(1-t)}^n}\pa{I_1+I_2}.\label{eq:def_loopsumzero_I12}
	\end{align}
    In the third step, we use inequalities like \eqref{example_triangular_2_edge_to_1} to reduce the loop \smash{$\mathscr{D}_{s,\ba}^{(n)}$} to \smash{$\mathscr{D}_{s,(a_1,a_{i+1},\ldots,a_n)}^{(n-i+1)}$},  producing an additional factor $\ell_s^{i-1}$. Furthermore, we obtain the loop factors {$\mathscr{D}_{t,(x_1,\ldots,x_{i+1})}^{(i+2)}(a_1)$} and \smash{$\mathscr{D}_{t,(x_1,\ldots,x_i)}^{(i+1)}(a_1)$} by applying the following inequality:
		\begin{equation}\label{bound_one_point_to_loop}
			\begin{aligned}
				\prod_{i=1}^{n}\mathscr{D}_{u}^2(|x_i-a|)\lesssim \mathscr{D}_{u,\bx}^{\pa{n+1}}(a),\quad \bx=(x_1,\ldots, x_n).
			\end{aligned}
		\end{equation}
    This inequality follows easily by induction from the elementary estimate $\mathscr{D}_{u}(|x-a|)\mathscr{D}_{u}(|y-a|)\lesssim \mathscr{D}_{u}(|x-y|)$. 
To conclude \eqref{eq:reduceto_Theta_diff}, it therefore suffices to prove
	\begin{equation}\label{eq:reduceto_Theta_diffi12}
		I_1+I_2\prec \ell_s^{n-i+1}\cdot\mathscr{D}_{t,\bx}^{(n)}.
	\end{equation}
	Note that the factor $\prod_{j=1}^{i}\mathscr{D}_t(|x_{j+1}-x_j|)$, which contains the variables $x_2,\ldots,x_i$, appears on both sides of this inequality and may therefore be cancelled. Consequently, without loss of generality (and up to minor notational changes), it suffices to prove \eqref{eq:reduceto_Theta_diffi12} in the case $i=1$. In this case the bound reduces to
	\begin{align}
		&\sum_{\ba=(a_1,\ldots,a_n)\in\ZL^n}\braket{a_{2}-a_1}_WR_t(|x_{2}-a_{2}|)\mathscr{D}_{t}(|x_1-a_1|)\mathscr{D}_{t}(|x_1-x_2|)\mathscr{D}_{t}(|x_2-a_1|)\prod_{2< j\leq n}\mathscr{D}_{t}^2(|x_j-a_j|)\cdot \mathscr{D}_{s,\ba}^{(n)}\nonumber\\
		&+\sum_{\ba=(a_1,\ldots,a_n)\in\ZL^n}\braket{a_{2}-a_1}_WR_t(|x_{2}-a_{1}|)\prod_{1\leq j\leq n}\mathscr{D}_{t}^2(|x_j-a_j|)\cdot \mathscr{D}_{s,\ba}^{(n)}\prec \ell_s^{n}\cdot\mathscr{D}_{t,\bx}^{(n)}.\label{target_after_difference}
	\end{align}

	We present the proof of \eqref{target_after_difference} only for the case $\alpha\in[1,2]$; the case $\alpha>2$ can be handled by a similar but slightly simpler argument. We first consider the first term on the LHS of \eqref{target_after_difference}. When $n=2$, we have
	\begin{align*}
		&~\sum_{a_1,a_2}\braket{a_{2}-a_1}_WR_t(|x_{2}-a_{2}|)\mathscr{D}_{t}(|x_1-a_1|)\mathscr{D}_{t}(|x_1-x_2|)\mathscr{D}_{t}(|x_2-a_1|)\cdot \mathscr{D}_{s}^2(|a_1-a_2|)\\
		\lesssim &~\pa{\frac{W}{\ell_t}}^{\alpha-1}\frac{\ell_s}{W}\mathscr{D}_{t}^2(|x_1-x_2|)\sum_{a_1,a_2}\pa{\frac{|x_2-a_2|}{W}+1}^{\alpha-2}\pa{\frac{|x_2-a_2|}{\ell_t}+1}^{1-\alpha}\pa{\frac{|a_1-a_2|}{\ell_s}+1}^{-\alpha}\\
		\prec &~\pa{\frac{W}{\ell_t}}^{\alpha-1}\frac{\ell_s^2}{W}\mathscr{D}_{t}^2(|x_1-x_2|)\sum_{a_2}\pa{\frac{|x_2-a_2|}{W}+1}^{\alpha-2}\pa{\frac{|x_2-a_2|}{\ell_t}+1}^{1-\alpha}\prec \ell_s^2\mathscr{D}_t^2(|x_1-x_2|).
	\end{align*}
	Next consider the case $n\ge 3$. We first sum over the index $a_1$ and use the bound
	\begin{align*}
		&~\sum_{a_1}\braket{a_{2}-a_1}_W\mathscr{D}_{t}(|x_1-a_1|)\mathscr{D}_{t}(|x_2-a_1|) \mathscr{D}_{s}(|a_1-a_2|)\mathscr{D}_{s}(|a_1-a_n|)\\
		\lesssim &~\ell_s \sum_{a_1}\mathscr{D}_{t}(|x_1-a_1|) \mathscr{D}_{s}(|a_1-a_n|)\pa{\frac{|x_2-a_1|}{\ell_t}+1}^{-(1+\alpha)/2}\pa{\frac{|a_1-a_2|}{\ell_s}+1}^{(1-\alpha)/2}\\
		\lesssim &~\ell_s \sum_{a_1}\mathscr{D}_{t}(|x_1-a_1|) \mathscr{D}_{s}(|a_1-a_n|)\pa{\frac{|x_2-a_2|}{\ell_t}+1}^{(1-\alpha)/2}\prec \ell_s^2\mathscr{D}_{t}(|x_1-a_n|)\pa{\frac{|x_2-a_2|}{\ell_t}+1}^{(1-\alpha)/2}.
	\end{align*}
	After this step, we apply inequalities such as \eqref{triangular_bound_B_t_scrD_r_1_r_2} successively for the indices $a_n,a_{n-1},\ldots,a_3$, which yields
	\begin{align}
		&\nonumber \sum_{\ba=(a_1,\ldots,a_n)\in\ZL^n}\braket{a_{2}-a_1}_WR_t(|x_{2}-a_{2}|)\mathscr{D}_{t}(|x_1-a_1|)\mathscr{D}_{t}(|x_1-x_2|)\mathscr{D}_{t}(|x_2-a_1|)\prod_{2< j\leq n}\mathscr{D}_{t}^2(|x_j-a_j|)\cdot \mathscr{D}_{s,\ba}^{(n)}\\
		&\prec \ell_s^{n}\sum_{a_2}R_t(|x_2-a_2|)\pa{\frac{|x_2-a_2|}{\ell_t}+1}^{(1-\alpha)/2}\mathscr{D}_{t,(x_3,\ldots,x_n,x_1)}^{(n)}(a_2,x_2)\prec \ell_s^{n}\mathscr{D}_{t,\bx}^{(n)},\label{target_after_difference1}
	\end{align}
	where in the second step we used the bound
	\begin{equation*}
		\begin{aligned}
			\sum_{a_2}R_t(|x_2-a_2|)\pa{{|x_2-a_2|}/{\ell_t}+1}^{(1-\alpha)/2}\mathscr{D}_{t}(|x_3-a_2|)\prec \mathscr{D}_t(|x_2-x_3|),
		\end{aligned}
	\end{equation*}
	which follows from an argument similar to those used in
	\eqref{triangular_convolution_example_1} and \eqref{triangular_convolution_example_2}.

We now consider the second term on the LHS of \eqref{target_after_difference}. When $n=2$, we have
	\begin{align*}
		&~\sum_{a_1,a_2}\braket{a_{2}-a_1}_WR_t(|x_{2}-a_{1}|)\mathscr{D}_{t}^2(|x_1-a_1|)\mathscr{D}_{t}^2(|x_2-a_2|) \mathscr{D}_{s}^2(|a_1-a_2|)\\
		\lesssim &~\ell_s\sum_{a_1,a_2}R_t(|x_{2}-a_{1}|)\mathscr{D}_{t}^2(|x_1-a_1|)\mathscr{D}_{t}^2(|x_2-a_2|) \pa{\frac{|a_1-a_2|}{\ell_s}+1}^{-\alpha}\\
		\prec &~\ell_s^2\sum_{a_1}R_t(|x_{2}-a_{1}|)\mathscr{D}_{t}^2(|x_1-a_1|) \pa{\frac{|a_1-x_2|}{\ell_t}+1}^{-\alpha}\prec \ell_s^2\mathscr{D}_{t}^2(|x_1-x_2|).
	\end{align*}
	The second and third steps again follow from arguments similar to those used in
	\eqref{triangular_convolution_example_1} and \eqref{triangular_convolution_example_2}.
	When $n\geq 3$, we first sum over the index $a_2$ and use
	\begin{align*}
		&~\sum_{a_2}\braket{a_2-a_1}_W\mathscr{D}_{t}^2(|x_2-a_2|)\mathscr{D}_{s}(|a_2-a_1|)\mathscr{D}_{s}(|a_2-a_3|)\\
		\lesssim &~\ell_s \sum_{a_2}\mathscr{D}_{t}^2(|x_2-a_2|)\pa{\frac{|a_2-a_1|}{\ell_s}+1}^{(1-\alpha)/2}\mathscr{D}_{s}(|a_2-a_3|)\\
		\lesssim &~\ell_s \pa{\frac{|x_2-a_1|}{\ell_t}+1}^{(1-\alpha)/2}\sum_{a_2}\mathscr{D}_{t}(|x_2-a_2|)\mathscr{D}_{s}(|a_2-a_3|)\prec \ell_s^2 \pa{\frac{|x_2-a_1|}{\ell_t}+1}^{(1-\alpha)/2}\mathscr{D}_t(|x_2-a_3|).
	\end{align*}
	Applying inequalities like \eqref{triangular_bound_B_t_scrD_r_1_r_2} successively for the indices $a_3,\ldots,a_n$, we obtain
	\begin{align}
		&\sum_{\ba=(a_1,\ldots,a_n)\in\ZL^n}\braket{a_{2}-a_1}_WR_t(|x_{2}-a_{1}|)\prod_{1\leq j\leq n}\mathscr{D}_{t}^2(|x_j-a_j|)\cdot \mathscr{D}_{s,\ba}^{(n)}\nonumber\\
		&\prec \ell_s^n\sum_{a_1}\mathscr{D}_t^2(|x_1-a_1|)R_t(|x_2-a_1|)\pa{\frac{|x_2-a_1|}{\ell_t}+1}^{(1-\alpha)/2}\mathscr{D}_{t,(x_3,\ldots,x_n)}^{(n-1)}(x_2,a_1)\nonumber\\
		&\lesssim \ell_s^n\sum_{a_1}\mathscr{D}_t(|x_1-a_1|)R_t(|x_2-a_1|)\pa{\frac{|x_2-a_1|}{\ell_t}+1}^{(1-\alpha)/2}\mathscr{D}_{t,(x_3,\ldots,x_n)}^{(n-1)}(x_2,x_1)\prec \ell_s^n\mathscr{D}_{t,\bx}^{(n)},\label{target_after_difference2}
	\end{align}
	where the last step again follows from an argument similar to that used for \eqref{triangular_convolution_example_1} and \eqref{triangular_convolution_example_2}.
	Combining \eqref{target_after_difference1} and \eqref{target_after_difference2}, we obtain \eqref{target_after_difference}, which completes the proof of \eqref{evolution_kernel_bound_3}.

	\section{Relations between \texorpdfstring{$\Xi$}{Xi}-parameters}\label{sec_Xi_inequalities}

	\subsection{Proof of Lemma \ref{lemma_Xi_reshape}}\label{sec_proof_of_Xi_bounding_relations}

The bounds \eqref{eq_cC_d_to_loops} and \eqref{eq_cC_od_to_cC_d} were essentially proved in the proof of \cite[Lemma A.2]{Band1D}, where only the $\max$-norm was tracked. Although the extension is conceptually straightforward, keeping track of all decay factors turns out to be technically involved. For completeness, we present the details here for the reader's convenience. It suffices to prove the estimates \eqref{eq_cC_d_to_loops}--\eqref{reshape_chains_into_loops} at a fixed time $u\in\qa{s,t}$; a standard $N^{-C}$-net argument then extends them uniformly to all $u\in[s,t]$.

For $n=2$, Lemma \ref{2_loop_to_1_chain} already yields 
\begin{equation}
	\Xi_{u,2}^{\cC_{\txt{d}}}\prec \Xi_{u,2}^{\cL}.\label{eq:Xi2_Cd}
\end{equation}
Next, assume that $n\geq 3$. To bound the chains by loops, we employ resolvent identities together with some standard large-deviation estimates; see, for example, \cite[Lemma 3.3]{ERDOS20121435}. For any $i\in\ZL$, let \smash{$H_u^{(i)}$} denote the $(N-1)\times(N-1)$ minor of $H_u$ obtained by removing the $i$-th row and column. The corresponding resolvent minor is the $(N-1)\times (N-1)$ matrix \smash{\(G_{u}^{\pa{i}}\p{\sigma}:=\p{H^{\pa{i}}_u-z_u\p{\sigma}}^{-1}.\)}
	For notational convenience, we regard \smash{$G_u^{(i)}(\sigma)$} as an $N\times N$ matrix by setting \smash{$G_{u,kl}^{(i)}(\sigma)=0$} whenever $k=i$ or $l=i$. The following resolvent identities follow directly from the Schur complement formula: for $i,j,k\in\ZL$,
	\begin{align}\label{resolvent_identities}        &G_{u,ij}=G_{u,ii}\pB{\delta_{ij}-\sum_{k}H_{u,ik}G_{u,kj}^{\pa{i}}},\quad  &G_{u,jk}^{\pa{i}}=G_{u,jk}-\frac{G_{u,ji}G_{u,ik}}{G_{u,ii}}=\pa{G_u-\frac{G_u\Delta^{i}G_u}{G_{u,ii}}}_{jk},
	\end{align}
	where the matrix $\Delta^{i}$ is defined by $\Delta^{i}_{xy}:=\delta_{xi}\delta_{yi}$ for $x,y\in \ZL$.
	With these notations, we further introduce the following variants of the $n$-$\cC$-chain \smash{$\cC_{u,\bsigma,\bx}^{(n)}(a,b)$} for any position $i\in\ZL$: 
	\begin{align*}
		\cC_{u,\bsigma,\bx}^{\pa{n},(i]}(a,b):=&\pb{G_u\p{\sigma_1}S^{\pa{x_1}}G_u^{\pa{i}}\p{\sigma_2}S^{\pa{x_2}}\cdots G_u^{\pa{i}}\p{\sigma_{n-1}} S^{\pa{x_{n-1}}}G_u^{\pa{i}}\p{\sigma_n}}_{ab},\\
		\cC_{u,\bsigma,\bx}^{\pa{n},[i)}(a,b):=&\pb{G_u^{\pa{i}}\p{\sigma_1}S^{\pa{x_1}}G_u^{\pa{i}}\p{\sigma_2}S^{\pa{x_2}}\cdots G_u^{\pa{i}}\p{\sigma_{n-1}} S^{\pa{x_{n-1}}}G_u\p{\sigma_n}}_{ab},\\
		\cC_{u,\bsigma,\bx}^{\pa{n},\qa{i}}(a,b):=&\pb{G_u^{\pa{i}}\p{\sigma_1}S^{\pa{x_1}}G_u^{\pa{i}}\p{\sigma_2}S^{\pa{x_2}}\cdots G_u^{\pa{i}}\p{\sigma_{n-1}}S^{\pa{x_{n-1}}}G_u^{\pa{i}}\p{\sigma_n}}_{ab},\\
		\cC_{u,\bsigma,\bx}^{\pa{n},\pa{i}}(a,b):=&\pb{G_u\p{\sigma_1}S^{\pa{x_1}}G_u^{\pa{i}}\p{\sigma_2}S^{\pa{x_2}}\cdots G_u^{\pa{i}}\p{\sigma_{n-1}}S^{\pa{x_{n-1}}}G_u\p{\sigma_n}}_{ab}.
	\end{align*}

	We first consider the diagonal case with $n=2k$. For any diagonal $n$-chain $\cC_{u,\bsigma,\bx}^{\pa{n}}(x)$ with $\bsigma=(\sigma_1,\ldots,\sigma_n)\in\ha{-,+}^n$ and $(\bx,x)=(x_1,\ldots,x_{n-1},x)\in\ZL^n$, we approximate it by 
	\begin{equation*}
		\begin{aligned}            &\cC_{u,\bsigma,\bx}^{\pa{n},\pa{x}}(x,x)=G_{u,xx}\p{\sigma_1}G_{u,xx}\p{\sigma_n}\sum_{a,b\in\ZL}H_{u,xa}\cC_{u,\bsigma,\bx}^{\pa{n},\qa{x}}\p{a,b}H_{u,bx},
		\end{aligned}
	\end{equation*}
	where we used the first resolvent identity in \eqref{resolvent_identities}. 
	Since the entries in the $x$-th row and column of $H_u$ are independent of \smash{$G_u^{\pa{x}}$}, we may apply the large deviation estimate from \cite[Lemma 3.3]{ERDOS20121435} to obtain 
	\begin{equation}\label{concentration_x_to_xx}
		\begin{aligned}                    \cC_{u,\bsigma,\bx}^{\pa{n},\pa{x}}(x,x)=&~uG_{u,xx}\p{\sigma_1}G_{u,xx}\p{\sigma_n}\sum_{a\in\ZL}\cC_{u,\bsigma,\bx}^{\pa{n},\qa{x}}\p{a,a}S_{ax}\\
			&~+\OO_\prec \pbb{\sum_{a,b\in\ZL}S_{xa}S_{xb} \cC_{u,\bsigma,\bx}^{\pa{n},\qa{x}}\pa{a,b}\pa{\cC_{u,\bsigma,\bx}^{\pa{n},\qa{x}}}^*\p{b,a}}^{1/2}.
		\end{aligned}
	\end{equation}
	Note that the leading term contains an $n$-loop in which all resolvents $G_u$ are replaced by $G_u^{(x)}$, while the error term corresponds to a $(2n)$-loop with all resolvents replaced by \smash{$G_u^{(x)}$}. Using the inequality $\mathrm{Tr}(AB)\le \mathrm{Tr}(A)\mathrm{Tr}(B)$ for any positive semidefinite matrices $A$ and $B$, we can bound the $(2n)$-loop in the error term by the product of two $n$-loops formed with \smash{$G_u^{(x)}$}.
	Therefore, to complete the proof, it remains to bound the $n$-loops formed with \smash{$G_u^{(x)}$}, as well as the difference \smash{$\cC_{u,\bsigma,\bx}^{\pa{n}}(x,x)-\cC_{u,\bsigma,\bx}^{\pa{n},\pa{x}}(x,x)$}.

	We first estimate the loop $\sum_{a}\cC_{u,\bsigma,\bx}^{\pa{n},\qa{x}}S_{ax}=\avgb{G_u^{\pa{x}}\p{\sigma_1}S^{\pa{x_1}}\cdots G_u^{\pa{x}}\p{\sigma_n}S^{\pa{x_n}}}$. Applying the second identity in \eqref{resolvent_identities} to the entries \smash{$G_u^{\pa{x}}$} in the loop, we can express this quantity as an $n$-$\cL$-loop together with several subleading terms. In these terms, the loop is decomposed into several $\cC$-chains, say $\ell$ of them (depending on how many \smash{$G_u^{\pa{x}}$} are replaced by $G_u\Delta^xG_u$), with lengths $n_1,\ldots, n_\ell\ge 2$, together with some powers of $1/G_{u,xx}$. Applying the entrywise local law in \eqref{step_1_supercritical_weak_local_law} to the factors $G_{u,xx}$ in the denominator, using the definitions in \eqref{eq:Xi_Cd_Cod}, and employing the inequality $\mathscr{D}_u(|a-x|)\mathscr{D}_u(|x-b|)\lesssim \mathscr{D}_u(|a-b|)$ to recover \smash{$\mathscr{D}_{u,(\bx,x)}^{\pa{n}}$}, we obtain
	\begin{align}            &\sum_{a\in\ZL}\cC_{u,\bsigma,\bx}^{\pa{n},\qa{x}}(a,a)S_{ax}=\cL_{u,\bsigma,(\bx,x)}^{\pa{n}}+B_u^{n}(0)\mathscr{D}_{u,(\bx,x)}^{\pa{n}}\cdot \OO_\prec \pbb{\sum_{\ell\ge 1}\sum_{n_1,\ldots,n_\ell\geq 2: n_1+\cdots +n_\ell=n+\ell}\Xi_{u,n_1}^{\cC_{\rd}}\cdots \Xi_{u,n_\ell}^{\cC_{\rd}}}\nonumber\\
		&=\cL_{u,\bsigma,(\bx,x)}^{\pa{n}}+B_u^{n}(0)\mathscr{D}_{u,(\bx,x)}^{\pa{n}}\cdot \OO_\prec\pbb{\sum_{\ell\ge 2}\sum_{2\leq n_1,\ldots,n_\ell\leq n-2: n_1+\cdots +n_\ell=n+\ell}\Xi_{u,n_1}^{\cC_{\rd}}\cdots \Xi_{u,n_\ell}^{\cC_{\rd}}}\nonumber\\
		&\quad +B_u^{n}(0)\mathscr{D}_{u,\pa{\bx,x}}^{\pa{n}}\cdot\OO_\prec\pbb{\Xi_{u,n+1}^{\cC_{\rd}}+\Xi_{u,2}^{\cC_{\rd}}\Xi_{u,n}^{\cC_{\rd}}+\Xi_{u,3}^{\cC_{\rd}}\Xi_{u,n-1}^{\cC_{\rd}}+\pb{\Xi_{u,2}^{\cC_{\rd}}}^2\Xi_{u,n-1}^{\cC_{\rd}}}\nonumber\\
		&=\cL_{u,\bsigma,(\bx,x)}^{\pa{n}}+B_u^{n}(0)\mathscr{D}_{u,(\bx,x)}^{\pa{n}}\cdot \OO_\prec\pbb{\sum_{\ell\ge 2}\sum_{2\leq n_1,\ldots,n_\ell\leq n-2: n_1+\cdots +n_\ell=n+\ell}\Xi_{u,n_1}^{\cC_{\rd}}\cdots \Xi_{u,n_\ell}^{\cC_{\rd}}+B_u^{-1/2}(0)\pa{\Xi_{u,2}^{\cL}}^{1/2}\Xi_{u,n}^{\cC_{\rd}}}\nonumber\\
		&\quad +B_u^{n}(0)\mathscr{D}_{u,\pa{\bx,x}}^{\pa{n}}\cdot\opr{\Xi_{u,2}^{\cC_{\rd}}\Xi_{u,n}^{\cC_{\rd}}+ \Xi_{u,4}^{\cC_{\rd}}\Xi_{u,n-2}^{\cC_{\rd}} +\pb{\Xi_{u,2}^{\cC_{\rd}}}^{3}\Xi_{u,n-2}^{\cC_{\rd}}}.\label{eq:n-loop-x}
	\end{align}
	In the second step, we singled out all terms containing a chain of length at least $n-1$. In the third step, we applied the bounds in \eqref{loop_bound_cut_long_to_short} to obtain \smash{$\Xi_{u,n+1}^{\cC_{\rd}}\lesssim B_u^{-1/2}(0)\pa{\Xi_{u,2}^{\cL}}^{1/2}\Xi_{u,n}^{\cC_{\rd}}$}, as well as 
	\begin{align*}
		\Xi_{u,3}^{\cC_{\rd}}\Xi_{u,n-1}^{\cC_{\rd}}&\lesssim \pa{\Xi_{u,2}^{\cC_{\rd}}\Xi_{u,4}^{\cC_{\rd}}}^{1/2}\pa{\Xi_{u,n-2}^{\cC_{\rd}}\Xi_{u,n}^{\cC_{\rd}}}^{1/2}\le \Xi_{u,4}^{\cC_{\rd}}\Xi_{u,n-2}^{\cC_{\rd}}+\Xi_{u,2}^{\cC_{\rd}}\Xi_{u,n}^{\cC_{\rd}},\\
		\pa{\Xi_{u,2}^{\cC_{\rd}}}^2\Xi_{u,n-1}^{\cC_{\rd}}&\lesssim \pa{\Xi_{u,2}^{\cC_{\rd}}}^2\pa{\Xi_{u,n-2}^{\cC_{\rd}}\Xi_{u,n}^{\cC_{\rd}}}^{1/2} \lesssim\Xi_{u,2}^{\cC_{\rd}}\Xi_{u,n}^{\cC_{\rd}}+\pa{\Xi_{u,2}^{\cC_{\rd}}}^3\Xi_{u,n-2}^{\cC_{\rd}}.\end{align*}
	As explained above, the second term on the RHS of \eqref{concentration_x_to_xx} can be bounded in a similar manner. Consequently, we obtain
	\begin{align}
		\pb{B_u^{n-1}(0)\mathscr{D}_{u,(\bx,x)}^{\pa{n}}}^{-1} \cC_{u,\bsigma,\bx}^{\pa{n},\pa{x}}(x,x) \prec &~\Xi_{u,n}^{\cL} + B_u(0){\Xi_{u,2}^{\cC_{\rd}}\Xi_{u,n}^{\cC_{\rd}}}+B_u^{1/2}(0)\pa{\Xi_{u,2}^{\cL}}^{1/2}\Xi_{u,n}^{\cC_{\rd}}\nonumber\\
		&  +B_u(0)\sum_{\ell\ge 2}\sum_{2\leq n_1,\ldots,n_\ell\leq n-2: n_1+\cdots +n_\ell=n+\ell}\Xi_{u,n_1}^{\cC_{\rd}}\cdots \Xi_{u,n_\ell}^{\cC_{\rd}}.\label{chain_to_loop_x_chain_bound}
	\end{align}
    
	Next, we control the difference $\cC_{u,\bsigma,\bx}^{\pa{n}}(x,x)-\cC_{u,\bsigma,\bx}^{\pa{n},\pa{x}}(x,x)$.
	Proceeding as above in \eqref{eq:n-loop-x}, we apply the second
	identity in \eqref{resolvent_identities} to replace each occurrence of \smash{$G_u^{\pa{x}}$} by $G_u$ in \smash{$\cC_{u,\bsigma,\bx}^{\pa{n},\pa{x}}(x,x)$}. This produces several terms in which the $n$-chain is decomposed into $\ell$ diagonal $\cC$-chains of lengths $n_1,\ldots, n_\ell\ge 2$, whenever $(\ell-1)$ factors \smash{$G_u^{\pa{x}}$} are replaced by $G_u\Delta^xG_u$. Applying the entrywise local law in \eqref{step_1_supercritical_weak_local_law} to the factors $G_{u,xx}$ in the denominator and using the definitions in \eqref{eq:Xi_Cd_Cod}, we obtain that for any small constant $\e>0$,
	\begin{align}
		\frac{\cC_{u,\bsigma,\bx}^{\pa{n}}(x,x)-\cC_{u,\bsigma,\bx}^{\pa{n},\pa{x}}(x,x)}{B_u^{n-1}(0)\mathscr{D}_{u,(\bx,x)}^{\pa{n}}}  \prec \sum_{\ell\ge 2}\sum_{2\le n_1,\ldots,n_\ell\le n-2:n_1+\cdots+n_\ell=n+\ell-1}\Xi_{u,n_1}^{\cC_{\rd}}\cdots \Xi_{u,n_\ell}^{\cC_{\rd}}+ \Xi_{u,2}^{\cC_{\rd}} \Xi_{u,n-1}^{\cC_{\rd}} \nonumber\\
		\lesssim \sum_{\ell\ge 2}\sum_{2\le n_1,\ldots,n_\ell\le n- 2:n_1+\cdots+n_\ell=n+\ell-1}\Xi_{u,n_1}^{\cC_{\rd}}\cdots \Xi_{u,n_\ell}^{\cC_{\rd}} + B_u^\e(0)\Xi_{u,n}^{\cC_{\rd}} + B_u^{-\e}(0) \pb{\Xi_{u,2}^{\cC_{\rd}}}^2 \Xi_{u,n-2}^{\cC_{\rd}},\label{eq:n-loop-x_diff}
	\end{align}
	where in the second step, we used \smash{$\Xi_{u,n-1}^{\cC_{\rd}}\lesssim \p{\Xi_{u,n-2}^{\cC_{\rd}} \Xi_{u,n}^{\cC_{\rd}}}^{1/2}$}, which follows from \eqref{loop_bound_cut_long_to_short}, together with the AM-GM inequality. Combining this estimate with \eqref{chain_to_loop_x_chain_bound} yields a self-consistent bound for $\Xi_{u,n}^{\cC_{\rd}}$. If \smash{$\Xi_{u,2}^{\cL}\le B_t^{-c}(0)$} for some small constant $c\in(0,1/100)$, then the self-consistent bound implies
	\begin{align}         \Xi_{u,n}^{\cC_{\rd}}\prec\Xi_{u,n}^{\cL}         & +B_u(0)\sum_{\ell\ge 2}\sum_{2\leq n_1,\ldots,n_\ell\leq n-2: n_1+\cdots +n_\ell=n+\ell}\Xi_{u,n_1}^{\cC_{\rd}}\cdots \Xi_{u,n_\ell}^{\cC_{\rd}} \nonumber\\
		&+\sum_{\ell\ge 2}\sum_{2\le n_1,\ldots,n_\ell\le n- 2:n_1+\cdots+n_\ell=n+\ell-1}\Xi_{u,n_1}^{\cC_{\rd}}\cdots \Xi_{u,n_\ell}^{\cC_{\rd}} .\label{cC_to_loops_fixed_time_bound}        \end{align}
	Here we used that $\e$ can be chosen arbitrarily small and that \smash{$\Xi_{u,2}^{\cC_{\rd}}\prec \Xi_{u,2}^{\cL}\le B_t^{-c}(0)$} by \eqref{eq:Xi2_Cd}. Finally, to obtain the stopping-time bound \eqref{eq_cC_d_to_loops} uniformly for $u\in[s,t]$ from the fixed-time bound above, we apply a standard $N^{-C}$-net argument together with a perturbation estimate.

	For \eqref{eq_cC_od_to_cC_d}, consider an off-diagonal $n$-chain \smash{$\cC_{u,\bsigma,\bx}^{\pa{n}}(x,y)$} with $n\geq 2$. We approximate it by \smash{$\cC_{u,\bsigma,\bx}^{\pa{n},(x]}(x,y)$} and control the difference between the two quantities. First, we apply the second identity in \eqref{resolvent_identities} to the entries \smash{$G_u^{\pa{x}}$} in \smash{$\cC_{u,\bsigma,\bx}^{\pa{n},(x]}(x,y)$}. This expresses the chain as an $n$-$\cC$-chain together with several subleading terms.
	In these terms, the $n$-chain is decomposed into an off-diagonal $\cC$-chain of length $n_1\ge 1$ and $(\ell-1)$ diagonal $\cC$-chains of lengths $n_2,\ldots, n_\ell\ge 2$, corresponding to the case where $(\ell-1)$ factors \smash{$G_u^{\pa{x}}$} are replaced by $G_u\Delta^xG_u$. This yields
	\begin{align}
		&\cC_{u,\bsigma,\bx}^{\pa{n}}(x,y)-\cC_{u,\bsigma,\bx}^{\pa{n},(x]}(x,y)\label{chain_x_to_0_d}\\
		&\prec {B_u^{n-1/2}(0)\mathscr{D}_{u,\bx}^{\pa{n}}\pa{x,y}\sum_{\ell\ge 2}\sum_{n_1\geq 1,n_2,\ldots,n_\ell\geq 2:n_1+\cdots+n_\ell=n+\ell-1}\Xi_{u,n_1}^{\cC_{\txt{od}}}\Xi_{u,n_2}^{\cC_{\rd}}\cdots \Xi_{u,n_\ell}^{\cC_{\rd}}}.\nonumber
	\end{align}
	Next, we control $\cC_{u,\bsigma,\bx}^{\pa{n},(x]}(x,y)$. Applying the first resolvent identity in \eqref{resolvent_identities} together with the large deviation estimate from \cite[Lemma 3.3]{ERDOS20121435}, we obtain
	\begin{equation*}
		\begin{aligned}
			\cC_{u,\bsigma,\bx}^{\pa{n},(x]}(x,y)=-G_{u,xx}\pa{HG_u^{\pa{x}}\p{\sigma_1}S^{\pa{x_1}}\cdots S^{\pa{x_{n-1}}}G_u^{\pa{x}}\p{\sigma_n}}_{xy}\prec \pbb{\sum_{a}S_{xa}\pa{\cC_{u,\bsigma,\bx}^{\pa{n},\qa{x}}}^*(y,a)\cdot \cC_{u,\bsigma,\bx}^{\pa{n},\qa{x}}(a,y) }^{1/2}.
        \end{aligned}
	\end{equation*}
	When $a=x$, the chain $\cC_{u,\bsigma,\bx}^{\pa{n},\qa{x}}(a,y)=0$ by the fact that $G_{xb}^{(x)}=0$ for $b\in\ZL$. When $a\ne x$, we again apply the second identity in \eqref{resolvent_identities} to the \smash{$G_u^{\pa{x}}$} entries in \smash{$\cC_{u,\bsigma,\bx}^{\pa{n},\qa{x}}(a,y)$}. This expresses the chain as an $n$-$\cC$-chain together with several subleading terms. In these terms, the $n$-chain is decomposed into two off-diagonal $\cC$-chains of lengths $n_1,n_2\ge 1$ and $(\ell-2)$ diagonal $\cC$-chains of lengths $n_3,\ldots, n_{\ell}\ge 2$, corresponding to the case where $(\ell-1)$ factors \smash{$G_u^{\pa{x}}$} are replaced by $G_u\Delta^xG_u$. Consequently, for $x\ne a,$
	\begin{align}
		&\cC_{u,\bsigma,\bx}^{\pa{n},\qa{x}}(a,y)-\cC_{u,\bsigma,\bx}^{\pa{n}}(a,y)\prec B_u^{n}(0)\mathscr{D}_{u,\bx}^{\pa{n}}\pa{a,y} \cdot \Xi_{u,1}^{\cC_{\txt{od}}}\Xi_{u,n}^{\cC_{\txt{od}}} \label{chain_x_to_0_od}\\
		& + {B_u^{n}(0)\mathscr{D}_{u,\bx}^{\pa{n}}\pa{a,y}\sum_{\ell\ge 2}\sum_{1\le n_1,n_2\le n-1,n_3,\ldots,n_\ell\geq 2: n_1+\cdots+n_\ell=n+\ell-1}\Xi_{u,n_1}^{\cC_{\txt{od}}}\Xi_{u,n_2}^{\cC_{\txt{od}}}\Xi_{u,n_3}^{\cC_{\rd}}\cdots \Xi_{u,n_\ell}^{\cC_{\rd}}}.\nonumber
	\end{align}
	Using \eqref{chain_x_to_0_od}, together with a simple AM-GM inequality, we deduce that for any small constant $\e>0$,
	\begin{align*}
		\qa{B_u^{n-1/2}(0)\mathscr{D}_{u,\bx}^{\pa{n}}\pa{x,y}}^{-1}  &\cC_{u,\bsigma,\bx}^{\pa{n},(x]}(x,y) \prec   \pa{\Xi_{u,2n}^{\cC_{\txt{d}}} }^{1/2} +B_u^{\e}(0)\Xi_{u,n}^{\cC_{\txt{od}}}+B_u^{1/2-\e}(0) \Xi_{u,1}^{\cC_{\txt{od}}}\Xi_{u,n}^{\cC_{\txt{od}}}\\
		+&B_u^{1/2-\e}(0)\sum_{\ell\ge 2}\sum_{1\le n_1,n_2\leq n-1,n_3,\ldots,n_\ell\geq 2: n_1+\cdots+n_\ell=n+\ell-1}\Xi_{u,n_1}^{\cC_{\txt{od}}}\Xi_{u,n_2}^{\cC_{\txt{od}}}\Xi_{u,n_3}^{\cC_{\rd}}\cdots \Xi_{u,n_\ell}^{\cC_{\rd}}.
	\end{align*}
	Combining this with \eqref{chain_x_to_0_d} yields a self-consistent estimate for $\Xi_{u,n}^{\cC_{\txt{od}}}$.  If \smash{$\Xi_{u,2}^{\cL}\le B_t^{-c}(0)$} for some small constant $c\in(0,1/100)$,  then \eqref{upper_bound_2_to_1} implies that \smash{$\Xi_{u,1}^{\cC_{\txt{od}}}\prec B_t^{-c/2}(0)$}. Solving the resulting self-consistent estimate for \smash{$\Xi_{u,n}^{\cC_{\txt{od}}}$}, and noting that $\e>0$ can be chosen arbitrarily small, we obtain
	\begin{align*}
		\Xi_{u,n}^{\cC_{\txt{od}}} \prec  &~ \pa{\Xi_{u,2n}^{\cC_{\txt{d}}} }^{1/2} +\sum_{\ell\ge 2}\sum_{n_1\geq 1,n_2,\ldots,n_\ell\geq 2:n_1+\cdots+n_\ell=n+\ell-1}\Xi_{u,n_1}^{\cC_{\txt{od}}}\Xi_{u,n_2}^{\cC_{\rd}}\cdots \Xi_{u,n_\ell}^{\cC_{\rd}}\\
		&~+B_u^{1/2}(0)\sum_{\ell\ge 2}\sum_{1\le n_1,n_2\leq n-1,n_3,\ldots,n_\ell\geq 2: n_1+\cdots+n_\ell=n+\ell-1}\Xi_{u,n_1}^{\cC_{\txt{od}}}\Xi_{u,n_2}^{\cC_{\txt{od}}}\Xi_{u,n_3}^{\cC_{\rd}}\cdots \Xi_{u,n_\ell}^{\cC_{\rd}}.
	\end{align*}
	Finally, \eqref{eq_cC_od_to_cC_d} follows from this estimate by a standard $N^{-C}$-net and perturbation argument.

	Next, we prove the reshaping estimates \eqref{reshape_2_chains_to_loops} and \eqref{reshape_chains_into_loops}. For \eqref{reshape_2_chains_to_loops}, applying the first resolvent identity in \eqref{resolvent_identities} to the diagonal $T$-variable $T_{u,xy}^{\bsigma}$, we obtain
	\begin{equation*}
		\begin{aligned}            T_{u,xy}^{\bsigma}&=\pb{G_u(\sigma_1)S^{\pa{x}}G_u\p{\sigma_2}}_{yy}=G_{u,yy}\p{\sigma_1}G_{u,yy}\p{\sigma_2}\qa{\pa{H_uG_u^{\pa{y}}(\sigma_1)S^{\pa{x}}G_u^{\pa{y}}\p{\sigma_2}H_u}_{yy}+S_{xy}}.
		\end{aligned}
	\end{equation*}
	Applying the large deviation estimate from \cite[Lemma 3.3]{ERDOS20121435}, we obtain 
	\begin{equation*}
		\begin{aligned}            T_{u,xy}^{\bsigma}=&~G_{u,yy}\p{\sigma_1}G_{u,yy}\p{\sigma_2}\qa{u\avgb{G_u^{\pa{y}}(\sigma_1)S^{\pa{x}}G_u^{\pa{y}}\p{\sigma_2}S^{\pa{y}}}+S_{xy}}\\
			&~+\opr{\avgb{G_u^{\pa{y}}(\sigma_1)S^{\pa{x}}G_u^{\pa{y}}\p{\sigma_2}S^{\pa{y}}G_u^{\pa{y}}(-\sigma_2)S^{\pa{x}}G_u^{\pa{y}}\pa{-\sigma_1}S^{\pa{y}}}}^{1/2}.
		\end{aligned}
	\end{equation*}
	We now bound the terms on the RHS. First, by \eqref{upper_bound_2_to_1} and \eqref{entrywise_bound_2_to_1}, we have
	\begin{equation*}
		\begin{aligned}
			\max_{\sigma\in\ha{-,+}}\|G_u(\sigma)-m(\sigma)\|_{\max}\prec B_u^{1/2}(0)\pa{\Xi_{u,2}^{\cL}}^{1/2}.
		\end{aligned}
	\end{equation*}
	Second, similarly to the arguments above, applying the second identity in \eqref{resolvent_identities} to each entry of $G^{(y)}_u$ yields
	\begin{equation*}
		\begin{aligned}
			\avgb{G_u^{\pa{y}}(\sigma_1)S^{\pa{x}}G_u^{\pa{y}}\p{\sigma_2}S^{\pa{y}}}-\cL_{u,xy}^{\bsigma}\prec {B_u^{2}(0)\mathscr{D}_{u}^2(|x-y|)\qa{\Xi_{u,3}^{\cC_{\rd}}+ \pb{\Xi_{u,2}^{\cC_{\rd}}}^2}}.
		\end{aligned}
	\end{equation*}
	Finally, the error term can be bounded in the same way as in \eqref{eq:n-loop-x}: 
	\begin{equation*}
		\begin{aligned}
			&~\avgb{G_u^{\pa{y}}(\sigma_1)S^{\pa{x}}G_u^{\pa{y}}\pa{\sigma_2}S^{\pa{y}}G_u^{\pa{y}}(-\sigma_2)S^{\pa{x}}G_u^{\pa{y}}\pa{-\sigma_1}S^{\pa{y}}}\\
			\prec&~B_u^{3}(0)\mathscr{D}_{u}^4(|x-y|)\pbb{\Xi_{u,4}^{\cL}+B_u(0)\sum_{n_1,\ldots,n_\ell\geq 2:n_1+\cdots+n_\ell=4+\ell}\Xi_{u,n_1}^{\cC_{\rd}}\cdots \Xi_{u,n_\ell}^{\cC_{\rd}} }.
		\end{aligned}
	\end{equation*}
	Combining the above estimates, we obtain
	\begin{equation*}
		\begin{aligned}    		T_{u,xy}^{\bsigma}&=m\p{\sigma_1}m\p{\sigma_2}\pa{u\cL_{u,xy}^{\bsigma}+S_{xy}} +B_u^{3/2}(0)\mathscr{D}_{u}^2(|x-y|)\cdot \opr{\pa{\Xi_{u,2}^{\cL}}^{3/2}+B_u^{1/2}(0) \qa{\Xi_{u,3}^{\cC_{\rd}}+ \pb{\Xi_{u,2}^{\cC_{\rd}}}^2}}\\
			&\quad +B_u^{3/2}(0)\mathscr{D}_{u}^2(|x-y|)\cdot\OO_\prec\pbb{\Xi_{u,4}^{\cL}+B_u(0)\sum_{n_1,\ldots,n_\ell\geq 2:n_1+\cdots+n_\ell=4+\ell}\Xi_{u,n_1}^{\cC_{\rd}}\cdots \Xi_{u,n_\ell}^{\cC_{\rd}}}^{1/2}.
		\end{aligned}
	\end{equation*}
	This concludes the proof of \eqref{reshape_2_chains_to_loops}.

    The bound \eqref{reshape_chains_into_loops} follows from the argument used to prove \eqref{eq_cC_d_to_loops}, together with the deterministic reshaping relation for the $\cK$-loops stated in the following lemma.
    \begin{lemma}\label{lemma_reshape_cK_chains_into_loops}
     Let $\bx=(x_1,\ldots,x_{n-1},x)\in\ZL^n$ and $\bsigma=(\sigma_1,\ldots,\sigma_n)\in\ha{-,+}^n$ with $n\geq 3$. For any $t\in[0,\tf]$, we use the abbreviation
        \begin{equation}\label{abbreviation_cK_chains_and_loops}
            \cK^{\star}_{p:q}(t):=\cK_{t,(\sigma_p,\ldots,\sigma_q),(x_p,\ldots,x_{q-1},x)}^{\star,(q-p+1)},\qquad \forall 1\leq p<q\leq n\quad \txt{and} \quad \star\in\ha{\cC,\cL}.
        \end{equation}
    Then, we have
        \begin{equation}\label{eq_interlacing_reshape}
            t\cK_{1:n}^{\cL}(t)=\sum_{r=1}^{n-1}\sum_{1=i_0<i_1<\cdots<i_{r-1}<i_r=n}(-1)^{r+1}m(\sigma_{i_0})^{-1}\prod_{j=0}^{r-1}\cK_{i_j:i_{j+1}}^{\cC}(t)m(\sigma_{i_{j+1}})^{-1}
        \end{equation}
     Equivalently, we have the reshaping relation
        \begin{equation}\label{eq_reshape_cK_chains_into_loops}
            \cK_{1:n}^{\cC}(t)=tm(\sigma_1)m(\sigma_n)\cK_{1:n}^{\cL}+\sum_{r=1}^{n-2}\sum_{1<i_1<\cdots<i_r<i_{r+1}=n}(-1)^{r+1}\cK_{1:i_1}^{\cC}(t)\prod_{j=1}^{r}\pa{m(\sigma_{i_j})^{-1}\cK_{i_j:i_{j+1}}^{\cC}(t)}.
        \end{equation}
    \end{lemma}

  We postpone the proof of this lemma until after completing the proof of \eqref{reshape_chains_into_loops}. Let $n\geq 3$, and consider a diagonal $n$-chain \smash{$\cC_{u,\bsigma,\bx}^{(n)}(x,x)$} with $\bsigma=(\sigma_1,\ldots,\sigma_n)\in\ha{-,+}^n$ and $(\bx,x)=(x_1,\ldots,x_{n-1},x)\in\ZL^n$. First, applying a large deviation estimate yields \eqref{concentration_x_to_xx}. We then apply the second identity in \eqref{resolvent_identities} to each \smash{$G_u^{(x)}$} on the RHS to recover the $G_u$-chain and loop structure, as in the argument around \eqref{eq:n-loop-x}.  Using also the simple inequality $\mathscr{D}_{u}(|a-x|)\mathscr{D}_{u}(|x-b|)\lesssim \mathscr{D}_{u}(|a-b|)$ to recover the loop-decay structure, together with the assumed bounds on the \(\Xi\)-parameters, we obtain \begin{equation} \label{eq:C_reshape_1}
    \cC_{u,\bsigma,\bx}^{\pa{n},\pa{x}}(x,x)-um(\sigma_1)m(\sigma_n)\cK_{u,\bsigma,(\bx,x)}^{\cL}  \prec  B_u^{n-1/2}(0)\mathscr{D}_{u,(\bx,x)}^{(n)} .
	\end{equation}
Similarly, following the argument around \eqref{eq:n-loop-x_diff}, we again apply the second identity in \eqref{resolvent_identities} to each \(G_u^{(x)}\) to recover the chain and loop structure, use the same decay inequality above, and invoke the assumed bounds on the \(\Xi\)-parameters. This yields
    \begin{align}\label{eq:C_reshape_2}
        \cC_{u,\bsigma,\bx}^{(n)}(x,x)-\cC_{u,\bsigma,\bx}^{(n),(x)}(x,x)=&~\sum_{r=1}^{n-2}\sum_{1<i_1<\cdots<i_r<i_{r+1}=n}(-1)^{r+1}\cK_{1:i_1}^{\cC}(u)\prod_{j=1}^{r}\pa{m(\sigma_{i_j})^{-1}\cK_{i_j:i_{j+1}}^{\cC}(u)} \nonumber\\
        &~+\OO_{\prec}\pa{B_u^{n-1/2}(0)\mathscr{D}_{u,(\bx,x)}^{(n)}}.
    \end{align}
    Combining \eqref{eq:C_reshape_1} and \eqref{eq:C_reshape_2}, we conclude \eqref{reshape_chains_into_loops} by \eqref{eq_reshape_cK_chains_into_loops}.

	For the proof of \Cref{lemma_reshape_cK_chains_into_loops} and for the later proof of Lemma \ref{lemma_n_cK_tree_bound}, we collect several basic properties of the $\cK$-loops in the following lemma.
	
	\begin{lemma}\label{lemma_basic_property_wh_cK}
		With the notation of Definition \ref{def_n_cK_loop}, the following properties hold.
		\begin{itemize}
			\item {\bf Recursive relation}: For $\bsigma=(\sigma_1,\ldots,\sigma_{n+1})\in\ha{-,+}^{n+1}$ and $\bx=(x_1,\ldots,x_{n+1})\in\ZL^{n+1}$ with $n\geq 1$, we have
			\begin{align}
				\qquad\qquad   &\wh \cK_{t,\bsigma,\bx}^{\pa{n+1}}=m\p{\sigma_1}\wh \cK_{t,\pa{\sigma_2,\ldots,\sigma_{n+1}},\pa{x_2,\ldots,x_n,x_1}}^{\pa{n}}\pa{1-tm\p{\sigma_1}m\p{\sigma_{n+1}}S}^{-1}_{x_1x_{n+1}}\label{recursive_relation_n_cK}\\
				&+tm\p{\sigma_1}\sum_{k=2}^{n}\sum_{x,y}\wh \cK_{t,\pa{\sigma_1,\ldots,\sigma_k},\pa{x_1,\ldots,x_{k-1},y}}^{\pa{k}}S_{xy}\cdot\wh \cK_{t,\pa{\sigma_k,\ldots,\sigma_{n+1}},\pa{x_k,\ldots,x_n,x}}^{\pa{n-k+2}}\pa{1-tm\p{\sigma_1}m\p{\sigma_{n+1}}S}^{-1}_{xx_{n+1}}.\nonumber
			\end{align}
			\item {\bf Shift invariance}: For any cyclic shift $\tau_k$ acting on $\pa{s_1,\dots,s_n}$ by $\tau_{k}\pa{s_1,\ldots,s_n}=\pa{s_{k+1},\ldots,s_{k+n}}$, with the convention $s_i \equiv s_j$ whenever $i=j \mod n$, we have for any $k\in\qq{n}$:
			\begin{equation}\label{shift_invariance}
				\wh \cK_{t,\bsigma,\bx}^{\pa{n}}=\wh \cK_{t,\tau_k{\bsigma},\tau_k\bx}^{\pa{n}}.
			\end{equation}
			\item {\bf Reduction to matrices}: For $\bsigma\in\ha{-,+}^n$, $i\neq j\in\qq{n}$, and $t\in\qa{0,\tf}$, the matrix
			\begin{equation}\label{i_j_symmetry_wh_cK}            A_{x_i,x_j}:=\sum_{x_k\in\ZL:\, k\in\qq{n}\setminus\ha{i,j}}\wh \cK_{t,\bsigma,(x_1,\ldots,x_n)}^{(n)}
			\end{equation}
			is a function of the variance profile matrix $S$, and it is symmetric. Moreover, the quantity 
			$$\sum_{x_k\in\ZL:\, k\in\qq{n}\setminus\ha{i}}\wh \cK_{t,\bsigma,(x_1,\ldots,x_n)}^{(n)}$$
			 does not depend on $x_i$, and hence is a constant.
		\end{itemize}
	\end{lemma}
	\begin{proof}
		For the recursive relation \eqref{recursive_relation_n_cK}, an argument similar to the proof of Lemma A.1 in \cite{fan2025blockreductionmethodrandom} shows that the quantities defined recursively through \eqref{recursive_relation_n_cK} also satisfy the evolution equation \eqref{evolution_wh_cK_loops}. Moreover, at $t=0$, starting from
		\smash{$\wh \cK_{0,\sigma,x}^{(1)}=m(\sigma)$}, the recursion \eqref{recursive_relation_n_cK} yields \smash{$\wh \cK_{0,(\sigma_1,\ldots,\sigma_n),(x_1,\ldots,x_n)}^{(n)}=m(\sigma_1)\cdots m(\sigma_n)\delta_{x_1x_2}\cdots \delta_{x_nx_1}$} for any fixed $n$. Therefore, by the uniqueness of the solution to the system of differential equations \eqref{evolution_wh_cK_loops}, we obtain \eqref{recursive_relation_n_cK}.
		The shift invariance \eqref{shift_invariance} follows immediately from the tree representation in \cite[Lemma 4.16]{truong2025localizationlengthfinitevolumerandom}. Finally, the third property follows directly from an induction argument based on \eqref{recursive_relation_n_cK}.\footnote{In this induction argument, we may assume $i=1$ by the shift invariance \eqref{shift_invariance}.}
	\end{proof}

    \begin{proof}[\bf Proof of \Cref{lemma_reshape_cK_chains_into_loops}]We employ an induction argument on $n\geq 3$ to prove the identity \eqref{eq_interlacing_reshape}, which is equivalent to the reshaping relation \eqref{eq_reshape_cK_chains_into_loops}. For $n\geq 2$, consider $(x_1,\ldots,x_n,x)\in\ZL^{n+1}$ and $(\sigma_1,\ldots,\sigma_{n+1})\in\ha{-,+}^{n+1}$, and adopt the notation introduced in \eqref{abbreviation_cK_chains_and_loops}. Regarding both sides of \eqref{recursive_relation_n_cK} as vectors indexed by $x_{n+1}$, we multiply both sides by the matrix $(1-tm(\sigma_1)m(\sigma_{n+1})S)$. Then, by the definition of the $\cK$-chains and loops in \Cref{def_n_cK_loop}, we obtain  
\begin{equation}\label{eq_first_reshape}
\cK_{1:n+1}^{\cC}-tm(\sigma_1)m(\sigma_{n+1})\cK_{1:n+1}^{\cL}=m(\sigma_1)S_{x_1x}\cK_{2:n+1}^{\cC}+  tm(\sigma_1)\sum_{k=2}^{n}\cK_{1:k}^{\cL}\cK_{k:n+1}^{\cC}.
\end{equation}
Setting \(n=2\) in \eqref{eq_first_reshape}, we obtain
\begin{equation*}
\begin{aligned}
\cK_{1:3}^{\cC}=tm(\sigma_1)m(\sigma_{3})\cK_{1:3}^{\cL}+m(\sigma_1)S_{x_1x}\cK_{2:3}^{\cC}+tm(\sigma_1)\cK_{1:2}^{\cL}\cK_{2:3}^{\cC} =tm(\sigma_1)m(\sigma_{3})\cK_{1:3}^{\cL}+m(\sigma_2)^{-1}\cK_{1:2}^{\cC}\cK_{2:3}^{\cC},
\end{aligned}
\end{equation*}
where in the second step we used the identity 
\(S_{x_1x}+t\cK_{1:2}^{\cL}=m(\sigma_1)^{-1}m(\sigma_2)^{-1}\cK_{1:2}^{\cC}.\)
Hence, the \(n=3\) case of \eqref{eq_reshape_cK_chains_into_loops} follows.
Now assume that \eqref{eq_interlacing_reshape} (and hence \eqref{eq_reshape_cK_chains_into_loops}) holds for all \(3\leq k\leq n\). Applying the induction hypothesis to \(\cK_{1:k}^{\cL}\) for each \(3\leq k\leq n\), and using the above identity for \(k=2\), we obtain
\begin{equation*}
        \begin{aligned}
            &~-m(\sigma_1)m(\sigma_{n+1})^{-1}S_{x_1x}\cK_{2:n+1}^{\cC}-tm(\sigma_1)m(\sigma_{n+1})^{-1}\sum_{k=2}^{n}\cK_{1:k}^{\cL}\cK_{k:n+1}^{\cC}\\
            =&~\sum_{k=2}^{n}\sum_{r=1}^{k-1}(-1)^{r}\sum_{1=i_0<i_1<\cdots<i_{r-1}<i_r=k}\prod_{j=0}^{r-1}\qa{\cK_{i_j:i_{j+1}}^{\cC}m(\sigma_{i_{j+1}})^{-1}}\cdot \cK_{k:n+1}^{\cC}m(\sigma_{n+1})^{-1}\\
            =&~\sum_{r=1}^{n-1}(-1)^{r}\sum_{1=i_0<i_1<\cdots<i_{r}<i_{r+1}=n+1}\prod_{j=0}^{r}\qa{\cK_{i_j:i_{j+1}}^{\cC}m(\sigma_{i_{j+1}})^{-1}}\\
            =&~\sum_{r=1}^{n}(-1)^{r+1}\sum_{1=i_0<i_1<\cdots<i_{r-1}<i_{r}=n+1}\prod_{j=0}^{r-1}\qa{\cK_{i_j:i_{j+1}}^{\cC}m(\sigma_{i_{j+1}})^{-1}}-\cK_{1:n+1}^{\cC}m(\sigma_{n+1})^{-1}.
        \end{aligned}
    \end{equation*}
Substituting this identity into \eqref{eq_first_reshape} yields the case \(n+1\) of \eqref{eq_interlacing_reshape}. This completes the induction and proves \eqref{eq_interlacing_reshape}.
    \end{proof}

	\subsection{Proof of Lemma \ref{lemma_dynamical_inequality}}\label{subsec:pf_lemma_dynamical_inequality}

   We prove only the case $n\geq 3$; when $n=2$, the estimate \eqref{eq_dynamical_inequality_n_2} follows from a similar, but considerably simpler, argument.
    We begin by considering the evolution equation \eqref{evolution_n_cL_cK_integrated} in the non-alternating case, i.e., when $\sigma_i=\sigma_{i+1}$ for some $i\in\qq{n}$.
	
	\begin{proof}[\bf Proof of Lemma \ref{lemma_dynamical_inequality}: non-alternating case]
		We illustrate the argument using the quadratic error term \eqref{eq:quad_error_n} in \eqref{evolution_n_cL_cK_integrated} as an example.
		Each term is the product of an $(\cL-\cK)$-loop and a $(\cC-\cK^{\cC})$-chain, and at least one of them has length at most $n-1$. Without loss of generality, we assume that the length of the loop is less than or equal to that of the chain.
		We then bound the $(\cL-\cK)$-loops using the $\Xi^{\cL-\cK}$ parameters and the $(\cC-\cK^{\cC})$-chains using the $\Xi^{\cC_{\rd}}$ parameters together with Lemma \ref{lemma_n_cK_bound}. Combining these bounds with the convolution inequality \eqref{example_triangular_loop_recovery} to recover the loop structure of \smash{$\mathscr{D}_{u,\bx}^{\pa{n}}$}, we obtain
		\begin{equation*}
			\begin{aligned}
				\mathcal{E}^{(n)}_{u, \boldsymbol{\sigma}, \bx}\prec \frac{1}{1-u}B_u^{n}(0)\mathscr{D}_{u,\bx}^{\pa{n}}\cdot \max_{k=2}^{n-1} \wh \Xi_{u,k}^{\cL-\cK}\wh \Xi_{u,n+2-k}^{\cC_{\rd}}.
			\end{aligned}
		\end{equation*}
		Combining this bound with the evolution kernel bound \eqref{evolution_kernel_bound_2}, we obtain 
		\begin{align*}            \int_{s}^{r}\pa{\bU_{u,r,\bsigma}^{\pa{n}}\circ\mathcal{E}^{(n)}_{u, \boldsymbol{\sigma}}}_{\bx}\dd u&\prec B_r^n(0)\mathscr{D}_{r,\bx}^{\pa{n}}\cdot \int_{s}^r\frac{1}{1-u}\max_{k=2}^{n-1} \wh \Xi_{u,k}^{\cL-\cK}\wh \Xi_{u,n+2-k}^{\cC_{\rd}}\,\rd u\\
			&\prec B_r^n(0)\mathscr{D}_{r,\bx}^{\pa{n}}\cdot \sup_{u\in\qa{s,r}} \max_{k=2}^{n-1} \wh \Xi_{u,k}^{\cL-\cK}\wh \Xi_{u,n+2-k}^{\cC_{\rd}}.
		\end{align*}
		A nearly identical argument bounds the second term on the RHS of \eqref{evolution_n_cL_cK_integrated} by \smash{$B_r^n(0)\mathscr{D}_{r,\bx}^{\pa{n}}\cdot \sup_{u\in\qa{s,r}} \max_{k=2}^{n-1} \wh \Xi_{u,k}^{\cL-\cK}$}. Combining the induction hypothesis \eqref{flow_n_loops_local_law_supercritical} at time $s$ with the evolution kernel bound \eqref{evolution_kernel_bound_2}, we bound the first term on the RHS of \eqref{evolution_n_cL_cK_integrated} by \smash{$B_r^n(0)\mathscr{D}_{r,\bx}^{\pa{n}}$}. 
		Next, applying the averaged local law in \eqref{step_2_sharp_local_law_supercritical} and bounding the $(n+1)$-chain in \eqref{eq:LW_long} using the $\Xi^{\cC_{\rd}}$ parameter, we estimate the light-weight term in \eqref{evolution_n_cL_cK_integrated} by \smash{$B_r^n(0)\mathscr{D}_{r,\bx}^{\pa{n}}\cdot \sup_{u\in\qa{s,r}} \wh \Xi_{u,n+1}^{\cC_{\rd}}$}. 
		
		It remains to bound the martingale term in \eqref{evolution_n_cL_cK_integrated}. Applying Lemma \ref{lemma_bound_martingale_to_cB_supercritical} and estimating the sum of the $(2n+2)$-chain \eqref{def_tensor_cB_supecritical} over $x$ using the $\Xi^{\cC_{\rd}}$ parameter together with the convolution inequality \eqref{example_triangular_loop_recovery}, we obtain
		\begin{equation}
			\label{eq:martingale_choice_1}
			\int_{s}^{r}\pa{\bU_{u,r,\bsigma}^{\pa{n}}\circ
				\dd\mathcal{B}^{(n)}_{u, \boldsymbol{\sigma}}}_{\bx} \prec B_r^n(0)\mathscr{D}_{r,\bx}^{\pa{n}}\cdot \sup_{u\in\qa{s,r}}\pa{\wh \Xi_{u,2n+2}^{\cC_{\rd}}}^{1/2}, 
		\end{equation}
		which corresponds to the choice \eqref{eq:choose_Lambda1}. 
		To obtain \eqref{eq:choose_Lambda2}, we recall \eqref{step_2_cU_martingale_to_cU_X}. From the estimates \eqref{step_2_sharp_local_law_supercritical} and \eqref{step_2_loops_bound_supercritical} established in Step~2, together with \eqref{step_2_cC_rd_bound_supercritical} and \eqref{step_2_cC_od_bound_supercritical} (which follow from \Cref{lemma_Xi_reshape}), we have
		\begin{equation*}        \Xi_{u,4}^{\cL}+\Xi_{u,2}^{\cC_{\rd}}+\Xi_{u,2}^{\cC_{\txt{od}}}\prec 1.
		\end{equation*}
		Substituting this estimate into Lemma \ref{lemma_n_loop_contraction} yields 
		\begin{equation*}
			\begin{aligned}
				X_{u,\bsigma,\bx}^{\cB,\pa{n}}\prec \frac{1}{\sqrt{1-u}}\pa{\Xi_{u,n-1}^{\cC_{\rd}} \Xi_{u,n}^{\cC_{\rd}}}^{1/2}\cdot B_u^{n-3/4}(0)\mathscr{D}_{u,\bx}^{(n)}.
			\end{aligned}
		\end{equation*}
		Applying the evolution kernel bound \eqref{evolution_kernel_bound_1} (since \eqref{evolution_kernel_bound_2} is not used here, the following bound also applies in the alternating case), and combining this estimate with \eqref{step_2_cU_martingale_to_cU_X} and \eqref{eq:Xi_by_whXi}, we obtain
		\begin{align}
			\int_{s}^{r}\pa{\bU_{u,r,\bsigma}^{\pa{n}}\circ
				\dd\mathcal{B}^{(n)}_{u, \boldsymbol{\sigma}}}_{\bx}&\prec \qa{\sup_{u\in\qa{s,r}}\Xi_{u,n-1}^{\cC_{\rd}}\Xi_{u,n}^{\cC_{\rd}}\cdot\int_{s}^{r}\frac{1}{1-u}\pa{\frac{1-u}{1-r}B_r^{n-3/4}(0)\mathscr{D}_{r,\bx}^{(n)}}^2\,\rd u}^{1/2} \nonumber\\
			&\prec\frac{1-s}{1-r}\sup_{u\in\qa{s,r}}\pa{\wh \Xi_{u,n-1}^{\cC_{\rd}}\wh \Xi_{u,n}^{\cC_{\rd}}}^{1/2}\cdot B_r^{n-3/4}(0)\mathscr{D}_{r,\bx}^{(n)}.\label{eq:martingale_weak_n}
		\end{align}
		This yields the parameter choice \eqref{eq:choose_Lambda2} and completes the estimate of the non-alternating loops by the RHS of \eqref{eq_dynamical_inequality}, after rescaling by \smash{$B_r^n(0)\mathscr{D}_{r,\bx}^{\pa{n}}$}.
	\end{proof}
	
	For the alternating case with $\sigma_i\neq \sigma_{i+1}$ for all $i\in\qq{n}$, we employ the following sum-zero operators.
	\begin{definition}\label{def_sum_zero_operator_supercritical}
		For any $n\geq 2$ and $u\in\qa{s,t}$, the \emph{partial sum operator} $P_u$ acts on an $n$-tensor $\cA^{(n)}$ as
		\begin{equation*}
			\begin{aligned}
				\pB{P_u\circ\cA^{(n)}}_{\bx}:=(1-u)^{n-1}\sum_{a_2,\ldots,a_n\in\ZL}\cA_{x_1a_2\cdots a_n}^{(n)}\prod_{i=2}^{n}\Theta_{u,x_1x_i},\qquad \forall \bx=(x_1,\ldots,x_n)\in\ZL^n.
			\end{aligned}
		\end{equation*}
		The \emph{sum-zero operator} $Q_u$ is then defined by $Q_u:=1-P_u$. 
	\end{definition}
	The following basic properties of the above operators can be verified directly.     
    \begin{lemma}\label{lemma:convolve_sumzero}
		For any $n\geq 2$, $u\in\qa{s,t}$, and $n$-tensor $\cA^{(n)}$, we have  
		\begin{equation}\label{eq:first_evolution_sumzero}
			\begin{aligned}
				\sum_{x_2,\ldots,x_n\in\ZL}\pa{Q_u\circ \cA^{(n)}}_{x_1x_2\cdots x_n}=0,\quad \forall x_1\in\ZL,
			\end{aligned}
		\end{equation}
		\begin{equation}\label{P_u_Q_u_bound}
			\begin{aligned}
				\normb{P_u\circ \cA^{(n)}}_{\mathscr{D}_u^{-1}}\prec \normb{\cA^{(n)}}_{\mathscr{D}_u^{-1}},\quad\normb{Q_u\circ \cA^{(n)}}_{\mathscr{D}_u^{-1}}\prec \normb{\cA^{(n)}}_{\mathscr{D}_u^{-1}}.
			\end{aligned}
		\end{equation}
	\end{lemma}        
	\begin{proof}      
		The sum-zero property \eqref{eq:first_evolution_sumzero} follows directly from the identity $\sum_{x_i\in\ZL}\Theta_{u,x_1x_i}=(1-u)^{-1}$. To prove \eqref{P_u_Q_u_bound}, using \eqref{Theta_bound_opposite_charge}, \eqref{example_triangular_2_edge_to_1}, and \eqref{bound_one_point_to_loop}, we obtain         \begin{equation}\label{P_u_action_bound}
			\begin{aligned}
				\pa{P_u\circ \cA^{(n)}}_{\bx}&=(1-u)^{n-1}\sum_{a_2,\ldots,a_n\in\ZL}\cA_{x_1a_2\cdots a_n}^{(n)}\prod_{i=2}^{n}\Theta_{u,x_1x_i}\\
				&\prec\ell_u^{-n+1}\normb{\cA^{(n)}}_{\mathscr{D}_u^{-1}}\sum_{a_2,\ldots,a_n\in\ZL}\mathscr{D}_{u,(x_1,a_2,\ldots,a_n)}^{\pa{n}}\cdot\prod_{i=2}^{n}\mathscr{D}_{u}^2(|x_1-x_i|)\\
				&\prec\normb{\cA^{(n)}}_{\mathscr{D}_u^{-1}}\cdot\prod_{i=2}^{n}\mathscr{D}_{u}^2(|x_1-x_i|)\lesssim\normb{\cA^{(n)}}_{\mathscr{D}_u^{-1}}\cdot\mathscr{D}_{u,\bx}^{(n)}.
			\end{aligned}
		\end{equation}
		This proves the first bound in \eqref{P_u_Q_u_bound}; the second bound follows immediately since $Q_u=1-P_u$. \end{proof}

	We are now ready to complete the proof of Lemma \ref{lemma_dynamical_inequality} by treating the alternating case.
	
	\begin{proof}[\bf Proof of Lemma \ref{lemma_dynamical_inequality}]
		In the case $\sigma_i\neq \sigma_{i+1}$ for all $i\in\qq{n}$, we decompose
		\begin{equation*}
			\begin{aligned}
				(\cL-\cK)_{r,\bsigma,\bx}^{(n)}=\qa{P_r\circ (\cL-\cK)_{r,\bsigma}^{(n)}}_{\bx}+\qa{Q_r\circ (\cL-\cK)_{r,\bsigma}^{(n)}}_{\bx}.
			\end{aligned}
		\end{equation*}
		For the first term, we apply Ward's identities \eqref{eq_Ward0} and \eqref{WI_calK} at a vertex $x_k$, and bound the resulting $(n-1)$-loops with the $\Xi^{\cL-\cK}$-parameters: 
		\begin{align}\label{eq:PrL_alternating0} \sum_{x_k} (\cL-\cK)_{r,\bsigma,\bx}^{(n)} \prec \frac{1}{1-r}B_r^{n-1}(0)\wh \Xi_{r,n-1}^{\cL-\cK} \cdot \mathscr{D}_{r,\bx^{(k)}}^{(n-1)},\quad \text{with}\quad \bx^{(k)}:=(x_1,\ldots, x_{k-1},x_{k+1},\ldots, x_n). 
		\end{align}
		Applying this bound with $k=n$, and using the convolution inequality \eqref{example_triangular_2_edge_to_1} to sum over $x_2,\ldots,x_{n-1}$ in \smash{$\mathscr{D}_{r,\bx^{(n)}}^{(n-1)}$}, together with \eqref{Theta_bound_opposite_charge} and \eqref{bound_one_point_to_loop}, we obtain
		\begin{align}\label{eq:PrL_alternating}
			\qa{P_r\circ (\cL-\cK)_{r,\bsigma}^{(n)}}_{\bx}&\prec \p{1-r}^{n-2}B_r^{2n-2}(0)\wh \Xi_{r,n-1}^{\cL-\cK}\cdot \prod_{i=2}^{n}\mathscr{D}_{r}^2(|x_1-x_i|) \cdot \sum_{a_2,\ldots,a_{n-1}}\mathscr{D}_{r,(x_1,a_2,\ldots,a_{n-1})}^{(n-1)} \nonumber\\
			&\lesssim B_r^{n}(0)\wh \Xi_{r,n-1}^{\cL-\cK}\cdot \mathscr{D}_{r,\bx}^{(n)}.
		\end{align}
		
		We next analyze $\q{Q_r\circ (\cL-\cK)_{r,\bsigma}^{(n)}}_{\bx}$, which satisfies the following equation by direct computation and Duhamel's principle (see also equation (5.91) in \cite{Band1D}):
		\begin{align}
			&\qa{Q_r\circ (\cL-\cK)_{r,\bsigma}^{(n)}}_{\bx}=\pa{\bU_{s,r,\bsigma}^{\pa{n}}\circ Q_s\circ \pa{\cL-\cK}_{s,\bsigma}^{\pa{n}}}_{\bx}+\sum_{l_{\cK}=3}^{n}\int_{s}^{r}\pa{\bU_{u,r,\bsigma}^{\pa{n}}\circ Q_u\circ\qa{\mathfrak{D}_{l_{\cK}}(\cL-\cK)}_{u,\bsigma}^{\pa{n}}}_{\bx}\, \rd u\nonumber\\
			&+ \int_{s}^{r}\pa{\bU_{u,r,\bsigma}^{\pa{n}}\circ Q_u\circ\mathcal{E}^{(n)}_{u, \boldsymbol{\sigma}}}_{\bx}\dd u +\int_{s}^{r}\pa{\bU_{u,r,\bsigma}^{\pa{n}}\circ Q_u\circ\mathcal{W}^{(n)}_{u, \boldsymbol{\sigma}}}_{\bx}
			\dd u + \int_{s}^{r}\pa{\bU_{u,r,\bsigma}^{\pa{n}}\circ Q_u\circ
				\dd\mathcal{B}^{(n)}_{u, \boldsymbol{\sigma}}}_{\bx}
			\nonumber\\
			&+\int_{s}^{r}\pa{\bU_{u,r,\bsigma}^{\pa{n}}\circ \qb{Q_u,\bTheta_{u,\bsigma}^{(n)}}\circ(\cL-\cK)_{u,\bsigma}^{\pa{n}}}_{\bx}\, \rd u-\int_{s}^{r}\pa{\bU_{u,r,\bsigma}^{\pa{n}}\circ\pa{\partial_{u}P_u}\circ(\cL-\cK)_{u,\bsigma}^{\pa{n}}}_{\bx}\, \rd u,\label{eq:derv_Qr_L-K}
		\end{align}
		where the operator $\bTheta_{u,\bsigma}^{(n)}$ is defined in \eqref{eq:bTheta_n}, $[Q_u,\bTheta_{u,\bsigma}^{(n)}]$ denotes the commutator, and $\partial_u P_u$ is defined by     \begin{equation*}
			\begin{aligned}
				\pB{(\partial_u P_u)\circ\cA^{(n)}}_{\bx}:=\sum_{a_2,\ldots,a_n\in\ZL}\cA_{x_1a_2\cdots a_n}^{(n)}\partial_u\pa{(1-u)^{n-1}\prod_{i=2}^{n}\Theta_{u,x_1x_i}},\quad \forall \bx=(x_1,\ldots,x_n)\in\ZL^n.
			\end{aligned}
		\end{equation*}
		Using the sum-zero property \eqref{eq:first_evolution_sumzero}, the bound \eqref{P_u_Q_u_bound}, and the evolution kernel bound \eqref{evolution_kernel_bound_3}, the first four terms on the RHS of \eqref{eq:derv_Qr_L-K} can be treated exactly as in the non-alternating case. For the last two terms, one checks that both \smash{$\qb{Q_u,\bTheta_{u,\bsigma}^{(n)}}\circ(\cL-\cK)_{u,\bsigma}^{\pa{n}}$} and \smash{$\pa{\partial_{u}P_u}\circ(\cL-\cK)_{u,\bsigma}^{\pa{n}}$} satisfy the sum-zero property \eqref{sum_zero_property_supercritical}. Moreover, we estimate them as follows. Combining the following bound \[\partial_u\Theta_{u,xy}=\sum_a \Theta_{u,xa}\Theta_{u,ay} \prec \sum_a B_u(|x-a|)B_u(|a-y|) \prec (1-u)^{-1}B_{u}(|x-y|),\]
		with the argument leading to \eqref{eq:PrL_alternating}, we obtain
		\begin{equation}\label{eq:term6}
			\pa{\pa{\partial_{u}P_u}\circ(\cL-\cK)_{u,\bsigma}^{\pa{n}}}_{\bx}\prec \frac{1}{1-u}B_u^{n}(0)\wh \Xi_{u,n-1}^{\cL-\cK}\cdot \mathscr{D}_{u,\bx}^{(n)}.
		\end{equation}
		Similarly, using \eqref{Theta_bound_opposite_charge} and \eqref{eq:PrL_alternating}, we obtain
		\begin{align*}
			\pa{\bTheta_{u,\bsigma}^{(n)}\circ P_u \circ(\cL-\cK)_{u,\bsigma}^{\pa{n}}}_{\bx}  \prec  B_u^{n+1}(0)\wh \Xi_{u,n-1}^{\cL-\cK}\sum_{i=1}^n\sum_{a_i}\mathscr{D}_u^2(|x_i-a_i|) \mathscr{D}_{u,\bx^{(i)}(a_i)}^{(n)} \prec \frac{1}{1-u}B_u^{n}(0)\wh \Xi_{u,n-1}^{\cL-\cK}\cdot \mathscr{D}_{u,\bx}^{(n)},
		\end{align*}
		where in the second step we applied the convolution inequality \eqref{example_triangular_loop_recovery}. On the other hand, we have
		\begin{align*}
			&\pa{P_u \circ \bTheta_{u,\bsigma}^{(n)}\circ (\cL-\cK)_{u,\bsigma}^{\pa{n}}}_{\bx}=(1-u)^{n-1}\prod_{i=2}^{n}\Theta_{u,x_1x_i} \cdot \pbb{\sum_{\ba}\Theta_{u,x_1a_1}(\cL-\cK)_{u,\bsigma,\ba}^{\pa{n}} + \frac{n-1}{1-u}\sum_{\ba'}(\cL-\cK)_{u,\bsigma,\ba^{(1)}(x_1)}^{\pa{n}}}\\
			&\prec \frac{ B_u^{n}(0)\wh \Xi_{u,n-1}^{\cL-\cK}}{(1-u)\ell_u^{n-1}} \prod_{i=2}^{n}\mathscr{D}_{u}^2(|x_1-x_i|) \cdot \pbb{\sum_{a_1,\ldots, a_{n-1}}\mathscr{D}_u^2(|x_1-a_1|) \mathscr{D}_{u,\ba^{(n)}}^{(n-1)}+ \ell_u \sum_{a_2,\ldots, a_{n-1}} \mathscr{D}_{u,(x_1,a_2,\ldots,a_{n-1})}^{(n-1)}}\\
			&\prec \frac{1}{1-u}B_u^{n}(0)\wh \Xi_{u,n-1}^{\cL-\cK}\cdot \mathscr{D}_{u,\bx}^{(n)}.
		\end{align*}
		Here we denote $\ba:=(a_1,\ldots,a_n)\in \ZL^n$, $\ba':=(a_2,\ldots,a_n)\in \ZL^{n-1}$, $\ba^{(n)}:=(a_1,\ldots,a_{n-1})\in \ZL^{n-1}$, and $\ba^{(1)}(x_1):=(x_1,a_2,\ldots,a_n)\in \ZL^n$. In the second step, we used \eqref{Theta_bound_opposite_charge} and \eqref{eq:PrL_alternating0}; in the last step, we applied the convolution inequalities \eqref{example_triangular_loop_recovery} and \eqref{example_triangular_2_edge_to_1}, together with \eqref{bound_one_point_to_loop}.
		Combining the above two bounds yields 
		\begin{align}\label{eq:term5}
			\pa{\qb{Q_u,\bTheta_{u,\bsigma}^{(n)}}\circ(\cL-\cK)_{u,\bsigma}^{\pa{n}}}_{\bx} \prec \frac{1}{1-u}B_u^{n}(0)\wh \Xi_{u,n-1}^{\cL-\cK}\cdot \mathscr{D}_{u,\bx}^{(n)}.
		\end{align}
		Together with \eqref{eq:term6}, applying the evolution kernel bound \eqref{evolution_kernel_bound_3} and integrating over $u$, we conclude that the last two terms in \eqref{eq:derv_Qr_L-K} are bounded by \smash{$\OO_\prec(B_r^{n}(0)\cdot\sup_{u\in[s,r]}\wh \Xi_{u,n-1}^{\cL-\cK}\cdot \mathscr{D}_{r,\bx}^{(n)})$}.

		It remains to handle the martingale term in \eqref{eq:derv_Qr_L-K}. For the choice \eqref{eq:choose_Lambda1}, we again apply Lemma \ref{lemma_bound_martingale_to_cB_supercritical} and express the resulting quadratic variation as a $(2n)$-tensor satisfying the sum-zero property (see equations (5.103) and (5.104) in \cite{Band1D}). This tensor can then be estimated in the same way as in the non-alternating case. Applying the evolution kernel bound \eqref{evolution_kernel_bound_3}, integrating over $u$, and using Markov's inequality, we bound the martingale term as in \eqref{eq:martingale_choice_1}. We omit the details. 
		For the choice \eqref{eq:choose_Lambda2}, we write 
		\begin{equation*}
			\begin{aligned}
				\int_{s}^{r}\pa{\bU_{u,r,\bsigma}^{\pa{n}}\circ Q_u\circ
					\dd\mathcal{B}^{(n)}_{u, \boldsymbol{\sigma}}}_{\bx}=\int_{s}^{r}\pa{\bU_{u,r,\bsigma}^{\pa{n}}\circ
					\dd\mathcal{B}^{(n)}_{u, \boldsymbol{\sigma}}}_{\bx}-\int_{s}^{r}\pa{\bU_{u,r,\bsigma}^{\pa{n}}\circ P_u\circ
					\dd\mathcal{B}^{(n)}_{u, \boldsymbol{\sigma}}}_{\bx}.
			\end{aligned}
		\end{equation*}
		The first term is bounded as in \eqref{eq:martingale_weak_n} (recall the remark preceding that equation). For the second term, applying Lemma \ref{lemma_bound_martingale_to_cB_supercritical} together with a CS inequality as in \eqref{step_2_cU_martingale_to_cU_X}, we obtain that
		\begin{equation*}
			\begin{aligned}            \E\absa{\int_{s}^{r}\pa{\bU_{u,r,\bsigma}^{\pa{n}}\circ P_u\circ
						\dd\mathcal{B}^{(n)}_{u, \boldsymbol{\sigma}}}_{\bx}}^p\lesssim \E\pa{\int_{s}^{r} \pa{\cU_{u,r}^{\otimes n}\circ P_u\circ X_{u,\bsigma}^{\cB,\pa{n}}}_{\bx}^2\,\rd u }^{p/2}
			\end{aligned}
		\end{equation*}
		for any fixed $p\in 2\N$.  
		Applying the bounds \eqref{P_u_Q_u_bound} and \eqref{evolution_kernel_bound_1}, integrating over $u$, and using Markov's inequality, we again obtain a bound of the same form as in \eqref{eq:martingale_weak_n}.

		In summary, we bound the RHS of \eqref{eq:derv_Qr_L-K} by
		\begin{equation*}
			\begin{aligned}
				\qa{Q_r\circ (\cL-\cK)_{r,\bsigma}^{(n)}}_{\bx}\prec B_r^n(0)\mathscr{D}_{r,\bx}^{\pa{n}}\cdot\ha{ \sup_{u\in\qa{s,r}}\pa{\max_{k=2}^{n-1}\wh \Xi_{u,k}^{\cL-\cK}+\max_{k=2}^{n-1} \wh \Xi_{u,k}^{\cL-\cK}\wh \Xi_{u,n+2-k}^{\cC_{\rd}}+\wh \Xi_{u,n+1}^{\cC_{\rd}}}+\Lambda_{r,n}^{1/2}}.
			\end{aligned}
		\end{equation*}
		Together with \eqref{eq:PrL_alternating}, this yields the desired estimate for the alternating loops by the RHS of \eqref{eq_dynamical_inequality}, and hence completes the proof. \end{proof}

	\begin{proof}[\bf Proof of \eqref{eq:EL-K_samesigma} when $\sigma_1\ne \sigma_2$]
Using the sum-zero operator from Definition \ref{def_sum_zero_operator_supercritical}, we decompose
	\begin{equation*}
		\begin{aligned}
			\E(\mathcal{L} - \mathcal{K})^{(2)}_{r, \boldsymbol{\sigma}, \bx}= \qb{P_r\circ\E(\mathcal{L} - \mathcal{K})^{(2)}_{r, \boldsymbol{\sigma}}}_{\bx}+\qb{Q_r\circ \E(\mathcal{L} - \mathcal{K})^{(2)}_{r, \boldsymbol{\sigma}}}_{\bx}.
		\end{aligned}
	\end{equation*}
	For the first term, applying Ward's identity and \eqref{expected_single_resolvent_local_law_supercritical}, we obtain
	\begin{equation}\label{eq:PEL-K_samesigma}
		\qb{P_r\circ\E(\mathcal{L} - \mathcal{K})^{(2)}_{r, \boldsymbol{\sigma}}}_{\bx}\prec B_r^{3}(0)\mathscr{D}_{r,\bx}^{(2)}.
	\end{equation}
	For the second term, we analyze its evolution equation by taking expectations in \eqref{eq:derv_Qr_L-K} with $n=2$. The second and martingale terms vanish. The first, third, and fourth terms can be bounded by \smash{$\OO_\prec(B_r^{5/2}(0)\mathscr{D}_{r,\bx}^{(2)})$} using arguments analogous to the case $\sigma_1=\sigma_2$, together with the evolution kernel bound \eqref{evolution_kernel_bound_3}.
	The last two terms can be estimated as in \eqref{eq:term6} and \eqref{eq:term5}, and improve by an additional factor $B_u(0)$ after taking expectations, by \eqref{expected_single_resolvent_local_law_supercritical}:
	\begin{equation}\label{eq:exp_term56}
		\absa{\pa{\qb{Q_u,\bTheta_{u,\bsigma}^{(n)}}\circ \E(\cL-\cK)_{u,\bsigma}^{\pa{2}}}_{\bx}}+    \absa{\pa{\pa{\partial_{u}P_u}\circ\E(\cL-\cK)_{u,\bsigma}^{\pa{2}}}_{\bx}}\prec \frac{1}{1-u}B_u^{3}(0)\mathscr{D}_{u,\bx}^{(2)}.
	\end{equation}   
	Applying the evolution kernel bound \eqref{evolution_kernel_bound_3} and integrating over $u$, these contributions are bounded by \smash{$\OO_\prec(B_r^{3}(0)\mathscr{D}_{r,\bx}^{(2)})$}. Altogether, we obtain
	\begin{equation*}
		\begin{aligned}
			\qb{Q_r\circ\E(\mathcal{L} - \mathcal{K})^{(2)}_{r, \boldsymbol{\sigma}}}_{\bx}\prec B_r^{5/2}(0)\mathscr{D}_{r,\bx}^{(2)}.
		\end{aligned}
	\end{equation*}
	Combining this with \eqref{eq:PEL-K_samesigma}, we recover \eqref{eq:EL-K_samesigma} in the case $\sigma_1\ne \sigma_2$. 
	\end{proof}

	\section{Proof of \texorpdfstring{$\cK$}{cK}-loop bounds}\label{sec_proof_of_lemma_n_cK_bound} 
	
	In this section, we present the proof of Lemma \ref{lemma_n_cK_tree_bound}. Similar estimates to \eqref{eq_n_cK_tree_bound} have been established in \cite{Band1D,Band2D,truong2025localizationlengthfinitevolumerandom,dubova2025delocalizationnonmeanfieldrandommatrices} through the analysis of a \emph{tree representation} formula for $\cK$-loops, and in \cite{erdos2025zigzagstrategyrandomband,fan2025blockreductionmethodrandom} via a dynamical approach. In the present work, we extend the dynamical method developed in our previous work \cite{fan2025blockreductionmethodrandom} to establish the bound \eqref{eq_n_cK_tree_bound}. We remark that the tree representation method can also be used to prove Lemma \ref{lemma_n_cK_tree_bound}, and it may even yield a stronger $\lesssim$-type bound instead of the $\prec$-type bound obtained here. 
	From a technical perspective, although our dynamical method is essentially equivalent to the tree representation analysis, it avoids the introduction of several intricate graphical notations as well as additional technical arguments required in our setting. In particular, it allows us to treat the power-law variance profile without assuming translation invariance of the variance profile or imposing conditions on the second-order difference of the $\Theta$-propagator, as in \cite{Band1D,Band2D,truong2025localizationlengthfinitevolumerandom,dubova2025delocalizationnonmeanfieldrandommatrices}.

	\subsection{Triangular operations on non-crossing trees}
	
	A key technical ingredient in the proof of Lemma \ref{lemma_n_cK_tree_bound} is the following bound derived from the triangle inequality:
	\begin{equation}\label{triangular_operation_bound_local}
		\begin{aligned}
			\mathscr{D}_t^2(|x-a|)\mathscr{D}_t^2(|y-a|)\lesssim \mathscr{D}_t^2(|x-y|)\mathscr{D}_t^2(|x-a|)+\mathscr{D}_t^2(|x-y|)\mathscr{D}_t^2(|y-a|).
		\end{aligned}
	\end{equation}
	Graphically, this inequality induces a \emph{switching operation} on non-crossing trees, defined as follows. Suppose the vertices  $(x_1,\ldots,x_n)$ are arranged cyclically on a circle and $T\in\txt{NCT}(\bx)$ is a non-crossing tree on $\bx$. 
	Assume that $x_j$ and $x_k$ are consecutive neighbors of $x_i$ in the cyclic order around the circle; that is, $(x_i,x_j),(x_i,x_k)\in E(T)$, and no other neighbor of $x_i$ lies on the arc between $x_j$ and $x_k$ that does not contain $x_i$. As illustrated in \Cref{tikz_triangular_operation}, the switching operation at the vertex $x_i$ with respect to the edges $\{(x_i,x_j),(x_i,x_k)\}$ generates two new non-crossing trees $T_1$ and $T_2$, obtained by replacing the edges $\{(x_i,x_j),(x_i,x_k)\}$ with $\{(x_i,x_k),(x_j,x_k)\}$ and $\{(x_i,x_j),(x_j,x_k)\}$, respectively.
	
	\begin{figure}[htbp]
		\centering
		\vspace{-0.2cm} 
		\hspace{-2.0cm}
		\begin{minipage}[c]{0.45\textwidth}
			\vspace{-0.5cm}
			\caption{Graph of the triangular bound \eqref{triangular_operation_bound_local}: }
			\label{tikz_triangular_operation}
		\end{minipage}
		\hspace{-0.2cm}
		\begin{minipage}[c]{0.55\textwidth}
			\centering
			\begin{tikzpicture}[
				scale=0.4,         fermion/.style={postaction={decorate}, decoration={markings, mark=at position 0.55 with {\arrow[scale=1.6]{Latex}}}},
				photon/.style={decorate, decoration={snake, segment length=5pt, amplitude=1.2pt}},
				dot/.style={circle, fill=black, inner sep=1.5pt}
				]
				
				\def\Cir{2.0}
				\def\angleI{65}           \def\angleK{-10}          \def\angleJ{225}  
				\begin{scope}[shift={(0,0)}]
					\draw[thick, fill=black!15] (0,0) circle (\Cir);
					
					\node[dot, label=\angleI:{$x_i$}] (xi1) at (\angleI:\Cir) {};
					\node[dot, label=\angleK:{$x_k$}] (xk1) at (\angleK:\Cir) {};
					\node[dot, label=\angleJ:{$x_j$}] (xj1) at (\angleJ:\Cir) {};
					
					\draw[thick] (xj1) -- (xi1);             \draw[thick] (xi1) -- (xk1);                    \end{scope}
				
				\node at (4.2, 0) {\Large $\lesssim$};
				
				\begin{scope}[shift={(8,0)}]
					\draw[thick, fill=black!15] (0,0) circle (\Cir);
					
					\node[dot, label=\angleI:{$x_i$}] (xi2) at (\angleI:\Cir) {};
					\node[dot, label=\angleK:{$x_k$}] (xk2) at (\angleK:\Cir) {};
					\node[dot, label=\angleJ:{$x_j$}] (xj2) at (\angleJ:\Cir) {};
					
					\draw[thick] (xi2) -- (xk2);             \draw[thick] (xj2) -- (xk2);                    \end{scope}
				
				\node at (12, 0) {\Large $+$};
				
				\begin{scope}[shift={(16,0)}]
					\draw[thick, fill=black!15] (0,0) circle (\Cir);
					
					\node[dot, label=\angleI:{$x_i$}] (xi3) at (\angleI:\Cir) {};
					\node[dot, label=\angleK:{$x_k$}] (xk3) at (\angleK:\Cir) {};
					\node[dot, label=\angleJ:{$x_j$}] (xj3) at (\angleJ:\Cir) {};
					
					\draw[thick] (xj3) -- (xi3);             \draw[thick] (xj3) -- (xk3);                    \end{scope}
				
			\end{tikzpicture}
		\end{minipage}
		\vspace{-0.2cm}
	\end{figure}
	
	The bound \eqref{triangular_operation_bound_local} then implies that
	\begin{equation*}
		\begin{aligned}
			\mathcal{D}_t(T)\lesssim \mathcal{D}_t(T_1)+\mathcal{D}_t(T_2).
		\end{aligned}
	\end{equation*}
	Clearly, each switching operation at a vertex $x_i$ reduces the degree of $x_i$ by one. Hence, the operation can be applied repeatedly at any vertex whose degree is at least $2$. Therefore, for any vertex $x\in V(T)$ with $T\in\txt{NCT}(\bx)$, repeated switching operations at $x$ yield
	\begin{equation*}
		\begin{aligned}
			\mathcal{D}_t(T)\lesssim \sum_{T'\in\txt{NCT}(\bx): \txt{deg}_{T'}(x)=1}\mathcal{D}_t(T'),
		\end{aligned}
	\end{equation*}
	where $\txt{deg}_{T'}(x)$ denotes the degree of $x$ in the tree $T'$. This observation leads to the following lemma.
	\begin{lemma}
		Fix $n\ge 2$.  For any $x\in\bx$ with $\bx\in\ZL^n$, the tree-shaped decay factor $\mathscr{T}_{t,\bx}^{(n)}$ is bounded by the contribution of trees in which $x$ has degree one, i.e., 
		\begin{equation}\label{eq_tree_decay_bounded_bound_deg_1}
			\mathscr{T}_{t,\bx}^{(n)}\lesssim \sum_{T\in\txt{NCT}(\bx): \txt{deg}_{T}(x)=1}\mathcal{D}_t(T).     \end{equation}
	\end{lemma}
	As a corollary, the tree-shaped decay factor satisfies the following convolution inequalities.
	
	\begin{lemma}\label{lemma_convolution_tree_decay}
		For any $\bx=(x_1,\ldots,x_k)\in\ZL^k$, $\by\in\ZL^n$ with fixed $k,n\geq 1$, we have 
		\begin{align}\label{eq:convolution_T_kn}
			& \sum_{a\in\ZL}\mathscr{T}_{t,(\bx,a)}^{(k+1)}\cdot\mathscr{T}_{t,(a,\by)}^{(n+1)}\lesssim \ell_t\cdot \mathscr{T}_{t,(\bx,\by)}^{(k+n)},\\
			&  \sum_{a\in\ZL}\mathscr{T}_{t,(\bx,a)}^{(k+1)}\lesssim \ell_t\cdot \mathscr{T}_{t,\bx}^{(k)}.\label{eq:convolution_T_0}
		\end{align}
	\end{lemma}
	\begin{proof}
		The bound \eqref{eq:convolution_T_0} follows immediately from \eqref{eq_tree_decay_bounded_bound_deg_1}. For \eqref{eq:convolution_T_kn}, we write
		\begin{equation*}
			\begin{aligned}
				\sum_{a\in\ZL}\mathscr{T}_{t,(\bx,a)}^{(k+1)}\cdot\mathscr{T}_{t,(a,\by)}^{(n+1)}\lesssim\sum_{a\in\ZL}\sum_{T_1\in\txt{NCT}(\bx,a):\deg_{T_1}(a)=1}\sum_{T_2\in\txt{NCT}(a,\by):\deg_{T_2}(a)=1}\mathcal{D}_{t}(T_1)\mathcal{D}_{t}(T_2).
			\end{aligned}
		\end{equation*}
		Applying the switching operation at $a$ on the RHS and then summing over $a$ yields \eqref{eq:convolution_T_kn}.
	\end{proof}
	Using the switching operation, we also obtain the following lemma, which states that for any $i\in\qq{n}$, the tree-shaped decay factor \smash{$\mathscr{T}_{t,\bx}^{(n)}$} can be bounded by the contribution of trees containing the edge $(x_i,x_{i+1})$.
	\begin{lemma}
		For any $\bx=(x_1,\ldots,x_n)\in \ZL^n$ and $i\in\qq{n}$ with $n\geq 2$, we have
		\begin{equation}\label{scrT_bounded_by_x_i_x_i_p_1}
			\mathscr{T}_{t,\bx}^{(n)}\lesssim \sum_{T\in\txt{NCT}(\bx): (x_i,x_{i+1})\in E(T)}\mathcal{D}_t(T).     \end{equation}
	\end{lemma}
	\begin{proof}
		Let $T\in\txt{NCT}(\bx)$. Suppose that $x_i$ is connected to a vertex $a$ along the geodesic path from $x_i$ to $x_{i+1}$ in $T$. If $a=x_{i+1}$, then $(x_i,x_{i+1})\in E(T)$ and the claim holds.
		Otherwise, since $T$ is non-crossing, we may apply a switching operation at $a$ to obtain two new non-crossing trees in which the graph distance between $x_i$ and $x_{i+1}$ is reduced by at least one compared with that in $T$. Iterating this procedure yields
		\begin{equation*}
			\mathcal{D}_t(T)\lesssim \sum_{T'\in\txt{NCT}(\bx): (x_i,x_{i+1})\in E(T')}\mathcal{D}_t(T'),
		\end{equation*}
		which implies \eqref{scrT_bounded_by_x_i_x_i_p_1} by summing over all $T\in\txt{NCT}(\bx)$.
	\end{proof}
	
	\subsection{Proof of Lemma \ref{lemma_n_cK_tree_bound}}
	
	We are now ready to prove the tree-shaped decay bound \eqref{eq_n_cK_tree_bound} by induction on $n$. For $n\in\{1,2\}$, the bound follows directly from Assumption \ref{assumption_input_bound}. Suppose now that $n\geq 3$ and that \eqref{eq_n_cK_tree_bound} holds for all $\cK^{\cL}$- and $\cK^{\cC}$-loops of lengths $1,\ldots, n-1$. It suffices to prove the following estimate for the $\cK\equiv \cK^{\cL}$-loops of length $n$: 
    \begin{equation}\label{target_n_cK_loop_tree_bound}
		\begin{aligned}
			\absa{\cK_{t,\bsigma,\bx}^{(n)}}\prec B_t^{n-1}(0)\cdot \mathscr{T}_{t,\bx}^{(n)}.
		\end{aligned}
	\end{equation}
    The corresponding bound for the \(\cK^{\cC}\)-loops follows from the reshaping relation \eqref{eq_reshape_cK_chains_into_loops}, the inequality $\mathscr{D}_{t}(|a-x|)\mathscr{D}_{t}(|x-b|)\lesssim \mathscr{D}_{t}(|a-b|)$, and the induction hypothesis.

	At time $t=0$, by the definition of the $\cK$-loops, for any $\bsigma=(\sigma_1,\ldots,\sigma_n)\in\ha{-,+}^n$ and $\bx=(x_1,\ldots,x_n)\in\ZL^n$, we have 
	\begin{equation}\label{eq:initialK0}
		\cK_{0,\bsigma,\bx}^{(n)}=m(\sigma_1)\cdots m(\sigma_n)\sum_{x\in\ZL}\prod_{i=1}^{n}S_{x_ix}\lesssim B_0^{n-1}(0) \mathscr{T}_{0,\bx}^{(n)},
	\end{equation}
	where the last inequality follows from a basic calculus estimate.  Now consider an arbitrary $t\in\qa{0,\tf}$. Applying Duhamel's principle to \eqref{eq_evolution_n_cK_loops} over the interval $[0,t]$, we obtain  
	\begin{equation}\label{eq:initialKt}        \cK_{t,\bsigma,\bx}^{(n)}=\pa{\bU_{0,t,\bsigma}^{(n)}\circ \cK_{0,\bsigma}^{(n)}}_{\bx}+\sum_{1\leq k<l\leq n: 2\le l-k\le n-2} \int_{0}^{t} \pa{\bU_{s,t,\bsigma}^{(n)}\circ \cA_{s,\bsigma}^{(k,l)}}_{\bx} \rd s,  
	\end{equation}
	where the tensor $\cA_{s,\bsigma}^{(k,l)}$ is defined by
	\begin{align*}
		\cA_{s,\bsigma,\ba}^{(k,l)}:=   \sum_{a\in\ZL}\pa{\pa{\cG_L}_{k,l}^{\pa{a}}\circ \cK_{s,\bsigma,\ba}^{\cL,\pa{n}}} \cdot \pa{\pa{\cG_R}_{k,l}^{\pa{a}}\circ \cK_{s,\bsigma,\ba}^{\cC,\pa{n}}},\quad \ba\in \ZL^n.
	\end{align*}
	Using the induction hypothesis together with the convolution inequality \eqref{eq:convolution_T_kn}, we obtain
	\begin{align}\label{eq:U2klt}
		\cA_{s,\bsigma,\ba}^{(k,l)} \prec \ell_s B_s^{n}(0)\mathscr{T}_{s,\ba}^{(n)} =(1-s)^{-1} \cdot B_s^{n-1}(0)\mathscr{T}_{s,\ba}^{(n)}.
	\end{align}
	
To estimate the RHS of \eqref{eq:initialKt}, we require suitable bounds on the evolution kernel $\bU_{s,t,\bsigma}^{(n)}$.  
	For any $0\leq s<t\leq \tf$ and $\bsigma=(\sigma_1,\ldots,\sigma_n)\in\ha{-,+}^n$, we decompose the evolution kernel as in \eqref{first_decomposition_cU_to_Theta_I}. 
	Moreover, in the case of alternating charge, i.e., when $\bsigma=(\sigma_1,\ldots,\sigma_n)$ satisfies $\sigma_i\neq \sigma_{i+1}$ for all $i\in\qq{n}$, we further decompose \smash{$\boldsymbol{\Theta}_{t,\bsigma}^{\qq{n},(n)}$} as in \eqref{decomposition_bTheta_to_vartheta}:
	\begin{align}
		&\boldsymbol{\Theta}_{t,\bsigma}^{\qq{n},(n)}=\vartheta_{t,\qq{n}}^{(n)}+\sum_{i=1}^{n-1}\pa{\vartheta_{t,\qq{i}}^{(n)}-\vartheta_{t,\qq{i+1}}^{(n)}}\label{second_decomposition_Theta_I}\\
		&=\vartheta_{t,\qq{n}}^{(n)}+\sum_{i=1}^{n-1}\pbb{\vartheta_{t,\qq{n}\setminus\ha{i+1}}^{(n)}-\vartheta_{t,\qq{n}}^{(n)}+\sum_{j=i+2}^{n}\pa{\vartheta_{t,\qq{i}\cup \qq{j+1,n}}^{(n)}-\vartheta_{t,\qq{i+1}\cup\qq{j+1,n}}^{(n)}-\vartheta_{t,\qq{i}\cup \qq{j,n}}^{(n)}+\vartheta_{t,\qq{i+1}\cup\qq{j,n}}^{(n)}}}.\nonumber
	\end{align}
	where we adopt the convention $\qq{n+1,n}\equiv \emptyset$. The structure of the differences in this decomposition is chosen carefully so that each summand indexed by $j$ involves two first-order differences of the $\Theta$-propagators, as in \eqref{Theta_difference_bound}. Moreover, the two ``distinguished edges'' attached to $x_{i+1}$ and $x_j$ are not ``entangled'' with the other edges attached to $x_1,\ldots, x_i$ and $x_{j+1},\ldots, x_n$, ensuring that the non-crossing structure of the trees will be preserved.
	Finally, for any $n$-tensor \smash{$\cA^{(n)}$}, we define the norm with the tree-shaped decay removed by
	\begin{equation*}
		\begin{aligned}        \normb{\cA^{(n)}}_{\mathscr{T}_t^{-1}}:=\max_{\bx\in\ZL^n}\absb{\cA^{(n)}_{\bx}/\mathscr{T}_{t,\bx}^{(n)}}.     \end{aligned}
	\end{equation*}
	The following lemma collects several estimates for the action of the evolution kernels in this tree-decay-normalized norm. Its proof is deferred to \Cref{sec_n_cK_technical_lemmas}.
	
	\begin{lemma}\label{lemma_evolution_kernel_tree_decay_bound}
		Fix $0\leq s<t\leq \tf$, and let $\cA^{(n)}$ be an arbitrary $n$-tensor. 
		\begin{enumerate}
			\item If $I\subsetneq \qq{n}$, then 
			\begin{equation}\label{evolution_kernel_tree_decay_bound_1}
				\normB{\boldsymbol{\Theta}_{t,\bsigma}^{I,(n)}\circ \cA^{(n)}}_{\mathscr{T}_t^{-1}}\prec \qa{\frac{\ell_s}{\ell_t(1-t)}}^{|I|} \normb{ \cA^{(n)}}_{\mathscr{T}_s^{-1}}.
			\end{equation}
			\item If $I=\qq{n}$ and $\bsigma$ is non-alternating, then
			\begin{equation}\label{evolution_kernel_tree_decay_bound_2}
				\quad   \norma{\boldsymbol{\Theta}_{t,\bsigma}^{\qq{n},(n)}\circ \cA^{(n)}}_{\mathscr{T}_t^{-1}}\prec \qa{\frac{\ell_s}{\ell_t(1-t)}}^{n-1} \normb{ \cA^{(n)}}_{\mathscr{T}_s^{-1}}\le \frac{1}{1-s} \qa{\frac{\ell_s}{\ell_t(1-t)}}^{n-1}\normb{ \cA^{(n)}}_{\mathscr{T}_s^{-1}}.
			\end{equation}
			\item Suppose that $\ell_t<N$ and $n\geq 4$. For any $i\in\qq{n-1}$ and $j\in\qq{i+2,n}$, we have 
			\begin{equation}\label{evolution_kernel_tree_decay_bound_3}
				\begin{aligned}
					\norma{\pa{\vartheta_{t,\qq{i}\cup \qq{j+1,n}}^{(n)}-\vartheta_{t,\qq{i+1}\cup\qq{j+1,n}}^{(n)}-\vartheta_{t,\qq{i}\cup \qq{j,n}}^{(n)}+\vartheta_{t,\qq{i+1}\cup\qq{j,n}}^{(n)}}\circ \cA^{(n)}}_{\mathscr{T}_t^{-1}}\\
					\prec \frac{1}{1-s} \qa{\frac{\ell_s}{\ell_t(1-t)}}^{n-1} \normb{ \cA^{(n)}}_{\mathscr{T}_s^{-1}}.
				\end{aligned}
			\end{equation}
			\item Suppose that $\ell_t<N$ and $n\geq 4$, and that $\cA^{(n)}$ is a $(1,i+1)$-symmetric $n$-tensor for some $i\in\qq{n-1}$ in the following sense. For any \(p\neq q\in\qq{n}\), an \(n\)-tensor \(\cA^{(n)}\) is said to be \((p,q)\)-symmetric if the induced matrix        
			\begin{equation*}
				\begin{aligned}
					(x_p,x_q)\mapsto \sum_{k\in\qq{n}\setminus \ha{p,q}}\sum_{x_k\in\ZL}\cA_{(x_1,\ldots,x_n)}^{(n)}
				\end{aligned}
			\end{equation*}  
			is symmetric. 
			Then, we have that
			\begin{equation}\label{evolution_kernel_tree_decay_bound_4}
				\norma{\pa{\vartheta_{t,\qq{n}\setminus\ha{i+1}}^{(n)}-\vartheta_{t,\qq{n}}^{(n)}}\circ \cA^{(n)}}_{\mathscr{T}_t^{-1}}\prec \frac{1}{1-s} \qa{\frac{\ell_s}{\ell_t(1-t)}}^{n-1} \normb{ \cA^{(n)}}_{\mathscr{T}_s^{-1}}.
			\end{equation}
			
			\item If $\ell_s=\ell_t=N$, then for any $i\in\qq{n-1}$,  
			\begin{equation}\label{evolution_kernel_tree_decay_bound_5}
				\norma{\pa{\vartheta_{t,\qq{i}}^{(n)}-\vartheta_{t,\qq{i+1}}^{(n)}}\circ \cA^{(n)}}_{\max}\prec \frac{1}{1-s} \pa{\frac{1}{1-t}}^{n-1}\cdot \normb{ \cA^{(n)}}_{\max}.
			\end{equation}
		\end{enumerate}
	\end{lemma}

	With these evolution kernel estimates at hand, we first control the equation \eqref{eq:initialKt} for non-alternating $\bsigma$. Combining the decomposition \eqref{first_decomposition_cU_to_Theta_I} with \eqref{evolution_kernel_tree_decay_bound_1} and \eqref{evolution_kernel_tree_decay_bound_2}, we obtain
	\[ \norma{\bU_{s,t,\bsigma}^{(n)}\circ \cA^{(n)}}_{\mathscr{T}_t^{-1}}\prec \qa{\frac{\ell_s(1-s)}{\ell_t(1-t)}}^{n-1}\normb{ \cA^{(n)}}_{\mathscr{T}_s^{-1}}.\]
	Applying \eqref{eq:initialK0} and \eqref{eq:U2klt} together with this evolution kernel estimate, and integrating over $s$, we bound \eqref{eq:initialKt} as in \eqref{target_n_cK_loop_tree_bound}.
	It therefore remains to bound \smash{$\cK_{t,\bsigma,\bx}^{(n)}$} for alternating $\bsigma$, where $n\geq 4$ is even. First consider $t$ such that $\ell_t<N$. Applying the decomposition \eqref{first_decomposition_cU_to_Theta_I} to \eqref{eq:initialKt} and using the estimates \eqref{eq:initialK0}, \eqref{eq:U2klt}, and \eqref{evolution_kernel_tree_decay_bound_1}, we obtain
	\begin{align}
		\cK_{t,\bsigma,\bx}^{(n)} &=t^n\pa{\boldsymbol{\Theta}_{t,\bsigma}^{\qq{n},(n)}\circ \cK_{0,\bsigma}^{(n)}}_{\bx} + \sum_{1\leq k<l\leq n: 2\le l-k\le n-2}\int_{0}^{t} (t-s)^n\pa{\boldsymbol{\Theta}_{t,\bsigma}^{\qq{n},(n)}\circ \cA_{s,\bsigma}^{(k,l)}}_{\bx}  \rd s +\opr{B_t^{n-1}(0)\mathscr{T}_{t,\bx}^{(n)}}\nonumber\\
		& = t^nm(\sigma_1)\cdots m(\sigma_n)\sum_{\ba=(a_1,\ldots,a_n)\in\ZL^n}\prod_{i=1}^{n}\Theta_{t,x_ia_1}\cdot \sum_{a\in\ZL}\prod_{i=1}^{n}S_{a_ia}\nonumber\\
		&\quad + \sum_{1\leq k<l\leq n: 2\le l-k\le n-2}\int_{0}^{t} (t-s)^n \sum_{\ba=(a_1,\ldots,a_n)\in\ZL^n}\prod_{i=1}^{n}\Theta_{t,x_ia_1} \cal A^{(k,l)}_{s,\bsigma,\ba} \rd s+\opr{B_t^{n-1}(0)\mathscr{T}_{t,\bx}^{(n)}}\nonumber\\
		&= \sum_{a_1\in\ZL}\prod_{i=1}^{n}\Theta_{t,x_ia_1}\cdot \cR_{t,\bsigma}^{(n)}+\opr{B_t^{n-1}(0)\mathscr{T}_{t,\bx}^{(n)}}.\label{eq:Kt_scalar}
	\end{align}
	In the second step, we additionally used the decomposition \eqref{second_decomposition_Theta_I} together with the estimates \eqref{evolution_kernel_tree_decay_bound_3} and \eqref{evolution_kernel_tree_decay_bound_4}. The symmetry required for the application of \eqref{evolution_kernel_tree_decay_bound_4} follows from the third property in Lemma \ref{lemma_basic_property_wh_cK}.
	In the third step, we introduce the scalar
	\[\cR_{t,\bsigma}^{(n)} := t^n\prod_{i=1}^{n}m(\sigma_i) +  \sum_{1\leq k<l\leq n: 2\le l-k\le n-2}\int_{0}^{t} (t-s)^n \sum_{a_2,\ldots,a_n} \cal A^{(k,l)}_{s,\bsigma,\ba} \rd s, \] 
	which is independent of $a_1$, again by the third property in Lemma \ref{lemma_basic_property_wh_cK}. 
	Next, we sum \eqref{eq:Kt_scalar} over $x_2,\ldots,x_n$. Applying Ward's identity from Lemma \ref{lemma_loops_Wards_identities} to the summation of \smash{$\cK_{t,\bsigma,\bx}^{(n)}$} over $x_n$, we control it by $(n-1)$-$\cK$-loops together with a factor $(1-t)^{-1}$. The resulting $(n-1)$-$\cK$-loops can then be bounded using the induction hypothesis. Summing over the remaining vertices and repeatedly applying \eqref{eq:convolution_T_0}, we obtain
	\begin{equation*}
		\begin{aligned}
			(1-t)^{-n}\cR_{t,\bsigma}^{(n)}\prec \sum_{x_2,\ldots,x_n}B_t^{n-1}(0)\mathscr{T}_{t,\bx}^{(n)} + (1-t)^{-1} \sum_{x_2,\ldots, x_{n-1}}B_t^{n-2}(0)\mathscr{T}_{t,(x_1,\ldots,x_{n-1})}^{(n-1)} \lesssim (1-t)^{-(n-1)} .
		\end{aligned}
	\end{equation*}
	This implies the following sum-zero property, analogous to that established in \cite[Lemma 3.10]{Band1D} for 1D regular random band matrices:
	\begin{equation*}
		\cR_{t,\bsigma}^{(n)}\prec 1-t.
	\end{equation*}
	Substituting this estimate into \eqref{eq:Kt_scalar} and using \eqref{Theta_bound_opposite_charge}, we obtain
	\begin{equation}\label{eq:alternating_sigma}
		\cK_{t,\bsigma,\bx}^{(n)}\prec (1-t)B_t^{n}(0) \sum_{a_1\in\ZL}\prod_{i=1}^{n}\mathscr{D}_{t}^2(|x_i-a_1|)+B_t^{n-1}(0)\mathscr{T}_{t,\bx}^{(n)}\prec B_t^{n-1}(0)\mathscr{T}_{t,\bx}^{(n)},
	\end{equation}
	where in the second step we used \smash{$\sum_{a\in\ZL}\prod_{i=1}^{n}\mathscr{D}_{t}^2(|x_i-a|) \le \sum_{a\in\ZL}\mathscr{T}_{t,(\bx,a)}^{(n+1)}$} together with \eqref{eq:convolution_T_0}. Combining this with the corresponding bound for the non-alternating case proves \eqref{target_n_cK_loop_tree_bound} for $t\in [0,1-(W/N)^{\al\wedge 2})$, where $\ell_t<N$. By continuity, the bound \eqref{target_n_cK_loop_tree_bound} also holds at $t_*:=1-(W/N)^{\alpha\wedge 2}$ whenever $t_*\le \tf$.

	For any $t\in\qa{t_*,\tf}$, we apply Duhamel's principle to \eqref{eq_evolution_n_cK_loops} on the interval $[t_*,t]$, which yields an equation analogous to \eqref{eq:initialKt}. The resulting terms can be controlled using the bound \eqref{target_n_cK_loop_tree_bound} at $t=t_*$ and the estimate \eqref{eq:U2klt} from the induction hypothesis, together with \eqref{evolution_kernel_tree_decay_bound_1}, \eqref{evolution_kernel_tree_decay_bound_2}, and \eqref{evolution_kernel_tree_decay_bound_5}, by arguments similar to those used in the regime $t\in[0,t_*)$. In fact, the analysis in this regime is simpler, since \smash{$\mathscr{T}_{t,\bx}^{(n)}\asymp 1$} whenever $\ell_t=N$;  hence no spatial decay factors need to be tracked. We therefore omit further details. This completes the proof of Lemma \ref{lemma_n_cK_tree_bound}.

	\subsection{Proof of \Cref{lemma_evolution_kernel_tree_decay_bound}}\label{sec_n_cK_technical_lemmas}

	For the proof of Lemma \ref{lemma_evolution_kernel_tree_decay_bound}, a key idea is to apply the switching operation (together with several analogous transformations) to reduce the problem to a number of simple cases that can be verified directly. To facilitate this reduction, we repeatedly use the following estimates arising from switching operations.
	
	\begin{lemma}\label{eq:graph_reduction}
		For any $0\leq s\leq t\leq \tf$ and $\bx=(x_1,\ldots,x_n)$ with $n\geq 2$, we have, for any $I\subsetneq \qq{n}$,
		\begin{equation}\label{eq_reduction_n_point}
			\sum_{\ba=(a_1,\ldots,a_{n})\in\ZL^n}\prod_{i\in I}\mathscr{D}_t^2(|x_i-a_i|)\cdot \prod_{i\notin I}\delta_{x_ia_i}\cdot\mathscr{T}_{s,\ba}^{(n)}\lesssim \ell_s^{|I|} \mathscr{T}_{t,\bx}^{(n)}.
		\end{equation}
		Moreover, if $I\subseteq \qq{n}\setminus\ha{1,n}$, the bound can be improved to
		\begin{equation}\label{eq_reduction_2_point}
			\sum_{\ba=(a_1,\ldots,a_{n})\in\ZL^n}\prod_{i\in I}\mathscr{D}_t^2(|x_i-a_i|)\cdot \prod_{i\notin I}\delta_{x_ia_i}\cdot\mathscr{T}_{s,\ba}^{(n)}\lesssim \ell_s^{|I|}\mathscr{D}_s^2(|x_1-x_n|)\sum_{k=1}^{n-1}\mathscr{T}_{t,(x_1,\ldots,x_k)}^{(k)}\mathscr{T}_{t,(x_{k+1},\ldots,x_n)}^{(n-k)}.
		\end{equation}
		In words, the RHS represents the sum of $\mathcal{D}_t(T)$ over all non-crossing trees $T$ on $\bx$ containing the edge $(x_1,x_n)\in E(T)$, with the factor $\mathscr{D}_t^2(|x_1-x_n|)$ replaced by $\mathscr{D}_s^2(|x_1-x_n|)$.
	\end{lemma}
	\begin{proof}
		We first prove \eqref{eq_reduction_n_point} by induction on $n$. For $n=2$, the estimate \eqref{eq_reduction_n_point} follows directly from the convolution inequality 
		\begin{equation}\label{convolution_D_t_2}
			\sum_{a\in\ZL}\mathscr{D}_s^2(|x-a|)\mathscr{D}_t^2(|a-y|)\lesssim \ell_s\mathscr{D}_t^2(|x-y|).
		\end{equation}
		Assume that $n\geq 3$ and that \eqref{eq_reduction_n_point} holds for $2,\ldots,n-1$. Without loss of generality, we may assume $1\notin I$. By \eqref{eq_tree_decay_bounded_bound_deg_1}, it suffices to show that for any $T\in\NCT(\ba)$ with $\deg_T(a_1)=1$,
		\begin{equation}\label{eq:sum_D2_DT}
			\sum_{\ba=(a_1,\ldots,a_{n})\in\ZL^n}\prod_{i\in I}\mathscr{D}_t^2(|x_i-a_i|)\cdot \prod_{i\notin I}\delta_{x_ia_i}\cdot \mathcal{D}_s(T)\lesssim \ell_s^{|I|} \mathscr{T}_{t,\bx}^{(n)}.
		\end{equation}
		Suppose that $(a_1,a_l)\in E(T)$ for some $l\in\qq{2,n}$. For ease of presentation, we only consider the case $2<l<n$ here, since the endpoint case $l\in\ha{2,n}$ follows from a similar argument. Let $T_1$ and $T_2$ denote the restrictions of $T$ to the vertex sets $(a_2,\ldots,a_l)$ and $(a_l,\ldots,a_n)$, respectively. Then, $T_1$ and $T_2$ are non-crossing trees with at most $(n-1)$ vertices. If $l\in I$, applying the induction hypothesis \eqref{eq_reduction_n_point} yields
		\begin{align}
			&~\sum_{\ba}\prod_{i\in I}\mathscr{D}_t^2(|x_i-a_i|)\cdot \prod_{i\notin I}\delta_{x_ia_i}\cdot \mathcal{D}_s(T)\nonumber\\
			=&~\sum_{\ba}\mathscr{D}_s^2(|a_1-a_l|)\cdot \prod_{i\in I}\mathscr{D}_t^2(|x_i-a_i|)\cdot \prod_{i\notin I}\delta_{x_ia_i}\cdot {\mathcal{D}_s(T_1)\mathcal{D}_s(T_2)}\nonumber\\
			\lesssim& ~\ell_s^{k(T_1)+k(T_2) -2}\sum_{a_l\in\ZL}\mathscr{D}_s^2(|x_1-a_l|)\mathscr{D}_t^2(|x_l-a_l|)\mathscr{T}_{t,(x_2,\ldots,x_{l-1},a_l)}^{(l-1)}\mathscr{T}_{t,(a_l,x_{l+1},\ldots,x_n)}^{(n-l+1)}, \label{eq:noncrossing1}        \end{align}
		where $k(T_1):=|\h{2,\ldots, l}\cap I|$ and $k(T_2):=|\h{l,\ldots, n}\cap I|$ so that $k(T_1)+k(T_2)=|I|+1$. To control the remaining sum over $a_l$, choose arbitrary trees $T_1'\in \txt{NCT}(x_2,\ldots,x_{l-1},a_l)$, $T_2'\in \txt{NCT}(a_l,x_{l+1},\ldots,x_n)$, together with the edges $(x_l,a_l)$. By repeatedly applying switching operations at $a_l$, we can bound \eqref{eq:noncrossing1} by
		\begin{align}\label{eq:noncrossing2}
			(\txt{\ref{eq:noncrossing1}})\lesssim \ell_s^{|I|-1}\sum_{a_l\in\ZL}\mathscr{D}_s^2(|x_1-a_l|)\mathscr{D}_t^2(|x_l-a_l|) \sum_{T_1',T_2'}^\star \mathcal D_t(T_1')\mathcal D_t(T_2') \lesssim \ell_s ^{|I|}\mathscr{T}_{t,\bx}^{(n)},
		\end{align}
		where we denote
		\[\sum_{T_1',T_2'}^\star \equiv \sum_{T_1'\in \txt{NCT}(x_2,\ldots, x_{l-1},a_l): \deg_{T_1'}(a_l)=1}\sum_{T_2'\in \txt{NCT}(a_l,x_{l+1},\ldots, x_{n}): \deg_{T_2'}(a_l)=1}.\] 
		In the second step of \eqref{eq:noncrossing2}, we used the following argument. Given $T_1'$ and $T_2'$, suppose that $a_l$ is connected to $x_{k_1}$ in $T_1'$ and to $x_{k_2}$ in $T_2'$, where $k_1\in \h{2,\ldots, l-1}$ and $k_2\in \h{l+1,\ldots, n}$.
		Using the estimate
		\[ \sum_{a_l}\mathscr{D}_s^2(|x_1-a_l|)\mathscr{D}_t^2(|x_l-a_l|)\mathscr{D}_t^2(|x_{k_1}-a_l|)\mathscr{D}_t^2(|x_{k_2}-a_l|) \lesssim \ell_s \mathscr{T}^{(4)}_{t,(x_1,x_{k_1},x_l,x_{k_2})},\]
		which follows from elementary calculus, and observing that $T\cup (T_1'\setminus\{(x_{k_1},a_l)\})\cup (T_2'\setminus\{(x_{k_2},a_l)\})$ forms a non-crossing tree for any $T\in\txt{NCT}(x_1,x_{k_1},x_l,x_{k_2})$, we obtain the second step of \eqref{eq:noncrossing2}. This proves \eqref{eq_reduction_n_point} in the case $l\in I$. 
		If $l\notin I$, applying the induction hypothesis and arguing similarly (in fact, more simply) yields 
		\begin{equation*}
			\begin{aligned}
				&\sum_{\ba}\prod_{i\in I}\mathscr{D}_t^2(|x_i-a_i|)\cdot \prod_{i\notin I}\delta_{x_ia_i}\cdot \mathcal{D}_s(T) \lesssim \ell_s^{k(T_1)+k(T_2)}\mathscr{D}_s^2(|x_1-x_l|)\mathscr{T}_{t,(x_2,\ldots,x_l)}^{(l-1)}\mathscr{T}_{t,(x_{l},\ldots,x_n)}^{(n-l+1)}\lesssim \ell_s^{|I|}\mathscr{T}_{t,\bx}^{(n)},
			\end{aligned}
		\end{equation*}
		where we used $k(T_1)+k(T_2)=|I|$ and $\mathscr{D}_s^2(|x_1-x_l|)\le \mathscr{D}_t^2(|x_1-x_l|)$. This completes the induction and establishes \eqref{eq_reduction_n_point}.

		The second bound \eqref{eq_reduction_2_point} follows readily from the weaker bound \eqref{eq_reduction_n_point}. Indeed, by \eqref{scrT_bounded_by_x_i_x_i_p_1}, we have
		\begin{equation}
			\begin{aligned}
				\mathscr{T}_{s,\ba}^{(n)}\lesssim \sum_{T\in\NCT(\ba):(a_1,a_n)\in E(T)} \cD_s(T)=\mathscr{D}_s^2(|a_1-a_n|)\sum_{k=1}^{n-1}\mathscr{T}_{s,(a_1,\ldots,a_k)}^{(k)}\mathscr{T}_{s,(a_{k+1},\ldots,a_n)}^{(n-k)}.
			\end{aligned}
		\end{equation}
		Substituting this bound into the LHS of \eqref{eq_reduction_2_point} and applying \eqref{eq_reduction_n_point} yields the desired estimate.
	\end{proof}
	
	Using \eqref{eq_reduction_n_point}, we can immediately derive the bound \eqref{evolution_kernel_tree_decay_bound_1}. For  $I\subsetneq\qq{n}$,  
	\begin{equation*}
		\begin{aligned}
			\pa{\bTheta_{t,\bsigma}^{I,(n)}\circ \cA^{(n)}}_{\bx}&\prec B_t^{|I|}(0)\normb{\cA^{(n)}}_{\mathscr{T}_{s}^{-1}}\sum_{\ba\in\ZL^n}\prod_{i\in I}\mathscr{D}_t^2(|x_i-a_i|)\prod_{i\notin I}\delta_{x_ia_i}\cdot \mathscr{T}_{s,\ba}^{(n)}\lesssim \qa{\frac{\ell_s}{\ell_t(1-t)}}^{|I|}\normb{\cA^{(n)}}_{\mathscr{T}_{s}^{-1}}\mathscr{T}_{t,\bx}^{(n)}.
		\end{aligned}
	\end{equation*}
	For the proof of \eqref{evolution_kernel_tree_decay_bound_2}, assume without loss of generality that $\sigma_1=\sigma_2$. Using \eqref{eq_reduction_n_point}, we obtain
	\begin{equation*}
		\begin{aligned}
			&\pa{\bTheta_{t,\bsigma}^{\qq{n},(n)}\circ \cA^{(n)}}_{\bx}\prec B_t^{n-1}(0)\normb{\cA^{(n)}}_{\mathscr{T}_{s}^{-1}}\sum_{\ba\in\ZL^n}B_0(|x_1-a_1|)\prod_{i=2}^n\mathscr{D}_t^2(|x_i-a_i|)\cdot \mathscr{T}_{s,\ba}^{(n)}\\
			&\lesssim \qa{\frac{\ell_s}{\ell_t(1-t)}}^{n-1}\normb{\cA^{(n)}}_{\mathscr{T}_{s}^{-1}}\sum_{a_1\in\ZL}B_0(|x_1-a_1|)\mathscr{T}_{t,(a_1,x_2,\ldots,x_n)}^{(n)}\lesssim \qa{\frac{\ell_s}{\ell_t(1-t)}}^{n-1}\normb{\cA^{(n)}}_{\mathscr{T}_{s}^{-1}}\mathscr{T}_{t,\bx}^{(n)}.
		\end{aligned}
	\end{equation*}
	The bound \eqref{evolution_kernel_tree_decay_bound_5} follows directly from the estimate $\norma{\Theta_t}_{\max}\prec \q{N(1-t)}^{-1}$ by \eqref{Theta_bound_opposite_charge}, together with the following estimate for any $x_1,x_2,y_1,y_2\in\ZL$:
	\begin{equation}
		\begin{aligned}
			\absa{\Theta_{t,x_1y_1}-\Theta_{t,x_2y_2}}\prec \upcirc{B_t}(|x_1-y_1|)+\upcirc{B_t}(|x_2-y_2|)\prec \frac{1}{N} \pa{\frac{N}{W}}^{\alpha\wedge 2}\le \frac{1}{N(1-s)},
		\end{aligned}
	\end{equation}
	which follows from \eqref{zero_mode_removed_bound} for  $s\in[0,\tf]$ with $\ell_t=N$. 
	It remains to prove the estimates  \eqref{evolution_kernel_tree_decay_bound_3} and \eqref{evolution_kernel_tree_decay_bound_4}. 
	For brevity, we present the argument only for $\alpha\in\qa{1,2}$, since the case $\alpha\in(2,\infty)$ follows from an analogous (and in fact simpler) argument. Without loss of generality, we assume throughout the proof that \(\norm{\cA^{(n)}}_{\mathscr{T}_s^{-1}}\le 1\). 
	
	\medskip
	\noindent{\bf Proof of \eqref{evolution_kernel_tree_decay_bound_3}.}
	To prove \eqref{evolution_kernel_tree_decay_bound_3}, we write
	\begin{align}
		&\qa{\pa{\vartheta_{t,\qq{i}\cup \qq{j+1,n}}^{(n)}-\vartheta_{t,\qq{i+1}\cup\qq{j+1,n}}^{(n)}-\vartheta_{t,\qq{i}\cup \qq{j,n}}^{(n)}+\vartheta_{t,\qq{i+1}\cup\qq{j,n}}^{(n)}}\circ \cA^{(n)}}_{\bx}\nonumber\\
		=&\sum_{\ba\in\ZL^n}\prod_{k\in\qq{i}\cup \qq{j+1,n}}\Theta_{t,x_ka_1}\cdot\prod_{k\in \qq{i+2, j-1}}\Theta_{t,x_ka_k}\cdot  (\Theta_{t,x_{i+1}a_{i+1}}-\Theta_{t,x_{i+1}a_1}) (\Theta_{t,x_ja_j}-\Theta_{t,x_ja_1})\cdot\cA^{(n)}_{\ba}.\label{eq:double_diff_ij}
	\end{align}
	Applying the bounds \eqref{Theta_bound_opposite_charge} and \eqref{Theta_difference_bound}, we obtain
	\begin{align}
		(\text{\ref{eq:double_diff_ij}})\prec B_t^n(0)\sum_{\ba \in\ZL^n}\prod_{k\in\qq{i}\cup \qq{j+1,n}}\mathscr{D}^2_{t}(|x_k-a_1|)\cdot \prod_{k\in \qq{i+2, j-1}}\mathscr{D}^2_{t}(|x_k-a_k|)\cdot  D_t(x_{i+1},x_j,a_{i+1},a_1,a_j) \mathscr{T}^{(n)}_{s,\ba},\nonumber
	\end{align}
	where we define
	\begin{equation*}
		\begin{aligned}
			&D_t(x_{i+1},x_j,a_{i+1},a_1,a_j):= {\braket{a_{i+1}-a_1}_W\braket{a_j-a_1}_W} R_t(|x_{i+1}-a_{i+1}|\vee |x_{i+1}-a_1|)R_t(|x_{j}-a_{j}|\vee |x_{j}-a_1|) \sum_{k=1}^{4}I_t^{(k)},
		\end{aligned}
	\end{equation*}
	and the terms $I_t^{(k)}\equiv I_t^{(k)}(x_{i+1},x_j,a_{i+1},a_1,a_j)$, $k\in\{1,2,3,4\}$, are given by
	\begin{equation*}
		\begin{aligned}
			&I_t^{(1)}:=\mathscr{D}^2_t(|x_{i+1}-a_{i+1}|)\mathscr{D}^2_t(|x_{j}-a_{j}|),\qquad\, I_t^{(2)}:=\mathscr{D}^2_t(|x_{i+1}-a_{i+1}|)\mathscr{D}^2_t(|x_{j}-a_{1}|),\\
			&I_t^{(3)}:=\mathscr{D}^2_t(|x_{i+1}-a_{1}|)\mathscr{D}^2_t(|x_{j}-a_{j}|),\qquad\quad I_t^{(4)}:=\mathscr{D}^2_t(|x_{i+1}-a_{1}|)\mathscr{D}^2_t(|x_{j}-a_{1}|).
		\end{aligned}
	\end{equation*}
	We repeatedly apply the switching operation to the factor $\prod_{k\in\qq{i}\cup \qq{j+1,n}}\mathscr{D}^2_{t}(|x_k-a_1|)$, which allows us to bound it by a sum of non-crossing trees with exactly one edge incident to $a_1$, say $\mathscr{D}^2_{t}(|x_{k_0}-a_1|)$. Furthermore, applying \eqref{eq:convolution_T_0}, we obtain
	\begin{equation*}
		\begin{aligned}
			\sum_{a_k\in\ZL:k\in \qq{2,i}\cup \qq{j+1,n}} \mathscr{T}^{(n)}_{s,\ba}\lesssim \ell_s^{n+i-j-1}\mathscr{T}_{s,\ba\pa{\{1\}\cup \qq{i+1,j}}}^{(j-i+1)} ,
		\end{aligned}
	\end{equation*}
	where $\ba\pa{\{1\}\cup \qq{i+1,j}}$ denotes the cycle of vertices obtained from $\ba$ by removing the vertices $a_k$ with $k\in \qq{2,i}\cup \qq{j+1,n}$ while preserving the order of the remaining elements. Consequently,
	\begin{align*}
		(\text{\ref{eq:double_diff_ij}})\prec &~ \ell_s^{n+i-j-1}B_t^n(0)\sum_{k_0\in \qq{i}\cup \qq{j+1,n}} \sum_{\substack{T\in \text{NCT}(x_1,\ldots, x_i, x_{j+1},\ldots, x_n,a_1): (x_{k_0},a_1)\in E(T)}}\sum_{a_k \in\ZL: k\in \{1\}\cup \qq{i+1,j}}  \frac{\mathcal D_t(T)}{\mathscr{D}_t^2(|x_{k_0}-a_1|)} \\
		&~\times \mathscr{D}^2_t(|x_{k_0}-a_1|) \prod_{k\in \qq{i+2, j-1}}\mathscr{D}^2_{t}(|x_k-a_k|)\cdot D_t(x_{i+1},x_j,a_{i+1},a_1,a_j)\cdot\mathscr{T}_{s,\ba\pa{\{1\}\cup \qq{i+1,j}}}^{(j-i+1)},
	\end{align*}
	Note that the factor ${\mathcal D_t(T)}/{\mathscr{D}_t^2(|x_{k_0}-a_1|)}$ does not depend on $a_1$, and therefore remains unaffected by the subsequent estimates. Hence, it suffices to consider the special case $k_0=1$, $i=1$, and $j=n$, for which we need to establish
	\begin{align*}
		\sum_{\ba \in\ZL^n} \mathscr{D}^2_{t}(|x_1-a_1|)\prod_{3\le k\le n-1}\mathscr{D}^2_{t}(|x_k-a_k|)\cdot D_t(x_{2},x_n,a_{2},a_1,a_n)\mathscr{T}^{(n)}_{s,\ba}\prec \ell_s^n\frac{\ell_t(1-t)}{\ell_s(1-s)}\cdot \mathscr{T}_{t,\bx}^{(n)}. \end{align*}
	Together with \eqref{eq_tree_decay_bounded_bound_deg_1}, it therefore suffices to show that for each $T\in\NCT(\ba)$ with $\deg_T(a_1)=1$,
	\begin{align}
		&\sum_{\ba\in\ZL^n} \mathscr{D}^2_{t}(|x_1-a_1|)\prod_{3\leq k\leq n-1}\mathscr{D}^2_{t}(|x_k-a_k|)\cdot D_t(x_{2},x_n,a_{2},a_1,a_n)\mathcal{D}_s(T)\prec \ell_s^n\frac{\ell_t(1-t)}{\ell_s(1-s)}\cdot  \mathscr{T}_{t,\bx}^{(n)}.\label{target_BBDscrD}
	\end{align}
	
	To prove \eqref{target_BBDscrD}, we first establish the following bound for the partial sums over $a_3,\ldots, a_{n-1}$:
	\begin{align}
		&\sum_{a_3,\ldots,a_{n-1}}\prod_{3\leq k\leq n-1}\mathscr{D}_{t}^2(|x_k-a_k|)\cdot \mathcal{D}_s\pa{T}\lesssim \ell_s^{n-3} \mathscr{T}_{s,(a_1,a_2,a_n)}^{(3)}\sum_{k=2}^{n-1} \mathscr{T}_{t,(a_2,x_3,\ldots,x_k)}^{(k-1)}\mathscr{T}_{t,(x_{k+1},\ldots,x_{n-1},a_n)}^{(n-k)}\nonumber\\
		&\qquad+\ell_s^{n-3} \mathscr{D}_s^2(|a_1-a_2|)\mathscr{D}_s^2(|a_1-a_n|)\sum_{2\leq i<j\le n-1}\mathscr{T}_{t,(a_2,x_3,\ldots,x_i)}^{(i-1)}\mathscr{T}_{t,(a_1,x_{i+1},\ldots,x_j)}^{(j-i+1)}\mathscr{T}_{t,(x_{j+1},\ldots,x_{n-1},a_n)}^{(n-j)}.\label{bound_a_3_to_n_1}
	\end{align}
	We justify this bound as follows. Let $a_r$ denote the unique neighbor of $a_1$ in $T$. If $r\in\ha{2,n}$, then applying \eqref{eq_reduction_2_point} directly yields the first term on the RHS of \eqref{bound_a_3_to_n_1}.
	If $r\in\qq{3,n-1}$, we apply \eqref{eq_reduction_2_point} separately to the restrictions of $T$ to the vertex sets $(a_2,\ldots,a_r)$ and $(a_r,\ldots,a_{n})$. This yields
	\begin{align}
		\sum_{a_3,\ldots,a_{n-1}}\prod_{3\leq k\leq n-1}&\mathscr{D}_{t}^2(|x_k-a_k|)\cdot \mathcal{D}_s\pa{T}\lesssim \ell_s^{n-4}\sum_{a_r}\mathscr{D}_s^2(|a_r-a_2|)\mathscr{D}_s^2(|a_1-a_r|)\mathscr{D}_s^2(|a_r-a_n|)\mathscr{D}_t^2(|x_r-a_r|)\nonumber\\
		&\times \sum_{i=2}^{r-1} \sum_{j=r}^{n-1} \mathscr{T}_{t,(a_2,x_3,\ldots,x_i)}^{(i-1)}\mathscr{T}_{t,(x_{i+1},\ldots,x_{r-1},a_r)}^{(r-i)}\mathscr{T}_{t,(a_r,x_{r+1},\ldots, x_j)}^{(j-r+1)}\mathscr{T}_{t,(x_{j+1},\ldots,x_{n-1},a_n)}^{(n-j)}.\label{bound_a_3_to_n_1_2}
	\end{align}
	Next, applying two switching operations to the three $\mathscr{D}_s^2$-edges incident to $a_r$, we bound their product by
	\begin{align*}
		\mathscr{D}_s^2(|a_1-a_2|)\mathscr{D}_s^2(|a_1-a_n|)\mathscr{D}_s^2(|a_1-a_r|) + \mathscr{T}^{(3)}_{s,(a_1,a_2,a_n)}\qa{\mathscr{D}_s^2(|a_2-a_r|)+\mathscr{D}_s^2(|a_n-a_r|)}. 
	\end{align*} 
	Substituting this bound into \eqref{bound_a_3_to_n_1_2} and estimating the resulting sum over $a_r$ using \eqref{eq_reduction_n_point}, we obtain the RHS of \eqref{bound_a_3_to_n_1}.

	With \eqref{bound_a_3_to_n_1}, to establish \eqref{target_BBDscrD} it remains to prove that
	\begin{align}
		\sum_{a_1,a_{2},a_n}&~ \mathcal{D}_{t,x_1a_1} D_t(x_{2},x_n,a_{2},a_1,a_n)\bigg(\mathscr{T}_{s,(a_1,a_2,a_n)}^{(3)}\sum_{k=2}^{n-1} \mathscr{T}_{t,(a_2,x_3,\ldots,x_k)}^{(k-1)}\mathscr{T}_{t,(x_{k+1},\ldots,x_{n-1},a_n)}^{(n-k)} \label{bound_a_3_to_n_1_3}\\
		+&~\mathcal{D}_{s,a_1a_2}\mathcal{D}_{s,a_1a_n}\sum_{2\leq i<j\le n-1}\mathscr{T}_{t,(a_2,x_3,\ldots,x_i)}^{(i-1)}\mathscr{T}_{t,(a_1,x_{i+1},\ldots,x_j)}^{(j-i+1)}\mathscr{T}_{t,(x_{j+1},\ldots,x_{n-1},a_n)}^{(n-j)}\bigg)\prec \ell_s^3\frac{\ell_t(1-t)}{\ell_s(1-s)}  \mathscr{T}_{t,\bx}^{(n)},\nonumber
	\end{align}
	where we abbreviate $\mathscr{D}^2_t(|a-b|)\equiv \mathcal D_{t,ab}$. This notation is convenient because $\mathcal D_{t,ab}=\cal D_t\p{a,b}$ in the sense of \eqref{eq:mathcalDT}, where $(a,b)$ denotes the tree consisting of the single edge $(a,b)$. We again bound each non-crossing $\mathscr{T}_t$-tree containing exactly one $a_k$ with $k\in\ha{1,2,n}$ using \eqref{eq_tree_decay_bounded_bound_deg_1}. Consequently, the estimate \eqref{bound_a_3_to_n_1_3} reduces to proving the following bounds:
	\begin{align}\label{bound_DT1T2}
		\sum_{a_1,a_{2},a_n} \mathcal{D}_{t,x_1a_1} \cdot \mathcal D_{t}(T_2)\mathcal D_{t}(T_n) \cdot D_t(x_{2},x_n,a_{2},a_1,a_n)\cdot \mathscr{T}_{s,(a_1,a_2,a_n)}^{(3)}\prec \ell_s^3\frac{\ell_t(1-t)}{\ell_s(1-s)}  \mathscr{T}_{t,\bx}^{(n)},
	\end{align}
	for any $T_2\in \text{NCT}(a_2,x_3,\ldots,x_k)$ with $\deg_{T_2}(a_2)=1$ and $T_n\in \text{NCT}(x_{k+1},\ldots,x_{n-1},a_n)$ with $\deg_{T_n}(a_n)=1$; 
	\begin{align}\label{bound_DT1T2T3}
		\sum_{a_1,a_{2},a_n} \mathcal{D}_{t,x_1a_1} \cdot \mathcal D_{s,a_1a_2}\mathcal D_{s,a_1a_n} \cdot \mathcal D_{t}(T_2)\mathcal D_{t}(T_1)\mathcal D_{t}(T_n) \cdot D_t(x_{2},x_n,a_{2},a_1,a_n)\prec \ell_s^3\frac{\ell_t(1-t)}{\ell_s(1-s)}  \mathscr{T}_{t,\bx}^{(n)},
	\end{align}
	for any $T_2\in \text{NCT}(a_2,x_3,\ldots,x_i)$ with $\deg_{T_2}(a_2)=1$, $T_n\in \text{NCT}(x_{j+1},\ldots, x_{n-1},a_n)$ with $\deg_{T_n}(a_n)=1$, and $T_1\in \text{NCT}(a_1,x_{i+1},\ldots,x_j)$ with $\deg_{T_1}(a_1)=1$. Extracting the edges in the trees $T_2,T_1,T_n$ that are incident to the vertices $a_2,a_1,a_n$ (note that if a tree consists of a single vertex, then no such edge exists), we reduce the proofs of \eqref{bound_DT1T2} and \eqref{bound_DT1T2T3} to establishing the following inequalities:
	\begin{align}
		&\sum_{a,b,c}D_t(x,z,a,b,c)\mathcal{D}_{s,ab}\mathcal{D}_{s,bc} \mathcal{D}_{t,yb}\mathcal{D}_{t,y'b}\prec \ell_s^{2+\al}\ell_t^{1-\alpha}\mathscr{T}_{t,(x,y,z,y')}^{(4)},\label{reshape_bound_1}\\
		&\sum_{a,b,c}D_t(x,z,a,b,c)\mathcal{D}_{s,ab}\mathcal{D}_{s,bc}\mathcal{D}_{t,yb}\mathcal{D}_{t,y'b}\mathcal{D}_{t,x'a}\prec \ell_s^{2+\al}\ell_t^{1-\alpha}\mathscr{T}_{t,(x,y,z,y',x')}^{(5)}, \label{reshape_bound_2}\\
		&\sum_{a,b,c}D_t(x,z,a,b,c)\mathcal{D}_{s,ab}\mathcal{D}_{s,bc}\mathcal{D}_{t,yb}\mathcal{D}_{t,y'b}\mathcal{D}_{t,z'c}\prec \ell_s^{2+\al}\ell_t^{1-\alpha}\mathscr{T}_{t,(x,y,z,z',y')}^{(5)}, \label{reshape_bound_3}\\
		&\sum_{a,b,c}D_t(x,z,a,b,c)\mathcal{D}_{s,ab}\mathcal{D}_{s,bc}\mathcal{D}_{t,yb}\mathcal{D}_{t,y'b} \mathcal{D}_{t,x'a}\mathcal{D}_{t,z'c}\prec \ell_s^{2+\al}\ell_t^{1-\alpha}\mathscr{T}_{t,(x,y,z,z',y',x')}^{(6)}, \label{reshape_bound_4}\\
		&\sum_{a,b,c}D_t(x,z,a,b,c)\mathscr{T}_{s,(a,b,c)}^{(3)}\mathcal{D}_{t,yb}\prec \ell_s^{2+\al}\ell_t^{1-\alpha}\mathscr{T}_{t,(x,y,z)}^{(3)},\label{reshape_bound_5}\\
		&\sum_{a,b,c}D_t(x,z,a,b,c)\mathscr{T}_{s,(a,b,c)}^{(3)}\mathcal{D}_{t,yb}\mathcal{D}_{t,x'a}\prec \ell_s^{2+\al}\ell_t^{1-\alpha} \mathscr{T}_{t,(x,y,z,x')}^{(4)},\label{reshape_bound_6}\\
		&\sum_{a,b,c}D_t(x,z,a,b,c)\mathscr{T}_{s,(a,b,c)}^{(3)}\mathcal{D}_{t,yb}\mathcal{D}_{t,z'c}\prec \ell_s^{2+\al}\ell_t^{1-\alpha}\mathscr{T}_{t,(x,y,z,z')}^{(4)},\label{reshape_bound_7}\\
		&\sum_{a,b,c}D_t(x,z,a,b,c)\mathscr{T}_{s,(a,b,c)}^{(3)}\mathcal{D}_{t,yb}\mathcal{D}_{t,x'a}\mathcal{D}_{t,z'c}\prec \ell_s^{2+\al}\ell_t^{1-\alpha} \mathscr{T}_{t,(x,y,z,z',x')}^{(5)}.\label{reshape_bound_8}
	\end{align}
	Here the indices $a,b,c,x,y,z$ correspond respectively to $a_2,a_1,a_n,x_2,x_1,x_n$, while the vertices adjacent to $a_2,a_1,a_n$ in the trees $T_2,T_1,T_n$ correspond to $x',y',z'$, respectively. Combining the tree-decay bounds from \eqref{reshape_bound_1}--\eqref{reshape_bound_8} with the relation \(\ell_t=W(1-t)^{-1/\alpha}\) in the regime \(\ell_t<N\), we obtain \eqref{bound_DT1T2} and \eqref{bound_DT1T2T3}, thereby completing the proof of \eqref{evolution_kernel_tree_decay_bound_3}. 
The proofs of the bounds \eqref{reshape_bound_1}--\eqref{reshape_bound_8} rely on case-by-case estimates based on straightforward but somewhat tedious calculus arguments. Moreover, the bounds \eqref{reshape_bound_1}--\eqref{reshape_bound_3} and \eqref{reshape_bound_5}--\eqref{reshape_bound_7} follow immediately from the bounds \eqref{reshape_bound_4} and \eqref{reshape_bound_8}, respectively. 
    Therefore, instead of presenting the full details, we provide a detailed proof of the most involved estimate, \eqref{reshape_bound_8}, in \Cref{sec_proof_of_reshape_bounds}, which illustrates the main ideas and several useful techniques. The proofs of the remaining estimates are similar and are therefore omitted.

	\medskip
	\noindent{\bf Proof of the bound \eqref{evolution_kernel_tree_decay_bound_4}.}   
	By the definition \eqref{eq:vartheta_tI}, we have
	\begin{align}
		&\qa{\pa{\vartheta_{t,\qq{n}\setminus\ha{i+1}}^{(n)}-\vartheta_{t,\qq{n}}^{(n)}}\circ \cA^{(n)}}_{\bx}=\sum_{\ba=(a_1,\ldots,a_n)\in\ZL^n}\prod_{k\in\qq{n}\setminus \ha{i+1}}\Theta_{t,x_ka_1}\pa{\Theta_{t,x_{i+1}a_{i+1}}-\Theta_{t,x_{i+1}a_1}}\cdot \cA^{(n)}_{\ba}\nonumber\\
		&\qquad=\sum_{a_1,a_{i+1}\in\ZL}\prod_{k\in\qq{n}\setminus \ha{i+1}}\Theta_{t,x_ka_1}\cdot \pa{\Theta_{t,x_{i+1}a_{i+1}}-\Theta_{t,x_{i+1}a_1}}\cdot\sum_{a_l\in\ZL:l\in\qq{n}\setminus\ha{1,i+1}} \cA^{(n)}_{\ba}.\label{evolution_kernel_tree_decay_bound_4pf}
	\end{align}
	The key observation in the proof of \eqref{evolution_kernel_tree_decay_bound_4} is that the $(1,i+1)$-symmetry of $\cA^{(n)}$ produces an additional difference factor. More precisely, we can rewrite
	\begin{equation*}
		\begin{aligned}
			(\text{\ref{evolution_kernel_tree_decay_bound_4pf}})&=\sum_{a_1,a_{i+1}}\prod_{k\in\qq{n}\setminus \ha{i+1}}\Theta_{t,x_ka_{i+1}}\cdot \pa{\Theta_{t,x_{i+1}a_1}-\Theta_{t,x_{i+1}a_{i+1}}}\cdot \cA^{(2)}_{a_1a_{i+1}}\\
			&=\frac{1}{2}\sum_{a_1,a_{i+1}}\pbb{\prod_{k\in\qq{n}\setminus \ha{i+1}}\Theta_{t,x_ka_1}-\prod_{k\in\qq{n}\setminus \ha{i+1}}\Theta_{t,x_ka_{i+1}}}\pa{\Theta_{t,x_{i+1}a_{i+1}}-\Theta_{t,x_{i+1}a_1}}\cdot\cA^{(2)}_{a_1a_{i+1}},
		\end{aligned}
	\end{equation*}
	where we abbreviate $\cA^{(2)}_{a_1a_{i+1}}\equiv\sum_{a_l\in\ZL:l\in\qq{n}\setminus\ha{1,i+1}} \cA^{(n)}_{\ba}$. To use the estimate \eqref{Theta_difference_bound}, we further decompose\footnote{As in \eqref{second_decomposition_Theta_I}, this decomposition is chosen so that no two edges become “entangled”, thereby preserving the non-crossing structure of the trees.}
	\begin{align*}
		&~\prod_{k\in\qq{n}\setminus \ha{i+1}}\Theta_{t,x_ka_1}-\prod_{k\in\qq{n}\setminus \ha{i+1}}\Theta_{t,x_ka_{i+1}}\\
		=&~\sum_{j=1}^{i}\pbb{\prod_{k\in\qq{n}\setminus \qq{j+1,i+1}}\Theta_{t,x_ka_1}\cdot \prod_{k=j+1}^{ i}\Theta_{t,x_ka_{i+1}}-\prod_{k\in\qq{n}\setminus \qq{j,i+1}}\Theta_{t,x_ka_1}\cdot \prod_{k=j}^{i}\Theta_{t,x_ka_{i+1}}}\\
		&~+\sum_{j=i+2}^{n}\pbb{\prod_{k=i+2}^{j}\Theta_{t,x_ka_1}\cdot \prod_{k\in\qq{n}\setminus\qq{i+1,j}}\Theta_{t,x_ka_{i+1}}-\prod_{k=i+2}^{j-1} \Theta_{t,x_ka_1}\cdot \prod_{k\in\qq{n}\setminus\qq{i+1,j-1}}\Theta_{t,x_ka_{i+1}}}\\
		=&~\sum_{j=1}^{i}\prod_{k\in\qq{n}\setminus \qq{j,i+1}}\Theta_{t,x_ka_1}\cdot \pa{\Theta_{t,x_ja_1}-\Theta_{t,x_ja_{i+1}}}\cdot \prod_{k=j+1}^i \Theta_{t,x_ka_{i+1}}\\
		&~+\sum_{j=i+2}^n \prod_{k=i+2}^{j-1}\Theta_{t,x_ka_1}\cdot \pa{\Theta_{t,x_ja_1}-\Theta_{t,x_ja_{i+1}}}\cdot \prod_{k\in\qq{n}\setminus\qq{i+1,\ldots,j}}\Theta_{t,x_ka_{i+1}}.
	\end{align*}
	Combining the above identities yields
	\begin{align}
		&\qa{\pa{\vartheta_{t,\qq{n}\setminus\ha{i+1}}^{(n)}-\vartheta_{t,\qq{n}}^{(n)}}\circ \cA^{(n)}}_{\bx}\label{1_difference_to_2}\\
		&=\frac{1}{2}\sum_{a_1,a_{i+1}}\sum_{j=1}^i  \prod_{k\in\qq{n}\setminus \qq{j,i+1}}\Theta_{t,x_ka_1}\cdot \pa{\Theta_{t,x_ja_1}-\Theta_{t,x_ja_{i+1}}}\cdot\prod_{k=j+1}^i \Theta_{t,x_ka_{i+1}}\cdot \pa{\Theta_{t,x_{i+1}a_{i+1}}-\Theta_{t,x_{i+1}a_1}}\cdot\cA^{(2)}_{a_1a_{i+1}}\nonumber\\
		&+\frac{1}{2}\sum_{a_1,a_{i+1}}\sum_{j=i+2}^n \prod_{k=i+2}^{j-1}\Theta_{t,x_ka_1}\cdot \pa{\Theta_{t,x_ja_1}-\Theta_{t,x_ja_{i+1}}}\cdot \prod_{k\in\qq{n}\setminus\qq{i+1,j}}\Theta_{t,x_ka_{i+1}}\cdot \pa{\Theta_{t,x_{i+1}a_{i+1}}-\Theta_{t,x_{i+1}a_1}}\cdot\cA^{(2)}_{a_1a_{i+1}}.\nonumber
	\end{align}
	Recalling that $\norm{\cA^{(n)}}_{\mathscr{T}_s^{-1}}\le 1$, we use \eqref{eq:convolution_T_0} to obtain
	\begin{equation*}
		\begin{aligned}
			\cA^{(2)}_{a_1a_{i+1}}\equiv \sum_{a_l\in\ZL:l\in\qq{n}\setminus\ha{1,i+1}} \cA^{(n)}_{\ba}\lesssim \ell_s^{n-2}\mathscr{D}_s^2(|a_1-a_{i+1}|).
		\end{aligned}
	\end{equation*}
	Applying this estimate together with the bounds \eqref{Theta_bound_opposite_charge} and \eqref{Theta_difference_bound}, we obtain the following bound for \eqref{1_difference_to_2}:
	\begin{equation*}
		\begin{aligned}
			&\quad \qa{\pa{\vartheta_{t,\qq{n}\setminus\ha{i+1}}^{(n)}-\vartheta_{t,\qq{n}}^{(n)}}\circ \cA^{(n)}}_{\bx}\\
			&\prec \ell_s^{n-2}B_t^n(0)\sum_{a_1,a_{i+1}}\sum_{j=1}^i \prod_{k\in\qq{n}\setminus \qq{j,i+1}}\mathscr{D}^2_{t}(|x_k-a_1|)\cdot \prod_{k=j+1}^i \mathscr{D}^2_{t}(|x_k-a_{i+1}|)\cdot  D_t(x_{i+1},x_j,a_{i+1},a_1,a_{i+1})\mathscr{D}_s^2(|a_1-a_{i+1}|)\\
			&+\ell_s^{n-2}B_t^n(0)\sum_{a_1,a_{i+1}}\sum_{j=i+2}^n \prod_{k=i+2}^{j-1}\mathscr{D}^2_{t}(|x_k-a_1|)\cdot \prod_{k\in\qq{n}\setminus\qq{i+1,j}}\mathscr{D}^2_{t}(|x_k-a_{i+1}|) \cdot D_t(x_{i+1},x_j,a_{i+1},a_1,a_{i+1})\mathscr{D}_s^2(|a_1-a_{i+1}|).
		\end{aligned}
	\end{equation*}
	Next, we apply the switching operations to the star-type factors
	\begin{equation*}
		\begin{aligned}
			\prod_{k\in\qq{n}\setminus \qq{j,i+1}}\mathscr{D}^2_{t}(|x_k-a_1|),\ \ \ \prod_{k=j+1}^i \mathscr{D}^2_{t}(|x_k-a_{i+1}|),\ \ \ \prod_{k=i+2}^{j-1}\mathscr{D}^2_{t}(|x_k-a_1|),\ \ \ \prod_{k\in\qq{n}\setminus\qq{i+1,j}}\mathscr{D}^2_{t}(|x_k-a_{i+1}|),
		\end{aligned}
	\end{equation*}
	bounding each of them by a sum of non-crossing trees, each containing at most one edge incident to $a_1$ or $a_{i+1}$. Then, using a reduction analogous to that leading to \eqref{reshape_bound_1}--\eqref{reshape_bound_8}, the proof of the bound \eqref{evolution_kernel_tree_decay_bound_4} is reduced to establishing the following inequalities:
	\begin{align}
		&\sum_{a,b}D_t(x,y,a,b,a)\mathcal{D}_{s,ab} \lesssim \ell_s^{1+\alpha}\ell_t^{1-\al}\mathscr{T}_{t,(x,y)}^{(2)},         \label{2_reshape_bound_1}\\
		&\sum_{a,b}D_t(x,y,a,b,a)\mathcal{D}_{s,ab}\mathcal{D}_{t,x'a}\lesssim \ell_s^{1+\alpha}\ell_t^{1-\al}\mathscr{T}_{t,(x',x,y)}^{(3)},\label{2_reshape_bound_2}\\
		&\sum_{a,b}D_t(x,y,a,b,a)\mathcal{D}_{s,ab} \mathcal{D}_{t,y'b} \lesssim \ell_s^{1+\alpha}\ell_t^{1-\alpha}\mathscr{T}_{t,(x,y',y)}^{(3)},\label{2_reshape_bound_3}\\
		&\sum_{a,b}D_t(x,y,a,b,a)\mathcal{D}_{s,ab} \mathcal{D}_{t,x'a} \mathcal{D}_{t,y'b}\lesssim  \ell_s^{1+\alpha}\ell_t^{1-\alpha}\mathscr{T}_{t,(x',x,y',y)}^{(4)}.\label{2_reshape_bound_4}        
	\end{align}
    The bounds \eqref{2_reshape_bound_1}--\eqref{2_reshape_bound_3} follow immediately from \eqref{2_reshape_bound_4}, which itself follows from straightforward calculus estimates and is considerably simpler than the argument for \eqref{reshape_bound_8}. We therefore omit the details. This completes the proof of \eqref{evolution_kernel_tree_decay_bound_4}.

	\subsection{Proof of the bound \eqref{reshape_bound_8}}\label{sec_proof_of_reshape_bounds}
	
	For simplicity of presentation, we introduce the abbreviation
	\[ \braket{a-b}_t:=|a-b|+\ell_t,\quad \forall t\in[0,\tf],\]
	which will be used throughout the remainder of the proof. To establish \eqref{reshape_bound_8}, we first note that, by \eqref{triangular_operation_bound_local},
	\begin{equation*}    \mathscr{T}_{s,(a,b,c)}^{(3)}\lesssim \mathcal{D}_{s,ba}\mathcal{D}_{s,ac}+\mathcal{D}_{s,bc}\mathcal{D}_{s,ac}.
	\end{equation*}
	By symmetry, it suffices to prove \eqref{reshape_bound_8} with \smash{$\mathscr{T}_{s,(a,b,c)}^{(3)}$} replaced by the second term $\mathcal{D}_{s,bc}\mathcal{D}_{s,ac}$ on the RHS. To cancel the singular factor $\braket{b-a}_W\braket{b-c}_W$ appearing in $D_t(x,z,a,b,c)$, we employ the estimate
	\begin{equation*}
		\begin{aligned}
			\mathcal{D}_{s,bc}\mathcal{D}_{s,ac} \lesssim \pa{\frac{\braket{a-b}_s}{\ell_s}}^{-1}\qa{\pa{\frac{\braket{b-c}_s}{\ell_s}}^{-\alpha}\pa{\frac{\braket{a-c}_s}{\ell_s}}^{-1-\alpha}  + \pa{\frac{\braket{b-c}_s}{\ell_s}}^{-1-\alpha}\pa{\frac{\braket{a-c}_s}{\ell_s}}^{-\alpha}}.
		\end{aligned}
	\end{equation*}
	Using this bound, we obtain
	\begin{align}
		&D_t(x,z,a,b,c)\mathcal{D}_{s,bc}\mathcal{D}_{s,ac} \lesssim  \ell_s^2 R_t(|x-a|\vee |x-b|)R_t(|z-c|\vee |z-b|)\sum_{k=1}^4\sum_{l=1}^{2}\pa{I_t^{(k)}J_s^{(l)}}(x,z,a,b,c) \nonumber\\
		&\qquad \lesssim \ell_s^2 R_t^2(0)\cdot \mathscr{R}_t(|x-a|\vee |x-b|)\mathscr{R}_t(|z-c|\vee |z-b|)\sum_{k=1}^4\sum_{l=1}^{2}\pa{I_t^{(k)}J_s^{(l)}}(x,z,a,b,c),\label{D_t_scrD_to_I_J}
	\end{align}
	where $R_t(0)=W^{\al-2}\ell_t^{1-\al}$ and $\mathscr{R}_t(r):=(r/W+1)^{\al-2}(r/\ell_t+1)^{1-\al}$ for $r\ge 0$. Furthermore, we have introduced the functions
	\begin{equation*}
		\begin{aligned}
			J_s^{(1)}(x,z,a,b,c):=&\pa{\frac{\braket{b-c}_s}{\ell_s}}^{1-\alpha}\pa{\frac{\braket{a-c}_s}{\ell_s}}^{-1-\alpha},\quad J_s^{(2)}(x,z,a,b,c):=\pa{\frac{\braket{b-c}_s}{\ell_s}}^{-\alpha}\pa{\frac{\braket{a-c}_s}{\ell_s}}^{-\alpha}.
		\end{aligned}
	\end{equation*}
	The remaining estimates are divided into eight cases according to $k=1,2,3,4$ and $l=1,2$. For brevity, we present detailed arguments only for the two representative cases $(k,l)=(1,1)$ and $(k,l)=(4,1)$; the remaining cases can be handled similarly.
	The key distinction among the different $(k,l)$ cases lies in how the two $\mathscr{R}_t$-factors are treated. In particular, there are three basic ways to bound the factor $\mathscr{R}_t(|x-a|\vee |x-b|)$:
	(1) $\mathscr R_t(|x-a|\vee |x-b|)\lesssim \mathscr R_t(|x-a|)$, (2) $\mathscr R_t(|x-a|\vee |x-b|)\lesssim \mathscr R_t(|x-b|)$, and (3) $\mathscr R_t(|x-a|\vee |x-b|)\lesssim \mathscr R_t(|a-b|)$. Each choice leads to a different graphical structure in the resulting summations. Once the appropriate bounds for the $\mathscr{R}_t$-factors are fixed, the remaining estimates follow by essentially the same arguments in all cases. Accordingly, for the cases other than $(k,l)=(1,1)$ and $(k,l)=(4,1)$, we only indicate how the $\mathscr{R}_t$-factors are controlled and omit the repetitive details.

	To obtain the factor $\ell_s^{2+\al}\ell_t^{1-\al}$ on the RHS of \eqref{reshape_bound_8}, in addition to the factor $\ell_s^2R_t^2(0)$, we need to extract an extra factor \(\ell_s\cdot (W^{2-\alpha}\ell_s^{\alpha-1})\cdot (W^{2-\al}\ell_t^{\alpha-1})\). Among the three summations over $a,b,c\in \ZL$, one summation will produce the factor $\ell_s$, while the remaining two summations will yield the other two factors via the bound
	\begin{equation}\label{eq:x-aa-y}
		\sum_{a\in\ZL}\pa{\frac{|x-a|}{W}+1}^{\alpha-2}\pa{\frac{|a-y|}{\ell}+1}^{1-\alpha}\prec W^{2-\alpha}\ell^{\alpha-1},\quad \forall W\le \ell\le N, \ x,y\in \ZL.
	\end{equation}
	Therefore, when assigning the $\mathscr{R}_t$-factors, we will deliberately arrange the sums so that expressions of the form \eqref{eq:x-aa-y} appear with $\ell\in\ha{\ell_s,\ell_t}$.
	
	\medskip
	\noindent
	{\bf Case 1:} $(k,l)=(1,1)$. In this case, we bound the $\mathscr{R}_t$-factors as $\mathscr R_t(|x-a|\vee |x-b|)\lesssim \mathscr R_t(|x-b|)$ and $\mathscr R_t(|z-c|\vee |z-b|)\lesssim \mathscr R_t(|b-c|)$, and it suffices to prove
	\begin{align}
		&\sum_{a,b,c\in\ZL}\mathcal{D}_{t,xa} \mathcal{D}_{t,x'a} \mathcal{D}_{t,zc} \mathcal{D}_{t,z'c} \mathcal{D}_{s,ac} J(x,y,b,c) \lesssim W^{4-2\alpha}\ell_s^{\alpha}\ell_t^{\alpha-1}\mathscr{T}_{t,(x,y,z,z',x')}^{(5)},\label{eq:kl11} 
	\end{align}
	where we abbreviate 
	\[J(x,y,b,c):=\mathscr R_t(|x-b|)\mathscr R_t(|b-c|) \cal D_{t,yb}\pa{{\braket{b-c}_s}/{\ell_s}}^{1-\alpha}.\] 
	Applying switching operations to the $\cal D_t$-edges at $a$ and $c$, we find that the LHS of \eqref{eq:kl11} is bounded by 
	\begin{align*}
		\mathcal{D}_{t,xx'}\mathcal{D}_{t,zz'} \sum_{a,b,c\in\ZL}\pa{\cal D_{t,xa}+\cal D_{t,x'a}}\pa{ \cal D_{t,zc} + \cal D_{t,z'c}}\mathcal{D}_{s,ac} J(x,y,b,c) .
	\end{align*}    
	By symmetry, it suffices to control the following contribution, since the remaining terms can be estimated in the same way:
	\begin{equation*}
		\begin{aligned}
			I_{1}:=\mathcal{D}_{t,xx'}\mathcal{D}_{t,zz'}\sum_{a,b,c\in\ZL}\cal D_{t,x'a} \cal D_{t,z'c} \cal D_{s,ac} J(x,y,b,c).
		\end{aligned}
	\end{equation*}

	Using the convolution inequality \eqref{convolution_D_t_2} to bound the summation over $a$, we obtain 
	\begin{align*}
		I_{1}&\lesssim \ell_s\mathcal{D}_{t,xx'}\mathcal{D}_{t,zz'}\sum_{b,c\in\ZL}\cal D_{t,x'c}  \cal D_{t,z'c} J(x,y,b,c) \\
		&\lesssim \ell_s\mathcal{D}_{t,xx'}\mathcal{D}_{t,zz'}\mathcal{D}_{t,x'z'}\sum_{b,c\in\ZL}\pa{\cal D_{t,x'c} + \cal D_{t,z'c}}J(x,y,b,c),
	\end{align*}
	where in the second step we applied another switching operation at $c$. By symmetry, it suffices to control the term containing the factor $\cal D_{t,x'c}$, which we denote by $I_2$; the remaining term can be treated in the same way. Using an argument similar to that for \eqref{triangular_convolution_example_1} and \eqref{triangular_convolution_example_2}, and decomposing the summation region according to whether $|x'-c|\ge |x'-b|/2$ or $|x'-c|< |x'-b|/2$, we obtain the following convolution inequality:
	\begin{equation}\label{three_scales_convolution_example_1}
		\begin{aligned}
			\sum_{c}\pa{\frac{\braket{x'-c}_t}{\ell_t}}^{-1-\alpha}\pa{\frac{\braket{b-c}_W}{W}}^{\alpha-2}\pa{\frac{\braket{b-c}_t}{\ell_t}}^{1-\alpha}\pa{\frac{\braket{b-c}_s}{\ell_s}}^{1-\alpha}\prec W^{2-\alpha}\ell_s^{\alpha-1}\pa{\frac{\braket{x'-b}_t}{\ell_t}}^{-\alpha}.
		\end{aligned}
	\end{equation}
	Using this inequality, we bound $I_2$ as 
	\begin{align}
		I_2&\prec W^{2-\alpha}\ell_s^{\alpha}\mathcal{D}_{t,xx'}\mathcal{D}_{t,zz'}\mathcal{D}_{t,x'z'} \sum_{b\in\ZL}\pa{\frac{\braket{x'-b}_t}{\ell_t}}^{-\alpha} \cal D_{t,yb} \mathscr{R}_t(|x-b|) .\label{eq:kl11_I2}
	\end{align}
	Next, applying the following estimate derived from the triangle inequality,
	\begin{equation}\label{eq:triangle_2form}
		\pa{\frac{\braket{x'-b}_t}{\ell_t}}^{-\alpha} \cal D_{t,yb}  \lesssim  \pa{\frac{\braket{x'-y}_t}{\ell_t}}^{-\alpha} \cal D_{t,yb} + \pa{\frac{\braket{x'-b}_t}{\ell_t}}^{-\alpha} \cal D_{t,x'y}, \end{equation}
	we further bound \eqref{eq:kl11_I2} by
	\begin{align*}
		I_2&\prec W^{2-\alpha}\ell_s^{\alpha}\mathcal{D}_{t,xx'}\mathcal{D}_{t,zz'}\mathcal{D}_{t,x'z'} \pa{\frac{\braket{x'-y}_t}{\ell_t}}^{-\alpha} \sum_{b\in\ZL} \cal D_{t,yb}\mathscr{R}_t(|x-b|)\\
		&+ W^{2-\alpha}\ell_s^{\alpha}\mathcal{D}_{t,xx'}\mathcal{D}_{t,zz'}\mathcal{D}_{t,x'z'}\mathcal{D}_{t,x'y} \sum_{b\in\ZL}\pa{\frac{\braket{x'-b}_t}{\ell_t}}^{-\alpha} \mathscr{R}_t(|x-b|).
	\end{align*}
	Using again arguments similar to those in \eqref{triangular_convolution_example_1} and \eqref{triangular_convolution_example_2}, together with the bound \eqref{eq:x-aa-y}, we obtain
	\begin{equation}\label{three_scales_convolution_example_2}
		\sum_{b\in\ZL} \cal D_{t,yb}\mathscr{R}_t(|x-b|)\prec W^{2-\al}\ell_t^{\al-1} \pa{\frac{\braket{x-y}_t}{\ell_t}}^{-1},\quad  \sum_{b\in\ZL}\pa{\frac{\braket{x'-b}_t}{\ell_t}}^{-\alpha} \mathscr{R}_t(|x-b|) \prec W^{2-\al}\ell_t^{\al-1}.
	\end{equation}
	Combining these estimates yields 
	\begin{align*}
		I_2 &\prec W^{4-2\alpha}\ell_s^{\alpha}\ell_t^{\al-1}\mathcal{D}_{t,xx'}\mathcal{D}_{t,zz'}\mathcal{D}_{t,x'z'} \qa{\pa{\frac{\braket{x'-y}_t}{\ell_t}}^{-\alpha}\pa{\frac{\braket{x-y}_t}{\ell_t}}^{-1}+\mathcal D_{t,x'y}}\\
		&\lesssim W^{4-2\alpha}\ell_s^{\alpha}\ell_t^{\al-1}\mathcal{D}_{t,xx'}\mathcal{D}_{t,zz'}\mathcal{D}_{t,x'z'} \pa{\mathcal D_{t,xy}+\mathcal D_{t,x'y}},
	\end{align*}
	where in the second step we used Young's inequality. 
	This completes the estimate for the case $(k,l)=(1,1)$.

	\medskip
	\noindent
	{\bf Case 2:} $(k,l)=(4,1)$. In this case, bounding the $\mathscr R_t$-factors as $\mathscr R_t(|x-a|\vee |x-b|)\lesssim \mathscr R_t(|x-a|)$ and $\mathscr R_t(|z-c|\vee |z-b|)\lesssim \mathscr R_t(|b-c|)$, it suffices to prove 
	\begin{align*}
		&\sum_{a,b,c}\mathcal{D}_{t,xb}\mathcal{D}_{t,yb}\mathcal{D}_{t,zb} \mathcal{D}_{t,x'a}  \mathcal{D}_{t,z'c} \mathcal{D}_{s,ac} \mathscr R_t(|x-a|)\mathscr R_t(|b-c|)\pa{\frac{\braket{b-c}_s}{\ell_s}}^{1-\alpha} \lesssim W^{4-2\alpha}\ell_s^{\alpha}\ell_t^{\alpha-1}\mathscr{T}_{t,(x,y,z,z',x')}^{(5)}. 
	\end{align*}
	The remainder of the proof is similar to the previous case $(k,l)=(1,1)$. More precisely, we first apply switching operations to the three $\cal D_t$-edges incident to the vertex $b$, replacing them by non-crossing trees on $(x,y,z)$ together with a single edge incident to $b$. As an example, we only estimate the following term, since the remaining terms can be handled analogously:
	\begin{align*}
		&I_1=\mathcal{D}_{t,xy}\mathcal{D}_{t,xz}\sum_{a,b,c\in\ZL}\mathcal{D}_{t,x'a}  \mathcal{D}_{t,z'c} \mathcal{D}_{s,ac} \mathscr R_t(|x-a|)\cdot \cal D_{t,xb}\mathscr R_t(|b-c|)\pa{\frac{\braket{b-c}_s}{\ell_s}}^{1-\alpha}.
	\end{align*}
	Bounding the summation over $b$ and applying the convolution inequality \eqref{three_scales_convolution_example_1}, we obtain
	\begin{align*}
		I_1 &\lesssim W^{2-\alpha}\ell_s^{\alpha-1}\mathcal{D}_{t,xy}\mathcal{D}_{t,xz}\sum_{a,c\in\ZL}\mathcal{D}_{t,x'a}  \mathcal{D}_{t,z'c} \mathcal{D}_{s,ac} \mathscr R_t(|x-a|)\pa{\frac{\braket{x-c}_t}{\ell_t}}^{-\alpha}\\
		&\lesssim W^{2-\alpha}\ell_s^{\alpha-1}\mathcal{D}_{t,xy}\mathcal{D}_{t,xz}\sum_{a,c\in\ZL}\mathcal{D}_{t,x'a} \mathcal{D}_{s,ac} \mathscr R_t(|x-a|)  \cdot \mathcal{D}_{t,z'c} \pa{\frac{\braket{x-z'}_t}{\ell_t}}^{-\alpha} \\
		&\quad + W^{2-\alpha}\ell_s^{\alpha-1}\mathcal{D}_{t,xy}\mathcal{D}_{t,xz}\sum_{a,c\in\ZL}\mathcal{D}_{t,x'a} \mathcal{D}_{s,ac} \mathscr R_t(|x-a|)  \cdot \mathcal{D}_{t,xz'} \pa{\frac{\braket{x-c}_t}{\ell_t}}^{-\alpha} ,
	\end{align*}
	where in the second step we applied a bound analogous to \eqref{eq:triangle_2form} to the term $\mathcal{D}_{t,z'c} \pa{{\braket{x-c}_t}/{\ell_t}}^{-\alpha}$. Next, for the summation over $c$, we apply \eqref{convolution_D_t_2} to the first term on the RHS and the following convolution inequality to the second term:
	\begin{equation*}
		\sum_{c\in\ZL}\pa{\frac{\braket{x-c}_t}{\ell_t}}^{-\alpha}\pa{\frac{\braket{a-c}_s}{\ell_s}}^{-1-\alpha}\lesssim \ell_s \pa{\frac{\braket{x-a}_t}{\ell_t}}^{-\alpha}.
	\end{equation*} 
	Consequently, we obtain
	\begin{align*}       
		I_1 & 
        \lesssim W^{2-\alpha}\ell_s^{\alpha}\mathcal{D}_{t,xy}\mathcal{D}_{t,xz}\pa{\frac{\braket{x-z'}_t}{\ell_t}}^{-\alpha}\sum_{a\in\ZL}\mathcal{D}_{t,x'a} \mathcal{D}_{t,z'a} \mathscr R_t(|x-a|) 
         \\
		&+ W^{2-\alpha}\ell_s^{\alpha}\mathcal{D}_{t,xy}\mathcal{D}_{t,xz}\mathcal{D}_{t,xz'} \sum_{a\in\ZL}\mathcal{D}_{t,x'a} \mathscr R_t(|x-a|)\pa{\frac{\braket{x-a}_t}{\ell_t}}^{-\alpha} \\
		&\lesssim W^{4-2\alpha}\ell_s^{\alpha}\ell_t^{\alpha-1}\mathcal{D}_{t,xy}\mathcal{D}_{t,xz}\ha{\mathcal D_{t,x'z'}\pa{\frac{\braket{x-z'}_t}{\ell_t}}^{-\alpha}\qa{\pa{\frac{\braket{x-x'}_t}{\ell_t}}^{-1}+\pa{\frac{\braket{x-z'}_t}{\ell_t}}^{-1}}+ \mathcal{D}_{t,xz'}\mathcal{D}_{t,xx'}}.
	\end{align*}
    For the first term, we apply the switching operation $\mathcal{D}_{t,x'a} \mathcal{D}_{t,z'a}\lesssim \mathcal{D}_{t,x'z'}\pa{\mathcal{D}_{t,x'a} + \mathcal{D}_{t,z'a}}$ together with the convolution inequality \eqref{three_scales_convolution_example_2} applied to $\sum_{a}\pa{\mathcal{D}_{t,x'a}+ \mathcal{D}_{t,z'a}} \mathscr R_t(|x-a|) $. 
    For the second term, we use the convolution estimate
	\begin{align*}
		\sum_{a}\pa{\frac{\braket{x'-a}_t}{\ell_t}}^{-1-\alpha} \pa{\frac{\abs{x-a}}{W}+1}^{\alpha-2}\pa{\frac{\braket{x-a}_t}{\ell_t}}^{1-2\alpha}\lesssim W^{2-\al}\ell_t^{\al-1}\mathcal{D}_{t,xx'},
	\end{align*}
	which follows from an argument similar to that used for \eqref{triangular_convolution_example_1} and \eqref{triangular_convolution_example_2}. Combining these bounds yields
	\begin{align*}       
		I_1 &\lesssim  W^{4-2\alpha}\ell_s^{\alpha}\ell_t^{\alpha-1}\mathcal{D}_{t,xy}\mathcal{D}_{t,xz}\pa{\mathcal D_{t,x'z'}\mathcal D_{t,xz'}+\mathcal D_{t,x'z'}\mathcal D_{t,xx'}+  \mathcal{D}_{t,xz'}\mathcal{D}_{t,xx'}},
	\end{align*}
	where we used Young's inequality in the derivation. This completes the estimate for the case $(k,l)=(4,1)$.

	\medskip
	\noindent
	{\bf Other cases.} For the cases other than $(k,l)\in\{(1,1),(4,1)\}$, we bound the $\mathscr R_t$-factors $\mathscr R_t(|x-a|\vee|x-b|)\mathscr R_t(|z-c|\vee |z-b|)$ as follows. If $k\in\ha{1,2}$, we bound them by $\mathscr R_t(|x-b|)\mathscr R_t(|b-c|)$; if $k\in\{3,4\}$, we bound them by $\mathscr R_t(|x-a|)\mathscr R_t(|b-c|)$. The remainder of the argument then proceeds in the same way as in the cases $(k,l)=(1,1)$ and $(k,l)=(4,1)$. We therefore omit the details.


\begin{thebibliography}{100}
	
	\bibitem{PRL_Anderson}
	E.~Abrahams, P.~W. Anderson, D.~C. Licciardello, and T.~V. Ramakrishnan.
	\newblock Scaling theory of localization: Absence of quantum diffusion in two
	dimensions.
	\newblock {\em Phys. Rev. Lett.}, 42:673--676, Mar 1979.
	
	\bibitem{PhysRevE.104.064202}
	R.~Agrawal, A.~Pandey, and S.~Puri.
	\newblock Enhancement in breaking of time-reversal invariance in the quantum
	kicked rotor.
	\newblock {\em Phys. Rev. E}, 104:064202, Dec 2021.
	
	\bibitem{Aizenman_book}
	M.~Aizenman and S.~Warzel.
	\newblock {\em Random Operators: Disorder Effects on Quantum Spectra and
		Dynamics}, volume 168 of {\em Graduate Studies in Mathematics}.
	\newblock American Mathematical Society, Providence, RI, 2015.
	
	\bibitem{QuadraticAOEK}
	O.~Ajanki, L.~Erd{\H{o}}s, and T.~Kr{\"u}ger.
	\newblock Quadratic vector equations on complex upper half-plane.
	\newblock {\em Memoirs of the American Mathematical Society}, 261(1261), 2019.
	
	\bibitem{ajanki2017universality}
	O.~H. Ajanki, L.~Erd{\H{o}}s, and T.~Kr{\"u}ger.
	\newblock Universality for general {Wigner}-type matrices.
	\newblock {\em Probability Theory and Related Fields}, 169(3):667--727, Dec
	2017.
	
	\bibitem{ALTSHULER1997487}
	B.~L. Altshuler and L.~S. Levitov.
	\newblock Weak chaos in a quantum {Kepler} problem.
	\newblock {\em Physics Reports}, 288(1--6):487--512, 1997.
	\newblock I.M. Lifshitz and Condensed Matter Theory.
	
	\bibitem{Anderson1958Absence}
	P.~W. Anderson.
	\newblock Absence of diffusion in certain random lattices.
	\newblock {\em Phys. Rev.}, 109:1492--1505, Mar 1958.
	
	\bibitem{BaoErd2015}
	Z.~Bao and L.~Erd{\H{o}}s.
	\newblock Delocalization for a class of random block band matrices.
	\newblock {\em Probability Theory and Related Fields}, 167(3):673--776, Apr
	2017.
	
	\bibitem{PhysRevB.98.134205}
	S.~Bera, G.~De~Tomasi, I.~M. Khaymovich, and A.~Scardicchio.
	\newblock Return probability for the {Anderson} model on the random regular
	graph.
	\newblock {\em Phys. Rev. B}, 98:134205, Oct 2018.
	
	\bibitem{bogomolny2018power}
	E.~Bogomolny and M.~Sieber.
	\newblock Power-law random banded matrices and ultrametric matrices:
	Eigenvector distribution in the intermediate regime.
	\newblock {\em Phys. Rev. E}, 98:042116, Oct 2018.
	
	\bibitem{BORGONOVI1999317}
	F.~Borgonovi, P.~Conti, D.~Rebuzzi, B.~Hu, and B.~Li.
	\newblock Cantori and dynamical localization in the {Bunimovich} stadium.
	\newblock {\em Physica D: Nonlinear Phenomena}, 131(1--4):317--326, 1999.
	\newblock Classical Chaos and its Quantum Manifestations.
	
	\bibitem{Borland1963TheNature}
	R.~E. Borland.
	\newblock The nature of the electronic states in disordered one-dimensional
	systems.
	\newblock {\em Proceedings of the Royal Society of London. Series A.
		Mathematical and Physical Sciences}, 274(1359):529--545, 1963.
	
	\bibitem{bourgade2017universality}
	P.~Bourgade, L.~Erd{\H{o}}s, H.-T. Yau, and J.~Yin.
	\newblock Universality for a class of random band matrices.
	\newblock {\em Advances in Theoretical and Mathematical Physics},
	21(3):739--800, 2017.
	
	\bibitem{bourgade2019random}
	P.~Bourgade, F.~Yang, H.-T. Yau, and J.~Yin.
	\newblock Random band matrices in the delocalized phase, {II}: Generalized
	resolvent estimates.
	\newblock {\em Journal of Statistical Physics}, 174(6):1189--1221, Mar 2019.
	
	\bibitem{bourgade2020random}
	P.~Bourgade, H.-T. Yau, and J.~Yin.
	\newblock Random band matrices in the delocalized phase {I}: Quantum unique
	ergodicity and universality.
	\newblock {\em Communications on Pure and Applied Mathematics},
	73(7):1526--1596, 2020.
	
	\bibitem{buijsman2026power}
	W.~Buijsman, M.~Haque, and I.~M. Khaymovich.
	\newblock Power-law banded random matrix ensemble as a model for quantum
	many-body {Hamiltonians}.
	\newblock {\em Phys. Rev. E}, 113:034116, Mar 2026.
	
	\bibitem{carrera2021multifractal}
	M.~Carrera-N{\'u}{\~n}ez, A.~M. Mart{\'\i}nez-Arg{\"u}ello, and J.~A.
	M{\'e}ndez-Berm{\'u}dez.
	\newblock Multifractal dimensions and statistical properties of critical
	ensembles characterized by the three classical {Wigner--Dyson} symmetry
	classes.
	\newblock {\em Physica A: Statistical Mechanics and its Applications},
	573:125965, 2021.
	
	\bibitem{ConJ-Ref2}
	G.~Casati, I.~Guarneri, F.~Izrailev, and R.~Scharf.
	\newblock Scaling behavior of localization in quantum chaos.
	\newblock {\em Phys. Rev. Lett.}, 64:5--8, Jan 1990.
	
	\bibitem{scalingabndCGMLIF1990PRL}
	G.~Casati, L.~Molinari, and F.~Izrailev.
	\newblock Scaling properties of band random matrices.
	\newblock {\em Phys. Rev. Lett.}, 64:1851--1854, Apr 1990.
	
	\bibitem{Chen2022}
	N.~Chen and C.~K. Smart.
	\newblock Random band matrix localization by scalar fluctuations.
	\newblock {\em arXiv:2206.06439}, 2022.
	
	\bibitem{CLG2023critical}
	W.~Chen, G.~Lemari\'e, and J.~Gong.
	\newblock Critical dynamics of long-range quantum disordered systems.
	\newblock {\em Phys. Rev. E}, 108:054127, Nov 2023.
	
	\bibitem{CEDS_EJP}
	G.~Cipolloni, L.~Erd{\H o}s, and D.~Schr{\"o}der.
	\newblock Optimal multi-resolvent local laws for {Wigner} matrices.
	\newblock {\em Electronic Journal of Probability}, 27(117):1--38, 2022.
	
	\bibitem{CES_Forum}
	G.~Cipolloni, L.~Erd{\H{o}}s, and D.~Schröder.
	\newblock Rank-uniform local law for {Wigner} matrices.
	\newblock {\em Forum of Mathematics, Sigma}, 10:e96, 2022.
	
	\bibitem{Cipolloni2024}
	G.~Cipolloni, R.~Peled, J.~Schenker, and J.~Shapiro.
	\newblock Dynamical localization for random band matrices up to {$W\ll
		N^{1/4}$}.
	\newblock {\em Communications in Mathematical Physics}, 405(3):82, Mar 2024.
	
	\bibitem{cohen2024complexity}
	K.~Cohen, Y.~Oz, and D.-l. Zhong.
	\newblock Complexity measure diagnostics of ergodic to many-body localization
	transition.
	\newblock {\em Phys. Rev. B}, 110:L180101, Nov 2024.
	
	\bibitem{Cuevas2001Anomalously}
	E.~Cuevas, V.~Gasparian, and M.~Ortu\~no.
	\newblock Anomalously large critical regions in power-law random matrix
	ensembles.
	\newblock {\em Phys. Rev. Lett.}, 87:056601, Jul 2001.
	
	\bibitem{cuevas2001fluctuations}
	E.~Cuevas, M.~Ortu\~no, V.~Gasparian, and A.~P\'erez-Garrido.
	\newblock Fluctuations of the correlation dimension at metal-insulator
	transitions.
	\newblock {\em Phys. Rev. Lett.}, 88:016401, Dec 2001.
	
	\bibitem{Localization1_2}
	R.~Drogin.
	\newblock Localization of one-dimensional random band matrices.
	\newblock {\em arXiv:2508.05802}, 2025.
	
	\bibitem{dubova2025delocalizationnonmeanfieldrandommatrices}
	S.~Dubova, F.~Yang, H.-T. Yau, and J.~Yin.
	\newblock Delocalization of non-mean-field random matrices in dimensions {$d\ge
		3$}.
	\newblock {\em arXiv:2507.20274}, 2025.
	
	\bibitem{DY}
	S.~Dubova and K.~Yang.
	\newblock Quantum diffusion and delocalization in one-dimensional band matrices
	via the flow method.
	\newblock {\em Journal of Mathematical Physics}, 67(3):033301, 2026.
	
	\bibitem{Band2D}
	S.~Dubova, K.~Yang, H.-T. Yau, and J.~Yin.
	\newblock Delocalization of two-dimensional random band matrices.
	\newblock {\em arXiv:2503.07606}, 2025.
	
	\bibitem{erdHos2011quantum}
	L.~Erd{\H{o}}s and A.~Knowles.
	\newblock Quantum diffusion and delocalization for band matrices with general
	distribution.
	\newblock {\em Annales Henri Poincar{\'e}}, 12(7):1227--1319, Nov 2011.
	
	\bibitem{erdHos2011quantum1}
	L.~Erd{\H{o}}s and A.~Knowles.
	\newblock Quantum diffusion and eigenfunction delocalization in a random band
	matrix model.
	\newblock {\em Communications in Mathematical Physics}, 303(2):509--554, Apr
	2011.
	
	\bibitem{Average_fluc}
	L.~Erd{\H{o}}s, A.~Knowles, and H.-T. Yau.
	\newblock Averaging fluctuations in resolvents of random band matrices.
	\newblock {\em Annales Henri Poincar{\'e}}, 14(8):1837--1926, Dec 2013.
	
	\bibitem{delocal}
	L.~Erd{\H{o}}s, A.~Knowles, H.-T. Yau, and J.~Yin.
	\newblock Delocalization and diffusion profile for random band matrices.
	\newblock {\em Communications in Mathematical Physics}, 323(1):367--416, Oct
	2013.
	
	\bibitem{EKYYEJP18}
	L.~Erd{\H{o}}s, A.~Knowles, H.-T. Yau, and J.~Yin.
	\newblock The local semicircle law for a general class of random matrices.
	\newblock {\em Electronic Journal of Probability}, 18(59):1--58, 2013.
	
	\bibitem{erdHos2024eigenstate}
	L.~Erd{\H{o}}s and V.~Riabov.
	\newblock Eigenstate thermalization hypothesis for {Wigner}-type matrices.
	\newblock {\em Communications in Mathematical Physics}, 405(12):282, Nov 2024.
	
	\bibitem{erdos2025zigzagstrategyrandomband}
	L.~Erd{\H{o}}s and V.~Riabov.
	\newblock The zigzag strategy for random band matrices.
	\newblock {\em arXiv:2506.06441}, 2025.
	
	\bibitem{ERDOS20121435}
	L.~Erd{\H{o}}s, H.-T. Yau, and J.~Yin.
	\newblock Rigidity of eigenvalues of generalized {Wigner} matrices.
	\newblock {\em Advances in Mathematics}, 229(3):1435--1515, 2012.
	
	\bibitem{evers2000fluctuations}
	F.~Evers and A.~D. Mirlin.
	\newblock Fluctuations of the inverse participation ratio at the {Anderson}
	transition.
	\newblock {\em Phys. Rev. Lett.}, 84:3690--3693, Apr 2000.
	
	\bibitem{evers2008anderson}
	F.~Evers and A.~D. Mirlin.
	\newblock {Anderson} transitions.
	\newblock {\em Rev. Mod. Phys.}, 80:1355--1417, Oct 2008.
	
	\bibitem{faez2011critical}
	S.~Faez, A.~Lagendijk, and A.~Ossipov.
	\newblock Critical scaling of polarization waves on a heterogeneous chain of
	resonators.
	\newblock {\em Phys. Rev. B}, 83:075121, Feb 2011.
	
	\bibitem{fan2025blockreductionmethodrandom}
	J.~Fan, F.~Yang, and J.~Yin.
	\newblock A block reduction method for random band matrices with general
	variance profiles.
	\newblock {\em arXiv:2507.11945}, 2025.
	
	\bibitem{PhysRevLett.66.986}
	M.~Feingold, D.~M. Leitner, and M.~Wilkinson.
	\newblock Spectral statistics in semiclassical random-matrix ensembles.
	\newblock {\em Phys. Rev. Lett.}, 66:986--989, Feb 1991.
	
	\bibitem{ScalingPropertyBandMatrixFYMA1991PRL}
	Y.~V. Fyodorov and A.~D. Mirlin.
	\newblock Scaling properties of localization in random band matrices: A
	{\ensuremath{\sigma}}-model approach.
	\newblock {\em Phys. Rev. Lett.}, 67:2405--2409, Oct 1991.
	
	\bibitem{Fyodorov_2009}
	Y.~V. Fyodorov, A.~Ossipov, and A.~Rodriguez.
	\newblock The {Anderson} localization transition and eigenfunction
	multifractality in an ensemble of ultrametric random matrices.
	\newblock {\em Journal of Statistical Mechanics: Theory and Experiment},
	2009(12):L12001, 2009.
	
	\bibitem{PhysRevE.73.026213}
	A.~M. Garc\'{\i}a-Garc\'{\i}a.
	\newblock Power spectrum characterization of the {Anderson} transition.
	\newblock {\em Phys. Rev. E}, 73:026213, Feb 2006.
	
	\bibitem{GARCIAGARCIA2004361}
	A.~M. García-García and K.~Takahashi.
	\newblock Long range disorder and {Anderson} transition in systems with chiral
	symmetry.
	\newblock {\em Nuclear Physics B}, 700(1--3):361--384, 2004.
	
	\bibitem{HeMa2018}
	Y.~He and M.~Marcozzi.
	\newblock Diffusion profile for random band matrices: A short proof.
	\newblock {\em Journal of Statistical Physics}, 177(4):666--716, Nov 2019.
	
	\bibitem{hopjan2023scale}
	M.~Hopjan and L.~Vidmar.
	\newblock Scale-invariant critical dynamics at eigenstate transitions.
	\newblock {\em Phys. Rev. Res.}, 5:043301, Dec 2023.
	
	\bibitem{ishii1973localization}
	K.~Ishii.
	\newblock Localization of eigenstates and transport phenomena in the
	one-dimensional disordered system.
	\newblock {\em Progress of Theoretical Physics Supplement}, 53:77--138, 1973.
	
	\bibitem{Kirsch2007}
	W.~Kirsch.
	\newblock An invitation to random {Schroedinger} operators.
	\newblock {\em arXiv:0709.3707}, 2007.
	
	\bibitem{PhysRevLett.79.1913}
	V.~E. Kravtsov and K.~A. Muttalib.
	\newblock New class of random matrix ensembles with multifractal eigenvectors.
	\newblock {\em Phys. Rev. Lett.}, 79:1913--1916, Sep 1997.
	
	\bibitem{Kravtsov2000Energy}
	V.~E. Kravtsov and A.~M. Tsvelik.
	\newblock Energy level dynamics in systems with weakly multifractal
	eigenstates: Equivalence to one-dimensional correlated fermions at low
	temperatures.
	\newblock {\em Phys. Rev. B}, 62:9888--9891, Oct 2000.
	
	\bibitem{L.S.Levitov_1989}
	L.~S. Levitov.
	\newblock Absence of localization of vibrational modes due to dipole-dipole
	interaction.
	\newblock {\em Europhysics Letters}, 9(1):83--86, 1989.
	
	\bibitem{PhysRevLett.64.547}
	L.~S. Levitov.
	\newblock Delocalization of vibrational modes caused by electric dipole
	interaction.
	\newblock {\em Phys. Rev. Lett.}, 64:547--550, Jan 1990.
	
	\bibitem{levitov1999critical}
	L.~S. Levitov.
	\newblock Critical {Hamiltonians} with long range hopping.
	\newblock {\em Annalen der Physik}, 511(7-9):697--706, 1999.
	
	\bibitem{liu2025edgestatisticsrandomband}
	D.-Z. Liu and G.~Zou.
	\newblock Edge statistics for random band matrices.
	\newblock {\em arXiv:2401.00492v2}, 2025.
	
	\bibitem{liu2025edgeuniversalityinhomogeneousrandom}
	D.-Z. Liu and G.~Zou.
	\newblock Edge universality for inhomogeneous random matrices.
	\newblock {\em arXiv:2508.17838}, 2025.
	
	\bibitem{lydzba2024normal}
	P.~\L{}yd\ifmmode~\dot{z}\else \.{z}\fi{}ba, R.~\ifmmode \acute{S}\else
	\'{S}\fi{}wi\ifmmode~\mbox{\k{e}}\else \k{e}\fi{}tek, M.~Mierzejewski,
	M.~Rigol, and L.~Vidmar.
	\newblock Normal weak eigenstate thermalization.
	\newblock {\em Phys. Rev. B}, 110:104202, Sep 2024.
	
	\bibitem{lydzba2020eigenstate}
	P.~\L{}yd\ifmmode~\dot{z}\else \.{z}\fi{}ba, M.~Rigol, and L.~Vidmar.
	\newblock Eigenstate entanglement entropy in random quadratic {Hamiltonians}.
	\newblock {\em Phys. Rev. Lett.}, 125:180604, Oct 2020.
	
	\bibitem{martinez2023coherent}
	M.~Martinez, G.~Lemarié, B.~Georgeot, C.~Miniatura, and O.~Giraud.
	\newblock Coherent forward scattering as a robust probe of multifractality in
	critical disordered media.
	\newblock {\em SciPost Physics}, 14(3):057, 2023.
	
	\bibitem{martinez2023scattering}
	A.~M. Mart\'{\i}nez-Arg\"uello, M.~Carrera-N\'u\~nez, and J.~A.
	M\'endez-Berm\'udez.
	\newblock Scattering and transport properties of the three classical
	{Wigner-Dyson} ensembles at the {Anderson} transition.
	\newblock {\em Phys. Rev. E}, 107:024139, Feb 2023.
	
	\bibitem{mendez2014generalized}
	J.~A. M{\'e}ndez-Berm{\'u}dez, A.~Alcazar-L{\'o}pez, and I.~Varga.
	\newblock On the generalized dimensions of multifractal eigenstates.
	\newblock {\em Journal of Statistical Mechanics: Theory and Experiment},
	2014(11):P11012, 2014.
	
	\bibitem{mendez2010scattering}
	J.~A. M\'endez-Berm\'udez, V.~A. Gopar, and I.~Varga.
	\newblock Scattering and transport statistics at the metal-insulator
	transition: A numerical study of the power-law banded random-matrix model.
	\newblock {\em Phys. Rev. B}, 82:125106, Sep 2010.
	
	\bibitem{PhysRevB.72.064108}
	J.~A. M\'endez-Berm\'udez and T.~Kottos.
	\newblock Probing the eigenfunction fractality using {Wigner} delay times.
	\newblock {\em Phys. Rev. B}, 72:064108, Aug 2005.
	
	\bibitem{mendez2006scattering}
	J.~A. M\'endez-Berm\'udez and I.~Varga.
	\newblock Scattering at the {Anderson} transition: Power-law banded random
	matrix model.
	\newblock {\em Phys. Rev. B}, 74:125114, Sep 2006.
	
	\bibitem{Mildenberger2006BoundaryMI}
	A.~Mildenberger, A.~R. Subramaniam, R.~Narayanan, F.~Evers, I.~A. Gruzberg, and
	A.~D. Mirlin.
	\newblock Boundary multifractality in critical one-dimensional systems with
	long-range hopping.
	\newblock {\em Phys. Rev. B}, 75:094204, Mar 2007.
	
	\bibitem{mirlin2000multifractality}
	A.~D. Mirlin and F.~Evers.
	\newblock Multifractality and critical fluctuations at the {Anderson}
	transition.
	\newblock {\em Phys. Rev. B}, 62:7920--7933, Sep 2000.
	
	\bibitem{Power-RBM}
	A.~D. Mirlin, Y.~V. Fyodorov, F.-M. Dittes, J.~Quezada, and T.~H. Seligman.
	\newblock Transition from localized to extended eigenstates in the ensemble of
	power-law random banded matrices.
	\newblock {\em Phys. Rev. E}, 54:3221--3230, Oct 1996.
	
	\bibitem{monthus2009statistical}
	C.~Monthus and T.~Garel.
	\newblock Statistical properties of two-particle transmission at an {Anderson}
	transition.
	\newblock {\em Journal of Physics A: Mathematical and Theoretical},
	42(47):475007, 2009.
	
	\bibitem{monthus2009statistics}
	C.~Monthus and T.~Garel.
	\newblock Statistics of the two-point transmission at {Anderson} localization
	transitions.
	\newblock {\em Phys. Rev. B}, 79:205120, May 2009.
	
	\bibitem{mott1961theory}
	N.~F. Mott and W.~D. Twose.
	\newblock The theory of impurity conduction.
	\newblock {\em Advances in Physics}, 10(38):107--163, 1961.
	
	\bibitem{ndawana2003energy}
	M.~L. Ndawana and V.~E. Kravtsov.
	\newblock Energy level statistics of a critical random matrix ensemble.
	\newblock {\em Journal of Physics A: Mathematical and General},
	36(12):3639--3645, 2003.
	
	\bibitem{paley2005statistical}
	C.~J. Paley, S.~N. Taraskin, and S.~R. Elliott.
	\newblock Statistical properties of the critical eigenstates in power-law
	random banded matrices across the band.
	\newblock {\em Phys. Rev. B}, 72:033105, Jul 2005.
	
	\bibitem{PhysRevB.57.10232}
	D.~A. Parshin and H.~R. Schober.
	\newblock Multifractal structure of eigenstates in the {Anderson} model with
	long-range off-diagonal disorder.
	\newblock {\em Phys. Rev. B}, 57:10232--10235, May 1998.
	
	\bibitem{PelSchShaSod}
	R.~Peled, J.~Schenker, M.~Shamis, and S.~Sodin.
	\newblock On the {Wegner} orbital model.
	\newblock {\em International Mathematics Research Notices}, 2019(4):1030--1058,
	2019.
	\newblock First published online 11 July 2017.
	
	\bibitem{Pono1997Coherent}
	I.~V. Ponomarev and P.~G. Silvestrov.
	\newblock Coherent propagation of interacting particles in a random potential:
	The mechanism of enhancement.
	\newblock {\em Phys. Rev. B}, 56:3742--3759, Aug 1997.
	
	\bibitem{rao2022power}
	W.-J. Rao.
	\newblock Power-law random banded matrix ensemble as the effective model for
	many-body localization transition.
	\newblock {\em The European Physical Journal Plus}, 137(3):398, Mar 2022.
	
	\bibitem{santra2025complexitytransitionschaoticquantum}
	G.~C. Santra, A.~Windey, S.~Bandyopadhyay, A.~Legramandi, and P.~Hauke.
	\newblock Complexity transitions in chaotic quantum systems: Nonstabilizerness,
	entanglement, and fractal dimension in {SYK} and random matrix models.
	\newblock {\em arXiv:2505.09707}, 2025.
	
	\bibitem{Sch2009}
	J.~Schenker.
	\newblock Eigenvector localization for random band matrices with power law band
	width.
	\newblock {\em Communications in Mathematical Physics}, 290(3):1065--1097, Sep
	2009.
	
	\bibitem{sen2025dissipativespectralformfactor}
	S.~Sen, S.~Kumar, A.~Sarkar, and M.~Kulkarni.
	\newblock Dissipative spectral form factor for the elliptic {Ginibre} unitary
	ensemble and applications.
	\newblock {\em Phys. Rev. B}, 114:014203, Jul 2026.
	
	\bibitem{Sodin2010Spectral}
	S.~Sodin.
	\newblock The spectral edge of some random band matrices.
	\newblock {\em Annals of Mathematics}, 172(3):2223--2251, 2010.
	
	\bibitem{thomson2020dynamics}
	S.~J. Thomson and M.~Schir{\'o}.
	\newblock Dynamics of disordered quantum systems using flow equations.
	\newblock {\em The European Physical Journal B}, 93(2):22, Feb 2020.
	
	\bibitem{truong2025localizationlengthfinitevolumerandom}
	S.~K. Truong, F.~Yang, and J.~Yin.
	\newblock On the localization length of finite-volume random block
	{Schr\"odinger} operators.
	\newblock {\em arXiv:2503.11382}, 2025.
	
	\bibitem{vallejo2024reducing}
	I.~Vallejo-Fabila, A.~K. Das, D.~A. Zarate-Herrada, A.~S. Matsoukas-Roubeas,
	E.~J. Torres-Herrera, and L.~F. Santos.
	\newblock Reducing dynamical fluctuations and enforcing self-averaging by
	opening many-body quantum systems.
	\newblock {\em Phys. Rev. B}, 110:075138, Aug 2024.
	
	\bibitem{varga2002fluctuation}
	I.~Varga.
	\newblock Fluctuation of correlation dimension and inverse participation number
	at the {Anderson} transition.
	\newblock {\em Phys. Rev. B}, 66:094201, Sep 2002.
	
	\bibitem{varga2000critical}
	I.~Varga and D.~Braun.
	\newblock Critical statistics in a power-law random-banded matrix ensemble.
	\newblock {\em Phys. Rev. B}, 61:R11859(R)--R11862(R), May 2000.
	
	\bibitem{varga2008entanglement}
	I.~Varga and J.~A. Méndez-Bermúdez.
	\newblock Entanglement in disordered systems at criticality.
	\newblock {\em physica status solidi (c)}, 5(3):867--870, 2008.
	
	\bibitem{vega2019multifractality}
	D.~A. Vega-Oliveros, J.~A. M\'endez-Berm\'udez, and F.~A. Rodrigues.
	\newblock Multifractality in random networks with power-law decaying bond
	strengths.
	\newblock {\em Phys. Rev. E}, 99:042303, Apr 2019.
	
	\bibitem{10.1214/18-EJP197}
	P.~von Soosten and S.~Warzel.
	\newblock The phase transition in the ultrametric ensemble and local stability
	of {Dyson Brownian} motion.
	\newblock {\em Electronic Journal of Probability}, 23(70):1--24, 2018.
	
	\bibitem{vonSoostenWarzel2019}
	P.~von Soosten and S.~Warzel.
	\newblock Delocalization and continuous spectrum for ultrametric random
	operators.
	\newblock {\em Annales Henri Poincar{\'e}}, 20(9):2877--2898, Sep 2019.
	
	\bibitem{Sooster2019}
	P.~von Soosten and S.~Warzel.
	\newblock Non-ergodic delocalization in the {Rosenzweig--Porter} model.
	\newblock {\em Letters in Mathematical Physics}, 109(4):905--922, Apr 2019.
	
	\bibitem{10.1214/19-ECP278}
	P.~von Soosten and S.~Warzel.
	\newblock Random characteristics for {Wigner} matrices.
	\newblock {\em Electronic Communications in Probability}, 24(75):1--12, 2019.
	
	\bibitem{Wigner}
	E.~P. Wigner.
	\newblock Characteristic vectors of bordered matrices with infinite dimensions.
	\newblock {\em Annals of Mathematics}, 62(3):548--564, 1955.
	
	\bibitem{MWilkinson_1991}
	M.~Wilkinson, M.~Feingold, and D.~M. Leitner.
	\newblock Localization and spectral statistics in a banded random matrix
	ensemble.
	\newblock {\em Journal of Physics A: Mathematical and General}, 24(1):175--182,
	1991.
	
	\bibitem{xu2024bulk}
	C.~Xu, F.~Yang, H.-T. Yau, and J.~Yin.
	\newblock Bulk universality and quantum unique ergodicity for random band
	matrices in high dimensions.
	\newblock {\em The Annals of Probability}, 52(3):765--837, 2024.
	
	\bibitem{YYYTexpansion2022CMP}
	F.~Yang, H.-T. Yau, and J.~Yin.
	\newblock Delocalization and quantum diffusion of random band matrices in high
	dimensions {II}: {T}-expansion.
	\newblock {\em Communications in Mathematical Physics}, 396(2):527--622, Dec
	2022.
	
	\bibitem{yang2021delocalization}
	F.~Yang, H.-T. Yau, and J.~Yin.
	\newblock Delocalization and quantum diffusion of random band matrices in high
	dimensions {I}: Self-energy renormalization.
	\newblock {\em Communications in Mathematical Physics}, 407(7):151, Jun 2026.
	
	\bibitem{Band1D_III}
	F.~Yang and J.~Yin.
	\newblock Random band matrices in the delocalized phase, {III}: averaging
	fluctuations.
	\newblock {\em Probability Theory and Related Fields}, 179(1):451--540, Feb
	2021.
	
	\bibitem{yang2025delocalizationgeneralclassrandom}
	F.~Yang and J.~Yin.
	\newblock Delocalization of a general class of random block {Schr\"odinger}
	operators.
	\newblock {\em arXiv:2501.08608}, 2025.
	
	\bibitem{Bandedge}
	F.~Yang and J.~Yin.
	\newblock Delocalization of random band matrices at the edge.
	\newblock {\em arXiv:2505.11993}, 2025.
	
	\bibitem{Band1D}
	H.-T. Yau and J.~Yin.
	\newblock Delocalization of one-dimensional random band matrices.
	\newblock {\em arXiv:2501.01718}, 2025.
	
	\bibitem{C.Yeung_1987}
	C.~Yeung and Y.~Oono.
	\newblock A conjecture on nonlocal random tight-binding models.
	\newblock {\em Europhysics Letters}, 4(9):1061--1065, 1987.
	
	\bibitem{yevtushenko2003virial}
	O.~Yevtushenko and V.~E. Kravtsov.
	\newblock Virial expansion for almost diagonal random matrices.
	\newblock {\em Journal of Physics A: Mathematical and General},
	36(30):8265--8289, 2003.
	
	\bibitem{zarateherrada2025dynamicaldetectionextendednonergodic}
	D.~A. Zarate-Herrada, I.~Vallejo-Fabila, L.~F. Santos, and E.~J.
	Torres-Herrera.
	\newblock Dynamical detection of extended nonergodic states in many-body
	quantum systems.
	\newblock {\em arXiv:2505.23910}, 2025.
	
	\bibitem{zhang2024magnetic}
	X.~Zhang and M.~S. Foster.
	\newblock Magnetic instability and spin-glass order beyond the {Anderson-Mott}
	transition in interacting power-law random banded matrix fermions.
	\newblock {\em Phys. Rev. B}, 110:155137, Oct 2024.
	
\end{thebibliography}
\end{document}